\documentclass[11pt,reqno]{amsart}

\pdfoutput=1

\usepackage{LPPSS1_Macros}

\title[Asymptotic stability of the degree-one vortex in the abelian YMH model]{Asymptotic stability of the degree-one vortex in the abelian Yang-Mills-Higgs model}

\author[J. L\"uhrmann]{Jonas L\"uhrmann}
\address{Department of Mathematics and Computer Science, University of Cologne, Cologne, Germany}
\email{jonas.luehrmann@uni-koeln.de}

\author[J. M. Palacios]{Jos\'e M. Palacios}
\address{Institute of Mathematics, \'Ecole Polytechnique F\'ed\'erale de Lausanne (EPFL), Lausanne, Switzerland}
\email{jose.palaciosarmesto@epfl.ch}

\author[F. Pusateri]{Fabio Pusateri}
\address{Department of Mathematics, University of Toronto, Toronto, Ontario, Canada}
\email{fabiop@mail.math.toronto.edu}

\author[W. Schlag]{Wilhelm Schlag}
\address{Department of Mathematics, Yale University, New Haven, Connecticut, USA}
\email{wilhelm.schlag@yale.edu}

\author[S. Shahshahani]{Sohrab Shahshahani}
\address{Department of Mathematics and Statistics, University of Massachusetts Amherst, Amherst, Massachusetts, USA}
\email{sshahshahani@umass.edu}

\thanks{
J.L.\  was partially supported by NSF CAREER grant DMS-2235233.
F.P.\ is supported in part by a start-up grant from the University of Toronto and NSERC grants RGPIN-2018-0648 and RGPIN-2025-06419. 
J.M.P.\ was partially supported by NSERC grants RGPIN-2018-0648, and is also supported by the Swiss National Science Foundation grant 225701. 
W.S.\ is partially supported by the United States NSF DMS-2350356. S.S. was supported by the Simons Foundation grant 639284.
}

\begin{document}

\begin{abstract}
We prove asymptotic stability of the degree-one vortex in the
$(1+2)$-dimensional abelian Yang-Mills-Higgs model at the self-dual
coupling, for small equivariant perturbations in weighted Sobolev spaces.
The abelian Yang-Mills-Higgs model is a classical relativistic field
theory on $(1+2)$-dimensional Minkowski space, describing a complex-valued
field coupled to an electromagnetic potential and admitting topological
solitons known as vortices.
  
This paper is the final and main part of a three-paper series. Under the orthogonal gauge condition used here, 
perturbations of the vortex are governed by a system of nonlinear Klein-Gordon equations for the dynamical variables, coupled to
an elliptic equation for the temporal component of the electromagnetic potential. The linearized operator has continuous spectrum $[1,\infty)$
and a single positive gap eigenvalue (internal mode) of multiplicity two, whose spectral properties, associated distorted Fourier theory, and linear decay estimates are developed in the two companion papers.

The main difficulty is the long-time analysis of the coupled radiation--internal-mode dynamics. In two space dimensions the dispersive
decay of the Klein-Gordon radiation is relatively weak, while the internal mode decays only on the long time scale dictated by nonlinear radiation damping.
At the same time, the Klein-Gordon equations for the radiation contain non-spatially localized variable coefficient quadratic interactions,
which cannot be treated perturbatively and require a normal form analysis.
We prove decay of the radiation by combining a good-bad decomposition, 
a flat-sharp decomposition, and a space-time resonance analysis carried out
relative to the flat Klein-Gordon flow. The passage between the flat analysis 
and the Klein-Gordon flow with potential is achieved through integrated local energy decay and transference estimates 
derived from the distorted Fourier theory. We expect this part of the nonlinear analysis to be of independent interest. 
\end{abstract}

\maketitle 
  
\setcounter{tocdepth}{2}

\tableofcontents

\section{Introduction}

This paper is the third and final part of a three-part series, cf. \cites{LPPSS1}{LPPSS2}, devoted to the asymptotic stability problem for equivariant perturbations of the degree-one vortex in the $(1+2)$-dimensional abelian Yang–Mills–Higgs model at the critical (self–dual) coupling.

\subsection{The abelian Yang-Mills-Higgs model}

The abelian Yang-Mills-Higgs model on $(1+2)$-dimensional Minkowski space is a classical relativistic field theory for a complex-valued field 
\begin{equation}
    \phi \colon \bbR^{1+2} \to \bbC 
\end{equation}
coupled to a real-valued connection $1$-form
\begin{equation}
    A = A_0 \, \ud t + A_1 \, \ud x_1 + A_2 \, \ud x_2, \quad A_\mu \colon \bbR^{1+2} \to \bbR, \quad 0 \leq \mu \leq 2.
\end{equation}
For a coupling constant $\lambda > 0$, the Lagrangian action functional is 
\begin{equation}
    \calL_\lambda\bigl[ \phi, A \bigr] := \int_{\bbR^{1+2}} \biggl( \frac14 F_{\mu \nu} F^{\mu \nu} + \frac12 \bfD_\mu \phi \overline{\bfD^\mu \phi} + \frac{\lambda}{8} \bigl( 1 - |\phi|^2 \bigr)^2 \biggr) \, \ud x \, \ud t.
\end{equation}
Here
\begin{equation}
    \begin{aligned}
        \bfD_\mu := \nabla_\mu - i A_\mu, \qquad F_{\mu \nu} := \nabla_\mu A_\nu - \nabla_\nu A_\mu, \quad 0 \leq \mu, \nu \leq 2,
    \end{aligned}
\end{equation}
denote the covariant derivative and the curvature of the connection.
Throughout, we use the Einstein summation convention with respect to the Minkowski metric $\mathrm{diag}[-1,1,1]$ on $\bbR^{1+2}$.

The associated Euler-Lagrange equations are 
\begin{equation} \label{equ:EL_equations}
    \left\{ \begin{aligned}
        \bfD^\mu \bfD_\mu \phi + \frac{\lambda}{2} \bigl( 1 - |\phi|^2 \bigr) \phi &= 0, \\ 
        \nabla^\mu F_{\mu \nu} + \Im \bigl( \overline{\phi} \bfD_\nu \phi \bigr) &= 0,
    \end{aligned} \right.
\end{equation}
and will be referred to as the (hyperbolic) abelian Yang-Mills-Higgs equations.
Defining the magnetic field $B$ and the electric field components $E_1, E_2$ by
\begin{equation}
    B := F_{12}, \qquad E_1 := F_{01}, \qquad E_2 := F_{02},
\end{equation}
we can write the equations \eqref{equ:EL_equations} in rectangular coordinates as
\begin{equation} \label{equ:EL_equations_rectangular} 
    \left\{ \begin{aligned}
        -\bfD_0^2 \phi + \bfD^j \bfD_j \phi + \frac{\lambda}{2} \bigl( 1 - |\phi|^2 \bigr) \phi &= 0, \\ 
        \partial_0 E_1 + \partial_2 B - \Im \bigl( \overline{\phi} \bfD_1 \phi \bigr) &= 0, \\
        \partial_0 E_2 - \partial_1 B - \Im \bigl( \overline{\phi} \bfD_2 \phi \bigr) &= 0, \\
        \partial_1 E_1 + \partial_2 E_2 - \Im \bigl( \overline{\phi} \bfD_0 \phi \bigr) &= 0.
    \end{aligned} \right.
\end{equation}
The system is invariant under the gauge transformations
\begin{equation} \label{equ:gauge_invariance}
    (\phi, A) \mapsto \bigl( e^{i\gamma} \phi, A + \ud \gamma \bigr)
\end{equation}
for any sufficiently regular function $\gamma \colon \bbR^{1+2} \to \bbR$. It is also invariant under space-time translations and under Lorentz transformations.
The (formally conserved) energy functional is given by
\begin{equation} \label{equ:energy_functional_general_lambda}
    \mathcal{E}_\lambda\bigl[\phi, A\bigr] := \frac12 \int_{\bbR^2} \biggl( B^2 + E_1^2+E_2^2 + |\bfD_0\phi|^2 + |\bfD_1\phi|^2 + |\bfD_2\phi|^2 + \frac{\lambda}{4} \bigl(1-|\phi|^2\bigr)^2 \biggr) \, \ud x.
\end{equation}
Global existence of sufficiently regular finite energy solutions to \eqref{equ:EL_equations} in the Lorenz gauge was established in \cite{BurMon1}.

The Higgs potential $\bbC \ni \phi \mapsto \frac{\lambda}{8} (1-|\phi|^2)^2$ has vacuum manifold $\bbS^1$.
Correspondingly, finite-energy configurations approach the vacuum at spatial infinity in the sense that, formally,
\begin{equation} \label{equ:intro_finite_energy_behavior_at_infinity}
    |\phi| \to 1, \quad |\bfD \phi| \to 0 \quad \text{as } |x| \to \infty.
\end{equation} 
Thus, the scalar field determines, at least formally, a map from the circle at infinity $\bbS^1_\infty$ to the vacuum manifold $\bbS^1$.
One can therefore associate to such a configuration an integer-valued winding number, or topological degree, $\mathrm{deg}(\phi) \in \bbZ$.
Under suitable convergence assumptions at infinity, the magnetic flux is quantized according to
\begin{equation} \label{equ:intro_magnetic_flux}
    \frac{1}{2\pi} \int_{\bbR^2} B \, \ud x = \mathrm{deg}(\phi) \in \bbZ.
\end{equation}
We refer to \cite[Chapter II]{bookJT}, \cite[Section 3]{Zhao26}, and the references therein for a rigorous discussion of the topological degree of finite energy abelian Yang-Mills-Higgs fields.

Static finite energy configurations in fixed topological degree classes are expected to play a central role in the long-time dynamics of the abelian Yang-Mills-Higgs equations \eqref{equ:EL_equations}.
The vacuum states $(\phi, A) = (e^{i\kappa},0)$, $\kappa \in [0,2\pi)$, have trivial topological degree; their asymptotic stability was addressed in \cite{Tsutsumi03}.
In every non-trivial topological degree class $n \in \bbZ \backslash \{0\}$ there exist vortex solutions, which in Coulomb gauge take the form 
\begin{equation} \label{equ:intro_vortex_solutions_general_lambda}
    (\phi_n, A_n) = \bigl( e^{in\theta} U_n(r), a_n(r) \, \ud \theta \bigr),
\end{equation}
where the radial profiles satisfy
\begin{equation}
    U_n(0) = a_n(0) = 0, \quad U_n(r) \to 1, \quad a_n(r) \to n \quad \text{as } r \to \infty.
\end{equation}
The existence of such vortices was proved in \cite{Plo1,BeCh}.
\cite{GusSigVort,GusVort} established that vortices of degree $\pm 1$ are (orbitally) stable for all values of the coupling constant $\lambda > 0$, while higher-degree vortices $|n| \geq 2$ are orbitally stable for $0 < \lambda < 1$ and unstable for $\lambda > 1$. 
The asymptotic stability problem for the degree-one vortex has remained open.

The main result of this three-paper series is the asymptotic stability of the degree-one vortex under the abelian Yang-Mills-Higgs equations \eqref{equ:EL_equations} in the self-dual case $\lambda = 1$ for equivariant perturbations that are small in a weighted Sobolev norm, see Theorem~\ref{thm:main} in Subsection~\ref{subsec:main_result} below for the precise statement of the result.
Throughout the remainder of this paper we restrict to the self-dual coupling $\lambda = 1$.

\subsection{The self-dual case}

For the self-dual coupling $\lambda = 1$, the Euler-Lagrange equations become 
\begin{equation} \label{equ:EL_equations_selfdual}
    \left\{ \begin{aligned}
        \bfD^\mu \bfD_\mu \phi + \frac12 \bigl( 1 - |\phi|^2 \bigr) \phi &= 0, \\ 
        \nabla^\mu F_{\mu \nu} + \Im \bigl( \overline{\phi} \bfD_\nu \phi \bigr) &= 0.
    \end{aligned} \right.
\end{equation}
This value of the coupling constant is distinguished by the Bogomolny structure of the static energy functional \cite{Bogo}.
Indeed, by completing the square and using an integration by parts identity, one obtains 
\begin{equation}\label{eq:staticenergyBogomolny1}
    \begin{aligned}
        E_{\mathrm{static}}\bigl[\phi, A \bigr] &:= \int_{\bbR^2} \biggl( \frac12 B^2 + \frac12 |\bfD_1\phi|^2 + \frac12 |\bfD_2\phi|^2 + \frac18 \bigl( 1 - |\phi|^2 \bigr)^2 \biggr) \, \ud x \\ 
        &= \frac12 \int_{\bbR^2} \biggl( \Bigl( B \mp \frac12 \bigl( 1 - |\phi|^2 \bigr) \Bigr)^2 + |(\bfD_1 \pm i \bfD_2) \phi|^2 \biggr) \, \ud x \pm \frac12 \int_{\bbR^2} B \, \ud x. 
    \end{aligned}
\end{equation}
The choice of sign corresponds to the sign of the topological degree.
In view of the magnetic flux quantization \eqref{equ:intro_magnetic_flux}, the topological degree $n \in \bbZ$ of a finite energy field thus provides a lower bound on the static energy $E_{\mathrm{static}}\bigl[\phi,A\bigr] \geq \pi |n|$.
Energy minimizers in a fixed degree class must therefore solve the first-order Bogomolny equations. For positive degree these are
\begin{equation} \label{equ:intro_1st_order_Bogomolny}
    B = \frac12 \bigl( 1 - |\phi|^2 \bigr), \quad (\bfD_1 + i \bfD_2) \phi = 0.
\end{equation}
The moduli space of solutions to \eqref{equ:intro_1st_order_Bogomolny} was thoroughly studied in \cite{Taubes80_1,Taubes80_2}, see also \cite[Chapter III]{bookJT}. In particular, all static solutions to \eqref{equ:EL_equations_selfdual} satisfy \eqref{equ:intro_1st_order_Bogomolny}.

In this work we focus on the static equivariant degree-one vortex centered at the origin
\begin{equation} \label{equ:intro_1vortex_selfdual}
    \underline{\Phi}(r,\theta) = e^{i\theta} U(r), \quad \underline{A}(r,\theta) = a_\theta(r) \, \ud \theta.
\end{equation}
The finite energy condition implies
\begin{equation}
    U(r) \to 1, \quad a_\theta(r) \to 1 \quad \text{as } r \to \infty.
\end{equation}
We denote by $(U, a_\theta) \in \calC^\infty(\bbR_+) \times \calC^\infty(\bbR_+)$ the unique solution to the corresponding system of ordinary differential equations
\begin{equation} \label{equ:intro_onevortex_ODEs}
    \left\{ \begin{aligned}
        \frac{1}{r} \partial_r a_\theta - \frac12 \bigl( 1 - U^2 \bigr) &= 0, \\
        \partial_r U - \frac{1-a_\theta}{r} U &= 0,
    \end{aligned} \right.
\end{equation}
that is, the Bogomolny equations in polar coordinates,
subject to the boundary conditions 
\begin{equation} \label{equ:intro_onevortex_ODEs_boundary_conditions}
    \begin{aligned}
        U(0) = 0, \quad \lim_{r \to \infty} U(r) = 1, \quad U'(r) > 0, \\ 
        a_\theta(0) = 0, \quad \lim_{r\to\infty} a_\theta(r) = 1, \quad a_\theta'(r) > 0.
    \end{aligned}
\end{equation}
We refer to \cite[Proposition 5.1]{LPPSS1} for a simple ODE proof of uniqueness of the system~\eqref{equ:intro_onevortex_ODEs} under the conditions \eqref{equ:intro_onevortex_ODEs_boundary_conditions}. 
In the sequel, we refer both to the pair $(\underline{\Phi}, \underline{A})$ and to the radial profiles $(U, a_\theta)$ as the degree-one vortex.

\begin{figure}[ht]
    \centering
    \includegraphics[scale=0.95]{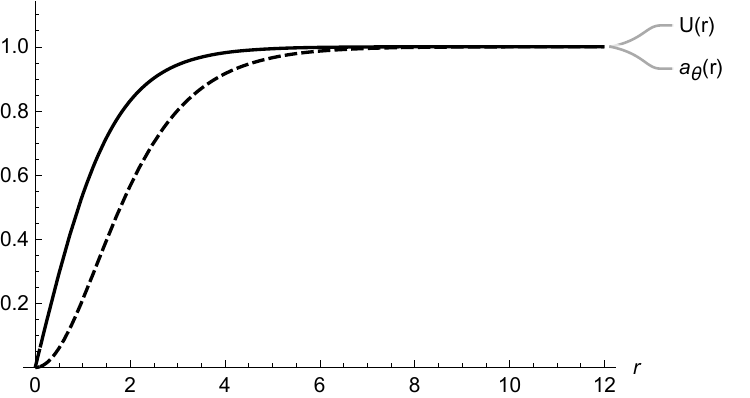}
    \caption{Plots of the radial profiles $U(r)$ and $a_\theta(r)$ of the degree-one vortex.}
    \label{figvortex}
\end{figure}

\medskip
\subsection{Hamiltonian formulation}\label{subsec:Hamiltonintro}

We recast the self-dual abelian Yang-Mills-Higgs equations \eqref{equ:EL_equations_selfdual} in Hamiltonian form. 
This formulation will be used to derive the evolution equations for equivariant perturbations of the degree-one vortex and will also be important in parts of the nonlinear analysis later on. 
As in Maxwell theory, the temporal component of the connection acts as a Lagrange multiplier enforcing Gauss’ law; see, for instance, \cite{DiracQMlec1,waldgrbook1}, and also \cite{Ruback1988} for a related derivation in a different gauge.

Recall the Lagrangian action functional
\begin{equation}
    \calL\bigl[ \phi, A \bigr] := \int_{\bbR^{1+2}} \mathscrbf{L}(\phi;A) \, \ud x \, \ud t
\end{equation}
with Lagrangian density
\begin{equation}
    \mathscrbf{L}(\phi;A) := \dfrac{1}{4} F_{\mu\nu}F^{\mu \nu} + \dfrac{1}{2}\bfD_\mu \phi \overline{\bfD^\mu \phi} + \dfrac{1}{8}\big( 1 - \vert \phi\vert^2 \big)^2.
\end{equation}
The associated Euler-Lagrange equations \eqref{equ:EL_equations_selfdual} read in rectangular coordinates 
\begin{equation} \label{equ:EL_equations_selfdual_rectangular}
    \left\{ \begin{aligned}
        -\bfD_0^2 \phi + \bfD^j \bfD_j \phi + \frac12 \bigl( 1 - |\phi|^2 \bigr) \phi &= 0, \\ 
        \partial_0 E_1 + \partial_2 B - \Im \bigl( \overline{\phi} \bfD_1 \phi \bigr) &= 0, \\
        \partial_0 E_2 - \partial_1 B - \Im \bigl( \overline{\phi} \bfD_2 \phi \bigr) &= 0, \\
        \partial_1 E_1 + \partial_2 E_2 - \Im \bigl( \overline{\phi} \bfD_0 \phi \bigr) &= 0.
    \end{aligned} \right.
\end{equation}

We define the position variables $\pos = (\pos_\phi, \overline{\pos}_\phi, \pos_{A_1},\pos_{A_2})$ by
\begin{equation} \label{equ:intro_def_pos_phi}
    \pos_\phi := \phi, \qquad \overline{\pos}_\phi := \overline{\phi}, \qquad \pos_{A_j} := A_j, \qquad j = 1, 2,
\end{equation}
and the momentum variables $\mom = (\mom_\phi,\overline{\mom}_\phi, \mom_{A_1},\mom_{A_2})$ by
\begin{equation} \label{equ:intro_def_mom_phi}
    \mom_\phi := \dfrac{\partial \mathscrbf{L}}{ \partial( \partial_0\overline{\pos}_\phi)} = - \dfrac{1}{2} \bfD_0 \phi, \quad \overline{\mom}_\phi := \dfrac{\partial \mathscrbf{L}}{ \partial( \partial_0\pos_\phi)} = - \dfrac{1}{2} \overline{\bfD_0 \phi}, \quad \mom_{A_j} := \dfrac{\partial \mathscrbf{L}}{\partial(\partial_0 \pos_{A_j})} = -E_j, \quad j = 1,2.
\end{equation}
Since the Lagrangian does not contain $\partial_0 A_0$, the temporal component $A_0$ has no conjugate momentum variable and is retained as a Lagrange multiplier.
The Hamiltonian is defined by
\begin{equation}\label{equ:introH}
    \mathcal{H} := \int_{\R^2} \mathscrbf{H} \, \ud x
\end{equation}
with Hamiltonian density 
\begin{equation} \label{equ:introHdensity}
    \begin{aligned}
        \mathscrbf{H} &:= \mathscrbf{L}\bigl(\pos_\phi, \pos_A; \partial_0\pos_\phi, \partial_0\overline{\pos}_\phi, \partial_0 \pos_A, A_0\bigr) - \mom_\phi \partial_0 \overline{\pos}_\phi - \overline{\mom}_\phi \partial_0 \pos_\phi - \sum_{j=1}^2 \mom_{A_j}\partial_0\pos_{A_j} \\ 
        &= \dfrac{1}{2} B^2 + \dfrac{1}{2}\vert \mom_A\vert^2 + 2\vert \mom_\phi\vert^2  + \dfrac{1}{2}\big\vert (\nabla_x - i \pos_A)  \pos_\phi \big\vert^2 
        + \dfrac{1}{8} \big(1-\vert \pos_\phi\vert^2 \big)^2 - 2 A_0 \im (\overline{\pos}_\phi \mom_\phi) - \mom_A \cdot \nabla_x A_0.
    \end{aligned}
\end{equation}
In passing to the second line, we eliminated the velocities in favor of the momenta $\partial_0 \phi = - 2 \mom_\phi + i A_0 \pos_\phi$ and $\partial_0 A_j = \partial_j A_0 - \mom_{A_j}$, while keeping $A_0$ as an independent Lagrange multiplier.

Then the Euler-Lagrange equations \eqref{equ:EL_equations_selfdual_rectangular} are equivalent to the constrained Hamiltonian system 
\begin{align} \label{eq:con_ham_1}
    \partial_t \pos = - \dfrac{\delta \mathcal{H}}{\delta \overline{\mom}}, \qquad \partial_t \mom = \dfrac{\delta \mathcal{H}}{\delta \overline{\pos}}, \qquad \dfrac{\delta \mathcal{H}}{\delta A_0}=0
\end{align}
for the conjugate variables $\pos = (\pos_\phi,\overline{\pos}_\phi, \pos_{A_1},\pos_{A_2})$ and $\mom = (\mom_\phi,\overline{\mom}_\phi, \mom_{A_1},\mom_{A_2})$.
Variation with respect to $A_0$ yields precisely Gauss’ law, namely the final equation in~\eqref{equ:EL_equations_selfdual_rectangular}.

\subsection{Equivariant perturbations of the degree-one vortex}

We now consider equivariant perturbations of the degree-one vortex
\begin{equation} \label{eq:hamper}
    \begin{aligned}
        \phi(t,r,\theta) & =  \underline{\Phi}(r,\theta)+\varphi(t,r,\theta) = e^{i\theta} \big( U(r) + \alpha(t,r) + i\beta(t,r) \big), \\
        A(t,r,\theta) & =  \underline{A}(r,\theta) + \eta(t,r,\theta) = a_\theta(r) \, \ud \theta  + \eta_\theta(t,r) \, \ud \theta + \eta_r(t,r) \, \ud r + \eta_0(t,r) \, \ud t.
    \end{aligned}
\end{equation}
Using the constrained Hamiltonian formulation \eqref{eq:con_ham_1} of the self-dual abelian Yang-Mills-Higgs system of equations, we derive the evolution equations for the perturbation variables. To this end we define the Hamiltonian relative to the vortex by 
\begin{equation} \label{eq:Hamper}
\begin{aligned}
&\mathcal{H}_0\bigl[\varphi,\overline{\varphi},\eta,\mom_\varphi,\overline{\mom}_\varphi,\mom_\eta,\eta_0\bigr] \\
    &\quad := \mathcal{H}\big[\underline{\Phi} + \varphi, \overline{\underline{\Phi}} + \overline{\varphi}, \underline{A}+\eta,\mom_{\varphi},\overline{\mom}_{\varphi},\mom_{\eta},\eta_0 \big] - \mathcal{H}\big[ \underline{\Phi},\overline{\underline{\Phi}},\underline{A},0,0,0,0 \big],
\end{aligned}
\end{equation}
where $\mom_\varphi := \mom_\phi$ and $\mom_\eta := \mom_A$ as in \eqref{equ:intro_def_mom_phi}.
For the spatial components of the perturbation connection, we introduce the notation 
\begin{align} \label{eq:hamperzetamu}
    \zeta := \frac{\eta_\theta}{r}, \qquad \mu := - \eta_r.
\end{align}
Then, if we define the positions
\begin{equation}\label{eq:con_ham_4.0}
    \pos_\alpha = 2\alpha, \qquad \pos_\beta = 2\beta,\qquad \pos_\zeta=\zeta,\qquad  \pos_\mu=\mu,
\end{equation} 
and momenta
\begin{equation}\label{eq:con_ham_4_alt0}
   \mom_\alpha   = - \partial_t\alpha- \eta_0 \beta ,\qquad \mom_\beta  = - \partial_t\beta+ \eta_0(U+\alpha), \qquad  \mom_\zeta = -\partial_t \zeta \qquad \mom_\mu = - \partial_t \mu  - \partial_r \eta_0,
\end{equation}
we obtain the (non-canonical) Hamiltonian system of equations
\begin{equation} \label{eq:con_ham_4_alt}
\begin{aligned}
    \partial_t \pos_\ast &= -2\frac{\delta\calH_0}{\delta\mom_\ast}, \qquad \partial_t\mom_\ast = 2\frac{\delta\calH_0}{\delta \pos_\ast}, \qquad \ast\in\{\alpha,\beta\}, \\
    \partial_t \pos_\ast &= -\frac{\delta\calH_0}{\delta\mom_\ast}, \qquad \, \, \, \partial_t\mom_\ast = \frac{\delta\calH_0}{\delta \pos_\ast}, \qquad \, \, \, \ast\in\{\zeta,\mu\},
\end{aligned}
\end{equation}
along with the constraint equation
\begin{equation} \label{eq:con_ham_4_alt_constraint}
    \frac{\delta \calH_0}{\delta \eta_0} = 0.
\end{equation}
Here, we view $\calH_0$ as a function of $(\pos_\ast,\mom_\ast,\eta_0)$ with $\ast\in\{\alpha,\beta,\zeta,\mu\}$.
The derivation of \eqref{eq:con_ham_4_alt}--\eqref{eq:con_ham_4_alt_constraint} together with an explicit expansion of $\calH_0$ is given Subsection~\ref{Ssec:Ham}.

\medskip 

Next, we impose the orthogonal gauge condition (or background gauge condition) on the perturbation
\begin{equation} \label{equ:intro_orthogonal_gauge_condition}
    \nabla^j \eta_j = \Im \bigl( \overline{\underline{\Phi}} \phi \bigr).
\end{equation}
Since the vortex connection $a_\theta(r) \, \ud \theta$ is divergence-free, this is equivalent to $\nabla^j A_j = \Im \bigl( \overline{\underline{\Phi}} \phi \bigr)$. We also refer to this condition as Stuart’s gauge \cite{Stu}. 
Using the above notation, it reads in polar coordinates 
\begin{equation} \label{equ:intro_Stuart_gauge_polar}
    \partial^\ast \mu = U\beta, \qquad \partial^\ast := -\partial_r-\frac{1}{r}.
\end{equation}
Starting from the constrained Hamiltonian equations \eqref{eq:con_ham_4_alt}, \eqref{eq:con_ham_4_alt_constraint}, imposing the gauge condition~\eqref{equ:intro_Stuart_gauge_polar}, and expanding $\calH_0$, see the calculations in Lemma~\ref{lemHamQ}, we obtain that the dynamical variables $\alpha, \zeta, \beta, \mu$ satisfy the system of equations 
\begin{equation} \label{equ:sys3}
    \begin{aligned}
        \pt^2 \begin{pmatrix} \alpha \\ \zeta \\ \beta \\ \mu \end{pmatrix} 
        + 
        \begin{pmatrix}
            L_1 & -2U' & 0 & 0 \\
            -2U' & L_2 & 0 & 0 \\ 
            0 & 0 & L_1 & -2U' \\
            0 & 0 & -2U' & L_2 
        \end{pmatrix}
        \begin{pmatrix} \alpha \\ \zeta \\ \beta \\ \mu \end{pmatrix} 
        = 
        \begin{pmatrix}
            \calN_\alpha \\ \calN_\zeta \\ \calN_\beta \\ \calN_\mu 
        \end{pmatrix},
    \end{aligned}
\end{equation}
coupled to the elliptic equation for the temporal component
\begin{equation} \label{equ:eta0_equation}
    \big(-\Delta + U^2\big) \eta_0 = \calN_0.
\end{equation}
Here we use the notation 
\begin{align} \label{eq:bL1L2}
    L_1 := -\Delta + \frac{(1-a_\theta)^2}{r^2} - \frac{\pr a_\theta}{r} + U^2, \quad L_2 := -\Delta + \frac{1}{r^2} + U^2. 
\end{align}
Moreover, introducing for convenience
\begin{equation}
    b := \frac{1-a_\theta}{r},
\end{equation}
the nonlinearities are given by
\begin{equation} \label{sysN1}
    \begin{aligned}
        \calN_\alpha &:= - \frac32 U \alpha^2 + \frac12 U \beta^2 - \frac12 (\alpha^2+\beta^2)\alpha - 2 \mu \pr\beta + 2 b \zeta \alpha - U (\zeta^2 + \mu^2) - (\zeta^2 + \mu^2) \alpha \\
        &\quad \quad - 2\eta_0 \pt \beta - (\pt \eta_0) \beta + U \eta_0^2 + \eta_0^2 \alpha, \\
        \calN_\zeta &:= b (\alpha^2 + \beta^2) - 2 U \zeta \alpha - \zeta (\alpha^2 +\beta^2), \\   
        \calN_\beta &:= - 2U \alpha\beta - \frac12 (\alpha^2+\beta^2)\beta + 2\mu \pr \alpha + 2 b \zeta \beta - (\zeta^2+\mu^2) \beta \\
        &\quad \quad + 2\eta_0 \pt \alpha + U (\pt \eta_0) + (\pt \eta_0) \alpha + \eta_0^2 \beta, \\ 
        \calN_\mu &:= -\pr \pt \eta_0 - \alpha \pr \beta + \beta \pr \alpha - 2U \alpha \mu - (\alpha^2+\beta^2) \mu,
    \end{aligned}
\end{equation}
and
\begin{equation}
    \calN_0 := -\beta \pt \alpha + \alpha \pt \beta - 2 U \eta_0 \alpha - \eta_0 (\alpha^2 + \beta^2).
\end{equation}
For later estimates, we also require an elliptic equation for $\partial_t \eta_0$.
Differentiating \eqref{equ:eta0_equation} in time and inserting the evolution equations \eqref{equ:sys3} gives
\begin{equation} \label{equ:pt_eta0_equation}
    \begin{aligned}
        \big(-\Delta + U^2\big) \pt \eta_0 
        &= \beta L_1\alpha - 2U' \beta \zeta - \alpha L_1 \beta + 2U' \alpha \mu + \calC_0
    \end{aligned}
\end{equation}
with 
\begin{equation}
    \begin{aligned}
        \calC_0 &:= - \beta \calN_\alpha + \alpha \calN_\beta + \pt \bigl( - 2U\eta_0\alpha - \eta_0 (\alpha^2 + \beta^2) \bigr).
    \end{aligned}
\end{equation}

\subsection{Spectral decomposition}

For the linearized operator in \eqref{equ:sys3} we introduce the short-hand notation
\begin{equation}\label{sysN1compact}
    \bfM := \begin{pmatrix} \bfL & 0 \\ 0 & \bfL \end{pmatrix}, 
    \quad 
    \bfL := \begin{pmatrix} L_1 & -2U' \\ -2U' & L_2 \end{pmatrix}.
\end{equation}
Then the system \eqref{equ:sys3} can be written more succinctly as
\begin{equation} \label{syscompact}
    \bigl( \pt^2 + \bfM \bigr) \bmr = \bmN
\end{equation}
with 
\begin{equation} \label{syscompact2}
    \bmr := \bigl(\alpha, \zeta, \beta, \mu\bigr)^T, \quad \bmN := \bigl(\calN_\alpha, \calN_\zeta, \calN_\beta, \calN_\mu\bigr)^T.
\end{equation}
In our companion work \cite{LPPSS1}, we determine the spectral properties of the linearized operator $\bfM$ with the aid of rigorous numerical computations; the relevant results are summarized in Proposition~\ref{prop:spectral_fourier}. In particular, we prove that $\bfM$ has purely absolutely continuous spectrum $[1,\infty)$ and that the threshold energy $1$ is neither an eigenvalue nor a resonance. Moreover, $\bfM$ has a unique positive gap eigenvalue (internal mode) $\lambda^2\in(0,1)$ of multiplicity two with numerical value $\lambda^2\approx0.7774$. The associated eigenspace is spanned by
\begin{equation}\label{eq:bmYsdef1}
    \bmY_1 := \bigl( Y_1, Y_2, 0, 0\bigr)^T, \qquad \bmY_2 := \bigl(0, 0, Y_1, Y_2\bigr)^T, \qquad \bfM \bmY_j = \lambda^2\bmY_j,
\end{equation}
where $Y_1 := \lambda^{-1} U \psi$, $Y_2 := - \lambda^{-1} \psi'$, and where $\psi$ is the $L^2_{r \ud r}$-normalized\footnote{This means $\int_0^\infty \psi(r)^2r\, dr=1$ without a $2\pi$ factor. That $\|\bmY_1\|_{L^2_{r \ud r}} = \|\bmY_2\|_{L^2_{r \ud r}} = 1$ follows from the eigenvalue equation and integration by parts. } radial ground state of $-\Delta + U^2$ with energy $\lambda^2$.
We refer to \cite{AIMCQ,AIBPMCNOQ} for related numerical work.

\medskip
\begin{figure}[ht] 
\centering
\begin{tikzpicture} 
\coordinate[label=below:\small {$0$}] (zero) at ($(0:0)+(270:2.1)$);
\coordinate[label=below:{\small $1$}] (one) at ($(0:6)+(270:2.1)$);
\coordinate[label=below:\textcolor{blue!50!black}{{\footnotesize \shortstack{internal mode\\$\lambda^2 \approx 0.7774$} }}] (im) at ($(0:4.668)+(270:2.1)$); 
\draw[->,thick] ($(180:1.0)+(270:2.1)$) --++ ($(0:11.5)$);
\draw[line width=5pt,color=orange,draw opacity=0.45] (one) --++ ($(0:4.5)$);
\draw[thick] ($(zero)+(270:0.08)$) -- ($(zero)+(90:0.08)$);
\draw[thick] ($(one)+(270:0.08)$) -- ($(one)+(90:0.08)$);
\fill[blue!50!black] (im) circle (2.5pt);
\end{tikzpicture}
\caption{Spectral features of the linearized operator $\bfM$: The orange band indicates the continuous spectrum. The blue dot represents the unique positive gap eigenvalue (internal mode).} \label{fig:spectral_features}
\end{figure}
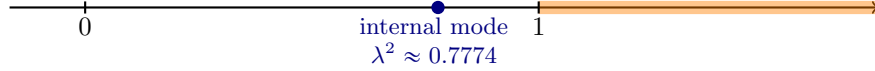

Correspondingly, we enact the spectral decomposition
\begin{equation} \label{equ:spectral_decomposition}
    \bmr(t) = \bmu(t) + z_1(t) \bmY_1 + z_2(t) \bmY_2
\end{equation}
with
\begin{equation}\label{equ:spectral_decomposition2}
    \bmu(t) := \bigl( u_1(t), u_2(t), u_3(t), u_4(t) \bigr)^T 
\end{equation}
and 
\begin{equation}
    z_j(t) := \langle \bmY_j, \bmr(t) \rangle := \int_0^\infty \bmY_j(r) \cdot \bmr(t,r) \, r \ud r, \quad j = 1, 2.
\end{equation}
We also introduce the orthogonal projection to the continuous spectral subspace of $\bigl(L^2_{r\ud r}(\bbR_+)\bigr)^4$ relative to the linearized operator $\bfM$,
\begin{equation} \label{equ:definition_bfP_c}
    \bfP_c \bmv := \bmv - \langle \bmY_1, \bmv \rangle \bmY_1 - \langle \bmY_2, \bmv \rangle \bmY_2.
\end{equation}
From \eqref{equ:sys3} and \eqref{equ:eta0_equation} we obtain a system of nonlinear Klein-Gordon equations for the radiation term $\bmu(t)$ and a system of nonlinear oscillators for the internal mode components $z_1(t)$, $z_2(t)$,
\begin{equation} \label{equ:system_evol_equations_u_z}
    \left\{ \begin{aligned}
        (\pt^2 + \bfM) \bmu &= \bfP_c \bmN, \\ 
        (\pt^2 + \lambda^2) z_1 &= \langle \bmY_1, \bmN \rangle, \\ 
        (\pt^2 + \lambda^2) z_2 &= \langle \bmY_2, \bmN \rangle,
    \end{aligned} \right.
\end{equation}
coupled to an elliptic equation for the temporal component $\eta_0$,
\begin{equation}\label{equ:system_evol_equations_u_z_part2}
    (-\Delta + U^2) \eta_0 = \calN_0.
\end{equation}
In terms of the variables $\bmu$ and $z_1, z_2$, the gauge condition \eqref{equ:intro_Stuart_gauge_polar} reads
\begin{equation}\label{eq:Stuartgaugeintro2}
    \frac{1}{r}\partial_r\big( r (u_4+z_2Y_2) \big)+Uu_3=0.
\end{equation}

\subsection{Main result} \label{subsec:main_result}

We are now in the position to provide a precise statement of our asymptotic stability result for equivariant perturbations of the degree-one vortex under the evolution of the self-dual abelian Yang-Mills-Higgs system of equations.

In order to measure the decay of the internal mode components, we introduce the quantity
\begin{equation} \label{eq:zeddefintro1}
    \mathfrak{z}(t):=\frac{1}{4}\sum_{j=1}^2\big(|z_j(t)|^2+\lambda^{-2}|\partial_tz_j(t)|^2\big).
\end{equation}
Moreover, we denote by $0 < \delta \ll 1$ a small absolute constant that will be used to measure losses from various endpoint estimates. We pick a large integer $N\gg1$ to denote the Sobolev regularity of the data and we assume $N \geq 100\delta^{-1}$. 

\begin{theorem} \label{thm:main}
There exists a small absolute constant $0 < \varepsilon_0 \ll 1$ with the following property:
Assume the initial data $(\bmu(0),\partial_t\bmu(0))$ and $(z_j(0),\partial_t z_j(0))$, $j=1,2$, satisfy the gauge condition~\eqref{eq:Stuartgaugeintro2}.
\setlength{\leftmargini}{1.5em}
\begin{itemize}
\item[(i)] Suppose for some $0 < \varepsilon \leq \varepsilon_0$ that
\begin{equation} \label{equ:new_assumption1}
    0 \leq \sqrt{\zed(0)} \leq \varepsilon,
\end{equation}
and that
\begin{equation} \label{equ:new_assumption2}
    \bigl\| \jap{x} \bigl( e^{i\theta} \bigl( \bmu(0), \pt \bmu(0) \bigr) \bigr) \bigr\|_{H^2_x(\bbR^2) \times H^1_x(\bbR^2)} + \bigl\| e^{i\theta} \bigl( \bmu(0), \pt \bmu(0) \bigr) \bigr\|_{H^N_x(\bbR^2) \times H^{N-1}_x(\bbR^2)} \leq \varepsilon.
\end{equation}
Then the solution to \eqref{equ:system_evol_equations_u_z}-\eqref{equ:system_evol_equations_u_z_part2} with this initial data exists for all times and the gauge condition \eqref{eq:Stuartgaugeintro2} is satisfied. Moreover, the following decay estimates hold 
\begin{equation} \label{eq:decay_zed_mainthm}
    \zed(t) \leq \frac54 \frac{\varepsilon^2}{1 + \Gamma_1 \varepsilon^2 t}, \quad 0 \leq t < \infty, 
\end{equation}
for some constant $\Gamma_1 > 0$ and
\begin{equation} \label{eq:decay_radiation_mainthm} 
    \| e^{i\theta}\bmu(t) \|_{L^\infty_x(\bbR^2)} \lesssim \begin{cases} 
                            \min\bigl\{ \varepsilon, \varepsilon^{1-10\delta} t^{-\frac12-4\delta} \bigr\},  &0 \leq t \leq \varepsilon^{-2}, \\
                            t^{-1+\delta},  &\varepsilon^{-2} \leq t < \infty.
                         \end{cases} 
\end{equation}
\item[(ii)]
If, in addition to \eqref{equ:new_assumption1} and \eqref{equ:new_assumption2}, we assume that
\begin{equation} \label{equ:new_assumption_frakz_equal}
    \sqrt{\frakz(0)} = \varepsilon,
\end{equation}
then we have 
\begin{equation} \label{equ:decay_zed_lower_bound_mainthm}
    \frac34 \frac{\varepsilon^2}{1 + \Gamma_0 \varepsilon^2 t} \leq \zed(t) \leq \frac54 \frac{\varepsilon^2}{1 + \Gamma_1 \varepsilon^2 t}, \quad 0 \leq t < \infty, 
\end{equation}
for some constants $\Gamma_0 > \Gamma_1 > 0$.
\end{itemize}
\end{theorem}

\begin{remark}
    The constants $\Gamma_0 > \Gamma_1 > 0$ in the statement of Theorem~\ref{thm:main} are the Fermi Golden Rule constants defined in \eqref{eq:Gamma01def} in the proof of Proposition~\ref{propODEdamp}.
    The fact that $\Gamma_1 > 0$ is verified in \cite[Proposition~4.3]{LPPSS1} using rigorous numerics.
\end{remark}

The core of the proof of Theorem~\ref{thm:main} is the long-time analysis of
the coupled system consisting of the nonlinear Klein--Gordon equations
\eqref{equ:system_evol_equations_u_z} for the radiation term and the nonlinear
oscillator equations for the internal mode components. In two space dimensions,
the dispersive decay of Klein--Gordon waves is weak enough that the quadratic
terms in the equation for the radiation term cannot be treated perturbatively. We therefore
develop a version of the space-time resonance method with potential in the
half-line setting. A key point is that much of the normal form analysis is
carried out relative to the flat Klein--Gordon flow, while transference
estimates derived from the distorted Fourier theory allow us to pass back to
the Klein--Gordon flow generated by the linearized operator $\bfM$.

The decay of the internal mode components is caused by their resonant nonlinear
coupling to the continuous spectrum, a mechanism known as nonlinear radiation
damping. This mechanism naturally produces the time scale
$t\sim\varepsilon^{-2}$ and explains the two regimes appearing in the decay
estimates \eqref{eq:decay_zed_mainthm} and \eqref{eq:decay_radiation_mainthm}.
In the course of the proof we obtain more refined information on both the
radiation term and the internal mode components; see the bootstrap assumptions
in Subsection~\ref{subsec:bootstrap_setup}.

We refer to Section~\ref{sec:overview_of_proof} for a detailed overview of the main ideas entering the proof of Theorem~\ref{thm:main}.

\subsection{References} 

It is customary to distinguish between gauged (or magnetic) vortex models, in which the scalar field is coupled to an electromagnetic field, and
non-gauged (or global) vortex models; see the standard monographs \cite{bookJT,bookManSut}. 
In both settings, one may further distinguish between the static theory and the study of the dynamics.

The static theory of magnetic vortices in the abelian Yang-Mills-Higgs model at the self-dual coupling was largely developed in \cite{Taubes80_1,Taubes80_2,bookJT}. See also \cite{Plo1,BeCh} for the equivariant setting. In the self-dual case vortices are energy minimizers within their degree class. Using the Bogomolny structure, \cite{Taubes80_1,Taubes80_2,bookJT} proved the existence and uniqueness of vortices for any finite prescribed vortex set, that is, the zeros of the scalar field $\phi$. Low energy dynamics as geodesic motion on the moduli space of self-dual vortices was predicted in \cite{Manton82}. See also \cite{bookManSut}.  This geodesic approximation was verified near the self-dual regime and on a large time interval in \cite{Stu}. 
Away from the critical coupling vortices are not necessary energy minimizers in their energy class and may not be spectrally stable. For a detailed spectral analysis and results on the dynamics of vortices away from the self-dual coupling see \cite{GusSig2,GusVort,GusVort0}. 
For the hyperbolic abelian Yang-Mills-Higgs equations \cite{BurMon1} proved global existence in the Lorenz gauge for all values of the coupling constant. Asymptotic stability of the vacuum states $(\phi,A) = (e^{i\kappa},0)$, $\kappa \in [0,2\pi)$, was proved in \cite{Tsutsumi03}. 
For the related Schr\"odinger flow, we refer to \cite{DemStu09}, and for the corresponding heat flow dynamics to \cite{GusVort0,GusVort,Zhao26,StuDem}. 
Related questions concerning vortex dynamics on Riemann surfaces have been studied in \cite{Stuart99,Dem13,LanSha1}. 

The static solutions to the (non-gauged) complex Ginzburg-Landau equation in the plane are also known as vortices, but they have infinite energy. 
Existence and uniqueness of equivariant vortices were proved in \cite{CEQ94,HerHer94}. Their dynamics have been studied extensively for the associated nonlinear Schr\"odinger equation, also known as the Gross-Pitaevskii equation. 
Global existence for finite-energy data was established in \cite{Gerard}, while the final-state problem for small perturbations of the vacuum was considered in \cite{GNT1}. Orbital stability of the degree-one vortex was proved in \cite{GPS21}; see also \cite{Wein,dPFK}. Further results
on non-gauged vortex dynamics can be found in \cite{Neu,OvSig1,OvSigNon1,OvSigNon2,OvSig2,delPJM}. 
Detailed analyses of the corresponding linearized operator, including the construction of the associated distorted Fourier transform, were carried out in
\cite{PaPu,CGP25,LSS25}.

From a broader perspective, the present work is concerned with the dynamics of topological solitons and more specifically  with asymptotic stability
problems in settings where relatively weak dispersion and low-power nonlinearities interact with additional spectral obstructions, such as internal modes or threshold resonances.
The literature on asymptotic stability of solitons is vast, and we do not attempt a comprehensive review here.
General overviews and further references can be found in the survey articles \cite{Martel_ICM,CuccMaeda20_Survey,Germain24_Review,KMM17_Review,Tao09}.
For a sample of related results in dimension one on the asymptotic stability of kinks in scalar field theories on the line and of ground states of focusing nonlinear Klein-Gordon equations on the line we refer to \cite{GermPusZhang22,GP20,KK11_2,DelMas20,CLL,KoYu23,ChenLuhr,KMM17,KM24,CuccMaeda_kink_23,AMP20,LSch,KriSchNLS,BP95,KNS12,KMM19,CuccMaedaMurgScrob23,LSch23,KaPu22,LL1,PalaciosPusateri24} and the references therein. 
For a sample of related results in higher dimensions see \cite{GKT07,GKT08,BejTat14,BPT24,BKT13,KriSch07,KLS14,ZeLi22,LMZ18,LOS17,LLOS23,NakSch11,NakSch12,Schlag09} and the references therein. 
Important progress has also been made on the stability of stationary solutions to quasilinear wave equations, particularly in general relativity. 
Although a review of this literature lies beyond the scope of the present paper, we mention \cite{DKSW16}, which is especially relevant to the discussion in
Section~\ref{sec:overview_of_proof}, and refer to \cite{SJOhICM1} for a survey of some of the developments in the quasilinear setting.

A central aspect of this work is the decay, under the full nonlinear evolution, of the oscillations associated with the internal mode. 
The nonlinear mechanism responsible for this decay is known as radiation damping. The use of radiation damping to establish the decay
of internal modes goes back to \cite{BP95,SW}; see also \cite{Sig}.
Since these pioneering works, the mechanism has been established in a variety of models possessing internal modes; see, for example,
\cite{BusSul03,BamCuc,SW,BP95,TsaiYau02,GANGS07,GW08,CucMiz08,LP22,CuccMaeda_kink_23,KM24,KMM17} and the references therein, as well as the survey article \cite{CuccMaeda20_Survey}.      

\subsection{Organization of the paper}

Section \ref{sec:overview_of_proof} provides an overview of the proof of Theorem \ref{thm:main}. We first explain the main difficulties and the mechanisms used to overcome them through a scalar model problem, and then describe how these ideas are implemented for the full system.

In Section \ref{sec:setup} we introduce the good-bad decomposition, which separates the spatially localized interactions involving the internal modes from the purely radiative interactions, and the flat-sharp decomposition, which isolates the component amenable to analysis using the standard Fourier transform. We also introduce the corresponding decomposition of the temporal component, classify the nonlinearities in Corollary \ref{cor:types_of_nonlinearities}, and formulate the bootstrap assumptions used throughout the proof.

Section \ref{sec:linear_theory} collects the linear inputs established in our companion works. Proposition \ref{prop:spectral_fourier} records the spectral properties of the linearized operator $\bfM$ and the associated distorted Fourier transform. We then state the dispersive and integrated local energy decay estimates for the flat and the perturbed Klein-Gordon evolutions (see in particular Propositions \ref{prop:dispersive_estimate_with_potential} and \ref{prop:ILED}). Finally, Proposition \ref{prop:transference} relates weighted norms of profiles defined using the flat Fourier transform to solutions of the perturbed Klein--Gordon equation associated with $\bfM$.

In Section \ref{sec:working_title_eta0_bounds} we estimate the temporal component. Proposition \ref{prop:eta0c_bounds} establishes pointwise decay and high Sobolev bounds for its good component, while Proposition \ref{prop:eta0d_bounds} proves spatially weighted pointwise and high Sobolev estimates for its bad component. In particular, the latter estimates quantify the spatial localization inherited from the internal modes.

Section \ref{sec:working_title_flat_bounds} contains the main nonlinear analysis for the auxiliary flat components of the good radiation. We perform quadratic normal form transformations for the non-localized interactions and incorporate the resulting terminal boundary corrections into the final data. Proposition \ref{prop:weighted_energies_g} proves the weighted estimate for the corresponding renormalized profiles, while Proposition \ref{prop:dispersive_decay_g} gives their improved dispersive decay. Uniform high Sobolev and integrated local energy estimates are proved in Propositions \ref{prop:HN_g} and \ref{prop:HN_ILED_g}. At the end of the section, these auxiliary estimates on $\bbR^2$ are transferred to the first angular modes defining the flat component of the good radiation.

In Section \ref{sec:working_title_sharp_bounds} we use the localization of the potential, together with the linear decay and transference estimates, to prove dispersive, high Sobolev, integrated local energy, and weighted profile bounds for the sharp component (see Propositions \ref{prop:dispersive_decay_g_sharp}, \ref{prop:HN_g_sharp}, and \ref{prop:weighted_energies_g_sharp}). Section \ref{sec:working_title_w_bounds} establishes the corresponding bounds for the bad component (see Propositions \ref{prop:dispersive_decay_h}, \ref{prop:HN_h}, and \ref{prop:weighted_energies_h}). The spatial localization of its internal mode forcing yields the required pointwise decay, while the same forcing produces the characteristic $t^{1/2}$ growth of its weighted profile norm.

The internal mode analysis is carried out in Sections \ref{sec:FGR1} and \ref{secIM2}. In Section \ref{sec:FGR1} we derive the effective system of ODEs for the two internal mode profiles and use the constrained Hamiltonian structure to identify the dissipative form generated by the Fermi Golden Rule. Proposition \ref{propODEdamp} reduces the decay law for $\vert Z_1\vert^2+\vert Z_2\vert^2$ to suitable estimates for the remainders in this system. These estimates are established in Section \ref{secIM2} (see Proposition \ref{propODERbounds}). The effective ODE analysis in Section \ref{sec:FGR1}, combined with the remainder estimates in Section \ref{secIM2}, establishes the nonlinear radiation damping of the internal modes.

Finally, in Section \ref{sec:conclusion_proof_thm_main} we combine the estimates for the temporal component, the flat, sharp, and bad radiation components, and the internal modes. We improve all the bootstrap assumptions, close the continuation argument, and complete the proof of Theorem \ref{thm:main}.

\medskip
\subsection{Notation and conventions} \label{subsec:notation}

We collect several notational conventions and basic definitions that will be used throughout this work.

\medskip 
\noindent {\it Standard conventions.}
We write $a \lesssim b$ when $a \leq C b$ for some absolute constant $C>0$ independent of $a$ and $b$; $a \simeq b$ means that $a\lesssim b$ and $b\lesssim a$. When $a$ and $b$ are expressions depending on our main variables or parameters, the inequalities are assumed to hold uniformly over these. 
We denote by $a+$ (respectively $a-$) a number $b > a$ (respectively $b < a$) that can be chosen arbitrarily close to $a$. 
For arbitrary $b$ and non-negative $a$, we use the short-hand notation $b = O(a)$ if $|b| \leq C a$.
Throughout, we use the Japanese bracket notation
\begin{equation}
    \jap{x} := (1+|x|^2)^{\frac12}, \quad \jap{\xi} := (1+|\xi|^2)^{\frac12}, \quad x, \xi \in \bbR^2,
\end{equation}
as well as $\jap{t} := (1+t^2)^{\frac12}$, $t \in \bbR$.
We sometimes use a dot symbol `$\cdot$' to distinguish the bounds on various quantities involved in product estimates.
Moreover, we use standard notation for Lebesgue and Sobolev norms such as $L^p$ and $W^{s,p}$ with $H^s = W^{s,2}$.

\medskip 
\noindent {\it Constants.} 
We denote by $C_0 \geq 1$ the constant in the bootstrap assumptions \eqref{equ:bootstrap_assumption_weighted_f}, \eqref{equ:bootstrap_assumption_Linfty_f}, \eqref{equ:bootstrap_assumption_HN_f}, \eqref{equ:bootstrap_assumption_HN_ILED_f}.
Throughout, $C_1 \geq 1$ indicates an absolute constant whose size may still change from place to place, while $C_2 \equiv C_2(C_0) \geq 1$ denotes a constant whose size depends on $C_0$.

\medskip 
\noindent {\it Radial $L^2$ inner products.}
For real-valued functions $f, g \in L^2_{r \ud r}(\mathbb R_+)$ and $\mathbb R^4$-valued functions $\bmY, \bmZ \in \bigl(L^2_{r \ud r}(\mathbb R_+)\bigr)^4$, we use the notation
\begin{equation}
    \langle f, g \rangle := \int_0^\infty f(r) g(r) \, r \ud r, \quad \langle \bmY, \bmZ \rangle := \int_0^\infty \bmY(r) \cdot \bmZ(r) \, r \ud r,
\end{equation}
where $\cdot$ denotes the Euclidean inner product on $\mathbb R^4$.

\medskip 
\noindent {\it Flat Fourier transform.}
Our conventions for the Fourier transform of a Schwartz function on $\bbR^n$ are
\begin{equation} \label{eq:FT}
 \begin{aligned}
  \widehat{\calF}[f](\xi) &= \hatf(\xi) = \frac{1}{(2\pi)^{\frac{n}{2}}} \int_{\bbR^n} e^{-ix\cdot\xi} f(x) \, \ud x, \\
  \widehat{\calF}^{-1}[f](x) &= \check{f}(x) = \frac{1}{(2\pi)^{\frac{n}{2}}}  \int_{\bbR^n}  e^{ix\cdot\xi} f(\xi) \, \ud \xi.
 \end{aligned}
\end{equation}
Then the convolution laws are given by
\begin{equation} \label{eq:convolution_laws}
 \widehat{\calF}[f \ast g] = (2\pi)^{\frac{n}{2}} \hatf \hatg, \quad \widehat{\calF}[f g] = \frac{1}{(2\pi)^{\frac{n}{2}}} \hatf \ast \hatg
\end{equation}
for $f, g \in \calS(\bbR^n)$.
We denote by $\langle D\rangle^s$ and $|D|^s$, $s\in\mathbb R$, the Fourier multiplier operators with symbols $\langle\xi\rangle^s$ and
$|\xi|^s$, respectively.
Moreover, for $\kappa \in \{\pm\}$ we use the short-hand notation $\widehat{f}^{\,\kappa}(\xi)$, where
\begin{equation}
    \widehat{f}^{+}(\xi) := \widehat{f}(\xi), \quad \widehat{f}^{-}(\xi) := \overline{\widehat{f}(\xi)}.
\end{equation}

\medskip
\noindent
{\it Cut-offs and Littlewood-Paley projections}.
We denote by $\chi \in C_c^\infty(\bbR^2)$ a smooth, non-negative, radially symmetric bump function with the property that $\chi(\xi) = 1$ for $|\xi| \leq 1$ and $\chi(\xi) = 0$ for $|\xi| \geq 2$. Then we set
\begin{equation}
    \varphi_0(\xi) := \chi(\xi), \quad \xi \in \bbR^2,
\end{equation}
and for every integer $n \geq 1$,
\begin{equation}
    \varphi_n(\xi) := \chi\Bigl(\frac{\xi}{2^n}\Bigr) - \chi\Bigl(\frac{\xi}{2^{n-1}}\Bigr), \quad \xi \in \bbR^2.
\end{equation}
Moreover, for any integer $M\in \bbZ$ we introduce the cut-offs 
\begin{equation}
    \varphi_{\leq M}(\xi) := \chi\Bigl(\frac{\xi}{2^M}\Bigr), \quad \varphi_{>M}(\xi) := 1 - \varphi_{\leq M}(\xi).
\end{equation}
Occasionally, we will also need to resolve low frequencies $\lesssim 1$ dyadically. To this end we use the following cut-off functions for any integer $n \in \bbZ$,
\begin{equation}
    \psi_n(\xi) := \chi\Bigl( \frac{\xi}{2^n} \Bigr) - \chi\Bigl( \frac{\xi}{2^{n-1}} \Bigr).
\end{equation}
Finally, we denote by $\widetilde{\varphi}_n(\xi)$ etc. suitable fattened versions of the cut-off functions $\varphi_n(\xi)$ with the property that $\widetilde{\varphi}_n(\xi) \varphi_n(\xi) = \varphi_n(\xi)$.

The standard Littlewood-Paley projection operators $P_n$ for integers $n \geq 0$ are defined by
\begin{equation}
    \widehat{P_n f}(\xi) := \varphi_n(\xi) \hatf(\xi).
\end{equation}

\section{Overview of the Proof of Theorem~\ref{thm:main}} \label{sec:overview_of_proof}

In this section we give an overview of the proof of Theorem~\ref{thm:main}, emphasizing the main difficulties and the mechanisms used to overcome them.
In Subsection~\ref{subsec:scalar_model_problem} we organize the discussion around a scalar model problem that captures many of the essential features of the full problem, including the coupling between radiation and the internal mode, non-spatially localized quadratic interactions, and derivative nonlinearities. In Subsection~\ref{subsec:implementation_full_system} we then explain how these ideas are implemented for the full nonlinear system and discuss several additional aspects of the proof of Theorem~\ref{thm:main}.

\subsection{The scalar model problem and the functional framework} \label{subsec:scalar_model_problem}
Let
\[
    H:=-\Delta+1+V,
\]
where $\Delta$ denotes the Laplacian on $\bbR^2$ and $V$ is a smooth radial potential satisfying
$|V(x)|\lesssim\jap{x}^{-2}$.  We assume that
\[
    \sigma_{\mathrm c}(H)=[1,\infty), \qquad \sigma_{\mathrm d}(H)=\{\lambda^2\}, \qquad H\fy=\lambda^2\fy, \qquad \|\fy\|_{L^2_x}=1, \qquad 2\lambda>1,
\]
and that the threshold $1$ is neither an eigenvalue nor a resonance.  
The projection onto the continuous spectral subspace of $L^2(\bbR^2)$ relative to $H$ is denoted by $P_c$.
We consider the scalar model problem
\begin{equation} \label{eq:overview_single_model}
    \left\{
    \begin{aligned}
        (\pt^2+H)u
        &=
        P_c(\fy^2) z^2
        +P_c\bigl(
        \underline{\Phi}_{\real}u^2
        +u\partial_1u 
        \bigr),\\
        \pt^2z+\lambda^2z
        &=
        2z\langle u,\fy^2\rangle,
    \end{aligned}
    \right.
\end{equation}
where the non-spatially localized coefficient $\underline{\Phi}_{\real}$ is given by 
\[
    \underline{\Phi}_{\real}(x)
    :=U(|x|)\frac{x_1}{|x|},
    \qquad
    U(r)\rightarrow1
    \quad\text{as }r\rightarrow\infty .
\]
We refer to $u=P_cu$ as the radiation and to $z$ as the internal mode component.
Our goal is to prove decay for $u$ and $z$ under the same smallness assumptions on the initial data as in the statement of Theorem~\ref{thm:main}. 

We briefly describe the role of the terms retained in \eqref{eq:overview_single_model}. The source $P_c (\varphi^2) z^2$ generates the leading radiation produced by the internal mode and dictates the weighted and pointwise bounds for the radiation. The term $2z\angles{u}{\varphi^2}$ is the leading feedback in the equation for the internal mode component; reinserting into it the radiation generated by $P_c(\varphi^2) z^2$ produces the Fermi Golden Rule damping. The term $\underline{\Phi}_{\real} u^2$ represents a variable coefficient non-localized quadratic term. The most delicate nonlinear analysis concerns the treatment of this term. The term $u\partial_1u$ is included because it is a top-order derivative quadratic term.

We note that the actual abelian Yang-Mills-Higgs equation presented in the introduction preserves degree-one equivariance, while the model problem \eqref{eq:overview_single_model} does not.
If $\bmu = e^{i\theta}u$ denotes a function of degree one, then a quadratic product of $\bmu$ and $\overline{\bmu}$ has degree zero or two.  By contrast, $\underline{\Phi} |\bmu|^2$ has degree one, while a product of $\bmu$ and $\partial_j\bmu$ can also have a degree-one component.  Here the scalar term $\underline{\Phi}_{\real}u^2$ models $\Re(\underline{\Phi}|\bmu|^2)$, so $u^2$ represents the degree-zero quantity $|\bmu|^2$ in the full system.
Also, the full system contains additional spatially localized quadratic and non-localized cubic terms, but their treatment
is not harder than the nonlinear terms in \eqref{eq:overview_single_model}, so we omit them from this overview.

The localized and non-localized nonlinearities in the radiation equation require different estimates, which motivates our first decomposition.
In part inspired by \cite{LP22}, we introduce the {\it good--bad decomposition}
\[u=v+w,\]
where
\begin{equation} \label{eq:overview_good_bad_model}
    \left\{
    \begin{aligned}
        (\pt^2+H)v &= P_c\calN_c(u),
        \quad  & \quad  (v,\pt v)|_{t=0}&=(u,\pt u)|_{t=0},\\
        (\pt^2+H)w &= P_c\calN_d(z),
        \quad  & \quad  w(0)&=\pt w(0)=0,
    \end{aligned}
    \right.
\end{equation}
with
\begin{equation} \label{eq:overview_single_model_nonlinearity}
    \calN_d(z) := P_c(\fy^2) z^2,
    \qquad
    \calN_c(u) := \underline{\Phi}_{\real}u^2 + u\partial_1u.
\end{equation}
The \emph{good part} $v$ carries the initial data for the radiation and is driven by the radiative part of the nonlinearity, namely $\calN_c(u)$,
whereas the \emph{bad part} $w$ has zero initial data and is forced by the localized internal mode source term $\calN_d(z)$.
The good-bad terminology derives from the weighted, rather than pointwise, estimates satisfied by the two components: the bad part satisfies much weaker weighted estimates, but still obeys the required dispersive bounds, see Subsection~\ref{subsubsec:overview_bad_component}.

Establishing decay of $v$ is the most involved part of the nonlinear analysis for \eqref{eq:overview_single_model}. 
Our approach is based on an additional \emph{flat-sharp decomposition} that combines the space-time resonance method for the flat Klein-Gordon equation with integrated local energy decay (ILED) estimates and transference relations for $\partial_t^2+H$, see Subsection~\ref{subsec:overview_flat_sharp}. 

To describe the functional framework we introduce the flat profile of the radiation $u(t)$ and the profile of the internal mode $z(t)$,
\begin{equation}
 \begin{aligned}
  f(t)&:=e^{-it\langle D\rangle}(2i\langle D\rangle)^{-1}
  (\partial_t+i\langle D\rangle)u(t),\\
  Z(t)&:=e^{-it\lambda}(2i\lambda)^{-1}(\partial_t+i\lambda)z(t).
 \end{aligned}
 \label{eq:overview_profiles}
\end{equation}
Our schematic bootstrap setup is
\begin{equation} \label{eq:overview_large_time_framework}
 |Z(t)|\leq 2C_0t^{-\frac{1}{2}},
 \quad
 \bigl\|\langle\xi\rangle^2\langle\nabla_\xi\rangle\widehat{f}(t)\bigr\|_{L_\xi^2}
 \leq 2C_0 t^{\frac12},
 \quad
 \bigl\|e^{it\langle D\rangle}f(t)\bigr\|_{L_x^\infty}
 \leq 2C_0 t^{-1+\delta}, \quad 0 \leq t \leq T,
\end{equation}
together with a uniform high Sobolev norm bound and an integrated local energy decay estimate
\begin{equation} \label{eq:overview_large_time_framework_ILED}
\sup_{0\leq t\leq T}\|f(t)\|_{H^N_x}+\sum_{1 \leq k \leq N} \, \bigl\| \jx^{-\frac32-\kappa} |D|^k e^{it\jD} f(t)\bigr\|_{L^2_t([0,T]; L^2_x)} \leq 2C_0 \varepsilon.
\end{equation}
Here $T$ is the final bootstrap time, $C_0$ is a large bootstrap constant, $0<\delta\ll1$ and  $N\gg\delta^{-1}\gg1$ are as in Theorem~\ref{thm:main}, and $0<\kappa\ll1$ is a fixed small number used to define the ILED norm. The second norm in \eqref{eq:overview_large_time_framework} is our main weighted norm. To be precise, the estimates \eqref{eq:overview_large_time_framework} need to be modified at short times $t\leq \varepsilon^{-2}$, see Subsection~\ref{subsubsec:FGRmodel}. The estimates needed at earlier times require additional work but do not present major new conceptual difficulties. Their precise form is recorded in the main bootstrap assumptions in Subsection~\ref{subsec:bootstrap_setup}. 
In the remainder of this overview we focus on the regime $t\geq \varepsilon^{-2}$, explaining how these bounds naturally arise and how the bootstrap is closed.

For the moment, let us observe how the bound on the weighted norm in \eqref{eq:overview_large_time_framework} does not guarantee 
the $L^\infty_x$ decay stated there. 
Indeed, the standard weighted dispersive estimate gives, 
\begin{equation} \label{eq:overview_weighted_dispersive}
    \bigl\|e^{it\jD}h\bigr\|_{L^\infty_x}
    \lesssim
    \jap{t}^{-1+\delta}
    \Bigl(    
    \bigl\|\jxi^2\nabla_\xi\widehat h\bigr\|_{L^2_\xi}
    +
    \bigl\|\jxi^2\widehat h\bigr\|_{L^2_\xi}
    \Bigr).
\end{equation}
Part of our strategy is to bootstrap independently the weighted norm and the dispersive decay in different ways for the good part $v$ and the bad part $w$. 

The use of the flat Fourier transform also seems incompatible with equation \eqref{eq:overview_single_model}, where the linear operator is not the flat Klein-Gordon operator. The flat-sharp decomposition introduced below is designed to connect the flat and perturbed Klein-Gordon evolutions. A key role in making this connection is played by the transference estimates described in Subsection~\ref{subsubsec:transference_relations_overview} below. The transference estimates use the relation between $\nabla_\xi\widehat f$ and $\Omega_{0j}u$, where $\Omega_{0j}:=x^j\partial_t+t\partial_j$ denotes a Lorentz boost vector field. For future reference in the overview, we record this relation. For any Schwartz function $\psi(t,\cdot)$,
\begin{equation} \label{eq:overview_Omega_Omegahat1}
    \widehat{\Omega_{0j} \psi}(t)(\xi) = i \widehat{\Omega}_{0j}\widehat{\psi}(t)(\xi), \qquad \mathrm{where} \qquad \widehat{\Omega}_{0j}:=t\xi_j+\partial_t\partial_{\xi^j}.
\end{equation}
Now if $\Psi(t)
    :=
    e^{-it\jD}(2i\jD)^{-1}
    (\partial_t+i\jD)\psi(t)$ denotes the profile of $\psi$, then by direct computation
\begin{equation}
\label{eq:overview_profile_boost_identity}
    \begin{aligned}
    \jxi^2\partial_{\xi_j}\widehat \Psi(t)
    ={}&
    -\xi_j\widehat \Psi(t)
    +
    \frac12e^{-it\jxi}
    \partial_{\xi_j}
    \bigl(
        \partial_t^2+\jxi^2
    \bigr)\widehat \psi(t) -
    \frac12e^{-it\jxi}
    (\partial_t+i\jxi)
    \widehat\Omega_{0j}\widehat \psi(t),
    \end{aligned}
\end{equation}
and consequently,
\begin{equation}
\label{equ:overview_transference_flat}
    \begin{aligned}
    \bigl\|\jxi^2\partial_{\xi_j}\widehat \Psi(t)\bigr\|_2
    \lesssim{}&
    \|\Psi(t)\|_{H^1_x}
    +
    \bigl\|
        \jx(\partial_t^2+\jD^2)\psi(t)
    \bigr\|_{L^2_x} + 
    \bigl\|
        (\partial_t+i\jD)\Omega_{0j}\psi(t)
    \bigr\|_{L^2_x}.
    \end{aligned}
\end{equation}
See \cite[Section~4]{LPPSS2} for more details on how \eqref{eq:overview_profile_boost_identity} and \eqref{equ:overview_transference_flat} are used to prove the transference estimate in Proposition~\ref{prop:transference}. We also refer to the introduction of \cite{LPPSS2} for a discussion of how in principle the transference relations can be obtained at the level of $\nabla_\xi$ of the Fourier transforms of the profiles, and without using the Lorentz boost vector fields.

\medskip 

We emphasize the modular nature of the flat-sharp approach developed here for proving decay. It consists of two main components:
\begin{itemize}
    \item[(i)] a robust nonlinear analysis for the flat Klein-Gordon equation,
    \item[(ii)] use of ILED estimates and transference relations for the passage between the flat analysis and the Klein-Gordon evolution with potential.
\end{itemize}
Related instances of this general paradigm can be found, for example, in \cite{RodSch04,DR10,Tat08}. In principle, each of the two components can be implemented by a variety of methods.
For the flat nonlinear analysis, we find the space-time resonance method particularly effective, as it provides a streamlined framework for implementing normal forms for quadratic nonlinearities.
To establish the linear theory, especially ILED estimates and transference relations, we use the distorted Fourier transform. 

For the implementation of the flat-sharp approach we also find it instructive to draw analogies with the classical vector field method. In Subsection~\ref{subsec:overview_vf} we informally describe how in that context the flat-sharp decomposition would be implemented for a simpler cubic nonlinearity. The interested reader may find this independent account clarifying before reading the discussion in Subsections~\ref{subsec:overview_flat_sharp}-\ref{subsubsec:overview_cubic_commutation}. Subsection~\ref{subsec:overview_vf} can also be safely skipped by the reader who is eager to see the immediate implementation for \eqref{eq:overview_single_model} and the full problem, where we work entirely on the Fourier side.

\subsubsection{Radiation damping and the decay of the internal mode} \label{subsubsec:FGRmodel}

The main difficulty in proving decay of the internal mode $z$ is that the fundamental solutions $e^{\pm i\lambda t}$ of the linear equation $z''+\lambda^2z=0$ are purely oscillatory and hence exhibit no decay.  Thus, any decay of $z$ must come from its nonlinear coupling to the radiation,
through the mechanism usually called radiation damping. The analysis of this
mechanism goes back to \cite{BP95,SW}. In the model \eqref{eq:overview_single_model}, the source
$P_c(\fy^2) z^2$ generates the leading internal-mode contribution to the
radiation, which is precisely the bad part $w$.  Since $w$ has zero initial data,
\begin{equation} \label{eq:overview_fgr_radiation}
    w(t)
    =
    \int_0^t
    \frac{\sin((t-s)\sqrt H)}{\sqrt H}
    P_c(\fy^2) z(s)^2\,\ud s .
\end{equation}
Since $u=v+w$, reinserting \eqref{eq:overview_fgr_radiation} into
$2z\langle u,\fy^2\rangle$ in the equation for $z$ produces the leading
cubic feedback. After separating the oscillations of $z$ by means of the
profile $Z$, one obtains an equation of the form
\begin{equation} \label{eq:overview_scalar_fgr_ode}
    \pt Z
    = - \Gamma |Z|^2 Z - i\Lambda |Z|^2 Z +\calR_Z,
    \qquad
    \Gamma\geq 0, \qquad \Lambda \in \bbR,
\end{equation}
where $\calR_Z$ are remainder terms that involve $Z$ and $u$.
The coefficient $\Gamma$ comes from the $2\lambda$-oscillation in
$z^2$ appearing in the feedback mentioned  above. It is given explicitly in \eqref{eq:overview_scalar_fgr_coefficient}.
Since $2\lambda>1,$ and $\sigma_{\mathrm c}(\sqrt H)=[1,\infty)$, this oscillation can excite the continuous spectrum and transfer energy from
the internal mode to the radiation. The strength of this coupling is measured by
\begin{equation} \label{eq:overview_scalar_fgr_coefficient}
    \Gamma
    =
    c_\lambda\,
    \Im\,\Big\langle
    P_c(\fy^2), \,
    (H-(2\lambda)^2+i0)^{-1}P_c(\fy^2)
    \Big\rangle ,
\end{equation}
where $(2\lambda)^2$ is the corresponding spectral value of $H$, and $c_\lambda>0$ is an explicit constant.  The Fermi Golden Rule condition states that $\Gamma$ does not vanish. It makes the cubic term $-\Gamma|Z|^2Z$ in \eqref{eq:overview_scalar_fgr_ode} dissipative, thereby producing the nonlinear damping of the internal mode. The other cubic term $-i\Lambda |Z|^2 Z$ does not affect the decay of the amplitude $|Z(t)|$, only the oscillatory behavior of $Z(t)$.
For our actual problem \eqref{equ:system_evol_equations_u_z} the corresponding non-vanishing condition is verified with the aid of rigorous numerics in our companion paper \cite{LPPSS1}, see Subsection~\ref{sssec:im}.

In particular, if $\calR_Z=0$, then
\begin{align}\label{eq:overview_scalar_fgr_decay}
    \pt |Z|^2=-2\Gamma|Z|^4, \qquad
    |Z(t)|^2
    =
    \frac{|Z(0)|^2}
    {1+2\Gamma|Z(0)|^2t}.
\end{align}
The proof of decay for $z$ therefore reduces to showing that under the bootstrap assumptions \eqref{eq:overview_large_time_framework}--\eqref{eq:overview_large_time_framework_ILED} the contribution of $\calR_Z$ to the ODE \eqref{eq:overview_scalar_fgr_ode} is perturbative. The relation \eqref{eq:overview_scalar_fgr_decay} also explains the need to separate the time scales $0 \leq t \leq \varepsilon^{-2}$ and $t\geq \varepsilon^{-2}$. Indeed, if $|Z(0)|=\varepsilon$, then according to~\eqref{eq:overview_scalar_fgr_decay}, $|Z(t)|$ is comparable to $\varepsilon$ for $t\leq \varepsilon^{-2}$ and decays at the rate $t^{-\frac{1}{2}}$ for $t\geq \varepsilon^{-2}$. If instead we only assume that $|Z(0)|\leq \varepsilon$, then the first estimate in \eqref{eq:overview_large_time_framework} should only be viewed as an upper bound.

The non-negativity of $\Gamma$ is evident from the explicit formula \eqref{eq:overview_scalar_fgr_coefficient}. Remarkably, it is by now understood that this general structure is present for Hamiltonian equations, see for instance \cite{BamCuc}. Our full system \eqref{equ:system_evol_equations_u_z} does not exactly fit into the framework of  \cite{BamCuc}, or other related existing works, because it is a system, rather than a scalar equation, and more importantly it is not of Hamiltonian form in the classical sense. Instead, in view of the gauge freedom, to write the system in Hamiltonian form one has to use the formalism of constrained Hamiltonians as in Subsection~\ref{subsec:Hamiltonintro}. Nevertheless, we show that this structure is enough to derive a condition that is analogous to the non-negativity of $\Gamma$ in \eqref{eq:overview_scalar_fgr_ode}. For our explicit nonlinearity, we have chosen to demonstrate this by direct computation rather than an abstract derivation, see Section~\ref{sec:FGR1}.

\subsubsection{Weighted growth and pointwise decay of the bad part} \label{subsubsec:overview_bad_component}

The need to treat the good and bad parts of the radiation differently was already observed in a quadratic model problem in three dimensions in \cite{LP22}. 
The bad part $w$ is driven by the localized source $P_c(\fy^2) z^2$.  Although, as shown below, its weighted profile is forced to grow like $t^{1/2}$, the localization of the source still yields 
sufficient, almost pointwise integrable, decay. 
Indeed, using the decay of $Z(s)$ suggested by \eqref{eq:overview_scalar_fgr_decay}, the Duhamel formula \eqref{eq:overview_fgr_radiation} and the $L^1$-$L^\infty$ dispersive estimate for $e^{\pm it\sqrt H}$ give
\begin{equation} \label{eq:overview_bad_pointwise}
    \begin{aligned}
        \|w(t)\|_{L^\infty_x}
        &\lesssim
        \int_0^t
        \frac{1}{\jap{t-s}}
        \frac{\varepsilon^2}{1+\varepsilon^2s}\,\ud s 
        \lesssim t^{-1+\delta},
        \qquad t\gtrsim\varepsilon^{-2}.
    \end{aligned}
\end{equation}
The proof of the $L^1$-$L^\infty$ decay estimate and other linear estimates for the perturbed propagator $e^{\pm it\sqrt H}$ use the 
construction of the matrix distorted Fourier transform associated with the operator \eqref{sysN1compact}. 
These are carried out in our companion paper \cite{LPPSS2}.

For the weighted estimate we introduce the flat profiles of $v$ and $w$ as
\begin{equation}\label{eq:overviewvwflatprofiles}
  g(t):=e^{-it\langle D\rangle}(2i\langle D\rangle)^{-1}
  (\partial_t+i\langle D\rangle)v(t),\qquad h(t):=e^{-it\langle D\rangle}(2i\langle D\rangle)^{-1}
  (\partial_t+i\langle D\rangle)w(t).
\end{equation}
 To estimate  $\nabla_\xi\widehat{f}$ we need to control $\nabla_\xi \widehat{g}$ and $\nabla_\xi \widehat{h}$. The contribution of the good profile $g$ is discussed below in Subsection~\ref{subsec:overview_flat_sharp}. For the bad profile $\nabla_\xi \widehat{h}$, in this subsection we ignore the difference between $\sqrt{H}$ and $\jap{D}$. This issue is addressed by the transference relations which are discussed in Subsections~\ref{subsec:overview_flat_sharp} and~\ref{subsubsec:transference_relations_overview}. Then, when $\nabla_\xi$ falls on the linear propagator $e^{-is\jap{\xi}}$ in the Duhamel formula for the profile, it produces a factor $s$. Using (dual) ILED, the leading term is therefore controlled schematically by
\begin{equation}
 \bigl\|\langle\xi\rangle^2\nabla_\xi\widehat{h}(t)\bigr\|_{L_\xi^2}
 \lesssim
 \bigl\|s|Z(s)|^2\bigr\|_{L_s^2([0,t])}
 +\{\text{lower-order terms}\}
 \lesssim t^{\frac12}.
 \label{eq:overview_bad_weighted}
\end{equation}
ILED estimates turn out to be quite useful at various stages of the proof
and streamline the treatment of many spatially localized nonlinear terms compared to previous works, e.g., \cite{LP22,SW}. The rigorous version of \eqref{eq:overview_bad_weighted} using transference relations is obtained in a similar manner as how we treat the sharp profile in Subsection~\ref{subsec:overview_flat_sharp} below.  The absence of a small coefficient in \eqref{eq:overview_bad_weighted} comes from the fact that
\[
 s|Z(s)|^2\simeq1,
 \qquad s\gtrsim\varepsilon^{-2}.
\]
However, since the implicit constant is independent of $C_0$, estimate \eqref{eq:overview_bad_weighted} improves on the weighted bootstrap estimate in \eqref{eq:overview_large_time_framework} if $C_0$ is chosen sufficiently large. 
 
The (optimal) bound \eqref{eq:overview_bad_weighted} explains why we need to propagate a growing weighted estimate for the full radiation as in \eqref{eq:overview_large_time_framework}; at the same time, pointwise decay of $w$ can be proved directly using the spatial localization of its source.  Instead, for the good part $v$, whose forcing is not spatially localized, the proof of pointwise decay requires a more elaborate argument
and new ideas, which we describe next.

\subsubsection{The flat-sharp decomposition} \label{subsec:overview_flat_sharp}

We now explain how the good component $v$ is treated using the flat-sharp decomposition and, in particular, how we are able to prove that it decays consistently with \eqref{eq:overview_large_time_framework}. 
Let 
\begin{align}\label{eq:overview_fs}
    v = P_c v^\flat + P_c v^\sharp,
\end{align}
where the flat component $v^\flat$ 
solves the flat Klein-Gordon equation
\begin{equation} \label{eq:overview_flat_equation}
    (\pt^2+\jD^2)v^\flat
    =
    \calN_c(u)
    =
    \underline{\Phi}_{\real}u^2
    +u\partial_1u,
\end{equation}
and the sharp component $v^\sharp$ 
solves
\begin{equation} \label{eq:overview_sharp_equation}
    (\pt^2+H) P_c v^\sharp = -P_c (Vv^\flat).
\end{equation}
We want the decay of $v^\flat$ to be determined by the source term in \eqref{eq:overview_flat_equation}, which is better than the linear decay of free Klein-Gordon waves. This is achieved by a judicious choice of final data at the terminal time $T$ of our bootstrap interval. More precisely, we impose final data for $v^\flat$ so that the terminal boundary term produced by the quadratic normal form that we are going to perform cancels; see \eqref{eq:overview_renormalized_profile} in Subsection \ref{subsec:overview_normal_forms} for more details on this.

In order to recover the bootstrap assumptions \eqref{eq:overview_large_time_framework} on the flat profile $f$, 
it is useful to also introduce the flat profiles associated with the flat and the sharp component of the good part of the radiation,
\begin{equation} \label{eq:gflatsharpdefoverview}
   g^\ast(t) := e^{-it\jD}(2i\jD)^{-1} (\pt+i\jD)P_cv^\ast(t),\qquad  \ast\in\{\flat,\sharp\}.
\end{equation}
With our suitable choice of data at $t=T$,
the final-state formulation associated to \eqref{eq:overview_flat_equation} yields faster than linear decay for the flat component, namely
\begin{equation} \label{eq:overview_flat_outputs}
    \bigl\|e^{it\jD}g^\flat(t)\bigr\|_{L^\infty_x}
    \lesssim t^{-\frac32+4\delta},
    \qquad
    \bigl\|\jxi^2\jap{\nabla_\xi}
    \widehat{g}^\flat_{\mathrm{ren}}(t)\bigr\|_{L^2_\xi}
    \lesssim t^{-\frac12+3\delta}.
\end{equation}
Here, $g^\flat_{\mathrm{ren}}$ is the quadratic normal-form renormalization of $g^\flat$, see \eqref{eq:overview_renormalized_profile}.
We refer to Subsection \ref{subsec:overview_normal_forms} for a detailed discussion of the variable coefficient normal form step and to Subsection~\ref{subsubsec:overview_cubic_commutation} for details on how to obtain \eqref{eq:overview_flat_outputs}.

For $g^\sharp$ note that the sharp equation is governed by the perturbed operator $H$. 
The linear transference estimate will allow us to control the weighted norm of this flat profile in terms of localized space-time norms of the forcing in equation~\eqref{eq:overview_sharp_equation}.
More precisely, our main transference estimate applied to \eqref{eq:overview_sharp_equation} gives, schematically,
\begin{equation} \label{eq:overview_sharp_weighted_schematic}
    \begin{aligned}
    \bigl\|\jxi^2\nabla_\xi\widehat g^\sharp(t)\bigr\|_{L^2_\xi}
    \lesssim{}&
    \{\text{weighted initial data}\}
    +
    \|g^\sharp(t)\|_{H^2_x}
    +
    \jt
    \bigl\|
        \jD^2e^{it\jD}g^\sharp(t)
    \bigr\|_{L^{\infty-}_x}
    \\
    &+
    \sum_{j=0,1}
    \bigl\|s\,\jx V\partial_s^jv^\flat
    \bigr\|_{L^2_s([0,t];L^2_x)}
    +
    \bigl\|s\,\jx\jD(Vv^\flat)
    \bigr\|_{L^2_s([0,t];L^2_x)}
    + \cdots
    \end{aligned}
\end{equation}
where `$\cdots$' are similar or more favorable terms;
see Proposition \ref{prop:transference} for the exact statement.
The localization of $V$, together with the bounds \eqref{eq:overview_flat_outputs} for $v^\flat$, control all
the forcing norms on the right-hand side of \eqref{eq:overview_sharp_weighted_schematic} and yield
\begin{equation}
\label{eq:overview_sharp_weighted_final}
    \bigl\|\jxi^2\jap{\nabla_\xi}
    \widehat g^\sharp(t)\bigr\|_2
    \lesssim t^{5\delta},
\end{equation}
which is an improvement over the bootstrap assumption~\eqref{eq:overview_large_time_framework}. 

The basic idea of the transference relation leading to \eqref{eq:overview_sharp_weighted_schematic} is as follows. First, by \eqref{eq:overview_profile_boost_identity} and \eqref{equ:overview_transference_flat} we reduce the task of estimating $\nabla_\xi\widehat{g}^\sharp$ to that of estimating $\Omega_{0j}v^\sharp$. Now suppose there exists an operator $\widetilde\Omega$ such that 
    \begin{equation}\label{eq:overviewtildeOmegacommute1_alt}
        \bigl[ \widetilde\Omega,\partial_t^2+H \bigr] = 0,
    \end{equation}
    and such that for any functions $g$, $G$ and for $\gamma\geq0$,
    \begin{equation}\label{eq:vsharpOmegaoverview1_alt}
    \|\Omega_{0j}g\|_{L^2_x}\lesssim \| \widetilde\Omega g\|_{L^2_x},
\end{equation}
and 
\begin{equation}\label{eq:vsharpOmegaoverview2_alt}
     \|\jap{x}^{\gamma}\widetilde\Omega G\|_{L^2_x}\lesssim t\|\jap{x}^{\gamma}G\|_{H^1_x}+\|\jap{x}^{1+\gamma} \partial_tG\|_{L^2_x}.
\end{equation} 
Note that the estimate \eqref{eq:vsharpOmegaoverview2_alt} 
is satisfied if $\widetilde\Omega$ is replaced by $\Omega_{0j}$. In view of \eqref{eq:vsharpOmegaoverview1_alt}, it suffices to estimate $\widetilde{\Omega}v^\sharp$. In view of \eqref{eq:overviewtildeOmegacommute1_alt} and ILED for \eqref{eq:overview_sharp_equation}, we then need to estimate $\widetilde{\Omega}P_c(Vv^\flat)$ in the dual ILED norm. The desired estimate \eqref{eq:overview_sharp_weighted_schematic} then follows by using \eqref{eq:vsharpOmegaoverview2_alt}.

For the pointwise bound on $e^{it\jap{D}} g^\sharp(t)$, by the same argument as in~\eqref{eq:overview_bad_pointwise},
\begin{equation} \label{eq:overview_sharp_duhamel}
    \begin{aligned}
        \|P_cv^\sharp(t)\|_{L^\infty_x}
        &\lesssim
        \{\text{norms of data}\} + 
        \int_0^t
        \jap{t-s}^{-1}
        \|Vv^\flat(s)\|_{W^{3,1}_x}\,\ud s
        \lesssim t^{-1+\delta}.
    \end{aligned}
\end{equation}
Note that the initial data for $v^\sharp$ include the contribution of $v^\flat$. The corresponding weighted norms for these are already estimated in the analysis of $g^\flat$.

\subsubsection{Transference relations} \label{subsubsec:transference_relations_overview}

The existence of the operator $\widetilde\Omega$ satisfying \eqref{eq:overviewtildeOmegacommute1_alt}--\eqref{eq:vsharpOmegaoverview2_alt} is one of the main inputs of our companion paper \cite{LPPSS2}. The idea behind the existence of $\widetilde\Omega$ is presented in Subsection~1.3 in \cite{LPPSS2} and the full details, leading to Proposition~\ref{prop:transference}, are provided in Section~4 of \cite{LPPSS2}. For completeness, here we briefly recall this. Recall the relation~\eqref{eq:overview_Omega_Omegahat1}. Taking the Fourier transform of the relation $[\Box+1,\Omega_{0j}]=0$ shows that
\begin{equation}\label{eq:OmegaFouriercommoverview1}
    \bigl[ \partial_t^2+\jap{\xi}^2, \widehat\Omega_{0j} \bigr] = 0.
\end{equation}
Taking this as the motivation, we hope to relate $\Omega_{0j}$ to an operator similar to $\widehat\Omega_{0j}$ but defined with respect to the distorted Fourier transform $\widetilde\calF$ for $H$. The latter is the spectral transformation diagonalizing the self-adjoint operator $H$, that is, $\widetilde\calF[H g](k)= \jap{k}^2\widetilde\calF[g](k)$, where we have used $k$ to denote the independent distorted Fourier variable.\footnote{To be precise, this procedure can be carried out mode by mode and $k$ is the spectral parameter. Alternatively one can use the Lippmann-Schwinger approach. In our actual problem under equivariant symmetry we are naturally lead to an operator restricted to the first angular mode.} For simplicity of exposition we present the remainder of the argument in spatial dimension one. Here the distorted Fourier transform $\widetilde\calF$ and its inverse $\widetilde\calF^{-1}$ can be written as
\begin{equation}\label{eq:dfToverviewdef1}
    \tilde g(k)
    :=\widetilde{\calF}[g](k)=\int_{\bbR} \phi(x,k)g(x)\,\ud x,\qquad \widetilde{\calF}^{-1}[\tilde g](k)=\int_{\bbR} \overline{\phi(x,k)}\tilde g(x)\,\ud k,
\end{equation}
and satisfy $\|g\|_{L^2_x(\bbR)}=\|\tilde g\|_{L^2_k(\bbR)}$. The distorted Fourier basis functions $\phi(x,k)$ satisfy $H\phi(x,k)=\jap{k}^2\phi(x,k)$. By the variation of parameters formula
\begin{equation}
    \phi(x,k)=e^{ixk}s_+(x,k)+e^{-ixk}s_-(x,k),
\end{equation}
where for suitable functions $q_\pm(k)$,
\begin{equation}
    s_+(x,k)=q_+(k)-i\int_{-\infty}^xe^{-iyk}V(y)\phi(y,k)\,\ud y,\qquad s_-(x,k)=q_-(k)-i\int_x^\infty e^{iyk}V(y)\phi(y,k)\,\ud y.
\end{equation}
It follows that
\begin{align*}
    &\partial_x\phi(x,k)=ike^{ixk}s_+(x,k)-ike^{-ixk}s_-(x,k),\\
    &x\phi(x,k)=-i\partial_k(e^{ixk}s_+(x,k)-e^{-ixk}s_-(x,k))+i(e^{ixk}\partial_ks_+(x,k)-e^{-ixk}\partial_ks_-(x,k)).
\end{align*}
Combining these relations with \eqref{eq:dfToverviewdef1} and the inversion formula $g(x)=\widetilde\calF^{-1}[g](x)$, and integrating by parts once in $k$ gives
\begin{align}\label{eq:transferenceintro1}
    \begin{split}
    \Omega u(t,x)&=i\int_{-\infty}^\infty\big(e^{ikx}s_+(x,k)-e^{-ixk}s_-(x,k)\big)\widetilde\Omega \tilde{u}(t,k)\,\ud k\\
    &\quad+ i\int_{-\infty}^\infty\big(e^{ikx}\partial_ks_+(x,k)-e^{-ixk}\partial_ks_-(x,k)\big)\partial_t\tilde{u}(t,k)\,\ud k.
    \end{split}
\end{align}
Here $\Omega= t\partial_x+x\partial_t$  and $\widetilde{\Omega}:=\partial_t\partial_k+tk$.  We can then hope to  bound the $L^2_x$ norm of the first integral in \eqref{eq:transferenceintro1} by the $L^2_k$ norm of $\widetilde\Omega\tilde{u}$, and treat the second integral in \eqref{eq:transferenceintro1} as an error term. Moreover, by analogy with the usual Fourier transform, where differentiation in the Fourier variable corresponds to multiplication on the physical side and vice versa, we expect that the estimate in \eqref{eq:vsharpOmegaoverview2_alt} is satisfied for this choice of $\widetilde{\Omega}$. The commutativity property \eqref{eq:overviewtildeOmegacommute1_alt} is a direct consequence of the algebraic relation \eqref{eq:OmegaFouriercommoverview1}. 

In the actual two-dimensional equivariant problem, the Bessel functions $J_n(rk)$ and $Y_n(rk)$ replace the exponentials $e^{\pm ixk}$. Moreover, the Lorentz boosts change angular degree: when applied to an equivariant function of degree $n$, they generate degrees $n-1$ and $n+1$.  Correspondingly, the distorted boost operators take the form
\[
    \widetilde\Omega_n^\pm
    :=
    tk+\partial_t
    \left(
        \partial_k\pm\frac{n}{k}
    \right).
\]
The matrix structure, the interaction of different Bessel orders, and the critical $\jap{r}^{-2}$ decay of the potential make the error analysis in the two-dimensional transference identity considerably more delicate.  The required mapping estimates, together with the resulting transference theorem, are proved in \cite[Subsection~1.3 and Section~4]{LPPSS2}.

\subsubsection{Variable-coefficient quadratic normal forms} \label{subsec:overview_normal_forms}

We now describe the main nonlinear analysis for the flat component satisfying \eqref{eq:overview_flat_equation}. 
Our treatment of the quadratic terms is based on a variable-coefficient normal-form method.
We recall that the use of the (classical) normal form method for flat Klein-Gordon equations with constant-coefficient quadratic nonlinearities goes back to \cite{Shatah1985}. 
For the more recent development of variable-coefficient normal forms we refer to \cite{GP20,LLS20,DelMas20} and references therein. 
In particular, in this work we draw on ideas from \cite{GP20} for the implementation of variable coefficient normal forms with a principal value singularity.

The two non-localized quadratic terms retained explicitly in the scalar model \eqref{eq:overview_flat_equation} present different difficulties. Upon taking the Fourier transform, the vortex coefficient breaks the convolution structure, and its Fourier transform introduces an additional frequency multiplier and a principal value singularity.
By contrast, the derivative quadratic nonlinearity retains the constant-coefficient quadratic Klein-Gordon phase, which never vanishes, but its multiplier carries one input derivative. 
In this overview we disregard the treatment of this latter
as it is more standard.

Let $f$ and $g^\flat$ denote the profiles of $u$ and $v^\flat$, as
in \eqref{eq:overview_profiles} and \eqref{eq:gflatsharpdefoverview}. We also use the notation $f^+:=f$ and $f^{-}:=\overline{f}$. For the contribution of the vortex $\underline{\Phi}_{\real}u^2$ to $g^\flat$, a representative integral is
\begin{equation} \label{eq:overview_type1_profile}
    \begin{aligned}
        \widehat{\calQ}^{(1)}(t,\xi)
        =
        \sum_{\kappa_1,\kappa_2\in\{\pm\}}
        \int_t^T\int_{\bbR^2}\int_{\bbR^2}
        \frac{e^{is\Psi^{(1)}_{\kappa_1\kappa_2}
        (\xi,\eta,\sigma)}}{2i\jxi}
        \widehat{\underline{\Phi}}_{\real}(\sigma) \widehat f^{\kappa_1}(s,\eta)
        \widehat f^{\kappa_2}(s,\xi-\eta-\sigma)
        \,\ud\sigma\,\ud\eta\,\ud s,
    \end{aligned}
\end{equation}
where
\begin{equation} \label{eq:overview_type1_phase}
    \Psi^{(1)}_{\kappa_1\kappa_2}(\xi,\eta,\sigma)
    :=
    -\jxi+\kappa_1\jap{\eta}
    +\kappa_2\jap{\xi-\eta-\sigma}.
\end{equation}
The main feature of this term is the singular Fourier transform of the vortex coefficient,
\begin{equation} \label{eq:overview_vortex_transform}
    \widehat{\underline{\Phi}}_{\real}(\sigma)
    =
    -i\,\varphi_{\leq1}(\sigma)\,
    \mathrm{p.v.}\frac{\sigma_1}{|\sigma|^3}
    +F_{\real}(\sigma),
\end{equation}
where the inverse Fourier transform of $F_{\real}$ is spatially localized, see Lemma~\ref{lem:FTvortex_re_im}.

We note that at $\sigma=0$, the phase in \eqref{eq:overview_type1_phase} becomes the usual quadratic Klein-Gordon phase
\[
    -\jxi+\kappa_1\jap{\eta}
    +\kappa_2\jap{\xi-\eta},
\]
which does not vanish. In a sufficiently small 
neighborhood of $\sigma=0$, the phase then has a positive lower bound for bounded frequencies.  
We can therefore divide by the phase and integrate by parts in time 
in a suitable region containing the principal-value singularity. 
Outside this suitably chosen neighborhood of $\sigma=0$, 
the phase may vanish, giving rise to time resonances. But there the principal value multiplier is regular in the sense that it furnishes spatial localization, which permits the use of ILED-type estimates.

We now discuss another relevant aspect connected to our normal form procedures, that is, the choice of final data at time $t=T$
for the flat equation \eqref{eq:overview_flat_equation}.
The normal form identity can be schematically written as
\begin{align}\label{eq:overview_time_normal_form}
        \int_t^T \int_{\bbR^2}\int_{\bbR^2}e^{is\Psi}F(s,\xi,\eta,\sigma) \, 
        \,\ud\sigma\,\ud\eta\,\ud s
        &=
        \left[
        \int_{\bbR^2}\int_{\bbR^2}
        \frac{e^{is\Psi}}{i\Psi}
        F(s,\xi,\eta,\sigma) \,\ud\sigma\,\ud\eta
        \right]_{s=t}^{s=T}
        \\
        &\quad - \int_t^T \int_{\bbR^2}\int_{\bbR^2}
        \frac{e^{is\Psi}}{i\Psi}
        \partial_s F(s,\xi,\eta,\sigma) \,\ud\sigma\,\ud\eta\,\ud s.
\end{align}
Let $\calB[f,f](s)$ denote the total quadratic boundary correction at
time $s$, obtained by summing, over all sign choices, the boundary
operators arising from every quadratic term to which
\eqref{eq:overview_time_normal_form} is applied.  After subtracting the
boundary correction at time $t$, the normal-form identity gives,
schematically,
\[
    g^\flat(t)-\calB[f,f](t)
    =
    g^\flat(T)-\calB[f,f](T)
    +
    \{\text{bulk terms
    }\}.
\]
The terminal difference
$g^\flat(T)-\calB[f,f](T)$ is independent of $t$,
although it remains a function of the frequency.  It would therefore
generate a homogeneous flat Klein-Gordon evolution, for which one has
only the standard linear decay estimate.  To eliminate this term and
retain the improved decay coming from the time-integrated terms
(see the estimates anticipated in \eqref{eq:overview_flat_outputs}), we choose
the final data of $v^\flat$ so that 
\begin{equation} \label{eq:overview_flat_final_data}
    g^\flat(T)=\calB[f,f](T).
\end{equation}
With this choice, the renormalized profile
\begin{align} \label{eq:overview_renormalized_profile}
    g^\flat_{\mathrm{ren}}
    :=
    g^\flat-\calB[f,f] \qquad \hbox{satisfies} \qquad
    g^\flat_{\mathrm{ren}}(T)=0,
\end{align}
and the equation for $g^\flat_{\mathrm{ren}}$, integrated backward from $T$, is driven by the cubic terms arising when $\partial_s$ in \eqref{eq:overview_time_normal_form} falls on one of the inputs.
The weighted estimate is proved for $g^\flat_{\mathrm{ren}}$, while $\calB[f,f]$ is estimated separately. 

\medskip
\subsubsection{Weighted estimates for the cubic terms}
\label{subsubsec:overview_cubic_commutation}
Substituting the profile equation for $\partial_sf$ into the bulk term in \eqref{eq:overview_time_normal_form} produces cubic integrals.
Oversimplifying, these terms can be thought of as the Duhamel integrals associated to terms of the form, say, $\underline{\Phi}_{\real} \cdot u \cdot u \cdot u$.
The main difficulty in the weighted estimate 
is when $\nabla_\xi$ hits the phase, producing a growing factor $is\nabla_\xi\Psi$. A direct cubic estimate placing two inputs in $L^\infty$ and the remaining
input in $L^2$ leaves
\[
    s \, \| u(s) \|_{L^\infty}^2 \lesssim s^{-1+2\delta},
\]
which is not integrable in time; therefore, oscillation must still be exploited. 
We then note that, for the cubic phase with one extra `coefficient frequency' $\sigma$ (cf. the quadratic phase in \eqref{eq:overview_type1_profile}-\eqref{eq:overview_type1_phase}),
that is,
\begin{align}
        \Psi_{\iota_1\iota_2\iota_3}
        := &
        -\jxi+\iota_1\jap{\xi_1}
        +\iota_2\jap{\xi_2} + 
        \iota_3 \jap{\xi-\xi_1-\xi_2-\sigma},
    \end{align}
we have the exact `commutation identity'
\begin{equation} \label{eq:overview_cubic_commutation}
    \begin{aligned}
        &
        \jxi\nabla_\xi\Psi_{\iota_1\iota_2\iota_3} =
        -\sigma 
        -\left(\iota_1\jap{\xi_1}\nabla_{\xi_1}
        +\iota_2\jap{\xi_2}\nabla_{\xi_2}
        \right) \Psi_{\iota_1\iota_2\iota_3}
        -
        \iota_3 \frac{\xi-\xi_1-\xi_2-\sigma}
        {\jap{\xi-\xi_1-\xi_2-\sigma}}
        \Psi_{\iota_1\iota_2\iota_3}.
    \end{aligned}
\end{equation}
After integration by parts in $\xi_1$ and $\xi_2$, the first two input-frequency derivatives in \eqref{eq:overview_cubic_commutation} are transferred to weighted input profiles, thereby removing the explicit factor $s$ from
these terms. A representative cubic contribution arising from 
this argument, when all inputs are given by the bad component $w$, can be estimated schematically by 
\[
    \int_t^T
    \|w(s)\|_{L^\infty_x}^2
    \bigl\|
        \jxi^2\nabla_\xi \widehat{h}(s)
    \bigr\|_{L^2_\xi}
    \,\ud s
    \lesssim
    \int_t^T
    s^{-2+2\delta}s^{1/2}\,\ud s
    \lesssim
    t^{-1/2+2\delta},
\]
consistent with the claimed \eqref{eq:overview_flat_outputs}.
The term containing the cubic phase on the right-hand side of \eqref{eq:overview_cubic_commutation} is converted into a time derivative using
\[
    i\Psi_{\iota_1\iota_2\iota_3}
    e^{is\Psi_{\iota_1\iota_2\iota_3}}
    =
    \partial_s
    \left(
    e^{is\Psi_{\iota_1\iota_2\iota_3}}
    \right).
\]
Integration by parts in $s$ then handles these contributions.
Finally, the extra factor of $\sigma$ in \eqref{eq:overview_cubic_commutation}, which is paired with the unfavorable growing factor of $s$, 
is used to mitigate the $\mathrm{p.v.}(\sigma_1/|\sigma|^3)$ singularity produced by the vortex $\underline{\Phi}_{\real}$ according to \eqref{eq:overview_vortex_transform};
the resulting kernel is almost in $L^2$ and this
suffices to gain enough localization and decay to overcome the $s$ factor and close the bootstrap for the good component.

Having obtained \eqref{eq:overview_flat_outputs},
the weighted estimate \eqref{eq:overview_sharp_weighted_final} 
on the sharp component then follows,
and these, together with \eqref{eq:overview_sharp_duhamel}
close the main bootstrap on the norms in \eqref{eq:overview_large_time_framework}.

\subsubsection{An analogy: The flat-sharp approach using the classical vector field method} \label{subsec:overview_vf}

In the context of the classical vector field method, and in view of \eqref{eq:overview_profile_boost_identity} and \eqref{equ:overview_transference_flat}, estimating $ \bigl\|\langle\xi\rangle^2\langle\nabla_\xi\rangle\widehat{f}(t)\bigr\|_{L_\xi^2}$ is analogous to estimating $\|\Omega u\|_{H^2_x}$, where $\Omega$ denotes any of the Lorentz boost
vector fields $\Omega_{0j}$. Here we informally describe the idea of the flat-sharp  approach by explaining how to estimate $\|\Omega v\|_{L^2_x}$ under the bootstrap assumptions \eqref{eq:overview_large_time_framework}--\eqref{eq:overview_large_time_framework_ILED}, and for a simpler cubic nonlinearity. 

To further simplify the discussion we assume here that $H$ has no discrete spectrum. That is, we provisionally replace the first equation in~\eqref{eq:overview_good_bad_model} by the simpler equation
\begin{equation}\label{eq:overview_good_bad_model_alt}
  (\partial_t^2+H)v= u^3,
  \qquad (v,\partial_tv)|_{t=0}=(u,\partial_tu)|_{t=0},\qquad \sigma_{\mathrm{d}}(H)= \text{\O}.
\end{equation}
If we could commute $\Omega$ with \eqref{eq:overview_good_bad_model_alt}, then our bootstrap assumptions \eqref{eq:overview_large_time_framework}--\eqref{eq:overview_large_time_framework_ILED} on $u$ would easily allow us to estimate $\|\Omega v\|_{L^2_x}$. To deal with the fact that $\partial_t^2+H$ and $\Omega$ do not commute, we introduce the flat-sharp decomposition $v=v^\flat+v^\sharp$. The flat part $v^\flat$ solves the same Klein-Gordon equation as $v$ but with the operator replaced by the flat Klein-Gordon operator $\Box+1$, and zero data at the final bootstrap time $T$:\footnote{For quadratic nonlinearities, the use of normal forms requires a different choice of final data. See Subsection~\ref{subsec:overview_normal_forms}.}
\begin{equation}\label{eq:flatsharpmodel2_alt}
    (\Box+1)v^\flat = u^3,\qquad v^\flat(T)=\partial_tv^\flat(T)=0.
\end{equation}
It follows that
\begin{equation}\label{eq:flatsharpmodel3_alt}
    (\partial_t^2+H)v^\sharp = -V v^\flat ,\qquad v^\sharp(0)=v(0)-v^\flat(0),\qquad \partial_tv^\sharp(0)=\partial_tv(0)-\partial_tv^\flat(0).
\end{equation}
The idea is that since $\Omega$ commutes with the equation for $v^\flat$, we can prove decay estimates for the flat component $v^\flat$. Moreover, since $v^\flat$ has vanishing final data, the decay estimates for $v^\flat$ are dictated by the source term in \eqref{eq:flatsharpmodel2_alt}, which is an improvement over the linear decay rate of free Klein-Gordon waves. For instance,  
\begin{equation}\label{eq:flatOmegaboundoverview1_alt}
    \|\Omega v^\flat(t)\|_{L^2_x}\lesssim \int_t^T \|u(s)\|_{L^\infty_x}^2\|\Omega u(s)\|_{L^2_x}\,\ud s
    \lesssim t^{-\frac{1}{2}+2\delta},
\end{equation}
where in analogy with the weighted estimate in  \eqref{eq:overview_large_time_framework}, we have assumed that $\|\Omega u(s)\|_{L^2_x}\lesssim s^{1/2}$.  As explained in Subsection~\ref{subsubsec:overview_bad_component}, the origin of this growth is the slow decay of the internal mode component. 
In analogy with \eqref{eq:overview_weighted_dispersive}, the estimate \eqref{eq:flatOmegaboundoverview1_alt} corresponds to the  pointwise bound $$\|v^\flat(t)\|_{L^\infty_x}\lesssim \jap{t}^{-3/2+3\delta}.$$ 

Turning to \eqref{eq:flatsharpmodel3_alt}, in view of the spatial decay of $V$ and the time decay of $v^\flat$, the right-hand side of \eqref{eq:flatsharpmodel3_alt} is spatially localized and satisfies 
$\jap{t}^{-3/2+3\delta}$ pointwise decay.  To overcome the non-commutativity of $\Omega$ and $\partial_t^2+H$, suppose 
\begin{itemize}
    \item integrated local energy decay holds for $\partial_t^2+H$,
    \item there exists an operator $\widetilde\Omega$ satisfying \eqref{eq:overviewtildeOmegacommute1_alt}--\eqref{eq:vsharpOmegaoverview2_alt}.
\end{itemize}
Estimate \eqref{eq:vsharpOmegaoverview1_alt} is the transference relation analogous to the more precise Proposition~\ref{prop:transference}.  Applying \eqref{eq:overviewtildeOmegacommute1_alt}, \eqref{eq:vsharpOmegaoverview1_alt}, and \eqref{eq:vsharpOmegaoverview2_alt} to \eqref{eq:flatsharpmodel3_alt} we see that if the potential $V$ decays sufficiently fast, then the right-hand side can be estimated in the dual ILED norm by
\begin{equation}\label{eq:flatOmegaboundoverview3_alt}
    \bigl\| \jx \widetilde{\Omega} \bigl( V v^\flat \bigr) \bigr\|_{L^2_s([0,t]; L^2_x(\bbR^2))} \lesssim \jap{t}^{\delta'}
\end{equation}
for some $\delta'>0$. This is much better than the growth $t^{1/2}$ of the weighted estimate in our bootstrap assumption~\eqref{eq:overview_large_time_framework}. 
The construction of $\widetilde{\Omega}$ is discussed in Subsection~\ref{subsubsec:transference_relations_overview}.

\subsection{Implementation for the full system} \label{subsec:implementation_full_system}

We now explain how some of the preceding arguments are implemented for the full system and how the three papers in the series enter the proof of Theorem~\ref{thm:main}.

\subsubsection{The bootstrap framework}
The full problem is written in symmetry-reduced equivariant variables as in \eqref{syscompact}-\eqref{syscompact2}.
To use the standard Fourier transform on $\bbR^2$ 
to estimate the radiation $\bmu$,
we lift the equations to the plane by conjugating the radiation
equation by $e^{i\theta}$ and write
\begin{align} \label{eq:imp_eithetau}
    e^{i\theta}\bmu
    =
    \bmu_\real+i\bmu_\imag,
\end{align}
where $\bmu_\real$ and $\bmu_\imag$ are real-valued vectors on $\bbR^2$.
We also rewrite radial derivatives in terms of rectangular derivatives and express the nonlinearities in terms of the components of $\bmu_\real$ and $\bmu_\imag$.
The flat Klein-Gordon profiles associated to \eqref{eq:imp_eithetau} are
\begin{equation}
    \bmf_\ast(t)
    :=
    e^{-it\jD}(2i\jD)^{-1}
    (\pt+i\jD)\bmu_\ast(t),
    \qquad
    \ast\in\{\real,\imag\},
\end{equation}
so that
\[
    \bmu_\ast(t)
    =
    2\Re\bigl(e^{it\jD}\bmf_\ast(t)\bigr).
\]
For the two internal mode components, we similarly define their profiles by
\begin{equation}
    Z_j(t)
    :=
    e^{-it\lambda}(2i\lambda)^{-1}
    (\pt+i\lambda)z_j(t),
    \quad
    j=1,2,
    \qquad
\end{equation}
and set
\begin{equation}
    \frakz(t):=|Z_1(t)|^2+|Z_2(t)|^2.
\end{equation}
The bootstrap assumptions concern $\frakz$ and the weighted, dispersive,
high Sobolev, and integrated local energy norms of $\bmf_\real$ and
$\bmf_\imag$. See the schematic \eqref{eq:overview_large_time_framework},
and the actual assumptions \eqref{equ:bootstrap_assumption_z} and
\eqref{equ:bootstrap_assumption_weighted_f}-\eqref{equ:bootstrap_assumption_HN_ILED_f}.
The piecewise form of the a priori bounds reflects the different
behavior of the quadratic internal mode forcing of the bad part before and after the damping time $t=\varepsilon^{-2}$.

\subsubsection{Spectral and linear input from the companion works \cites{LPPSS1}{LPPSS2}}

Our companion papers \cites{LPPSS1}{LPPSS2} establish the spectral results and linear estimates that underpin the nonlinear analysis in the present work. In this subsection, we review the aspects of those results that are most relevant here.

In our first companion paper \cite{LPPSS1} we exploit the Bogomolny structure
\begin{equation}
\label{eq:overviewBogomolny1}
    \bfL=\calB^\ast\calB,
\end{equation}
with diagonal super-symmetric partner operator
\begin{equation}
\label{eq:overviewBogomolny2}
    \calL
    :=
    \calB\calB^\ast
    =
    \pmat{\calL_1&0\\0&\calL_2},
\end{equation}
where $\calL_1$ and $\calL_2$ are scalar, singular, radial Schr\"odinger operators; see \cite[Subsection~3.1]{LPPSS1}.  
Recall the block-diagonal structure of the linearized operator $\bfM$ in terms of $\bfL$ from \eqref{sysN1compact}. This allows us to reduce the spectral analysis
for the matrix operator $\bfL$ to that of the scalar operators $\calL_1$ and $\calL_2$.
Using this reduction, together with analytic arguments and rigorous
numerics, in \cite{LPPSS1} we prove that the continuous spectrum of $\bfM$ is
purely absolutely continuous and equals $[1,\infty)$, that the threshold
$1$ is neither an eigenvalue nor a resonance, and that $\bfM$ has a unique
gap eigenvalue $\lambda^2\in(0,1)$ of multiplicity two.  Its eigenspace
gives the two internal modes, and $2\lambda>1$, so their second harmonic
lies in the continuous spectrum.  More precisely,
\[
    \lambda^2\in[0.777471875,\,0.777473750].
\]
In the same work we give a computer-assisted proof, with rigorous error analysis, of the fact that the relevant Fermi Golden Rule non-degeneracy holds. These spectral conclusions are summarized in Proposition~\ref{prop:spectral_fourier}.

We now briefly describe the numerical part of~\cite{LPPSS1}. The calculation
begins with quantitative estimates for the vortex profile $(U,a)$. Near
$r=0$, the vortex and the spectral solutions are represented by Frobenius
series, with explicit bounds for the terms left out of the finite sums. For
large $r$, comparison arguments and Volterra equations give bounds for the
vortex and for the solutions determined by their prescribed behavior at
infinity. On the bounded interval between these two regions, interval
arithmetic produces intervals containing the exact solutions and their
derivatives. The solutions constructed from the two ends are then compared
at a common radius. For $\calL_2$, a threshold Wronskian bounded away from
zero excludes a resonance, while opposite signs of the eigenvalue Wronskian
at the endpoints of the interval displayed above prove the existence of an
eigenvalue inside that interval. An eigenvalue-counting estimate, whose
integral is also bounded rigorously, rules out a second eigenvalue. The
corresponding assertions for $\calL_1$ are proved analytically with a numerically assisted check that a potential $V_1>1$.

The verification of the Fermi Golden Rule reduces to evaluating explicit radial integrals
involving the internal mode, the vortex, and the scalar generalized
eigenfunctions. The contributions near the origin are estimated from the
Frobenius series, the integrals over the bounded region are treated by
interval arithmetic, and the contributions at infinity are controlled by
the decay estimates and the Volterra equations. This proves
\begin{equation}\label{eq:FGRnumoverview1}
    D_1<0,\qquad D_2<0,\qquad D_{12}<0,
\end{equation}
for the coefficients appearing in \eqref{eq:overview_full_internal_mode_system} below. 
These inequalities are the non-degeneracy input used to obtain the dissipative terms in the equations for the internal modes.

In our second companion paper \cite{LPPSS2} we develop the distorted Fourier
theory for $\bfM$, which provides a spectral representation of the
perturbed Klein-Gordon evolution on the continuous spectral subspace.  The
reduction
\eqref{eq:overviewBogomolny1}-\eqref{eq:overviewBogomolny2} relates this
construction to the distorted Fourier theories of $\calL_1$ and
$\calL_2$. As a consequence, we establish dispersive, integrated local energy, mapping,
and transference estimates for the matrix operator $\bfM$ that are summarized in
Subsections~\ref{sec:linear_decay_estimates}--\ref{sec:mapping_properties},
see in particular Proposition~\ref{prop:transference}.  
Note that the final results that are used in this article are stated without reference to the distorted Fourier transform. 
Here we only point out the critical inverse square decay of the potential in $\bfM$ that makes the development of the distorted Fourier theory in \cite{LPPSS2} especially challenging.

\subsubsection{The radiation and temporal component}\label{sssec:ILED}

Recall the system \eqref{equ:system_evol_equations_u_z}.
The good-bad decomposition, analogous to \eqref{eq:overview_good_bad_model}, takes the form
\begin{equation}
\label{eq:overview_full_good_bad_decomposition}
    \bmu=\bmv+\bmw,
    \qquad
    (\pt^2+\bfM)\bmv=\bfP_c\bmN_c,
    \qquad
    (\pt^2+\bfM)\bmw=\bfP_c\bmN_d.
\end{equation}
As already described, the good part $\bmv$ inherits the initial radiation data, while
the bad part $\bmw$ has zero initial data. In the full system, every term in $\bmN_d$ contains either an
internal mode factor or the bad temporal component $\eta_{0,d}$, which will be defined below, and that can be treated using spatial localization.  
Apart from the leading source with two internal mode factors considered in \eqref{eq:overview_fgr_radiation}, the remaining terms
contain either a radiation factor or $\eta_{0,d}$, both of which decay faster than the internal modes. They are therefore controlled
perturbatively and do not affect the $t^{1/2}$ weighted growth or the
pointwise decay of $\bmw$.  

The good part $\bmv$ of the radiation is decomposed further into a flat and a sharp component 
\begin{equation}
    \bmv
    =
    \bfP_c\bmv^\flat+\bfP_c\bmv^\sharp.
\end{equation}
After conjugation by $e^{i\theta}$, the real and imaginary parts of
$\bmv^\flat$ satisfy flat Klein-Gordon equations driven by the
corresponding parts of $e^{i\theta}\bmN_c$.  Their final profiles are
chosen to equal the quadratic boundary corrections produced by the
normal-form analysis. The sharp component then satisfies
\begin{equation}
    (\pt^2+\bfM)\bfP_c\bmv^\sharp
    =
    -\bfP_c(\bfV\bmv^\flat),
    \qquad
    |\bfV(r)|\lesssim\jap{r}^{-2}.
\end{equation}
The analysis of non-localized quadratic terms is then confined to the flat
component, while the sharp component is estimated using the transference
and integrated local energy estimates for its localized forcing.

The temporal component is decomposed at the level of its elliptic equation
\eqref{equ:eta0_equation}:
\begin{equation}
    \eta_0=\eta_{0,c}+\eta_{0,d},
    \qquad
    (-\Delta+U^2)\eta_{0,c}=\calN_{0,c},
    \qquad
    (-\Delta+U^2)\eta_{0,d}=\calN_{0,d};
\end{equation}
see \eqref{equ:eta0c_equation} and \eqref{equ:eta0d_equation}.
Here $\calN_{0,c}$ contains no internal mode variables or factors
$\eta_{0,d}$, while $\calN_{0,d}$ contains the remaining terms.  
Except for the terms carrying $\eta_{0,d}$  all other terms in $\calN_{0,d}$ are spatially localized. The terms with $\eta_{0,d}$ in $\calN_{0,d}$ are perturbative in the elliptic estimate for $\eta_{0,d}.$ The off-diagonal exponential
decay of the Green kernel of $(-\Delta+U^2)^{-1}$ then implies that
$\eta_{0,d}$ is also localized.
Unlike the hyperbolic decomposition, this
splitting is performed at each fixed time and involves no choice of initial
data.

We further decompose the good part of the temporal component as
\begin{equation}
    \eta_{0,c}
    =
    \eta_{0,c}^\flat+\eta_{0,c}^\sharp, \qquad \ \hbox{where} \ \qquad
    (-\Delta+1)\eta_{0,c}^\flat=\calN_{0,c},
    \qquad \
    (-\Delta+U^2)\eta_{0,c}^\sharp
    =
    (1-U^2)\eta_{0,c}^\flat.
\end{equation}
See \eqref{equ:eta0c_flat_elliptic_equ} and
\eqref{equ:eta0c_sharp_elliptic_equ}.  Here $-\Delta$ is the full
Laplacian on $\bbR^2$, with $\eta_{0,c}$ viewed as a radial function.
Since
\begin{equation}
    \widehat{\eta}_{0,c}^\flat(\xi)
    =
    \jxi^{-2}\widehat{\calN}_{0,c}(\xi),
\end{equation}
the flat component of the good part of $\eta_0$ is suitable for analysis in the usual Fourier
variables, while the sharp component of the good part is generated by a localized coefficient.  
The quadratic terms involving $\pt\eta_{0,c}^\flat$ on the right-hand sides of \eqref{sysN1}
can then be treated by the same normal form analysis as the other non-localized quadratic terms,
while the other contributions are localized or at least cubic.

Altogether,
\begin{equation}
    \bmu
    =
    \bfP_c\bmv^\flat
    +
    \bfP_c\bmv^\sharp
    +
    \bfP_c\bmw.
\end{equation}
The estimates for these three components yield, for
$t\geq\varepsilon^{-2}$ and $\ast\in\{\real,\imag\}$,
\begin{equation}
\label{eq:overview_full_radiation_outputs}
    \bigl\|
        \jxi^2\jap{\nabla_\xi}\widehat{\bmf}_\ast(t)
    \bigr\|_{L^2_\xi}
    \lesssim
    t^{\frac12},
    \qquad
    \bigl\|e^{it\jD}\bmf_\ast(t)\bigr\|_{L^\infty_x}
    \lesssim
    t^{-1+\delta}.
\end{equation}
The bad part dictates the weighted growth, but its pointwise decay is proved directly from its localized source, as in \eqref{eq:overview_bad_pointwise}. 
The good part satisfies stronger weighted estimates, which give the required decay of its flat and sharp components as in Subsection~\ref{subsec:overview_flat_sharp}.  
Combining the three components also gives the high Sobolev and integrated local energy bounds, as well as the finer full bootstrap estimates for shorter times $t \leq \varepsilon^{-2}$.

\subsubsection{The internal modes and closure of the bootstrap}\label{sssec:im}

As discussed in Subsection~\ref{subsubsec:FGRmodel}, Section~\ref{sec:FGR1} uses the constrained Hamiltonian formulation to
derive the analogue of the scalar ODE
\eqref{eq:overview_scalar_fgr_ode} for the two internal mode profiles $Z_j := (2i\lambda)^{-1} e^{-it\lambda} (\pt+i\lambda)z_j$, $j=1,2$. This leads to a coupled system of ODEs, which is initially analyzed assuming bounds
on the remainder terms that are proved later in
Section~\ref{secIM2}.  
More precisely, the profiles $Z_1,Z_2$ satisfy
\begin{equation}
\label{eq:overview_full_internal_mode_system}
\left\{
\begin{aligned}
\dot Z_1={}&
    (D_1+iH_1)|Z_1|^2Z_1
    +2(D_{12}+iH_{12})|Z_2|^2Z_1
    +(C_{12}+iC_{12}')\overline{Z_1}Z_2^2
    +R_1,
\\
\dot Z_2={}&
    (D_2+iH_2)|Z_2|^2Z_2
    +2(D_{12}+iH_{12})|Z_1|^2Z_2
    +(C_{12}+iC_{12}')\overline{Z_2}Z_1^2
    +R_2.
\end{aligned}
\right.
\end{equation}
Here $R_1$ and $R_2$ denote remainders analogous to $\calR_Z$ in
\eqref{eq:overview_scalar_fgr_ode}. The coefficients  $D_1,D_2,D_{12}$, and $C_{12}$ together play the role of the scalar damping coefficient $\Gamma$ in \eqref{eq:overview_scalar_fgr_coefficient}, while $H_1,H_2,H_{12}$, and $C_{12}'$ play the role of the phase-shifting coefficient $\Lambda$ in \eqref{eq:overview_scalar_fgr_ode}.  

Since $\bfM$ is matrix-valued and the internal eigenspace is spanned by $\bmY_1,\bmY_2$, the above coefficients involve self- and cross-interactions 
(see Lemmas~\ref{lem1} and~\ref{lem3} and the formulas
\eqref{FGRODE14}--\eqref{FGRODE14'}) and we can verify that
\begin{equation}
    D_1,D_2,D_{12}\leq0,
    \qquad
    |C_{12}|\leq\sqrt{D_1D_2}.
\end{equation}
The contribution of the imaginary coefficients drop out from the equation for $\frakz=|Z_1|^2+|Z_2|^2$, and the real coefficients give 
\begin{equation}
\label{eq:overview_full_internal_mode_inequality}
    -\Gamma_0\frakz^2(t)
    \leq
    \frac{\ud}{\ud t}\frakz(t)
    -
    2\Re\bigl(R_1\overline{Z_1}\bigr)
    -
    2\Re\bigl(R_2\overline{Z_2}\bigr)
    \leq
    -\Gamma_1\frakz^2(t),
    \qquad
    0 \leq \Gamma_1 \leq \Gamma_0,
\end{equation}
where
\begin{align*}
    \Gamma_0 &:= 4 \max \bigl\{ |D_1|, |D_2|, |D_{12}| \bigr\} \geq 0, \\
    \Gamma_1 &:= 2 \min \bigl\{ |D_{12}|,\sqrt{D_1D_2} \bigr\} \min \big\{ \sqrt{D_1/D_2}, \sqrt{D_2/D_1} \bigr\} \geq 0.
\end{align*}
The Fermi Golden Rule non-degeneracy, verified in \cite[Proposition~4.3]{LPPSS1} using rigorous numerics, gives $\Gamma_1>0$, see \eqref{eq:FGRnumoverview1}.

Section~\ref{secIM2} establishes the remainder bounds needed to apply
Proposition~\ref{propODEdamp}; see Proposition~\ref{propODERbounds}.
Consequently,
\begin{equation}
    \frakz(t) \lesssim
    \frac{\varepsilon^2}{1+\Gamma_1\varepsilon^2t},
\end{equation}
with a comparable lower bound as in \eqref{equ:decay_zed_lower_bound_mainthm} in the case $\frz(0) = \varepsilon^2$,
as stated in Theorem \ref{thm:main}. Finally, this decay controls the internal mode forcing in the radiation equations and, together with the radiation estimates obtained above, allows us to improve all the bootstrap assumptions and complete the proof of Theorem \ref{thm:main}.

\section{Setting up the Analysis} \label{sec:setup}

In this section we prepare the system of evolution equations \eqref{equ:system_evol_equations_u_z}, \eqref{equ:system_evol_equations_u_z_part2} for the radiation term $\bmu(t)$, the temporal component $\eta_0(t)$, and the internal mode components $z_1(t)$, $z_2(t)$ for the nonlinear analysis in the remainder of this paper.

\subsection{Good-bad decomposition} \label{subsec:good_bad_decomposition}

We will adopt a space-time resonances method approach to establish decay for the radiation term.
To this end we first decompose the radiation term as well as the temporal component into a ``good'' part and a ``bad'' part, following \cite{LP22}.
The ``bad'' part is driven solely by nonlinear terms containing at least one internal mode component, which are thus spatially localized, whereas the ``good'' part receives feedback only from nonlinear terms involving radiation term components. As a consequence, we will see that the ``good'' part is amenable to more favorable weighted energy estimates, which motivates the terminology. Instead, the weighted energy bounds for the ``bad'' part are strongly growing. Nevertheless, we can obtain acceptable dispersive decay estimates for the ``bad'' part by exploiting the spatial localization of its nonlinear terms.

We introduce the notation
\begin{equation}
    \alpha = \alpha_d + u_1, \quad \zeta=\zeta_d+u_2, \quad \beta=\beta_d+u_3, \quad \mu=\mu_d+u_4
\end{equation}
with
\begin{align*}
    \alpha_d := z_1 Y_1, \quad \zeta_d := z_1 Y_2, \quad \beta_d := z_2 Y_1, \quad \mu_d := z_2 Y_2.
\end{align*}
We first lay out the good-bad decomposition $\calN_0 = \calN_{0,c} + \calN_{0,d}$ for the nonlinearities in the elliptic equation for $\eta_0$.
To this end we define the ``good'' part of the temporal component $\eta_{0,c}$ via the nonlinear equation
\begin{equation} \label{equ:eta0c_equation}
    (-\Delta + U^2) \eta_{0,c} = \calN_{0,c}
\end{equation}
with
\begin{equation}
    \calN_{0,c} := u_1 \partial_t u_3 - u_3 \partial_t u_1 - 2U\eta_{0,c} u_1 - \eta_{0,c} (u_1^2+u_3^2). 
\end{equation}
Then the ``bad'' part of the temporal component $\eta_{0,d} := \eta_0 - \eta_{0,c}$ satisfies
\begin{equation} \label{equ:eta0d_equation}
    (-\Delta + U^2) \eta_{0,d} = \calN_{0,d},
\end{equation}
where
\begin{align*}
    \calN_{0,d} &:= \alpha \partial_t \beta_d - \beta\partial_t\alpha_d + \alpha_d \partial_t u_3 - \beta_d\partial_t u_1 -2U \eta_{0,d} \alpha - 2U\eta_{0,c} \alpha_{d} \\
     &\quad \quad - \eta_{0,d} \bigl(\alpha^2+\beta^2\bigr) -\eta_{0,c} \bigl(\alpha_d^2 + 2\alpha_d u_1 + \beta_d^2 + 2\beta_d u_3\bigr).
\end{align*}
Importantly, the ``bad'' part of the temporal component $\eta_{0,d}(t)$ is spatially localized as a consequence of both the localization of all nonlinear terms involving an internal mode component and the off-diagonal exponential decay of the Green's function of $-\Delta+U^2$; see Proposition~\ref{prop:eta0d_bounds}.

Now we are in the position to lay out the good-bad decomposition for the nonlinearities of the radiation term
\begin{equation}
    \calN_{\ast} = \calN_{\ast, c} + \calN_{\ast, d}, \quad \ast \in \{ \alpha, \zeta, \beta, \mu \},
\end{equation}
where the ``good'' parts of the nonlinearities are given by
\begin{align}\label{equ:N_cd}
    \begin{split}
        \calN_{\alpha,c} &:= - \frac32 U u_1^2 + \frac12 U u_3^2 - \frac12 (u_1^2+u_3^2)u_1 - 2 u_4 \pr u_3 + 2 b u_2 u_1 - U (u_2^2 + u_4^2) - (u_2^2 + u_4^2) u_1 \\
        &\quad \quad - 2\eta_{0,c} \pt u_3 - (\pt \eta_{0,c}) u_3 + U \eta_{0,c}^2 + \eta_{0,c}^2 u_1, \\
        \calN_{\zeta,c} & := b (u_1^2 + u_3^2) - 2 U u_2 u_1 - u_2 (u_1^2 +u_3^2), \\   
        \calN_{\beta,c} & := - 2U u_1u_3 - \frac12 (u_1^2+u_3^2)u_3 + 2u_4 \pr u_1 + 2 b u_2 u_3 - (u_2^2+u_4^2) u_3 \\
        &\quad \quad + 2\eta_{0,c} \pt u_1 + U (\pt \eta_{0,c}) + (\pt \eta_{0,c}) u_1 + \eta_{0,c}^2 u_3, \\ 
        \calN_{\mu,c} & := -\pr \pt \eta_{0,c} - u_1 \pr u_3 + u_3 \pr u_1 - 2U u_1 u_4 - (u_1^2+u_3^2) u_4,
    \end{split}
\end{align}
and the ``bad'' parts are then defined as
\begin{align*}
    \calN_{\ast,d} := \calN_{\ast}-\calN_{\ast,c},  \quad \ast \in \{ \alpha, \zeta, \beta, \mu \}.
\end{align*}
We write 
\begin{equation} \label{equ:decomposition_bmN_cont_disc}
    \bmN = \bmN_c + \bmN_d, \quad \bmN_c := \bigl( \calN_{\alpha,c}, \calN_{\zeta,c}, \calN_{\beta,c}, \calN_{\mu,c} \bigr)^T, \quad \bmN_d := \bigl( \calN_{\alpha,d}, \calN_{\zeta,d}, \calN_{\beta,d}, \calN_{\mu,d} \bigr)^T.
\end{equation}
Correspondingly, we arrive at the following good-bad decomposition for the radiation term 
\begin{equation} \label{equ:good_bad_decomposition}
    \bmu(t) = \bmv(t) + \bmw(t),
\end{equation}
where the ``good'' part $\bmv(t)$ satisfies
\begin{equation} \label{equ:evol_equ_bmv}
    \left\{ \begin{aligned}
        &(\pt^2 + \bfM) \bmv = \bfP_c \bmN_c \\
        &\bmv(0) = \bmu(0), \quad \pt \bmv(0) = \pt \bmu(0),
    \end{aligned} \right.
\end{equation}
while the ``bad'' part $\bmw(t)$ experiences the feedback of the internal mode components with zero initial data
\begin{equation} \label{equ:evol_equ_bmw}
    \left\{ \begin{aligned}
        &(\pt^2 + \bfM) \bmw = \bfP_c \bmN_d \\
        &\bmw(0) = 0, \quad \pt \bmw(0) = 0.
    \end{aligned} \right.
\end{equation}

\subsection{Flat-sharp decomposition} \label{subsec:flat_sharp_decomposition}

\subsubsection{Motivation}

The analysis of the good part $\bmv(t)$ of the radiation term is based on the space-time resonances method.
The main difficulty lies in the normal form analysis of the non-localized quadratic interactions appearing in the nonlinearity $\bmN_c$.
While the distorted Fourier transform associated with $\bfM$ is well suited for linear estimates, it is less convenient for the nonlinear analysis in this paper, in particular for the normal form analysis. For that reason we isolate a component of $\bmv(t)$ whose evolution is governed by the flat Klein-Gordon operator and can therefore be studied using flat Fourier variables on $\bbR^2$. The remaining contribution records the effect of the decaying potential in $\bfM$.

This leads to a second decomposition of the good part $\bmv(t)$ of the radiation term, which we refer to as its \emph{flat-sharp decomposition}. Roughly speaking, the flat component captures the asymptotically free Klein-Gordon dynamics responsible for the leading-order nonlinear interactions, whereas the sharp component measures the discrepancy between the true evolution generated by $\bfM$ and the flat model. Since this discrepancy is caused by a decaying potential, the forcing term for the sharp component is spatially localized. In order to recover acceptable bounds for the sharp component, it is important that the flat component enjoys faster than linear decay. To this end we introduce the flat component below as a solution to a suitably chosen final state problem.

To implement this idea, we first separate the free Klein-Gordon part of the linearized operator. Recalling the definition of $\bfM$, we write
\begin{equation}
    \bfM = \Bigl(-\Delta + \frac{1}{r^2} + 1 \Bigr) \mathbb{I}_{4 \times 4} + \bfV
\end{equation}
with a symmetric matrix potential
\begin{equation} \label{equ:definition_bfV}
    \bfV := \begin{pmatrix} V_1 & -2U' & 0 & 0 \\ -2U' & V_2 & 0 & 0 \\ 0 & 0 & V_1 & -2U' \\ 0 & 0 & -2U' & V_2 \end{pmatrix}
\end{equation}
where the scalar potentials $V_1(r)$, $V_2(r)$ are given by
\begin{equation} \label{equ:definition_V1_V2_for_bfV}
    \begin{aligned}
        V_1 := - \frac{1}{r^2} + \frac{(1-a_\theta)^2}{r^2} - \frac{\pr a_\theta}{r} + U^2 - 1, \quad V_2 := U^2 - 1.
    \end{aligned}
\end{equation}
In particular, we have $|\bfV(r)| \lesssim \jap{r}^{-2}$. Thus, $\bfM$ may be viewed as the flat Klein-Gordon operator on the first angular momentum sector plus a decaying potential.

\subsubsection{Auxiliary flat Klein-Gordon evolutions}

The flat Fourier transform is easiest to use in rectangular coordinates.
Since the components of the radiation variable belong to the first angular momentum sector, it is natural to pass from the equivariant variables to ordinary scalar functions on $\bbR^2$ by multiplying by $e^{i\theta}$, and to then separate real and imaginary parts. From the evolution equation \eqref{equ:evol_equ_bmv} for the good part of the radiation term, we obtain
\begin{equation} \label{equ:motivating_bmv_ast_auxiliary_solutions}
    \bigl( \pt^2 - \Delta + 1 + \bfV \bigr) \bigl( e^{i\theta} \bmv \bigr) = e^{i\theta} \bfP_c \bmN_c.
\end{equation}
The preceding discussion motivates the introduction of auxiliary flat Klein-Gordon evolutions driven by the same good nonlinearity, but without the localized potential and without the projection $\bfP_c$. More precisely, we define real-valued functions $\bmv_\real^\flat(t,x)$ and $\bmv_\imag^\flat(t,x)$ by
\begin{equation} \label{equ:evol_equ_bmv_real_flat}
    \bigl( \pt^2 -\Delta + 1 \bigr) \bmv_\real^\flat = \Re\bigl( e^{i\theta} \bmN_c \bigr), \quad (t,x) \in [0,T]\times\bbR^2,
\end{equation}
and
\begin{equation} \label{equ:evol_equ_bmv_imag_flat}
    \bigl( \pt^2 -\Delta + 1 \bigr) \bmv_\imag^\flat = \Im\bigl( e^{i\theta} \bmN_c \bigr), \quad (t,x) \in [0,T]\times\bbR^2.
\end{equation}
Here $[0,T]$ denotes the bootstrap time interval on which the main estimates are carried out.

A key feature of the forthcoming normal form analysis is that the quadratic interactions generate boundary terms at the final time $T$. If one were to impose vanishing final data, these boundary contributions would reappear in the representation formula for the flat profiles and lead to losses in the weighted energy estimates. We therefore incorporate them directly into the final data and prescribe
\begin{equation} \label{equ:data_at_T_bmv_real_flat}
    \bmv_\real^\flat(T) = 2 \Re\bigl( e^{iT\jD} \bmB_\real[\bmf,\bmf](T) \bigr), \quad \pt \bmv_\real^\flat(T) = -2\jD \Im\bigl( e^{iT\jD} \bmB_\real[\bmf,\bmf](T) \bigr),
\end{equation}
and
\begin{equation} \label{equ:data_at_T_bmv_imag_flat}
    \bmv_\imag^\flat(T) = 2 \Re\bigl( e^{iT\jD} \bmB_\imag[\bmf,\bmf](T) \bigr), \quad \pt \bmv_\imag^\flat(T) = -2\jD \Im\bigl( e^{iT\jD} \bmB_\imag[\bmf,\bmf](T) \bigr).
\end{equation}
The quadratic expressions $\bmB_\real[\bmf,\bmf](T)$ and $\bmB_\imag[\bmf,\bmf](T)$ are the boundary terms produced by the normal form transformations in Section~\ref{sec:working_title_flat_bounds}. Their expressions are specified precisely in \eqref{equ:definition_bmB_boundary_terms} in Subsection~\ref{subsec:quadratic_boundary_terms} and are not ostensibly of degree one. 
Here $\bmf$ collectively denotes the flat profiles $\bmf_\real$ and $\bmf_\imag$ defined in in \eqref{equ:definition_bmf_flat_profile} below. 
Since all bootstrap assumptions for the radiation term on the interval $[0,T]$ are formulated in terms $\bmf_\real$ and $\bmf_\imag$, see \eqref{equ:bootstrap_assumption_weighted_f}-\eqref{equ:bootstrap_assumption_HN_ILED_f} below, the quantities $\bmB_\real[\bmf,\bmf](T)$ and $\bmB_\imag[\bmf,\bmf](T)$ are well-defined under these assumptions.

\subsubsection{Definition of the flat and sharp components}

Since the final states \eqref{equ:data_at_T_bmv_real_flat} and \eqref{equ:data_at_T_bmv_imag_flat} for the auxiliary Klein-Gordon evolutions \eqref{equ:evol_equ_bmv_real_flat} and \eqref{equ:evol_equ_bmv_imag_flat} may not be restricted to degree one, the solutions $\bmv_\real^\flat$ and $\bmv_\imag^\flat$ are not confined to the first angular momentum sector. However, the original radiation variable $\bmv$ belongs to that sector. Consequently, only the first angular mode of the auxiliary flat evolution is relevant for comparison with the true dynamics.

We write the angular decompositions
\begin{equation}
    \bmv_\real^\flat(t,x) = \sum_{n\in\bbZ} \bmv_{\real,n}^\flat(t,r) e^{in\theta}, \qquad \bmv_{\real,n}^\flat(t,r) := \frac{1}{2\pi} \int_0^{2\pi} e^{-in\theta} \bmv_\real^\flat(t,x) \, \ud \theta,
\end{equation}
and
\begin{equation}
    \bmv_\imag^\flat(t,x) = \sum_{n\in\bbZ} \bmv_{\imag,n}^\flat(t,r) e^{in\theta}, \qquad \bmv_{\imag,n}^\flat(t,r) := \frac{1}{2\pi} \int_0^{2\pi} e^{-in\theta} \bmv_\imag^\flat(t,x) \, \ud \theta.
\end{equation}
We then define the flat component of the good part of the radiation by
\begin{equation}
    \bmv^\flat(t,r) := \bmv_{\real,1}^\flat(t,r) + i \bmv_{\imag,1}^\flat(t,r).
\end{equation}
The associated sharp component is defined by
\begin{equation}
    \bmv^\sharp(t,r) := \bmv(t,r) - \bmv^\flat(t,r).
\end{equation}

By construction, $\bmv^\flat(t,r)$ solves the flat Klein-Gordon equation restricted to the first angular momentum sector,
\begin{equation} \label{equ:evol_equ_bmv_flat}
    \Bigl( \pt^2 -\Delta + \frac{1}{r^2} + 1 \Bigr) \bmv^\flat = \bmN_c, \quad (t,r) \in [0,T] \times (0,\infty),
\end{equation}
with final data $\bigl( \bmv^\flat(T), \pt \bmv^\flat(T) \bigr)$ prescribed at time $t=T$ through \eqref{equ:data_at_T_bmv_real_flat}-\eqref{equ:data_at_T_bmv_imag_flat}. This is the reason for the terminology ``flat component.''

Since $\bmv = \bfP_c \bmv$, we have
\begin{equation}
    \bmv(t) = \bfP_c \bmv^\flat(t) + \bfP_c \bmv^\sharp(t).
\end{equation}
Combining \eqref{equ:evol_equ_bmv} and \eqref{equ:evol_equ_bmv_flat}, we find that $\bfP_c \bmv^\sharp$ satisfies
\begin{equation} \label{equ:evol_equ_bmv_sharp}
    \left\{ \begin{aligned}
        &(\pt^2 + \bfM) \bfP_c \bmv^\sharp = - \bfP_c \bigl( \bfV \bmv^\flat \bigr), \quad (t,r) \in [0,T] \times (0,\infty), \\
        &\bfP_c \bmv^\sharp(0) = \bfP_c \bmu(0) - \bfP_c \bmv^\flat(0), \quad \bfP_c (\pt \bmv^\sharp)(0) = \bfP_c (\pt \bmu)(0) - \bfP_c (\pt \bmv^\flat)(0).
    \end{aligned} \right.
\end{equation}
Thus, the sharp component is driven exclusively by the spatially localized term $\bfP_c \bigl( \bfV \bmv^\flat \bigr)$. In particular, it contains no genuinely non-localized quadratic interactions. The role of $\bmv^\sharp$ is therefore fundamentally different from that of $\bmv^\flat$: the latter carries the resonant dynamics and requires a detailed normal form analysis, whereas the former can be estimated using linear dispersive theory together with the spatial decay of the coefficients of $\bfV$ and faster than linear dispersive decay of the flat component~$\bmv^\flat$.

\subsubsection{Flat-sharp decomposition of the temporal component}

We also introduce an analogous decomposition for the good part of the temporal component. The elliptic operator $-\Delta+U^2$ differs from $-\Delta+1$ by the localized coefficient $U^2-1$. We therefore define the flat contribution $\eta_{0,c}^\flat$ by
\begin{equation} \label{equ:eta0c_flat_elliptic_equ}
    (-\Delta + 1) \eta_{0,c}^\flat = \calN_{0,c},
\end{equation}
and set
\begin{equation}
    \eta_{0,c}^\sharp := \eta_{0,c} - \eta_{0,c}^\flat.
\end{equation}
Since $\eta_{0,c}$ solves \eqref{equ:eta0c_equation}, it follows that the sharp component satisfies
\begin{equation} \label{equ:eta0c_sharp_elliptic_equ}
    (-\Delta + U^2) \eta_{0,c}^\sharp = (1-U^2) \eta_{0,c}^\flat.
\end{equation}
Thus the sharp component of the good part of the temporal component is generated solely by the spatially localized term $(1-U^2) \eta_{0,c}^\flat$. This mirrors the flat-sharp decomposition of the good part of the radiation term and allows the estimates for $\eta_{0,c}$ to be reduced to estimates for the flat elliptic variable $\eta_{0,c}^\flat$ plus perturbative bounds for $\eta_{0,c}^\sharp$.

\subsection{Structure of the nonlinearities} \label{subsec:structure_nonlinearities}

The purpose of this subsection is to classify the nonlinearities
$\Re\bigl(e^{i\theta}\bmN_c\bigr)$ and $\Im\bigl(e^{i\theta}\bmN_c\bigr)$ appearing in the flat Klein--Gordon
equations \eqref{equ:evol_equ_bmv_real_flat} and \eqref{equ:evol_equ_bmv_imag_flat} into a finite number of types that
reflect the structure of the subsequent analysis. The classification is
organized according to two guiding principles: whether an interaction is
spatially localized or non-localized, and whether it is genuinely
quadratic or at least cubic. Since the temporal component $\eta_{0}$ is itself (at least) quadratic in the radiation variables, nonlinear
terms involving $\eta_{0}$ require a separate discussion.

The nonlinearities $\Re\bigl( e^{i\theta} \bmN_c \bigr)$ and $\Im\bigl( e^{i\theta} \bmN_c \bigr)$ feature the components of $e^{i\theta} \bmu(t)$ as well as the temporal component $\eta_{0,c}(t)$.
Recalling that $\underline{\Phi}(r,\theta) = e^{i\theta} U(r)$, we for instance have 
\begin{equation}
    \begin{aligned}
        e^{i\theta} \calN_{\alpha,c} &= -\frac32 \underline{\Phi} \bigl| e^{i\theta} u_1 \bigr|^2 + \frac12 \underline{\Phi} \bigl| e^{i\theta} u_3 \bigr|^2 - \frac12 \bigl( \bigl| e^{i\theta} u_1 \bigr|^2 + \bigl| e^{i\theta} u_3 \bigr|^2 \bigr) \bigl( e^{i\theta} u_1 \bigr) \\
        &\quad - \bigl( e^{i\theta} u_4 \bigr) (\partial_1 - i\partial_2) \bigl( e^{i\theta} u_3 \bigr) - \overline{e^{i\theta} u_4} (\partial_1 + i \partial_2) \bigl( e^{i\theta} u_3 \bigr) \\
        &\quad + (1-a_\theta) \bigl( e^{i\theta} u_2 \bigr) (\partial_1 - i\partial_2) \bigl( e^{i\theta} u_1 \bigr) - (1-a_\theta) \bigl( e^{-i\theta} u_2 \bigr) (\partial_1 + i \partial_2) \bigl( e^{i\theta} u_1 \bigr) \\
        &\quad - \underline{\Phi} \bigl( \bigl| e^{i\theta} u_2 \bigr|^2 + \bigl| e^{i\theta} u_4 \bigr|^2 \bigr) - \bigl( \bigl| e^{i\theta} u_2 \bigr|^2 + \bigl| e^{i\theta} u_4 \bigr|^2 \bigr) \bigl( e^{i\theta} u_1 \bigr) \\ 
        &\quad - 2 \eta_{0,c} \pt \bigl( e^{i\theta} u_3 \bigr) - \pt \eta_{0,c} \bigl( e^{i\theta} u_3 \bigr) + \underline{\Phi} \eta_{0,c}^2 + \eta_{0,c}^2 \bigl( e^{i\theta} u_1 \bigr).
    \end{aligned}
\end{equation}
Here in order to rewrite the nonlinearities $-2 e^{i\theta} u_4 \partial_r u_3$ and $2 b e^{i\theta} u_2 u_1$, we used the simple fact that $\frac{1}{r} e^{i\theta} u_1 = \frac{-i}{r} \partial_\theta (e^{i\theta} u_1)$ as well as the identities
\begin{equation}
    \begin{aligned}
        \partial_r = \frac12 e^{i\theta} (\partial_1 - i\partial_2) + \frac12 e^{-i\theta} (\partial_1 + i\partial_2), \quad \frac{1}{r} \partial_\theta = \frac{i}{2} e^{i\theta} (\partial_1-i\partial_2) - \frac{i}{2} e^{-i\theta} (\partial_1+i\partial_2).
    \end{aligned}
\end{equation}

Next, we decompose $e^{i\theta} \bmu(t)$ into its real and imaginary parts
\begin{equation}
    \begin{aligned}
        e^{i\theta} \bmu(t) = \bmu_\real(t) + i \bmu_\imag(t), \quad \bmu_\real(t) := \Re \bigl( e^{i\theta} \bmu(t) \bigr), \quad \bmu_\imag(t) := \Im \bigl( e^{i\theta} \bmu(t) \bigr),
    \end{aligned}
\end{equation}
and we also decompose the vortex profile into its real and imaginary parts 
\begin{equation}
    \underline{\Phi} = \underline{\Phi}_\real + i \underline{\Phi}_\imag, \quad \underline{\Phi}_\real := \Re \bigl( e^{i\theta} U(r) \bigr), \quad \underline{\Phi}_\imag := \Im \bigl( e^{i\theta} U(r) \bigr).
\end{equation}
This leads to the following expressions for the first components of the (vectorial) nonlinearities $\Re\bigl( e^{i\theta} \bmN_c \bigr)$ and $\Im\bigl( e^{i\theta} \bmN_c \bigr)$,
\begin{equation}
    \begin{aligned}
        \Re\bigl( e^{i\theta} \calN_{\alpha,c} \bigr) &= -\frac32 \underline{\Phi}_\real \bigl( u_{1,\real}^2 + u_{1,\imag}^2 \bigr) + \frac12 \underline{\Phi}_\real \bigl( u_{3,\real}^2 + u_{3,\imag}^2 \bigr) - \frac12 \bigl( u_{1,\real}^2 + u_{1,\imag}^2 + u_{3,\real}^2 + u_{3,\imag}^2 \bigr) u_{1,\real} \\ 
        &\quad - 2 u_{4,\real} \partial_1 u_{3,\real} - 2 u_{4,\imag} \partial_2 u_{3,\real} + 2 (1-a_\theta) u_{2,\real} \partial_2 u_{1,\imag} - 2 (1-a_\theta) u_{2,\imag} \partial_1 u_{1,\imag} \\
        &\quad - \underline{\Phi}_\real \bigl( u_{2,\real}^2 + u_{2,\imag}^2 + u_{4,\real}^2 + u_{4,\imag}^2 \bigr) - \bigl( u_{2,\real}^2 + u_{2,\imag}^2 + u_{4,\real}^2 + u_{4,\imag}^2 \bigr) u_{1,\real} \\ 
        &\quad - 2 \eta_{0,c} \pt u_{3,\real} - \pt \eta_{0,c} u_{3,\real} + \underline{\Phi}_\real \eta_{0,c}^2 + \eta_{0,c}^2 u_{1,\real}, \\
        \Im\bigl( e^{i\theta} \calN_{\alpha,c} \bigr) &= -\frac32 \underline{\Phi}_\imag \bigl( u_{1,\real}^2 + u_{1,\imag}^2 \bigr) +\frac12 \underline{\Phi}_\imag \bigl( u_{3,\real}^2 + u_{3,\imag}^2 \bigr) - \frac12 \bigl( u_{1,\real}^2 + u_{1,\imag}^2 + u_{3,\real}^2 + u_{3,\imag}^2 \bigr) u_{1,\imag} \\
        &\quad - 2 u_{4,\real} \partial_1 u_{3,\imag}- 2 u_{4,\imag} \partial_2 u_{3,\imag} - 2(1-a_\theta)u_{2,\real}\partial_2u_{1,\real} + 2(1-a_\theta)u_{2,\imag}\partial_1u_{1,\real} \\
        &\quad - \underline{\Phi}_\imag \bigl( u_{2,\real}^2 + u_{2,\imag}^2 + u_{4,\real}^2 + u_{4,\imag}^2 \bigr) - \bigl( u_{2,\real}^2 + u_{2,\imag}^2 + u_{4,\real}^2 + u_{4,\imag}^2 \bigr)u_{1,\imag} \\
        &\quad -2\eta_{0,c}\pt u_{3,\imag} - \pt\eta_{0,c}u_{3,\imag} +\underline{\Phi}_\imag\eta_{0,c}^2+\eta_{0,c}^2u_{1,\imag}.
    \end{aligned}
\end{equation}
For the second components of $\Re\bigl( e^{i\theta} \bmN_c \bigr)$ and $\Im\bigl( e^{i\theta} \bmN_c \bigr)$ we find 
\begin{equation}
\begin{aligned}
    \Re\bigl( e^{i\theta} \calN_{\zeta,c} \bigr) &= (1-a_\theta)\bigl(u_{1,\real}\partial_2u_{1,\imag}-u_{1,\imag}\partial_1u_{1,\imag}\bigr) + (1-a_\theta)\bigl(u_{3,\real}\partial_2u_{3,\imag} -u_{3,\imag}\partial_1u_{3,\imag}\bigr) \\
    &\quad -2\underline{\Phi}_\real \bigl(u_{2,\real}u_{1,\real}+u_{2,\imag}u_{1,\imag}\bigr) - \bigl(u_{1,\real}^2+u_{1,\imag}^2+u_{3,\real}^2+u_{3,\imag}^2\bigr)u_{2,\real}, \\ 
    \Im\bigl( e^{i\theta} \calN_{\zeta,c} \bigr) &= (1-a_\theta) \bigl( -u_{1,\real}\partial_2u_{1,\real} + u_{1,\imag}\partial_1u_{1,\real} \bigr) + (1-a_\theta)\bigl( -u_{3,\real}\partial_2u_{3,\real} + u_{3,\imag}\partial_1u_{3,\real} \bigr) \\
    &\quad -2\underline{\Phi}_\imag \bigl(u_{2,\real}u_{1,\real}+u_{2,\imag}u_{1,\imag}\bigr) -\bigl(u_{1,\real}^2+u_{1,\imag}^2+u_{3,\real}^2+u_{3,\imag}^2\bigr)u_{2,\imag},
\end{aligned}
\end{equation}
and for the third components of $\Re\bigl( e^{i\theta} \bmN_c \bigr)$ and $\Im\bigl( e^{i\theta} \bmN_c \bigr)$ we obtain
\begin{equation}
\begin{aligned}
    \Re\bigl( e^{i\theta} \calN_{\beta,c} \bigr) &= -2\underline{\Phi}_\real \bigl(u_{1,\real}u_{3,\real}+u_{1,\imag}u_{3,\imag}\bigr) -\frac12\bigl(u_{1,\real}^2+u_{1,\imag}^2+u_{3,\real}^2+u_{3,\imag}^2\bigr)u_{3,\real} \\
    &\quad +2u_{4,\real}\partial_1u_{1,\real} +2u_{4,\imag}\partial_2u_{1,\real} +2(1-a_\theta) \bigl( u_{2,\real}\partial_2u_{3,\imag} - u_{2,\imag}\partial_1u_{3,\imag} \bigr) \\
    &\quad -\bigl(u_{2,\real}^2+u_{2,\imag}^2+u_{4,\real}^2+u_{4,\imag}^2\bigr)u_{3,\real} \\
    &\quad +2\eta_{0,c}\pt u_{1,\real}+\underline{\Phi}_\real\pt\eta_{0,c}+\pt\eta_{0,c}u_{1,\real}+\eta_{0,c}^2u_{3,\real}, \\
    \Im\bigl( e^{i\theta} \calN_{\beta,c} \bigr) &= -2\underline{\Phi}_\imag \bigl(u_{1,\real}u_{3,\real}+u_{1,\imag}u_{3,\imag}\bigr) - \frac12\bigl(u_{1,\real}^2+u_{1,\imag}^2+u_{3,\real}^2+u_{3,\imag}^2\bigr)u_{3,\imag} \\
    &\quad +2u_{4,\real}\partial_1u_{1,\imag}+2u_{4,\imag}\partial_2u_{1,\imag} -2(1-a_\theta) \bigl( u_{2,\real}\partial_2u_{3,\real} - u_{2,\imag}\partial_1u_{3,\real} \bigr) \\
    &\quad -\bigl(u_{2,\real}^2+u_{2,\imag}^2+u_{4,\real}^2+u_{4,\imag}^2\bigr)u_{3,\imag} \\
    &\quad +2\eta_{0,c}\pt u_{1,\imag} +\underline{\Phi}_\imag\pt\eta_{0,c}+\pt\eta_{0,c}u_{1,\imag}+\eta_{0,c}^2u_{3,\imag}.
\end{aligned}
\end{equation}
Finally, the fourth components of $\Re\bigl( e^{i\theta} \bmN_c \bigr)$ and $\Im\bigl( e^{i\theta} \bmN_c \bigr)$ are given by
\begin{equation}
\begin{aligned}
    \Re\bigl( e^{i\theta} \calN_{\mu,c} \bigr) &= -\partial_1\pt\eta_{0,c} -u_{1,\real}\partial_1u_{3,\real} -u_{1,\imag}\partial_2u_{3,\real} + u_{3,\real}\partial_1u_{1,\real} +u_{3,\imag}\partial_2u_{1,\real} \\
    &\quad -2\underline{\Phi}_\real \bigl(u_{1,\real}u_{4,\real}+u_{1,\imag}u_{4,\imag}\bigr) -\bigl(u_{1,\real}^2+u_{1,\imag}^2+u_{3,\real}^2+u_{3,\imag}^2\bigr)u_{4,\real}, \\
    \Im\bigl( e^{i\theta} \calN_{\mu,c} \bigr) &= -\partial_2\pt\eta_{0,c} -u_{1,\real}\partial_1u_{3,\imag} -u_{1,\imag}\partial_2u_{3,\imag} + u_{3,\real}\partial_1u_{1,\imag} +u_{3,\imag}\partial_2u_{1,\imag} \\
    &\quad -2\underline{\Phi}_\imag \bigl(u_{1,\real}u_{4,\real}+u_{1,\imag}u_{4,\imag}\bigr) -\bigl(u_{1,\real}^2+u_{1,\imag}^2+u_{3,\real}^2+u_{3,\imag}^2\bigr)u_{4,\imag}.
\end{aligned}
\end{equation}

Collecting the terms in the preceding formulas, we obtain the following
schematic classification of the nonlinearities that enter the auxiliary flat
Klein--Gordon evolutions \eqref{equ:evol_equ_bmv_real_flat}, \eqref{equ:evol_equ_bmv_imag_flat}. 
In the statement below, numerical constants and signs are suppressed, and $\partial_j$, $j=1,2$, denotes a Cartesian coordinate derivative.

\begin{corollary} \label{cor:types_of_nonlinearities}
For each $\bullet \in \{\alpha,\zeta,\beta,\mu\}$, the components $\Re\bigl(e^{i\theta}\calN_{\bullet,c}\bigr)$ and $\Im\bigl(e^{i\theta}\calN_{\bullet,c}\bigr)$ are finite real linear combinations of nonlinear terms of the following ten types:
\begin{enumerate} 
    \item $\underline{\Phi}_\ast u_{k,\ast_k} u_{l,\ast_l}$, $\ast, \ast_k, \ast_l \in \{\real, \imag\}$, $k, l \in \{ 1, 2, 3, 4\}$,
    \item $u_{k, \ast_k} \partial_j u_{l, \ast_l}$, $j \in \{1,2\}$, $\ast_k, \ast_l \in \{\real, \imag\}$, $k, l \in \{1, 2, 3, 4 \}$,    
    \item $u_{k,\ast_k} u_{l, \ast_l} u_{m, \ast_m}$, $\ast_k, \ast_l, \ast_m \in \{\real, \imag\}$, $k, l, m \in \{1, 2, 3, 4 \}$,
    \item $(1-a_\theta) u_{k,\ast_k} \partial_j u_{l, \ast_l}$, $j \in \{1,2\}$, $\ast_k, \ast_l \in \{\real, \imag\}$, $k, l \in \{1, 2, 3, 4 \}$,
    \item $\partial_j \partial_t \eta_{0,c}$, $j \in \{1, 2\}$,
    \item $\underline{\Phi}_\ast (\partial_t \eta_{0,c})$, $\ast \in \{\real, \imag\}$,    
    \item $\eta_{0,c} \partial_t u_{k,\ast_k}$, $\ast_k \in \{\real, \imag\}$, $k \in \{1, 2, 3, 4\}$,
    \item $(\partial_t \eta_{0,c}) u_{k, \ast_k}$, $\ast_k \in \{\real, \imag\}$, $k \in \{1, 2, 3, 4\}$, 
    \item $\underline{\Phi}_\ast \eta_{0,c}^2$, $\ast \in \{\real, \imag\}$,
    \item $\eta_{0,c}^2 u_{k, \ast_k}$, $\ast_k \in \{\real, \imag\}$, $k \in \{1, 2, 3, 4\}$.
\end{enumerate}
\end{corollary}

The distinction between the ten classes is important for the analysis in
Section~\ref{sec:working_title_flat_bounds}. The only genuinely
non-spatially localized quadratic interactions are those of types
\emph{(1)}, \emph{(2)}, \emph{(5)}, and \emph{(6)}. In types \emph{(5)}
and \emph{(6)} this refers to the leading quadratic contribution to
$\partial_t\eta_{0,c}$ obtained from the elliptic equation for the temporal
component. These four classes are therefore the terms that require the
normal form analysis in Section~\ref{sec:working_title_flat_bounds}. The
remaining classes are perturbative in a more direct way: type \emph{(3)}
is cubic, type \emph{(4)} carries the spatially decaying coefficient
$1-a_\theta$, and types \emph{(7)}--\emph{(10)} are at least cubic once the
quadratic character of $\eta_{0,c}$ is taken into account.

Finally, we also consider the nonlinearities in the elliptic equations for $\eta_{0,c}^\flat$ and $\pt \eta_{0,c}^\flat$. The former can be rewritten as 
\begin{equation} \label{equ:eta0c_flat_equation_rewritten}
    \begin{aligned}
        (-\Delta + 1) \eta_{0,c}^\flat &= -2 \underline{\Phi} \eta_{0,c} \overline{\bigl(e^{i\theta} u_1\bigr)} - \overline{\bigl(e^{i\theta} u_3\bigr)} \pt \bigl( e^{i\theta} u_1 \bigr) \\
        &\quad + \overline{\bigl(e^{i\theta} u_1\bigr)} \pt \bigl( e^{i\theta} u_3 \bigr) - \eta_{0,c} \bigl( \bigl| e^{i\theta} u_1 \bigr|^2 + \bigl| e^{i\theta} u_3 \bigr|^2 \bigr) \\
        &= -2 \eta_{0,c} \bigl( \underline{\Phi}_{\real} u_{1,\real} + \underline{\Phi}_{\imag} u_{1,\imag} \bigr) - \bigl( u_{3,\real} \pt u_{1,\real} + u_{3,\imag} \pt u_{1,\imag} \bigr) \\ 
        &\quad + u_{1,\real} \pt u_{3,\real} + u_{1,\imag} \pt u_{3,\imag} - \eta_{0,c} \bigl( u_{1,\real}^2 + u_{1,\imag}^2 + u_{3,\real}^2 + u_{3,\imag}^2 \bigr).
    \end{aligned}
\end{equation}
To determine the equation for $\pt \eta_{0,c}^\flat$, we reinsert the evolution equation for the radiation term in \eqref{equ:system_evol_equations_u_z}, which gives
\begin{equation}
    \begin{aligned}
        &(-\Delta + 1) \pt \eta_{0,c}^\flat \\
        &\quad = -u_3 \pt^2 u_1 + u_1 \pt^2 u_3 - 2 U \pt \eta_{0,c} u_1 - 2 U \eta_{0,c} \pt u_1 - \pt \eta_{0,c} \bigl( u_1^2 + u_3^2 \bigr) - \eta_{0,c} \bigl( 2 u_1 \pt u_1 + 2 u_3 \pt u_3 \bigr) \\ 
        &\quad = u_3 \Bigl( -\Delta + \frac{1}{r^2} \Bigr) u_1 - u_1 \Bigl( -\Delta + \frac{1}{r^2} \Bigr) u_3 - 2U' u_2 u_3 + 2U' u_1 u_4 - u_3 ( \bfP_c \bmN )_1 + u_1 (\bfP_c \bmN)_3 \\
        &\quad \quad - 2 U \pt \eta_{0,c} u_1 - 2 U \eta_{0,c} \pt u_1 - \pt \eta_{0,c} \bigl( u_1^2 + u_3^2 \bigr) - \eta_{0,c} \bigl( 2 u_1 \pt u_1 + 2 u_3 \pt u_3 \bigr).
    \end{aligned}
\end{equation}
We can then rewrite the preceding equation as 
\begin{equation} \label{equ:pteta0c_flat_equation_rewritten}
    \begin{aligned}
        &(-\Delta + 1) \pt \eta_{0,c}^\flat \\
        &\quad = u_{3,\real} (-\Delta u_{1,\real}) + u_{3,\imag} (-\Delta u_{1,\imag}) - u_{1,\real} (-\Delta u_{3,\real}) - u_{1,\imag} (-\Delta u_{3,\imag}) \\
        &\quad \quad -2 U' \bigl( u_{2,\real} u_{3,\real} + u_{2,\imag} u_{3,\imag} \bigr) + 2 U' \bigl( u_{1,\real} u_{4,\real} + u_{1,\imag} u_{4,\imag} \bigr) \\
        &\quad \quad -u_{3,\real} \Re\bigl( e^{i\theta} (\bfP_c \bmN)_1 \bigr) -u_{3,\imag} \Im\bigl(e^{i\theta} (\bfP_c \bmN)_1\bigr) \\
        &\quad \quad + u_{1,\real} \Re\bigl(e^{i\theta} (\bfP_c \bmN)_3\bigr) + u_{1,\imag} \Im\bigl(e^{i\theta} (\bfP_c \bmN)_3\bigr) \\
        &\quad \quad - 2 \underline{\Phi}_\real \pt \eta_{0,c} u_{1,\real} - 2 \underline{\Phi}_\imag \pt \eta_{0,c} u_{1,\imag} - 2 \underline{\Phi}_\real  \eta_{0,c} \pt u_{1,\real} - 2 \underline{\Phi}_\imag \eta_{0,c} \pt u_{1,\imag} \\ 
        &\quad \quad - \pt \eta_{0,c} \bigl( u_{1,\real}^2 + u_{1,\imag}^2 + u_{3,\real}^2 + u_{3,\imag}^2 \bigr) \\
        &\quad \quad - \eta_{0,c} \bigl( 2 u_{1,\real} \pt u_{1,\real} + 2 u_{1,\imag} \pt u_{1,\imag} + 2 u_{3,\real} \pt u_{3,\real} + 2 u_{3,\imag} \pt u_{3,\imag} \bigr). 
    \end{aligned}
\end{equation}
Observe that the first line on the right-hand side of \eqref{equ:pteta0c_flat_equation_rewritten} can be rewritten as
\begin{equation} \label{equ:pteta0c_flat_equation_leading_order_term_rewritten}
    \begin{aligned}
        &u_{3,\real} (-\Delta u_{1,\real}) + u_{3,\imag} (-\Delta u_{1,\imag}) - u_{1,\real} (-\Delta u_{3,\real}) - u_{1,\imag} (-\Delta u_{3,\imag}) \\
        &= \nabla \cdot \bigl( u_{1,\real} \nabla u_{3,\real} - u_{3,\real} \nabla u_{1,\real} \bigr) + \nabla \cdot \bigl( u_{1,\imag} \nabla u_{3,\imag} - u_{3,\imag} \nabla u_{1,\imag} \bigr).
    \end{aligned}
\end{equation}
This div-curl structure is relevant later when we prove energy estimates with $N$-derivatives.

\subsection{Bootstrap setup} \label{subsec:bootstrap_setup}

By the local and global existence result from \cite{BurMon1}, the solution to the system of evolution equations \eqref{equ:system_evol_equations_u_z}, \eqref{equ:system_evol_equations_u_z_part2} exists locally in time at $H^N_x(\bbR^2)$-regularity on some maximal interval of existence $[0,T_\ast)$ for some $0 < T_\ast \leq \infty$ under the hypotheses of Theorem~\ref{thm:main} and preserves the gauge condition \eqref{eq:Stuartgaugeintro2}. Moreover, if $\|(e^{i\theta} \bmu(t), e^{i\theta} \pt \bmu(t))\|_{H^N_x \times H^{N-1}_x}$ stays small for all $0 \leq t < T_\ast$, then we have $T_\ast = \infty$. We will recall in Section~\ref{sec:conclusion_proof_thm_main} how the results from \cite{BurMon1} apply in the context of the orthogonal gauge condition~\eqref{equ:intro_Stuart_gauge_polar} in this work. 

We establish global existence and the decay estimates \eqref{eq:decay_zed_mainthm} and \eqref{eq:decay_radiation_mainthm} in part (i) of the statement of Theorem~\ref{thm:main} by a bootstrap argument. The additional lower bound \eqref{equ:decay_zed_lower_bound_mainthm} for the internal mode in part (ii) does not require a separate bootstrap assumption. Once the bootstrap estimates used for part (i) have been established, the lower bound in \eqref{equ:decay_zed_lower_bound_mainthm} follows directly from part (ii) of Proposition~\ref{propODEdamp} under the additional assumption \eqref{equ:new_assumption_frakz_equal}.


In this subsection we lay out the precise bootstrap setup for the analysis in the remainder of this paper.
Fix $0 < T < T_\ast$. All estimates below are formulated on an arbitrary interval $[0,T]$, with constants independent of $T$.
In the estimates below we may assume $T \geq \varepsilon^{-2}$ (and thus $T_\ast > \varepsilon^{-2}$). If $T < \varepsilon^{-2}$, the long-time regimes are empty, and the corresponding estimates are obtained by restricting the short-time arguments to $[0,T]$.

The profiles of the internal mode components are denoted by
\begin{equation}\label{eq:Zjprofiles1}
    Z_j(t) := e^{-it\lambda} (2i\lambda)^{-1} (\pt + i \lambda) z_j(t), \quad j = 1, 2.
\end{equation}
Observe that $\frakz(t)$ introduced in \eqref{eq:zeddefintro1} satisfies
\begin{equation} 
    \mathfrak{z}(t) = |Z_1(t)|^2 + |Z_2(t)|^2, \quad t \geq 0.
\end{equation}
The majority of our nonlinear analysis focuses on the auxiliary evolutions for $\bmv_\real^\flat(t,x)$ and $\bmv_\imag^\flat(t,x)$ in Section~\ref{sec:working_title_flat_bounds} using the space-time resonances method based on the flat Fourier transform.
Its nonlinearities feature the radiation term in terms of the variables $\bmu_\real(t,x)$ and $\bmu_\imag(t,x)$.
It is therefore natural to introduce their associated flat profiles
\begin{equation} \label{equ:definition_bmf_flat_profile}
    \begin{aligned}
        \bmf_\real(t) &:= e^{-it\jD} (2i\jD)^{-1} (\pt+i\jD) \Re\bigl( e^{i\theta} \bmu(t) \bigr), \\
        \bmf_\imag(t) &:= e^{-it\jD} (2i\jD)^{-1} (\pt+i\jD) \Im\bigl( e^{i\theta} \bmu(t) \bigr),
    \end{aligned}
\end{equation}
and to formulate our bootstrap assumptions about the radiation term in terms of these flat profiles.

Recall from the statement of Theorem~\ref{thm:main} that $0 < \varepsilon_0 \ll 1$ denotes a small absolute constant whose size will be fixed sufficiently small in the course of the proof of Theorem~\ref{thm:main}. Moreover, we pick small absolute constants $0 < \delta \ll 1$ and $0 < \kappa \ll 1$. Additionally, we fix $N \in \bbN$ sufficiently large such that 
\begin{equation} \label{equ:assumption_N_delta}
    \frac{100}{N} \leq \delta.
\end{equation}
We denote by $C_0 \geq 1$ a sufficiently large absolute constant. 
Recall from \eqref{equ:new_assumption1}, \eqref{equ:new_assumption2} in part (i) of the statement of Theorem~\ref{thm:main} that for some $0 < \varepsilon \leq \varepsilon_0$ the initial conditions satisfy $\mathfrak{z}(0) \leq \varepsilon^2$ and 
\begin{equation} 
    \begin{aligned}
        &\bigl\| \jap{x} \bigl( e^{i\theta} \bigl( \bmu(0), \pt \bmu(0) \bigr) \bigr) \bigr\|_{H^2_x(\bbR^2) \times H^1_x(\bbR^2)} + \bigl\| e^{i\theta} \bigl( \bmu(0), \pt \bmu(0) \bigr) \bigr\|_{H^N_x(\bbR^2) \times H^{N-1}_x(\bbR^2)} \leq \varepsilon.
    \end{aligned}
\end{equation}
We impose the following bootstrap assumption for the upper decay bound of the internal mode components on the interval $[0,T]$,
\begin{equation} \label{equ:bootstrap_assumption_z}
    \frakz(t) \leq \frac{3}{2} \frakz_{\Gamma_1}(t), \quad 0 \leq t \leq T,
\end{equation}
where we use the notation
\begin{equation} \label{eq:frakzfirstdef1}
    \frakz_{\Gamma_1}(t) := \frac{\varepsilon^2}{1 + \Gamma_1 \varepsilon^2 t}
\end{equation}
with $\Gamma_1 > 0$ the Fermi Golden Rule constant defined in \eqref{eq:Gamma01def}.
These are the only bootstrap assumptions involving the internal mode components. In particular, under the additional hypothesis $\frakz(0) = \varepsilon^2$ in 
\eqref{equ:new_assumption_frakz_equal} from part (ii) of the statement of Theorem~\ref{thm:main}, we do not bootstrap a lower bound for $\frakz(t)$. The lower bound asserted in \eqref{equ:decay_zed_lower_bound_mainthm} will instead follow directly from part (ii) of Proposition~\ref{propODEdamp}.


Moreover, we make the following bootstrap assumptions about the flat profiles $\bmf_\ast(t)$, $\ast~\in~\{\real, \imag\}$, and their evolutions on a given bootstrap interval $[0,T]$,
\begin{align}
    \bigl\| \jxi^2 \jap{\nabla_\xi} \widehat{\bmf}_\ast(t,\xi) \bigr\|_{L^2_\xi(\bbR^2)} &\leq 2C_0 \left\{ \begin{aligned}
                            &\varepsilon^{1-20\delta}, &&0 \leq t \leq 1, \\
                            &\varepsilon^{1-20\delta} t^{\frac12}, &&1 \leq t \leq \varepsilon^{-1-20\delta}, \\
                            &\varepsilon^2 t^{\frac32},  &&\varepsilon^{-1-20\delta} \leq t \leq \varepsilon^{-2}, \\
                            &t^{\frac12},  &&\varepsilon^{-2} \leq t \leq T,
                         \end{aligned} \right. \label{equ:bootstrap_assumption_weighted_f} \\ 
    \bigl\| e^{it\jD} \bmf_\ast(t)\bigr\|_{L^\infty_x(\bbR^2)} &\leq 2C_0 \left\{ \begin{aligned}
                            &\varepsilon^{1-10\delta}, &&0 \leq t \leq 1, \\    
                            &\varepsilon^{1-10\delta} t^{-\frac12-4\delta},  &&1 \leq t \leq \varepsilon^{-2}, \\
                            &t^{-1+\delta},  &&\varepsilon^{-2} \leq t \leq T,
                         \end{aligned} \right. \label{equ:bootstrap_assumption_Linfty_f} \\ 
    \bigl\| \bmf_\ast(t)\bigr\|_{H^N_x(\bbR^2)} &\leq 2C_0 \varepsilon, \quad 0 \leq t \leq T, \label{equ:bootstrap_assumption_HN_f} \\
    \sum_{1 \leq k \leq N} \, \bigl\| \jx^{-\frac32-\kappa} |D|^k e^{it\jD} \bmf_\ast(t)\bigr\|_{L^2_t([0,T]; L^2_x(\bbR^2))} &\leq 2C_0 \varepsilon. \label{equ:bootstrap_assumption_HN_ILED_f} 
\end{align}

Let us briefly indicate the role of the bootstrap norms. The bound \eqref{equ:bootstrap_assumption_z} encodes the expected Fermi Golden Rule decay of the internal mode components.
The weighted Fourier bound \eqref{equ:bootstrap_assumption_weighted_f} is the main vector-field type norm used in the space-time resonance analysis of the flat component of the good part of the radiation term. Its strong growth is caused by the quadratic internal mode source terms in the evolution equation for the radiation term.
The dispersive estimate \eqref{equ:bootstrap_assumption_Linfty_f} is the pointwise decay bound needed in all nonlinear product estimates.
The uniform-in-time high Sobolev norm bound \eqref{equ:bootstrap_assumption_HN_f} allows us to absorb all the losses of derivatives that come up in our arguments and helps to streamline some of the normal form analysis; for example, we can estimate all the multi-linear symbols that appear in the analysis without too much effort, and do not need to worry about the derivatives in the nonlinear terms. 
The integrated local energy decay estimate at high Sobolev regularity \eqref{equ:bootstrap_assumption_HN_ILED_f} is primarily used in the derivation of the high Sobolev norm bound \eqref{equ:bootstrap_assumption_HN_f}.

\begin{rem}
Under the assumptions \eqref{equ:bootstrap_assumption_Linfty_f} and 
\eqref{equ:bootstrap_assumption_HN_f}
it follows that for all $1 \leq \ell \leq 40$ and any $p \geq N/20$,
\begin{align}\label{equ:GNSell}
\bigl\| \jD^\ell e^{it\jD} \bmf_\ast(t) \bigr\|_{L^p_x(\bbR^2)} 
 & \lesssim_{C_0} \left\{ \begin{aligned}
&\varepsilon^{1-10\delta} \jt^{-\frac12-\frac72\delta},  &&0 \leq t \leq \varepsilon^{-2},
\\
 &t^{-1+\frac32\delta},  
 &&\varepsilon^{-2} \leq t \leq T.
\end{aligned} \right.
\end{align}
The inequality \eqref{equ:GNSell} follows from the Gagliardo-Nirenberg-Sobolev interpolation inequality
\begin{equation} \label{equ:GNS}
 \bigl\| \jap{D}^k g \bigr\|_{L^p_x(\bbR^2)} 
 \lesssim \bigl\| \jap{D}^N g \bigr\|_{L^2_x(\bbR^2)}^{\frac{k-\frac{2}{p}}{N-1}} 
 \|g\|_{L^\infty_x(\bbR^2)}^{1- \frac{k-\frac{2}{p}}{N-1}}, \quad 1 \leq k < N-1, 
\end{equation}
the bootstrap assumptions \eqref{equ:bootstrap_assumption_Linfty_f} 
and \eqref{equ:bootstrap_assumption_HN_f},
and our choice of parameters ($100/N \leq \delta \leq 1/100$). 
\end{rem}

\subsection{Decomposition of the flat profiles} \label{subsec:flat_profiles}

We now decompose the flat profiles $\bmf_\real(t)$ and $\bmf_\imag(t)$ according to the good-bad decomposition introduced in Subsection~\ref{subsec:good_bad_decomposition} and the flat-sharp decomposition introduced in Subsection~\ref{subsec:flat_sharp_decomposition}.
Combining the good-bad splitting $\bmu(t) = \bmv(t) + \bmw(t)$ with the flat-sharp splitting $\bmv(t) = \bmv^\flat(t) + \bmv^\sharp(t)$, and recalling that $\bmu(t) = \bfP_c \bmu(t)$, we obtain
\begin{equation} \label{equ:decomposition_radiation_term_flat_profiles_definitions}
    \bmu(t) = \bfP_c \bmv^\flat(t) + \bfP_c \bmv^\sharp(t) + \bfP_c \bmw(t),
\end{equation}
whence 
\begin{align} 
        \bmu_\real(t) &= \Re\bigl( e^{i\theta} \bfP_c \bmv^\flat(t) \bigr) + \Re\bigl( e^{i\theta} \bfP_c \bmv^\sharp(t) \bigr) + \Re\bigl( e^{i\theta} \bfP_c \bmw(t) \bigr), \label{equ:decomposition_bmu_real_flat_profiles_definitions} \\ 
        \bmu_\imag(t) &= \Im\bigl( e^{i\theta} \bfP_c \bmv^\flat(t) \bigr) + \Im\bigl( e^{i\theta} \bfP_c \bmv^\sharp(t) \bigr) + \Im\bigl( e^{i\theta} \bfP_c \bmw(t) \bigr). \label{equ:decomposition_bmu_imag_flat_profiles_definitions}
\end{align}
In order to improve the bootstrap bounds \eqref{equ:bootstrap_assumption_weighted_f}, \eqref{equ:bootstrap_assumption_Linfty_f}, \eqref{equ:bootstrap_assumption_HN_f}, \eqref{equ:bootstrap_assumption_HN_ILED_f} for the flat profiles $\bmf_\real(t)$ and $\bmf_\imag(t)$ associated with the radiation term variables $\bmu_\real(t)$ and $\bmu_\imag(t)$, it is natural to derive improved bounds for the flat profiles associated with each of the three components on the right-hand sides of \eqref{equ:decomposition_bmu_real_flat_profiles_definitions} and \eqref{equ:decomposition_bmu_imag_flat_profiles_definitions}.
Correspondingly, we introduce the flat profiles associated with $\Re\bigl( e^{i\theta} \bfP_c \bmv^\sharp(t)\bigr)$ and $\Im\bigl( e^{i\theta} \bfP_c \bmv^\sharp(t)\bigr)$ by
\begin{equation} \label{equ:definition_bmg_sharp_profile}
    \begin{aligned}
        \bmg^\sharp_\real(t) &:= e^{-it\jD} (2i\jD)^{-1} (\pt+i\jD) \Re \bigl( e^{i\theta} \bfP_c \bmv^\sharp(t) \bigr), \\
        \bmg^\sharp_\imag(t) &:= e^{-it\jD} (2i\jD)^{-1} (\pt+i\jD) \Im \bigl( e^{i\theta} \bfP_c \bmv^\sharp(t) \bigr),
    \end{aligned}
\end{equation}
and we denote the flat profiles associated with $\Re\bigl(e^{i\theta} \bfP_c \bmw(t)\bigr)$ and $\Im\bigl( e^{i\theta} \bfP_c \bmw(t)\bigr)$ by
\begin{equation} \label{equ:definition_bmh_profile}
    \begin{aligned}
        \bmh_\real(t) &:= e^{-it\jD} (2i\jD)^{-1} (\pt+i\jD) \Re \bigl( e^{i\theta} \bfP_c \bmw(t) \bigr), \\
        \bmh_\imag(t) &:= e^{-it\jD} (2i\jD)^{-1} (\pt+i\jD) \Im \bigl( e^{i\theta} \bfP_c \bmw(t) \bigr).   
    \end{aligned}
\end{equation}
This leads to the following decomposition of the flat profiles $\bmf_\real(t)$ and $\bmf_\imag(t)$ associated with the radiation term variables $\bmu_\real(t)$ and $\bmu_\imag(t)$,
\begin{align}
    \bmf_\real(t) &= e^{-it\jD} (2i\jD)^{-1} (\pt + i\jD) \, \Re\bigl( e^{i\theta} \bfP_c \bmv^\flat(t) \bigr) + \bmg^\sharp_\real(t) + \bmh_\real(t), \label{equ:decomposition_bmf_real_flat_profiles_definitions} \\
    \bmf_\imag(t) &= e^{-it\jD} (2i\jD)^{-1} (\pt + i\jD) \, \Im\bigl( e^{i\theta} \bfP_c \bmv^\flat(t) \bigr) + \bmg^\sharp_\imag(t) + \bmh_\imag(t). \label{equ:decomposition_bmf_imag_flat_profiles_definitions}
\end{align}
We do not estimate the first terms on the right-hand sides of \eqref{equ:decomposition_bmf_real_flat_profiles_definitions} and of \eqref{equ:decomposition_bmf_imag_flat_profiles_definitions} directly. Instead, we estimate the flat profiles associated with the auxiliary evolutions $\bmv_\real^\flat(t,x)$ and $\bmv_\imag^\flat(t,x)$,
\begin{equation} \label{equ:definition_bmg_flat_profile}
    \begin{aligned}
        \bmg^\flat_\real(t) &:= e^{-it\jD} (2i\jD)^{-1} (\pt+i\jD) \bmv^\flat_\real(t), \\
        \bmg^\flat_\imag(t) &:= e^{-it\jD} (2i\jD)^{-1} (\pt+i\jD) \bmv^\flat_\imag(t),
    \end{aligned}
\end{equation}
and then transfer these bounds to the flat profiles associated with $\Re\bigl( e^{i\theta} \bfP_c \bmv^\flat(t) \bigr)$ and $\Im\bigl( e^{i\theta} \bfP_c \bmv^\flat(t) \bigr)$ using the transfer estimates from Subsection~\ref{subsec:transfer_estimates}.

We record that the evolution equations for the profiles $\bmf_\real(t)$ and $\bmf_\imag(t)$ are given by
\begin{equation} \label{equ:evol_equation_flat_profiles_bmf_real_imag}
    \begin{aligned}
        \pt \bmf_\real(t) &= e^{-it\jD} (2i\jD)^{-1} \Bigl( - \bfV \bmu_\real(t) + \Re \bigl( e^{i\theta} \bfP_c \bfN_c(t) \bigr) + \Re \bigl( e^{i\theta} \bfP_c \bfN_d(t) \bigr) \Bigr), \\
        \pt \bmf_\imag(t) &= e^{-it\jD} (2i\jD)^{-1} \Bigl( - \bfV \bmu_\imag(t) + \Im \bigl( e^{i\theta} \bfP_c\bfN_c(t) \bigr) + \Im \bigl( e^{i\theta} \bfP_c \bfN_d(t) \bigr) \Bigr).
    \end{aligned}
\end{equation}
Moreover, from \eqref{equ:evol_equ_bmv_real_flat}, \eqref{equ:evol_equ_bmv_imag_flat} we obtain that the evolution equations of the profiles $\bmg^\flat_\real(t)$ and $\bmg^\flat_\imag(t)$ are given by
\begin{equation} \label{equ:evol_equ_bmg_real_imag_flat}
    \begin{aligned}
        \pt \bmg^\flat_\real(t) &= e^{-it\jD} (2i\jD)^{-1} \Re\bigl(e^{i\theta} \bmN_c(t) \bigr), \\
        \pt \bmg^\flat_\imag(t) &= e^{-it\jD} (2i\jD)^{-1} \Im\bigl(e^{i\theta} \bmN_c(t) \bigr). 
    \end{aligned}
\end{equation}

\section{Linear Theory} \label{sec:linear_theory}

In this section we collect the linear theory used for the nonlinear analysis in this work.
We first recall the spectral properties of the linearized operator $\bfM$
and the associated distorted Fourier transform. 
We then record the linear decay estimates for the free Klein-Gordon evolution and for the evolution generated by $\bfM$. 
Next, we state a transference estimate that controls weighted norms of flat profiles for solutions to inhomogeneous Klein--Gordon equations
associated with $\bfM$ whose forcing terms are spatially localized. 
Finally, we record several auxiliary mapping properties of the linearized operator.
All proofs are contained in our companion works \cites{LPPSS1}{LPPSS2}, to which we give precise references below.

\subsection{Spectral and distorted Fourier theory} \label{sec:spectral_fourier}

The relevant spectral properties of $\bfM$ and the description of the associated distorted Fourier transform are summarized in the next proposition. We use $\calM_2(\bbC)$ to denote the space of $2\times2$ complex-valued matrices, and $\bfP_c$ to denote the orthogonal projection onto the continuous spectral subspace of $(L^2_{\mathrm{rad}}(\R^2))^4$ relative to $\bfM$.

\begin{proposition} \label{prop:spectral_fourier}
Let $\bfM$ be the self-adjoint operator defined by \eqref{sysN1compact} on $(L^2_{\mathrm{rad}}(\R^2))^4$. Then:
\begin{itemize}
\item[(i)] The restriction of $\bfM$ to $\operatorname{ran}\bfP_c$ has purely absolutely continuous spectrum $[1,\infty)$. Moreover, the threshold
energy $1$ is neither an eigenvalue nor a resonance of $\bfM$.
\item[(ii)] $\bfM$ has a unique gap eigenvalue of multiplicity two (an internal mode), which we denote by $\lambda^2\in(0,1)$ (numerically $\lambda^2\approx 0.7774$). The associated eigenspace is spanned by
\begin{equation}
    \bmY_1 := \bigl( Y_1, Y_2, 0, 0\bigr)^T, \qquad \bmY_2 := \bigl(0, 0, Y_1, Y_2\bigr)^T, \qquad \bfM \bmY_j = \lambda^2\bmY_j,
\end{equation}
where $Y_1 := \lambda^{-1} U \psi$, $Y_2 := - \lambda^{-1} \psi'$, and where $\psi$ is the $L^2_{r \ud r}$-normalized\footnote{This means $\int_0^\infty \psi(r)^2r\, dr=1$ without a $2\pi$ factor. That $\|\bmY_1\|_{L^2_{r \ud r}} = \|\bmY_2\|_{L^2_{r \ud r}} = 1$ follows from the eigenvalue equation and integration by parts. } radial ground state of $-\Delta + U^2$ with energy $\lambda^2$.
\item [(iii)] There exists a  function $E \colon [0,\infty)\times[0,\infty)\to\calM_2(\bbC)$ such that with $\bfE := {\textstyle \pmat{E&0\\ 0&E}}$, and the distorted Fourier transform $\wtilcalF$ relative to $\bfM$ and its inverse $\wtilcalF^{-1}$ defined as
\begin{equation}
    \wtilcalF[ \bmg ](\freq) \equiv \widetilde{\bmg}(\freq) := \int_0^\infty \overline{\bfE(r,\freq)}^{t}  \bmg(r) r\, \ud r, 
\end{equation}
and
\begin{equation}
    \wtilcalF^{-1}[\bmh](r) := \int_0^\infty \bfE(r,\freq)  \bmh(\freq) \freq\, \ud \freq,
\end{equation}
where $\bmg \in L^2_{r \ud r}(\bbR_+;\bbR^4)$ and $\bmh\in L^2_{k \ud k}(\bbR_+;\bbR^4)$, the following properties are satisfied:
\begin{equation}
    \wtilcalF[ \bfM \bmg ](\freq) = \jap{k}^2 \widetilde{\bmg}(\freq),\qquad \bfP_c = \wtilcalF^{-1} \circ \wtilcalF, \qquad  \bigl\| \widetilde{\bmg} \bigr\|_{L^2_{\freq \ud \freq}(\bbR_+;\bbR^4)} = \bigl\| \bfP_c\,\bmg \bigr\|_{L^2_{r \ud r}(\bbR_+;\bbR^4)}.
\end{equation}
\end{itemize}
\end{proposition}
\begin{proof}
Items (i) and (ii) follow from \cite[Theorem~1.1]{LPPSS1}, while item (iii) is proved in \cite[Theorem~1.2]{LPPSS2}.
\end{proof}

\subsection{Linear decay estimates} \label{sec:linear_decay_estimates}

We next collect the linear decay estimates used in this work, first for the free Klein-Gordon evolution on $\bbR^2$ and then for the evolution generated by $\bfM$.

\subsubsection{Linear decay estimates for the free two-dimensional Klein--Gordon equation without symmetry restrictions}

We use the following estimates for the free Klein-Gordon evolution on $\bbR^2$. Although these estimates are standard, proofs are included in
\cite{LPPSS2} for the reader's convenience. We begin with the dispersive estimate.

\begin{lemma}[Dispersive estimate] \label{lem:flat_dispersive_estimate}
    For every $0 < \delta \ll 1$ there exists $C \geq 1$ such that for any Schwartz function $g \in \calS(\bbR^2)$ and for any $t \in \bbR$,
    \begin{equation}
        \bigl\| e^{it\jD} g \bigr\|_{L^\infty_x(\bbR^2)} \leq \frac{C}{\jt^{1-\delta}} \Bigl( \bigl\| \jxi^2 \nabla_\xi \widehat{g}(\xi) \bigr\|_{L^2_\xi(\bbR^2)} + \bigl\| \jxi^2 \widehat{g}(\xi) \bigr\|_{L^2_\xi(\bbR^2)} \Bigr).
    \end{equation}    
\end{lemma}
\begin{proof}
    See \cite[Lemma~A.4]{LPPSS2}.
\end{proof}

The next two lemmas give dual integrated local energy decay estimates, with Lemma~\ref{lem:flat_dual_ILED_low_freq} exploiting additional vanishing at low frequency. Lemma~\ref{lem:flat_ILED_high_Sobolev_norm} records the corresponding final-state inhomogeneous ILED estimate.

\begin{lemma}[Dual ILED] \label{lem:flat_dual_ILED}
    Let $\kappa > 0$. There exists $C \geq 1$ such that for any $0 \leq t \leq T$, 
    \begin{equation}
        \biggl\| \int_t^T e^{\pm i s \jxi} \widehat{G}(s,\xi) \, \ud s \biggr\|_{L^2_\xi(\bbR^2)} \leq C \bigl\| \jx^{1+\kappa} G(s,x) \bigr\|_{L^2_s([t,T]; L^2_x(\bbR^2))}.
    \end{equation}
\end{lemma}
\begin{proof}
    See \cite[Lemma~A.2]{LPPSS2}.
\end{proof}

\begin{lemma}[Dual ILED with low-frequency vanishing] \label{lem:flat_dual_ILED_low_freq}
    Let $\kappa > 0$. There exists $C \geq 1$ such that for any $0 \leq t \leq T$, 
    \begin{equation}
        \biggl\| \int_t^T \frac{\xi}{\jxi} e^{\pm i s \jxi} \widehat{G}(s,\xi) \, \ud s \biggr\|_{L^2_\xi(\bbR^2)} \leq C \bigl\| \jx^{\frac12+\kappa} G(s,x) \bigr\|_{L^2_s([t,T]; L^2_x(\bbR^2))}.
    \end{equation}
\end{lemma}
\begin{proof}
    See \cite[Lemma~A.1]{LPPSS2}.
\end{proof}

\begin{lemma}[ILED with low-frequency vanishing] \label{lem:flat_ILED_high_Sobolev_norm}
    Let $\kappa > 0$. There exists $C \geq 1$ such that for any $T > 0$,
    \begin{equation} \label{equ:flat_ILED_high_Sobolev_final_time}
    \begin{aligned}
        &\biggl\| \jx^{-\frac12-\kappa} \int_t^T \frac{|D|}{\jD} e^{i(t-s)\jD} G(s) \, \ud s \biggr\|_{L^2_t([0,T]; L^2_x(\bbR^2))} \\
        &\quad \quad \quad \leq C \inf_{G = G_1 + G_2} \Bigl( \bigl\|G_1\bigr\|_{L^1_t([0,T]; L^2_x(\bbR^2))} + \bigl\|\jx^{\frac12+\kappa} G_2\bigr\|_{L^2_t([0,T]; L^2_x(\bbR^2))} \Bigr).
    \end{aligned}
    \end{equation}    
\end{lemma}
\begin{proof}
    See \cite[Lemma~A.3]{LPPSS2}.
\end{proof}

\subsubsection{Linear decay estimates for the Klein-Gordon evolution generated by $\bfM$}

We now recall the linear decay estimates for the Klein-Gordon evolution $e^{\pm it\sqrt{\bfM}} \bfP_c$. For the proofs we refer to \cite{LPPSS2}, noting that the corresponding results in \cite{LPPSS2} are stated for the evolution $e^{\pm it\sqrt{\bfL}}$. Because of the block diagonal structure of $\bfM$, see \eqref{sysN1compact}, the corresponding estimates for $e^{\pm it\sqrt{\bfM}}$ follow directly from the ones for $e^{\pm it\sqrt{\bfL}}$. We begin with the  dispersive decay estimates. 

\begin{proposition} \label{prop:dispersive_estimate_with_potential}
For $j = 0, 1$ and all $t \in \bbR$, it holds that
\begin{equation} \label{eq:Linftydispersiveboundnonsharp2}
    \bigl\| e^{\pm i t\sqrt{\bfM}} \bfM^{-\frac{j}{2}} \bfP_c \bmf \bigr\|_{L^\infty_r(\bbR_+)} \lesssim \jap{t}^{-1} \bigl\| e^{i\theta} \bmf \bigr\|_{W^{4-j,1}_x(\bbR^2)}.   
\end{equation}
Moreover, for any $p\in(2,\infty)$, for $j=0,1$, and all $t\in\bbR$,
    \begin{equation} \label{eq:Linftyweighteddispersivebound1}
        \bigl\| e^{\pm it\sqrt{\bfM}} \bfM^{-\frac{j}{2}} \bfP_c \bmf \bigr\|_{L^\infty_{r}(\bbR_+)} \lesssim \jap{t}^{-1+\frac{2}{p}} \bigl\| \jap{x}\jap{D}^{2-j} \bigl( e^{i\theta}\bmf \bigr) \bigr\|_{L^2_x(\bbR^2)}.
    \end{equation}
\end{proposition}
\begin{proof}
    See \cite[Corollary~3.3]{LPPSS2}.
\end{proof}

We will also use the following integrated local energy decay estimates and their dual versions.

\begin{proposition} \label{prop:ILED}
For any $\gamma_1>\frac{1}{2}$ and $\gamma_2\geq0$ with $\gamma_1+\gamma_2\geq1$ and $(\gamma_1,\gamma_2)\neq(1,0)$, we have 
\begin{equation} \label{equ:ILED}
    \begin{aligned} 
    \bigl\| \langle r\rangle^{-\gamma_1} e^{\pm i s \sqrt{\bfM}} \bfM^{-\frac{\gamma_2}{2}} (\bfM-\bfI)^{\frac{\gamma_2}{2}} \, \bfP_c \bmf \bigr\|_{L^2_s(\bbR; L^2_{r\ud r}(\bbR_+))} \lesssim_{\gamma_1,\gamma_2} \bigl\| \bfP_c \bmf \bigr\|_{L^2_{r\ud r}(\bbR_+)},
    \end{aligned}
\end{equation}
as well as for all $t > 0$,
\begin{equation} \label{equ:dual_ILED}
    \biggl\| \int_0^t e^{\pm is\sqrt{\bfM}} \bfM^{-\frac{\gamma_2}{2}} (\bfM-\bfI)^{\frac{\gamma_2}{2}} \bfP_c \bmF(s,\cdot) \, \ud s \biggr\|_{L^2_{r\ud r}(\bbR_+)} \lesssim_{\gamma_1,\gamma_2} \bigl\| \jap{r}^{\gamma_1} \bfP_c \bmF \bigr\|_{L^2_s([0,t]; L^2_{r\ud r}(\bbR_+))}.
\end{equation}
Moreover, for any $T > 0$ and any $\kappa > 0$ it holds that
\begin{equation} \label{equ:inhomog_ILED}
    \biggl\| \jap{r}^{-\frac32-\kappa} \int_0^t e^{\pm i (t-s) \sqrt{\bfM}} \bfM^{-\frac12} \bfP_c \bmF(s,\cdot) \, \ud s \biggr\|_{L^2_t([0,T]; L^2_{r\ud r}(\bbR_+))} \lesssim \bigl\| \jap{r}^{1+\kappa} \bfP_c \bmF \bigr\|_{L^2_t([0,T]; L^2_{r\ud r}(\bbR_+))}.
\end{equation}
\end{proposition}
\begin{proof}
    See \cite[Lemmas~3.7 and~3.8]{LPPSS2}.
\end{proof}

Finally, the analysis of certain remainder terms in the ODEs for the internal modes requires the following two pointwise-in-time local decay estimates.

\begin{lem}\label{lem:linpaperlocaldecaycutoff1}
Let $\chi_0+\chi_1=1$ be a smooth partition of unity on $[0,\infty)$, where $\chi_0$ is supported in $[0,2]$ and $\chi_1$ is supported in $[1,\infty)$.
Then for any $\bmf$ satisfying $\jap{r}\bmf\in L^2_{r\ud r}$,
\begin{align}
\big\Vert \langle r\rangle^{-1}e^{it\sqrt{\bfM}}\chi_1\big(\sqrt{\bfM-\bfI}\big)\bfP_c\mathbf{f}\big\Vert_{L^2_{r\ud r}}&\lesssim \dfrac{1}{t}\Vert \langle r\rangle \mathbf{f}\Vert_{L^2_{r\ud r}},  \quad t\geq1. 
\end{align}
\end{lem}
\begin{proof}
This follows from \cite[Lemma~3.5]{LPPSS2} with $k_0$ there set to one.
\end{proof}

\begin{lemma}\label{lem:regresbd1}
    Let $\bmf$ satisfy $\jap{r}^2\bmf\in L^2_{r\ud r}$, and define 
    \begin{equation}
        \bigl(\sqrt{\bfM}-2\lambda+i0^+\bigr)^{-1} \bmf := \lim_{\epsilon\to 0^+} \bigl(\sqrt{\bfM}-2\lambda+i\epsilon\bigr)^{-1}\bmf. 
    \end{equation}
    Then for $q = 0, 1$ and any $t \geq 0$,
    \begin{equation*}
        \|\jap{r}^{-2}\partial_r^qe^{it\sqrt{\bfM}}(\sqrt{\bfM}-2\lambda+i0^+)^{-1}\bmf\|_{L^2_{r\ud r}}\lesssim \frac{1}{\jap{t}}\|\jap{r}^2\bmf\|_{L^2_{r\ud r}}.
    \end{equation*}
\end{lemma}
\begin{proof}
    See \cite[Lemma~3.9]{LPPSS2}.
\end{proof}

\subsection{Transference estimate} \label{sec:transference}

The following estimate controls weighted energies of the flat profiles of solutions to Klein--Gordon equations generated by $\bfM$ with spatially
localized forcing. We apply it in this work to the sharp component of the good part of the radiation term and to the bad part of the radiation term.

\begin{proposition} \label{prop:transference}
    Suppose 
    \begin{equation} \label{equ:weighted_estimate_black_box_evol_equ}
        \left\{ \begin{aligned}
            &(\pt^2 + \bfM) \bfP_c \bmv = \bfP_c \bmF, \quad (t,r) \in [0,T] \times [0,\infty), \\
            &\bigl( \bfP_c \bmv(0), \bfP_c \pt \bmv(0) \bigr) = (\bma, \bmb).
        \end{aligned} \right.
    \end{equation}
    Define the flat profiles
    \begin{equation} \label{equ:weighted_estimate_black_box_def_profile}
        \begin{aligned}
            \bmg_\real(t) &:= e^{-it\jD} (2i\jD)^{-1} (\pt + i\jD) \Re\bigl( e^{i\theta} \bfP_c \bmv(t) \bigr), \\ 
            \bmg_\imag(t) &:= e^{-it\jD} (2i\jD)^{-1} (\pt + i\jD) \Im\bigl( e^{i\theta} \bfP_c \bmv(t) \bigr).
        \end{aligned}
    \end{equation}
    Then we have for $0 \leq t \leq T$,
    \begin{equation} \label{equ:weighted_estimate_black_box_bound}
        \begin{aligned}
            &\sum_{\ast \in \{\real, \imag\}} \bigl\| \jxi^2 \nabla_\xi \widehat{\bmg}_\ast(t,\xi) \bigr\|_{L^2_\xi} \\
            &\quad \quad \lesssim \bigl\| \jx e^{i\theta} \bma \bigr\|_{H^2_x} + \bigl\| \jx e^{i\theta} \bmb \bigr\|_{H^1_x} + \bigl\| \jx e^{i\theta} \bmF(0,r) \bigr\|_{L^2_x} + \bigl\| \jx e^{i\theta} (\pt \bmF)(0,r) \bigr\|_{L^2_x} \\ 
            &\quad \quad \quad + \bigl\| \jx e^{i\theta} \bmF(t,r) \bigr\|_{L^2_x} + \sum_{\ast \in \{\real, \imag\}} \Bigl( \bigl\| \bmg_\ast(t) \bigr\|_{H^2_x} + \jt \bigl\| \jD^2 e^{it\jD} \bmg_\ast(t) \bigr\|_{L^{\infty-}_x} \Bigr) \\ 
            &\quad \quad \quad + \sum_{j=0,1} \bigl\| s \cdot \jx \bigl( e^{i\theta} \partial_s^j \bmF(s,r) \bigr) \bigr\|_{L^2_s([0,t]; L^2_x)} + \bigl\| s \cdot \jx \jD \bigl( e^{i\theta} \bmF(s,r) \bigr) \bigr\|_{L^2_s([0,t]; L^2_x)} \\
            &\quad \quad \quad + \sum_{j=1,2} \bigl\| \jap{x} \bigl( e^{i\theta} \partial_s^j \bmF(s,r) \bigr) \bigr\|_{L^1_s([0,t]; L^2_x)} + \bigl\| \jx \jD \bigl( e^{i\theta} \bmF(s,r) \bigr) \bigr\|_{L^1_s([0,t]; L^2_x)}.
        \end{aligned}
    \end{equation}
\end{proposition}
\begin{proof}
    See \cite[Theorem~1.3]{LPPSS2}.
\end{proof}

\subsection{Auxiliary estimates} \label{sec:mapping_properties}

In the analysis of the evolution equations for the sharp component of the good part $\bmv^\sharp(t)$ and for the bad part $\bmw(t)$, it will occasionally be technically convenient to work with the operator
\begin{equation}
    \bfM_{\mathrm{nr}} := \bigl(-\Delta + 1) \bbI_{4 \times 4} + \bfV
\end{equation}
with $\bfV$ defined in \eqref{equ:definition_bfV}. The latter is related to the linearized operator $\bfM$ via the identity
\begin{equation} \label{equ:relation_bfM_bfMnr}
    e^{i\theta} \bigl( \bfM f \bigr)(r) = \bfM_{\mathrm{nr}} \bigl( e^{i\theta} f(r) \bigr).
\end{equation}
We shall make use of the following mapping properties of the operator $\bfM_{\mathrm{nr}}$.

\begin{lemma} \label{lem:mapping_properties_bfMnr}
    For any integer $k \geq 1$, any $\gamma \in \bbR$, and any function $\bmh$ of the form $\bmh(r,\theta)=e^{i\theta}\underline{\bmh}(r)$,
    \begin{align}
        \bigl\| \jx^\gamma \bfM_{\nr}^{\frac{k}{2}} \bmh \bigr\|_{L^2_x(\bbR^2)} &\lesssim \bigl\| \jx^\gamma \bmh \bigr\|_{H^k_x(\bbR^2)}, \label{equ:bfMnr_mapping_Hk} \\ 
        \bigl\| \jx^\gamma \bmh \bigr\|_{H^k_x(\bbR^2)} &\lesssim \bigl\| \jx^\gamma \bmh \bigr\|_{L^2_x(\bbR^2)} + \bigl\| \jx^\gamma \bfM_{\nr}^{\frac{k}{2}} \bmh \bigr\|_{L^2_x(\bbR^2)}, \label{equ:bfMnr_mapping_Mk} \\ 
        \bigl\| \jD^{-1} \bfM_{\nr}^{\frac12} \bmh \bigr\|_{L^\infty_x(\bbR^2)} &\lesssim \bigl\| \bmh \bigr\|_{L^\infty_x(\bbR^2)}. \label{equ:bfMnr_mapping_jDinverse_Linfty}
    \end{align}
\end{lemma}
\begin{proof}
    The estimates \eqref{equ:bfMnr_mapping_Hk}, \eqref{equ:bfMnr_mapping_Mk} are proved in \cite[Lemma~A.6]{LPPSS2}, while the last estimate \eqref{equ:bfMnr_mapping_jDinverse_Linfty} is proved in \cite[Lemma~A.5]{LPPSS2}. 
\end{proof}

\section{Bounds for the Temporal Component} \label{sec:working_title_eta0_bounds}

In this section we establish pointwise decay estimates and high Sobolev norm
estimates for the temporal component $\eta_0$ and its time derivatives
that will be used in the subsequent nonlinear analysis. We treat
separately the decomposition
\begin{equation}
    \eta_0 = \eta_{0,c}^{\flat} + \eta_{0,c}^{\sharp} + \eta_{0,d}
\end{equation}
introduced in Subsections~\ref{subsec:good_bad_decomposition} and~\ref{subsec:flat_sharp_decomposition}. 
The flat component $\eta_{0,c}^{\flat}$ of the good part of the temporal component is estimated directly from the flat elliptic
equations \eqref{equ:eta0c_flat_equation_rewritten} and
\eqref{equ:pteta0c_flat_equation_rewritten}. 
The sharp component $\eta_{0,c}^{\sharp}$ of the good part of the temporal component is then controlled in terms of $\eta_{0,c}^{\flat}$ through the localized elliptic relation \eqref{equ:eta0c_sharp_elliptic_equ},
while the bad part $\eta_{0,d}$ of the temporal component is estimated in spatially weighted norms by exploiting
the spatial localization of the internal mode components and the off-diagonal exponential decay of the Green's function of $-\Delta+U^2$.

We first record the estimates for the flat component $\eta_{0,c}^\flat(t)$ of the good part of the temporal component.

\begin{proposition} \label{prop:eta0c_bounds}
    Suppose the bootstrap assumptions \eqref{equ:bootstrap_assumption_z}, \eqref{equ:bootstrap_assumption_weighted_f}, \eqref{equ:bootstrap_assumption_Linfty_f}, \eqref{equ:bootstrap_assumption_HN_f}, \eqref{equ:bootstrap_assumption_HN_ILED_f} are in place.  
    Then we have the decay estimates
    \begin{equation} \label{equ:eta0c_Linfty_bound}
        \bigl\| \jD^2 \eta_{0,c}^\flat(t) \bigr\|_{L^\infty_x(\bbR^2)} \lesssim \left\{ \begin{aligned}
                            &\varepsilon^{2-20\delta} \jt^{-1+7\delta}, \quad &&0 \leq t \leq \varepsilon^{-2}, \\ 
                            &t^{-2+3\delta},  &&\varepsilon^{-2} \leq t \leq T,
                         \end{aligned} \right.  
    \end{equation}
    and for $1 \leq l \leq 3$,
    \begin{equation} \label{equ:pteta0c_Linfty_bound}
        \bigl\| \jD^2 \pt^l \eta_{0,c}^\flat(t) \bigr\|_{L^\infty_x(\bbR^2)} \lesssim \left\{ \begin{aligned}
                            &\varepsilon^{2-20\delta} \jt^{-1+7\delta}, \quad &&0 \leq t \leq \varepsilon^{-2}, \\
                            &t^{-2+3\delta},  &&\varepsilon^{-2} \leq t \leq T.
                         \end{aligned} \right.  
    \end{equation}
    Moreover, the following high Sobolev norm bounds hold
    \begin{equation} \label{equ:eta0c_HN_bound}
        \bigl\| \jD^{N+1} \eta_{0,c}^\flat(t) \bigr\|_{L^2_x(\bbR^2)} \lesssim \left\{ \begin{aligned}
                            &\varepsilon^{2-10\delta} \jt^{-\frac12-3\delta}, \quad &&0 \leq t \leq \varepsilon^{-2}, \\
                            &\varepsilon t^{-1+2\delta},  && \varepsilon^{-2} \leq t \leq T,
                         \end{aligned} \right.  
    \end{equation}
    and 
    \begin{equation} \label{equ:pteta0c_HN_bound}
        \bigl\| \jD^{N} \pt \eta_{0,c}^\flat(t) \bigr\|_{L^2_x(\bbR^2)} \lesssim \left\{ \begin{aligned}
                            &\varepsilon^{2-10\delta} \jt^{-\frac12-3\delta}, \quad &&0 \leq t \leq \varepsilon^{-2}, \\
                            &\varepsilon t^{-1+2\delta},  &&\varepsilon^{-2} \leq t \leq T.
                         \end{aligned} \right.  
    \end{equation}
\end{proposition}

We next state the corresponding estimates for the bad part $\eta_{0,d}$ of the temporal component. 
By construction, this component is generated by interactions involving at least one internal mode component and is therefore spatially
localized; accordingly, its bounds are formulated in weighted norms.

\begin{proposition} \label{prop:eta0d_bounds}
    Suppose the bootstrap assumptions \eqref{equ:bootstrap_assumption_z}, \eqref{equ:bootstrap_assumption_weighted_f}, \eqref{equ:bootstrap_assumption_Linfty_f}, \eqref{equ:bootstrap_assumption_HN_f}, \eqref{equ:bootstrap_assumption_HN_ILED_f} are in place.
    Then we have the weighted decay estimates 
    \begin{equation} \label{equ:eta0d_Linfty_bound}
        \bigl\| \jx^{10} \jD^2 \eta_{0,d}(t) \bigr\|_{L^\infty_x(\bbR^2)} \lesssim \left\{ \begin{aligned}
                            &\varepsilon^2 + \varepsilon^{2-10\delta} \jt^{-\frac12-3\delta}, \quad &&0 \leq t \leq \varepsilon^{-2}, \\
                            &t^{-1},  &&\varepsilon^{-2} \leq t \leq T,
                         \end{aligned} \right.  
    \end{equation}
    and for $1 \leq l \leq 3$,
    \begin{equation} \label{equ:pteta0d_Linfty_bound}
        \bigl\| \jx^{10} \jD^{2} (\pt^l \eta_{0,d})(t) \bigr\|_{L^\infty_x(\bbR^2)} \lesssim \left\{ \begin{aligned}
                            &\varepsilon^{2-10\delta} \jt^{-\frac12-3\delta}, \quad &&0 \leq t \leq \varepsilon^{-2}, \\
                            &t^{-\frac32+2\delta},  &&\varepsilon^{-2} \leq t \leq T.
                         \end{aligned} \right.  
    \end{equation}
    Moreover, the following weighted high Sobolev norm bounds hold
    \begin{equation} \label{equ:eta0d_HN_bound}
        \bigl\| \jx^{10} \jD^{N+1} \eta_{0,d}(t) \bigr\|_{L^2_x(\bbR^2)} \lesssim \varepsilon^2, \quad 0 \leq t \leq T,
    \end{equation}
    and 
    \begin{equation} \label{equ:pteta0d_HN_bound}
        \bigl\| \jx^{10} \jD^{N} (\pt \eta_{0,d})(t) \bigr\|_{L^2_x(\bbR^2)} \lesssim \varepsilon^2, \quad 0 \leq t \leq T.
    \end{equation}
\end{proposition}

\begin{remark}
    The spatial weight $\jx^{10}$ in the weighted estimates in Proposition~\ref{prop:eta0d_bounds} provides sufficient spatial localization to close all nonlinear estimates involving the bad part of the temporal component $\eta_{0,d}(t)$ and its time derivatives $\pt^l \eta_{0,d}(t)$, $1 \leq l \leq 3$. In principle, we could establish analogous weighted estimates for spatial weights $\jx^M$ for any integer $M \geq 1$. However, the size of the initial data has to be chosen depending on the integer $M$, so we have to fix a weight that suffices to close all nonlinear estimates.
\end{remark}

To pass from the flat component to the full good part of the temporal component, we use the elliptic equation
\eqref{equ:eta0c_sharp_elliptic_equ}. We record the resulting comparison estimates in the following lemma.
\begin{lemma} \label{lem:eta0c_sharp_in_terms_of_flat}
    For $0 \leq t \leq T$ we have 
    \begin{align}
        \eta_{0,c}^\sharp(t) &= (-\Delta+U^2)^{-1} \bigl( (1-U^2) \eta_{0,c}^\flat(t) \bigr), \label{equ:eta0c_sharp_in_terms_of_flat} \\ 
        \pt^l \eta_{0,c}^\sharp(t) &= (-\Delta+U^2)^{-1} \bigl( (1-U^2) \pt^l \eta_{0,c}^\flat(t) \bigr), \quad 1 \leq l \leq 3. \label{equ:pteta0c_sharp_in_terms_of_flat}
    \end{align}
    It follows that for $0 \leq t \leq T$,
    \begin{align}
        \bigl\| \jD^2 \eta_{0,c}^\sharp(t) \bigr\|_{L^\infty_x(\bbR^2)} &\lesssim \bigl\| \eta_{0,c}^\flat(t) \bigr\|_{L^\infty_x(\bbR^2)}, \label{equ:eta0c_sharp_Linfty_bound} \\ 
        \bigl\| \jD^2 \pt^l \eta_{0,c}^\sharp(t) \bigr\|_{L^\infty_x(\bbR^2)} &\lesssim \bigl\| \pt^l \eta_{0,c}^\flat(t) \bigr\|_{L^\infty_x(\bbR^2)}, \quad 1 \leq l \leq 3, \label{equ:pteta0c_sharp_Linfty_bound} \\
        \bigl\| \jD^k \eta_{0,c}^\sharp(t) \bigr\|_{L^2_x(\bbR^2)} &\lesssim \bigl\| \jD^{k-2} \eta_{0,c}^\flat(t) \bigr\|_{L^2_x(\bbR^2)}, \quad k \geq 2, \label{equ:eta0c_sharp_Hk_bound} \\
        \bigl\| \jD^k \pt^l \eta_{0,c}^\sharp(t) \bigr\|_{L^2_x(\bbR^2)} &\lesssim \bigl\| \jD^{k-2} \pt^l \eta_{0,c}^\flat(t) \bigr\|_{L^2_x(\bbR^2)}, \quad k \geq 2, \quad 1 \leq l \leq 3. \label{equ:pteta0c_sharp_Hk_bound}
    \end{align}
\end{lemma}
\begin{proof}
    The identities \eqref{equ:eta0c_sharp_in_terms_of_flat} and \eqref{equ:pteta0c_sharp_in_terms_of_flat} follow directly from the defining elliptic equation \eqref{equ:eta0c_sharp_elliptic_equ} for the sharp component of the good part of the temporal component.
    The asserted estimates \eqref{equ:eta0c_sharp_Linfty_bound}, \eqref{equ:pteta0c_sharp_Linfty_bound}, \eqref{equ:eta0c_sharp_Hk_bound}, \eqref{equ:pteta0c_sharp_Hk_bound} are then an immediate consequence of the off-diagonal exponential decay of the Green's function of $(-\Delta+U^2)^{-1}$ and the fact that the sharp component gains two derivatives with respect to the flat component by elliptic regularity.
\end{proof}

With these comparison estimates in hand, we turn to the proof of Proposition~\ref{prop:eta0c_bounds}.

\begin{proof}[Proof of Proposition~\ref{prop:eta0c_bounds}]
    We begin with the decay estimate \eqref{equ:eta0c_Linfty_bound} for $\eta_{0,c}^\flat(t)$.
    From \eqref{equ:eta0c_flat_equation_rewritten} we obtain for $0 \leq t \leq T$,
    \begin{equation}
        \begin{aligned}
            &\bigl\| \jD^2 \eta_{0,c}^\flat(t) \bigr\|_{L^\infty_x} \\
            &= \Bigl\| \jD^2 (-\Delta+1)^{-1} \Bigl( -2 \eta_{0,c} \bigl( \underline{\Phi}_{\real} u_{1,\real} + \underline{\Phi}_{\imag} u_{1,\imag} \bigr) - \bigl( u_{3,\real} \pt u_{1,\real} + u_{3,\imag} \pt u_{1,\imag} \bigr) \\ 
            &\quad \quad \quad \quad \quad \quad \quad \quad \quad \quad + u_{1,\real} \pt u_{3,\real} + u_{1,\imag} \pt u_{3,\imag} - \eta_{0,c} \bigl( u_{1,\real}^2 + u_{1,\imag}^2 + u_{3,\real}^2 + u_{3,\imag}^2 \bigr) \Bigr) \Bigr\|_{L^\infty_x} \\ 
            &\lesssim \Bigl( \bigl\|\eta_{0,c}^\flat(t)\bigr\|_{L^\infty_x} + \bigl\|\eta_{0,c}^\sharp(t)\bigr\|_{L^\infty_x} \Bigr) \Bigl( \|\underline{\Phi}_\real\|_{L^\infty_x} \bigl\| e^{it\jD} \bmf_\real(t) \bigr\|_{L^\infty_x} + \|\underline{\Phi}_\imag\|_{L^\infty_x} \bigl\| e^{it\jD} \bmf_\imag(t) \bigr\|_{L^\infty_x} \Bigr) \\ 
            &\quad + \bigl\| e^{it\jD} \bmf_\real(t) \bigr\|_{L^\infty_x} \bigl\| \jD e^{it\jD} \bmf_\real(t) \bigr\|_{L^\infty_x} + \bigl\| e^{it\jD} \bmf_\imag(t) \bigr\|_{L^\infty_x} \bigl\| \jD e^{it\jD} \bmf_\imag(t) \bigr\|_{L^\infty_x} \\ 
            &\quad + \Bigl( \bigl\|\eta_{0,c}^\flat(t)\bigr\|_{L^\infty_x} + \bigl\|\eta_{0,c}^\sharp(t)\bigr\|_{L^\infty_x} \Bigr) \Bigl( \bigl\| e^{it\jD} \bmf_\real(t) \bigr\|_{L^\infty_x}^2 + \bigl\| e^{it\jD} \bmf_\imag(t) \bigr\|_{L^\infty_x}^2 \Bigr).
        \end{aligned}
    \end{equation}
    Note that by \eqref{equ:eta0c_sharp_Linfty_bound} we have $\bigl\| \eta_{0,c}^\sharp(t) \bigr\|_{L^\infty_x} \lesssim \bigl\| \eta_{0,c}^\flat(t) \bigr\|_{L^\infty_x}$.
    Hence, upon reabsorbing the terms with a factor $\bigl\| \eta_{0,c}^\flat(t) \bigr\|_{L^\infty_x}$ into the left-hand side of the preceding inequality,
    and using the Gagliardo-Nirenberg-Sobolev inequality~\eqref{equ:GNS}, we find that
    \begin{equation}
        \begin{aligned}
            \bigl\| \jD^2 \eta_{0,c}^\flat(t) \bigr\|_{L^\infty_x} 
            &\lesssim \bigl\| e^{it\jD} \bmf_\real(t) \bigr\|_{L^\infty_x} \bigl\| \jD e^{it\jD} \bmf_\real(t) \bigr\|_{L^\infty_x} + \bigl\{ \text{similar or better terms}\bigr\}. 
        \end{aligned}
    \end{equation}
    Invoking \eqref{equ:GNSell} we find that this is bounded by $t^{-2+3\delta}$ for long times $\varepsilon^{-2} \leq t \leq T$ and by $\varepsilon^{2-20\delta} \jt^{-1+7\delta}$ for short times $0 \leq t \leq \varepsilon^{-2}$.

    The derivation of the decay estimate \eqref{equ:pteta0c_Linfty_bound} for $\pt^l \eta_{0,c}^\flat(t)$, $1 \leq l \leq 3$, proceeds similarly, using \eqref{equ:pteta0c_flat_equation_rewritten} and \eqref{equ:pteta0c_flat_equation_leading_order_term_rewritten}. We note that in the case of $l=3$, on the right-hand side of \eqref{equ:pteta0c_flat_equation_rewritten} we have to reinsert the evolution equation for the components of $\pt^2 \bmu_\real$ and of $\pt^2 \bmu_\imag$ obtained from \eqref{equ:system_evol_equations_u_z}. Since the reinserted nonlinearities are of higher order, this only produces similar or better terms.

    Next, we turn to the high Sobolev norm bound \eqref{equ:eta0c_HN_bound} for $\eta_{0,c}^\flat(t)$.
    From \eqref{equ:eta0c_flat_equation_rewritten} we obtain for $0 \leq t \leq T$, using standard product estimates,
    \begin{equation}
        \begin{aligned}
            &\bigl\| \jD^{N+1} \eta_{0,c}^\flat(t) \bigr\|_{L^2_x} \\
            &\lesssim \Bigl\| \jD^{N-1} \Bigl( \eta_{0,c} \bigl( \underline{\Phi}_{\real} u_{1,\real} + \underline{\Phi}_{\imag} u_{1,\imag} \bigr) \Bigr) \Bigr\|_{L^2_x} + \Bigl\| \jD^{N-1} \bigl( u_{3,\real} \pt u_{1,\real} + u_{3,\imag} \pt u_{1,\imag} \bigr) \Bigr\|_{L^2_x} \\ 
            &\quad \quad + \Bigl\| \jD^{N-1} \bigl( u_{1,\real} \pt u_{3,\real} + u_{1,\imag} \pt u_{3,\imag} \bigr) \Bigr\|_{L^2_x} + \Bigl\| \jD^{N-1} \Bigl( \eta_{0,c} \bigl( u_{1,\real}^2 + u_{1,\imag}^2 + u_{3,\real}^2 + u_{3,\imag}^2 \bigr) \Bigr) \Bigr\|_{L^2_x} \\             
            &\lesssim \Bigl( \bigl\| \jD^{N-1} \eta_{0,c}^\flat(t) \bigr\|_{L^2_x} + \bigl\| \jD^{N-1} \eta_{0,c}^\sharp(t) \bigr\|_{L^2_x} \Bigr) \Bigl( \| \underline{\Phi}_{\real} \|_{L^\infty_x} \bigl\| e^{it\jD} \bmf_\real(t) \bigr\|_{L^\infty_x} + \| \underline{\Phi}_{\imag} \|_{L^\infty_x} \bigl\| e^{it\jD} \bmf_\imag(t) \bigr\|_{L^\infty_x} \Bigr) \\ 
            &\quad + \Bigl( \bigl\| \eta_{0,c}^\flat(t) \bigr\|_{L^\infty_x} + \bigl\| \eta_{0,c}^\sharp(t) \bigr\|_{L^\infty_x} \Bigr) \\
            &\quad \quad \quad \quad \times \Bigl( \bigl\| \jD^{N-1} \underline{\Phi}_{\real} \bigr\|_{L^\infty_x} \bigl\| \jD^{N-1} \bmf_\real(t) \bigr\|_{L^2_x} + \bigl\| \jD^{N-1} \underline{\Phi}_{\imag} \|_{L^\infty_x} \bigl\| \jD^{N-1} \bmf_\imag(t) \bigr\|_{L^2_x} \Bigr) \\ 
            &\quad + \bigl\| \jD^N \bmf_\real(t) \bigr\|_{L^2_x} \bigl\| \jD e^{it\jD} \bmf_\real(t) \bigr\|_{L^\infty_x} + \bigl\| \jD^N \bmf_\imag(t) \bigr\|_{L^2_x} \bigl\| \jD e^{it\jD} \bmf_\imag(t) \bigr\|_{L^\infty_x} \\
            &\quad + \Bigl( \bigl\| \jD^{N-1} \eta_{0,c}^\flat(t) \bigr\|_{L^2_x} + \bigl\| \jD^{N-1} \eta_{0,c}^\sharp(t) \bigr\|_{L^2_x} \Bigr) \Bigl( \bigl\| e^{it\jD} \bmf_\real(t) \bigr\|_{L^\infty_x}^2 + \bigl\| e^{it\jD} \bmf_\imag(t) \bigr\|_{L^\infty_x}^2 \Bigr) \\ 
            &\quad + \Bigl( \bigl\| \eta_{0,c}^\flat(t) \bigr\|_{L^\infty_x} + \bigl\| \eta_{0,c}^\sharp(t) \bigr\|_{L^\infty_x} \Bigr) \\
            &\quad \quad \quad \quad \times \Bigl( \bigl\| \jD^{N-1} \bmf_\real(t) \bigr\|_{L^2_x} \bigl\| e^{it\jD} \bmf_\real(t) \bigr\|_{L^\infty_x} + \bigl\| \jD^{N-1} \bmf_\imag(t) \bigr\|_{L^2_x} \bigl\| e^{it\jD} \bmf_\imag(t) \bigr\|_{L^\infty_x} \Bigr).
        \end{aligned}
    \end{equation}
    Recall from \eqref{equ:eta0c_sharp_Linfty_bound} and \eqref{equ:eta0c_sharp_Hk_bound} that we have $\bigl\|\eta_{0,c}^\sharp(t)\bigr\|_{L^\infty_x} \lesssim \bigl\|\eta_{0,c}^\flat(t)\bigr\|_{L^\infty_x}$ and $\bigl\| \jD^{N-1} \eta_{0,c}^\sharp(t) \bigr\|_{L^2_x} \lesssim \bigl\| \jD^{N-3} \eta_{0,c}^\flat(t)\bigr\|_{L^2_x}$.
    Thus, upon reabsorbing the terms with a factor $\bigl\| \jD^{N-3} \eta_{0,c}^\flat(t) \bigr\|_{L^2_x}$ into the left-hand side of the preceding inequality, and using the previously established $L^\infty_x$-bounds \eqref{equ:eta0c_Linfty_bound}, 
    we conclude
    \begin{equation}
        \begin{aligned}
            \bigl\| \jD^{N+1} \eta_{0,c}^\flat(t) \bigr\|_{L^2_x}
            &\lesssim \bigl\| \jD^N \bmf_\real(t) \bigr\|_{L^2_x} \bigl\| \jD e^{it\jD} \bmf_\real(t) \bigr\|_{L^\infty_x} + \bigl\{ \text{similar or better terms}\bigr\}.
        \end{aligned}
    \end{equation}
    Hence, using \eqref{equ:GNSell}, we infer that the last line is bounded by $\varepsilon t^{-1+2\delta}$ for long times $\varepsilon^{-2} \leq t \leq T$ and by $\varepsilon^{2-10\delta} \jt^{-\frac12-3\delta}$ for short times $0 \leq t \leq \varepsilon^{-2}$.
    
    The proof of the $H^{N}_x$-bound \eqref{equ:pteta0c_HN_bound} for $\pt \eta_{0,c}^\flat(t)$ proceeds similarly, using \eqref{equ:pteta0c_flat_equation_rewritten} and \eqref{equ:pteta0c_flat_equation_leading_order_term_rewritten}.
    This finishes the proof of Proposition~\ref{prop:eta0c_bounds}.
\end{proof}

We now turn to the weighted estimates for the bad part of the temporal component and prove Proposition~\ref{prop:eta0d_bounds}.

\begin{proof}[Proof of Proposition~\ref{prop:eta0d_bounds}]
    We begin with the weighted decay estimate \eqref{equ:eta0d_Linfty_bound} for $\eta_{0,d}(t)$.
    From \eqref{equ:eta0d_equation} we obtain that
    \begin{equation} \label{equ:eta0d_greens_fct_expression}
        \begin{aligned}
            \eta_{0,d} &= - (-\Delta+U^2)^{-1} \Bigl[ Y_1^2 \bigl( (\pt z_1) z_2 - z_1 (\pt z_2) \bigr) + Y_1 \bigl( z_2 (\pt u_1) - z_1 (\pt u_3) \bigr) \\ 
            &\qquad \qquad \qquad \qquad \quad + Y_1 \bigl( (\pt z_1) u_3 - (\pt z_2) u_1 \bigr) + 2 Y_1 U \eta_{0,c} z_1 + 2 U \eta_{0,d} (z_1 Y_1 + u_1) \\
            &\qquad \qquad \qquad \qquad \quad + \eta_{0,c} \Bigl( Y_1^2 \bigl( z_1^2 + z_2^2 \bigr) + Y_1 \bigl( 2 z_1 u_1 + 2 z_2 u_3 \bigr) \Bigr) \\
            &\qquad \qquad \qquad \qquad \quad + \eta_{0,d} \Bigl( \bigl( z_1 Y_1 + u_1 \bigr)^2 + \bigl( z_2 Y_1 + u_3 \bigr)^2 \Bigr) \Bigr].
        \end{aligned}
    \end{equation}
    Using the off-diagonal exponential decay of the Green's function of $(-\Delta + U^2)^{-1}$, it follows that for $0 \leq t \leq T$,
    \begin{equation}
        \begin{aligned}
            &\bigl\| \jx^{10} \jD^2 \eta_{0,d}(t) \bigr\|_{L^\infty_x} \\
            &\lesssim \bigl\| \jx^{10} Y_1^2 \bigr\|_{L^\infty_x} \, \frakz(t) + \bigl\| \jx^{10} Y_1 \bigr\|_{L^\infty_x} \, \frakz(t)^{\frac12} \Bigl( \bigl\| \jD e^{it\jD} \bmf_\real(t) \bigr\|_{L^\infty_x} + \bigl\| \jD e^{it\jD} \bmf_\imag(t) \bigr\|_{L^\infty_x} \Bigr) \\ 
            &\quad + \bigl\| \jx^{10} Y_1 \bigr\|_{L^\infty_x} \Bigl( \bigl\| \eta_{0,c}^{\flat}(t) \bigr\|_{L^\infty_x} + \bigl\| \eta_{0,c}^{\sharp}(t) \bigr\|_{L^\infty_x} \Bigr) \, \frakz(t)^{\frac12} + \bigl\| \jx^{10} Y_1^2 \bigr\|_{L^\infty_x} \Bigl( \bigl\| \eta_{0,c}^{\flat}(t) \bigr\|_{L^\infty_x} + \bigl\| \eta_{0,c}^{\sharp}(t) \bigr\|_{L^\infty_x} \Bigr) \, \frakz(t) \\
            &\quad + \bigl\| \jx^{10} \eta_{0,d}(t) \bigr\|_{L^\infty_x} \Bigl( \frakz(t)^{\frac12} + \bigl\| e^{it\jD} \bmf_\real(t) \bigr\|_{L^\infty_x} + \bigl\| e^{it\jD} \bmf_\imag(t) \bigr\|_{L^\infty_x} \Bigr) \\ 
            &\quad  + \bigl\| \jx^{10} Y_1 \bigr\|_{L^\infty_x} \Bigl( \bigl\| \eta_{0,c}^{\flat}(t) \bigr\|_{L^\infty_x} + \bigl\| \eta_{0,c}^{\sharp}(t) \bigr\|_{L^\infty_x} \Bigr) \, \frakz(t)^{\frac12} \Bigl( \bigl\| e^{it\jD} \bmf_\real(t) \bigr\|_{L^\infty_x} + \bigl\| e^{it\jD} \bmf_\imag(t) \bigr\|_{L^\infty_x} \Bigr) \\
            &\quad + \bigl\| \jx^{10} \eta_{0,d}(t) \bigr\|_{L^\infty_x} \Bigl( \frakz(t) + \bigl\| e^{it\jD} \bmf_\real(t) \bigr\|_{L^\infty_x}^2 + \bigl\| e^{it\jD} \bmf_\imag(t) \bigr\|_{L^\infty_x}^2 \Bigr).
        \end{aligned}
    \end{equation}
    Reabsorbing the terms with a factor $\bigl\| \jx^{10} \eta_{0,d}(t) \bigr\|_{L^\infty_x}$ into the left-hand side of the preceding inequality, we find that
    \begin{equation}
        \begin{aligned}
            \bigl\| \jx^{10} \jD^2 \eta_{0,d}(t) \bigr\|_{L^\infty_x} 
            &\lesssim \frakz(t) + \frakz(t)^{\frac12} \bigl\| \jD e^{it\jD} \bmf_\real(t) \bigr\|_{L^\infty_x} + \bigl\{\text{similar or better terms}\bigr\}. 
        \end{aligned}
    \end{equation}
    Invoking \eqref{equ:GNSell}, we find that the last line is bounded by $t^{-1}$ for long times $\varepsilon^{-2} \leq t \leq T$ and by $\varepsilon^2 + \varepsilon^{2-10\delta} \jt^{-\frac12-3\delta}$ for short times $0 \leq t \leq \varepsilon^{-2}$.

    Next, we turn to the weighted high Sobolev norm bound \eqref{equ:eta0d_HN_bound} for $\eta_{0,d}(t)$.
    Again using the off-diagonal exponential decay of the Green's function of $(-\Delta + U^2)^{-1}$, we conclude from \eqref{equ:eta0d_greens_fct_expression} for $0 \leq t \leq T$ that 
    \begin{equation}
        \begin{aligned}
            &\bigl\| \jx^{10} \jD^{N+1} \eta_{0,d}(t) \bigr\|_{L^2_x} \\
            &\lesssim \bigl\| \jx^{10} Y_1^2 \bigr\|_{H^{N-1}_x} \, \frakz(t) + \bigl\| \jx^{10} Y_1 \bigr\|_{W^{N-1,\infty}_x} \, \frakz(t)^{\frac12} \bigl( \|\bmf_\real(t)\|_{H^N_x} + \|\bmf_\imag(t)\|_{H^N_x} \bigr) \\
            &\quad + \bigl\| \jx^{10} Y_1 \bigr\|_{W^{N-1,\infty}_x} \Bigl( \bigl\|\eta_{0,c}^\flat(t)\bigr\|_{H^{N-1}_x} + \bigl\|\eta_{0,c}^\sharp(t)\bigr\|_{H^{N-1}_x} \Bigr) \, \frakz(t)^{\frac12} \\ 
            &\quad + \bigl\| \jx^{10} \jD^{N-1} \eta_{0,d}(t) \bigr\|_{L^2_x} \Bigl( \frakz(t)^{\frac12} + \bigl\| \jD^{N-1} e^{it\jD} \bmf_{\real}(t) \bigr\|_{L^\infty_x} + \bigl\| \jD^{N-1} e^{it\jD} \bmf_{\imag}(t) \bigr\|_{L^\infty_x} \Bigr) \\ 
            &\quad + \Bigl( \bigl\|\eta_{0,c}^\flat(t)\bigr\|_{H^{N-1}_x} + \bigl\|\eta_{0,c}^\sharp(t)\bigr\|_{H^{N-1}_x} \Bigr) \bigl\|\jx^{10} Y_1^2\bigr\|_{W^{N-1,\infty}_x} \, \frakz(t) \\ 
            &\quad + \Bigl( \bigl\|\eta_{0,c}^\flat(t)\bigr\|_{H^{N-1}_x} + \bigl\|\eta_{0,c}^\sharp(t)\bigr\|_{H^{N-1}_x} \Bigr) \bigl\|\jx^{10} Y_1 \bigr\|_{W^{N-1,\infty}_x} \, \frakz(t)^{\frac12} \Bigl( \bigl\| \bmf_{\real}(t) \bigr\|_{H^N_x} + \bigl\| \bmf_{\imag}(t) \bigr\|_{H^N_x} \Bigr) \\ 
            &\quad + \bigl\| \jx^{10} \jD^{N-1} \eta_{0,d}(t) \bigr\|_{L^\infty_x} \Bigl( \frakz(t) + \|\bmf_\real(t)\|_{H^{N}_x}^2 + \|\bmf_\imag(t)\|_{H^N_x}^2 \Bigr).
        \end{aligned}
    \end{equation}
    Hence, upon reabsorbing the term with a factor of $\bigl\| \jx^{10} \jD^{N-1} \eta_{0,d}(t) \bigr\|_{L^2_x}$ into the left-hand side of the preceding inequality, we find that
    \begin{equation}
        \begin{aligned}
            \bigl\| \jx^{10} \jD^{N+1} \eta_{0,d}(t) \bigr\|_{L^2_x} \lesssim \frakz(t) + \frakz(t)^{\frac12} \|\bmf_\real(t)\|_{H^N_x} + \bigl\{ \text{similar or better terms}\bigr\}.
        \end{aligned}
    \end{equation}   
    By \eqref{equ:bootstrap_assumption_z}, \eqref{equ:bootstrap_assumption_HN_f}, the last line is bounded by $\varepsilon^2$ for all times $0 \leq t \leq T$.

    It remains to establish the weighted decay estimates \eqref{equ:pteta0d_Linfty_bound} for $\pt^l \eta_{0,d}$, $1 \leq l \leq 3$, and the weighted high Sobolev norm bound \eqref{equ:pteta0d_HN_bound} for $\pt \eta_{0,d}$. From \eqref{equ:eta0d_equation} we obtain that $\pt \eta_{0,d}$ satisfies the elliptic equation
    \begin{equation}
        \begin{aligned}
            &(-\Delta + U^2) \pt \eta_{0,d} \\
            &= Y_1^2 \bigl( z_1 \pt^2 z_2 - z_2 \pt^2 z_1 \bigr) + Y_1 \bigl( u_1 (\pt^2 z_2) - u_3 (\pt^2 z_1) \bigr) + Y_1 \bigl( z_1 \pt^2 u_3 - z_2 \pt^2 u_1 \bigr) \\ 
            &\quad - 2U (\pt \eta_{0,d}) \bigl( z_1 Y_1 + u_1 \bigr) - 2 U \eta_{0,d} \bigl( (\pt z_1) Y_1 + \pt u_1 \bigr) \\ 
            &\quad - 2U \bigl( \pt \eta_{0,c}^\flat + \pt \eta_{0,c}^\sharp \bigr) z_1 Y_1 - 2U \bigl( \eta_{0,c}^\flat + \eta_{0,c}^\sharp \bigr) (\pt z_1) Y_1 \\ 
            &\quad - (\pt \eta_{0,d}) \Bigl( \bigl( z_1 Y_1 + u_1 \bigr)^2 + \bigl( z_2 Y_1 + u_3 \bigr)^2 \Bigr) \\
            &\quad - \eta_{0,d} \Bigl( 2 \bigl( z_1 Y_1 + u_1 \bigr) \bigl( (\pt z_1) Y_1 + (\pt u_1) \bigr) + 2 \bigl( z_2 Y_1 + u_3 \bigr) \bigl( (\pt z_2) Y_1 + (\pt u_3) \bigr) \Bigr) \\ 
            &\quad - (\pt \eta_{0,c}) \Bigl( \bigl(z_1^2 + z_2^2\bigr) Y_1 + 2 Y_1 \bigl( z_1 u_1 + z_2 u_3 \bigr) \Bigr) \\ 
            &\quad - \eta_{0,c} \Bigl( \bigl( 2 z_1 (\pt z_1) + 2 z_2 (\pt z_2) \bigr) Y_1 + 2 Y_1 \bigl( (\pt z_1) u_1 + z_1 (\pt u_1) + (\pt z_2) u_3 + z_2 (\pt u_3) \bigr) \Bigr).
        \end{aligned}
    \end{equation}
    Inserting the evolution equations \eqref{equ:system_evol_equations_u_z} for $\bmu$, $z_1$, $z_2$ in the first line on the right-hand side, we find that
    \begin{equation} \label{equ:pteta0d_equation}
        \begin{aligned}
            &(-\Delta + U^2) \pt \eta_{0,d} \\
            &= Y_1^2 \bigl( z_1 \langle \bmY_2, \bmN \rangle - z_2 \langle \bmY_1, \bmN \rangle \bigr) + Y_1 \bigl( -\lambda^2 u_1 z_2 + u_1 \langle \bmY_2, \bmN \rangle + \lambda^2 u_3 z_1 - u_3 \langle \bmY_1, \bmN \rangle \bigr) \\
            &\quad - \bigl( \cos(\theta) Y_1 \bigr) z_1 \bigl( -\Delta + V_1 + 1 \bigr) u_{3,\real} - \bigl( \sin(\theta) Y_1 \bigr) z_1 \bigl( -\Delta + V_1 + 1 \bigr) u_{3,\imag} \\
            &\quad + \bigl( \cos(\theta) Y_1 \bigr) z_2 \bigl( -\Delta + V_1 + 1 \bigr) u_{1,\real} + \bigl( \sin(\theta) Y_1 \bigr) z_2 \bigl( -\Delta + V_1 + 1 \bigr) u_{1,\imag} \\ 
            &\quad + 2U' Y_1 \bigl( z_1 u_4 - z_2 u_2 \bigr) + Y_1 \bigl( z_1 (\bfP_c \bmN)_3 - z_2 (\bfP_c \bmN)_1 \bigr) \\ 
            &\quad - 2U (\pt \eta_{0,d}) \bigl( z_1 Y_1 + u_1 \bigr) - 2 U \eta_{0,d} \bigl( (\pt z_1) Y_1 + \pt u_1 \bigr) \\ 
            &\quad - 2U \bigl( \pt \eta_{0,c}^\flat + \pt \eta_{0,c}^\sharp \bigr) z_1 Y_1 - 2U \bigl( \eta_{0,c}^\flat + \eta_{0,c}^\sharp \bigr) (\pt z_1) Y_1 \\ 
            &\quad - (\pt \eta_{0,d}) \Bigl( \bigl( z_1 Y_1 + u_1 \bigr)^2 + \bigl( z_2 Y_1 + u_3 \bigr)^2 \Bigr) \\
            &\quad - \eta_{0,d} \Bigl( 2 \bigl( z_1 Y_1 + u_1 \bigr) \bigl( (\pt z_1) Y_1 + (\pt u_1) \bigr) + 2 \bigl( z_2 Y_1 + u_3 \bigr) \bigl( (\pt z_2) Y_1 + (\pt u_3) \bigr) \Bigr) \\ 
            &\quad - \bigl( (\pt \eta_{0,c}^\flat) + (\pt \eta_{0,c}^\sharp) \bigr) \Bigl( \bigl(z_1^2 + z_2^2\bigr) Y_1 + 2 Y_1 \bigl( z_1 u_1 + z_2 u_3 \bigr) \Bigr) \\ 
            &\quad - \bigl( \eta_{0,c}^\flat + \eta_{0,c}^\sharp \bigr) \Bigl( \bigl( 2 z_1 (\pt z_1) + 2 z_2 (\pt z_2) \bigr) Y_1 + 2 Y_1 \bigl( (\pt z_1) u_1 + z_1 (\pt u_1) + (\pt z_2) u_3 + z_2 (\pt u_3) \bigr) \Bigr).
        \end{aligned}
    \end{equation}
    Using the off-diagonal exponential decay of the Green's function of $(-\Delta+U^2)^{-1}$, the decay estimate $|\langle \bmY_j, \bmN(t) \rangle| \lesssim \frakz(t) + \bigl\{\text{better terms}\bigr\}$, as well as Lemma~\ref{lem:eta0c_sharp_in_terms_of_flat}, it follows that
    \begin{equation}
        \begin{aligned}
            &\bigl\| \jx^{10} \jD^2 (\pt \eta_{0,d})(t) \bigr\|_{L^\infty_x} \\
            &\lesssim \bigl\| \jx^{10} Y_1^2 \bigr\|_{L^\infty_x} \, \frakz(t)^{\frac12} \Bigl( \frakz(t) + \bigl\{ \text{better terms} \bigr\} \Bigr) \\ 
            &\quad + \bigl\| \jx^{10} Y_1 \bigr\|_{L^\infty_x} \Bigl( \frakz(t)^{\frac12} \bigl( \bigl\| e^{it\jD} \bmf_\real(t) \bigr\|_{L^\infty_x} + \bigl\| e^{it\jD} \bmf_\imag(t) \bigr\|_{L^\infty_x} \bigr) + \bigl\{ \text{better terms} \bigr\} \Bigr) \\
            &\quad + \bigl\| \jx^{10} Y_1 \bigr\|_{L^\infty_x} \, \frakz(t)^{\frac12} \Bigl( \bigl\| \jD^2 e^{it\jD} \bmf_\real(t) \bigr\|_{L^\infty_x} + \bigl\| \jD^2 e^{it\jD} \bmf_\imag(t) \bigr\|_{L^\infty_x} \Bigr) \\ 
            &\quad + \bigl\| \jx^{10} Y_1 \bigr\|_{L^\infty_x} \frakz(t)^{\frac12} \Bigl( \bigl\| e^{it\jD} \bmf_\real(t) \bigr\|_{L^\infty_x} + \bigl\|  e^{it\jD} \bmf_\imag(t) \bigr\|_{L^\infty_x} \Bigr) \\ 
            &\quad + \bigl\| \jx^{10} Y_1 \bigr\|_{L^\infty_x} \frakz(t)^{\frac12} \Bigl( \frakz(t) + \bigl\{ \text{better terms} \bigr\} \Bigr) \\ 
            &\quad + \bigl\| \jx^{10} (\pt \eta_{0,d})(t) \bigr\|_{L^\infty_x} \Bigl( \frakz(t)^{\frac12} + \bigl\| e^{it\jD} \bmf_\real(t) \bigr\|_{L^\infty_x} + \bigl\|  e^{it\jD} \bmf_\imag(t) \bigr\|_{L^\infty_x} \Bigr) \\ 
            &\quad + \bigl\| \jx^{10} \eta_{0,d}(t) \bigr\|_{L^\infty_x} \Bigl( \frakz(t)^{\frac12} + \bigl\| \jD e^{it\jD} \bmf_\real(t) \bigr\|_{L^\infty_x} + \bigl\| \jD e^{it\jD} \bmf_\imag(t) \bigr\|_{L^\infty_x} \Bigr) \\ 
            &\quad + \bigl\| \jx^{10} Y_1 \bigr\|_{L^\infty_x} \, \frakz(t)^{\frac12} \Bigl( \bigl\| \pt \eta_{0,c}^\flat(t) \bigr\|_{L^\infty_x} + \bigl\|  \eta_{0,c}^\flat(t) \bigr\|_{L^\infty_x} \Bigr) \\ 
            &\quad + \bigl\| \jx^{10} (\pt \eta_{0,d})(t) \bigr\|_{L^\infty_x} \Bigl( \frakz(t) + \bigl\| e^{it\jD} \bmf_\real(t) \bigr\|_{L^\infty_x}^2 + \bigl\|  e^{it\jD} \bmf_\imag(t) \bigr\|_{L^\infty_x}^2 \Bigr) \\ 
            &\quad + \bigl\| \jx^{10} \eta_{0,d}(t) \bigr\|_{L^\infty_x} \Bigl( \frakz(t) + \bigl\| \jD e^{it\jD} \bmf_\real(t) \bigr\|_{L^\infty_x}^2 + \bigl\|  \jD e^{it\jD} \bmf_\imag(t) \bigr\|_{L^\infty_x}^2 \Bigr) \\ 
            &\quad + \bigl\| \jx^{10} Y_1 \bigr\|_{L^\infty_x} \Bigl( \bigl\| \eta_{0,c}^\flat(t) \bigr\|_{L^\infty_x} + \bigl\| (\pt \eta_{0,c}^\flat)(t) \bigr\|_{L^\infty_x} \Bigr) \Bigl( \frakz(t) + \frakz(t)^{\frac12} \bigl\| \jD e^{it\jD} \bmf(t) \bigr\|_{L^\infty_x} \Bigr).
        \end{aligned}
    \end{equation}
    Thus, upon reabsorbing the terms with a factor of $\bigl\| \jx^{10} (\pt \eta_{0,d})(t) \bigr\|_{L^\infty_x}$ into the left-hand side of the preceding inequality, using the $L^\infty_x$ bound for $\eta_{0,d}(t)$ from \eqref{equ:eta0d_Linfty_bound} as well as the bounds from Proposition~\ref{prop:eta0c_bounds}, we conclude that
    \begin{equation}
        \begin{aligned}
            \bigl\| \jx^{10} \jD^2 (\pt \eta_{0,d})(t) \bigr\|_{L^\infty_x} 
            &\lesssim \frakz(t)^{\frac32} + \frakz(t)^{\frac12} \bigl\| \jD^2 e^{it\jD} \bmf_\real(t) \bigr\|_{L^\infty_x} + \bigl\{ \text{similar or better terms} \bigr\}.
        \end{aligned}
    \end{equation}
    By \eqref{equ:GNSell} this is bounded by $t^{-\frac32 + 2\delta}$ for long times $\varepsilon^{-2} \leq t \leq T$ and by $\varepsilon^{2-10\delta} \jt^{-\frac12-3\delta}$ for short times $0 \leq t \leq \varepsilon^{-2}$.    
    The derivation of the corresponding weighted $L^\infty_x$ decay estimates for $(\pt^l \eta_{0,d})(t)$, $l = 2, 3$, proceeds analogously, starting from the equation \eqref{equ:pteta0d_equation}. We note that upon differentiating \eqref{equ:pteta0d_equation} in time once or twice, we again have to insert the evolution equations for $\pt^2 z_j$ and for the components of $\pt^2 \bmu_\ast$, $\ast \in \{\real, \imag\}$, obtained from \eqref{equ:system_evol_equations_u_z}, but this only produces similar or better terms.

    It remains to establish the weighted high Sobolev norm bound \eqref{equ:pteta0d_HN_bound} for $(\pt \eta_{0,d})(t)$. 
    From \eqref{equ:pteta0d_equation} we infer, using the off-diagonal exponential decay of the Green's function of $(-\Delta+U^2)^{-1}$, the decay estimate $|\langle \bmY_j, \bmN(t) \rangle| \lesssim \frakz(t) + \bigl\{\text{better terms}\bigr\}$, as well as Lemma~\ref{lem:eta0c_sharp_in_terms_of_flat},
    \begin{equation}
        \begin{aligned}
            &\bigl\| \jx^{10} \jD^N (\pt \eta_{0,d})(t) \bigr\|_{L^2_x} \\
            &\lesssim \bigl\| \jx^{10} Y_1^2 \bigr\|_{H^{N-2}_x} \, \frakz(t)^{\frac12} \Bigl( \frakz(t) + \bigl\{\text{better terms}\bigr\} \Bigr) \\ 
            &\quad + \bigl\| \jx^{10} Y_1 \bigr\|_{W^{N-2,\infty}_x} \bigl( \|\bmf_\real(t)\|_{H^{N-2}_x} + \|\bmf_\imag(t)\|_{H^{N-2}_x} \bigr) \Bigl( \frakz(t) + \bigl\{\text{better terms}\bigr\} \Bigr) \\
            &\quad + \bigl\| \jx^{10} e^{i\theta} Y_1 \bigr\|_{W^{N-2,\infty}_x} \, \frakz(t)^{\frac12} \bigl( \|\bmf_\real(t)\|_{H^{N}_x} + \|\bmf_\imag(t)\|_{H^{N}_x} \bigr) \\ 
            &\quad + \bigl\| \jx^{10} Y_1 \bigr\|_{W^{N-2,\infty}_x} \, \frakz(t)^{\frac12} \bigl( \|\bmf_\real(t)\|_{H^{N-2}_x} + \|\bmf_\imag(t)\|_{H^{N-2}_x} \bigr) \\
            &\quad + \bigl\| \jx^{10} Y_1 \bigr\|_{W^{N-2,\infty}_x} \Bigl( \frakz(t)^{\frac32} + \bigl\{\text{similar or better terms}\bigr\} \Bigr) \\ 
            &\quad + \bigl\| \jx^{10} \jD^{N-2} (\pt \eta_{0,d})(t) \bigr\|_{L^2_x} \Bigl( \frakz(t)^{\frac12} + \|\bmf_\real(t)\|_{H^{N}_x} + \|\bmf_\imag(t)\|_{H^{N}_x} \Bigr) \\
            &\quad + \bigl\| \jx^{10} \jD^{N} \eta_{0,d}(t) \bigr\|_{L^2_x} \Bigl( \frakz(t) + \|\bmf_\real(t)\|_{H^{N}_x} + \|\bmf_\imag(t)\|_{H^{N}_x} \Bigr) \\ 
            &\quad + \bigl\| \jx^{10} Y_1 \bigr\|_{W^{N-2,\infty}_x} \frakz(t)^{\frac12} \Bigl( \bigl\| (\pt \eta_{0,c}^\flat)(t) \bigr\|_{H^N_x} + \bigl\| \eta_{0,c}^\flat(t)\bigr\|_{H^N_x} \Bigr) \\ 
            &\quad + \bigl\| \jx^{10} \jD^{N-2} (\pt \eta_{0,d})(t) \bigr\|_{L^2_x} \Bigl( \frakz(t) + \|\bmf_\real(t)\|_{H^{N}_x}^2 + \|\bmf_\imag(t)\|_{H^{N}_x}^2 \Bigr) \\ 
            &\quad + \bigl\| \jx^{10} \jD^{N-2} \eta_{0,d}(t) \bigr\|_{L^2_x} \Bigl( \frakz(t) + \|\bmf_\real(t)\|_{H^{N}_x}^2 + \|\bmf_\imag(t)\|_{H^{N}_x}^2 \Bigr) \\ 
            &\quad + \bigl\| \jx^{10} Y_1 \bigr\|_{W^{N-2,\infty}_x} \Bigl( \bigl\| \eta_{0,c}^\flat(t) \bigr\|_{H^N_x} + \bigl\| (\pt \eta_{0,c}^\flat)(t) \bigr\|_{H^N_x} \Bigr) \Bigl( \frakz(t) + \frakz(t)^{\frac12} \bigl( \|\bmf_\real(t)\|_{H^{N}_x} + \|\bmf_\imag(t)\|_{H^{N}_x} \bigr) \Bigr).
        \end{aligned}
    \end{equation}
    Once more reabsorbing the terms with a factor of $\bigl\| \jx^{10} \jD^{N-2} (\pt \eta_{0,d})(t) \bigr\|_{L^2_x}$ into the left-hand side of the preceding inequality, using the bound \eqref{equ:eta0d_HN_bound} as well as the bounds from Proposition~\ref{prop:eta0c_bounds}, we find that
    \begin{equation}
        \begin{aligned}
            \bigl\| \jx^{10} \jD^N (\pt \eta_{0,d})(t) \bigr\|_{L^2_x} 
            &\lesssim \frakz(t)^{\frac12} \|\bmf_\real(t)\|_{H^N_x} + \bigl\{\text{similar or better terms}\bigr\}.
        \end{aligned}
    \end{equation}
    By \eqref{equ:bootstrap_assumption_z}, \eqref{equ:bootstrap_assumption_HN_f}, the last line is bounded by $\varepsilon^2$ for all times $0 \leq t \leq T$.
    This finishes the proof of Proposition~\ref{prop:eta0d_bounds}.    
\end{proof}

\section{Bounds for the Flat Component of the Good Part} \label{sec:working_title_flat_bounds}

The purpose of this section is to derive weighted energy, dispersive, high Sobolev, and ILED bounds for the auxiliary flat profiles $\bmg^\flat_\real(t)$ and $\bmg^\flat_\imag(t)$ introduced in \eqref{equ:definition_bmg_flat_profile}. Recall that these profiles are associated with the real-valued flat Klein--Gordon evolutions $\bmv^\flat_\real(t)$ and $\bmv^\flat_\imag(t)$, which solve \eqref{equ:evol_equ_bmv_real_flat} and \eqref{equ:evol_equ_bmv_imag_flat} on the bootstrap interval $[0,T]$.
The estimates obtained here are auxiliary estimates outside symmetry. They are transferred in Subsection~\ref{subsec:transfer_estimates}
to the first angular modes of the flat component of the good part of the radiation term,
which enter the decompositions of the full flat profiles $\bmf_\real(t)$ and $\bmf_\imag(t)$ in
\eqref{equ:decomposition_bmf_real_flat_profiles_definitions} and
\eqref{equ:decomposition_bmf_imag_flat_profiles_definitions}.

We remind the reader that the reason for imposing final data on the flat component is the subsequent
analysis of the sharp component. By \eqref{equ:evol_equ_bmv_sharp}, the
sharp component is forced by $-\bfP_c(\bfV \bmv^\flat)$, where $\bfV$ is spatially localized. To close the estimates for $\bmv^\sharp$, it is important that $\bmv^\flat$ exhibit improved decay coming directly from the nonlinear contributions rather than merely behave like a generic linear Klein--Gordon wave. Our final-state formulation is designed to furnish this nonlinear decay mechanism.

\subsection{Summary of results and overview} \label{subsec:summary_gflat_bounds}

We begin by summarizing the main results of this section.
The weighted energy estimate in the first proposition below is stated for the normal-form-renormalized profile
\[
    \bmg^\flat_\ast(t)-\bmB_\ast[\bmf,\bmf](t),
    \qquad \ast\in\{\real,\imag\},
\]
rather than for $\bmg^\flat_\ast(t)$ itself. The quadratic expression
$\bmB_\ast[\bmf,\bmf](t)$ collects all the boundary terms generated by the normal form
transformations of the non-spatially localized quadratic interactions. Its
precise definition is given in \eqref{equ:definition_bmB_boundary_terms}.
The point of the final-state prescription
\eqref{equ:data_at_T_bmv_real_flat}--\eqref{equ:data_at_T_bmv_imag_flat}
is that the renormalized profile, not the raw profile, has zero final data
at time $T$.

\begin{proposition} \label{prop:weighted_energies_g}
    Suppose the bootstrap assumptions \eqref{equ:bootstrap_assumption_z}, \eqref{equ:bootstrap_assumption_weighted_f}, \eqref{equ:bootstrap_assumption_Linfty_f}, \eqref{equ:bootstrap_assumption_HN_f}, \eqref{equ:bootstrap_assumption_HN_ILED_f} are in place.  
    Then we have for $\ast \in \{\real, \imag\}$,
    \begin{equation}
        \Bigl\| \jxi^2 \jap{\nabla_\xi} \Bigl( \widehat{\bmg}^\flat_{\ast}(t,\xi) - \widehat{\calF}\bigl[ \bmB_\ast[\bmf, \bmf](t) \bigr](\xi) \Bigr) \Bigr\|_{L^2_\xi(\bbR^2)} \leq C_2 \left\{ \begin{aligned}
                            &\varepsilon^{1-6\delta}, \quad &&0 \leq t \leq \varepsilon^{-2}, \\
                            &t^{-\frac12+3\delta},  &&\varepsilon^{-2} \leq t \leq T,
                         \end{aligned} \right.  
    \end{equation}
    and    
    \begin{equation}
        \Bigl\| \jxi^2 \jap{\nabla_\xi} \widehat{\calF}\bigl[ \bmB_\ast[\bmf, \bmf](t) \bigr](\xi) \Bigr\|_{L^2_\xi(\bbR^2)} \leq C_2 \left\{ \begin{aligned}
                            &\varepsilon^{2-30\delta} \max\bigl\{ \jt^{\frac12-4\delta}, \varepsilon^{1+20\delta} t^{1-3\delta} \bigr\}, \quad &&0 \leq t \leq \varepsilon^{-2}, \\
                            &\varepsilon t^{2\delta},  &&\varepsilon^{-2} \leq t \leq T,
                         \end{aligned} \right.   
    \end{equation}    
    where the boundary terms $\bmB_\ast[\bmf, \bmf](t)$ are defined in \eqref{equ:definition_bmB_boundary_terms} and where $C_2 \equiv C_2(C_0)$ depends on the size of the constant $C_0$ in the bootstrap assumptions \eqref{equ:bootstrap_assumption_weighted_f}, \eqref{equ:bootstrap_assumption_Linfty_f}, \eqref{equ:bootstrap_assumption_HN_f}, \eqref{equ:bootstrap_assumption_HN_ILED_f}.
\end{proposition}

The next proposition establishes improved dispersive decay estimates for the evolutions of the auxiliary flat profiles $\bmg^\flat_\real(t)$ and $\bmg^\flat_\imag(t)$.

\begin{proposition} \label{prop:dispersive_decay_g}
    Suppose the bootstrap assumptions \eqref{equ:bootstrap_assumption_z}, \eqref{equ:bootstrap_assumption_weighted_f}, \eqref{equ:bootstrap_assumption_Linfty_f}, \eqref{equ:bootstrap_assumption_HN_f}, \eqref{equ:bootstrap_assumption_HN_ILED_f} are in place.  
    Then we have for $\ast \in \{\real, \imag\}$,
    \begin{equation} \label{equ:dispersive_decay_g}
        \bigl\| e^{it\jD} \bmg^\flat_\ast(t)\bigr\|_{L^\infty_x(\bbR^2)} 
        \leq C_2 \left\{ \begin{aligned}
                            &\varepsilon^{1-6\delta} \jt^{-1+\delta}, \quad  &&0 \leq t \leq \varepsilon^{-2}, \\
                            &t^{-\frac32+4\delta},  &&\varepsilon^{-2} \leq t \leq T,
                         \end{aligned} \right.
    \end{equation}    
    where $C_2 \equiv C_2(C_0)$ depends on the size of the constant $C_0$ in the bootstrap assumptions \eqref{equ:bootstrap_assumption_weighted_f}, \eqref{equ:bootstrap_assumption_Linfty_f}, \eqref{equ:bootstrap_assumption_HN_f}, \eqref{equ:bootstrap_assumption_HN_ILED_f}.
\end{proposition}

We also obtain a uniform-in-time bound on a high Sobolev norm of the auxiliary flat profiles.

\begin{proposition} \label{prop:HN_g}
    Suppose the bootstrap assumptions \eqref{equ:bootstrap_assumption_z}, \eqref{equ:bootstrap_assumption_weighted_f}, \eqref{equ:bootstrap_assumption_Linfty_f}, \eqref{equ:bootstrap_assumption_HN_f}, \eqref{equ:bootstrap_assumption_HN_ILED_f} are in place. 
    Then we have for $\ast \in \{\real, \imag\}$ and $0 \leq t \leq T$,
    \begin{equation}
        \bigl\| \bmg^\flat_\ast(t)\bigr\|_{H^N_x(\bbR^2)} \leq C_2 \varepsilon^{\frac32-2\delta},
    \end{equation}
    where $C_2 \equiv C_2(C_0)$ depends on the size of the constant $C_0$ in the bootstrap assumptions \eqref{equ:bootstrap_assumption_weighted_f}, \eqref{equ:bootstrap_assumption_Linfty_f}, \eqref{equ:bootstrap_assumption_HN_f}, \eqref{equ:bootstrap_assumption_HN_ILED_f}.
\end{proposition}

Finally, we establish an integrated local energy decay estimate at high Sobolev regularity for the evolutions of the auxiliary flat profiles.

\begin{proposition} \label{prop:HN_ILED_g}
    Suppose the bootstrap assumptions \eqref{equ:bootstrap_assumption_z}, \eqref{equ:bootstrap_assumption_weighted_f}, \eqref{equ:bootstrap_assumption_Linfty_f}, \eqref{equ:bootstrap_assumption_HN_f}, \eqref{equ:bootstrap_assumption_HN_ILED_f} are in place.  
    Then we have for some $0 < \kappa \ll 1$ and for $\ast \in \{\real, \imag\}$,
    \begin{equation} \label{equ:HN_ILED_g}
        \sum_{1 \leq k \leq N} \bigl\| \jx^{-\frac12-\kappa} |D|^k e^{it\jD} \bmg^\flat_\ast(t)\bigr\|_{L^2_t([0,T]; L^2_x(\bbR^2))} \leq C_2 \varepsilon^{\frac32-2\delta},
    \end{equation} 
    where $C_2 \equiv C_2(C_0)$ depends on the size of the constant $C_0$ in the bootstrap assumptions \eqref{equ:bootstrap_assumption_weighted_f}, \eqref{equ:bootstrap_assumption_Linfty_f}, \eqref{equ:bootstrap_assumption_HN_f}, \eqref{equ:bootstrap_assumption_HN_ILED_f}.
\end{proposition}

The main difficulty in the derivation of the above bounds comes from the quadratic interactions without spatial localization in the nonlinearities $\Re\bigl(e^{i\theta} \bmN_c\bigr)$ and $\Im\bigl(e^{i\theta} \bmN_c\bigr)$ of the auxiliary flat Klein-Gordon equations \eqref{equ:evol_equ_bmv_real_flat} and \eqref{equ:evol_equ_bmv_imag_flat}. 
According to the classification in Corollary~\ref{cor:types_of_nonlinearities}, the only such contributions are those of types (1), (2), (5), and (6). 
These interactions cannot be estimated directly in the profile equations and must instead be treated by normal form transformations, i.e., by integrating by parts in time in the corresponding Duhamel integrals. 
The normal form transformations convert the non-spatially localized quadratic interactions into effectively cubic or at least spatially localized quadratic terms, but they also generate quadratic boundary terms at the two endpoints of the interval $[t,T]$.
This endpoint structure explains the choice of the final data in
\eqref{equ:data_at_T_bmv_real_flat}--\eqref{equ:data_at_T_bmv_imag_flat}.
Ideally, one would like to impose zero final data for the flat component at
time $t=T$. 
However, if one imposed the raw condition $\bmg^\flat_\ast(T)=0$, $\ast\in\{\real,\imag\}$, then the normal form identity would leave a terminal quadratic boundary contribution $\bmB_\ast[\bmf,\bmf](T)$. 
The available weighted energy bounds for these quadratic boundary expressions are not compatible with the desired improved dispersive decay estimate for the evolution of the renormalized profile.
Instead, we incorporate the terminal quadratic boundary term into the final state data. Equivalently, the quantity which has zero final data at time $T$
is not the raw profile $\bmg^\flat_\ast$, but the normal-form-renormalized
profile
\begin{equation}
    \bmg^\flat_\ast(t)-\bmB_\ast[\bmf,\bmf](t), 
    \qquad \ast\in\{\real,\imag\}.
\end{equation}

We now spell this out for the profile $\bmg^\flat_\real(t)$. The other profile $\bmg^\flat_\imag(t)$ can be treated exactly in the same way, and we omit the details. 
For times $0 \leq t \leq T$, we obtain from \eqref{equ:evol_equ_bmg_real_imag_flat} that
\begin{equation}
    \begin{aligned}
        \bmg^\flat_\real(t)
        &= \bmg^\flat_\real(T) - \int_t^T \partial_s \bmg^\flat_\real(s) \, \ud s \\
        &= \bmg^\flat_\real(T)
        - \int_t^T e^{-is\jD} (2i\jD)^{-1}
        \Re\bigl( e^{i\theta} \bmN_c(s) \bigr) \, \ud s.
    \end{aligned}
\end{equation}
After applying the normal form transformations to all non-spatially localized
quadratic terms in the preceding Duhamel integral, we obtain the schematic
identity
\begin{equation}
    \begin{aligned}
        &- \int_t^T e^{-is\jD} (2i\jD)^{-1}
        \Re\bigl( e^{i\theta} \bmN_c(s) \bigr) \, \ud s \\
        &\quad =
        - \bmB_\real[\bmf,\bmf](T)
        + \bmB_\real[\bmf,\bmf](t)
        + \int_t^T e^{-is\jD} (2i\jD)^{-1}
        \bmN_{\real,c}^{\text{renorm}}(s) \, \ud s.
    \end{aligned}
\end{equation}
Here $\bmB_\real[\bmf,\bmf]$ denotes the sum of the quadratic boundary
operators produced by the normal form transformations, while
$\bmN_{\real,c}^{\text{renorm}}$ denotes the corresponding renormalized
remainder, consisting of effectively cubic or at least spatially localized quadratic terms. The
precise expression for $\bmB_\real[\bmf,\bmf]$ is recorded in \eqref{equ:definition_bmB_boundary_terms} in 
Subsection~\ref{subsec:quadratic_boundary_terms} below.

The final data \eqref{equ:data_at_T_bmv_real_flat} are chosen so that
\begin{equation}
    \bmg^\flat_\real(T)=\bmB_\real[\bmf,\bmf](T).
\end{equation}
Indeed, if $h=\bmB_\real[\bmf,\bmf](T)$ and
\begin{equation}
    \bmv^\flat_\real(T)
    =2\Re\bigl(e^{iT\jD}h\bigr), \qquad
    \pt\bmv^\flat_\real(T)
    =-2\jD\Im\bigl(e^{iT\jD}h\bigr),
\end{equation}
then, by the definition \eqref{equ:definition_bmg_flat_profile},
\begin{equation}
    e^{-iT\jD}(2i\jD)^{-1}(\pt+i\jD)\bmv^\flat_\real(T)=h.
\end{equation}
Thus, \eqref{equ:data_at_T_bmv_real_flat} is equivalent to
$\bmg^\flat_\real(T)=\bmB_\real[\bmf,\bmf](T)$. Consequently, the terminal
quadratic boundary contribution cancels identically, and we arrive at the
final-state representation for the renormalized profile 
\begin{equation} \label{equ:final_state_representation_renormalized_profile}
    \begin{aligned}
        \bmg^\flat_\real(t)-\bmB_\real[\bmf,\bmf](t)
        =
        \int_t^T e^{-is\jD} (2i\jD)^{-1}
        \bmN_{\real,c}^{\text{renorm}}(s) \, \ud s.
    \end{aligned}
\end{equation}

We stress that the cancellation in \eqref{equ:final_state_representation_renormalized_profile} is essential for the weighted energy estimate in Proposition~\ref{prop:weighted_energies_g} and for the dispersive estimate in Proposition~\ref{prop:dispersive_decay_g}.
Instead, in the proofs of the high Sobolev norm bounds in Proposition~\ref{prop:HN_g} and of the ILED bounds in Proposition~\ref{prop:HN_ILED_g}, we also use normal form transformations for all non-spatially localized quadratic terms, but the resulting boundary terms can be estimated directly also at the terminal time $T$.

\medskip 

We close this subsection by describing the organization of the rest of Section~\ref{sec:working_title_flat_bounds}. 
Since the estimates for $\bmg_\real^\flat(t)$ and $\bmg_\imag^\flat(t)$ are completely analogous, in what follows we only focus on the estimates for $\bmg_\real^\flat(t)$.
In Subsection~\ref{subsec:quadratic_boundary_terms} we specify the precise expressions for the boundary terms $\bmB_\real[\bmf,\bmf]$ and $\bmB_\imag[\bmf,\bmf]$.
In Subsection~\ref{subsec:preparations} we collect the multilinear estimates and phase bounds used throughout the space-time resonance analysis of the nonlinearities in the remainder of this section.
In Subsection~\ref{subsec:proofs_propositions_gflat} we establish the proofs of Propositions~\ref{prop:weighted_energies_g}, \ref{prop:HN_g}, and
\ref{prop:HN_ILED_g} by systematically estimating, type by type, the contributions listed in Corollary~\ref{cor:types_of_nonlinearities}. 
The dispersive bound of Proposition~\ref{prop:dispersive_decay_g} is then derived in Subsection~\ref{subsec:proof_dispersive_decay_g}, and the transfer estimates are stated in Subsection~\ref{subsec:transfer_estimates}.

\subsection{Quadratic boundary terms} \label{subsec:quadratic_boundary_terms}

In this subsection we record the precise expressions for the (vectorial) quadratic boundary terms $\bmB_\real[\bmf, \bmf](t)$ and $\bmB_\imag[\bmf, \bmf](t)$ that arise in the remainder of this section from the normal form analysis of the quadratic nonlinearities without spatial localization.
For $\ast \in \bigl\{ \real, \imag \bigr\}$ and times $0 \leq t \leq T$ these boundary terms are given by
\begin{equation} \label{equ:definition_bmB_boundary_terms}
    \begin{aligned}
        \bmB_\ast[\bmf, \bmf](t) = \Bigl( B_{\ast, 1}[\bmf, \bmf](t), B_{\ast, 2}[\bmf, \bmf](t), B_{\ast, 3}[\bmf, \bmf](t), B_{\ast, 4}[\bmf, \bmf](t) \Bigr)^{T}.
    \end{aligned}
\end{equation}
The four components correspond to the four components $\alpha,\zeta,\beta,\mu$ of the vectorial radiation term. The explicit formulas for the components are specified below in \eqref{equ:definition_Breal1}, \eqref{equ:definition_Bimag1}, \eqref{equ:definition_Brealimag2}, \eqref{equ:definition_Brealimag3}, and \eqref{equ:definition_Brealimag4}.

Only the classes (1), (2), (5), and (6) in Corollary~\ref{cor:types_of_nonlinearities} contribute to
$\bmB_\ast[\bmf,\bmf]$. The remaining classes either carry spatial localization or are at least cubic, and therefore do not require a normal form analysis.
Correspondingly, we introduce four families of quadratic boundary operators:
\begin{equation} \label{equ:families_of_boundary_operators}
\calB_{\kappa_1\kappa_2}^{\pvdots,j}, \qquad
\calB_{\kappa_1\kappa_2}^{\delta_0,j}, \qquad
\calB_{\kappa_1\kappa_2}^{\delta_0,j,\eta_0}, \qquad
\calB_{\kappa_1\kappa_2}^{\pvdots,j,\eta_0}, 
\end{equation}
where $j \in \{1, 2\}$ and $\kappa_1, \kappa_2 \in \{\pm\}$.
The first two correspond to the non-localized quadratic terms of types
(1) and (2). The latter two arise from the leading quadratic part of
$\partial_t\eta_{0,c}$ in the terms of types (5) and (6).
All of these quadratic boundary operators involve the phase function
\begin{equation}
    \Psi_{\kappa_1\kappa_2}(\xi,\xi_1,\xi_2) := -\jxi + \kappa_1 \jap{\xi_1} + \kappa_2 \jap{\xi_2}, \quad \kappa_1, \kappa_2 \in \{ \pm \}, \quad \xi, \xi_1, \xi_2 \in \bbR^2.
\end{equation}
The bilinear operators below are evaluated at an arbitrary time $s\in[0,T]$. In the examples following the definitions, $t$ denotes the lower endpoint of the Duhamel integral, so that $0\leq t\leq s\leq T$.

For the contributions of nonlinear terms of type (1) we use the bilinear operators
\begin{equation}
    \begin{aligned}
        &\calB_{\kappa_1\kappa_2}^{\pvdots, j} \bigl[ g, h \bigr](s,x) \\
        &:= \widehat{\calF}^{-1}\biggl[ (2\pi)^2 \int_{\bbR^2} \int_{\bbR^2} (2i\jxi)^{-1} \frac{e^{is \Psi_{\kappa_1\kappa_2}(\xi,\xi_1,\xi_2)}}{i\Psi_{\kappa_1\kappa_2}(\xi,\xi_1,\xi_2)}  \widehat{g}^{\kappa_1}(s,\xi_1) \widehat{h}^{\kappa_2}(s,\xi_2)  \\
        &\qquad \qquad \times  \varphi_{\leq -M}\bigl( R(\xi_1,\xi_2) (\xi-\xi_1-\xi_2) \bigr) (-i) \varphi_{\leq 1}(\xi-\xi_1-\xi_2) \, \pvdots \, \frac{(\xi-\xi_1-\xi_2)_j}{|\xi-\xi_1-\xi_2|^3} \, \ud \xi_1 \, \ud \xi_2 \biggr](x),
    \end{aligned}
\end{equation}
where the cut-offs are defined as in Subsection~\ref{subsec:contributions_type1} below and where the factor $(-i)$ in front of the cut-off $\varphi_{\leq 1}$ is the constant in front of the leading order term in the expression \eqref{equ:FTvortex_re_im} for the flat Fourier transform of the components of the degree-one vortex $\underline{\Phi}$ established in Lemma~\ref{lem:FTvortex_re_im}. 
The bilinear expressions $\calB_{\kappa_1\kappa_2}^{\pvdots, j} \bigl[ g, h \bigr](s,x)$ arise from the singular parts of the quadratic expressions of type (1) discussed in Subsection~\ref{subsec:contributions_type1}.
For instance, the nonlinearity $\underline{\Phi}_\real u_{1,\real}^2$ comes with the singular part
\begin{equation}
    \begin{aligned}
        \calQ^{(1)}_{1,1,\mathrm{sing}}(t) &:= \sum_{\kappa_1, \kappa_2 \in \{\pm\}} \widehat{\calF}^{-1}\biggl[ (2\pi)^2 \int_t^T \int_{\bbR^2} \int_{\bbR^2} (2i\jxi)^{-1} e^{is\Psi_{\kappa_1 \kappa_2}(\xi,\xi_1,\xi_2)} \widehat{f}^{\kappa_1}_{1,\real}(s,\xi_1) \widehat{f}^{\kappa_2}_{1,\real}(s,\xi_2) \\
        &\qquad \times \varphi_{\leq -M}\bigl( R(\xi_1,\xi_2) (\xi-\xi_1-\xi_2) \bigr) (-i) \varphi_{\leq 1}(\xi-\xi_1-\xi_2) \, \pvdots \, \frac{(\xi-\xi_1-\xi_2)_1}{|\xi-\xi_1-\xi_2|^3} \, \ud \xi_1 \, \ud \xi_2 \, \ud s \biggr].
    \end{aligned}
\end{equation}
Integrating by parts in time we find
\begin{equation}
    \begin{aligned}
        &\calQ^{(1)}_{1,1,\mathrm{sing}}(t) = \sum_{\kappa_1, \kappa_2 \in\{\pm\}} \biggl( \calB_{\kappa_1\kappa_2}^{\pvdots, 1} \bigl[ f_{1,\real}, f_{1,\real} \bigr](T) - \calB_{\kappa_1\kappa_2}^{\pvdots, 1} \bigl[ f_{1,\real}, f_{1,\real} \bigr](t) \biggr) - \int_t^T \ldots \, \ud s.
    \end{aligned}
\end{equation}

For the contributions of nonlinear terms of type (2) discussed in Subsection~\ref{subsec:contributions_type2} we use the bilinear operators
\begin{equation}
    \begin{aligned}
        &\calB_{\kappa_1\kappa_2}^{\delta_0, j} \bigl[ g, h \bigr](s,x) := \widehat{\calF}^{-1}\biggl[ 2\pi \int_{\bbR^2} (2i\jxi)^{-1} (\xi_1)_j \frac{e^{is \Psi_{\kappa_1\kappa_2}(\xi,\xi-\xi_1,\xi_1)}}{\Psi_{\kappa_1\kappa_2}(\xi,\xi-\xi_1,\xi_1)}  \widehat{g}^{\kappa_1}(s,\xi-\xi_1) \widehat{h}^{\kappa_2}(s,\xi_1) \, \ud \xi_1 \biggr](x).
    \end{aligned}
\end{equation}
For example, for a non-localized quadratic nonlinearity $u_{4,\real} \partial_1 u_{3,\real}$ of type (2) we have 
\begin{equation}
    \begin{aligned}
        &\int_t^T e^{-is\jD} (2i\jD)^{-1} \bigl( u_{4,\real}(s) \partial_1 u_{3,\real}(s) \bigr) \, \ud s \\
        &\quad = \sum_{\kappa_1, \kappa_2 \in\{\pm\}} \biggl( \calB_{\kappa_1\kappa_2}^{\delta_0, 1} \bigl[ f_{4,\real}, f_{3,\real} \bigr](T) - \calB_{\kappa_1\kappa_2}^{\delta_0, 1} \bigl[ f_{4,\real}, f_{3,\real} \bigr](t) \biggr) - \int_t^T \ldots \, \ud s.
    \end{aligned}
\end{equation}

For the contributions of nonlinear terms of type (5) discussed in Subsection~\ref{subsec:contributions_type5} we define the bilinear operators
\begin{equation}
    \begin{aligned}
        &\calB_{\kappa_1\kappa_2}^{\delta_0, j, \eta_{0}} \bigl[ g, h \bigr](s,x) \\
        &\quad := \widehat{\calF}^{-1}\biggl[ 2\pi \int_{\bbR^2} (2i\jxi)^{-1} (\xi)_j \jxi^{-2} (\xi \cdot \xi_1) \frac{e^{is \Psi_{\kappa_1\kappa_2}(\xi,\xi-\xi_1,\xi_1)}}{\Psi_{\kappa_1\kappa_2}(\xi,\xi-\xi_1,\xi_1)}  \widehat{g}^{\kappa_1}(s,\xi-\xi_1) \widehat{h}^{\kappa_2}(s,\xi_1) \, \ud \xi_1 \biggr](x).
    \end{aligned}
\end{equation}
Only the leading flat quadratic part of nonlinear terms of type (5) of the form $-\partial_j \partial_t\eta_{0,c}$ contribute to these boundary operators. All terms involving the sharp temporal component, or cubic and higher-order expressions are treated
as perturbative remainders in Subsection~\ref{subsec:contributions_type5}. 
In view of \eqref{equ:pteta0c_flat_equation_rewritten}, \eqref{equ:pteta0c_flat_equation_leading_order_term_rewritten}, their flat components $-\partial_j \pt \eta_{0,c}^\flat$, $j \in \{1,2\}$, contribute the non-localized quadratic nonlinearities
\begin{equation}
    \begin{aligned}
        -\partial_j (-\Delta+1)^{-1} \nabla \cdot \bigl( u_{1,\real} \nabla u_{3,\real} - u_{3,\real} \nabla u_{1,\real} + u_{1,\imag} \nabla u_{3,\imag} - u_{3,\imag} \nabla u_{1,\imag} \bigr).
    \end{aligned}
\end{equation}
For instance, we have
\begin{equation}
    \begin{aligned}
        &\int_t^T e^{-is\jD} (2i\jD)^{-1} \Bigl( -\partial_j (-\Delta+1)^{-1} \nabla \cdot \bigl( u_{1,\real} \nabla u_{3,\real} \bigr) \Bigr) \, \ud s \\
        &\quad = \sum_{\kappa_1, \kappa_2 \in\{\pm\}} \biggl( \calB_{\kappa_1\kappa_2}^{\delta_0, j, \eta_0} \bigl[ f_{1,\real}, f_{3,\real} \bigr](T) - \calB_{\kappa_1\kappa_2}^{\delta_0, j, \eta_0} \bigl[ f_{1,\real}, f_{3,\real} \bigr](t) \biggr) - \int_t^T \ldots \, \ud s.
    \end{aligned}
\end{equation}

Finally, for the contributions of nonlinear terms of type (6) discussed in Subsection~\ref{subsec:contributions_type6} we introduce the bilinear operators 
\begin{equation}
    \begin{aligned}
        &\calB_{\kappa_1\kappa_2}^{\pvdots, j, \eta_{0}} \bigl[ g, h \bigr](s,x) \\
            &\quad := \widehat{\calF}^{-1}\biggl[ (2\pi)^2 \int_{\bbR^2} \int_{\bbR^2} (2i\jxi)^{-1} \jap{\xi_1+\xi_2}^{-2} \bigl( -(\xi_1+\xi_2)\cdot \xi_2 \bigr) \frac{e^{is \Psi_{\kappa_1\kappa_2}(\xi,\xi_1,\xi_2)}}{i\Psi_{\kappa_1\kappa_2}(\xi,\xi_1,\xi_2)} \\
            &\qquad \qquad \qquad \times  \widehat{g}^{\kappa_1}(s,\xi_1) \widehat{h}^{\kappa_2}(s,\xi_2) \varphi_{\leq -M}\bigl( R(\xi_1,\xi_2) (\xi-\xi_1-\xi_2) \bigr) (-i) \varphi_{\leq 1}(\xi-\xi_1-\xi_2) \\
            &\qquad \qquad \qquad \qquad \qquad \qquad \qquad \qquad \qquad \qquad \qquad \qquad \times \pvdots \, \frac{(\xi-\xi_1-\xi_2)_j}{|\xi-\xi_1-\xi_2|^3} \, \ud \xi_1 \, \ud \xi_2 \biggr](x).
    \end{aligned}
\end{equation}
These bilinear expressions originate precisely from the non-localized quadratic expressions contained in the singular parts of the terms $\underline{\Phi}_\real \pt \eta_{0,c}$, respectively $\underline{\Phi}_\imag \pt \eta_{0,c}$, by the same reasoning as above.

We now lay out the precise expressions for $\bmB_\real[\bmf, \bmf](t)$ and $\bmB_\imag[\bmf, \bmf](t)$.
We explain this in the context of the component $B_{\real, 1}[\bmf, \bmf](t)$, which collects all quadratic boundary terms coming from the nonlinearity $\Re\bigl( e^{i\theta} \calN_{\alpha,c}\bigr)$. In view of the exact expression for $\Re\bigl( e^{i\theta} \calN_{\alpha,c}\bigr)$ in Subsection~\ref{subsec:structure_nonlinearities}, the corresponding boundary terms of type $\calB_{\kappa_1\kappa_2}^{\delta_0, j}$ stem from
\begin{equation}
    \begin{aligned}
        - 2 u_{4,\real} \partial_1 u_{3,\real} - 2 u_{4,\imag} \partial_2 u_{3,\real},
    \end{aligned}
\end{equation}
while the corresponding boundary terms of type $\calB_{\kappa_1\kappa_2}^{\pvdots, j}$ stem from 
\begin{equation}
    \begin{aligned}
        -\frac32 \underline{\Phi}_\real \bigl( u_{1,\real}^2 + u_{1,\imag}^2 \bigr) + \frac12 \underline{\Phi}_\real \bigl( u_{3,\real}^2 + u_{3,\imag}^2 \bigr) - \underline{\Phi}_\real \bigl( u_{2,\real}^2 + u_{2,\imag}^2 + u_{4,\real}^2 + u_{4,\imag}^2 \bigr).
    \end{aligned}
\end{equation}
Thus, we obtain 
\begin{equation} \label{equ:definition_Breal1}
    \begin{aligned}
        B_{\real, 1}[\bmf, \bmf](t) &:= \sum_{\kappa_1, \kappa_2 \in\{\pm\}} \Bigl( -2 \calB_{\kappa_1\kappa_2}^{\delta_0, 1} \bigl[ f_{4,\real}, f_{3,\real} \bigr](t) -2 \calB_{\kappa_1\kappa_2}^{\delta_0, 2} \bigl[ f_{4,\imag}, f_{3,\real} \bigr](t) \Bigr) \\ 
        &\qquad + \sum_{\kappa_1, \kappa_2 \in\{\pm\}} \Bigl( - {\textstyle \frac32} \calB_{\kappa_1\kappa_2}^{\pvdots, 1} \bigl[ f_{1,\real}, f_{1,\real} \bigr](t) - {\textstyle \frac32} \calB_{\kappa_1\kappa_2}^{\pvdots, 1} \bigl[ f_{1,\imag}, f_{1,\imag} \bigr](t) \\
        &\qquad \qquad \qquad \qquad + {\textstyle \frac12} \calB_{\kappa_1\kappa_2}^{\pvdots, 1} \bigl[ f_{3,\real}, f_{3,\real} \bigr](t) + {\textstyle \frac12} \calB_{\kappa_1\kappa_2}^{\pvdots, 1} \bigl[ f_{3,\imag}, f_{3,\imag} \bigr](t) \\
        &\qquad \qquad \qquad \qquad - \calB_{\kappa_1\kappa_2}^{\pvdots, 1} \bigl[ f_{2,\real}, f_{2,\real} \bigr](t) - \calB_{\kappa_1\kappa_2}^{\pvdots, 1} \bigl[ f_{2,\imag}, f_{2,\imag} \bigr](t) \\
        &\qquad \qquad \qquad \qquad - \calB_{\kappa_1\kappa_2}^{\pvdots, 1} \bigl[ f_{4,\real}, f_{4,\real} \bigr](t) - \calB_{\kappa_1\kappa_2}^{\pvdots, 1} \bigl[ f_{4,\imag}, f_{4,\imag} \bigr](t) \Bigr),
    \end{aligned}
\end{equation}
and by analogous reasoning,
\begin{equation} \label{equ:definition_Bimag1}
    \begin{aligned}
        B_{\imag, 1}[\bmf, \bmf](t) 
        &:= \sum_{\kappa_1, \kappa_2 \in\{\pm\}} \Bigl( 
        -2 \calB_{\kappa_1\kappa_2}^{\delta_0, 1} \bigl[ f_{4,\real}, f_{3,\imag} \bigr](t) 
        -2 \calB_{\kappa_1\kappa_2}^{\delta_0, 2} \bigl[ f_{4,\imag}, f_{3,\imag} \bigr](t) \Bigr) \\ 
        &\qquad + \sum_{\kappa_1, \kappa_2 \in\{\pm\}} \Bigl( 
        - {\textstyle \frac32} \calB_{\kappa_1\kappa_2}^{\pvdots, 2} \bigl[ f_{1,\real}, f_{1,\real} \bigr](t) 
        - {\textstyle \frac32} \calB_{\kappa_1\kappa_2}^{\pvdots, 2} \bigl[ f_{1,\imag}, f_{1,\imag} \bigr](t) \\
        &\qquad \qquad \qquad \qquad 
        + {\textstyle \frac12} \calB_{\kappa_1\kappa_2}^{\pvdots, 2} \bigl[ f_{3,\real}, f_{3,\real} \bigr](t) 
        + {\textstyle \frac12} \calB_{\kappa_1\kappa_2}^{\pvdots, 2} \bigl[ f_{3,\imag}, f_{3,\imag} \bigr](t) \\
        &\qquad \qquad \qquad \qquad 
        - \calB_{\kappa_1\kappa_2}^{\pvdots, 2} \bigl[ f_{2,\real}, f_{2,\real} \bigr](t) 
        - \calB_{\kappa_1\kappa_2}^{\pvdots, 2} \bigl[ f_{2,\imag}, f_{2,\imag} \bigr](t) \\
        &\qquad \qquad \qquad \qquad 
        - \calB_{\kappa_1\kappa_2}^{\pvdots, 2} \bigl[ f_{4,\real}, f_{4,\real} \bigr](t) 
        - \calB_{\kappa_1\kappa_2}^{\pvdots, 2} \bigl[ f_{4,\imag}, f_{4,\imag} \bigr](t) \Bigr). 
    \end{aligned}
\end{equation}
Moreover, 
\begin{equation} \label{equ:definition_Brealimag2}
    \begin{aligned}
        B_{\real, 2}[\bmf, \bmf](t) 
        &:= \sum_{\kappa_1, \kappa_2 \in\{\pm\}} \Bigl(
        -2 \calB_{\kappa_1\kappa_2}^{\pvdots, 1} \bigl[ f_{2,\real}, f_{1,\real} \bigr](t)
        -2 \calB_{\kappa_1\kappa_2}^{\pvdots, 1} \bigl[ f_{2,\imag}, f_{1,\imag} \bigr](t)
        \Bigr), \\
        B_{\imag, 2}[\bmf, \bmf](t) 
        &:= \sum_{\kappa_1, \kappa_2 \in\{\pm\}} \Bigl(
        -2 \calB_{\kappa_1\kappa_2}^{\pvdots, 2} \bigl[ f_{2,\real}, f_{1,\real} \bigr](t)
        -2 \calB_{\kappa_1\kappa_2}^{\pvdots, 2} \bigl[ f_{2,\imag}, f_{1,\imag} \bigr](t)
        \Bigr), 
    \end{aligned}
\end{equation}
as well as 
\begin{equation} \label{equ:definition_Brealimag3}
    \begin{aligned}
        B_{\real, 3}[\bmf, \bmf](t)
        &:= \sum_{\kappa_1, \kappa_2 \in\{\pm\}} \Bigl(
        2 \calB_{\kappa_1\kappa_2}^{\delta_0, 1} \bigl[ f_{4,\real}, f_{1,\real} \bigr](t)
        +2 \calB_{\kappa_1\kappa_2}^{\delta_0, 2} \bigl[ f_{4,\imag}, f_{1,\real} \bigr](t) \Bigr) \\
        &\qquad
        + \sum_{\kappa_1, \kappa_2 \in\{\pm\}} \Bigl( -2 \calB_{\kappa_1\kappa_2}^{\pvdots, 1} \bigl[ f_{1,\real}, f_{3,\real} \bigr](t)
        -2 \calB_{\kappa_1\kappa_2}^{\pvdots, 1} \bigl[ f_{1,\imag}, f_{3,\imag} \bigr](t) \Bigr) \\
        &\qquad
        + \sum_{\kappa_1, \kappa_2 \in\{\pm\}} \Bigl( \calB_{\kappa_1\kappa_2}^{\pvdots, 1, \eta_0} \bigl[ f_{1,\real}, f_{3,\real} \bigr](t)
        - \calB_{\kappa_1\kappa_2}^{\pvdots, 1, \eta_0} \bigl[ f_{3,\real}, f_{1,\real} \bigr](t) \\
        &\qquad \qquad \qquad \qquad 
        + \calB_{\kappa_1\kappa_2}^{\pvdots, 1, \eta_0} \bigl[ f_{1,\imag}, f_{3,\imag} \bigr](t)
        - \calB_{\kappa_1\kappa_2}^{\pvdots, 1, \eta_0} \bigl[ f_{3,\imag}, f_{1,\imag} \bigr](t)
        \Bigr), \\
        B_{\imag, 3}[\bmf, \bmf](t)
        &:= \sum_{\kappa_1, \kappa_2 \in\{\pm\}} \Bigl(
        2 \calB_{\kappa_1\kappa_2}^{\delta_0, 1} \bigl[ f_{4,\real}, f_{1,\imag} \bigr](t)
        +2 \calB_{\kappa_1\kappa_2}^{\delta_0, 2} \bigl[ f_{4,\imag}, f_{1,\imag} \bigr](t) \Bigr) \\
        &\qquad
        + \sum_{\kappa_1, \kappa_2 \in\{\pm\}} \Bigl( -2 \calB_{\kappa_1\kappa_2}^{\pvdots, 2} \bigl[ f_{1,\real}, f_{3,\real} \bigr](t)
        -2 \calB_{\kappa_1\kappa_2}^{\pvdots, 2} \bigl[ f_{1,\imag}, f_{3,\imag} \bigr](t) \Bigr) \\
        &\qquad
        + \sum_{\kappa_1, \kappa_2 \in\{\pm\}} \Bigl( \calB_{\kappa_1\kappa_2}^{\pvdots, 2, \eta_0} \bigl[ f_{1,\real}, f_{3,\real} \bigr](t)
        - \calB_{\kappa_1\kappa_2}^{\pvdots, 2, \eta_0} \bigl[ f_{3,\real}, f_{1,\real} \bigr](t) \\
        &\qquad \qquad \qquad \qquad 
        + \calB_{\kappa_1\kappa_2}^{\pvdots, 2, \eta_0} \bigl[ f_{1,\imag}, f_{3,\imag} \bigr](t)
        - \calB_{\kappa_1\kappa_2}^{\pvdots, 2, \eta_0} \bigl[ f_{3,\imag}, f_{1,\imag} \bigr](t) 
        \Bigr).
    \end{aligned}
\end{equation}
Finally,
\begin{equation} \label{equ:definition_Brealimag4}
    \begin{aligned}
        B_{\real, 4}[\bmf, \bmf](t) 
        &:= \sum_{\kappa_1, \kappa_2 \in\{\pm\}} \Bigl(
        - \calB_{\kappa_1\kappa_2}^{\delta_0, 1} \bigl[ f_{1,\real}, f_{3,\real} \bigr](t)
        - \calB_{\kappa_1\kappa_2}^{\delta_0, 2} \bigl[ f_{1,\imag}, f_{3,\real} \bigr](t)   \\
        &\qquad \qquad \qquad \qquad 
        + \calB_{\kappa_1\kappa_2}^{\delta_0, 1} \bigl[ f_{3,\real}, f_{1,\real} \bigr](t)
        + \calB_{\kappa_1\kappa_2}^{\delta_0, 2} \bigl[ f_{3,\imag}, f_{1,\real} \bigr](t) \Bigr) \\
        &\qquad
        + \sum_{\kappa_1, \kappa_2 \in\{\pm\}} \Bigl( -2 \calB_{\kappa_1\kappa_2}^{\pvdots, 1} \bigl[ f_{1,\real}, f_{4,\real} \bigr](t)
        -2 \calB_{\kappa_1\kappa_2}^{\pvdots, 1} \bigl[ f_{1,\imag}, f_{4,\imag} \bigr](t) \Bigr) \\
        &\qquad
        + \sum_{\kappa_1, \kappa_2 \in\{\pm\}} \Bigl( \calB_{\kappa_1\kappa_2}^{\delta_0,1,\eta_0} \bigl[ f_{1,\real}, f_{3,\real} \bigr](t)
        - \calB_{\kappa_1\kappa_2}^{\delta_0,1,\eta_0} \bigl[ f_{3,\real}, f_{1,\real} \bigr](t) \\
        &\qquad \qquad \qquad \qquad 
        + \calB_{\kappa_1\kappa_2}^{\delta_0,1,\eta_0} \bigl[ f_{1,\imag}, f_{3,\imag} \bigr](t)
        - \calB_{\kappa_1\kappa_2}^{\delta_0,1,\eta_0} \bigl[ f_{3,\imag}, f_{1,\imag} \bigr](t) \Bigr), \\
        B_{\imag, 4}[\bmf, \bmf](t) 
        &:= \sum_{\kappa_1, \kappa_2 \in\{\pm\}} \Bigl(
        - \calB_{\kappa_1\kappa_2}^{\delta_0, 1} \bigl[ f_{1,\real}, f_{3,\imag} \bigr](t)
        - \calB_{\kappa_1\kappa_2}^{\delta_0, 2} \bigl[ f_{1,\imag}, f_{3,\imag} \bigr](t) \\
        &\qquad \qquad \qquad \qquad 
        + \calB_{\kappa_1\kappa_2}^{\delta_0, 1} \bigl[ f_{3,\real}, f_{1,\imag} \bigr](t)
        + \calB_{\kappa_1\kappa_2}^{\delta_0, 2} \bigl[ f_{3,\imag}, f_{1,\imag} \bigr](t) \Bigr) \\
        &\qquad
        + \sum_{\kappa_1, \kappa_2 \in\{\pm\}} \Bigl( -2 \calB_{\kappa_1\kappa_2}^{\pvdots, 2} \bigl[ f_{1,\real}, f_{4,\real} \bigr](t)
        -2 \calB_{\kappa_1\kappa_2}^{\pvdots, 2} \bigl[ f_{1,\imag}, f_{4,\imag} \bigr](t) \Bigr) \\
        &\qquad
        + \sum_{\kappa_1, \kappa_2 \in\{\pm\}} \Bigl( \calB_{\kappa_1\kappa_2}^{\delta_0,2,\eta_0} \bigl[ f_{1,\real}, f_{3,\real} \bigr](t)
        - \calB_{\kappa_1\kappa_2}^{\delta_0,2,\eta_0} \bigl[ f_{3,\real}, f_{1,\real} \bigr](t) \\
        &\qquad \qquad \qquad \qquad 
        + \calB_{\kappa_1\kappa_2}^{\delta_0,2,\eta_0} \bigl[ f_{1,\imag}, f_{3,\imag} \bigr](t)
        - \calB_{\kappa_1\kappa_2}^{\delta_0,2,\eta_0} \bigl[ f_{3,\imag}, f_{1,\imag} \bigr](t)
        \Bigr).
    \end{aligned}
\end{equation}

\subsection{Preparations} \label{subsec:preparations}

In this subsection we collect several multilinear estimates and bounds on a quadratic phase function that will be used frequently in the space-time resonance analysis in the remainder of this section. Moreover, we determine the flat Fourier transform of the degree-one vortex and derive auxiliary estimates for the time derivative of the flat profiles $\bmf_\real(t)$ and $\bmf_\imag(t)$. 

\subsubsection{Multilinear estimates}

We introduce two types of bilinear operators
\begin{equation}
    \begin{aligned}
        \calB_{\delta_0}\bigl[f, g\bigr](x) &:= \widehat{\calF}^{-1}_{\xi\mapsto x}\biggl[ \int_{\bbR^2} \fraka(\xi_1, \xi-\xi_1) \hatf(\xi_1) \hatg(\xi-\xi_1) \, \ud \xi_1 \biggr](x) \\ 
        &= \frac{1}{2\pi} \int_{\bbR^2} \int_{\bbR^2} e^{i x \cdot (\xi_1+\xi_2)} \fraka(\xi_1,\xi_2) \hatf(\xi_1) \hatg(\xi_2) \, \ud \xi_1 \, \ud \xi_2, \\
        \calB_{\pvdots}\bigl[f, g\bigr](x) &:= \widehat{\calF}^{-1}_{\xi\mapsto x}\biggl[ \int_{\bbR^2} \int_{\bbR^2}  \frakb(\xi_1, \xi-\xi_1-\sigma, \sigma) \hatf(\xi_1) \hatg(\xi-\xi_1-\sigma) \, \pvdots \frac{(\sigma)_j}{|\sigma|^3} \, \ud \sigma \, \ud \xi_1 \biggr](x) \\
        &= \frac{1}{2\pi} \int_{\bbR^2} \int_{\bbR^2} \int_{\bbR^2} e^{i x \cdot (\xi_1+\xi_2+\sigma)} \frakb(\xi_1,\xi_2,\sigma) \hatf(\xi_1) \hatg(\xi_2) \, \pvdots \frac{(\sigma)_j}{|\sigma|^3} \, \ud \sigma \, \ud \xi_1 \, \ud \xi_2,
    \end{aligned}
\end{equation}
where $j \in \{1, 2\}$ and $\sigma = \bigl( (\sigma)_1, (\sigma)_2 \bigr) \in \bbR^2$.
The following bilinear H\"older-type bounds will be used repeatedly in the nonlinear analysis. 

\begin{lemma} \label{lem:bilinear_estimates}
    If $1 \leq p, p_1, p_2 \leq \infty$ satisfy $\frac{1}{p} = \frac{1}{p_1} + \frac{1}{p_2}$, then
    \begin{align}
        \bigl\| \calB_{\delta_0}\bigl[f, g\bigr] \bigr\|_{L^p_x(\bbR^2)} &\lesssim \bigl\| \widehat{\calF}^{-1}[\fraka] \bigr\|_{L^1((\bbR^2)^2)} \|f\|_{L^{p_1}_x(\bbR^2)} \|g\|_{L^{p_2}_x(\bbR^2)}, \label{equ:bilinear_Hoelder_type_delta} \\
        \bigl\| \calB_{\pvdots}\bigl[f, g\bigr] \bigr\|_{L^p_x(\bbR^2)} &\lesssim \bigl\| \widehat{\calF}^{-1}[\frakb] \bigr\|_{L^1((\bbR^2)^3)} \|f\|_{L^{p_1}_x(\bbR^2)} \|g\|_{L^{p_2}_x(\bbR^2)}. \label{equ:bilinear_Hoelder_type_pv}
    \end{align}
\end{lemma}
\begin{proof}
    The first bilinear estimate \eqref{equ:bilinear_Hoelder_type_delta} follows from the representation
    \begin{equation}
        \begin{aligned}
            \calB_{\delta_0}\bigl[f, g\bigr](x) = \frac{1}{2\pi} \int_{\bbR^2} \int_{\bbR^2} \widehat{\calF}^{-1}[\fraka](y_1,y_2) f(x-y_1) g(x-y_2) \, \ud y_1 \, \ud y_2
        \end{aligned}
    \end{equation}
    by Minkowski's integral inequality, H\"older's inequality and the translation invariance of Lebesgue norms.
    For the proof of \eqref{equ:bilinear_Hoelder_type_pv} we recall the two-dimensional Riesz kernel identity
    \begin{equation}
        \calR_j(x) := \widehat{\calF}^{-1}\biggl[ \pvdots \, \frac{\xi_j}{|\xi|^3} \biggr](x) = i \frac{x_j}{|x|}, \quad j \in \{1,2\},
    \end{equation}
    in the sense of tempered distributions, where the normalization \eqref{eq:FT} is used.
    In particular, we have $\|\calR_j\|_{L^\infty_x(\bbR^2)} \lesssim 1$. Then \eqref{equ:bilinear_Hoelder_type_pv} follows from the representation
    \begin{equation}
        \begin{aligned}
            \calB_{\pvdots}\bigl[f, g\bigr](x) = \frac{1}{2\pi} \int_{\bbR^2} \int_{\bbR^2} \int_{\bbR^2} \widehat{\calF}^{-1}[\frakb](y_1,y_2,y_3) f(x-y_1) g(x-y_2) \calR_j(x-y_3) \, \ud y_1 \, \ud y_2 \, \ud y_3,
        \end{aligned}
    \end{equation}
    as above.
\end{proof}

Moreover, on many occasions we will face the following two types of trilinear operators
\begin{equation}
    \begin{aligned}
        &\calC_{\delta_0}\bigl[f, g, h\bigr](x) \\
        &\quad := \widehat{\calF}^{-1}_{\xi\mapsto x}\biggl[ \int_{\bbR^2} \int_{\bbR^2} \frakc(\xi_1, \xi_2, \xi-\xi_1-\xi_2) \hatf(\xi_1) \hatg(\xi_2) \hath(\xi-\xi_1-\xi_2) \, \ud \xi_1 \, \ud \xi_2 \biggr](x) \\ 
        &\quad = \frac{1}{2\pi} \int_{\bbR^2} \int_{\bbR^2} \int_{\bbR^2} e^{i x \cdot (\xi_1+\xi_2+\xi_3)} \frakc(\xi_1,\xi_2,\xi_3) \hatf(\xi_1) \hatg(\xi_2) \hath(\xi_3) \, \ud \xi_1 \, \ud \xi_2 \, \ud \xi_3, \\
        &\calC_{\pvdots}\bigl[f, g, h\bigr](x) \\
        &\quad := \widehat{\calF}^{-1}_{\xi\mapsto x}\biggl[ \int_{\bbR^2} \int_{\bbR^2} \int_{\bbR^2} \frakd(\xi_1, \xi_2, \xi-\xi_1-\xi_2-\sigma, \sigma) \hatf(\xi_1) \hatg(\xi_2) \hath(\xi-\xi_1-\xi_2-\sigma) \\
        &\qquad \qquad \qquad \qquad \qquad \qquad \qquad \qquad \qquad \qquad \qquad \qquad \qquad \qquad \times \, \pvdots \frac{(\sigma)_j}{|\sigma|^3} \, \ud \sigma \, \ud \xi_1 \, \ud \xi_2 \biggr](x) \\
        &\quad = \frac{1}{2\pi} \int_{\bbR^2} \int_{\bbR^2} \int_{\bbR^2} \int_{\bbR^2} e^{i x \cdot (\xi_1+\xi_2+\xi_3+\sigma)} \frakd(\xi_1,\xi_2,\xi_3,\sigma) \hatf(\xi_1) \hatg(\xi_2) \hath(\xi_3) \, \pvdots \frac{(\sigma)_j}{|\sigma|^3} \, \ud \sigma \, \ud \xi_1 \, \ud \xi_2 \, \ud \xi_3,
    \end{aligned}
\end{equation}
where $j \in \{1, 2\}$ and $\sigma = \bigl( (\sigma)_1, (\sigma)_2 \bigr) \in \bbR^2$.
Additionally, we will make use of the following two types of quadrilinear operators 
\begin{equation}
    \begin{aligned}
        &\calQ_{\delta_0}\bigl[f, g, h, q\bigr](x) \\
        &\quad := \widehat{\calF}^{-1}_{\xi\mapsto x}\biggl[ \int_{\bbR^2} \int_{\bbR^2} \int_{\bbR^2} \frakm(\xi_1, \xi_2, \xi_3, \xi-\xi_1-\xi_2-\xi_3) \\
        &\qquad \qquad \qquad \qquad \times \hatf(\xi_1) \hatg(\xi_2) \hath(\xi_3) \hatq(\xi-\xi_1-\xi_2-\xi_3) \, \ud \xi_1 \, \ud \xi_2 \, \ud \xi_3 \biggr](x) \\ 
        &\quad = \frac{1}{2\pi} \int_{\bbR^2} \int_{\bbR^2} \int_{\bbR^2} \int_{\bbR^2} e^{i x \cdot (\xi_1+\xi_2+\xi_3+\xi_4)} \frakm(\xi_1,\xi_2,\xi_3,\xi_4) \\
        &\qquad \qquad \qquad \qquad \times \hatf(\xi_1) \hatg(\xi_2) \hath(\xi_3) \hatq(\xi_4) \, \ud \xi_1 \, \ud \xi_2 \, \ud \xi_3 \, \ud \xi_4, \\
        &\calQ_{\pvdots}\bigl[f, g, h, q\bigr](x) \\
        &\quad := \widehat{\calF}^{-1}_{\xi\mapsto x}\biggl[ \int_{\bbR^2} \int_{\bbR^2} \int_{\bbR^2} \int_{\bbR^2} \frakn(\xi_1, \xi_2, \xi_3, \xi-\xi_1-\xi_2-\xi_3-\sigma, \sigma) \\
        &\qquad \qquad \qquad \qquad \times \hatf(\xi_1) \hatg(\xi_2) \hath(\xi_3) \hatq(\xi-\xi_1-\xi_2-\xi_3-\sigma) \, \pvdots \frac{(\sigma)_j}{|\sigma|^3} \, \ud \sigma \, \ud \xi_1 \, \ud \xi_2 \, \ud \xi_3 \biggr](x) \\
        &\quad = \frac{1}{2\pi} \int_{\bbR^2} \int_{\bbR^2} \int_{\bbR^2} \int_{\bbR^2} \int_{\bbR^2} e^{i x \cdot (\xi_1+\xi_2+\xi_3+\xi_4+\sigma)} \frakn(\xi_1,\xi_2,\xi_3,\xi_4,\sigma) \\
        &\qquad \qquad \qquad \qquad \times \hatf(\xi_1) \hatg(\xi_2) \hath(\xi_3) \hatq(\xi_4) \, \pvdots \frac{(\sigma)_j}{|\sigma|^3} \, \ud \sigma \, \ud \xi_1 \, \ud \xi_2 \, \ud \xi_3 \, \ud \xi_4,
    \end{aligned}
\end{equation}
where $j \in \{1, 2\}$ and $\sigma = \bigl( (\sigma)_1, (\sigma)_2 \bigr) \in \bbR^2$. 
Finally, we will also encounter the following quintilinear operator 
\begin{equation}
    \begin{aligned}
        &\calR_{\delta_0}\bigl[f, g, h, q, u \bigr](x) \\
        &\quad := \widehat{\calF}^{-1}_{\xi\mapsto x}\biggl[ \int_{\bbR^2} \int_{\bbR^2} \int_{\bbR^2} \int_{\bbR^2} \frakq(\xi_1, \xi_2, \xi_3, \xi_4, \xi-\xi_1-\xi_2-\xi_3-\xi_4) \\
        &\qquad \qquad \qquad \qquad \times \hatf(\xi_1) \hatg(\xi_2) \hath(\xi_3) \hatq(\xi_4) \hatu(\xi-\xi_1-\xi_2-\xi_3-\xi_4) \, \ud \xi_1 \, \ud \xi_2 \, \ud \xi_3 \, \ud \xi_4 \biggr](x) \\ 
        &\quad = \frac{1}{2\pi} \int_{\bbR^2} \int_{\bbR^2} \int_{\bbR^2} \int_{\bbR^2} \int_{\bbR^2} e^{i x \cdot (\xi_1+\xi_2+\xi_3+\xi_4+\xi_5)} \frakq(\xi_1,\xi_2,\xi_3,\xi_4,\xi_5) \\
        &\qquad \qquad \qquad \qquad \times \hatf(\xi_1) \hatg(\xi_2) \hath(\xi_3) \hatq(\xi_4) \hatu(\xi_5) \, \ud \xi_1 \, \ud \xi_2 \, \ud \xi_3 \, \ud \xi_4 \, \ud \xi_5.
    \end{aligned}
\end{equation}

In the next lemmas we record multilinear H\"older-type bounds for the preceding operators. Their proofs are completely analogous to the proof of Lemma~\ref{lem:bilinear_estimates} and are therefore omitted.

\begin{lemma} \label{lem:trilinear_estimates}
    If $1 \leq p, p_1, p_2, p_3 \leq \infty$ satisfy $\frac{1}{p} = \frac{1}{p_1} + \frac{1}{p_2} + \frac{1}{p_3}$, then
    \begin{equation}
        \bigl\| \calC_{\delta_0}\bigl[f, g, h\bigr] \bigr\|_{L^p_x(\bbR^2)} \lesssim \bigl\| \widehat{\calF}^{-1}[\frakc] \bigr\|_{L^1((\bbR^2)^3)} \|f\|_{L^{p_1}_x(\bbR^2)} \|g\|_{L^{p_2}_x(\bbR^2)} \|h\|_{L^{p_3}_x(\bbR^2)}
    \end{equation}
    and
    \begin{equation}
        \bigl\| \calC_{\pvdots}\bigl[f, g, h\bigr] \bigr\|_{L^p_x(\bbR^2)} \lesssim \bigl\| \widehat{\calF}^{-1}[\frakd] \bigr\|_{L^1((\bbR^2)^4)} \|f\|_{L^{p_1}_x(\bbR^2)} \|g\|_{L^{p_2}_x(\bbR^2)} \|h\|_{L^{p_3}_x(\bbR^2)}.
    \end{equation}
\end{lemma}

\begin{lemma} \label{lem:quadrilinear_estimates}
    If $1 \leq p, p_1, p_2, p_3, p_4 \leq \infty$ satisfy $\frac{1}{p} = \frac{1}{p_1} + \frac{1}{p_2} + \frac{1}{p_3} + \frac{1}{p_4}$, then
    \begin{equation}
        \bigl\| \calQ_{\delta_0}\bigl[f, g, h, q\bigr] \bigr\|_{L^p_x(\bbR^2)} \lesssim \bigl\| \widehat{\calF}^{-1}[\frakm] \bigr\|_{L^1((\bbR^2)^4)} \|f\|_{L^{p_1}_x(\bbR^2)} \|g\|_{L^{p_2}_x(\bbR^2)} \|h\|_{L^{p_3}_x(\bbR^2)} \|q\|_{L^{p_4}_x(\bbR^2)}
    \end{equation}
    and
    \begin{equation}
        \bigl\| \calQ_{\pvdots}\bigl[f, g, h, q\bigr] \bigr\|_{L^p_x(\bbR^2)} \lesssim \bigl\| \widehat{\calF}^{-1}[\frakn] \bigr\|_{L^1((\bbR^2)^5)} \|f\|_{L^{p_1}_x(\bbR^2)} \|g\|_{L^{p_2}_x(\bbR^2)} \|h\|_{L^{p_3}_x(\bbR^2)}\|q\|_{L^{p_4}_x(\bbR^2)}.
    \end{equation}
\end{lemma}

\begin{lemma} \label{lem:quintiilinear_estimates}
    If $1 \leq p, p_1, p_2, p_3, p_4, p_5 \leq \infty$ satisfy $\frac{1}{p} = \frac{1}{p_1} + \frac{1}{p_2} + \frac{1}{p_3} + \frac{1}{p_4} + \frac{1}{p_5}$, then
    \begin{equation}
        \bigl\| \calR_{\delta_0}\bigl[f, g, h, q, u\bigr] \bigr\|_{L^p_x(\bbR^2)} \lesssim \bigl\| \widehat{\calF}^{-1}[\frakq] \bigr\|_{L^1((\bbR^2)^5)} \|f\|_{L^{p_1}_x(\bbR^2)} \|g\|_{L^{p_2}_x(\bbR^2)} \|h\|_{L^{p_3}_x(\bbR^2)} \|q\|_{L^{p_4}_x(\bbR^2)} \|u\|_{L^{p_5}_x(\bbR^2)}.
    \end{equation}
\end{lemma}

\subsubsection{Bounds for the quadratic phase function}

The following quadratic phase function will come up in the space-time resonance analysis of the quadratic nonlinearities in the remainder of this section. For $\kappa_1, \kappa_2 \in \{\pm\}$ and $\xi, \xi_1, \xi_2 \in \bbR^2$ we set
\begin{equation} \label{equ:definition_Psi_quadratic}
    \Psi_{\kappa_1 \kappa_2}(\xi,\xi_1,\xi_2) := -\jxi + \kappa_1 \jap{\xi_1} + \kappa_2 \jap{\xi_2}.
\end{equation}
In the next lemma we establish lower bounds on the phase function \eqref{equ:definition_Psi_quadratic} that are needed for the normal form transformations of the quadratic nonlinearities. 

\begin{lemma}[Lower bound on phase function] \label{lem:lower_bound_quadratic_phase}
    For any $\kappa_1, \kappa_2 \in \{\pm\}$ and any $\xi_1, \xi_2 \in \bbR^2$ we have 
    \begin{equation} \label{equ:lower_bound_quadr_phase_delta}
        \bigl| \Psi_{\kappa_1 \kappa_2}(\xi_1+\xi_2,\xi_1,\xi_2) \bigr| \geq \frac12 \frac{1}{\min\bigl\{ \jxione, \jxitwo \bigr\}}.
    \end{equation}
    Moreover, we have for any $\kappa_1, \kappa_2 \in \{\pm\}$ and any $\xi_1, \xi_2 \in \bbR^2$ that if $\sigma \in \bbR^2$ satisfies 
    \begin{equation} \label{equ:lower_bound_quadr_phase_sigma_assumpt}
        |\sigma| \lesssim \frac{2^{-M}}{\min\bigl\{ \jxione, \jxitwo \bigr\}}    
    \end{equation}
    for some sufficiently large absolute constant $M \gg 1$, then 
    \begin{equation} \label{equ:lower_bound_quadr_phase_pv}
        \bigl| \Psi_{\kappa_1 \kappa_2}(\xi_1+\xi_2+\sigma,\xi_1,\xi_2) \bigr| \geq \frac14 \frac{1}{\min\bigl\{ \jxione, \jxitwo \bigr\}}.
    \end{equation}
\end{lemma}
\begin{proof}
    We begin with the proof of \eqref{equ:lower_bound_quadr_phase_delta}. First, we consider the case $(\kappa_1, \kappa_2) = (+,+)$.
    It is elementary to check that 
    \begin{equation} \label{equ:lower_bound_quadr_phase_proof1}
         \Psi_{++}(\xi_1+\xi_2, \xi_1,\xi_2) = -\jap{\xi_1+\xi_2} + \jxione + \jxitwo > 0.
    \end{equation}
    By inspection, we have $\Psi_{++}(\xi_1,\xi_1,0) = 1$ for all $\xi_1 \in \bbR^2$. 
    Fix some $\xi_2 \ne 0$, and assume that~\eqref{equ:lower_bound_quadr_phase_proof1} attains a minimum in the $\xi_1$-plane. The critical points in $\xi_1$ satisfy $\nabla_{\xi_1} \Psi_{++}(\xi_1+\xi_2,\xi_1,\xi_2)=0$, or equivalently,
    \begin{equation} \label{equ:lower_bound_quadr_phase_proof2}
        \frac{\xi_1}{\jxione} = \frac{\xi_1+\xi_2}{\jap{\xi_1+\xi_2}},
    \end{equation}
    which in particular implies that $2\xi_1 \cdot \xi_2 + |\xi_2|^2 = 0$. The solutions are given by the line $\xi_1 = -\frac12 \xi_2 + \lambda \xi_2^\perp$, $\lambda \in \bbR$ for a unit vector $\xi_2^\perp \perp \xi_2$. Plugging this expression into~\eqref{equ:lower_bound_quadr_phase_proof2}, we see that $\xi_2 = 0$, which is a contradiction.
    Hence, \eqref{equ:lower_bound_quadr_phase_proof1} does not attain a minimum in the $\xi_1$-plane, which means that the lower bound is found asymptotically as $|\xi_1| \to \infty$.  
    It is elementary to verify that
    \begin{equation}
        \liminf_{|\xi_1| \to \infty} \, \bigl( -\jap{\xi_1+\xi_2} + \jap{\xi_1} \bigr) \geq -|\xi_2|,
    \end{equation}
    whence 
    \begin{equation}
        \liminf_{|\xi_1| \to \infty} \, \bigl( -\jap{\xi_1+\xi_2} + \jxione + \jxitwo \bigr) \geq \jxitwo - |\xi_2| \geq \frac{1}{2\jxitwo}.
    \end{equation}
    By symmetry, we then conclude 
    \begin{equation}
        \Psi_{++}(\xi_1+\xi_2,\xi_1,\xi_2) \geq \frac12 \frac{1}{\min\bigl\{ \jxione, \jxitwo \bigr\}}.
    \end{equation}
    The cases $(\kappa_1, \kappa_2) = (+,-)$ and $(\kappa_1, \kappa_2) = (-,+)$ follow from the case $(\kappa_1, \kappa_2) = (+,+)$ by permutation of the variables, while in the case $(\kappa_1, \kappa_2) = (-,-)$ we have the trivial lower bound $|\Psi_{--}(\xi_1+\xi_2,\xi_1,\xi_2)| \geq 3$.
      
    Next, we turn to the proof of \eqref{equ:lower_bound_quadr_phase_pv}. 
    By direct computation, we find 
    \begin{equation}
        \begin{aligned}
            \Psi_{\kappa_1\kappa_2}(\xi_1+\xi_2+\sigma,\xi_1,\xi_2) &= -\jap{\xi_1+\xi_2+\sigma} + \kappa_1 \jap{\xi_1} + \kappa_2 \jap{\xi_2} \\
            &= -\jap{\xi_1+\xi_2} + \kappa_1 \jap{\xi_1} + \kappa_2 \jap{\xi_2} + \frac{|\xi_1+\xi_2|^2 - |\xi_1+\xi_2+\sigma|^2}{\jap{\xi_1+\xi_2} + \jap{\xi_1+\xi_2+\sigma}}.
        \end{aligned}
    \end{equation}
    Hence, by the previously established lower bound \eqref{equ:lower_bound_quadr_phase_delta} and by assumption \eqref{equ:lower_bound_quadr_phase_sigma_assumpt}, we obtain
    \begin{equation}
        \begin{aligned}
            \bigl| \Psi_{\kappa_1\kappa_2}(\xi_1+\xi_2+\sigma,\xi_1,\xi_2) \bigr| &\geq \bigl| -\jap{\xi_1+\xi_2} + \kappa_1 \jap{\xi_1} + \kappa_2 \jap{\xi_2} \bigr| - |\sigma| \frac{|\xi_1+\xi_2| + |\xi_1+\xi_2+\sigma|}{\jap{\xi_1+\xi_2} + \jap{\xi_1+\xi_2+\sigma}} \\
            &\geq \frac12 \frac{1}{\min\bigl\{ \jxione, \jxitwo \bigr\}} - \frac{C 2^{-M}}{\min\{\jxione, \jxitwo\}} \\
            &\geq \frac14 \frac{1}{\min\bigl\{ \jxione, \jxitwo \bigr\}},
        \end{aligned}
    \end{equation}
    as desired.
\end{proof}

Next, we obtain derivative bounds for the inverse of the phase function \eqref{equ:definition_Psi_quadratic}, which are key to establish H\"older-type bounds for the multilinear expressions resulting from the normal form transformations of the quadratic nonlinearities.

\begin{lemma}[Derivative bounds for the inverse of the phase function] \label{lem:derivative_bounds_quadratic_phase}
    For any $\kappa_1, \kappa_2 \in \{ \pm \}$, any $\xi_1, \xi_2 \in \bbR^2$, and any integers $a,b \geq 0$, we have 
    \begin{equation} \label{equ:derivative_bounds_inverse_phase_delta}
         \Bigl| \nabla_{\xi_1}^a \nabla_{\xi_2}^b \frac{1}{\Psi_{\kappa_1 \kappa_2}(\xi_1+\xi_2,\xi_1,\xi_2)} \Big| \lesssim \jxione^{-a} \jxitwo^{-b} \bigl( \min\{ \jxione, \jxitwo \} \bigr)^{1+2a+2b}.
    \end{equation}
    Moreover, we have for any $\kappa_1, \kappa_2 \in \{\pm\}$ and any $\xi_1, \xi_2 \in \bbR^2$ that if $\sigma \in \bbR^2$ satisfies
    \begin{equation} \label{equ:derivative_bounds_inverse_phase_sigma_assumption}
        |\sigma| \lesssim \frac{2^{-M}}{\min\bigl\{ \jxione, \jxitwo \bigr\}}    
    \end{equation}
    for some sufficiently large absolute constant $M \gg 1$, then for any integers $a, b, c \geq 0$,
    \begin{equation}  \label{equ:derivative_bounds_inverse_phase_pv}
        \Bigl| \nabla_{\xi_1}^a \nabla_{\xi_2}^b \nabla_\sigma^c \frac{1}{\Psi_{\kappa_1 \kappa_2}(\xi_1+\xi_2+\sigma,\xi_1,\xi_2)} \Big| \lesssim \jxione^{-a} \jxitwo^{-b} \bigl( \min\{ \jxione, \jxitwo \} \bigr)^{1+2a+2b+c}.
    \end{equation} 
\end{lemma}
\begin{proof}
The asserted bounds are immediate in the case $(\kappa_1, \kappa_2) = (-,-)$.
We now consider the most delicate case $(\kappa_1,\kappa_2)=(+,+)$. For any $\xi_1,\xi_2\in\mathbb R^2$ and
any integers $a,b\geq 0$ with $a+b>0$, we claim that
\begin{equation} \label{equ:derivative_bounds_phase_delta}
         \bigl| \nabla_{\xi_1}^a \nabla_{\xi_2}^b \Psi_{++}(\xi_1+\xi_2,\xi_1,\xi_2) \bigr| \lesssim \jxione^{-a} \jxitwo^{-b} \bigl( \min\{ \jxione, \jxitwo \} \bigr)^{a+b}.
    \end{equation}
Moreover, under the assumption \eqref{equ:derivative_bounds_inverse_phase_sigma_assumption},
for any integers $a,b,c\geq 0$ with $a+b+c>0$,
\begin{equation}  \label{equ:derivative_bounds_phase_pv}
        \Bigl| \nabla_{\xi_1}^a \nabla_{\xi_2}^b \nabla_\sigma^c \Psi_{++}(\xi_1+\xi_2+\sigma,\xi_1,\xi_2) \Big| \lesssim \jxione^{-a} \jxitwo^{-b} \bigl( \min\{ \jxione, \jxitwo \} \bigr)^{a+b}.
\end{equation}
These bounds follow by repeated differentiation, using the standard
symbol estimates for $\langle\cdot\rangle$ and the following elementary estimates for any $x,y \in \bbR^2$,
    \begin{align}
        \biggl| \jap{x} \biggl( \frac{x+y}{\jap{x+y}} - \frac{x}{\jap{x}} \biggr) \biggr| &\lesssim \min \bigl\{ \jap{x}, \jap{y} \bigr\}, \label{equ:elementary1} \\
        \biggl| \jap{x}^2 \biggl( -\frac{1}{\jap{x+y}} + \frac{1}{\jap{x}} \biggr) \biggr| &\lesssim \bigl( \min \bigl\{ \jap{x}, \jap{y} \bigr\} \bigr)^2, \label{equ:elementary2} \\
        \bigl| \jap{x} \jap{x+y}^{-1} \bigr| &\lesssim \min \bigl\{ \jap{x}, \jap{y} \bigr\}. \label{equ:elementary3}
    \end{align}
Then repeated differentiation of $\Psi_{++}^{-1}$, together with the preceding
bounds and Lemma~\ref{lem:lower_bound_quadratic_phase}, yields
\eqref{equ:derivative_bounds_inverse_phase_delta} and
\eqref{equ:derivative_bounds_inverse_phase_pv} in the case $(\kappa_1,\kappa_2) = (+,+)$.

It remains to consider the mixed-sign cases $(\kappa_1,\kappa_2) = (+,-)$ and $(\kappa_1,\kappa_2) = (-,+)$.
We distinguish between comparable input frequencies $\jap{\xi_1} \simeq \jap{\xi_2}$ and separated input frequencies $\jap{\xi_1} \ll \jap{\xi_2}$ or $\jap{\xi_1} \gg \jap{\xi_2}$.
In the comparable regime, the standard symbol bounds and Lemma~\ref{lem:lower_bound_quadratic_phase} suffice.
In the separated regime, if the sign attached to the larger input frequency is positive, the preceding difference estimates yield the same phase-derivative bounds as in the $(+,+)$ case; if it is negative, the phase is elliptic and satisfies
\begin{equation}
    \bigl| \Psi_{\kappa_1\kappa_2}(\xi_1+\xi_2,\xi_1,\xi_2) \bigr| \gtrsim \max\{\langle\xi_1\rangle,\langle\xi_2\rangle\}.
\end{equation}
The same argument applies to the shifted phase under the assumption \eqref{equ:derivative_bounds_inverse_phase_sigma_assumption}, completing the proof.
\end{proof}

\begin{remark}
    The derivative bounds for the phase function established in Lemma~\ref{lem:derivative_bounds_quadratic_phase} are far from optimal; compare, for instance, with the sharper estimates in \cite[Lemma 2.1]{Hayashi_Naumkin16_2DKG}. Nevertheless, the bounds obtained here are entirely sufficient for the multilinear H\"older-type estimates used throughout the nonlinear analysis of this section. Indeed, our bootstrap framework propagates a uniform-in-time bound on a sufficiently high Sobolev norm, allowing any low-frequency derivative losses arising from Lemma~\ref{lem:derivative_bounds_quadratic_phase} to be absorbed without difficulty.    
\end{remark}

\subsubsection{Flat Fourier transform of the degree-one vortex}

Finally, we determine the leading order behavior of the flat Fourier transforms of the real and imaginary components of the degree-one vortex $\underline{\Phi}(r,\theta) = e^{i\theta} U(r)$.

\begin{lemma}[Flat Fourier transforms of the real and imaginary vortex components] \label{lem:FTvortex_re_im}
Let
\[
        \underline{\Phi}_{\mathrm{re}}(r,\theta)
        :=
        \cos(\theta)U(r),
        \qquad
        \underline{\Phi}_{\mathrm{im}}(r,\theta)
        :=
        \sin(\theta)U(r).
\]
Then we have in the sense of tempered distributions for $\zeta = (\zeta_1, \zeta_2) \in \bbR^2$,
\begin{equation} \label{equ:FTvortex_re_im}
    \begin{aligned}
        \widehat{\calF}\bigl[ \underline{\Phi}_{\mathrm{re}} \bigr](\zeta)
        &=
        -i\,\varphi_{\leq 1}(\zeta) \,
        \mathrm{p.v.}\frac{\zeta_1}{|\zeta|^3}
        +
        F_{\mathrm{re}}(\zeta),  \\
        \widehat{\calF}\bigl[ \underline{\Phi}_{\mathrm{im}} \bigr](\zeta)
        &=  
        -i\,\varphi_{\leq 1}(\zeta) \,
        \mathrm{p.v.}\frac{\zeta_2}{|\zeta|^3}
        +
        F_{\mathrm{im}}(\zeta), 
    \end{aligned}
\end{equation}
where the remainder terms $F_{\mathrm{re}}, F_{\mathrm{im}} \in C^\infty(\mathbb R^2)$ satisfy 
\begin{equation} \label{equ:FTvortex_re_im_remainder}
        \bigl| \nabla_\zeta^k F_{\mathrm{re}}(\zeta) \bigr|
        +
        \bigl| \nabla_\zeta^k F_{\mathrm{im}}(\zeta) \bigr|
        \leq
        C_{k}\langle \zeta\rangle^{-4-k} \quad \text{for all } k \geq 0.
\end{equation}
In particular, $\langle x\rangle^m \check F_{\mathrm{re}}, \langle x\rangle^m \check F_{\mathrm{im}} \in L^\infty_x(\mathbb R^2)$ for every $m \geq 0$.
\end{lemma}

\begin{proof}
We identify \(x=(x_1,x_2)\in\mathbb R^2\) with \(x=re^{i\theta}\). Then
\[
        \underline{\Phi}_{\mathrm{re}}(x)
        =
        \frac{x_1}{r}U(r),
        \qquad
        \underline{\Phi}_{\mathrm{im}}(x)
        =
        \frac{x_2}{r}U(r).
\]
We first compute the exact leading singular Fourier transforms of \(x_j/r\), \(j=1,2\).
With the normalization \eqref{eq:FT} for the flat Fourier transform
we have the standard two-dimensional Riesz kernel identity $\widehat{r^{-1}}(\zeta) = |\zeta|^{-1}$.
Thus, we obtain in the sense of tempered distributions,
\begin{equation} \label{equ:FT_xj_over_r_exact}
        \widehat{\calF}\Bigl[ \frac{x_j}{r} \Bigr](\zeta)
        =
        i\partial_{\zeta_j}\bigl(|\zeta|^{-1}\bigr)
        =
        -i\,\mathrm{p.v.}\frac{\zeta_j}{|\zeta|^3},
        \qquad j=1,2.
\end{equation}

It remains to show that, after cutting off the singularity at low frequency, the rest is smooth and
rapidly decaying to the order claimed. Let \(\chi\in C_c^\infty(\mathbb R^2)\) be radial, with
\(\chi(x)=1\) for \(|x|\leq 1\). We use the standard vortex asymptotics
\[
        U(r)=\kappa r+O(r^3)
        \quad\text{as }r\to 0,
        \qquad
        U(r)=1-W(r)
        \quad\text{as }r\to\infty,
\]
where \(W\) and all its derivatives decay exponentially. For \(j=1,2\), define
\[
        S_j(x)
        :=
        \frac{x_j}{r}U(r)\chi(x)
        +
        \frac{x_j}{r}\bigl(U(r)-1\bigr) \bigl(1-\chi(x)\bigr).
\]
Then \(S_j\in \mathcal S(\mathbb R^2)\). Moreover,
\begin{equation} \label{equ:vortex_component_decomp}
        \frac{x_j}{r}U(r)
        =
        \frac{x_j}{r}(1-\chi(x))
        +
        S_j(x).
\end{equation}
Setting
\[
        K_j(\zeta)
        :=
        -i\,\mathrm{p.v.}\frac{\zeta_j}{|\zeta|^3},
\]
we obtain from \eqref{equ:FT_xj_over_r_exact} and \eqref{equ:vortex_component_decomp} that
\begin{equation} \label{equ:component_FT_preliminary_decomp}
        \widehat{\calF}\Bigl[\frac{x_j}{r}U(r)\Bigr](\zeta)
        =
        K_j(\zeta)
        -
        \frac1{2\pi}(K_j*\widehat\chi)(\zeta)
        +
        \widehat{S_j}(\zeta).
\end{equation}
Define
\begin{equation} \label{equ:definition_Fj_remainder}
        F_j(\zeta)
        :=
        \bigl( 1-\varphi_{\leq 1}(\zeta) \bigr) K_j(\zeta)
        -
        \frac1{2\pi}(K_j*\widehat\chi)(\zeta)
        +
        \widehat{S_j}(\zeta),
        \qquad j=1,2.
\end{equation}
Then \eqref{equ:component_FT_preliminary_decomp} becomes
\begin{equation} \label{equ:component_FT_final_decomp}
        \widehat{\calF}\Bigl[\frac{x_j}{r}U(r)\Bigr](\zeta)
        =
        \varphi_{\leq 1}(\zeta)K_j(\zeta)
        +
        F_j(\zeta).
\end{equation}
The term \(F_j\) is smooth near the origin because \((1-\varphi_{\leq 1})K_j\) vanishes near
\(\zeta=0\), \(K_j \ast \widehat\chi\) is the convolution of a tempered distribution with a Schwartz
function, and \(\widehat{S_j}\in\mathcal S(\mathbb R^2)\).

It remains to verify the decay of \(F_j\) at infinity. Since \(\chi\) is radial, hence even, we have
\begin{equation} \label{equ:chi_moments}
        \frac1{2\pi}\int_{\mathbb R^2}\widehat\chi(\eta) \, \ud \eta
        =
        \chi(0)
        =
        1,
        \qquad
        \int_{\mathbb R^2}\eta_m\widehat\chi(\eta) \, \ud \eta
        =
        0,
        \quad m=1,2.
\end{equation}
For \(|\zeta|\gg 1\), the cutoff \(\varphi_{\leq 1}(\zeta)\) vanishes, and hence
\[
        \bigl( 1-\varphi_{\leq 1}(\zeta) \bigr) K_j(\zeta)
        -
        \frac1{2\pi}(K_j \ast \widehat\chi)(\zeta)
        =
        K_j(\zeta)
        -
        \frac1{2\pi}
        \int_{\mathbb R^2}K_j(\zeta-\eta)\widehat\chi(\eta) \,\ud \eta.
\]
Using \eqref{equ:chi_moments} and Taylor expanding \(K_j(\zeta-\eta)\) around \(\zeta\), the zeroth
and first order terms cancel. Since $|\nabla_\zeta^\ell K_j(\zeta)| \lesssim_\ell |\zeta|^{-2-\ell}$ for $|\zeta| \geq 1$ and all $\ell \geq 0$,
the second order Taylor remainder gives
\begin{equation} \label{equ:Kj_cancellation_decay}
        \left|
        \nabla_\zeta^k
        \left(
        \bigl( 1-\varphi_{\leq 1} \bigr) K_j
        -
        \frac1{2\pi}K_j \ast \widehat\chi
        \right)(\zeta)
        \right|
        \lesssim_k
        \langle \zeta \rangle^{-4-k} \quad \text{for all } k \geq 0.
\end{equation}
The rapidly decaying term \(\widehat{S_j}\) satisfies stronger bounds. Therefore, $|\nabla_\zeta^k F_j(\zeta)| \lesssim_k \langle \zeta\rangle^{-4-k}$ and the remaining assertions in the statement of Lemma~\ref{lem:FTvortex_re_im} follow.
\end{proof}

\subsubsection{Bounds for the time derivative of the flat profiles $\bmf_\real(t)$ and $\bmf_\imag(t)$} 

The following lemma decomposes $\pt \bmf_\ast(t)$, $\ast \in \{\real, \imag\}$, into spatially localized and non-localized components and establishes separate high Sobolev and pointwise decay bounds for each.

\begin{lemma} \label{lem:splitting_pt_bmf}
    Suppose the bootstrap assumptions \eqref{equ:bootstrap_assumption_z}, \eqref{equ:bootstrap_assumption_weighted_f}, \eqref{equ:bootstrap_assumption_Linfty_f}, \eqref{equ:bootstrap_assumption_HN_f}, \eqref{equ:bootstrap_assumption_HN_ILED_f} are in place.  
    Then we have for $\ast \in \{\real, \imag\}$ and any $0 \leq t \leq T$,
    \begin{equation} \label{equ:splitting_pt_bmf}
        \partial_t \bmf_\ast(t) = e^{-it\jD} \bigl( \calQ_\ast^{\mathrm{l}}(t) + \calQ_\ast^{\mathrm{nl}}(t) \bigr), 
    \end{equation}
    where 
    \begin{align}
        \bigl\| \jx^2 \calQ_\ast^{\mathrm{l}}(t) \bigr\|_{H^N_x} &\lesssim \varepsilon, \label{equ:calQl_HN_bound} \\ 
        \bigl\| \jD \jx^2 \calQ_\ast^{\mathrm{l}}(t) \bigr\|_{L^\infty_x} &\lesssim \min\bigl\{ \varepsilon, t^{-1+\delta} \bigr\}, \label{equ:calQl_Linfty_bound}
    \end{align}
    and 
    \begin{align}
        \bigl\| \calQ_\ast^{\mathrm{nl}}(t) \bigr\|_{H^N_x} &\lesssim \varepsilon \cdot \min\bigl\{ \varepsilon, t^{-1+2\delta} \bigr\}, \label{equ:calQnl_HN_bound} \\ 
        \bigl\| \jD \calQ_\ast^{\mathrm{nl}}(t) \bigr\|_{L^\infty_x} &\lesssim \Bigl( \min\bigl\{ \varepsilon, t^{-1+2\delta} \bigr\} \Bigr)^2. \label{equ:calQnl_Linfty_bound}
    \end{align}    
\end{lemma}
\begin{proof}
    It suffices to consider the case $\ast = \real$. We decompose the right-hand side of the evolution equation \eqref{equ:evol_equation_flat_profiles_bmf_real_imag} for $\pt \bmf_\real$ into spatially localized and spatially non-localized terms. Recalling also the definition \eqref{equ:definition_bfP_c} of the projection $\bfP_c$, we correspondingly set
    \begin{align}
        \calQ_\real^{\mathrm{l}}(t) &:= (2i\jD)^{-1} \biggl( -\bfV \bmu_\real(t) + \Re\bigl( e^{i\theta} \bmN_d(t) \bigr) - \sum_{j=1}^2 \langle \bmY_j, \bmN(t) \rangle \Re\bigl( e^{i\theta} \bmY_j \bigr) \biggr), \label{equ:calQl_definition_pt_profile_lemma} \\         
        \calQ_\real^{\mathrm{nl}}(t) &:= (2i\jD)^{-1} \Bigl( \Re\bigl( e^{i\theta} \bmN_c(t) \bigr) \Bigr). \label{equ:calQnl_definition_pt_profile_lemma}
    \end{align}
    The high Sobolev norm bound \eqref{equ:calQl_HN_bound} for the first term on the right-hand side of \eqref{equ:calQl_definition_pt_profile_lemma} follows from the high Sobolev norm bootstrap assumption~\eqref{equ:bootstrap_assumption_HN_f} and the bound $\|\jap{r}^2 \bfV(r)\|_{W^{N,\infty}_x} \lesssim 1$, confer \eqref{equ:definition_bfV} and \eqref{equ:definition_V1_V2_for_bfV}.
    The estimate \eqref{equ:calQl_HN_bound} for the second term on the right-hand side of \eqref{equ:calQl_definition_pt_profile_lemma} is a consequence of the spatial localization properties of $e^{i\theta} \bmN_d(t)$, the fact that it can be viewed as an at least quadratic expression in $\bmu_\ast$ and $z_1, z_2$ with at most one derivative on $\bmu_\ast$, and the bootstrap assumptions \eqref{equ:bootstrap_assumption_Linfty_f}, \eqref{equ:bootstrap_assumption_HN_f} along with the bounds on the bad part $\eta_{0,d}(t)$ of the temporal component from Proposition~\ref{prop:eta0d_bounds}. For the last terms on the right-hand side of \eqref{equ:calQl_definition_pt_profile_lemma} it suffices to use the same a priori bounds together with the smoothness and spatial localization of $e^{i\theta} \bfY_j$.

    To obtain the $L^\infty_x$-bound \eqref{equ:calQl_Linfty_bound} for the first term on the right-hand side of \eqref{equ:calQl_definition_pt_profile_lemma}, we place $\jap{r}^2 \bfV(r)$ into $L^\infty_x(\bbR^2)$ and use the a priori decay estimate \eqref{equ:bootstrap_assumption_Linfty_f}, while for all other (at least quadratic) terms on the right-hand side of \eqref{equ:calQl_definition_pt_profile_lemma} it suffices to use their spatial localization along with the bootstrap assumptions \eqref{equ:bootstrap_assumption_Linfty_f}, \eqref{equ:bootstrap_assumption_HN_f} and the bounds from Proposition~\ref{prop:eta0d_bounds} for the bad part of the temporal component $\eta_{0,d}$.

    In order to prove \eqref{equ:calQnl_HN_bound} and \eqref{equ:calQnl_Linfty_bound} for the spatially non-localized terms \eqref{equ:calQnl_definition_pt_profile_lemma}, it suffices to establish the estimates
    \begin{equation}
        \bigl\| e^{i\theta} \bmN_c(t) \bigr\|_{H^{N-1}_x} \lesssim \varepsilon \cdot \min\bigl\{ \varepsilon, t^{-1+2\delta} \bigr\}, \quad \bigl\| e^{i\theta} \bmN_c(t) \bigr\|_{L^\infty_x} \lesssim \Bigl( \min\bigl\{ \varepsilon, t^{-1+2\delta} \bigr\} \Bigr)^2.
    \end{equation}
    These bounds follow from the observation that every term in $e^{i\theta} \bmN_c(t)$ is effectively at least quadratic in the components of $\bmu_\ast(t)$ and $\nabla \bmu_\ast(t)$. This is a consequence of the classification of the nonlinearities in Corollary~\ref{cor:types_of_nonlinearities}, together with the elliptic equations \eqref{equ:eta0c_flat_equation_rewritten} and \eqref{equ:pteta0c_flat_equation_rewritten} for $\eta_{0,c}^\flat(t)$ and $\pt \eta_{0,c}^\flat(t)$, respectively.
    The desired estimates then follow from the bootstrap assumptions \eqref{equ:bootstrap_assumption_Linfty_f}, \eqref{equ:bootstrap_assumption_HN_f}, the Gagliardo--Nirenberg--Sobolev interpolation estimate \eqref{equ:GNSell}, the bounds for $\eta^\flat_{0,c}(t)$ from Proposition~\ref{prop:eta0c_bounds}, and the bounds for $\eta^\sharp_{0,c}(t)$ from Lemma~\ref{lem:eta0c_sharp_in_terms_of_flat}.
\end{proof}

\subsection{Proofs of Propositions~\ref{prop:weighted_energies_g}, \ref{prop:HN_g}, and \ref{prop:HN_ILED_g}} \label{subsec:proofs_propositions_gflat}

In this subsection we prove Propositions~\ref{prop:weighted_energies_g}, \ref{prop:HN_g}, and \ref{prop:HN_ILED_g}. 
We estimate, class by class, the contributions of the ten types of nonlinearities listed in Corollary~\ref{cor:types_of_nonlinearities} to the evolution equations \eqref{equ:evol_equ_bmg_real_imag_flat} for the auxiliary flat profiles $\bmg^\flat_\real(t)$ and $\bmg^\flat_\imag(t)$. 
For each class, we establish the weighted energy estimates, the high Sobolev norm bounds, and the ILED bounds required in the three propositions. 
Together with the explanations from Subsection~\ref{subsec:summary_gflat_bounds}, this completes their proofs. 
Since the arguments for $\bmg^\flat_\real(t)$ and $\bmg^\flat_\imag(t)$ are analogous, we present the details only for $\bmg^\flat_\real(t)$.

\subsubsection{Contributions of type (1)} \label{subsec:contributions_type1}
We begin with the contributions of type (1), which are quadratic terms of the form $\underline{\Phi}_\ast u_{k,\ast_k}(t) u_{l,\ast_l}(t)$ with $\ast, \ast_k, \ast_l \in \{\real,\imag\}$ and $k, l \in \{1, 2, 3, 4\}$. In fact, it suffices to consider all terms of the form $\underline{\Phi}_{\real} u_{k,\real}(t) u_{l,\real}(t)$ with $k, l \in \{1,2,3,4\}$, because $u_{k,\real}(t)$ and $u_{k,\imag}(t)$ enjoy the same bootstrap bounds and because the treatment of $\underline{\Phi}_{\imag}$ is analogous to the treatment of $\underline{\Phi}_{\real}$. 
The contribution of a quadratic term $\underline{\Phi}_{\real} u_{k,\real}(t) u_{l,\real}(t)$ with $k, l \in \{1,2,3,4\}$ to the Fourier transform of the flat profile $\widehat{\bmg}_\real^\flat(t,\xi)$ is given by
\begin{equation}
    \begin{aligned}
        \widehat{\calQ}^{(1)}_{k,l}(t,\xi) &:= \int_t^T (2i\jxi)^{-1} e^{-is\jxi} \widehat{\calF}\bigl[ \underline{\Phi}_{\real} u_{k,\real}(s) u_{l,\real}(s) \bigr](\xi) \, \ud s \\
        &= (2\pi)^2 \sum_{\kappa_1, \kappa_2 \in \{\pm\}} \int_t^T \int_{\bbR^2} \int_{\bbR^2} (2i\jxi)^{-1} e^{is\Psi_{\kappa_1 \kappa_2}(\xi,\xi_1,\xi_2)} \widehat{f}^{\kappa_1}_{k,\real}(s,\xi_1) \widehat{f}^{\kappa_2}_{l,\real}(s,\xi_2) \\
        &\qquad \qquad \qquad \qquad \qquad \qquad \qquad \qquad \qquad \qquad \times \widehat{\underline{\Phi}}_\real(\xi-\xi_1-\xi_2) \, \ud \xi_1 \, \ud \xi_2 \, \ud s,
    \end{aligned}
\end{equation}
where the phase function $\Psi_{\kappa_1 \kappa_2}(\xi,\xi_1,\xi_2)$ is defined as in \eqref{equ:definition_Psi_quadratic}.
In the next three subsections we separately consider the weighted estimate, the high Sobolev norm bound, and the ILED estimate for this contribution.

\medskip 
\paragraph{{\it Weighted estimate for $\widehat{\calQ}^{(1)}_{k,l}(t,\xi)$}} \label{subsubsec:weighted_estimate_contribution_type1}
In view of the decomposition \eqref{equ:FTvortex_re_im} of the Fourier transform of the vortex $\widehat{\underline{\Phi}}_\real$ determined in Lemma~\ref{lem:FTvortex_re_im}, we decompose $\widehat{\calQ}^{(1)}_{k,l}(t,\xi)$ into a singular part, a regular part, and a remainder term
\begin{equation}
    \widehat{\calQ}_{k,l}^{(1)}(t,\xi) = \widehat{\calQ}^{(1)}_{k,l, \mathrm{sing}}(t,\xi) + \widehat{\calQ}^{(1)}_{k,l, \mathrm{reg}}(t,\xi) + \widehat{\calQ}^{(1)}_{k,l, \mathrm{rem}}(t,\xi),
\end{equation}
where the remainder term is given by
\begin{equation}
    \begin{aligned}
        \widehat{\calQ}^{(1)}_{k,l,\mathrm{rem}}(t,\xi) &:= (2\pi)^2 \sum_{\kappa_1, \kappa_2 \in \{\pm\}} \int_t^T \int_{\bbR^2} \int_{\bbR^2} (2i\jxi)^{-1} e^{is\Psi_{\kappa_1 \kappa_2}(\xi,\xi_1,\xi_2)} \widehat{f}^{\kappa_1}_{k,\real}(s,\xi_1) \widehat{f}^{\kappa_2}_{l,\real}(s,\xi_2) \\
        &\qquad \qquad \qquad \qquad \qquad \qquad \qquad \qquad \qquad \qquad \qquad \times F(\xi-\xi_1-\xi_2) \, \ud \xi_1 \, \ud \xi_2 \, \ud s
    \end{aligned}
\end{equation} 
with $F := F_\real \in \calC^\infty(\bbR^2)$ satisfying $|\nabla_\eta^\ell F(\eta)| \lesssim_\ell \jap{\eta}^{-4-\ell}$ for all $\ell \geq 0$ and $\jx^k \check{F} \in L^\infty_x(\bbR^2)$ for all $k \geq 0$ as in the statement of Lemma~\ref{lem:FTvortex_re_im}.
The singular and regular parts stem from the principal value component in $\widehat{\Phi}_\real$.
More precisely, we pick a radially symmetric smooth bump function $\varphi_{\leq 1}(\xi)$ with $\varphi_{\leq 1}(\xi) = 1$ for $|\xi| \leq 1$ and $\varphi_{\leq 1}(\xi) = 0$ for $|\xi| \geq 2$. Set $\varphi_{>1}(\xi) := 1 - \varphi_{\leq 1}(\xi)$. Moreover, we introduce the function 
\begin{equation}
    R(\xi_1, \xi_2) := \frac{\jxione \jxitwo}{\jxione + \jxitwo},
\end{equation}
which is a smoothed out version of $\min\{ \jxione, \jxitwo \}$. Finally, fix $M \gg 1$. Then we define (up to irrelevant constants) the regular part of $\widehat{\calQ}^{(1)}_{k,l}(t,\xi)$ by
\begin{equation}
    \begin{aligned}
        \widehat{\calQ}^{(1)}_{k,l,\mathrm{reg}}(t,\xi) &:= \sum_{\kappa_1, \kappa_2 \in \{\pm\}} \int_t^T \int_{\bbR^2} \int_{\bbR^2} (2i\jxi)^{-1} e^{is\Psi_{\kappa_1 \kappa_2}(\xi,\xi_1,\xi_2)} \widehat{f}^{\kappa_1}_{k,\real}(s,\xi_1) \widehat{f}^{\kappa_2}_{l,\real}(s,\xi_2) \\
        &\qquad \times \varphi_{>-M}\bigl( R(\xi_1,\xi_2) (\xi-\xi_1-\xi_2) \bigr) \varphi_{\leq 1}(\xi-\xi_1-\xi_2) \frac{(\xi-\xi_1-\xi_2)_1}{|\xi-\xi_1-\xi_2|^3} \, \ud \xi_1 \, \ud \xi_2 \, \ud s,
    \end{aligned}
\end{equation}
and the singular part of $\widehat{\calQ}^{(1)}_{k,l}(t,\xi)$ by
\begin{equation}
    \begin{aligned}
        \widehat{\calQ}^{(1)}_{k,l,\mathrm{sing}}(t,\xi) &:= \sum_{\kappa_1, \kappa_2 \in \{\pm\}} \int_t^T \int_{\bbR^2} \int_{\bbR^2} (2i\jxi)^{-1} e^{is\Psi_{\kappa_1 \kappa_2}(\xi,\xi_1,\xi_2)} \widehat{f}^{\kappa_1}_{k,\real}(s,\xi_1) \widehat{f}^{\kappa_2}_{l,\real}(s,\xi_2) \\
        &\qquad \times \varphi_{\leq -M}\bigl( R(\xi_1,\xi_2) (\xi-\xi_1-\xi_2) \bigr) \varphi_{\leq 1}(\xi-\xi_1-\xi_2) \pvdots \, \frac{(\xi-\xi_1-\xi_2)_1}{|\xi-\xi_1-\xi_2|^3} \, \ud \xi_1 \, \ud \xi_2 \, \ud s.
    \end{aligned}
\end{equation}

We first deal with the straightforward estimates for the contribution of the remainder term, and then turn to the contributions of the regular and the singular parts.

\medskip
\noindent \underline{Weighted estimate for the remainder term $\widehat{\calQ}^{(1)}_{k,l,\mathrm{rem}}(t,\xi)$.}
Here we use the dual integrated local energy decay estimate from Lemma~\ref{lem:flat_dual_ILED_low_freq} to obtain that
\begin{align}
    \nonumber
        &\bigl\| \jxi^2 \nabla_\xi \widehat{\calQ}^{(1)}_{k,l,\mathrm{rem}}(t,\xi) \bigr\|_{L^2_\xi(\bbR^2)} 
        \\ \nonumber
        &\lesssim \sum_{\kappa_1, \kappa_2 \in \{\pm\}} \, \biggl\| \int_t^T \frac{\xi}{\jxi} e^{-is\jxi} \cdot s \cdot \jxi \int_{\bbR^2} \int_{\bbR^2} e^{\kappa_1 is \jxione} \widehat{f}^{\kappa_1}_{k,\real}(s,\xi_1) e^{\kappa_2 is \jxitwo} \widehat{f}^{\kappa_2}_{l,\real}(s,\xi_2) 
        \\ \nonumber
        &\qquad \qquad \qquad \qquad \qquad \qquad \qquad \qquad \qquad \times F(\xi-\xi_1-\xi_2) \, \ud \xi_1 \, \ud \xi_2 \, \ud s \biggr\|_{L^2_\xi(\bbR^2)} + \bigl\{ \text{simpler terms} \bigr\} \\ 
        &\lesssim \sum_{\kappa_1, \kappa_2 \in \{\pm\}} \biggl\| s \cdot \Bigl\| \jx \jD \Bigl( \bigl( e^{\kappa_1 i s\jD} f^{\kappa_1}_{k,\real}(s) \bigr) \bigl( e^{\kappa_2 i s\jD} f^{\kappa_2}_{l,\real}(s) \bigr) \check{F} \Bigr) \Bigr\|_{L^2_x} \biggr\|_{L^2_s([t,T])} + \bigl\{ \text{simpler terms} \bigr\} 
        \\
        &\lesssim \Bigl\| s \cdot  \bigl\| \jD e^{is\jD} \bmf_{\real}(s) \bigr\|_{L^\infty_x} \bigl\| e^{is\jD} \bmf_{\real}(s) \bigr\|_{L^\infty_x} \bigl\| \jx \jD \check{F} \bigr\|_{L^2_x} \Bigr\|_{L^2_s([t,T])} + \bigl\{ \text{simpler terms} \bigr\}. \label{eq:GNSprototype}
\end{align}
Invoking \eqref{equ:GNSell}, we conclude that for times $\varepsilon^{-2} \leq t \leq T$ the first term on the right-hand side is bounded by
\begin{equation}\label{eq:GNSprototype1}
    \begin{aligned}
        \big\| s \cdot s^{-1+(3/2)\delta} 
        \cdot s^{-1+\delta} 
        \big\|_{L^2_s([t,T])} \lesssim t^{-\frac12+3\delta}.
    \end{aligned}
\end{equation}
When instead $0 \leq t \leq \varepsilon^{-2}$,
we split the $L^2_s$ integral into the intervals $[t,\varepsilon^{-2}]$ 
and $[\varepsilon^{-2},T]$, and, using also the previous estimate,
we obtain that \eqref{eq:GNSprototype} is bounded by
\begin{equation}\label{eq:GNSprototype2}
    \begin{aligned}
        &\Bigl\| s \cdot \varepsilon^{1-10\delta} \js^{-\frac12-3\delta}  
         \cdot \varepsilon^{1-10\delta} \js^{-\frac12-4\delta}
        \Bigr\|_{L^2_s([t,\varepsilon^{-2}])} 
        + 
        \varepsilon^{1-6\delta} 
        \lesssim \varepsilon^{1-6\delta}.
    \end{aligned}
\end{equation}
In what follows we will encounter several terms that resemble \eqref{eq:GNSprototype} and which lead to estimating
quantities of the form \eqref{eq:GNSprototype1} and \eqref{eq:GNSprototype2}; we can then apply similar arguments and will sometimes omit the details and refer to the above.

\medskip
\noindent \underline{Weighted estimate for the regular part $\widehat{\calQ}^{(1)}_{k,l,\mathrm{reg}}(t,\xi)$.}
To estimate the contribution of the regular part to the weighted energy bounds for $\widehat{\bmg}_\real^\flat(t,\xi)$, we compute 
\begin{equation}
    \begin{aligned}
        &\jxi^2 \nabla_\xi \widehat{\calQ}^{(1)}_{k,l,\mathrm{reg}}(t,\xi) \\ 
        &= \sum_{\kappa_1, \kappa_2 \in \{\pm\}} \frac12 \int_t^T s \cdot \frac{\xi}{\jxi} \cdot e^{-is\jxi} \jxi \int_{\bbR^2} \int_{\bbR^2}  e^{\kappa_1 is \jxione} \widehat{f}^{\kappa_1}_{k,\real}(s,\xi_1) e^{\kappa_2 is \jxitwo} \widehat{f}^{\kappa_2}_{l,\real}(s,\xi_2) \\
        &\qquad \qquad \times \varphi_{>-M}\bigl( R(\xi_1,\xi_2) (\xi-\xi_1-\xi_2) \bigr) \varphi_{\leq 1}(\xi-\xi_1-\xi_2) \frac{(\xi-\xi_1-\xi_2)_1}{|\xi-\xi_1-\xi_2|^3} \, \ud \xi_1 \, \ud \xi_2 \, \ud s \\ 
        &\quad + \sum_{\kappa_1, \kappa_2 \in \{\pm\}} \frac{1}{2i} \int_t^T e^{-is\jxi} \jxi \int_{\bbR^2} \int_{\bbR^2}  e^{\kappa_1 is \jxione} \widehat{f}^{\kappa_1}_{k,\real}(s,\xi_1) e^{\kappa_2 is \jxitwo} \widehat{f}^{\kappa_2}_{l,\real}(s,\xi_2) \\
        &\qquad \qquad \times \nabla_\xi \biggl( \varphi_{>-M}\bigl( R(\xi_1,\xi_2) (\xi-\xi_1-\xi_2) \bigr) \varphi_{\leq 1}(\xi-\xi_1-\xi_2) \frac{(\xi-\xi_1-\xi_2)_1}{|\xi-\xi_1-\xi_2|^3} \biggr) \, \ud \xi_1 \, \ud \xi_2 \, \ud s \\ 
        &\quad + \bigl\{ \text{simpler terms} \bigr\}.
    \end{aligned}
\end{equation}
Using the dual ILED estimate from Lemma~\ref{lem:flat_dual_ILED_low_freq} with $\kappa = \frac12$ for the first term on the right-hand side, we find
\begin{equation}
    \begin{aligned}
        &\bigl\| \jxi^2 \nabla_\xi \widehat{\calQ}^{(1)}_{k,l,\mathrm{reg}}(t,\xi) \bigr\|_{L^2_\xi} \\ 
        &\lesssim \sum_{\kappa_1, \kappa_2 \in \{\pm\}} \Bigl\| s \cdot \bigl\| \frakQ\bigl[ f^{\kappa_1}_{k,\real}(s), f^{\kappa_2}_{l,\real}(s) \bigr](x) \bigr\|_{L^2_x} \Bigr\|_{L^2_s([t,T])} + \sum_{\kappa_1, \kappa_2 \in \{\pm\}} \int_t^T \bigl\| \frakQ\bigl[ f^{\kappa_1}_{k,\real}(s), f^{\kappa_2}_{l,\real}(s) \bigr](x) \bigr\|_{L^2_x} \, \ud s \\ 
        &\quad \quad + \bigl\{ \text{simpler terms} \bigr\}.
    \end{aligned}
\end{equation}
with
\begin{equation}
    \begin{aligned}
        &\frakQ\bigl[ f^{\kappa_1}_{k,\real}(s), f^{\kappa_2}_{l,\real}(s) \bigr](x) \\
        &:= \widehat{\calF}^{-1}\biggl[ \jxi \int_{\bbR^2} \int_{\bbR^2}  e^{\kappa_1 is \jxione} \widehat{f}^{\kappa_1}_{k,\real}(s,\xi_1) e^{\kappa_2 is \jxitwo} \widehat{f}^{\kappa_2}_{l,\real}(s,\xi_2) \\
        &\qquad \qquad \times \nabla_\xi \biggl( \varphi_{>-M}\bigl( R(\xi_1,\xi_2) (\xi-\xi_1-\xi_2) \bigr) \varphi_{\leq 1}(\xi-\xi_1-\xi_2) \frac{(\xi-\xi_1-\xi_2)_1}{|\xi-\xi_1-\xi_2|^3} \biggr) \, \ud \xi_1 \, \ud \xi_2 \biggr](x).
    \end{aligned}
\end{equation}
Here we turned the $\jap{x}$ weight (with $\kappa = \frac12$) on the right-hand side of the dual ILED estimate from Lemma~\ref{lem:flat_dual_ILED_low_freq} into the $\nabla_\xi$ derivative in the definition of $\frakQ\bigl[ f^{\kappa_1}_{k,\real}(s), f^{\kappa_2}_{l,\real}(s) \bigr](x)$ in the preceding equation.
Upon changing variables and dyadically decomposing the frequencies, we can write $\frakQ\bigl[ f^{\kappa_1}_{k,\real}(s), f^{\kappa_2}_{l,\real}(s) \bigr](x)$ as a linear combination of terms
\begin{equation}
    \begin{aligned}
        \frakQ'\bigl[ f^{\kappa_1}_{k,\real}(s), f^{\kappa_2}_{l,\real}(s) \bigr](x) &:= \sum_{0 \leq n_2 \leq n_1} \sum_{-M-n_2 \leq n \leq 0} \int_{\bbR^2_\eta} \int_{\bbR^2_{\xi_1}} \int_{\bbR^2_{\xi_2}} e^{ix\cdot(\xi_1+\xi_2+\eta)} \frakm_{n_1,n_2,n}(\xi_1,\xi_2,\eta) \\
        &\qquad \qquad \times \jxione^3 e^{\kappa_1 is \jxione} \widehat{f}^{\kappa_1}_{k,\real}(s,\xi_1) \jxitwo^3 e^{\kappa_2 is \jxitwo} \widehat{f}^{\kappa_2}_{l,\real}(s,\xi_2) \varphi_{\leq 5}(\eta) \, \ud \xi_1 \, \ud \xi_2 \, \ud \eta 
    \end{aligned}
\end{equation}
with
\begin{equation}
    \begin{aligned}
        \frakm_{n_1,n_2,n}(\xi_1,\xi_2,\eta) &:= \varphi_{n_1}(\xi_1) \varphi_{n_2}(\xi_2) \psi_n(\eta) \jap{\xi_1+\xi_2+\eta} \jap{\xi_1}^{-3} \jxitwo^{-3} \\
        &\qquad \qquad \qquad\times \nabla_\eta \Bigl( \varphi_{>-M}\bigl( R(\xi_1,\xi_2) \eta \bigr) \varphi_{\leq 1}(\eta) \frac{\eta_1}{|\eta|^3} \Bigr),
    \end{aligned}
\end{equation}
where without loss of generality we assumed $n_1 \geq n_2$. 
Under the frequency constraints $0 \leq n_2 \leq n_1$ and $-M-\min\{n_1,n_2\} \leq n \leq 0$, the symbol $\frakm_{n_1,n_2,n}(\xi_1,\xi_2,\eta)$ satisfies
\begin{equation}
    \bigl| \nabla_{\xi_1}^a \nabla_{\xi_2}^b \nabla_\eta^{c} \frakm_{n_1,n_2,n}(\xi_1,\xi_2,\eta) \bigr| \lesssim 2^{-an_1} 2^{-bn_2} 2^{-cn} 2^{-2n_1} 2^{-3(n_2+n)},
\end{equation}
whence 
\begin{equation}
    \bigl\| \widehat{\calF}^{-1}\bigl[\frakm_{n_1,n_2,n}\bigr] \bigr\|_{L^1( (\bbR^2)^3 )} \lesssim 2^{-2n_1} 2^{-3(n_2+n)}.
\end{equation}
Thus, using the trilinear estimates from Lemma~\ref{lem:trilinear_estimates}, 
we obtain 
\begin{equation}
    \begin{aligned}
        &\bigl\| \frakQ'\bigl[ f^{\kappa_1}_{k,\real}(s), f^{\kappa_2}_{l,\real}(s) \bigr](x) \bigr\|_{L^2_x} \\
        &\lesssim \sum_{0 \leq n_2 \leq n_1} \sum_{-M-n_2 \leq n \leq 0} 2^{-2 n_1} 2^{-3(n_2+n)} \bigl\| \jD^3 e^{is\jD} \bmf_\real(s) \bigr\|_{L^\infty_x} \bigl\| \jD^3 e^{is\jD} \bmf_\real(s) \bigr\|_{L^\infty_x} \bigr\| \widehat{\calF}^{-1}\bigl[ \varphi_{\leq 5} \bigr] \bigr\|_{L^2_x} 
        \\ 
        &\lesssim \bigl\| \jD^3 e^{is\jD} \bmf_\real(s) \bigr\|_{L^\infty_x}^2.
    \end{aligned}
\end{equation}
Returning to the weighted energy estimate for the regular part, we arrive at the bound
\begin{equation}
    \begin{aligned}
        \bigl\| \jxi^2 \nabla_\xi \widehat{\calQ}^{(1)}_{k,l,\mathrm{reg}}(t,\xi) \bigr\|_{L^2_\xi}
        \lesssim \Bigl\| s \,
        \bigl\| \jD^3 e^{is\jD} \bmf_\real(s) \bigr\|_{L^\infty_x}^2
        \Bigr\|_{L^2_s([t,T])} 
        + \int_t^T \bigl\| \jD^3 e^{is\jD} \bmf_\real(s) \bigr\|_{L^\infty_x}^2
        \, \ud s 
        \\ 
        \quad \quad + \bigl\{ \text{simpler terms} \bigr\}.
    \end{aligned}
\end{equation}
The first term on the right-hand side yields the largest contribution. Proceeding exactly as for the weighted estimate of $\widehat{\calQ}_{k,l,\mathrm{rem}}^{(1)}(t,\xi)$,
see \eqref{eq:GNSprototype}, 
using \eqref{equ:GNSell} we obtain for long times $\varepsilon^{-2} \leq t \leq T$ that
\begin{equation}
    \bigl\| \jxi^2 \nabla_\xi \widehat{\calQ}^{(1)}_{k,l,\mathrm{reg}}(t,\xi) \bigr\|_{L^2_\xi} \lesssim t^{-\frac12+3\delta},
\end{equation}
while for short times $0 \leq t \leq \varepsilon^{-2}$ we have 
\begin{equation}
    \bigl\| \jxi^2 \nabla_\xi \widehat{\calQ}^{(1)}_{k,l,\mathrm{reg}}(t,\xi) \bigr\|_{L^2_\xi} \lesssim \varepsilon^{1-6\delta}.
\end{equation}
The two bounds above are consistent with the desired bounds stated in Proposition \ref{prop:weighted_energies_g}.

\medskip 
\noindent \underline{Weighted estimate for the singular part $\widehat{\calQ}^{(1)}_{k,l,\mathrm{sing}}(t,\xi)$.}
Integrating by parts in time, we obtain that
\begin{equation} \label{equ:contribution_type1_calQ1_integrated_by_parts}
    \begin{aligned}
        &\widehat{\calQ}^{(1)}_{k,l,\mathrm{sing}}(t,\xi) \\
        &= \sum_{\kappa_1, \kappa_2 \in \{\pm\}} \Bigl[ \widehat{\calB}_{\kappa_1\kappa_2}^{\pvdots,1}\bigl[ f_{k,\real}, f_{l,\real} \bigr](s, \xi) \Bigr]_{s=t}^{s=T} \\ 
        &\quad - \sum_{\kappa_1, \kappa_2 \in \{\pm\}} \int_t^T \int_{\bbR^2} \int_{\bbR^2} (2i\jxi)^{-1} \frac{1}{i\Psi_{\kappa_1\kappa_2}(\xi,\xi_1,\xi_2)} e^{is \Psi_{\kappa_1\kappa_2}(\xi,\xi_1,\xi_2)} \partial_s \widehat{f}^{\kappa_1}_{k,\real}(s,\xi_1) \widehat{f}^{\kappa_2}_{l,\real}(s,\xi_2) \\
        &\qquad \qquad \times \varphi_{\leq -M}\bigl( R(\xi_1,\xi_2) (\xi-\xi_1-\xi_2) \bigr) \varphi_{\leq 1}(\xi-\xi_1-\xi_2) \pvdots \, \frac{(\xi-\xi_1-\xi_2)_1}{|\xi-\xi_1-\xi_2|^3} \, \ud \xi_1 \, \ud \xi_2 \, \ud s \\ 
        &\quad + \bigl\{ \text{similar terms} \bigr\}
    \end{aligned}
\end{equation}
with phase function $\Psi_{\kappa_1\kappa_2}(\xi,\xi_1,\xi_2)$ defined as in \eqref{equ:definition_Psi_quadratic} and boundary terms 
\begin{equation} \label{equ:calQ1_contribution_type1_boundary_term}
    \begin{aligned}
        &\widehat{\calB}_{\kappa_1\kappa_2}^{\pvdots,1}\bigl[ f_{k,\real}, f_{l,\real} \bigr](s,\xi) \\
        &:= \int_{\bbR^2} \int_{\bbR^2} (2i\jxi)^{-1} \frac{1}{i\Psi_{\kappa_1\kappa_2}(\xi,\xi_1,\xi-\xi_1-\sigma)} e^{is \Psi_{\kappa_1\kappa_2}(\xi,\xi_1,\xi-\xi_1-\sigma)}  \\
        &\qquad \times \widehat{f}^{\kappa_1}_{k,\real}(s,\xi_1) \widehat{f}^{\kappa_2}_{l,\real}(s,\xi-\xi_1-\sigma) \varphi_{\leq -M}\bigl( R(\xi_1,\xi-\xi_1-\sigma) \sigma \bigr) \varphi_{\leq 1}(\sigma) \pvdots \, \frac{(\sigma)_1}{|\sigma|^3} \, \ud \xi_1 \, \ud \sigma.
    \end{aligned}
\end{equation}

We begin with the weighted estimate for the boundary terms. By direct computation, we obtain up to irrelevant constants
\begin{equation}
    \begin{aligned}
        &\jxi^2 \nabla_\xi \widehat{\calB}_{\kappa_1\kappa_2}^{\pvdots,1}\bigl[ f_{k,\real}, f_{l,\real} \bigr](s,\xi) \\
        &\simeq s \cdot \int_{\bbR^2} \int_{\bbR^2} \jxi \frac{\nabla_{\xi} \Psi_{\kappa_1\kappa_2}(\xi,\xi_1,\xi-\xi_1-\sigma)}{\Psi_{\kappa_1\kappa_2}(\xi,\xi_1,\xi-\xi_1-\sigma)} e^{is \Psi_{\kappa_1\kappa_2}(\xi,\xi_1,\xi-\xi_1-\sigma)}  \\
        &\qquad \times \widehat{f}^{\kappa_1}_{k,\real}(s,\xi_1) \widehat{f}^{\kappa_2}_{l,\real}(s,\xi-\xi_1-\sigma) \varphi_{\leq -M}\bigl( R(\xi_1,\xi-\xi_1-\sigma) \sigma \bigr) \varphi_{\leq 1}(\sigma) \pvdots \, \frac{(\sigma)_1}{|\sigma|^3} \, \ud \xi_1 \, \ud \sigma \\
        &\quad + \int_{\bbR^2} \int_{\bbR^2} \jxi \frac{1}{\Psi_{\kappa_1\kappa_2}(\xi,\xi_1,\xi-\xi_1-\sigma)} e^{is \Psi_{\kappa_1\kappa_2}(\xi,\xi_1,\xi-\xi_1-\sigma)}  \\
        &\qquad \times \widehat{f}^{\kappa_1}_{k,\real}(s,\xi_1) \nabla_\xi \widehat{f}^{\kappa_2}_{l,\real}(s,\xi-\xi_1-\sigma) \varphi_{\leq -M}\bigl( R(\xi_1,\xi-\xi_1-\sigma) \sigma \bigr) \varphi_{\leq 1}(\sigma) \pvdots \, \frac{(\sigma)_1}{|\sigma|^3} \, \ud \xi_1 \, \ud \sigma \\
        &\quad + \bigl\{ \text{similar and easier terms} \bigr\}.
    \end{aligned}
\end{equation}
It follows that
\begin{equation}
    \begin{aligned}
        &\bigl\| \jxi^2 \nabla_\xi \widehat{\calB}_{\kappa_1\kappa_2}^{\pvdots,1}\bigl[ f_{k,\real}, f_{l,\real} \bigr](s,\xi) \bigr\|_{L^2_\xi} \\
        &\lesssim s \cdot \sum_{n_1,n_2 \geq 0} \, \biggl\| \int_{\bbR^2} \int_{\bbR^2} \int_{\bbR^2} e^{ix\cdot(\xi_1+\xi_2+\sigma)} \frakm_{n_1,n_2}(\xi_1,\xi_2,\sigma) \varphi_{n_1}(\xi_1) e^{is\kappa_1 \jap{\xi_1}} \widehat{f}^{\kappa_1}_{k,\real}(s,\xi_1) \\
        &\qquad \qquad \qquad \qquad \qquad \qquad \times \varphi_{n_2}(\xi_2) e^{is\kappa_2\jap{\xi_2}} \widehat{f}^{\kappa_2}_{l,\real}(s,\xi_2) \pvdots \, \frac{(\sigma)_1}{|\sigma|^3} \, \ud \sigma \, \ud \xi_1 \, \ud \xi_2 \biggr\|_{L^2_x} \\ 
        &\quad + \sum_{n_1,n_2 \geq 0} \, \biggl\| \int_{\bbR^2} \int_{\bbR^2} \int_{\bbR^2} e^{ix\cdot(\xi_1+\xi_2+\sigma)} \frakn_{n_1,n_2}(\xi_1,\xi_2,\sigma) \varphi_{n_1}(\xi_1) e^{is\kappa_1 \jap{\xi_1}} \widehat{f}^{\kappa_1}_{k,\real}(s,\xi_1) \\
        &\qquad \qquad \qquad \qquad \qquad \qquad \times \varphi_{n_2}(\xi_2) e^{is\kappa_2\jap{\xi_2}} \nabla_{\xi_2} \widehat{f}^{\kappa_2}_{l,\real}(s,\xi_2) \pvdots \, \frac{(\sigma)_1}{|\sigma|^3} \, \ud \sigma \, \ud \xi_1 \, \ud \xi_2 \biggr\|_{L^2_x} \\
        &\quad + \bigl\{ \text{similar and easier terms} \bigr\},
    \end{aligned}
\end{equation}
where 
\begin{equation}
    \begin{aligned}
        &\frakm_{n_1,n_2}(\xi_1,\xi_2,\sigma) \\ 
        &\quad := \jap{\xi_1+\xi_2+\sigma} \frac{\nabla_{\xi_2} \Psi_{\kappa_1\kappa_2}(\xi_1+\xi_2+\sigma,\xi_1,\xi_2)}{\Psi_{\kappa_1\kappa_2}(\xi_1+\xi_2+\sigma,\xi_1,\xi_2)} \widetilde{\varphi}_{n_1}(\xi_1) \widetilde{\varphi}_{n_2}(\xi_2) \varphi_{\leq -M}\bigl( R(\xi_1,\xi_2) \sigma \bigr) \varphi_{\leq 1}(\sigma), \\ 
        &\frakn_{n_1,n_2}(\xi_1,\xi_2,\sigma) \\ 
        &\quad := \jap{\xi_1+\xi_2+\sigma} \frac{1}{\Psi_{\kappa_1\kappa_2}(\xi_1+\xi_2+\sigma,\xi_1,\xi_2)} \widetilde{\varphi}_{n_1}(\xi_1) \widetilde{\varphi}_{n_2}(\xi_2) \varphi_{\leq -M}\bigl( R(\xi_1,\xi_2) \sigma \bigr) \varphi_{\leq 1}(\sigma). 
    \end{aligned}
\end{equation}
By Lemma~\ref{lem:derivative_bounds_quadratic_phase} we have for integers $0 \leq a,b,c \leq 3$,
\begin{equation}
    \begin{aligned}
        &\bigl| \nabla_{\xi_1}^a \nabla_{\xi_2}^b \nabla_\sigma^c \frakm_{n_1,n_2}(\xi_1,\xi_2,\sigma) \bigr| + \bigl| \nabla_{\xi_1}^a \nabla_{\xi_2}^b \nabla_\sigma^c \frakn_{n_1,n_2}(\xi_1,\xi_2,\sigma) \bigr| \\
        &\qquad \qquad \qquad \qquad \qquad \qquad \qquad \qquad \lesssim 2^{-a n_1} 2^{-b n_2} 2^{\max\{n_1, n_2\}} 2^{16 \min\{n_1,n_2\}},    
    \end{aligned}
\end{equation}
whence 
\begin{equation}
    \begin{aligned}
        \bigl\| \widehat{\calF}^{-1}\bigl[ \frakm_{n_1,n_2} \bigr] \bigr\|_{L^1( (\bbR^2)^3 )} + \bigl\| \widehat{\calF}^{-1}\bigl[ \frakn_{n_1,n_2} \bigr] \bigr\|_{L^1( (\bbR^2)^3 )} \lesssim 2^{\max\{n_1, n_2\}} 2^{16 \min\{n_1,n_2\}}. 
    \end{aligned}
\end{equation}
Using the bilinear estimates from Lemma~\ref{lem:bilinear_estimates}, we conclude that
\begin{equation}
    \begin{aligned}
        &\bigl\| \jxi^2 \nabla_\xi \widehat{\calB}_{\kappa_1\kappa_2}^{\pvdots,1}\bigl[ f_{k,\real}, f_{l,\real} \bigr](s,\xi) \bigr\|_{L^2_\xi} \\
        &\lesssim s \cdot \sum_{0 \leq n_2 \leq n_1} \, 2^{\max\{n_1, n_2\}} 2^{16 \min\{n_1,n_2\}} \bigl\| P_{n_1} \bmf_{\real}(s) \bigr\|_{L^2_x} \bigl\| e^{is\jD} P_{n_2} \bmf_{\real}(s) \bigr\|_{L^\infty_x} \\  
        &\quad + \sum_{n_1,n_2 \geq 0} \, 2^{\max\{n_1, n_2\}} 2^{16 \min\{n_1,n_2\}} \bigl\| e^{is\jD} P_{n_1} \bmf_{\real}(s) \bigr\|_{L^\infty_x} \bigl\| \varphi_{n_2}(\xi) \nabla_{\xi} \widehat{\bmf}_{\real}(s,\xi) \bigr\|_{L^2_\xi} \\
        &\quad + \bigl\{ \text{similar and easier terms} \bigr\} \\
        &\lesssim s \cdot \|\bmf_{\real}(s)\|_{H^{18}_x} \bigl\| e^{is\jD} \bmf_{\real}(s) \bigr\|_{L^\infty_x} + \bigl\| \jD^{17} e^{is\jD} \bmf_{\real}(s) \bigr\|_{L^\infty_x} \bigl\| \jxi^2 \nabla_{\xi} \widehat{\bmf}_{\real}(s,\xi) \bigr\|_{L^2_\xi} + \ldots 
    \end{aligned}
\end{equation}
where in the first sum we could without loss of generality assume that $0 \leq n_2 \leq n_1$. 
Using the bootstrap assumptions \eqref{equ:bootstrap_assumption_weighted_f} and \eqref{equ:GNSell},
we obtain for times $\varepsilon^{-2} \leq s \leq T$ the bound
\begin{equation}
    \bigl\| \jxi^2 \nabla_\xi \widehat{\calB}_{\kappa_1\kappa_2}^{\pvdots,1}\bigl[ f_{k,\real}, f_{l,\real} \bigr](s,\xi) \bigr\|_{L^2_\xi} \lesssim \varepsilon s^\delta + s^{-\frac12+2\delta} + \bigl\{\text{stronger decaying terms}\bigr\},
\end{equation}
while for short times $0 \leq s \leq \varepsilon^{-2}$, we find that
\begin{equation}
    \begin{aligned}
    \bigl\| \jxi^2 \nabla_\xi \widehat{\calB}_{\kappa_1\kappa_2}^{\pvdots,1}\bigl[ f_{k,\real}, f_{l,\real} \bigr](s,\xi) \bigr\|_{L^2_\xi} &\lesssim \varepsilon^{2-10\delta} \js^{\frac12-4\delta} + \varepsilon^{2-30\delta} \js^{-3\delta} \max\bigl\{ 1, \varepsilon^{1+20\delta} s \bigr\} \\ 
    &\qquad \qquad \qquad \quad \quad \quad + \bigl\{\text{stronger decaying terms}\bigr\}.
    \end{aligned}
\end{equation}

Next, we turn to the weighted estimate for the second main term on the right-hand side of \eqref{equ:contribution_type1_calQ1_integrated_by_parts}.
Here we insert the Fourier transform of the time derivative of the flat profile $\partial_s \bmf_\real(s)$, using the decomposition from Lemma~\ref{lem:splitting_pt_bmf},
\begin{equation} \label{equ:contribution_type1_calQ1_sing_profile_equation_splitting_recalled}
    \begin{aligned}
        \partial_s \widehat{\bmf}_\real(s,\xi_1) &= e^{-is\jap{\xi_1}} \Bigl( \widehat{\calF}\bigl[ \calQ^{\mathrm{l}}_\real(s) \bigr](\xi_1) + \widehat{\calF}\bigl[ \calQ^{\mathrm{nl}}_\real(s) \bigr](\xi_1) \Bigr),
    \end{aligned}
\end{equation}
where we recall that $\calQ^{\mathrm{l}}_\real(s)$ are the spatially localized terms, while all the nonlinearities in 
\begin{equation*}
    \calQ^{\mathrm{nl}}_\real(s) := (2i\jD)^{-1} \Re\bigl( e^{i\theta} \bmN_c(s) \bigr)
\end{equation*}
are not spatially localized.
To estimate the contribution of the spatially localized terms $\calQ^{\mathrm{l}}_\real(s)$, it suffices to use the bounds \eqref{equ:calQl_HN_bound}, \eqref{equ:calQl_Linfty_bound} from Lemma~\ref{lem:splitting_pt_bmf}, whereas the contributions of the non-spatially localized terms $\calQ^{\mathrm{nl}}_\real(s)$ must be treated by exploiting their specific structure according to the classification in Corollary~\ref{cor:types_of_nonlinearities}.

We begin with the contributions from the spatially localized terms $\calQ^{\mathrm{l}}_\real(s)$.
Inserting the contribution of $\calQ^{\mathrm{l}}_\real(s)$ into \eqref{equ:contribution_type1_calQ1_integrated_by_parts}, we obtain for the corresponding weighted energy estimate
\begin{equation}
    \begin{aligned}
        &\biggl\| \jxi^2 \nabla_\xi \sum_{\kappa_1, \kappa_2 \in \{\pm\}} \int_t^T e^{-is\jxi} \int_{\bbR^2} \int_{\bbR^2} (2i\jxi)^{-1} \frac{1}{\Psi_{\kappa_1\kappa_2}(\xi,\xi-\xi_2-\sigma,\xi_2)} \\
        &\qquad \qquad \qquad \qquad \qquad \qquad \times \widehat{\calF}\bigl[ \bigl( \calQ_\real^{\mathrm{l}}(s) \bigr)_k \bigr]^{\kappa_1}(\xi-\xi_2-\sigma) e^{is\kappa_2\jap{\xi_2}} \widehat{f}^{\kappa_2}_{\ell,\real}(s,\xi_2) \\ 
        &\qquad \qquad \qquad \qquad \qquad \qquad \times \varphi_{\leq -M}\bigl( R(\xi-\xi_2-\sigma,\xi_2) \sigma \bigr) \varphi_{\leq 1}(\sigma) \, \pvdots \frac{(\sigma)_1}{|\sigma|^3} \, \ud \xi_2 \, \ud \sigma \, \ud s \biggr\|_{L^2_\xi} \\ 
        &\lesssim \sum_{\kappa_1, \kappa_2 \in \{\pm\}} \biggl\| \int_t^T s \cdot \frac{\xi}{\jxi} \cdot e^{-is\jxi} \jxi \int_{\bbR^2} \int_{\bbR^2} \frac{1}{\Psi_{\kappa_1\kappa_2}(\xi,\xi-\xi_2-\sigma,\xi_2)} \\
        &\qquad \qquad \qquad \qquad \qquad \qquad \times \widehat{\calF}\bigl[ \bigl( \calQ_\real^{\mathrm{l}}(s) \bigr)_k \bigr]^{\kappa_1}(\xi-\xi_2-\sigma) e^{is\kappa_2\jap{\xi_2}} \widehat{f}^{\kappa_2}_{\ell,\real}(s,\xi_2) \\ 
        &\qquad \qquad \qquad \qquad \qquad \qquad \times \varphi_{\leq -M}\bigl( R(\xi-\xi_2-\sigma,\xi_2) \sigma \bigr) \varphi_{\leq 1}(\sigma) \, \pvdots \frac{(\sigma)_1}{|\sigma|^3} \, \ud \xi_2 \, \ud \sigma \, \ud s \biggr\|_{L^2_\xi} \\
        &\quad + \sum_{\kappa_1, \kappa_2 \in \{\pm\}} \biggl\| \int_t^T e^{-is\jxi} \jxi \int_{\bbR^2} \int_{\bbR^2} \frac{1}{\Psi_{\kappa_1\kappa_2}(\xi,\xi-\xi_2-\sigma,\xi_2)} \\
        &\qquad \qquad \qquad \qquad \qquad \qquad \times \nabla \widehat{\calF}\bigl[ \bigl( \calQ_\real^{\mathrm{l}}(s) \bigr)_k \bigr]^{\kappa_1}(\xi-\xi_2-\sigma) e^{is\kappa_2\jap{\xi_2}} \widehat{f}^{\kappa_2}_{\ell,\real}(s,\xi_2) \\ 
        &\qquad \qquad \qquad \qquad \qquad \qquad \times \varphi_{\leq -M}\bigl( R(\xi-\xi_2-\sigma,\xi_2) \sigma \bigr) \varphi_{\leq 1}(\sigma) \, \pvdots \frac{(\sigma)_1}{|\sigma|^3} \, \ud \xi_2 \, \ud \sigma \, \ud s \biggr\|_{L^2_\xi} \\
        &\quad + \bigl\{\text{similar or easier terms}\bigr\}.
    \end{aligned}
\end{equation}
It suffices to provide the details for bounding the first term on the right-hand side of the preceding estimate. Using the dual ILED estimate from Lemma~\ref{lem:flat_dual_ILED_low_freq}, we obtain that it is bounded by
\begin{equation} \label{equ:contribution_type1_calQ1_inserted_calQl}
    \begin{aligned}        
        &\sum_{\kappa_1,\kappa_2 \in \{\pm\}} \, \Biggl\| s \cdot \biggl\| \jx \widehat{\calF}^{-1}_{\xi \mapsto x}\biggl[ \jxi \int_{\bbR^2} \int_{\bbR^2} \frac{1}{\Psi_{\kappa_1\kappa_2}(\xi,\xi-\xi_2-\sigma,\xi_2)} \\
        &\qquad \qquad \qquad \qquad \qquad \qquad \times \widehat{\calF}\bigl[ \bigl( \calQ_\real^{\mathrm{l}}(s) \bigr)_k \bigr]^{\kappa_1}(\xi-\xi_2-\sigma) e^{is\kappa_2\jap{\xi_2}} \widehat{f}^{\kappa_2}_{\ell,\real}(s,\xi_2) \\ 
        &\qquad \qquad \qquad \qquad \qquad \qquad \times \varphi_{\leq -M}\bigl( R(\xi-\xi_2-\sigma,\xi_2) \sigma \bigr) \varphi_{\leq 1}(\sigma) \, \pvdots \frac{(\sigma)_1}{|\sigma|^3} \, \ud \xi_2 \, \ud \sigma \biggr](x) \biggr\|_{L^2_x} \Biggr\|_{L^2_s([t,T])} \\ 
        &\lesssim \sum_{\kappa_1,\kappa_2 \in \{\pm\}} \sum_{n_1,n_2 \geq 0} \, \Biggl\| s \cdot \biggl\| \int_{\bbR^2} \int_{\bbR^2} \int_{\bbR^2} e^{ix\cdot(\xi_1+\xi_2+\sigma)} \frakm_{n_1,n_2}^{\kappa_1\kappa_2}(\xi_1,\xi_2,\sigma) \varphi_{n_1}(\xi_1) \nabla_{\xi_1} \widehat{\calF}\bigl[ \bigl( \calQ_\real^{\mathrm{l}}(s) \bigr)_k \bigr]^{\kappa_1}(\xi_1) \\
        &\qquad \qquad \qquad \qquad \qquad \qquad \qquad \qquad \times \varphi_{n_2}(\xi_2) e^{is\kappa_2\jap{\xi_2}} \widehat{f}^{\kappa_2}_{\ell,\real}(s,\xi_2) \, \pvdots \frac{(\sigma)_1}{|\sigma|^3} \, \ud \xi_1 \, \ud \xi_2 \, \ud \sigma \biggr](x) \biggr\|_{L^2_x} \Biggr\|_{L^2_s([t,T])} \\ 
        &\quad + \bigl\{\text{similar or easier terms}\bigr\}
    \end{aligned}
\end{equation}
with
\begin{equation}
    \begin{aligned}
        &\frakm_{n_1,n_2}^{\kappa_1\kappa_2}(\xi_1,\xi_2,\sigma) \\
        &\quad := \frac{1}{\Psi_{\kappa_1\kappa_2}(\xi_1+\xi_2+\sigma, \xi_1, \xi_2)} \jap{\xi_1+\xi_2+\sigma}  \widetilde{\varphi}_{n_1}(\xi_1) \widetilde{\varphi}_{n_2}(\xi_2) \varphi_{\leq -M}\bigl( R(\xi_1,\xi_2) \sigma \bigr) \varphi_{\leq 1}(\sigma).
    \end{aligned}
\end{equation}
By Lemma~\ref{lem:derivative_bounds_quadratic_phase} we have for integers $0 \leq a, b, c \leq 3$, 
\begin{equation}
    \bigl| \nabla_{\xi_1}^a \nabla_{\xi_2}^b \nabla_\sigma^c \frakm_{n_1,n_2}^{\kappa_1\kappa_2}(\xi_1,\xi_2,\sigma) \bigr| \lesssim 2^{-a n_1} 2^{-b n_2} 2^{\max\{n_1, n_2\}} 2^{16\min\{n_1,n_2\}},
\end{equation}
whence 
\begin{equation}
    \bigl\| \widehat{\calF}^{-1}\bigl[ \frakm_{n_1,n_2}^{\kappa_1\kappa_2} \bigr] \bigr\|_{L^1((\bbR^2)^3)} \lesssim 2^{\max\{n_1, n_2\}} 2^{16\min\{n_1,n_2\}}.
\end{equation}
Using the bilinear estimates from Lemma~\ref{lem:bilinear_estimates}, the Cauchy-Schwarz inequality, and the Gagliardo-Nirenberg-Sobolev interpolation inequality~\eqref{equ:GNS}, it follows that the first term on the right-hand side of \eqref{equ:contribution_type1_calQ1_inserted_calQl} can be crudely bounded by
\begin{equation}
    \begin{aligned}
        &\sum_{n_1,n_2 \geq 0} 2^{\max\{n_1, n_2\}} 2^{16\min\{n_1,n_2\}} \Bigl\| s \cdot \bigl\| P_{n_1} \bigl( \jap{x} \calQ_\real^{\mathrm{l}}(s) \bigr) \bigr\|_{L^2_x} \bigl\| P_{n_2} e^{is\jD} \bmf_\real(s) \bigr\|_{L^\infty_x} \Bigr\|_{L^2_s([t,T])} \\
        &\lesssim \Bigl\| s \cdot \bigl\| \jD^{16} \bigl( \jx \calQ_\real^{\mathrm{l}}(s) \bigr) \bigr\|_{L^2_x} \bigl\| \jD^{16} e^{is\jD} \bmf_\real(s) \bigr\|_{L^\infty_x} \Bigr\|_{L^2_s([t,T])} \\ 
        &\lesssim \Bigl\| s \cdot \bigl\| \jD^{16} \jx^2 \calQ_\real^{\mathrm{l}}(s) \bigr\|_{L^{\infty-}_x} \bigl\| \jD^{16} e^{is\jD} \bmf_\real(s) \bigr\|_{L^\infty_x} \Bigr\|_{L^2_s([t,T])} \\ 
        &\lesssim \Bigl\| s \cdot \bigl\| \jx^2 \calQ_\real^{\mathrm{l}}(s) \bigr\|_{H^N_x}^{\frac{16-\frac{2}{\infty-}}{N-1}} \bigl\| \jx^2 \calQ_\real^{\mathrm{l}}(s) \bigr\|_{L^\infty_x}^{1-\frac{16-\frac{2}{\infty-}}{N-1}} \bigl\|\bmf_\real(s)\bigr\|_{H^N_x}^{\frac{16}{N-1}} \bigl\| e^{is\jD} \bmf_\real(s) \bigr\|_{L^\infty_x}^{1-\frac{16}{N-1}} \Bigr\|_{L^2_s([t,T])}. 
    \end{aligned}
\end{equation}
Using the estimates \eqref{equ:calQl_HN_bound}, \eqref{equ:calQl_Linfty_bound} from Lemma~\ref{lem:splitting_pt_bmf} and the bootstrap assumptions \eqref{equ:bootstrap_assumption_Linfty_f}, \eqref{equ:bootstrap_assumption_HN_f}, 
we infer that for long times $\varepsilon^{-2} \leq t \leq T$, we obtain the estimate $t^{-\frac12+3\delta}$ for the last line. Proceeding analogously, for short times $0 \leq t \leq \varepsilon^{-2}$ we obtain the bound $\varepsilon^{1-6\delta}$.

Next, we turn to the contributions stemming from the non-spatially localized terms $\calQ^{\mathrm{nl}}_\real(s) := (2i\jD)^{-1} \Re\bigl( e^{i\theta} \bmN_c(s) \bigr)$ in \eqref{equ:contribution_type1_calQ1_sing_profile_equation_splitting_recalled}, which are of the types (1)--(10) as listed in Corollary~\ref{cor:types_of_nonlinearities}.
We begin with the contribution of a term of type (1), which without loss of generality we may assume to be of the form $\underline{\Phi}_\real u_{k_1,\real} u_{k_2,\real}$ for some $k_1, k_2 \in \{1, 2, 3, 4\}$. Inserting this term into the second main term on the right-hand side of \eqref{equ:contribution_type1_calQ1_integrated_by_parts} and expressing $u_{k_1,\real}$, $u_{k_2,\real}$ in terms of the evolution of their profiles, we are led to consider weighted estimates $\bigl\| \jxi \nabla_\xi \calI_1(t,\xi) \bigr\|_{L^2_\xi}$ for (slightly more general) terms of the form 
\begin{equation}
    \begin{aligned}
        \calI_1(t,\xi) &:= \int_t^T \int_{\bbR^2} \int_{\bbR^2} \int_{\bbR^2} \int_{\bbR^2} e^{is\Psi_{\iota_1 \iota_2 \iota_3}(\xi,\xi_1,\xi_2,\xi-\xi_1-\xi_2-\eta)} \frakm(\xi,\xi_1,\sigma)   \\
        &\qquad \times \widehat{f}^{\iota_1}_{k_1,\real}(s,\xi_1) \widehat{f}^{\iota_2}_{k_2,\real}(s,\xi_2) \widehat{f}^{\iota_3}_{k_3,\real}(s,\xi-\xi_1-\xi_2-\eta) \widehat{\underline{\Phi}}_{\real}(\eta-\sigma) \, \pvdots \frac{(\sigma)_1}{|\sigma|^3} \, \ud \xi_1 \, \ud \xi_2 \, \ud \eta \, \ud \sigma \, \ud s 
    \end{aligned}
\end{equation}
with 
\begin{equation}
    \begin{aligned}
        \frakm(\xi,\xi_1,\sigma) := \frac{1}{\Psi_{\kappa_1 \kappa_2}(\xi,\xi-\xi_1-\sigma,\xi_1)} \jap{\xi-\xi_1-\sigma}^{-1} \varphi_{\leq -M}\bigl( R(\xi-\xi_1-\sigma,\xi_1) \sigma\bigr) \varphi_{\leq 1}(\sigma)
    \end{aligned}
\end{equation}
for arbitrary $\iota_1, \iota_2, \iota_3 \in \{\pm\}$, $\kappa_1, \kappa_2 \in \{\pm\}$, $k_1, k_2, k_3 \in \{1,2,3,4\}$, where
\begin{equation} \label{equ:definition_Psi_cubic}
    \Psi_{\iota_1 \iota_2 \iota_3}(\xi,\xi_1,\xi_2,\xi_3) := -\jxi + \iota_1 \jap{\xi_1} + \iota_2 \jap{\xi_2} + \iota_3 \jap{\xi_3}.
\end{equation}
To this end we invoke the following commutation identity
\begin{equation} \label{equ:commutation_identity_cubic_pv}
    \begin{aligned}
        &\bigl( \jap{\xi} \nabla_\xi + \iota_1 \jap{\xi_1} \nabla_{\xi_1} + \iota_2 \jap{\xi_2} \nabla_{\xi_2} \bigr) \Psi_{\iota_1 \iota_2 \iota_3}(\xi,\xi_1,\xi_2,\xi-\xi_1-\xi_2-\eta) \\
        &\quad \quad \quad \quad = - \eta - \iota_3 \frac{\xi-\xi_1-\xi_2-\eta}{\jap{\xi-\xi_1-\xi_2-\eta}} \Psi_{\iota_1 \iota_2 \iota_3}(\xi,\xi_1,\xi_2,\xi-\xi_1-\xi_2-\eta).
    \end{aligned}
\end{equation}
Inserting \eqref{equ:commutation_identity_cubic_pv} and integrating by parts in the frequency variables $\xi_1$, $\xi_2$, we find that 
\begin{equation}
    \begin{aligned}
        &\jxi \nabla_\xi \calI_1(t,\xi) \\
        &= \int_t^T \int_{\bbR^2} \int_{\bbR^2} \int_{\bbR^2} \int_{\bbR^2} e^{is\Psi_{\iota_1 \iota_2 \iota_3}(\xi,\xi_1,\xi_2,\xi-\xi_1-\xi_2-\eta)} \frakm(\xi,\xi_1,\sigma)   \\
        &\qquad \quad \times \iota_1 \jap{\xi_1} \nabla_{\xi_1} \widehat{f}^{\iota_1}_{k_1,\real}(s,\xi_1) \widehat{f}^{\iota_2}_{k_2,\real}(s,\xi_2) \widehat{f}^{\iota_3}_{k_3,\real}(s,\xi-\xi_1-\xi_2-\eta) \widehat{\underline{\Phi}}_{\real}(\eta-\sigma) \, \pvdots \frac{(\sigma)_1}{|\sigma|^3} \, \ud \xi_1 \, \ud \xi_2 \, \ud \eta \, \ud \sigma \, \ud s \\ 
        &\quad + \int_t^T \int_{\bbR^2} \int_{\bbR^2} \int_{\bbR^2} \int_{\bbR^2} e^{is\Psi_{\iota_1 \iota_2 \iota_3}(\xi,\xi_1,\xi_2,\xi-\xi_1-\xi_2-\eta)} \frakm(\xi,\xi_1,\sigma)   \\
        &\qquad \quad \times \widehat{f}^{\iota_1}_{k_1,\real}(s,\xi_1) \iota_2 \jap{\xi_2} \nabla_{\xi_2} \widehat{f}^{\iota_2}_{k_2,\real}(s,\xi_2) \widehat{f}^{\iota_3}_{k_3,\real}(s,\xi-\xi_1-\xi_2-\eta) \widehat{\underline{\Phi}}_{\real}(\eta-\sigma) \, \pvdots \frac{(\sigma)_1}{|\sigma|^3} \, \ud \xi_1 \, \ud \xi_2 \, \ud \eta \, \ud \sigma \, \ud s \\ 
        &\quad + \int_t^T \int_{\bbR^2} \int_{\bbR^2} \int_{\bbR^2} \int_{\bbR^2} e^{is\Psi_{\iota_1 \iota_2 \iota_3}(\xi,\xi_1,\xi_2,\xi-\xi_1-\xi_2-\eta)} \frakm(\xi,\xi_1,\sigma)   \\
        &\qquad \quad \times \widehat{f}^{\iota_1}_{k_1,\real}(s,\xi_1) \widehat{f}^{\iota_2}_{k_2,\real}(s,\xi_2) \bigl( \jxi - \iota_1 \jap{\xi_1} - \iota_2 \jap{\xi_2} \bigr) \bigl( \nabla \widehat{f}^{\iota_3}_{k_3,\real}\bigr)(s,\xi-\xi_1-\xi_2-\eta) \widehat{\underline{\Phi}}_{\real}(\eta-\sigma) \\ 
        &\qquad \qquad \qquad \qquad \qquad \qquad \qquad \qquad \qquad \qquad \qquad \qquad \qquad \qquad \qquad \quad \times  \, \pvdots \frac{(\sigma)_1}{|\sigma|^3} \, \ud \xi_1 \, \ud \xi_2 \, \ud \eta \, \ud \sigma \, \ud s \\ 
        &\quad + \int_t^T \int_{\bbR^2} \int_{\bbR^2} \int_{\bbR^2} \int_{\bbR^2} i \cdot s \cdot (-\eta) \cdot e^{is\Psi_{\iota_1 \iota_2 \iota_3}(\xi,\xi_1,\xi_2,\xi-\xi_1-\xi_2-\eta)} \frakm(\xi,\xi_1,\sigma)   \\
        &\qquad \quad \times \widehat{f}^{\iota_1}_{k_1,\real}(s,\xi_1) \widehat{f}^{\iota_2}_{k_2,\real}(s,\xi_2) \widehat{f}^{\iota_3}_{k_3,\real}(s,\xi-\xi_1-\xi_2-\eta) \widehat{\underline{\Phi}}_{\real}(\eta-\sigma) \, \pvdots \frac{(\sigma)_1}{|\sigma|^3} \, \ud \xi_1 \, \ud \xi_2 \, \ud \eta \, \ud \sigma \, \ud s \\ 
        &\quad + \int_t^T \int_{\bbR^2} \int_{\bbR^2} \int_{\bbR^2} \int_{\bbR^2} i \cdot s \cdot (-\iota_3) \cdot \frac{\xi-\xi_1-\xi_2-\eta}{\jap{\xi-\xi_1-\xi_2-\eta}} \cdot \Psi_{\iota_1 \iota_2 \iota_3} \cdot e^{is\Psi_{\iota_1 \iota_2 \iota_3}(\xi,\xi_1,\xi_2,\xi-\xi_1-\xi_2-\eta)} \frakm(\xi,\xi_1,\sigma)   \\
        &\qquad \quad \times \widehat{f}^{\iota_1}_{k_1,\real}(s,\xi_1) \widehat{f}^{\iota_2}_{k_2,\real}(s,\xi_2) \widehat{f}^{\iota_3}_{k_3,\real}(s,\xi-\xi_1-\xi_2-\eta) \widehat{\underline{\Phi}}_{\real}(\eta-\sigma) \, \pvdots \frac{(\sigma)_1}{|\sigma|^3} \, \ud \xi_1 \, \ud \xi_2 \, \ud \eta \, \ud \sigma \, \ud s \\ 
        &\quad + \int_t^T \int_{\bbR^2} \int_{\bbR^2} \int_{\bbR^2} \int_{\bbR^2} e^{is\Psi_{\iota_1 \iota_2 \iota_3}(\xi,\xi_1,\xi_2,\xi-\xi_1-\xi_2-\eta)} \bigl( \jxi \nabla_\xi + \iota_1 \jap{\xi_1} \nabla_{\xi_1} + \iota_2 \jap{\xi_2} \nabla_{\xi_2} \bigr) \frakm(\xi,\xi_1,\sigma)   \\
        &\qquad \quad \times \widehat{f}^{\iota_1}_{k_1,\real}(s,\xi_1) \widehat{f}^{\iota_2}_{k_2,\real}(s,\xi_2) \widehat{f}^{\iota_3}_{k_3,\real}(s,\xi-\xi_1-\xi_2-\eta) \widehat{\underline{\Phi}}_{\real}(\eta-\sigma) \, \pvdots \frac{(\sigma)_1}{|\sigma|^3} \, \ud \xi_1 \, \ud \xi_2 \, \ud \eta \, \ud \sigma \, \ud s \\ 
        &=: \calJ_1(t,\xi) + \ldots + \calJ_6(t,\xi).
    \end{aligned}
\end{equation}
For the first term $\calJ_1(t,\xi)$, we find that 
\begin{equation}
    \begin{aligned}
        &\bigl\|\calJ_1(t,\xi)\bigr\|_{L^2_\xi} \\
        &\lesssim \sum_{n_1, n_2, n_3, n_4 \geq 0} \int_t^T \biggl\| \int_{\bbR^2} \int_{\bbR^2} \int_{\bbR^2} \int_{\bbR^2} \int_{\bbR^2} e^{ix\cdot(\xi_1+\xi_2+\xi_3+\xi_4+\sigma)} \frakn_{n_1,n_2,n_3,n_4}(\xi_1,\xi_2,\xi_3,\xi_4,\sigma) \\ 
        &\qquad \qquad \qquad \qquad \times \varphi_{n_1}(\xi_1) e^{is\iota_1\jap{\xi_1}} \jap{\xi_1} \nabla_{\xi_1} \widehat{f}^{\iota_1}_{k_1,\real}(s,\xi_1) \varphi_{n_2}(\xi_2) e^{is\iota_2\jap{\xi_2}} \widehat{f}^{\iota_2}_{k_2,\real}(s,\xi_2) \\
        &\qquad \qquad \qquad \qquad \times \varphi_{n_3}(\xi_3) e^{is\iota_3\jap{\xi_3}} \widehat{f}^{\iota_3}_{k_3,\real}(s,\xi_3) \varphi_{n_4}(\xi_4) \widehat{\underline{\Phi}}_{\real}(\xi_4) \, \pvdots \frac{(\sigma)_1}{|\sigma|^3} \, \ud \xi_1 \, \ud \xi_2 \, \ud \xi_3 \, \ud \xi_4 \, \ud \sigma \biggr\|_{L^2_x} \ud s
    \end{aligned}
\end{equation}
with 
\begin{equation} \label{equ:weighted_energy_estimate_type1_reinserted_type1_frakn}
    \begin{aligned}
        &\frakn_{n_1,n_2,n_3,n_4}(\xi_1,\xi_2,\xi_3,\xi_4,\sigma) \\
        &:= \frac{1}{\Psi_{\kappa_1 \kappa_2}(\xi_1+\xi_2+\xi_3+\xi_4+\sigma,\xi_2+\xi_3+\xi_4,\xi_1)} \jap{\xi_2+\xi_3+\xi_4}^{-1} \\
        &\quad \quad \quad \times \varphi_{\leq -M}\bigl( R(\xi_2+\xi_3+\xi_4, \xi_1) \sigma\bigr)  
        \widetilde{\varphi}_{n_1}(\xi_1) \widetilde{\varphi}_{n_2}(\xi_2) \widetilde{\varphi}_{n_3}(\xi_3) \widetilde{\varphi}_{n_4}(\xi_4) 
        \varphi_{\leq 1}(\sigma).
    \end{aligned}
\end{equation}
By Lemma~\ref{lem:derivative_bounds_quadratic_phase} we have for any integers $a, b, c, d, e \geq 0$ the bounds
\begin{equation}
    \begin{aligned}
        &\bigl| \nabla_{\xi_1}^a \nabla_{\xi_2}^b \nabla_{\xi_3}^c \nabla_{\xi_4}^d \nabla_{\sigma}^e \frakn_{n_1, n_2, n_3, n_4}(\xi_1,\xi_2,\xi_3,\xi_4,\sigma) \bigr| \\
        &\quad \lesssim \jap{\xi_1}^{-a} \jap{\xi_2+\xi_3+\xi_4}^{-b-c-d} \bigl( \min \bigl\{ \jap{\xi_1}, \jap{\xi_2+\xi_3+\xi_4} \bigr\} \bigr)^{1+2(a+b+c+d)+e} \\
        &\quad \lesssim \jap{\xi_1}^{-a} \jap{\xi_2}^{-b} \jap{\xi_3}^{-c} \jap{\xi_4}^{-d} \bigl( \min\bigl\{ \jap{\xi_2}, \jap{\xi_3+\xi_4} \bigr\} \bigr)^b \bigl( \min \bigl\{ \jap{\xi_3}, \jap{\xi_2+\xi_4} \bigr\} \bigr)^c  \\ 
        &\quad \quad \quad \times \bigl( \min\bigl\{ \jap{\xi_4}, \jap{\xi_2+\xi_3} \bigr\} \bigr)^d \bigl( \min \bigl\{ \jap{\xi_1}, \jap{\xi_2+\xi_3+\xi_4} \bigr\} \bigr)^{1+2(a+b+c+d)+e} \\
        &\quad \lesssim 2^{-an_1} 2^{-bn_2} 2^{-cn_3} 2^{-dn_4} \cdot 2^{bn_2} \cdot 2^{cn_3} \cdot 2^{d n_4} \cdot 2^{(1+2(a+b+c+d)+e) \max\{n_2,n_3,n_4\}}.
    \end{aligned}
\end{equation}
The latter imply for integers $0 \leq a, b, c, d, e \leq 3$ the crude, but sufficient bounds 
\begin{equation}
    \begin{aligned}
        \bigl| \nabla_{\xi_1}^a \nabla_{\xi_2}^b \nabla_{\xi_3}^c \nabla_{\xi_4}^d \nabla_{\sigma}^e \frakn_{n_1, n_2, n_3, n_4}(\xi_1,\xi_2,\xi_3,\xi_4,\sigma) \bigr| \lesssim 2^{-an_1} 2^{-bn_2} 2^{-cn_3} 2^{-dn_4} \cdot 2^{31 n_2} \cdot 2^{31 n_3} \cdot 2^{31 n_4},
    \end{aligned}
\end{equation}
whence 
\begin{equation} \label{equ:weighted_energy_estimate_type1_reinserted_type1_frakn_inverse_FT}
    \begin{aligned}
        \bigl\| \widehat{\calF}^{-1}\bigl[\frakn_{n_1, n_2, n_3, n_4}\bigr] \bigr\|_{L^1( (\bbR^2)^5 )} \lesssim 2^{31 n_2} \cdot 2^{31 n_3} \cdot 2^{31 n_4}.
    \end{aligned}
\end{equation}
By the quadrilinear estimates from Lemma~\ref{lem:quadrilinear_estimates} and the Gagliardo-Nirenberg-Sobolev inequality \eqref{equ:GNS}, we conclude that
\begin{equation}
    \begin{aligned}
        &\bigl\|\calJ_1(t,\xi)\bigr\|_{L^2_\xi} \\
        &\lesssim \sum_{n_1, n_2, n_3, n_4 \geq 0} \int_t^T 2^{31 n_2} 2^{31 n_3} 2^{31 n_4} \bigl\| \varphi_{n_1}(\xi_1) \jap{\xi_1} \nabla_{\xi_1} \widehat{\bmf}_\real(s,\xi_1) \bigr\|_{L^2_{\xi_1}} \bigl\| P_{n_2} e^{is\jD} \bmf_\real(s) \bigr\|_{L^\infty_x} \\
        &\qquad \qquad \qquad \qquad \qquad \qquad \qquad \qquad \qquad \qquad \qquad \times \bigl\| P_{n_3} e^{is\jD} \bmf_\real(s) \bigr\|_{L^\infty_x} \bigl\| P_{n_4} \Phi_{\real} \bigr\|_{L^\infty_x} \, \ud s \\ 
        &\lesssim \int_t^T \bigl\| \jap{\xi_1}^2 \nabla_{\xi_1} \widehat{\bmf}_\real(s) \bigr\|_{L^2_{\xi_1}} \bigl\| \jD^{32} e^{is\jD} \bmf_\real(s) \bigr\|_{L^\infty_x}^2 \bigl\| \jD^{32} \Phi_{\real} \bigr\|_{L^\infty_x} \, \ud s 
        .
    \end{aligned}
\end{equation}
Using the bootstrap assumptions \eqref{equ:bootstrap_assumption_weighted_f}, 
and \eqref{equ:GNSell}, we conclude for times $t \geq \varepsilon^{-2}$ that 
\begin{equation}\label{equ:J1A}
    \begin{aligned}
        \bigl\|\calJ_1(t,\xi)\bigr\|_{L^2_\xi} 
        & \lesssim \int_t^T s^{\frac12} 
        \cdot \bigl( s^{-1+\frac32\delta} \bigr)^2 \, \ud s \lesssim t^{-\frac12+3\delta}.
    \end{aligned}
\end{equation}
For times $\varepsilon^{-1-20\delta} \leq t \leq \varepsilon^{-2}$ we obtain
\begin{equation}\label{equ:J1B}
    \begin{aligned}
        \bigl\|\calJ_1(t,\xi)\bigr\|_{L^2_\xi} 
        &\lesssim \int_t^{\varepsilon^{-2}} \varepsilon^2 s^{\frac32} 
        \cdot 
        \bigl( \varepsilon^{1-10\delta} s^{-\frac12-\frac72\delta} \bigr)^{2} 
        \, \ud s +
        \int_{\varepsilon^{-2}}^T s^{\frac12} \cdot 
        \bigl( s^{-1+\frac32\delta} \bigr)^{2} \, \ud s 
        \lesssim 
        \varepsilon^{1-6\delta}.
    \end{aligned}
\end{equation}
Finally, for times $0 \leq t \leq \varepsilon^{-1-20\delta}$ we estimate the contribution with $s \in [t,\varepsilon^{-1-20\delta}]$ and use \eqref{equ:J1A}  and \eqref{equ:J1B} for the contributions with 
$s \geq \varepsilon^{-1-20\delta}$, to find that
\begin{equation}
    \begin{aligned}
        \bigl\|\calJ_1(t,\xi)\bigr\|_{L^2_\xi} 
        &\lesssim \int_t^{\varepsilon^{-1-20\delta}} \varepsilon^{1-20\delta} \js^{\frac12} \cdot 
        \bigl( \varepsilon^{1-10\delta} \js^{-\frac12-\frac72\delta} \bigr)^2 \, \ud s  + \varepsilon^{1-6\delta}
        &\lesssim \varepsilon^{\frac52-50\delta} 
        + \varepsilon^{1-6\delta} \lesssim \varepsilon^{1-6\delta}.
    \end{aligned}
\end{equation}

The terms $\calJ_2(t,\xi)$ and $\calJ_3(t,\xi)$ can be estimated similarly and we omit the details. We now turn to the term $\calJ_4(t,\xi)$, where we exploit the frequency factor $\eta$. Writing $\eta = (\eta-\sigma) + \sigma$, we decompose $\calJ_4(t,\xi)$ into 
\begin{equation}
    \begin{aligned}
        &\calJ_4(t,\xi) \\
        &= - \int_t^T \int_{\bbR^2} \int_{\bbR^2} \int_{\bbR^2} \int_{\bbR^2} i \cdot s \cdot e^{is\Psi_{\iota_1 \iota_2 \iota_3}(\xi,\xi_1,\xi_2,\xi-\xi_1-\xi_2-\eta)} \frakm(\xi,\xi_1,\sigma)   \\
        &\qquad \quad \times \widehat{f}^{\iota_1}_{k_1,\real}(s,\xi_1) \widehat{f}^{\iota_2}_{k_2,\real}(s,\xi_2) \widehat{f}^{\iota_3}_{k_3,\real}(s,\xi-\xi_1-\xi_2-\eta) (\eta-\sigma) \widehat{\underline{\Phi}}_{\real}(\eta-\sigma) \, \pvdots \frac{(\sigma)_1}{|\sigma|^3} \, \ud \xi_1 \, \ud \xi_2 \, \ud \eta \, \ud \sigma \, \ud s \\ 
        &\quad - \int_t^T \int_{\bbR^2} \int_{\bbR^2} \int_{\bbR^2} \int_{\bbR^2} i \cdot s \cdot e^{is\Psi_{\iota_1 \iota_2 \iota_3}(\xi,\xi_1,\xi_2,\xi-\xi_1-\xi_2-\eta)} \frakm(\xi,\xi_1,\sigma)   \\
        &\qquad \quad \times \widehat{f}^{\iota_1}_{k_1,\real}(s,\xi_1) \widehat{f}^{\iota_2}_{k_2,\real}(s,\xi_2) \widehat{f}^{\iota_3}_{k_3,\real}(s,\xi-\xi_1-\xi_2-\eta) \widehat{\underline{\Phi}}_{\real}(\eta-\sigma) \, \sigma \frac{(\sigma)_1}{|\sigma|^3} \, \ud \xi_1 \, \ud \xi_2 \, \ud \eta \, \ud \sigma \, \ud s \\ 
        &=: \calJ_{4,1}(t,\xi) + \calJ_{4,2}(t,\xi).
    \end{aligned}
\end{equation}
Then to estimate $\calJ_{4,1}(t,\xi)$, we observe that in view of the $r^{-1}$ decay of angular derivatives and the fact that the vortex $\underline{\Phi}(x) = e^{i\theta} U(r)$ is smooth \cite{Taubes80_1}, we have 
\begin{equation}
    \bigl| \widehat{\calF}^{-1}\bigl[ \xi \widehat{\underline{\Phi}}_\real(\xi) \bigr](x) \bigr| \lesssim \jx^{-1}, 
\end{equation}
which can be placed into $L^{2+}_x(\bbR^2)$. Thus, with $\frakn_{n_1,n_2,n_3,n_4}(\xi_1,\xi_2,\xi_3,\xi_4,\sigma)$ defined as in \eqref{equ:weighted_energy_estimate_type1_reinserted_type1_frakn} and using Lemma~\ref{lem:quadrilinear_estimates} with \eqref{equ:weighted_energy_estimate_type1_reinserted_type1_frakn_inverse_FT}, we obtain the crude but sufficient bound
\begin{equation}
    \begin{aligned}
        &\bigl\|\calJ_{4,1}(t,\xi)\bigr\|_{L^2_\xi} \\
        &\lesssim \sum_{n_1, n_2, n_3, n_4 \geq 0} \int_t^T s \cdot \biggl\| \int_{\bbR^2} \int_{\bbR^2} \int_{\bbR^2} \int_{\bbR^2} \int_{\bbR^2} e^{ix\cdot(\xi_1+\xi_2+\xi_3+\xi_4+\sigma)} \frakn_{n_1,n_2,n_3,n_4}(\xi_1,\xi_2,\xi_3,\xi_4,\sigma) \\ 
        &\qquad \qquad \qquad \qquad \times \varphi_{n_1}(\xi_1) e^{is\iota_1\jap{\xi_1}} \widehat{f}^{\iota_1}_{k_1,\real}(s,\xi_1) \varphi_{n_2}(\xi_2) e^{is\iota_2\jap{\xi_2}} \widehat{f}^{\iota_2}_{k_2,\real}(s,\xi_2) \\
        &\qquad \qquad \qquad \qquad \times \varphi_{n_3}(\xi_3) e^{is\iota_3\jap{\xi_3}} \widehat{f}^{\iota_3}_{k_3,\real}(s,\xi_3) \varphi_{n_4}(\xi_4) \xi_4 \widehat{\underline{\Phi}}_{\real}(\xi_4) \, \pvdots \frac{(\sigma)_1}{|\sigma|^3} \, \ud \xi_1 \, \ud \xi_2 \, \ud \xi_3 \, \ud \xi_4 \, \ud \sigma \biggr\|_{L^2_x} \ud s \\ 
        &\lesssim \sum_{n_1, n_2, n_3, n_4 \geq 0} \int_t^T s \cdot 2^{31 n_2} 2^{31 n_3} 2^{31 n_4} \bigl\| P_{n_1} e^{is\jD} \bmf_\real(s) \bigr\|_{L^{\infty-}_x} \bigl\| P_{n_2} e^{is\jD} \bmf_\real(s) \bigr\|_{L^\infty_x} \\
        &\qquad \qquad \qquad \qquad \qquad \qquad \qquad \qquad \times \bigl\| P_{n_3} e^{is\jD} \bmf_\real(s) \bigr\|_{L^\infty_x} \bigl\| P_{n_4} \widehat{\calF}^{-1}\bigl[ \xi_4 \widehat{\underline{\Phi}}_{\real}(\xi_4) \bigr] \bigr\|_{L^{2+}_x} \, \ud s \\ 
        &\lesssim \int_t^T s \cdot \bigl\| \jD e^{is\jD} \bmf_\real(s) \bigr\|_{L^{\infty-}_x} \bigl\| \jD^{32} e^{is\jD} \bmf_\real(s) \bigr\|_{L^\infty_x}^2 \bigl\| \jD^{32} \widehat{\calF}^{-1}\bigl[ \xi_4 \widehat{\underline{\Phi}}_{\real}(\xi_4) \bigr] \bigr\|_{L^{2+}_x} \, \ud s 
        . 
    \end{aligned}
\end{equation}
Using \eqref{equ:GNSell}, we conclude, for long times $\varepsilon^{-2} \leq t \leq T$, the crude but more than sufficient bound
\begin{equation}
    \begin{aligned}
        \bigl\|\calJ_{4,1}(t,\xi)\bigr\|_{L^2_\xi} &\lesssim \int_t^T s^{-2+4\delta} \, \ud s \lesssim t^{-1+4\delta}.
    \end{aligned}
\end{equation}
Instead for short times $0 \leq t \leq \varepsilon^{-2}$, using also the last bound above, we find 
\begin{equation}
    \begin{aligned}
        \bigl\|\calJ_{4,1}(t,\xi)\bigr\|_{L^2_\xi} &\lesssim \varepsilon^{3-30\delta} \int_t^{\varepsilon^{-2}} \js^{-\frac12} \, \ud s 
        + \varepsilon^{2-8\delta} \lesssim \varepsilon^{2-30\delta}.
    \end{aligned}
\end{equation}

For the term $\calJ_{4,2}(t,\xi)$ we first recall that for some absolute constant $c$ we have for $j, l \in \{1,2\}$,
\begin{equation}
    \begin{aligned}
        \widehat{\calF}^{-1}\biggl[ \frac{(\sigma)_j (\sigma)_l}{|\sigma|^3} \biggr](x) = c \biggl( \frac{\delta_{j l}}{|x|} - \frac{x_j x_l}{|x|^3} \biggr),
    \end{aligned}
\end{equation}
which can be placed into $L^{2-}_x(\bbR^2) + L^{2+}_x(\bbR^2)$. 
Thus, spending a little bit of derivative to go down from $L^2_x(\bbR^2)$ to $L^{2-}_x(\bbR^2)$, and then proceeding similarly as above, using the quintilinear estimate from Lemma~\ref{lem:quintiilinear_estimates}, we can crudely bound by
\begin{equation}
    \begin{aligned}
        &\bigl\|\calJ_{4,2}(t,\xi)\bigr\|_{L^2_\xi} \\
        &\lesssim \sum_{n_1, n_2, n_3, n_4 \geq 0} \int_t^T s \cdot \biggl\| \jD^{0+} \int_{\bbR^2} \int_{\bbR^2} \int_{\bbR^2} \int_{\bbR^2} \int_{\bbR^2} e^{ix\cdot(\xi_1+\xi_2+\xi_3+\xi_4+\sigma)} \frakn_{n_1,n_2,n_3,n_4}(\xi_1,\xi_2,\xi_3,\xi_4,\sigma) \\ 
        &\qquad \qquad \qquad \qquad \times \varphi_{n_1}(\xi_1) e^{is\iota_1\jap{\xi_1}} \widehat{f}^{\iota_1}_{k_1,\real}(s,\xi_1) \varphi_{n_2}(\xi_2) e^{is\iota_2\jap{\xi_2}} \widehat{f}^{\iota_2}_{k_2,\real}(s,\xi_2) \\
        &\qquad \qquad \qquad \qquad \times \varphi_{n_3}(\xi_3) e^{is\iota_3\jap{\xi_3}} \widehat{f}^{\iota_3}_{k_3,\real}(s,\xi_3) \varphi_{n_4}(\xi_4) \widehat{\underline{\Phi}}_{\real}(\xi_4) \sigma \frac{(\sigma)_1}{|\sigma|^3} \, \ud \xi_1 \, \ud \xi_2 \, \ud \xi_3 \, \ud \xi_4 \, \ud \sigma \biggr\|_{L^{2-}_x} \ud s \\ 
        &\lesssim \sum_{n_1, n_2, n_3, n_4 \geq 0} \int_t^T s \cdot 2^{n_1} 2^{32 n_2} 2^{32 n_3} 2^{32 n_4} \bigl\| P_{n_1} e^{is\jD} \bmf_\real(s) \bigr\|_{L^{\infty-}_x \cap L^\infty_x} \bigl\| P_{n_2} e^{is\jD} \bmf_\real(s) \bigr\|_{L^\infty_x} \\
        &\qquad \qquad \qquad \qquad \qquad \qquad \qquad \times \bigl\| P_{n_3} e^{is\jD} \bmf_\real(s) \bigr\|_{L^\infty_x} \bigl\| P_{n_4} \underline{\Phi}_\real \bigr\|_{L^\infty_x} \Bigl\| \widehat{\calF}^{-1} \Bigl[ \sigma \frac{(\sigma)_1}{|\sigma|^3} \Bigr] \Bigr\|_{L^{2-}_x + L^{2+}_x}  \, \ud s \\ 
        &\lesssim \int_t^T s \cdot \bigl\| \jD^2 e^{is\jD} \bmf_\real(s) \bigr\|_{L^{\infty-}_x \cap L^\infty_x} \bigl\| \jD^{33} e^{is\jD} \bmf_\real(s) \bigr\|_{L^\infty_x}^2 \, \ud s. 
    \end{aligned}
\end{equation}
The last line gives rise to favorable bounds as for the preceding term $\calJ_{4,1}(t,\xi)$. We omit the details.

Finally, in order to estimate the term $\calJ_5(t,\xi)$ we exploit the presence of the phase factor to first integrate by parts in time. Upon reinserting the evolution equations for the profiles, using the bounds from Lemma~\ref{lem:splitting_pt_bmf}, the resulting terms have a lot of room and we omit the details for the lengthy, but straightforward estimates. The last term $\calJ_6(t,\xi)$ is rather mild and we omit the details.

Next, we continue with the contribution of a term of type (2) from $\calQ_\real^{\mathrm{nl}}(s)$ in \eqref{equ:contribution_type1_calQ1_sing_profile_equation_splitting_recalled} when inserted back into the second main term on the right-hand side of \eqref{equ:contribution_type1_calQ1_integrated_by_parts}. 
Without loss of generality, we may assume a term of type (2) to be of the form $u_{k_1,\real} \partial_j u_{k_2,\real}$ for some $k_1, k_2 \in \{1, 2, 3, 4\}$ and $j \in \{1, 2\}$. 
Expressing $u_{k_1,\real}$, $u_{k_2,\real}$ in terms of the evolution of their profiles, we then have to consider weighted estimates $\bigl\| \jxi \nabla_\xi \calI_2(t,\xi) \bigr\|_{L^2_\xi}$ for general terms of the form 
\begin{equation}
    \begin{aligned}
        \calI_2(t,\xi) &:= \int_t^T \int_{\bbR^2} \int_{\bbR^2} \int_{\bbR^2} \int_{\bbR^2} e^{is\Psi_{\iota_1 \iota_2 \iota_3}(\xi,\xi_1,\xi_2,\xi-\xi_1-\xi_2-\sigma)} \frakm(\xi,\xi_1,\xi_2,\sigma)   \\
        &\quad \times \jap{\xi_1+\xi_2}^{-1} (\xi_2)_j \, \widehat{f}^{\iota_1}_{k_1,\real}(s,\xi_1) \widehat{f}^{\iota_2}_{k_2,\real}(s,\xi_2) \widehat{f}^{\iota_3}_{k_3,\real}(s,\xi-\xi_1-\xi_2-\sigma) \, \pvdots \frac{(\sigma)_1}{|\sigma|^3} \, \ud \xi_1 \, \ud \xi_2 \, \ud \sigma \, \ud s 
    \end{aligned}
\end{equation}
with 
\begin{equation}
    \begin{aligned}
        \frakm(\xi,\xi_1,\xi_2,\sigma) := \frac{1}{\Psi_{\kappa_1 \kappa_2}(\xi,\xi_1+\xi_2,\xi-\xi_1-\xi_2-\sigma)} \varphi_{\leq -M}\bigl( R(\xi_1+\xi_2,\xi-\xi_1-\xi_2-\sigma) \sigma\bigr) \varphi_{\leq 1}(\sigma)
    \end{aligned}
\end{equation}
for arbitrary $\iota_1, \iota_2, \iota_3 \in \{\pm\}$, $\kappa_1, \kappa_2 \in \{\pm\}$, $k_1, k_2, k_3 \in \{1,2,3,4\}$, and where $\Psi_{\iota_1 \iota_2 \iota_3}(\xi,\xi_1,\xi_2,\xi_3)$ is defined as in \eqref{equ:definition_Psi_cubic}.
To this end we proceed similarly as for the term $\calI_1(t,\xi)$ using the commutation identity \eqref{equ:commutation_identity_cubic_pv} (with $\eta$ replaced by $\sigma$).
In view of the additional frequency factor $(\xi_2)_j$, some additional care is necessary for those terms when a frequency derivative $\nabla_{\xi_2}$ falls onto the second input $\widehat{f}^{\iota_2}_{k_2,\real}(s,\xi_2)$. Dyadically localizing $|\xi_1| \sim 2^{n_1}$ and $|\xi_2| \sim 2^{n_2}$ for integers $n_1, n_2 \geq 0$, we have $|\jap{\xi_1+\xi_2}^{-1} (\xi_2)_j| \lesssim 1$ whenever $|n_1-n_2| \geq 5$. Only in the case of high-high interactions $|n_1 - n_2| < 5$, we may have $|\xi_1+\xi_2| \ll |\xi_2|$, but in that case additional frequency factors on the second input can be shifted to the first input. 
The lengthy details are similar to the treatment of $\calI_1(t,\xi)$, and are thus omitted.

We continue with the discussion of the contributions of terms of type (3) and (4) from $\calQ_\real^{\mathrm{nl}}(s)$ in \eqref{equ:contribution_type1_calQ1_sing_profile_equation_splitting_recalled} when inserted back into the second main term on the right-hand side of \eqref{equ:contribution_type1_calQ1_integrated_by_parts}. These are strictly easier to deal with than terms of type (1) or (2), because terms of type (3) are cubic, while terms of type (4) have additional spatial localization owing to the term $(1-a_\theta)$. We therefore omit the details.

Next, we consider the contribution of terms of type (5) of the form $\partial_j \partial_t \eta_{0,c}$ for $j \in \{1,2\}$ from $\calQ_\real^{\mathrm{nl}}(s)$ in \eqref{equ:contribution_type1_calQ1_sing_profile_equation_splitting_recalled} when inserted back into the second main term on the right-hand side of \eqref{equ:contribution_type1_calQ1_integrated_by_parts}. 
This leads to weighted energy estimates $\bigl\| \jxi \nabla_\xi \calI_5(t,\xi) \bigr\|_{L^2_\xi}$ for terms of the form
\begin{equation} \label{equ:weighted_energy_estimate_type1_reinserted_type5_calI5}
    \begin{aligned}
        &\calI_5(t,\xi) := \int_t^T \int_{\bbR^2} \int_{\bbR^2} \frac{1}{i\Psi_{\kappa_1\kappa_2}(\xi,\xi_1,\xi_2)} e^{-is\jxi} \bigl( \jap{\xi_1}^{-1} (\xi_1)_j \bigr) \widehat{\calF}\bigl[ \pt \eta_{0,c}^\flat(s) + \pt \eta_{0,c}^\sharp(s) \bigr]^{\kappa_1}(\xi_1) \\
        &\quad \times e^{is\kappa_2 \jap{\xi_2}} \hatf_{l,\real}^{\kappa_2}(s,\xi_2) \varphi_{\leq-M}\bigl( R(\xi_1,\xi_2) (\xi-\xi_1-\xi_2) \bigr) \varphi_{\leq 1}(\xi-\xi_1-\xi_2) \pvdots \frac{(\xi-\xi_1-\xi_2)_1}{|\xi-\xi_1-\xi_2|^3} \, \ud \xi_1 \, \ud \xi_2 \, \ud s.
    \end{aligned}
\end{equation}
In view of the spatial localization of $\partial_t \eta_{0,c}^\sharp(s) = (-\Delta + U^2)^{-1} \bigl( (1-U^2) \pt \eta_{0,c}^\flat(s) \bigr)$, see Lemma~\ref{lem:eta0c_sharp_in_terms_of_flat},
and the fast $L^\infty_x$-decay for $\pt \eta_{0,c}^\flat(s)$ obtained in Proposition~\ref{prop:eta0c_bounds}, the contribution of $\pt \eta_{0,c}^\sharp(s)$ in \eqref{equ:weighted_energy_estimate_type1_reinserted_type5_calI5} can be easily controlled using the dual ILED estimate from Lemma~\ref{lem:flat_dual_ILED_low_freq}. 
To estimate the contribution of $\pt \eta_{0,c}^\flat(s)$, we write, using \eqref{equ:pteta0c_flat_equation_rewritten} and \eqref{equ:pteta0c_flat_equation_leading_order_term_rewritten},
\begin{equation} \label{equ:weighted_energy_estimate_type1_reinserted_type5_pteta0c_flat}
    \begin{aligned}
        \pt \eta_{0,c}^\flat 
        &= (-\Delta + 1)^{-1} \Bigl( \nabla \cdot \bigl( u_{1,\real} \nabla u_{3,\real} - u_{3,\real} \nabla u_{1,\real} \bigr) + \nabla \cdot \bigl( u_{1,\imag} \nabla u_{3,\imag} - u_{3,\imag} \nabla u_{1,\imag} \bigr)  \\
        &\quad \quad \quad \quad \quad \quad \quad -2 U' \bigl( u_{2,\real} u_{3,\real} + u_{2,\imag} u_{3,\imag} \bigr) + 2 U' \bigl( u_{1,\real} u_{4,\real} + u_{1,\imag} u_{4,\imag} \bigr) \\
        &\quad \quad \quad \quad \quad \quad \quad -u_{3,\real} \Re\bigl( e^{i\theta} (\bfP_c \bmN)_1 \bigr) -u_{3,\imag} \Im\bigl(e^{i\theta} (\bfP_c \bmN)_1\bigr) \\
        &\quad \quad \quad \quad \quad \quad \quad + u_{1,\real} \Re\bigl(e^{i\theta} (\bfP_c \bmN)_3\bigr) + u_{1,\imag} \Im\bigl(e^{i\theta} (\bfP_c \bmN)_3\bigr) \\
        &\quad \quad \quad \quad \quad \quad \quad - 2 \underline{\Phi}_\real \pt \eta_{0,c} u_{1,\real} - 2 \underline{\Phi}_\imag \pt \eta_{0,c} u_{1,\imag} - 2 \underline{\Phi}_\real  \eta_{0,c} \pt u_{1,\real} - 2 \underline{\Phi}_\imag \eta_{0,c} \pt u_{1,\imag} \\ 
        &\quad \quad \quad \quad \quad \quad \quad - \pt \eta_{0,c} \bigl( u_{1,\real}^2 + u_{1,\imag}^2 + u_{3,\real}^2 + u_{3,\imag}^2 \bigr) \\
        &\quad \quad \quad \quad \quad \quad \quad - \eta_{0,c} \bigl( 2 u_{1,\real} \pt u_{1,\real} + 2 u_{1,\imag} \pt u_{1,\imag} + 2 u_{3,\real} \pt u_{3,\real} + 2 u_{3,\imag} \pt u_{3,\imag} \bigr) \Bigr).
    \end{aligned}
\end{equation}
The first line on the right-hand side of \eqref{equ:weighted_energy_estimate_type1_reinserted_type5_pteta0c_flat} when inserted back into \eqref{equ:weighted_energy_estimate_type1_reinserted_type5_calI5} produces trilinear terms, for which the corresponding weighted energy estimates can be obtained similarly as before using the commutation identity \eqref{equ:commutation_identity_cubic_pv}. 
A minor point worth commenting on is that the derivative in the nonlinearities such as $u_{1,\real} \nabla u_{3,\real}$ in the first line on the right-hand side of \eqref{equ:weighted_energy_estimate_type1_reinserted_type5_pteta0c_flat} gets compensated by the output $(-\Delta+1)^{-1} \nabla$ apart from the case of high-high interactions. But for high-high interactions derivatives can be transferred between the two inputs so that there is no danger of a pile-up of too many derivatives on an input onto which a frequency derivative falls. 

The quadratic nonlinearities in the second line on the right-hand side of \eqref{equ:weighted_energy_estimate_type1_reinserted_type5_pteta0c_flat} are spatially localized, thus their contributions when inserted back into \eqref{equ:weighted_energy_estimate_type1_reinserted_type5_calI5} can be easily estimated using the dual ILED estimate from Lemma~\ref{lem:flat_dual_ILED_low_freq}

The nonlinearities in the remaining lines on the right-hand side of \eqref{equ:weighted_energy_estimate_type1_reinserted_type5_pteta0c_flat} are at least of cubic order, and their contributions are thus easier to bound. We omit the lengthy details for these harmless contributions.

\medskip 

The contributions of terms of the type (6)--(10) from $\calQ_\real^{\mathrm{nl}}(s)$ in \eqref{equ:contribution_type1_calQ1_sing_profile_equation_splitting_recalled} when inserted back into the second main term on the right-hand side of \eqref{equ:contribution_type1_calQ1_integrated_by_parts} are similar or easier to the previously discussed cases. We leave the details to the reader.

\medskip
\paragraph{{\it High Sobolev norm estimate for $\widehat{\calQ}^{(1)}_{k,l}(t,\xi)$}} \label{subsubsec:calQ1_high_Sobolev}

For the high Sobolev norm estimate for $\widehat{\calQ}^{(1)}_{k,l}(t,\xi)$, we have to take into account the limited decay \eqref{equ:FTvortex_re_im_remainder} on the frequency side of the remainder term in the decomposition~\eqref{equ:FTvortex_re_im} of the flat Fourier transform of the components of the vortex $\underline{\Phi}$ obtained in Lemma~\ref{lem:FTvortex_re_im}. We therefore distinguish the case when at least two derivatives fall on the vortex, producing a sufficiently spatially localized term, and the case when all (apart from at most one) derivatives fall onto the inputs. Correspondingly, we write
\begin{equation} \label{equ:contribution1_HN_decomposition}
    \begin{aligned}
        &\bigl\| \jxi^N \widehat{\calQ}_{k,l}^{(1)}(t,\xi) \bigr\|_{L^2_\xi} \\
        &\lesssim \sum_{|\alpha| \leq N-1} \biggl\| \int_t^T e^{-is\jD} \partial^\alpha \bigl( \underline{\Phi}_{\real} u_{k,\real}(s) u_{l,\real}(s) \bigr) \, \ud s \biggr\|_{L^2_x} \\ 
        &\lesssim 
        \sum_{0 \leq |\beta| \leq 1} \sum_{|\gamma| \leq N-1-|\beta|} \biggl\| \int_t^T e^{-is\jD} \Bigl( \bigl( \partial^\beta \underline{\Phi}_{\real} \bigr) \partial^\gamma \bigl( u_{k,\real}(s) u_{l,\real}(s) \bigr) \Bigr) \, \ud s \biggr\|_{L^2_x} \\
        &\quad + \sum_{2 \leq |\beta| \leq N-1} \sum_{|\gamma| \leq N-1-|\beta|} \biggl\| \int_t^T e^{-is\jD} \Bigl( \bigl( \partial^\beta \underline{\Phi}_{\real} \bigr) \partial^\gamma \bigl( u_{k,\real}(s) u_{l,\real}(s) \bigr) \Bigr) \, \ud s \biggr\|_{L^2_x},
    \end{aligned}
\end{equation}
where $\alpha, \beta, \gamma \in \bbN_0^2$ are multi-indices.
Observe that $\bigl\| \jx^{|\beta|} \partial^\beta \underline{\Phi}_{\real} \bigr\|_{L^\infty_x} \lesssim 1$ for any multi-index $|\beta| \geq 1$.
Owing to this spatial localization, the estimate for the second term on the right-hand side of \eqref{equ:contribution1_HN_decomposition} is straightforward using the dual ILED estimate from Lemma~\ref{lem:flat_dual_ILED} with $0 < \kappa \ll 1$,
\begin{equation}
    \begin{aligned}
        &\sum_{2 \leq |\beta| \leq N-1} \sum_{|\gamma| \leq N-1-|\beta|} \biggl\| \int_t^T e^{-is\jD} \Bigl( \bigl( \partial^\beta \underline{\Phi}_{\real} \bigr) \partial^\gamma \bigl( u_{k,\real}(s) u_{l,\real}(s) \bigr) \Bigr) \, \ud s \biggr\|_{L^2_x} \\
        &\lesssim \sum_{2 \leq |\beta| \leq N-1} \sum_{|\gamma| \leq N-1-|\beta|} \Bigl\| \bigl\| \jx^{1+\kappa} \bigl( \partial^\beta \underline{\Phi}_{\real} \bigr) \partial^\gamma \bigl( u_{k,\real}(s) u_{l,\real}(s) \bigr) \bigr\|_{L^2_x} \Bigr\|_{L^2_s([t,T])} \\ 
        &\lesssim \sum_{2 \leq |\beta| \leq N-1} \bigl\|\jx^{1+\kappa} \bigl( \partial^\beta \underline{\Phi}_{\real} \bigr) \bigr\|_{L^\infty_x} 
        \Bigl\| \bigl\| e^{is\jD} \bmf_\real(s) \bigr\|_{L^\infty_x} \|\bmf_\real(s)\|_{H^{N-3}_x} \Bigr\|_{L^2_s([t,T])}.
    \end{aligned}
\end{equation}
Using \eqref{equ:bootstrap_assumption_Linfty_f}, \eqref{equ:bootstrap_assumption_HN_f} we obtain for large times $\varepsilon^{-2} \leq t \leq T$ that
\begin{equation}
    \begin{aligned}
        \sum_{2 \leq |\beta| \leq N-1} \sum_{|\gamma| \leq N-1-|\beta|} \biggl\| \int_t^T e^{-is\jD} \Bigl( \bigl( \partial^\beta \underline{\Phi}_{\real} \bigr) \partial^\gamma \bigl( u_{k,\real}(s) u_{l,\real}(s) \bigr) \Bigr) \, \ud s \biggr\|_{L^2_x} 
        &\lesssim \bigl\| s^{-1+\delta} \cdot \varepsilon \bigr\|_{L^2_s([t,T])} \\
        &\lesssim \varepsilon t^{-\frac12+\delta} \lesssim \varepsilon^{2-2\delta},
    \end{aligned}
\end{equation}
while we find for short times $0 \leq t \leq \varepsilon^{-2}$ that
\begin{equation}
    \begin{aligned}
        &\sum_{2 \leq |\beta| \leq N-1} \sum_{|\gamma| \leq N-1-|\beta|} \biggl\| \int_t^T e^{-is\jD} \Bigl( \bigl( \partial^\beta \underline{\Phi}_{\real} \bigr) \partial^\gamma \bigl( u_{k,\real}(s) u_{l,\real}(s) \bigr) \Bigr) \, \ud s \biggr\|_{L^2_x} \\
        &\lesssim \bigl\| \varepsilon^{1-10\delta} \js^{-\frac12-4\delta} \cdot \varepsilon \bigr\|_{L^2_s([t,\varepsilon^{-2}])} + \bigl\| s^{-1+\delta} \cdot \varepsilon \bigr\|_{L^2_s([\varepsilon^{-2},T])} \lesssim \varepsilon^{2-10\delta} + \varepsilon^{2-2\delta} \lesssim \varepsilon^{2-10\delta}.
    \end{aligned}
\end{equation}

For the first term on the right-hand side of \eqref{equ:contribution1_HN_decomposition}, we discuss the case $\beta = 0$ when all derivatives fall on the inputs. The other case $|\beta| = 1$ is analogous.
For any multi-index $\gamma = (\gamma_1,\gamma_2) \in \bbN_0^2$ with $|\gamma| \leq N-1$ this leads to the following terms on the Fourier side
\begin{equation} \label{equ:calQ1gamma_def}
    \begin{aligned}
        \widehat{\calQ}^{(1), \gamma}_{k,l}(t,\xi) &:= \sum_{\kappa_1, \kappa_2 \in \{\pm\}} \int_t^T \int_{\bbR^2} \int_{\bbR^2} e^{is\Psi_{\kappa_1 \kappa_2}(\xi,\xi_1,\xi_2)} \bigl( i (\xi_1+\xi_2) \bigr)^\gamma  \\
        &\qquad \qquad \qquad \qquad \times \widehat{f}^{\kappa_1}_{k,\real}(s,\xi_1) \widehat{f}^{\kappa_2}_{l,\real}(s,\xi_2) \widehat{\underline{\Phi}}_\real(\xi-\xi_1-\xi_2) \, \ud \xi_1 \, \ud \xi_2 \, \ud s,
    \end{aligned}
\end{equation}
where we use the notation $(\xi_1+\xi_2)^\gamma = (\xi_1+\xi_2)_1^{\gamma_1} (\xi_1+\xi_2)_2^{\gamma_2}$.
Proceeding as at the beginning of Subsection~\ref{subsubsec:weighted_estimate_contribution_type1}, we insert the decomposition of the flat Fourier transform of $\underline{\Phi}_\real$ from Lemma~\ref{lem:FTvortex_re_im}. Correspondingly, we obtain a decomposition of $\widehat{\calQ}^{(1), \gamma}_{k,l}(t,\xi)$ into a singular part, a regular part, and a remainder term
\begin{equation} \label{equ:calQ1gamma_decomposition}
    \widehat{\calQ}_{k,l}^{(1), \gamma}(t,\xi) = \widehat{\calQ}^{(1), \gamma}_{k,l, \mathrm{sing}}(t,\xi) + \widehat{\calQ}^{(1), \gamma}_{k,l, \mathrm{reg}}(t,\xi) + \widehat{\calQ}^{(1), \gamma}_{k,l, \mathrm{rem}}(t,\xi),
\end{equation}
where the singular part is given by
\begin{equation}
    \begin{aligned}
        \widehat{\calQ}^{(1),\gamma}_{k,l,\mathrm{sing}}(t,\xi) &:= \sum_{\kappa_1, \kappa_2 \in \{\pm\}} \int_t^T \int_{\bbR^2} \int_{\bbR^2} e^{is\Psi_{\kappa_1 \kappa_2}(\xi,\xi_1,\xi_2)} \bigl( i (\xi_1+\xi_2) \bigr)^\gamma \widehat{f}^{\kappa_1}_{k,\real}(s,\xi_1) \widehat{f}^{\kappa_2}_{l,\real}(s,\xi_2) \\
        &\qquad \times \varphi_{\leq -M}\bigl( R(\xi_1,\xi_2) (\xi-\xi_1-\xi_2) \bigr) \varphi_{\leq 1}(\xi-\xi_1-\xi_2) \, \pvdots \, \frac{(\xi-\xi_1-\xi_2)_1}{|\xi-\xi_1-\xi_2|^3} \, \ud \xi_1 \, \ud \xi_2 \, \ud s,
    \end{aligned}
\end{equation}
the regular part reads
\begin{equation}
    \begin{aligned}
        \widehat{\calQ}^{(1), \gamma}_{k,l,\mathrm{reg}}(t,\xi) &:= \sum_{\kappa_1, \kappa_2 \in \{\pm\}} \int_t^T \int_{\bbR^2} \int_{\bbR^2} e^{is\Psi_{\kappa_1 \kappa_2}(\xi,\xi_1,\xi_2)} \bigl( i (\xi_1+\xi_2) \bigr)^\gamma \widehat{f}^{\kappa_1}_{k,\real}(s,\xi_1) \widehat{f}^{\kappa_2}_{l,\real}(s,\xi_2) \\
        &\qquad \times \varphi_{>-M}\bigl( R(\xi_1,\xi_2) (\xi-\xi_1-\xi_2) \bigr) \varphi_{\leq 1}(\xi-\xi_1-\xi_2) \frac{(\xi-\xi_1-\xi_2)_1}{|\xi-\xi_1-\xi_2|^3} \, \ud \xi_1 \, \ud \xi_2 \, \ud s,
    \end{aligned}
\end{equation}
and the remainder term is
\begin{equation}
    \begin{aligned}
        \widehat{\calQ}^{(1), \gamma}_{k,l,\mathrm{rem}}(t,\xi) &:= \sum_{\kappa_1, \kappa_2 \in \{\pm\}} \int_t^T \int_{\bbR^2} \int_{\bbR^2} e^{is\Psi_{\kappa_1 \kappa_2}(\xi,\xi_1,\xi_2)} \bigl( i (\xi_1+\xi_2) \bigr)^\gamma \widehat{f}^{\kappa_1}_{k,\real}(s,\xi_1) \widehat{f}^{\kappa_2}_{l,\real}(s,\xi_2) \\
        &\qquad \qquad \qquad \qquad \qquad \qquad \qquad \qquad \qquad \qquad \qquad \times F(\xi-\xi_1-\xi_2) \, \ud \xi_1 \, \ud \xi_2 \, \ud s.
    \end{aligned}
\end{equation}
We remind the reader of the notation $R(\xi_1,\xi_2) := \jap{\xi_1} \jap{\xi_2} (\jap{\xi_1} + \jap{\xi_2})^{-1}$. Moreover, we recall from the statement of Lemma~\ref{lem:FTvortex_re_im} that $F := F_\real \in \calC^\infty(\bbR^2)$ satisfies $|\nabla_\eta^\ell F(\eta)| \lesssim_\ell \jap{\eta}^{-4-\ell}$ for all $\ell \geq 0$ and $\jx^k \check{F} \in L^\infty_x(\bbR^2)$ for all $k \geq 0$.

We begin with the straightforward estimates for the contribution of the remainder term and of the regular part, and subsequently turn to the contributions of the singular part.

\medskip
\noindent \underline{High Sobolev norm estimate for the remainder term $\widehat{\calQ}^{(1), \gamma}_{k,l,\mathrm{rem}}(t,\xi)$.}
Using the dual ILED estimate from Lemma~\ref{lem:flat_dual_ILED}, and recalling that $|\gamma| \leq N-1$, we obtain for some $0 < \kappa \ll 1$,
\begin{equation}
    \begin{aligned}
        \bigl\| \widehat{\calQ}^{(1), \gamma}_{k,l,\mathrm{rem}}(t,\xi) \bigr\|_{L^2_\xi} &\lesssim \Bigl\| \bigl\| \jx^{1+\kappa} \partial^\gamma \bigl( e^{is\kappa_1 \jD} f^{\kappa_1}_{k,\real}(s) e^{is\kappa_2 \jD} f^{\kappa_2}_{l,\real}(s) \bigr) \check{F} \bigr\|_{L^2_x} \Bigr\|_{L^2_s([t,T])} \\ 
        &\lesssim \Bigl\| \bigl\| \partial^\gamma \bigl( e^{is\kappa_1 \jD} f^{\kappa_1}_{k,\real}(s) e^{is\kappa_2 \jD} f^{\kappa_2}_{l,\real}(s) \bigr) \bigr\|_{L^2_x} \bigl\| \jx^{1+\kappa} \check{F} \bigr\|_{L^\infty_x} \Bigr\|_{L^2_s([t,T])} \\ 
        &\lesssim \Bigl\| \|\bmf(s)\|_{H^{N-1}_x} \bigl\|e^{is\jD} \bmf(s)\bigr\|_{L^\infty_x} \Bigr\|_{L^2_s([t,T])}.
    \end{aligned}
\end{equation}
Using \eqref{equ:bootstrap_assumption_Linfty_f}, \eqref{equ:bootstrap_assumption_HN_f}, we obtain for large times $\varepsilon^{-2} \leq t \leq T$ that
\begin{equation}
    \begin{aligned}
        \bigl\| \widehat{\calQ}^{(1), \gamma}_{k,l,\mathrm{rem}}(t,\xi) \bigr\|_{L^2_\xi} &\lesssim \bigl\| \varepsilon \cdot s^{-1+\delta} \bigr\|_{L^2_s([t,T])} \lesssim \varepsilon t^{-\frac12+\delta} \lesssim \varepsilon^{2-2\delta},
    \end{aligned}
\end{equation}
while for short times $0 \leq t \leq \varepsilon^{-2}$ we conclude 
\begin{equation}
    \begin{aligned}
        \bigl\| \widehat{\calQ}^{(1), \gamma}_{k,l,\mathrm{rem}}(t,\xi) \bigr\|_{L^2_\xi} &\lesssim \bigl\| \varepsilon \cdot \varepsilon^{1-10\delta} \js^{-\frac12-4\delta} \bigr\|_{L^2_s([t,\varepsilon^{-2}])} + \bigl\| \varepsilon \cdot s^{-1+\delta} \bigr\|_{L^2_s([\varepsilon^{-2},T])} \lesssim \varepsilon^{2-10\delta} + \varepsilon^{2-2\delta} \lesssim \varepsilon^{2-10\delta}.
    \end{aligned}
\end{equation}

\medskip
\noindent \underline{High Sobolev norm estimate for the regular part $\widehat{\calQ}^{(1), \gamma}_{k,l,\mathrm{reg}}(t,\xi)$.}
Using the dual ILED estimate from Lemma~\ref{lem:flat_dual_ILED}, here we have
\begin{equation} \label{equ:high_sobolev_calQ1_regular_part}
    \begin{aligned}
        \bigl\| \widehat{\calQ}^{(1),\gamma}_{k,l,\mathrm{reg}}(t,\xi) \bigr\|_{L^2_\xi} &\lesssim \sum_{\kappa_1, \kappa_2 \in \{\pm\} } \Bigl\| \bigl\| \jx^2 \frakQ^{\gamma,0}\bigl[ f^{\kappa_1}_{k,\real}(s), f^{\kappa_2}_{l,\real}(s) \bigr] \bigr\|_{L^2_x} \Bigr\|_{L^2_s([t,T])} \\ 
        &\lesssim \sum_{0 \leq m \leq 2} \sum_{\kappa_1, \kappa_2 \in \{\pm\}} \Bigl\| \bigl\| \frakQ^{\gamma,m} \bigl[ f^{\kappa_1}_{k,\real}(s), f^{\kappa_2}_{l,\real}(s) \bigr] \bigr\|_{L^2_x} \Bigr\|_{L^2_s([t,T])} 
    \end{aligned}
\end{equation}
with
\begin{equation}
    \begin{aligned}
        &\frakQ^{\gamma, m}\bigl[ f^{\kappa_1}_{k,\real}(s), f^{\kappa_2}_{l,\real}(s) \bigr](x) \\
        &:= \widehat{\calF}^{-1}\biggl[ \int_{\bbR^2} \int_{\bbR^2} \bigl( i (\xi_1+\xi_2) \bigr)^\gamma e^{\kappa_1 is\jxione} \widehat{f}^{\kappa_1}_{k,\real}(s,\xi_1) e^{\kappa_2 is\jxitwo} \widehat{f}^{\kappa_2}_{l,\real}(s,\xi_2) \\
        &\qquad \qquad \times \nabla_\xi^m \biggl( \varphi_{>-M}\bigl( R(\xi_1,\xi_2) (\xi-\xi_1-\xi_2) \bigr) \varphi_{\leq 1}(\xi-\xi_1-\xi_2) \frac{(\xi-\xi_1-\xi_2)_1}{|\xi-\xi_1-\xi_2|^3} \biggr) \, \ud \xi_1 \, \ud \xi_2 \biggr](x).
    \end{aligned}
\end{equation}
We provide the details for the estimates in the case of $m=2$, the other two cases leading to more favorable bounds.
Changing variables and dyadically decomposing the frequencies, we can write $\frakQ^{\gamma,2}\bigl[ f^{\kappa_1}_{k,\real}(s), f^{\kappa_2}_{l,\real}(s) \bigr](x)$ as a linear combination of terms of the form
\begin{equation}
    \begin{aligned}
        \frakQ^{\gamma,2}_{hl}\bigl[ f^{\kappa_1}_{k,\real}(s), f^{\kappa_2}_{l,\real}(s) \bigr](x) &:= \sum_{0 \leq n_2 \leq n_1} \sum_{-n_2 \leq n \leq 0} \int_{\bbR^2_\eta} \int_{\bbR^2_{\xi_1}} \int_{\bbR^2_{\xi_2}} e^{ix\cdot(\xi_1+\xi_2+\eta)} \frakm_{n_1,n_2,n}(\xi_1,\xi_2,\eta) \\
        &\qquad \quad \times \jxione^N e^{\kappa_1 is \jxione} \widehat{f}^{\kappa_1}_{k,\real}(s,\xi_1) \jxitwo^4 e^{\kappa_2 is \jxitwo} \widehat{f}^{\kappa_2}_{l,\real}(s,\xi_2) \varphi_{\leq 5}(\eta) \, \ud \xi_1 \, \ud \xi_2 \, \ud \eta 
    \end{aligned}
\end{equation}
with
\begin{equation}
    \begin{aligned}
        \frakm_{n_1,n_2,n}(\xi_1,\xi_2,\eta) &:= \varphi_{n_1}(\xi_1) \varphi_{n_2}(\xi_2) \psi_n(\eta) \bigl( i (\xi_1+\xi_2) \bigr)^\gamma \jap{\xi_1}^{-N} \jxitwo^{-4} \\
        &\qquad \qquad \qquad \qquad \times \nabla_\eta^2 \biggl( \varphi_{>-M}\bigl( R(\xi_1,\xi_2) \eta \bigr) \varphi_{\leq 1}(\eta) \frac{\eta_1}{|\eta|^3} \biggr),
    \end{aligned}
\end{equation}
where without loss of generality we assumed $n_1 \geq n_2$. Recall that $|\gamma| \leq N-1$.
Under the frequency constraints $0 \leq n_2 \leq n_1$ and $-\min\{n_1,n_2\} \leq n \leq 0$, the symbol $\frakm_{n_1,n_2,n}(\xi_1,\xi_2,\eta)$ satisfies
\begin{equation}
    \bigl| \nabla_{\xi_1}^a \nabla_{\xi_2}^b \nabla_\eta^{c} \frakm_{n_1,n_2,n}(\xi_1,\xi_2,\eta) \bigr| \lesssim 2^{-an_1} 2^{-bn_2} 2^{-cn} 2^{-n_1} 2^{-4(n_2+n)},
\end{equation}
whence 
\begin{equation}
    \bigl\| \widehat{\calF}^{-1}\bigl[\frakm_{n_1,n_2,n}\bigr] \bigr\|_{L^1( (\bbR^2)^3 )} \lesssim 2^{-n_1} 2^{-4(n_2+n)}.
\end{equation}
Thus, using the trilinear estimates from Lemma~\ref{lem:trilinear_estimates}, we obtain
\begin{equation}
    \begin{aligned}
        &\bigl\| \frakQ^{\gamma,2}_{hl}\bigl[ f^{\kappa_1}_{k,\real}(s), f^{\kappa_2}_{l,\real}(s) \bigr](x) \bigr\|_{L^2_x} \\
        &\lesssim \sum_{0 \leq n_2 \leq n_1} \sum_{-n_2 \leq n \leq 0} 2^{-n_1} 2^{-4(n_2+n)} \bigl\| \jD^N \bmf_\real(s) \bigr\|_{L^2_x} \bigl\| \jD^4 e^{is\jD} \bmf_\real(s) \bigr\|_{L^\infty_x} \bigr\| \widehat{\calF}^{-1}\bigl[ \varphi_{\leq 5} \bigr] \bigr\|_{L^\infty_x} \\
        &\lesssim 
        \bigl\|\bmf_\real(s) \bigr\|_{H^N_x} 
        \bigl\|  \jD^4 e^{is\jD} \bmf_\real(s) \bigr\|_{L^\infty_x}.
    \end{aligned}
\end{equation}
Returning to the high Sobolev norm estimate in \eqref{equ:high_sobolev_calQ1_regular_part}, invoking \eqref{equ:GNSell}
for long times $\varepsilon^{-2} \leq t \leq T$ we obtain the bound
\begin{equation}
    \begin{aligned}
        \bigl\| \widehat{\calQ}^{(1),\gamma}_{k,l,\mathrm{reg}}(t,\xi) \bigr\|_{L^2_\xi} 
        &\lesssim \bigl\| \varepsilon 
        \cdot s^{-1+\frac32\delta} \bigr\|_{L^2_s([t,T])} 
        \lesssim \varepsilon^{2-4\delta},
    \end{aligned}
\end{equation}
while for short times $0 \leq t \leq \varepsilon^{-2}$ we find 
\begin{equation}
    \begin{aligned}
        \bigl\| \widehat{\calQ}^{(1),\gamma}_{k,l,\mathrm{reg}}(t,\xi) \bigr\|_{L^2_\xi} 
        & \lesssim \Bigl\| \varepsilon \cdot
        \varepsilon^{1-10\delta} \js^{-\frac12-\frac72\delta} 
        \Bigr\|_{L^2_s([t,\varepsilon^{-2}])} 
        + \varepsilon^{2-4\delta} 
        \lesssim \varepsilon^{2-10\delta}.
    \end{aligned}
\end{equation}

\medskip
\noindent \underline{High Sobolev norm estimate for the singular part $\widehat{\calQ}^{(1), \gamma}_{k,l,\mathrm{sing}}(t,\xi)$.}
Here we first integrate by parts in time. This gives
\begin{equation} \label{equ:high_sobolev_calQ1_rewritten_IBP}
    \widehat{\calQ}^{(1),\gamma}_{k,l,\mathrm{sing}}(t,\xi) 
    = \sum_{\kappa_1, \kappa_2 \in \{\pm\}} \Bigl[ \widehat{\calB}^{\pvdots,1,\gamma}_{\kappa_1\kappa_2}\bigl[ f_{k,\real}, f_{l,\real} \bigr](s, \xi) \Bigr]_{s=t}^{s=T} + \widehat{\calC}^{(1), \gamma}_{k,l,\mathrm{sing}}(t,\xi),
\end{equation}
where the boundary terms are
\begin{equation} \label{equ:high_sobolev_calQ1_rewritten_boundary_terms}
    \begin{aligned}
        &\widehat{\calB}^{\pvdots,1,\gamma}_{\kappa_1\kappa_2}\bigl[ f_{k,\real}, f_{l,\real} \bigr](s,\xi) \\
        &:= \int_{\bbR^2} \int_{\bbR^2} \frac{1}{i\Psi_{\kappa_1\kappa_2}(\xi,\xi_1,\xi_2)} e^{is \Psi_{\kappa_1\kappa_2}(\xi,\xi_1,\xi_2)} \bigl( i (\xi_1+\xi_2) \bigr)^\gamma \widehat{f}^{\kappa_1}_{k,\real}(s,\xi_1) \widehat{f}^{\kappa_2}_{l,\real}(s,\xi_2)  \\
        &\qquad \qquad \times \varphi_{\leq -M}\bigl( R(\xi_1,\xi_2) (\xi-\xi_1-\xi_2) \bigr) \varphi_{\leq 1}(\xi-\xi_1-\xi_2) \, \pvdots \, \frac{(\xi-\xi_1-\xi_2)_1}{|\xi-\xi_1-\xi_2|^3} \, \ud \xi_1 \, \ud \xi_2 
    \end{aligned}
\end{equation}
with phase function $\Psi_{\kappa_1\kappa_2}(\xi,\xi_1,\xi_2)$ defined in \eqref{equ:definition_Psi_quadratic}, and where the bulk terms can be written as
\begin{equation} \label{equ:high_sobolev_calQ1_rewritten_bulk_terms}
    \widehat{\calC}^{(1), \gamma}_{k,l,\mathrm{sing}}(t,\xi) := - \sum_{\kappa_1, \kappa_2 \in \{\pm\}} \int_t^T \widehat{\calB}^{\pvdots,1,\gamma}_{\kappa_1\kappa_2}\bigl[ \partial_s f_{k,\real}, f_{l,\real} \bigr](s,\xi) \, \ud s + \bigl\{ \text{similar term} \bigr\}.
\end{equation}
Here ``similar term'' is the one where $\partial_s$ hits the other profile $f_{l,\real}(s)$, which can be handled analogously.

We first estimate the boundary terms \eqref{equ:high_sobolev_calQ1_rewritten_boundary_terms}. 
For $0 \leq t \leq s \leq T$ we have 
\begin{equation}
    \begin{aligned}
        &\bigl\| \widehat{\calB}^{\pvdots,1,\gamma}_{\kappa_1\kappa_2}\bigl[ f_{k,\real}, f_{l,\real} \bigr](s,\xi) \bigr\|_{L^2_\xi} \\
        &\lesssim \sum_{0 \leq n_2 \leq n_1} \, \biggl\| \int_{\bbR^2} \int_{\bbR^2} \int_{\bbR^2} e^{ix\cdot(\xi_1+\xi_2+\sigma)} \frakm_{n_1,n_2}(\xi_1,\xi_2,\sigma) \varphi_{n_1}(\xi_1) e^{is\kappa_1 \jap{\xi_1}} \widehat{f}^{\kappa_1}_{k,\real}(s,\xi_1) \\
        &\qquad \qquad \qquad \qquad \qquad \qquad \qquad \qquad \times \varphi_{n_2}(\xi_2) e^{is\kappa_2\jap{\xi_2}} \widehat{f}^{\kappa_2}_{l,\real}(s,\xi_2) \, \pvdots \, \frac{(\sigma)_1}{|\sigma|^3} \, \ud \sigma \, \ud \xi_1 \, \ud \xi_2 \biggr\|_{L^2_x} \\
        &\quad + \bigl\{ \text{similar terms} \bigr\},
    \end{aligned}
\end{equation}
where (we disregard $\gamma$ in the notation)
\begin{equation} \label{equ:high_sobolev_calQ1_rewritten_symbolm}
    \begin{aligned}
        &\frakm_{n_1,n_2}(\xi_1,\xi_2,\sigma) \\ 
        &\quad := \frac{1}{\Psi_{\kappa_1\kappa_2}(\xi_1+\xi_2+\sigma,\xi_1,\xi_2)} \bigl( i (\xi_1+\xi_2) \bigr)^\gamma \widetilde{\varphi}_{n_1}(\xi_1) \widetilde{\varphi}_{n_2}(\xi_2) \varphi_{\leq -M}\bigl( R(\xi_1,\xi_2) \sigma \bigr) \varphi_{\leq 1}(\sigma).
    \end{aligned}
\end{equation}
Using Lemma~\ref{lem:derivative_bounds_quadratic_phase}, we have for integers $0 \leq a,b,c \leq 3$,
\begin{equation}
    \begin{aligned}
        \bigl| \nabla_{\xi_1}^a \nabla_{\xi_2}^b \nabla_\sigma^c \frakm_{n_1,n_2}(\xi_1,\xi_2,\sigma) \bigr| \lesssim 2^{-a n_1} 2^{-b n_2} 2^{|\gamma| \max\{n_1, n_2\}} 2^{16 \min\{n_1, n_2\}}, \quad |\gamma| \leq N-1,
    \end{aligned}
\end{equation}
whence 
\begin{equation} \label{equ:calQ1_high_Sobolev_rewritten_symbolm_bound}
    \begin{aligned}
        \bigl\| \widehat{\calF}^{-1}\bigl[ \frakm_{n_1,n_2} \bigr] \bigr\|_{L^1( (\bbR^2)^3 )} \lesssim 2^{|\gamma| \max\{n_1, n_2\}} 2^{16 \min\{n_1, n_2\}}, \quad |\gamma| \leq N-1.
    \end{aligned}
\end{equation}
By the bilinear estimates from Lemma~\ref{lem:bilinear_estimates} we conclude 
\begin{equation} \label{equ:calQ1_high_Sobolev_rewritten_boundary1}
    \begin{aligned}
        &\bigl\| \widehat{\calB}^{\pvdots,1,\gamma}_{\kappa_1\kappa_2}\bigl[ f_{k,\real}, f_{l,\real} \bigr](s,\xi) \bigr\|_{L^2_\xi} \\
        &\lesssim \sum_{0 \leq n_2 \leq n_1} \, 2^{(N-1) n_1} 2^{16 n_2} \bigl\| P_{n_1} \bmf_{\real}(s) \bigr\|_{L^2_x} \bigl\| P_{n_2} e^{is\jD} \bmf_\real(s) \bigr\|_{L^\infty_x} + \bigl\{ \text{similar terms} \bigr\}, \\
        &\lesssim \bigl\| \bmf_{\real}(s) \bigr\|_{H^N_x} \bigl\| \jD^{16} e^{is\jD} \bmf_\real(s) \bigr\|_{L^\infty_x} + \bigl\{ \text{similar terms} \bigr\}, \\
    \end{aligned}
\end{equation}
Invoking \eqref{equ:bootstrap_assumption_HN_f} and \eqref{equ:GNSell}, it follows for long times $\varepsilon^{-2} \leq s \leq T$ that
\begin{equation} \label{equ:calQ1_high_Sobolev_rewritten_boundary2}
    \begin{aligned}
        \bigl\| \widehat{\calB}^{\pvdots,1,\gamma}_{\kappa_1\kappa_2}\bigl[ f_{k,\real}, f_{l,\real} \bigr](s,\xi) \bigr\|_{L^2_\xi} \lesssim \varepsilon \cdot t^{-1+2\delta} \lesssim \varepsilon^{3-4\delta},
    \end{aligned}
\end{equation}
while for short times $0 \leq t \leq s \leq \varepsilon^{-2}$ we find that
\begin{equation} \label{equ:calQ1_high_Sobolev_rewritten_boundary3}
    \begin{aligned}
        \bigl\| \widehat{\calB}^{\pvdots,1,\gamma}_{\kappa_1\kappa_2}\bigl[ f_{k,\real}, f_{l,\real} \bigr](s,\xi) \bigr\|_{L^2_\xi} \lesssim \varepsilon \cdot \varepsilon^{1-10\delta} \js^{-\frac12-3\delta} \lesssim \varepsilon^{2-10\delta}.
    \end{aligned}
\end{equation}

We can now proceed to estimate the bulk term in \eqref{equ:high_sobolev_calQ1_rewritten_bulk_terms}.
Using the splitting \eqref{equ:splitting_pt_bmf} from Lemma~\ref{lem:splitting_pt_bmf}, we see that the desired high Sobolev norm bound is a consequence of the following two estimates for $0 \leq t \leq T$,
\begin{align}
    \biggl\| \int_t^T \widehat{\calB}^{\pvdots,1,\gamma}_{\kappa_1\kappa_2}\Bigl[ \bigl( e^{-is\jD} \calQ^{\mathrm{l}}_\real \bigr)_k, f_{l,\real} \Bigr](s,\xi) \, \ud s \biggr\|_{L^2_\xi} &\lesssim \varepsilon^{\frac32-2\delta}, \label{equ:calQ1_high_Sobolev_rewritten_bulk_localized_contr} \\
    \biggl\| \int_t^T \widehat{\calB}^{\pvdots,1,\gamma}_{\kappa_1\kappa_2}\Bigl[ \bigl( e^{-is\jD} \calQ^{\mathrm{nl}}_\real \bigr)_k, f_{l,\real} \Bigr](s,\xi) \, \ud s \biggr\|_{L^2_\xi} &\lesssim \varepsilon^{2-5\delta}. \label{equ:calQ1_high_Sobolev_rewritten_bulk_nonlocalized_contr} 
\end{align}

{\it Proof of \eqref{equ:calQ1_high_Sobolev_rewritten_bulk_localized_contr}.}
Using the dual ILED estimate from Lemma~\ref{lem:flat_dual_ILED} and adopting the same notation for the symbol in \eqref{equ:high_sobolev_calQ1_rewritten_symbolm}, we estimate
\begin{equation}
    \begin{aligned}
        &\biggl\| \int_t^T \widehat{\calB}^{\pvdots,1,\gamma}_{\kappa_1\kappa_2}\Bigl[ \bigl( e^{-is\jD} \calQ^{\mathrm{l}}_\real \bigr)_k, f_{l,\real} \Bigr](s,\xi) \, \ud s \biggr\|_{L^2_\xi} \\
        &\lesssim \sum_{n_1, n_2 \geq 0} \, \biggl\| \jx^2 \int_{\bbR^2} \int_{\bbR^2} \int_{\bbR^2} e^{ix\cdot(\xi_1+\xi_2+\sigma)} \frakm_{n_1,n_2}(\xi_1,\xi_2,\sigma) \varphi_{n_1}(\xi_1) \widehat{\calF}\bigl[ \bigl( \calQ^{\mathrm{l}}_\real(s) \bigr)_k \bigr]^{\kappa_1}(\xi_1) \\
        &\qquad \qquad \qquad \qquad \qquad \qquad \qquad \qquad \times \varphi_{n_2}(\xi_2) e^{is\kappa_2\jap{\xi_2}} \widehat{f}^{\kappa_2}_{l,\real}(s,\xi_2) \, \pvdots \, \frac{(\sigma)_1}{|\sigma|^3} \, \ud \sigma \, \ud \xi_1 \, \ud \xi_2 \biggr\|_{L^2_x} \\
        &\lesssim \sum_{n_1, n_2 \geq 0} \, \biggl\| \int_{\bbR^2} \int_{\bbR^2} \int_{\bbR^2} e^{ix\cdot(\xi_1+\xi_2+\sigma)} \frakm_{n_1,n_2}(\xi_1,\xi_2,\sigma) \varphi_{n_1}(\xi_1) \nabla_{\xi_1}^2 \widehat{\calF}\bigl[ \bigl( \calQ^{\mathrm{l}}_\real(s) \bigr)_k \bigr]^{\kappa_1}(\xi_1) \\
        &\qquad \qquad \qquad \qquad \qquad \times \varphi_{n_2}(\xi_2) e^{is\kappa_2\jap{\xi_2}} \widehat{f}^{\kappa_2}_{l,\real}(s,\xi_2) \, \pvdots \, \frac{(\sigma)_1}{|\sigma|^3} \, \ud \sigma \, \ud \xi_1 \, \ud \xi_2 \biggr\|_{L^2_x} + \bigl\{\text{similar terms}\bigr\} \\ 
        &=: \sum_{n_1, n_2 \geq 0} C_{n_1,n_2}^{\mathrm{l},\gamma} + \bigl\{\text{similar terms}\bigr\}.
    \end{aligned}
\end{equation}
When $0 \leq n_2 \leq n_1$, using the bilinear estimate from Lemma~\ref{lem:bilinear_estimates} together with the symbol bound \eqref{equ:calQ1_high_Sobolev_rewritten_symbolm_bound}, the high Sobolev norm bound \eqref{equ:calQl_HN_bound} from Lemma~\ref{lem:splitting_pt_bmf} and the Gagliardo-Nirenberg-Sobolev interpolation estimate \eqref{equ:GNSell}, we obtain
\begin{equation}
    \begin{aligned}
        \sum_{0 \leq n_2 \leq n_1} C_{n_1,n_2}^{\mathrm{l},\gamma} &\lesssim \sum_{0 \leq n_2 \leq n_1} 2^{(N-1) n_1} 2^{16n_2} \Bigl\| \bigl\| P_{n_1} \jx^2 \calQ^{\mathrm{l}}_\real(s) \bigr\|_{L^2_x} \bigl\| P_{n_2} e^{is\jD} \bmf_\real(s) \bigr\|_{L^\infty_x} \Bigr\|_{L^2_s([t,T])} \\ 
        &\lesssim \Bigl\| \bigl\| \jx^2 \calQ^{\mathrm{l}}_\real(s) \bigr\|_{H^N_x} \bigl\| \jD^{16} e^{is\jD} \bmf_\real(s) \bigr\|_{L^\infty_x} \Bigr\|_{L^2_s([t,T])} \lesssim \varepsilon^{2-10\delta}.
    \end{aligned}
\end{equation}
When $0 \leq n_1 \leq n_2$, using again Lemma~\ref{lem:bilinear_estimates} together with \eqref{equ:calQ1_high_Sobolev_rewritten_symbolm_bound}, interpolating between the high Sobolev norm bound \eqref{equ:calQl_HN_bound} and the $L^\infty_x$-decay estimate \eqref{equ:calQl_Linfty_bound}, and invoking the bootstrap assumption \eqref{equ:bootstrap_assumption_HN_f}, we find 
\begin{equation}
    \begin{aligned}
        \sum_{0 \leq n_1 \leq n_2} C_{n_1,n_2}^{\mathrm{l},\gamma} &\lesssim \sum_{0 \leq n_1 \leq n_2} 2^{(N-1) n_2} 2^{16n_1} \Bigl\| \bigl\| P_{n_1} \jx^2 \calQ^{\mathrm{l}}_\real(s) \bigr\|_{L^\infty_x} \bigl\| P_{n_2} e^{is\jD} \bmf_\real(s) \bigr\|_{L^2_x} \Bigr\|_{L^2_s([t,T])} \\ 
        &\lesssim \Bigl\| \bigl\| \jD^{16} \jx^2 \calQ^{\mathrm{l}}_\real(s) \bigr\|_{L^\infty_x} \bigl\| \bmf_\real(s) \bigr\|_{H^N_x} \Bigr\|_{L^2_s([t,T])} \\ 
        &\lesssim \Bigl\| \bigl\| \jx^2 \calQ^{\mathrm{l}}_\real(s) \bigr\|_{L^\infty_x}^{1-\frac{16}{N-1}} \bigl\| \jD^{N} \jx^2 \calQ^{\mathrm{l}}_\real(s) \bigr\|_{L^2_x}^{\frac{16}{N-1}} \cdot \varepsilon \Bigr\|_{L^2_s([t,T])} \lesssim \varepsilon^{\frac32-2\delta},        
    \end{aligned}
\end{equation}
where the last inequality above follows by distinguishing the cases $0 \leq t \leq \varepsilon^{-1}$ and $\varepsilon^{-1} \leq t \leq T$. 

{\it Proof of \eqref{equ:calQ1_high_Sobolev_rewritten_bulk_nonlocalized_contr}.}
Using the bilinear estimates from Lemma~\ref{lem:bilinear_estimates} together with the symbol bound \eqref{equ:calQ1_high_Sobolev_rewritten_symbolm_bound}, we have for any $2 \leq p,q \leq \infty$ with $\frac{1}{p} + \frac{1}{q} = \frac12$ that
\begin{equation} \label{equ:calQ1_high_Sobolev_rewritten_calQnl_contribution}
    \begin{aligned}
        &\biggl\| \widehat{\calB}^{\pvdots,1,\gamma}_{\kappa_1\kappa_2}\Bigl[ \bigl( e^{-is\jD} \calQ^{\mathrm{nl}}_\real \bigr)_k, f_{l,\real} \Bigr](s,\xi) \, \ud s \biggr\|_{L^2_\xi} \\ 
        &\lesssim \sum_{n_1, n_2 \geq 0} 2^{(N-1) \max\{n_1,n_2\}} 2^{16 \min\{n_1,n_2\}} \bigl\| P_{n_1} \calQ^{\mathrm{nl}}_\real(s) \bigr\|_{L^p_x} \bigl\| P_{n_2} e^{is\jD} \bmf_\real(s) \bigr\|_{L^q_x} =: \sum_{n_1, n_2 \geq 0} C_{n_1,n_2}^{\mathrm{nl},\gamma}. 
    \end{aligned}
\end{equation}
In the case $0 \leq n_2 \leq n_1$, this gives 
\begin{equation}
    \begin{aligned}
        \sum_{0 \leq n_2 \leq n_1} C_{n_1,n_2}^{\mathrm{nl},\gamma} &\lesssim \bigl\| \calQ^{\mathrm{nl}}_\real(s) \bigr\|_{H^N_x} \bigl\| \jD^{16} e^{is\jD} \bmf_\real(s) \bigr\|_{L^\infty_x}.
    \end{aligned}
\end{equation}
Using the $H^N_x$-bound \eqref{equ:calQnl_HN_bound} and the interpolation estimate \eqref{equ:GNSell}, one can verify that the integral over $[t,T]$ of the preceding line is bounded by $C \varepsilon^{2-4\delta}$ by looking separately at the cases $t \lesssim \varepsilon^{-1}$ and $t \gtrsim \varepsilon^{-1}$.

When $0 \leq n_1 \leq n_2$, using again \eqref{equ:calQ1_high_Sobolev_rewritten_calQnl_contribution} along with interpolation,  we find
\begin{equation}
    \begin{aligned}
        \sum_{0 \leq n_1 \leq n_2} C_{n_1,n_2}^{\mathrm{nl},\gamma} &\lesssim \bigl\| \jD^{16} \calQ^{\mathrm{nl}}_\real(s) \bigr\|_{L^\infty_x} \bigl\| \bmf_\real(s) \bigr\|_{H^N_x} \lesssim \bigl\| \calQ^{\mathrm{nl}}_\real(s) \bigr\|_{L^\infty_x}^{1-\frac{16}{N-1}} \bigl\| \calQ^{\mathrm{nl}}_\real(s) \bigr\|_{H^N_x}^{\frac{16}{N-1}} \bigl\| \bmf_\real(s) \bigr\|_{H^N_x}.
    \end{aligned}
\end{equation}
Integrating the preceding line over $[t,T]$, and inserting the bounds \eqref{equ:calQnl_HN_bound}, \eqref{equ:calQnl_Linfty_bound} together with the bootstrap assumption \eqref{equ:bootstrap_assumption_HN_f}, we obtain an acceptable upper bound $C \varepsilon^{2-5\delta}$.

\medskip
\paragraph{{\it ILED estimate at high Sobolev norm for $\widehat{\calQ}^{(1)}_{k,l}(t,\xi)$}} \label{subsubsec:calQ1_ILED}

Just like for the high Sobolev norm bound considered in the preceding Subsection~\ref{subsubsec:calQ1_high_Sobolev}, for the proof of the ILED bound \eqref{equ:HN_ILED_g} we also have to take into account the limited decay~\eqref{equ:FTvortex_re_im_remainder} on the frequency side of the remainder term in the decomposition~\eqref{equ:FTvortex_re_im} of the flat Fourier transform of the components of the vortex $\underline{\Phi}$.
For this reason we again distinguish the case when at least two derivatives fall onto the vortex, producing a sufficiently spatially localized term, and the case when all (apart from at most one) derivatives fall onto the inputs. Fix $0 < \kappa \ll 1$. Correspondingly, we write
\begin{equation} \label{equ:contribution1_ILED_decomposition}
    \begin{aligned}
        &\sum_{1 \leq k \leq N} \Bigl\| \jx^{-\frac12-\kappa} |D|^k e^{it\jD} \calQ_{k,l}^{(1)}(t) \Bigr\|_{L^2_t([0,T]; L^2_x(\bbR^2))} \\  
        &\lesssim \sum_{0 \leq |\beta| \leq 1} \sum_{|\gamma| \leq N-1-|\beta|} \biggl\| \jx^{-\frac12-\kappa} \int_t^T \frac{|D|}{\jD} e^{i(t-s)\jD} \Bigl( \bigl( \partial^\beta \underline{\Phi}_{\real} \bigr) \partial^\gamma \bigl( u_{k,\real}(s) u_{l,\real}(s) \bigr) \Bigr) \, \ud s \biggr\|_{L^2_t([0,T]; L^2_x)} \\
        &\quad + \sum_{2 \leq |\beta| \leq N-1} \sum_{|\gamma| \leq N-1-|\beta|} \biggl\| \jx^{-\frac12-\kappa} \int_t^T \frac{|D|}{\jD} e^{i(t-s)\jD} \Bigl( \bigl( \partial^\beta \underline{\Phi}_{\real} \bigr) \partial^\gamma \bigl( u_{k,\real}(s) u_{l,\real}(s) \bigr) \Bigr) \, \ud s \biggr\|_{L^2_t([0,T]; L^2_x)},
    \end{aligned}
\end{equation}
where $\beta, \gamma \in \bbN_0^2$ are multi-indices and we recall that $\bigl\| \jx^{|\beta|} \partial^\beta \underline{\Phi}_{\real} \bigr\|_{L^\infty_x} \lesssim 1$ for any multi-index $|\beta| \geq 1$. Thanks to this spatial localization, using the ILED estimate from Lemma~\ref{lem:flat_ILED_high_Sobolev_norm}, the second term on the right-hand side of \eqref{equ:contribution1_ILED_decomposition} can be estimated by
\begin{equation}
    \begin{aligned}
        &\sum_{2 \leq |\beta| \leq N-1} \sum_{|\gamma| \leq N-1-|\beta|} \biggl\| \jx^{-\frac12-\kappa} \int_t^T \frac{|D|}{\jD} e^{i(t-s)\jD} \Bigl( \bigl( \partial^\beta \underline{\Phi}_{\real} \bigr) \partial^\gamma \bigl( u_{k,\real}(s) u_{l,\real}(s) \bigr) \Bigr) \, \ud s \biggr\|_{L^2_t([0,T]; L^2_x)} \\
        &\lesssim \sum_{2 \leq |\beta| \leq N-1} \sum_{|\gamma| \leq N-1-|\beta|} \Bigl\| \jx^{\frac12+\kappa} \bigl( \partial^\beta \underline{\Phi}_{\real} \bigr) \partial^\gamma \bigl( u_{k,\real}(s) u_{l,\real}(s) \bigr) \Bigr\|_{L^2_t([0,T];L^2_x)} \\
        &\lesssim \sum_{2 \leq |\beta| \leq N-1} \bigl\| \jx^{\frac12+\kappa} \bigl( \partial^\beta \underline{\Phi}_{\real} \bigr) \bigr\|_{L^\infty_x} \Bigl\| \bigl\| e^{is\jD} \bmf_\real(s) \bigr\|_{L^\infty_x} \|\bmf_\real(s)\|_{H^{N-3}_x} \Bigr\|_{L^2_s([0,T])}.
    \end{aligned}
\end{equation}
Invoking \eqref{equ:bootstrap_assumption_Linfty_f}, \eqref{equ:bootstrap_assumption_HN_f}, we obtain that the last line is bounded by $\varepsilon^{2-10\delta}$, just like for the second term on the right-hand side of \eqref{equ:contribution1_HN_decomposition} in the high Sobolev norm estimate for $\widehat{\calQ}^{(1)}_{k,l}(t,\xi)$, 

For the first term on the right-hand side of \eqref{equ:contribution1_ILED_decomposition}, we discuss the case $\beta = 0$ when all derivatives fall on the inputs, the other case $|\beta| = 1$ being analogous. 
Proceeding as in Subsection~\ref{subsubsec:calQ1_high_Sobolev} for the high Sobolev norm estimates, we first consider the Duhamel integrals on the Fourier side
\begin{equation} \label{equ:calQ1_ILED_pass_to_Fourier}
    \begin{aligned}
        &\sum_{0 \leq |\gamma| \leq N-1} \, \biggl\| \jx^{-\frac12-\kappa} \int_t^T \frac{|D|}{\jD} e^{i(t-s)\jD} \Bigl( \bigl( \partial^\beta \underline{\Phi}_{\real} \bigr) \partial^\gamma \bigl( u_{k,\real}(s) u_{l,\real}(s) \bigr) \Bigr) \, \ud s \biggr\|_{L^2_t([0,T]; L^2_x)} \\ 
        &= \sum_{0 \leq |\gamma| \leq N-1} \, \biggl\| \jx^{-\frac12-\kappa} \frac{|D|}{\jD} e^{it\jD} \widehat{\calF}^{-1}\Bigl[ \widehat{\calQ}^{(1),\gamma}_{k,l}(t,\xi) \Bigr](x) \biggr\|_{L^2_t([0,T]; L^2_x)} 
    \end{aligned}
\end{equation}
with $\widehat{\calQ}^{(1),\gamma}_{k,l}(t,\xi)$ defined exactly as in \eqref{equ:calQ1gamma_def}. For the convenience of the reader, we recall that
\begin{equation} 
    \begin{aligned}
        \widehat{\calQ}^{(1), \gamma}_{k,l}(t,\xi) &:= \sum_{\kappa_1, \kappa_2 \in \{\pm\}} \int_t^T \int_{\bbR^2} \int_{\bbR^2} e^{is(-\jxi + \kappa_1 \jxione + \kappa_2 \jxitwo)} \bigl( i (\xi_1+\xi_2) \bigr)^\gamma  \\
        &\qquad \qquad \qquad \qquad \times \widehat{f}^{\kappa_1}_{k,\real}(s,\xi_1) \widehat{f}^{\kappa_2}_{l,\real}(s,\xi_2) \widehat{\underline{\Phi}}_\real(\xi-\xi_1-\xi_2) \, \ud \xi_1 \, \ud \xi_2 \, \ud s.
    \end{aligned}
\end{equation}
Next, we insert the decomposition of the flat Fourier transform of $\underline{\Phi}_\real$ from Lemma~\ref{lem:FTvortex_re_im} to obtain the same decomposition \eqref{equ:calQ1gamma_decomposition} of $\widehat{\calQ}^{(1),\gamma}_{k,l}$ into a singular part, a regular part, and a remainder term. 
Both the regular part and the remainder term are spatially localized. Thus, using the ILED estimate from Lemma~\ref{lem:flat_ILED_high_Sobolev_norm} and placing the output in the weighted $L^2_t([0,T]; L^2_x(\bbR^2))$ space, we may then estimate their contributions to the right-hand side of \eqref{equ:calQ1_ILED_pass_to_Fourier} exactly as their corresponding contributions to the high Sobolev norm estimate in the preceding Subsection~\ref{subsubsec:calQ1_high_Sobolev}. In the expression for the singular part $\widehat{\calQ}_{k,l,\mathrm{sing}}^{(1),\gamma}(t,\xi)$, we first integrate by parts in time exactly as in \eqref{equ:high_sobolev_calQ1_rewritten_IBP}.  
In the bulk terms on the right-hand side of \eqref{equ:high_sobolev_calQ1_rewritten_IBP} we then have to insert the decomposition \eqref{equ:splitting_pt_bmf} of $\partial_s \widehat{f}^{\kappa_1}_{k,\real}(s,\xi_1)$ from Lemma~\ref{lem:splitting_pt_bmf} again. The resulting contributions to \eqref{equ:calQ1_ILED_pass_to_Fourier} can then be estimated exactly as their corresponding contributions to the high Sobolev norm bound from the preceding Subsection~\ref{subsubsec:calQ1_high_Sobolev}, except that one first invokes the ILED estimate from Lemma~\ref{lem:flat_ILED_high_Sobolev_norm}, placing non-spatially localized outputs in $L^1_t L^2_x([0,T]; L^2_x(\bbR^2))$ and spatially localized outputs in the weighted $L^2_t([0,T]; L^2_x(\bbR^2))$ norm.
It remains to discuss the boundary terms on the right-hand side of \eqref{equ:high_sobolev_calQ1_rewritten_IBP}. Their contributions to \eqref{equ:calQ1_ILED_pass_to_Fourier} can be estimated directly as follows
\begin{equation}
    \begin{aligned}
        &\sum_{0 \leq |\gamma| \leq N-1} \, \Biggl\| \jx^{-\frac12-\kappa} \frac{|D|}{\jD} e^{it\jD} \widehat{\calF}^{-1}\biggl[ \sum_{\kappa_1, \kappa_2 \in \{\pm\}} \Bigl[ \widehat{\calB}^{\pvdots,1,\gamma}_{\kappa_1\kappa_2}\bigl[ f_{k,\real}, f_{l,\real} \bigr](s, \xi) \Bigr]_{s=t}^{s=T} \biggr] \Biggr\|_{L^2_t([0,T]; L^2_x)} \\
        &\lesssim \sum_{0 \leq |\gamma| \leq N-1} \sum_{\kappa_1, \kappa_2 \in \{\pm\}} \, \biggl( \Bigl\| \bigl\| \widehat{\calB}^{\pvdots,1,\gamma}_{\kappa_1\kappa_2}\bigl[ f_{k,\real}, f_{l,\real} \bigr](T, \xi) \bigr\|_{L^2_\xi} \Bigr\|_{L^2_t([0,T])} + \Bigl\| \bigl\| \widehat{\calB}^{\pvdots,1,\gamma}_{\kappa_1\kappa_2}\bigl[ f_{k,\real}, f_{l,\real} \bigr](t, \xi) \bigr\|_{L^2_\xi} \Bigr\|_{L^2_t([0,T])} \biggr),
    \end{aligned}
\end{equation}
where we just dropped the weight and the propagator $\partial \jD^{-1} e^{it\jD}$ in $L^2_x$.
Recall that we consider the more difficult case $T \geq \varepsilon^{-2}$. Using \eqref{equ:calQ1_high_Sobolev_rewritten_boundary1}, \eqref{equ:calQ1_high_Sobolev_rewritten_boundary2}, \eqref{equ:calQ1_high_Sobolev_rewritten_boundary3}, we then conclude that the last line is bounded by 
\begin{equation}
    \begin{aligned}
        \Bigl\| \varepsilon \cdot T^{-1+2\delta} \Bigr\|_{L^2_t([0,T])} + \Bigl\| \varepsilon \cdot \varepsilon^{1-10\delta} \jt^{-\frac12-3\delta} \Bigr\|_{L^2_t([0,\varepsilon^{-2}])} + \Bigl\| \varepsilon \cdot t^{-1+2\delta} \Bigr\|_{L^2_t([\varepsilon^{-2},T])} \lesssim \varepsilon^{2-10\delta}.
    \end{aligned}
\end{equation}

\subsubsection{Contributions of type (2)} \label{subsec:contributions_type2}

Next, we consider the contributions of type (2), which are quadratic terms of the form $u_{k,\ast_k}(t) \partial_j u_{l,\ast_l}(t)$, with $j \in \{1,2\}$, $\ast_k, \ast_l \in \{\real,\imag\}$, and $k, l \in \{1, 2, 3, 4\}$. By symmetry, it suffices to discuss the quadratic terms $u_{k,\real}(t) \partial_1 u_{l,\real}(t)$ with $k, l \in \{1, 2, 3, 4\}$. Their contribution to the profile $\widehat{\bmg}_\real^\flat(t,\xi)$ is given by
\begin{equation}
    \begin{aligned}
        \widehat{\calQ}^{(2)}_{k,l}(t,\xi) 
        &:= \int_t^T (2i\jxi)^{-1} e^{-is\jxi} \widehat{\calF}\bigl[ u_{k,\real}(s) \partial_1 u_{l,\real}(s) \bigr](\xi) \, \ud s \\
        &= 2\pi \sum_{\kappa_1, \kappa_2 \in \{\pm\}} \int_t^T \int_{\bbR^2} (2i\jxi)^{-1} e^{is\Psi_{\kappa_1 \kappa_2}(\xi, \xi_1, \xi-\xi_1)} \widehat{f}^{\kappa_1}_{k,\real}(s,\xi_1) \, i(\xi-\xi_1)_1 \widehat{f}^{\kappa_2}_{l,\real}(s,\xi-\xi_1) \, \ud \xi_1 \, \ud s,
    \end{aligned}
\end{equation}
where the phase $\Psi_{\kappa_1 \kappa_2}(\xi, \xi_1, \xi-\xi_1)$ is defined as in \eqref{equ:definition_Psi_quadratic}.
Integrating by parts in time gives 
\begin{equation} \label{equ:contribution_type2_calQ2_integrated_by_parts}
    \begin{aligned}
        &\widehat{\calQ}^{(2)}_{k,l}(t,\xi) \\
        &= \sum_{\kappa_1, \kappa_2 \in \{\pm\}} \Bigl[ \widehat{\calB}_{\kappa_1 \kappa_2}^{\delta_0,1}\bigl[ f_{k,\real}, f_{l,\real} \bigr](s,\xi) \Bigr]_{s=t}^{s=T} \\ 
        &\quad - 2\pi \sum_{\kappa_1, \kappa_2 \in \{\pm\}} \int_t^T \int_{\bbR^2} (2i\jxi)^{-1} \frac{1}{\Psi_{\kappa_1 \kappa_2}(\xi, \xi_1, \xi-\xi_1)} e^{is\Psi_{\kappa_1 \kappa_2}(\xi, \xi_1, \xi-\xi_1)} \\
        &\qquad \qquad \qquad \qquad \qquad \qquad \qquad \qquad \times \partial_s \widehat{f}^{\kappa_1}_{k,\real}(s,\xi_1) \, (\xi-\xi_1)_1 \widehat{f}^{\kappa_2}_{l,\real}(s,\xi-\xi_1) \, \ud \xi_1 \, \ud s \\
        &\quad + \bigl\{\text{similar terms}\bigr\}
    \end{aligned}
\end{equation}
with boundary terms
\begin{equation}
    \begin{aligned}
        \widehat{\calB}_{\kappa_1 \kappa_2}^{\delta_0,1}\bigl[ f_{k,\real}, f_{l,\real} \bigr](s,\xi) := 2\pi \int_{\bbR^2} (2i\jxi)^{-1} \frac{e^{is\Psi_{\kappa_1 \kappa_2}(\xi, \xi_1, \xi-\xi_1)}}{\Psi_{\kappa_1 \kappa_2}(\xi, \xi_1, \xi-\xi_1)} \widehat{f}^{\kappa_1}_{k,\real}(s,\xi_1) \, (\xi-\xi_1)_1 \widehat{f}^{\kappa_2}_{l,\real}(s,\xi-\xi_1) \, \ud \xi_1. 
    \end{aligned}
\end{equation}
The weighted energy estimates for the contribution of $\widehat{\calQ}^{(2)}_{k,l}(t,\xi)$ are now obtained by an analogous normal-form strategy as in the treatment of the singular contribution $\widehat{\calQ}^{(1)}_{k,l,\mathrm{sing}}(t,\xi)$ in the preceding Subsection~\ref{subsubsec:weighted_estimate_contribution_type1} concerning contributions of type (1). 
The only new issue is the derivative on one input: after dyadic localization, some fixed-output estimates are not directly summable over the output frequency. In those cases we retain Littlewood–Paley orthogonality and square-sum the output pieces.
After reinserting the evolution equations for the profiles, the normal-form remainders in particular give rise to non-localized cubic terms. For their weighted energy estimates we use the commutation identity \eqref{equ:commutation_identity_cubic_pv} with $\eta=0$.
The high Sobolev norm estimate and the ILED estimate for the contribution of $\widehat{\calQ}^{(2)}_{k,l}(t,\xi)$ are proved in an analogous manner as the corresponding bounds for the singular part $\widehat{\calQ}^{(1),\gamma}_{k,l,\mathrm{sing}}(t,\xi)$ in the preceding Subsections~\ref{subsubsec:calQ1_high_Sobolev} and \ref{subsubsec:calQ1_ILED}, up to the additional dyadic square summation via Littlewood-Paley trichotomies.
In what follows we illustrate this mechanism at the example of the derivation of the weighted energy estimates for the contribution of $\widehat{\calQ}^{(2)}_{k,l}(t,\xi)$, while we leave the analogous details for the high Sobolev norm estimate and for the ILED estimate to the reader.

\medskip 

We first discuss the weighted estimate for the boundary terms. By direct computation, we find that up to irrelevant constants
\begin{equation}
    \begin{aligned}
        &\jxi^2 \nabla_\xi \widehat{\calB}_{\kappa_1\kappa_2}^{\delta_0,1}\bigl[ f_{k,\real}, f_{l,\real} \bigr](s,\xi) \\
        &\simeq s \cdot \jxi \int_{\bbR^2} \frac{\nabla_{\xi} \Psi_{\kappa_1\kappa_2}(\xi,\xi_1,\xi-\xi_1)}{\Psi_{\kappa_1\kappa_2}(\xi,\xi_1,\xi-\xi_1)} e^{is \Psi_{\kappa_1\kappa_2}(\xi,\xi_1,\xi-\xi_1)} \widehat{f}^{\kappa_1}_{k,\real}(s,\xi_1) \, (\xi-\xi_1)_1 \widehat{f}^{\kappa_2}_{l,\real}(s,\xi-\xi_1) \, \ud \xi_1 \\
        &\quad + \jxi \int_{\bbR^2} \frac{1}{\Psi_{\kappa_1\kappa_2}(\xi,\xi_1,\xi-\xi_1)} e^{is \Psi_{\kappa_1\kappa_2}(\xi,\xi_1,\xi-\xi_1)} \widehat{f}^{\kappa_1}_{k,\real}(s,\xi_1) \, (\xi-\xi_1)_1 \nabla_\xi \widehat{f}^{\kappa_2}_{l,\real}(s,\xi-\xi_1) \, \ud \xi_1 \\
        &\quad + \bigl\{ \text{similar and easier terms} \bigr\} \\
        &=: I(s,\xi) + II(s,\xi) + \bigl\{ \text{similar and easier terms} \bigr\}.
    \end{aligned}
\end{equation}
The second term $II(s,\xi)$ is tightest in terms of dyadic frequency summations and so we explain how to estimate it. 
We have 
\begin{equation} \label{equ:type2_boundary_squaresum}
    \begin{aligned}
        \bigl\| II(s,\xi) \bigr\|_{L^2_\xi} \lesssim \biggl( \sum_{n \geq 0} \, \bigl\| II_n(s,\xi) \bigr\|_{L^2_\xi}^2 \biggr)^{\frac12}
    \end{aligned}
\end{equation}
with 
\begin{equation}
    \begin{aligned}
        II_n(s,\xi) &:= \varphi_n(\xi) \jxi \int_{\bbR^2} \frac{e^{is \Psi_{\kappa_1\kappa_2}(\xi,\xi_1,\xi-\xi_1)}}{\Psi_{\kappa_1\kappa_2}(\xi,\xi_1,\xi-\xi_1)} \widehat{f}^{\kappa_1}_{k,\real}(s,\xi_1) \, (\xi-\xi_1)_1 \nabla_\xi \widehat{f}^{\kappa_2}_{l,\real}(s,\xi-\xi_1) \, \ud \xi_1.
    \end{aligned}
\end{equation}
By a standard Littlewood-Paley trichotomy we decompose each term $II_n(s,\xi)$ into low-high, high-low, and high-high interactions. Thus, we write
\begin{equation}
    \begin{aligned}
        II_n(s,\xi) = II_n^{lh}(s,\xi) + II_n^{hl}(s,\xi) + II_n^{hh}(s,\xi)
    \end{aligned}
\end{equation}
with 
\begin{equation}
    \begin{aligned}
        II_n^{lh}(s,\xi) &:= \sum_{|n_2-n| \leq 5} \sum_{0 \leq n_1 \leq n_2-5} \varphi_n(\xi) \jxi \int_{\bbR^2} \frac{e^{is \Psi_{\kappa_1\kappa_2}(\xi,\xi_1,\xi-\xi_1)}}{\Psi_{\kappa_1\kappa_2}(\xi,\xi_1,\xi-\xi_1)} \varphi_{n_1}(\xi_1) \widehat{f}^{\kappa_1}_{k,\real}(s,\xi_1) \\
        &\qquad \qquad \qquad \qquad \qquad \qquad \qquad \quad \times \varphi_{n_2}(\xi-\xi_1) (\xi-\xi_1)_1 \nabla_\xi \widehat{f}^{\kappa_2}_{l,\real}(s,\xi-\xi_1) \, \ud \xi_1, \\
        II_n^{hl}(s,\xi) &:= \sum_{|n_1-n| \leq 5} \sum_{0 \leq n_2 \leq n_1-5} \varphi_n(\xi) \jxi \int_{\bbR^2} \frac{e^{is \Psi_{\kappa_1\kappa_2}(\xi,\xi_1,\xi-\xi_1)}}{\Psi_{\kappa_1\kappa_2}(\xi,\xi_1,\xi-\xi_1)} \varphi_{n_1}(\xi_1) \widehat{f}^{\kappa_1}_{k,\real}(s,\xi_1) \\
        &\qquad \qquad \qquad \qquad \qquad \qquad \qquad \quad \times \varphi_{n_2}(\xi-\xi_1) (\xi-\xi_1)_1 \nabla_\xi \widehat{f}^{\kappa_2}_{l,\real}(s,\xi-\xi_1) \, \ud \xi_1, \\
        II_n^{hh}(s,\xi) &:= \sum_{n_1 \geq n-5} \sum_{|n_2-n_1| \leq 5} \varphi_n(\xi) \jxi \int_{\bbR^2} \frac{e^{is \Psi_{\kappa_1\kappa_2}(\xi,\xi_1,\xi-\xi_1)}}{\Psi_{\kappa_1\kappa_2}(\xi,\xi_1,\xi-\xi_1)} \varphi_{n_1}(\xi_1) \widehat{f}^{\kappa_1}_{k,\real}(s,\xi_1) \\
        &\qquad \qquad \qquad \qquad \qquad \qquad \qquad \quad \times \varphi_{n_2}(\xi-\xi_1) (\xi-\xi_1)_1 \nabla_\xi \widehat{f}^{\kappa_2}_{l,\real}(s,\xi-\xi_1) \, \ud \xi_1.   
    \end{aligned}
\end{equation}
Then we bound the low-high interactions by
\begin{equation}
    \begin{aligned}
        \bigl\| II_n^{lh}(s,\xi) \bigr\|_{L^2_\xi} &\lesssim \sum_{|n_2-n| \leq 5} \sum_{0 \leq n_1 \leq n_2-5} 2^n \biggl\| \int_{\bbR^2} \int_{\bbR^2} e^{ix\cdot(\xi_1+\xi_2)} \frakm_{n_1,n_2}(\xi_1,\xi_2) \varphi_{n_1}(\xi_1) e^{is\kappa_1 \jap{\xi_1}} \widehat{f}^{\kappa_1}_{k,\real}(s,\xi_1) \\ 
        &\qquad \qquad \qquad \qquad \qquad \qquad \qquad \qquad \times \varphi_{n_2}(\xi_2) (\xi_2)_1 e^{is\kappa_2 \jap{\xi_2}} \nabla_{\xi_2} \widehat{f}^{\kappa_2}_{l,\real}(s,\xi_2) \, \ud \xi_1 \, \ud \xi_2 \biggr\|_{L^2_x} 
    \end{aligned}
\end{equation}
with 
\begin{equation}
    \frakm_{n_1,n_2}(\xi_1,\xi_2) := \frac{1}{\Psi_{\kappa_1\kappa_2}(\xi_1+\xi_2,\xi_1,\xi_2)} \widetilde{\varphi}_{n_1}(\xi_1) \widetilde{\varphi}_{n_2}(\xi_2).
\end{equation}
By Lemma~\ref{lem:derivative_bounds_quadratic_phase} we have for integers $0 \leq a, b \leq 3$,
\begin{equation}
    \bigl| \nabla_{\xi_1}^a \nabla_{\xi_2}^b \frakm_{n_1, n_2}(\xi_1,\xi_2) \bigr| \lesssim 2^{-a n_1} 2^{-b n_2} 2^{13 n_1},
\end{equation}
whence 
\begin{equation}
    \bigl\| \widehat{\calF}^{-1}\bigl[ \frakm_{n_1,n_2} \bigr] \bigr\|_{L^1( (\bbR^2)^2 )} \lesssim 2^{13 n_1}.
\end{equation}
The bilinear estimates from Lemma~\ref{lem:bilinear_estimates} now give
\begin{equation}
    \begin{aligned}
        \bigl\| II_n^{lh}(s,\xi) \bigr\|_{L^2_\xi} &\lesssim \sum_{|n_2-n| \leq 5} \sum_{0 \leq n_1 \leq n_2-5} 2^n \cdot 2^{13 n_1} \bigl\| P_{n_1} e^{is\jD} \bmf_\real(s) \bigr\|_{L^\infty_x} \cdot 2^{n_2} \bigl\| \varphi_{n_2}(\xi_2) \nabla_{\xi_2} \widehat{\bmf}_\real(s,\xi_2) \bigr\|_{L^2_{\xi_2}} \\
        &\lesssim \sum_{|n_2-n| \leq 5} \bigl\| \jD^{14} e^{is\jD} \bmf_\real(s) \bigr\|_{L^\infty_x} \cdot 2^{2n_2} \bigl\| \varphi_{n_2}(\xi_2)  \nabla_{\xi_2} \widehat{\bmf}_\real(s,\xi_2) \bigr\|_{L^2_{\xi_2}}.
    \end{aligned}
\end{equation}
Proceeding analogously, we obtain for the high-low interactions
\begin{equation}
    \begin{aligned}
        \bigl\| II_n^{hl}(s,\xi) \bigr\|_{L^2_\xi} &\lesssim \sum_{|n_1-n| \leq 5} \sum_{0 \leq n_2 \leq n_1-5} 2^n \cdot 2^{13 n_2} \bigl\| P_{n_1} e^{is\jD} \bmf_\real(s) \bigr\|_{L^\infty_x} \cdot 2^{n_2} \bigl\| \varphi_{n_2}(\xi_2) \nabla_{\xi_2} \widehat{\bmf}_\real(s,\xi_2) \bigr\|_{L^2_{\xi_2}} \\
        &\lesssim 2^{-n} \bigl\| \jD^{15} e^{is\jD} \bmf_\real(s) \bigr\|_{L^\infty_x} \bigl\| \jap{\xi_2}^2 \nabla_{\xi_2} \widehat{\bmf}_\real(s,\xi_2) \bigr\|_{L^2_{\xi_2}},
    \end{aligned}
\end{equation}
and for the high-high interactions
\begin{equation}
    \begin{aligned}
        \bigl\| II_n^{hh}(s,\xi) \bigr\|_{L^2_\xi} &\lesssim \sum_{n_1 \geq n-5} \sum_{|n_2-n_1| \leq 5} 2^n \cdot 2^{13 n_1} \bigl\| P_{n_1} e^{is\jD} \bmf_\real(s) \bigr\|_{L^\infty_x} \cdot 2^{n_2} \bigl\| \varphi_{n_2}(\xi_2) \nabla_{\xi_2} \widehat{\bmf}_\real(s,\xi_2) \bigr\|_{L^2_{\xi_2}} \\
        &\lesssim 2^{-n} \bigl\| \jD^{15} e^{is\jD} \bmf_\real(s) \bigr\|_{L^\infty_x} \bigl\| \jap{\xi_2}^2 \nabla_{\xi_2} \widehat{\bmf}_\real(s,\xi_2) \bigr\|_{L^2_{\xi_2}}.
    \end{aligned}
\end{equation}
Inserting the preceding bounds back into \eqref{equ:type2_boundary_squaresum}, square-summing, 
we arrive at the acceptable bound
\begin{equation}
    \begin{aligned}
        \bigl\| II(s,\xi) \bigr\|_{L^2_{\xi}} &\lesssim \bigl\| \jD^{15} e^{is\jD} \bmf_\real(s) \bigr\|_{L^\infty_x} \bigl\| \jxi^2 \nabla_\xi \widehat{\bmf}_\real(s,\xi) \bigr\|_{L^2_\xi}
        .
    \end{aligned}
\end{equation}

Next, we discuss the weighted energy estimate for the second main term on the right-hand side of \eqref{equ:contribution_type2_calQ2_integrated_by_parts}, where we have to insert the time derivative of the flat profile $\partial_s \bmf_\real(s)$ again, using the decomposition from Lemma~\ref{lem:splitting_pt_bmf}. 
Here we only consider the tightest case when a contribution of type (2) from $\calQ_\real^{\mathrm{nl}}(s)$ in \eqref{equ:contribution_type1_calQ1_sing_profile_equation_splitting_recalled} gets reinserted, which without loss of generality we may assume to be of the form $u_{k_1,\real} \partial_1 u_{k_2,\real}$ for some $k_1, k_2 \in \{1,2,3,4\}$. 
We are thus led to consider weighted bounds $\bigl\| \jap{\xi} \nabla_\xi \calI_2(t,\xi) \bigr\|_{L^2_\xi}$ for (slightly more general) terms of the form
\begin{equation}
    \begin{aligned}
        \calI_2(t,\xi) &:= \int_t^T \int_{\bbR^2} \int_{\bbR^2} e^{is\Psi_{\iota_1 \iota_2 \iota_3}(\xi,\xi_1,\xi_2,\xi-\xi_1-\xi_2)} \frac{1}{\Psi_{\kappa_1 \kappa_2}(\xi,\xi_1+\xi_2,\xi-\xi_1-\xi_2)} \jap{\xi_1+\xi_2}^{-1} \\
        &\qquad \times \widehat{f}^{\iota_1}_{k_1,\real}(s,\xi_1) \, (\xi_2)_1 \widehat{f}^{\iota_2}_{k_2,\real}(s,\xi_2) \, (\xi-\xi_1-\xi_2)_1 \widehat{f}^{\iota_3}_{k_3,\real}(s,\xi-\xi_1-\xi_2) \, \ud \xi_1 \, \ud \xi_2 \, \ud s 
    \end{aligned}
\end{equation}
for arbitrary $\iota_1, \iota_2, \iota_3 \in \{\pm\}$, $\kappa_1, \kappa_2 \in \{\pm\}$, $k_1, k_2, k_3 \in \{1,2,3,4\}$, and with the phase $\Psi_{\iota_1 \iota_2 \iota_3}$ defined as in \eqref{equ:definition_Psi_cubic}.
Using the commutation identity
\begin{equation} \label{equ:commutation_identity_cubic_delta}
    \begin{aligned}
        &\bigl( \jap{\xi} \nabla_\xi + \iota_1 \jap{\xi_1} \nabla_{\xi_1} + \iota_2 \jap{\xi_2} \nabla_{\xi_2} \bigr) \Psi_{\iota_1 \iota_2 \iota_3}(\xi,\xi_1,\xi_2,\xi-\xi_1-\xi_2) \\
        &\qquad \qquad \qquad \qquad \qquad \qquad \qquad \quad = - \iota_3 \frac{\xi-\xi_1-\xi_2}{\jap{\xi-\xi_1-\xi_2}} \Psi_{\iota_1 \iota_2 \iota_3}(\xi,\xi_1,\xi_2,\xi-\xi_1-\xi_2),
    \end{aligned}
\end{equation}
and integrating by parts in the frequency variables $\xi_1, \xi_2$, we find that
\begin{equation}
    \begin{aligned}
        &\jxi \nabla_\xi \calI_2(t,\xi) \\ 
        &= \int_t^T \int_{\bbR^2} \int_{\bbR^2} e^{is\Psi_{\iota_1 \iota_2 \iota_3}(\xi,\xi_1,\xi_2,\xi-\xi_1-\xi_2)} \frac{1}{\Psi_{\kappa_1 \kappa_2}(\xi,\xi_1+\xi_2,\xi-\xi_1-\xi_2)} \jap{\xi_1+\xi_2}^{-1} \\
        &\qquad \times \iota_1 \jap{\xi_1} \nabla_{\xi_1} \widehat{f}^{\iota_1}_{k_1,\real}(s,\xi_1) \, (\xi_2)_1 \widehat{f}^{\iota_2}_{k_2,\real}(s,\xi_2) \, (\xi-\xi_1-\xi_2)_1 \widehat{f}^{\iota_3}_{k_3,\real}(s,\xi-\xi_1-\xi_2) \, \ud \xi_1 \, \ud \xi_2 \, \ud s \\ 
        &\quad + \int_t^T \int_{\bbR^2} \int_{\bbR^2} e^{is\Psi_{\iota_1 \iota_2 \iota_3}(\xi,\xi_1,\xi_2,\xi-\xi_1-\xi_2)} \frac{1}{\Psi_{\kappa_1 \kappa_2}(\xi,\xi_1+\xi_2,\xi-\xi_1-\xi_2)} \jap{\xi_1+\xi_2}^{-1} \\
        &\qquad \times \widehat{f}^{\iota_1}_{k_1,\real}(s,\xi_1) \, \iota_2 \jap{\xi_2} (\xi_2)_1 \nabla_{\xi_2} \widehat{f}^{\iota_2}_{k_2,\real}(s,\xi_2) \, (\xi-\xi_1-\xi_2)_1 \widehat{f}^{\iota_3}_{k_3,\real}(s,\xi-\xi_1-\xi_2) \, \ud \xi_1 \, \ud \xi_2 \, \ud s \\ 
        &\quad + \int_t^T \int_{\bbR^2} \int_{\bbR^2} e^{is\Psi_{\iota_1 \iota_2 \iota_3}(\xi,\xi_1,\xi_2,\xi-\xi_1-\xi_2)} \frac{1}{\Psi_{\kappa_1 \kappa_2}(\xi,\xi_1+\xi_2,\xi-\xi_1-\xi_2)} \jap{\xi_1+\xi_2}^{-1} \\
        &\qquad \quad \times \widehat{f}^{\iota_1}_{k_1,\real}(s,\xi_1) \, (\xi_2)_1 \widehat{f}^{\iota_2}_{k_2,\real}(s,\xi_2) \, \bigl(\jxi - \iota_1 \jap{\xi_1} - \iota_2 \jap{\xi_2}\bigr) \\ 
        &\qquad \quad \times (\xi-\xi_1-\xi_2)_1 (\nabla \widehat{f}^{\iota_3}_{k_3,\real})(s,\xi-\xi_1-\xi_2) \, \ud \xi_1 \, \ud \xi_2 \, \ud s \\ 
        &\quad + \int_t^T \int_{\bbR^2} \int_{\bbR^2} e^{is\Psi_{\iota_1 \iota_2 \iota_3}(\xi,\xi_1,\xi_2,\xi-\xi_1-\xi_2)} \frac{1}{\Psi_{\kappa_1 \kappa_2}(\xi,\xi_1+\xi_2,\xi-\xi_1-\xi_2)} \jap{\xi_1+\xi_2}^{-1} \\
        &\qquad \quad \times i \cdot s \cdot (-\iota_3) \cdot \frac{\xi-\xi_1-\xi_2}{\jap{\xi-\xi_1-\xi_2}} \Psi_{\iota_1 \iota_2 \iota_3} \\
        &\qquad \quad \times \widehat{f}^{\iota_1}_{k_1,\real}(s,\xi_1) \, (\xi_2)_1 \widehat{f}^{\iota_2}_{k_2,\real}(s,\xi_2) \, (\xi-\xi_1-\xi_2)_1 \widehat{f}^{\iota_3}_{k_3,\real}(s,\xi-\xi_1-\xi_2) \, \ud \xi_1 \, \ud \xi_2 \, \ud s \\ 
        &\quad + \bigl\{ \text{similar or easier terms} \bigr\} \\
        &=: \calJ_1(t,\xi) + \calJ_2(t,\xi) + \calJ_3(t,\xi) + \calJ_4(t,\xi) + \bigl\{ \text{similar or easier terms} \bigr\}.
    \end{aligned}
\end{equation}
We now estimate the $L^2_\xi$-norm of the above terms by dyadically decomposing the input and output frequencies and using the trilinear estimates from Lemma~\ref{lem:trilinear_estimates}, similarly to the derivation of the weighted bound for the term $\calI_1(t,\xi)$ in Subsection~\ref{subsubsec:weighted_estimate_contribution_type1}. The terms $\calJ_1(t,\xi)$ and $\calJ_2(t,\xi)$ have enough room to directly sum over the dyadic frequency pieces in view of the inverse frequency weight $\jap{\xi_1+\xi_2}^{-1}$. In the term $\calJ_4(t,\xi)$ we integrate by parts in time again using the additional phase factor $\Psi_{\iota_1 \iota_2 \iota_3}$. Only the term $\calJ_3(t,\xi)$ is tight and requires square summation over the output frequencies inside the time integral via suitable Littlewood-Paley trichotomies, similarly to the treatment of the weighted energy estimate of the boundary term above. We omit the lengthy but standard details.

\subsubsection{Contributions of type (3)} \label{equ:subsec:contributions_type3}
We leave the details of treating the contributions of all cubic terms of the form $u_{k,\ast_k}(t) u_{l,\ast_l}(t) u_{m,\ast_m}$
with $\ast_k, \ast_l, \ast_m \in \{\real,\imag\}$, and $k, l, m \in \{1, 2, 3, 4\}$ to the reader. Similar cubic terms result from contributions of type (1) and of type (2) after a normal form transformation and have already been treated in detail before.

\subsubsection{Contributions of type (4)} \label{equ:subsec:contributions_type4}
Next, we turn to the contributions of all spatially localized quadratic terms of the form 
$(1-a_\theta) u_{k,\ast_k} \partial_j u_{l, \ast_l}$ with $j \in \{1,2\}$, $\ast_k, \ast_l \in \{\real, \imag\}$, $k, l \in \{1, 2, 3, 4 \}$.
It suffices to discuss the details for terms of the form $(1-a_\theta) u_{k,\real} \partial_1 u_{l, \real}$ with $k, l \in \{1, 2, 3, 4 \}$.
Their contribution to the profile $\widehat{\bmg}_\real^\flat(t,\xi)$ is given by
\begin{equation}
    \begin{aligned}
        \widehat{\calQ}^{(4)}_{k,l}(t,\xi) 
        &:= \int_t^T (2i\jxi)^{-1} e^{-is\jxi} \widehat{\calF}\bigl[ (1-a_\theta) u_{k,\real} \partial_1 u_{l, \real} \bigr](\xi) \, \ud s.
    \end{aligned}
\end{equation}

\medskip
\paragraph{{\it Weighted estimate for $\widehat{\calQ}^{(4)}_{k,l}(t,\xi)$}} 
Using the dual ILED estimate from Lemma~\ref{lem:flat_dual_ILED_low_freq}, we have 
\begin{equation} \label{equ:calQ4_weighted_estimate}
    \begin{aligned}
        &\bigl\| \jxi^2 \nabla_\xi \widehat{\calQ}^{(4)}_{k,l}(t,\xi) \bigr\|_{L^2_\xi} \\
        &\lesssim \biggl\| \int_t^T s \cdot \frac{\xi}{\jxi} \cdot e^{-is\jxi} \jxi \widehat{\calF}\bigl[ (1-a_\theta) u_{k,\real}(s) \partial_1 u_{l, \real}(s) \bigr](\xi) \, \ud s \biggr\|_{L^2_\xi} \\ 
        &\quad + \biggl\| \int_t^T e^{-is\jxi} \jxi \nabla_\xi \widehat{\calF}\bigl[ (1-a_\theta) u_{k,\real}(s) \partial_1 u_{l, \real}(s) \bigr](\xi) \, \ud s \biggr\|_{L^2_\xi} + \bigl\{\text{similar or easier terms}\bigr\} \\ 
        &\lesssim \biggl\| s \cdot \bigl\| \jx \jD \bigl( (1-a_\theta) u_{k,\real}(s) \partial_1 u_{l,\real}(s) \bigr) \bigr\|_{L^2_x} \biggr\|_{L^2_s([t,T])} \\ 
        &\quad + \int_t^T \bigl\| \jx \jD \bigl( (1-a_\theta) u_{k,\real}(s) \partial_1 u_{l,\real}(s) \bigr) \bigr\|_{L^2_x} \, \ud s + \bigl\{\text{similar or easier terms}\bigr\}.
    \end{aligned}
\end{equation}
Then for the first term on the right-hand side we 
have
\begin{equation}
    \begin{aligned}
        &\biggl\| s \cdot \bigl\| \jx \jD \bigl( (1-a_\theta) u_{k,\real}(s) \partial_1 u_{l,\real}(s) \bigr) \bigr\|_{L^2_x} \biggr\|_{L^2_s([t,T])} \\
        &\lesssim \biggl\| s \cdot \bigl\| \jx (1-a_\theta) \bigr\|_{H^1_x} \cdot \bigl\| \jD e^{is\jD} \bmf_\real(s) \bigr\|_{L^\infty_x} \cdot \bigl\| \jD^2 e^{is\jD} \bmf_\real(s) \bigr\|_{L^\infty_x} \biggr\|_{L^2_s([t,T])} 
    \end{aligned}
\end{equation}
Invoking \eqref{equ:GNSell} we find, as for similar terms before (see for example \eqref{eq:GNSprototype}) that for long times $\varepsilon^{-2} \leq t \leq T$ the last line is bounded by $t^{-\frac12+3\delta}$, while for short times $0 \leq t \leq \varepsilon^{-2}$ we have the estimate $\varepsilon^{1-5\delta}$.
The bounds for the other terms on the right-hand side of \eqref{equ:calQ4_weighted_estimate} are even more favorable and we omit the details.

\medskip
\paragraph{{\it High Sobolev norm estimate for $\widehat{\calQ}^{(4)}_{k,l}(t,\xi)$}}
Using the dual ILED estimate from Lemma~\ref{lem:flat_dual_ILED} with $0 < \kappa \ll 1$, 
we obtain
\begin{equation}
    \begin{aligned}
        \bigl\| \jxi^N \widehat{\calQ}_{k,l}^{(4)}(t,\xi) \bigr\|_{L^2_\xi} 
        &\lesssim \biggl\| \int_t^T e^{-is\jxi} \jxi^{N-1} \widehat{\calF}\bigl[ (1-a_\theta) u_{k,\real}(s) \partial_1 u_{l,\real}(s) \bigr](\xi) \, \ud s \biggr\|_{L^2_\xi} \\ 
        &\lesssim \Bigl\| \bigl\| \jx^{1+\kappa} \jD^{N-1} \bigl( (1-a_\theta) u_{k,\real}(s) \partial_1 u_{l,\real}(s) \bigr) \bigr\|_{L^2_x} \Bigr\|_{L^2_s([t,T])} \\
        &\lesssim \Bigl\| \bigl\| \jx^{1+\kappa} (1-a_\theta) \bigr\|_{W^{N-1,\infty}_x} \bigl\| \jD e^{is\jD} \bmf_\real(s) \bigr\|_{L^\infty_x} \bigl\|\bmf_\real(s)\bigr\|_{H^N_x} \Bigr\|_{L^2_s([0,T])} .
    \end{aligned}
\end{equation}
From \eqref{equ:GNSell} we infer that the preceding line is bounded by $\varepsilon^{2-10\delta}$.

\medskip
\paragraph{{\it ILED estimate at high Sobolev norm for $\widehat{\calQ}^{(4)}_{k,l}(t,\xi)$}}
Using the ILED estimate \eqref{equ:flat_ILED_high_Sobolev_final_time} from Lemma~\ref{lem:flat_ILED_high_Sobolev_norm}, we obtain for some $0 < \kappa \ll 1$,
\begin{equation}
    \begin{aligned}
        &\sum_{1 \leq k \leq N} \bigl\| \jx^{-\frac12-\kappa} |D|^k e^{it\jD} \calQ_{k,l}^{(4)}(t) \bigr\|_{L^2_t([0,T]; L^2_x)} \\
        &= \sum_{1 \leq k \leq N} \, \biggl\| \jx^{-\frac12-\kappa} \int_t^T (2i\jD)^{-1} |D|^k e^{i(t-s)\jD} \bigl( (1-a_\theta) u_{k,\real}(s) \partial_1 u_{l,\real}(s) \bigr) \, \ud s \biggr\|_{L^2_t([0,T]; L^2_x)} \\
        &\lesssim \sum_{1 \leq k \leq N} \, \Bigl\| \jx^{\frac12+\kappa} |D|^{k-1} \bigl( (1-a_\theta) u_{k,\real}(s) \partial_1 u_{l,\real}(s) \bigr) \Bigr\|_{L^2_t([0,T]; L^2_x)} \\
        &\lesssim \bigl\| \jx^{\frac12+\kappa} (1-a_\theta) \bigr\|_{W^{N-1,\infty}_x} \Bigl\| \bigl\| \jD e^{is\jD} \bmf_\real(s) \bigr\|_{L^\infty_x} \|\bmf_\real(s)\|_{H^N_x} \Bigr\|_{L^2_s([0,T])}.
    \end{aligned}
\end{equation}
By the same argument as in the preceding high Sobolev norm estimate for $\widehat{\calQ}^{(4)}_{k,l}(t,\xi)$, using 
\eqref{eq:GNSprototype}, the last line is bounded by $\varepsilon^{2-10\delta}$.

\subsubsection{Contributions of type (5)} \label{subsec:contributions_type5}

Now we turn to contributions of the form $\partial_j \partial_t \eta_{0,c}$ with $j \in \{1,2\}$.
Decomposing $\partial_t \eta_{0,c} = \partial_t \eta_{0,c}^\flat + \partial_t \eta_{0,c}^\sharp$ into its flat and sharp components, and inserting \eqref{equ:pteta0c_flat_equation_rewritten}, \eqref{equ:pteta0c_flat_equation_leading_order_term_rewritten}, 
the corresponding contribution to $\widehat{\bmg}_\real^\flat(t,\xi)$ is given by
\begin{equation} \label{equ:contribution_type5_decomposition}
    \begin{aligned}
        &\int_t^T (2i\jxi)^{-1} e^{-is\jxi} \widehat{\calF}\bigl[ \partial_j \partial_t \eta_{0,c}(s) \bigr](\xi) \, \ud s \\
        &\simeq \int_t^T \bigl( \jxi^{-1} (\xi)_j \bigr) e^{-is\jxi} \widehat{\calF}\Bigl[ (-\Delta + U^2)^{-1} \bigl( (1-U^2) \partial_t \eta_{0,c}^\flat(s) \bigr) \Bigr](\xi) \, \ud s \\ 
        &\quad + \int_t^T \bigl( \jxi^{-1} (\xi)_j \bigr) e^{-is\jxi} \jxi^{-2} \xi \cdot \widehat{\calF}\Bigl[ u_{1,\real} \nabla u_{3,\real} - u_{3,\real} \nabla u_{1,\real} +  u_{1,\imag} \nabla u_{3,\imag} - u_{3,\imag} \nabla u_{1,\imag} \Bigr](\xi) \, \ud s \\        
        &\quad + \int_t^T \bigl( \jxi^{-1} (\xi)_j \bigr) e^{-is\jxi} \jxi^{-2} \widehat{\calF}\Bigl[ U' \bigl( u_{2,\real} u_{3,\real} + u_{2,\imag} u_{3,\imag} + u_{1,\real} u_{4,\real} + u_{1,\imag} u_{4,\imag} \bigr) \Bigr](\xi) \, \ud s \\
        &\quad + \int_t^T \bigl( \jxi^{-1} (\xi)_j \bigr) e^{-is\jxi} \jxi^{-2} \widehat{\calF}\Bigl[ \bigl\{\text{cubic and higher-order nonlinearities}\bigr\} \Bigr](\xi) \, \ud s. 
    \end{aligned}
\end{equation}
We now discuss schematically the weighted energy estimate, the high Sobolev norm estimate as well as the ILED estimate for the terms on the right-hand side of \eqref{equ:contribution_type5_decomposition}.
The first term on the right-hand side is spatially localized, moreover Proposition~\ref{prop:eta0c_bounds} provides sufficiently fast $L^\infty_x$-decay for $\pt \eta_{0,c}^\flat(s)$, and even provides decay for the $H^N_x$-norm of $\pt \eta_{0,c}^\flat(s)$. The weighted energy estimate for the contribution of the first term can thus be closed easily using the dual ILED estimate from Lemma~\ref{lem:flat_dual_ILED_low_freq}. For the high Sobolev norm estimate we also use the dual ILED estimate from Lemma~\ref{lem:flat_dual_ILED_low_freq}, while for the ILED bound we use Lemma~\ref{lem:flat_ILED_high_Sobolev_norm}.
The second term on the right-hand side of \eqref{equ:contribution_type5_decomposition} is essentially of the same type as contributions of type (2), which have already been treated above. We point out that the derivative on one of the inputs gets compensated by the inverse frequency weight $\jxi^{-2} \xi$ on the output apart from the case of high-high interactions, but in that case one can easily transfer derivatives between the two inputs.
The third term on the right-hand side of \eqref{equ:contribution_type5_decomposition} is spatially localized and quadratic, whence essentially of the same type as contributions of type (4), which have been discussed in detail above. 
Finally, the fourth term on the right-hand side of \eqref{equ:contribution_type5_decomposition} involves only terms of at least cubic order. We skip the details for these tedious, but harmless terms.

\subsubsection{Contributions of type (6)} \label{subsec:contributions_type6}
We continue with contributions of the form $\underline{\Phi}_{\ast} \pt \eta_{0,c}$ with $\ast \in \{\real,\imag\}$.
Without loss of generality we consider $\underline{\Phi}_\real \pt\eta_{0,c}$.
As in the preceding case of contributions of type (5), we decompose $\partial_t \eta_{0,c} = \partial_t \eta_{0,c}^\flat + \partial_t \eta_{0,c}^\sharp$ into its flat and sharp components, and insert \eqref{equ:pteta0c_flat_equation_rewritten}, \eqref{equ:pteta0c_flat_equation_leading_order_term_rewritten}.
The resulting contribution to $\widehat{\bmg}_\real^\flat(t,\xi)$ is then given by
\begin{equation} \label{equ:contribution_type6_decomposition}
    \begin{aligned}
        &\int_t^T (2i\jxi)^{-1} e^{-is\jxi} \widehat{\calF}\bigl[ \underline{\Phi}_\real \partial_t \eta_{0,c}(s) \bigr](\xi) \, \ud s \\
        &\simeq \int_t^T \jxi^{-1} e^{-is\jxi} \widehat{\calF}\Bigl[ \underline{\Phi}_\real (-\Delta + U^2)^{-1} \bigl( (1-U^2) \partial_t \eta_{0,c}^\flat(s) \bigr) \Bigr](\xi) \, \ud s \\         
        &\quad + \int_t^T \jxi^{-1} e^{-is\jxi} \widehat{\calF}\Bigl[ \underline{\Phi}_\real (-\Delta+1)^{-1} \nabla \cdot \bigl( u_{1,\real} \nabla u_{3,\real} - u_{3,\real} \nabla u_{1,\real} \\ 
        &\quad \quad \quad \quad \quad \quad \quad \quad \quad \quad \quad \quad \quad \quad \quad \quad \quad \quad \quad \quad + u_{1,\imag} \nabla u_{3,\imag} - u_{3,\imag} \nabla u_{1,\imag} \bigr) \Bigr](\xi) \, \ud s \\ 
        &\quad + \int_t^T \jxi^{-1} e^{-is\jxi} \widehat{\calF}\Bigl[ \underline{\Phi}_\real (-\Delta+1)^{-1} \Bigl( U' \bigl( u_{2,\real} u_{3,\real} + u_{2,\imag} u_{3,\imag} + u_{1,\real} u_{4,\real} + u_{1,\imag} u_{4,\imag} \bigr) \Bigr) \Bigr](\xi) \, \ud s \\
        &\quad + \int_t^T \jxi^{-1} e^{-is\jxi} \widehat{\calF}\Bigl[ \underline{\Phi}_\real (-\Delta+1)^{-1} \Bigl( \bigl\{\text{cubic and higher-order nonlinearities}\bigr\} \Bigr) \Bigr](\xi) \, \ud s. 
    \end{aligned}
\end{equation}
We briefly discuss the weighted energy estimate, the high Sobolev norm estimate as well as the ILED estimate for the terms on the right-hand side of \eqref{equ:contribution_type6_decomposition}. 
The first term on the right-hand side of \eqref{equ:contribution_type6_decomposition} can be treated exactly in the same manner as the first term on the right-hand side of \eqref{equ:contribution_type5_decomposition}.
The second term on the right-hand side of \eqref{equ:contribution_type6_decomposition} is effectively of the same type as contributions of type (1), which have already been treated above.
The third term on the right-hand side of \eqref{equ:contribution_type6_decomposition} is again spatially localized and quadratic, whence essentially of the same type as contributions of type (4), which have been discussed in detail above. Finally, we again skip the details for the tedious, but harmless higher-order nonlinearities in the fourth term on the right-hand side of \eqref{equ:contribution_type6_decomposition}.

\subsubsection{Contributions of type (7)--(10)} \label{equ:subsec:contributions_type7_to_10}
The contributions of type (7)--(10) can all be considered to be terms of at least cubic order and are thus harmless. We omit the details.

\subsection{Proof of Proposition~\ref{prop:dispersive_decay_g}} \label{subsec:proof_dispersive_decay_g}

It remains to deduce the dispersive decay estimate for the evolutions of the flat profiles $\bmg_\real^\flat(t)$ and $\bmg_\imag^\flat(t)$ asserted in Proposition~\ref{prop:dispersive_decay_g}.

\begin{proof}[Proof of Proposition~\ref{prop:dispersive_decay_g}]
    For $\ast \in \{\real,\imag\}$ and $0 \leq t \leq T$, we decompose 
    \begin{equation} \label{equ:dispersive_decay_g_decomposition}
        e^{it\jD} \bmg^\flat_\ast(t) = e^{it\jD} \bigl( \bmg^\flat_\ast(t) - \bmB_\ast[\bmf, \bmf](t) \bigr) + e^{it\jD} \bmB_\ast[\bmf, \bmf](t) 
    \end{equation}
    with $\bmB_\ast[\bmf, \bmf](t)$ defined in \eqref{equ:definition_bmB_boundary_terms}.
    For the first term on the right-hand side, the standard dispersive estimate from Lemma~\ref{lem:flat_dispersive_estimate} gives
    \begin{equation}
        \bigl\| e^{it\jD} \bigl( \bmg^\flat_\ast(t) - \bmB_\ast[\bmf, \bmf](t) \bigr) \bigr\|_{L^\infty_x(\bbR^2)} \lesssim \frac{1}{\jt^{1-\delta}} \Bigl\| \jxi^2 \jap{\nabla_\xi} \Bigl( \widehat{\bmg}^\flat_{\ast}(t,\xi) - \widehat{\calF}\bigl[ \bmB_\ast[\bmf, \bmf](t) \bigr](\xi) \Bigr) \Bigr\|_{L^2_\xi(\bbR^2)}.
    \end{equation}
    The weighted energy bounds from Proposition~\ref{prop:weighted_energies_g} then yield the asserted dispersive decay for the first term on the right-hand side of \eqref{equ:dispersive_decay_g_decomposition}.
    Instead, we estimate the second term on the right-hand side of \eqref{equ:dispersive_decay_g_decomposition} directly. 
    Recall that the components of $\bmB_\ast[\bmf, \bmf](t)$ encompass the four families of bilinear operators \eqref{equ:families_of_boundary_operators}.
    We illustrate how to estimate the boundary term \eqref{equ:calQ1_contribution_type1_boundary_term}, which belongs to the family $\calB_{\kappa_1 \kappa_2}^{\pvdots,j}$. The estimates for all other types of boundary terms are completely analogous. Here we have 
    \begin{equation}
        \begin{aligned}
            &\Bigl\| e^{it\jD} \widehat{\calF}^{-1}\Bigl[ \widehat{\calB}_{\kappa_1\kappa_2}^{\pvdots,1}\bigl[ f_{k,\real}, f_{l,\real} \bigr](t,\xi) \Bigr] \Bigr\|_{L^\infty_x} \\ 
            &\lesssim \sum_{n_1, n_2 \geq 0} \, \biggl\| \int_{\bbR^2} \int_{\bbR^2} \int_{\bbR^2} e^{ix\cdot(\xi_1+\xi_2+\sigma)} \frakm_{n_1,n_2}(\xi_1,\xi_2,\sigma) \varphi_{n_1}(\xi_1) e^{it\kappa_1 \jap{\xi_1}} \widehat{f}^{\kappa_1}_{k,\real}(t,\xi_1) \\
            &\qquad \qquad \qquad \qquad \qquad \qquad \qquad \times \varphi_{n_2}(\xi_2) e^{it\kappa_2\jap{\xi_2}} \widehat{f}^{\kappa_2}_{l,\real}(t,\xi_2) \, \pvdots \, \frac{(\sigma)_1}{|\sigma|^3} \, \ud \sigma \, \ud \xi_1 \, \ud \xi_2 \biggr\|_{L^\infty_x} 
        \end{aligned}
    \end{equation}
    with
\begin{equation}
    \begin{aligned}
        \frakm_{n_1,n_2}(\xi_1,\xi_2,\sigma) &:= \jap{\xi_1+\xi_2+\sigma}^{-1} \frac{1}{\Psi_{\kappa_1\kappa_2}(\xi_1+\xi_2+\sigma,\xi_1,\xi_2)} \\
        &\qquad \qquad \times \widetilde{\varphi}_{n_1}(\xi_1) \widetilde{\varphi}_{n_2}(\xi_2) \varphi_{\leq -M}\bigl( R(\xi_1,\xi_2) \sigma \bigr) \varphi_{\leq 1}(\sigma). 
    \end{aligned}
\end{equation}    
By Lemma~\ref{lem:derivative_bounds_quadratic_phase} we have for integers $0 \leq a,b,c \leq 3$,
\begin{equation}
    \begin{aligned}
        \bigl| \nabla_{\xi_1}^a \nabla_{\xi_2}^b \nabla_\sigma^c \frakm_{n_1,n_2}(\xi_1,\xi_2,\sigma) \bigr|  \lesssim 2^{-a n_1} 2^{-b n_2}  2^{16 \min\{n_1,n_2\}},    
    \end{aligned}
\end{equation}
whence 
\begin{equation}
    \begin{aligned}
        \bigl\| \widehat{\calF}^{-1}\bigl[ \frakm_{n_1,n_2} \bigr] \bigr\|_{L^1( (\bbR^2)^3 )} \lesssim 2^{16 \min\{n_1,n_2\}}. 
    \end{aligned}
\end{equation}
Using the bilinear estimates from Lemma~\ref{lem:bilinear_estimates} 
we obtain
\begin{equation}
        \begin{aligned}
            &\Bigl\| e^{it\jD} \widehat{\calF}^{-1}\Bigl[ \widehat{\calB}_{\kappa_1\kappa_2}^{\pvdots,1}\bigl[ f_{k,\real}, f_{l,\real} \bigr](t,\xi) \Bigr] \Bigr\|_{L^\infty_x} \\ 
            &\lesssim \sum_{n_1, n_2 \geq 0} 2^{16 \min\{n_1,n_2\}} \bigl\| P_{n_1} \bigl( e^{it\jD} \bmf_\real(t) \bigr) \bigr\|_{L^\infty_x} \bigl\| P_{n_2} \bigl( e^{it\jD} \bmf_\real(t) \bigr) \bigr\|_{L^\infty_x} \\ 
            &\lesssim \bigl\| \jD^{16} e^{it\jD} \bmf_\real(t) \bigr\|_{L^\infty_x} \bigl\| \jD e^{it\jD} \bmf_\real(t) \bigr\|_{L^\infty_x} 
            .
        \end{aligned}
\end{equation}   
By \eqref{equ:GNSell} the last line is bounded by $t^{-2+3\delta}$ for long times $\varepsilon^{-2} \leq t \leq T$ and by $\varepsilon^{2-20\delta} \jt^{-1-7\delta}$ for short times $0 \leq t \leq \varepsilon^{-2}$, which is faster than the dispersive decay of the first term on the right-hand side of \eqref{equ:dispersive_decay_g_decomposition}.
\end{proof}

\subsection{Transfer estimates} \label{subsec:transfer_estimates}

In this subsection we provide a collection of technical estimates that allow us to bound various quantities involving the flat component $\bmv^\flat(t)$ of the good part in terms of the solutions $\bmv^\flat_\ast(t)$, $\ast \in \{\real,\imag\}$, to \eqref{equ:evol_equ_bmv_real_flat}, \eqref{equ:evol_equ_bmv_imag_flat} and their profiles $\bmg_\ast^\flat(t)$.
For the latter, the key bounds have already been established in Propositions~\ref{prop:weighted_energies_g}, \ref{prop:dispersive_decay_g}, \ref{prop:HN_g}, \ref{prop:HN_ILED_g}, which constitute the main outcome of Section~\ref{sec:working_title_flat_bounds}.

We shall informally refer to the technical estimates in Lemma \ref{lem:transferring_bounds_vflat_variables} below as {\it transfer estimates}. Their purpose is to transfer the previously established bounds for the auxiliary evolutions $\bmv^\flat_\ast(t)$, $\ast \in \{\real,\imag\}$, and their profiles $\bmg_\ast^\flat(t)$ into corresponding bounds for quantities involving the flat component $\bmv^\flat(t)$. Such estimates will be used repeatedly in the next section, where we analyze the evolution equations for the sharp component $\bmv^\sharp(t)$, and in the last Section~\ref{sec:conclusion_proof_thm_main}, where we conclude the proof of Theorem~\ref{thm:main}.
As a preparation for the proof of Lemma~\ref{lem:transferring_bounds_vflat_variables} we recall the angular decompositions 
\begin{equation} \label{equ:transferring_bounds_angular_decomp}
        \bmv_\ast^\flat(t,x) = \sum_{n\in\bbZ} \bmv_{\ast,n}^\flat(t,r) e^{in\theta}, \quad \ast \in \{\real, \imag\}.        
\end{equation}
Moreover, we emphasize that $\jap{D} = \sqrt{1-\Delta}$ and $|D| = \sqrt{-\Delta}$ preserve the degree as well as the decomposition into angular modes. For instance, if $f(r)$ is a radially symmetric function, then
\begin{equation} \label{equ:transferring_bounds_identities}
    \begin{aligned}
            \jD \bigl( f(r) e^{in\theta} \bigr) &= e^{in\theta} \bigl( \jap{D_n} f \bigr)(r), \quad n \in \bbZ, \\
            |D| \bigl( f(r) e^{in\theta} \bigr) &= e^{in\theta} \bigl( |D_n| f \bigr)(r), \quad n \in \bbZ, \\
            e^{it\jD} \bigl( f(r) e^{in\theta} \bigr) &= e^{in\theta} \bigl( e^{it\jap{D_n}} f \bigr)(r), \quad n \in \bbZ,
    \end{aligned}
\end{equation}
with 
\begin{equation}
        \jap{D_n} := \sqrt{-\partial_r^2-r^{-1}\partial_r+n^2r^{-2}+1}, \qquad |D_n| := \sqrt{-\partial_r^2-r^{-1}\partial_r+n^2r^{-2}}, \quad n \in \bbZ.
\end{equation}
Analogous statements hold for the action of $\jD$, $|D|$, and $e^{it\jD}$ on $f(r) \sin(m\theta)$ and $f(r) \cos(m\theta)$ for $m \in \bbN_0$.

\begin{lemma} \label{lem:transferring_bounds_vflat_variables}
For any $1\leq p\leq \infty$, $k \geq 0$, and $l \in \{0,1\}$, we have 
\begin{equation} \label{equ:transferring_bounds_vflat_Lp}
    \begin{aligned}
    \bigl\| \jap{D}^k \bigl( e^{i\theta} \pt^l \bmv^\flat(t) \bigr) \bigr\|_{L^p_x(\bbR^2)} + \bigl\| \jap{D_1}^k \pt^l \bmv^\flat(t) \bigr\|_{L^p_{r\ud r}(\bbR_+)} &\lesssim \sum_{\ast\in\{\real,\imag\}} \bigl\| \jap{D}^k \pt^l \bmv^\flat_\ast(t) \bigr\|_{L^p_x(\bbR^2)} \\
    &\lesssim \sum_{\ast\in\{\real,\imag\}} \bigl\| \jD^{k+l} e^{it\jD} \bmg_\ast^\flat(t) \bigr\|_{L^p_x(\bbR^2)}.
    \end{aligned}
\end{equation}
Moreover, it holds that
\begin{equation} \label{equ:transferring_bounds_vflat_weighted_energy}
    \begin{aligned}
        &\bigl\| \jD^2 \jap{x} e^{-it\jap{D}} (2i\jap{D})^{-1} (\partial_t+i\jap{D}) \Re\bigl( e^{i\theta}\bmv^\flat(t) \bigr) \bigr\|_{L^2_x(\bbR^2)} \\
        &\quad + \bigl\| \jD^2 \jap{x} e^{-it\jap{D}} (2i\jap{D})^{-1} (\partial_t+i\jap{D}) \Im\bigl( e^{i\theta}\bmv^\flat(t) \bigr) \bigr\|_{L^2_x(\bbR^2)} \\
        &\lesssim \sum_{\ast \in \{\real,\imag\}} \bigl\| \jD^2 \jap{x} e^{-it\jap{D}} (2i\jap{D})^{-1} (\partial_t+i\jap{D}) \bmv^\flat_\ast(t) \bigr\|_{L^2_x(\bbR^2)} \\
        &\lesssim \sum_{\ast\in\{\real,\imag\}} \bigl\| \jD^2 \jx \bmg_\ast^\flat(t) \bigr\|_{L^2_x(\bbR^2)},
    \end{aligned}
\end{equation}
as well as
\begin{equation} \label{equ:transferring_bounds_vflat_HN}
    \begin{aligned}
        &\bigl\| (2i\jD)^{-1} (\pt + i\jD) \Re\bigl( e^{i\theta} \bmv^\flat(t) \bigr) \bigr\|_{H^N_x(\bbR^2)} + \bigl\| (2i\jD)^{-1} (\pt + i\jD) \Im\bigl( e^{i\theta} \bmv^\flat(t) \bigr) \bigr\|_{H^N_x(\bbR^2)} \\
        &\lesssim \sum_{\ast \in \{\real,\imag\}} \bigl\| (2i\jD)^{-1} (\pt + i\jD) \bmv_\ast^\flat(t) \bigr\|_{H^N_x(\bbR^2)} \lesssim \sum_{\ast\in\{\real,\imag\}} \bigl\| \bmg_\ast^\flat(t) \bigr\|_{H^N_x(\bbR^2)}.
    \end{aligned}
\end{equation}
Additionally, for any $k \geq 0$ and $\sigma \geq 0$, it holds that
\begin{equation} \label{equ:transferring_bounds_vflat_ILED}
    \begin{aligned}
        &\bigl\| \jx^{-\sigma} |D|^k (2i\jD)^{-1} (\pt + i\jD) \Re\bigl( e^{i\theta} \bmv^\flat(t) \bigr) \bigr\|_{L^2_t([0,T]; L^2_x(\bbR^2))} \\
        &\quad +  \bigl\| \jx^{-\sigma} |D|^k (2i\jD)^{-1} (\pt + i\jD) \Im\bigl( e^{i\theta} \bmv^\flat(t) \bigr) \bigr\|_{L^2_t([0,T]; L^2_x(\bbR^2))} \\
        &\lesssim \sum_{\ast \in \{\real,\imag\}} \bigl\| \jx^{-\sigma} |D|^k (2i\jD)^{-1} (\pt + i\jD) \bmv_\ast^\flat(t) \bigr\|_{L^2_t([0,T]; L^2_x(\bbR^2))} \\
        &\lesssim \sum_{\ast \in \{\real,\imag\}} \bigl\| \jx^{-\sigma} |D|^k e^{it\jD} \bmg_\ast^\flat(t) \bigr\|_{L^2_t([0,T]; L^2_x(\bbR^2))},
    \end{aligned}
\end{equation}
as well as
\begin{equation} \label{equ:transferring_bounds_V_bmvflat_term}
    \bigl\| \jD^k \jx^\sigma \bigl( e^{i\theta} \bfV \bmv^\flat(t) \bigr) \bigr\|_{L^2_x(\bbR^2)} \lesssim \sum_{\ast\in\{\real,\imag\}} \, \bigl\| \jD^k \jx^\sigma \bfV \bmv_\ast^\flat(t) \bigr\|_{L^2_x(\bbR^2)}.
\end{equation}
Finally, we have for any $k \geq 0$ and $\sigma \geq 0$ that 
\begin{equation} \label{equ:transferring_bounds_bfPc_vflat}
    \begin{aligned}
        \bigl\| \jD^k \jx^\sigma \bigl( e^{i\theta} \bfP_c \bmv^\flat(t) \bigr) \bigr\|_{L^2_x(\bbR^2)} &\lesssim \sum_{\ast\in\{\real,\imag\}} \bigl\| \jD^k \jx^\sigma e^{it\jD} \bmg_\ast^\flat(t) \bigr\|_{L^2_x(\bbR^2)}, \\
        \bigl\| \jD^k \jx^\sigma \bigl( e^{i\theta} \bfP_c \pt \bmv^\flat(t) \bigr) \bigr\|_{L^2_x(\bbR^2)} &\lesssim \sum_{\ast\in\{\real,\imag\}} \bigl\| \jD^{k+1} \jx^\sigma e^{it\jD} \bmg_\ast^\flat(t) \bigr\|_{L^2_x(\bbR^2)}.
    \end{aligned}
\end{equation}
\end{lemma}
\begin{proof}
    We begin with the $L^p$-estimate \eqref{equ:transferring_bounds_vflat_Lp} for integers $k \geq 0$ and $l \in \{0,1\}$.
    Since 
    \begin{equation*}
        \bigl\| \jap{D}^k \bigl( e^{i\theta} \pt^l \bmv^\flat(t) \bigr) \bigr\|_{L^p_x(\bbR^2)} = \bigl\| e^{i\theta} \jap{D_1}^k \bigl( \pt^l \bmv^\flat(t) \bigr) \bigr\|_{L^p_x(\bbR^2)} \lesssim \bigl\| \jap{D_1}^k \bigl( \pt^l \bmv^\flat(t) \bigr) \bigr\|_{L^p_{r\ud r}(\bbR_+)}, 
    \end{equation*}
    it suffices to estimate the latter term.
    From 
    \begin{equation*}
        \jap{D_1}^k \pt^l \bmv^\flat(t,r) = \frac{1}{2\pi} \int_{0}^{2\pi} e^{-i\theta} \bigl( \jD^k \pt^l \bmv^\flat_{\real}(t,r,\theta)+i\jD^k \pt^l\bmv^\flat_{\imag}(t,r,\theta) \bigr) \, \ud \theta,
    \end{equation*}
    we then obtain the asserted bound
    \begin{equation}
        \begin{aligned}
            \bigl\| \jap{D_1}^k \pt^l \bmv^\flat(t,r) \bigr\|_{L^p_{r\ud r}} \lesssim \sum_{\ast\in\{\real,\imag\}} \bigl\| \jD^k \pt^l \bmv^\flat_\ast \bigr\|_{L^1_{\ud\theta}L^p_{r\ud r}} &\lesssim \sum_{\ast\in\{\real,\imag\}} \bigl\| \jD^k \pt^l \bmv^\flat_\ast \bigr\|_{L^p_x(\bbR^2)} \\
            &\lesssim \sum_{\ast\in\{\real,\imag\}} \bigl\| \jD^{k+l} e^{it\jD} \bmg_\ast^\flat(t) \bigr\|_{L^p_x(\bbR^2)},
        \end{aligned}
    \end{equation}
    where the last estimate follows straight from the definition of the flat profiles $\bmg_\ast^\flat(t)$, $\ast \in \{\real,\imag\}$, in \eqref{equ:definition_bmg_flat_profile}.

    Let us now turn to the details of the proof of the weighted Sobolev estimate \eqref{equ:transferring_bounds_vflat_weighted_energy}. 
    From \eqref{equ:transferring_bounds_angular_decomp} and \eqref{equ:transferring_bounds_identities}, we infer
    \begin{equation}
        \jD^2 \jx e^{-it\jD} (2i\jD)^{-1} (\pt + i\jD) \bmv_\ast^\flat = \sum_{n\in\bbZ} e^{in\theta} \bigl( \jap{D_n}^2 \jap{r} e^{-it\jap{D_n}} (2i\jap{D_n})^{-1} (\pt+i\jap{D_n}) \bmv_{\ast,n}^\flat \bigr), 
    \end{equation}
    and thus by $L^2_{\ud \theta}$-orthogonality,
    \begin{equation}
        \begin{aligned}
        &\sum_{n\in\bbZ} \, \bigl\| \jap{D_n}^2 \jap{r} e^{-it\jap{D_n}} (2i\jap{D_n})^{-1} (\pt+i\jap{D_n}) \bmv_{\ast,n}^\flat \bigr\|_{L^2_{r\ud r}}^2 \\
        &\quad = \bigl\| \jD^2 \jx e^{-it\jD} (2i\jD)^{-1} (\pt + i\jD) \bmv_\ast^\flat \bigr\|_{L^2_x(\bbR^2)}^2.
        \end{aligned}
    \end{equation}
    In particular, the preceding identity implies
    \begin{equation} \label{equ:transferring_bounds_vflat_pm1_cos_sin}
        \begin{aligned}
            &\bigl\| \jD^2 \jx e^{-it\jD} (2i\jD)^{-1} (\pt + i\jD) \bigl( \cos(\theta) \bmv_{\ast,\pm 1}^\flat \bigr) \bigr\|_{L^2_x(\bbR^2)}^2 \\
            &\quad + \bigl\| \jD^2 \jx e^{-it\jD} (2i\jD)^{-1} (\pt + i\jD) \bigl( \sin(\theta) \bmv_{\ast,\pm 1}^\flat \bigr) \bigr\|_{L^2_x(\bbR^2)}^2 \\
            &= 2\pi \bigl\| \jap{D_1}^2 \jap{r} e^{-it\jap{D_1}} (2i\jap{D_1})^{-1} (\pt+i\jap{D_1}) \bmv_{\ast,\pm 1}^\flat \bigr\|_{L^2_{r\ud r}}^2 \\
            &\lesssim \bigl\| \jD^2 \jx e^{-it\jD} (2i\jD)^{-1} (\pt + i\jD) \bmv_\ast^\flat \bigr\|_{L^2_x(\bbR^2)}^2.
        \end{aligned}
    \end{equation}
    Since $\bmv_\real^\flat$ and $\bmv_\imag^\flat$ are real-valued solutions to \eqref{equ:evol_equ_bmv_real_flat}, respectively to \eqref{equ:evol_equ_bmv_imag_flat}, we have $\overline{\bmv^\flat_{\ast,1}}=\bmv^\flat_{\ast,-1}$. Hence, $\Re( e^{i\theta} \bmv^\flat ) = \Re \bigl( e^{i\theta}(\bmv^\flat_{\real,1}+i\bmv^\flat_{\imag,1}) \bigr)$ can be written as a linear combination of $\cos(\theta) \bmv^\flat_{\ast,\pm1}(t,r)$ and $\sin(\theta) \bmv^\flat_{\ast,\pm1}(t,r)$, and analogously for $\Im( e^{i\theta} \bmv^\flat)$.
    The asserted estimate \eqref{equ:transferring_bounds_vflat_weighted_energy} then follows immediately from \eqref{equ:transferring_bounds_vflat_pm1_cos_sin}.

    The proofs of the estimates \eqref{equ:transferring_bounds_vflat_HN}, \eqref{equ:transferring_bounds_vflat_ILED}, \eqref{equ:transferring_bounds_V_bmvflat_term} proceed analogously, and we omit the details.
    For the last asserted estimate \eqref{equ:transferring_bounds_bfPc_vflat} we note that upon inserting the definition of the projection to the continuous spectral subspace, we obtain
    \begin{equation}
        \begin{aligned}
            &\bigl\| \jD^k \jx^\sigma \bigl( e^{i\theta} \bfP_c \bmv^\flat(t) \bigr) \bigr\|_{L^2_x(\bbR^2)} \\
            &\leq \bigl\| \jD^k \jx^\sigma \bigl( e^{i\theta} \bmv^\flat(t) \bigr) \bigr\|_{L^2_x(\bbR^2)} + \sum_{j=1,2} \bigl\| \jD^k \jx^\sigma \bigl( \langle \bmY_j, \bmv^\flat(t) \rangle e^{i\theta} \bmY_j \bigr) \bigr\|_{L^2_x(\bbR^2)} \\
            &\lesssim \bigl\| \jD^k \jx^\sigma \bigl( e^{i\theta} \bmv^\flat(t) \bigr) \bigr\|_{L^2_x(\bbR^2)} + \bigl\| \bmv^\flat(t) \bigr\|_{L^2_{r\ud r}},
        \end{aligned}
    \end{equation}
    and similarly for the weighted Sobolev norm involving $\pt \bmv^\flat(t)$. At this point the proof proceeds analogously as before.
\end{proof}

\section{Bounds for the Sharp Component of the Good Part} \label{sec:working_title_sharp_bounds}

In this section we establish bounds for the flat profiles $\bmg_\ast^\sharp(t)$, $\ast \in \{\real, \imag\}$, defined in \eqref{equ:definition_bmg_sharp_profile}, which are associated with the sharp component $\bmv^\sharp(t)$ of the good part of the radiation term. 
We begin with dispersive decay estimates for the evolution $e^{it\jD} \bmg_\ast^\sharp(t)$.

\begin{proposition} \label{prop:dispersive_decay_g_sharp}
    Suppose the bootstrap assumptions \eqref{equ:bootstrap_assumption_z}, \eqref{equ:bootstrap_assumption_weighted_f}, \eqref{equ:bootstrap_assumption_Linfty_f}, \eqref{equ:bootstrap_assumption_HN_f}, \eqref{equ:bootstrap_assumption_HN_ILED_f} are in place.  
    Then we have for $\ast \in \{\real, \imag\}$,    
    \begin{equation} \label{equ:dispersive_decay_g_sharp}
        \bigl\|  e^{it\jD} \bmg^\sharp_\ast(t)\bigr\|_{L^\infty_x(\bbR^2)} 
        \leq C_2 \left\{ \begin{aligned}
                            &\varepsilon^{1-6\delta} \jt^{-1+3\delta}, \quad &&0 \leq t \leq \varepsilon^{-2}, \\    
                            &\varepsilon^{1-10\delta} t^{-1+\delta},  &&\varepsilon^{-2} \leq t \leq T,
                         \end{aligned} \right.
    \end{equation}    
    where $C_2 \equiv C_2(C_0)$ depends on the size of the constant $C_0$ in the bootstrap assumptions \eqref{equ:bootstrap_assumption_weighted_f}, \eqref{equ:bootstrap_assumption_Linfty_f}, \eqref{equ:bootstrap_assumption_HN_f}, \eqref{equ:bootstrap_assumption_HN_ILED_f}.
\end{proposition}
\begin{proof}
    We present the details for the dispersive estimate for $e^{it\jD} \bmg_\real^\sharp(t)$, noting that the case of $e^{it\jD} \bmg_\imag^\sharp(t)$ can be handled analogously.
    By the definition \eqref{equ:definition_bmg_sharp_profile} of the flat profile $\bmg^\sharp_\real(t)$, we have 
    \begin{equation} \label{equ:dispersive_decay_vsharp_decomposition}
        \begin{aligned}
            e^{it\jD} \bmg_\real^\sharp(t) = \frac12 \Re \bigl( e^{i\theta} \bfP_c \bmv^\sharp(t) \bigr) + \frac{1}{2i\jD} \Re \bigl( e^{i\theta} \bfP_c \pt \bmv^\sharp(t) \bigr) =: I(t) + II(t).
        \end{aligned}
    \end{equation}
    From the Duhamel formulation of the evolution equation \eqref{equ:evol_equ_bmv_sharp} for $\bmv^\sharp(t)$, we infer
    \begin{equation} \label{equ:dispersive_decay_vsharp_termI}
        \begin{aligned}
            \bigl\| I(t) \bigr\|_{L^\infty_x(\bbR^2)} \leq \bigl\| \bfP_c \bmv^\sharp(t) \bigr\|_{L^\infty_r} &\leq \bigl\| \cos\bigl(t\sqrt{\bfM}\bigr) \bfP_c \bmv^\sharp(0) \bigr\|_{L^\infty_r} + \bigl\| \sin\bigl(t\sqrt{\bfM}\bigr) \bfM^{-\frac12} \bfP_c \pt \bmv^\sharp(0) \bigr\|_{L^\infty_r} \\ 
            &\quad + \biggl\| \int_0^t \frac{\sin\bigl((t-s) \sqrt{\bfM}\bigr)}{\sqrt{\bfM}} \bfP_c \bigl( \bfV \bmv^\flat(s) \bigr) \, \ud s \biggr\|_{L^\infty_r}. 
        \end{aligned}
    \end{equation}
    For the first two terms on the right-hand side of \eqref{equ:dispersive_decay_vsharp_termI}, we use the dispersive decay estimate \eqref{eq:Linftyweighteddispersivebound1} from Proposition~\ref{prop:dispersive_estimate_with_potential} for some fixed $p \geq \frac{2}{\delta}$, the transfer estimate \eqref{equ:transferring_bounds_bfPc_vflat}, \eqref{equ:new_assumption2}, and Proposition~\ref{prop:weighted_energies_g}, to obtain for $0 \leq t \leq T$ that
    \begin{equation}
        \begin{aligned}
            &\bigl\| \cos\bigl(t\sqrt{\bfM}\bigr) \bfP_c \bmv^\sharp(0) \bigr\|_{L^\infty_r} + \bigl\| \sin\bigl(t\sqrt{\bfM}\bigr) \bfM^{-\frac12} \bfP_c \pt \bmv^\sharp(0) \bigr\|_{L^\infty_r} \\
            &\lesssim \jt^{-1+\delta} \Bigl( \bigl\| \jx \bigl( e^{i\theta} \bfP_c \bmv^\sharp(0) \bigr) \bigr\|_{H^2_x} + \bigl\| \jx \bigl( e^{i\theta} \bfP_c \pt \bmv^\sharp(0) \bigr) \bigr\|_{H^1_x} \Bigr) \\
            &\lesssim \jt^{-1+\delta} \Bigl( \bigl\| \jx \bigl( e^{i\theta} \bmu(0) \bigr) \bigr\|_{H^2_x} + \bigl\| \jx \bigl( e^{i\theta} \pt\bmu(0) \bigr) \bigr\|_{H^1_x} \\            
            &\quad \quad \quad \quad \quad \quad + \bigl\| \jx \bigl( e^{i\theta} \bfP_c \bmv^\flat(0) \bigr) \bigr\|_{H^2_x} + \bigl\| \jx \bigl( e^{i\theta} \bfP_c \pt \bmv^\flat(0) \bigr) \bigr\|_{H^1_x} \Bigr) \\ 
            &\lesssim \jt^{-1+\delta} \biggl( \bigl\| \jx \bigl( e^{i\theta} \bmu(0) \bigr) \bigr\|_{H^2_x} + \bigl\| \jx \bigl( e^{i\theta} \pt\bmu(0) \bigr) \bigr\|_{H^1_x} + \sum_{\ast\in\{\real,\imag\}} \bigl\| \jx \bmg_\ast^\flat(0) \bigr\|_{H^2_x} \biggr) \\    
            &\lesssim \jt^{-1+\delta} \bigl( \varepsilon + \varepsilon^{1-6\delta} \bigr) \lesssim \varepsilon^{1-6\delta} \jt^{-1+\delta}.
        \end{aligned}
    \end{equation}
    Moreover, using the dispersive decay estimate \eqref{eq:Linftydispersiveboundnonsharp2} from Proposition~\ref{prop:dispersive_estimate_with_potential}, the transfer estimate \eqref{equ:transferring_bounds_vflat_Lp}, and the Gagliardo-Nirenberg-Sobolev inequality \eqref{equ:GNS}, we may bound the last term on the right-hand side of \eqref{equ:dispersive_decay_vsharp_termI} by
    \begin{equation}
        \begin{aligned}
            &\biggl\| \int_0^t \frac{\sin\bigl((t-s) \sqrt{\bfM}\bigr)}{\sqrt{\bfM}} \bfP_c \bigl( \bfV \bmv^\flat(s) \bigr) \, \ud s \biggr\|_{L^\infty_r} 
            \lesssim \int_0^t \frac{1}{\jap{t-s}} \bigl\| e^{i\theta} \bfV \bmv^\flat(s) \bigr\|_{W^{3,1}_x} \, \ud s  \\
            &\lesssim \int_0^t \frac{1}{\jap{t-s}} \bigl\| \bfV \bigr\|_{W^{3,1+}_x} \bigl\| \jD^3 \bigl( e^{i\theta} \bmv^\flat(s) \bigr) \bigr\|_{L^{\infty-}_x} \, \ud s  
            \lesssim \int_0^t \frac{1}{\jap{t-s}} \sum_{\ast \in \{\real,\imag\}} \bigl\| \jD^3 \bmv_\ast^\flat(s) \bigr\|_{L^{\infty-}_x} \, \ud s \\
            &\lesssim \sum_{\ast\in\{\real,\imag\}} \int_0^t \frac{1}{\jap{t-s}} \bigl\| \bmg_\ast^\flat(s) \bigr\|_{H^N_x}^{\frac{3-\frac{2}{\infty-}}{N-1}} \bigl\| e^{is\jD} \bmg_\ast^\flat(s) \bigr\|_{L^\infty_x}^{1-\frac{3-\frac{2}{\infty-}}{N-1}} \, \ud s. 
        \end{aligned}
    \end{equation}
    Invoking Proposition~\ref{prop:dispersive_decay_g} and Proposition~\ref{prop:HN_g}, we obtain for short times $0 \leq t \leq \varepsilon^{-2}$ the crude but sufficient bound
    \begin{equation}\label{eq:vsharpforFGR1}
        \begin{aligned}
            \biggl\| \int_0^t \frac{\sin\bigl((t-s) \sqrt{\bfM}\bigr)}{\sqrt{\bfM}} \bfP_c \bigl( \bfV \bmv^\flat(s) \bigr) \, \ud s \biggr\|_{L^\infty_r} 
            &\lesssim \int_0^t \frac{1}{\jap{t-s}} \bigl( \varepsilon^{\frac32 - 2\delta} \bigr)^{\frac{3-\frac{2}{\infty-}}{N-1}} \bigl( \varepsilon^{1-6\delta} \js^{-1+\delta} \bigr)^{1-\frac{3-\frac{2}{\infty-}}{N-1}} \, \ud s \\ 
            &\lesssim \varepsilon^{1-6\delta} \int_0^t \frac{1}{\jap{t-s}} \frac{1}{\js^{1-2\delta}} \, \ud s \\
            &\lesssim \varepsilon^{1-6\delta} \jt^{-1+3\delta}.
        \end{aligned}
    \end{equation}
    Instead, for large times $\varepsilon^{-2} \leq t \leq T$ the bounds from Proposition~\ref{prop:dispersive_decay_g} and Proposition~\ref{prop:HN_g} imply the crude but sufficient bounds
    \begin{equation}\label{eq:vsharpforFGR2}
        \begin{aligned}
            \biggl\| \int_0^t \frac{\sin\bigl((t-s) \sqrt{\bfM}\bigr)}{\sqrt{\bfM}} \bfP_c \bigl( \bfV \bmv^\flat(s) \bigr) \, \ud s \biggr\|_{L^\infty_r} 
            &\lesssim \int_0^{\varepsilon^{-2}} \frac{1}{\jap{t-s}} \bigl( \varepsilon^{\frac32 - 2\delta} \bigr)^{\frac{3-\frac{2}{\infty-}}{N-1}} \bigl( \varepsilon^{1-6\delta} \js^{-1+\delta} \bigr)^{1-\frac{3-\frac{2}{\infty-}}{N-1}} \, \ud s \\
            &\quad + \int_{\varepsilon^{-2}}^t \frac{1}{\jap{t-s}} \bigl( \varepsilon^{\frac32 - 2\delta} \bigr)^{\frac{3-\frac{2}{\infty-}}{N-1}} \bigl( s^{-\frac32+4\delta} \bigr)^{1-\frac{3-\frac{2}{\infty-}}{N-1}} \, \ud s \\ 
            &\lesssim \varepsilon^{1-12\delta} t^{-1} \lesssim \varepsilon^{1-10\delta} t^{-1+\delta}.
        \end{aligned}
    \end{equation}    

    In order to estimate the second term $II(t)$ on the right-hand side of \eqref{equ:dispersive_decay_vsharp_decomposition}, we use the mapping property \eqref{equ:bfMnr_mapping_jDinverse_Linfty} and the relation \eqref{equ:relation_bfM_bfMnr} along with the Duhamel formula for $\bfP_c \pt \bmv^\sharp(t)$ to find that
    \begin{equation} \label{equ:dispersive_decay_vsharp_termII} 
        \begin{aligned}
            \bigl\| II(t) \bigr\|_{L^\infty_x(\bbR^2)} &= \bigl\| \jD^{-1} \bfM_{\mathrm{nr}}^{\frac12} \bfM_{\mathrm{nr}}^{-\frac12} \bigl( e^{i\theta} \bfP_c \pt \bmv^\sharp(t) \bigr) \bigr\|_{L^\infty_x(\bbR^2)} \\
            &\lesssim \bigl\| \bfM_{\mathrm{nr}}^{-\frac12} \bigl( e^{i\theta} \bfP_c \pt \bmv^\sharp(t) \bigr) \bigr\|_{L^\infty_x(\bbR^2)} \simeq \bigl\| \bfM^{-\frac12} \bfP_c \pt\bmv^\sharp(t) \bigr\|_{L^\infty_r} \\ 
            &\lesssim \bigl\| \sin\bigl(t\sqrt{\bfM}\bigr) \bfP_c \bmv^\sharp(0) \bigr\|_{L^\infty_r} + \bigl\| \cos\bigl(t\sqrt{\bfM}\bigr) \bfM^{-\frac12} \bfP_c \pt \bmv^\sharp(0) \bigr\|_{L^\infty_r} \\ 
            &\quad + \biggl\| \int_0^t \frac{\cos\bigl((t-s) \sqrt{\bfM}\bigr)}{\sqrt{\bfM}} \bfP_c \bigl( \bfV \bmv^\flat(s) \bigr) \, \ud s \biggr\|_{L^\infty_r}. 
        \end{aligned}
    \end{equation}
    In view of \eqref{equ:dispersive_decay_vsharp_termI}, it follows immediately that \eqref{equ:dispersive_decay_vsharp_termII} leads to the same decay estimates for $II(t)$ as those established above for $I(t)$.

    Combining the preceding estimates and observing that they are all dominated by the corresponding bounds \eqref{equ:dispersive_decay_g_sharp} in the statement of Proposition~7.1, we obtain the asserted decay estimate.
\end{proof}

Next, we establish high Sobolev norm bounds for the profiles $\bmg_\ast^\sharp(t)$, $\ast \in \{\real, \imag\}$, and ILED estimates for the evolutions $e^{it\jD} \bmg_\ast^\sharp(t)$, $\ast \in \{\real,\imag\}$.

\begin{proposition} \label{prop:HN_g_sharp}
    Suppose the bootstrap assumptions \eqref{equ:bootstrap_assumption_z}, \eqref{equ:bootstrap_assumption_weighted_f}, \eqref{equ:bootstrap_assumption_Linfty_f}, \eqref{equ:bootstrap_assumption_HN_f}, \eqref{equ:bootstrap_assumption_HN_ILED_f} are in place.  
    Then we have for $\ast \in \{\real, \imag\}$, 
    \begin{equation} \label{equ:HN_g_sharp}
        \sup_{0\leq t\leq T} \, \bigl\| \bmg^\sharp_\ast(t)\bigr\|_{H^N_x(\bbR^2)} + \sum_{1 \leq k \leq N} \bigl\| \jx^{-\frac32-\kappa} |D|^k e^{it\jD} \bmg^\sharp_\ast(t)\bigr\|_{L^2_t([0,T]; L^2_x(\bbR^2))} \leq C_1 \varepsilon,
    \end{equation}
    where $C_1 \geq 1$ is some absolute constant.
\end{proposition}
\begin{proof}
    For $\ast \in \{\real, \imag\}$ we have 
    \begin{equation} \label{equ:HN_gsharp_decomposition}
        \bigl\| \bmg_\ast^\sharp(t) \bigr\|_{H^N_x} \lesssim \bigl\| e^{i\theta} \bfP_c \bmv^\sharp(t) \bigr\|_{H^N_x} + \bigl\| e^{i\theta} \bfP_c \pt \bmv^\sharp(t) \bigr\|_{H^{N-1}_x}
    \end{equation}
    as well as 
    \begin{equation} \label{equ:ILED_gsharp_decomposition}
        \begin{aligned}
            &\sum_{1 \leq k \leq N} \, \bigl\| \jx^{-\frac32-\kappa} |D|^k e^{it\jD} \bmg^\sharp_\ast(t)\bigr\|_{L^2_t([0,T]; L^2_x(\bbR^2))} \lesssim \bigl\| \jx^{-\frac32-\kappa} e^{it\jD} \bmg^\sharp_\ast(t)\bigr\|_{L^2_t([0,T]; H^N_x(\bbR^2))} \\ 
            &\lesssim \bigl\| \jx^{-\frac32-\kappa} e^{i\theta} \bfP_c \bmv^\sharp(t) \bigr\|_{L^2_t([0,T]; H^N_x(\bbR^2))} + \bigl\| \jx^{-\frac32-\kappa} e^{i\theta} \bfP_c \partial_t \bmv^\sharp(t) \bigr\|_{L^2_t([0,T]; H^{N-1}_x(\bbR^2))}.
        \end{aligned}
    \end{equation}    
    We begin by estimating the first terms on the right-hand side of \eqref{equ:HN_gsharp_decomposition} and of \eqref{equ:ILED_gsharp_decomposition}.
    Using the mapping property \eqref{equ:bfMnr_mapping_Mk} and the relation \eqref{equ:relation_bfM_bfMnr},
    it follows that the first term on the right-hand side of \eqref{equ:HN_gsharp_decomposition} can be bounded by
    \begin{equation} \label{equ:HN_gsharp_termI_start}
        \begin{aligned}
            \bigl\| e^{i\theta} \bfP_c \bmv^\sharp(t) \bigr\|_{H^N_x} &\lesssim \bigl\| e^{i\theta} \bfP_c \bmv^\sharp(t) \bigr\|_{L^2_x} + \bigl\| \bfM_{\mathrm{nr}}^{\frac{N}{2}} \bigl( e^{i\theta} \bfP_c \bmv^\sharp(t) \bigr) \bigr\|_{L^2_x} \\ 
            &\lesssim \bigl\| \bfP_c \bmv^\sharp(t) \bigr\|_{L^2_{r \ud r}} + \bigl\| \bfM^{\frac{N}{2}} \bfP_c \bmv^\sharp(t) \bigr\|_{L^2_{r \ud r}}.
        \end{aligned}
    \end{equation}
    Similarly, by the mapping property \eqref{equ:bfMnr_mapping_Mk} we can estimate the first term on the right-hand side of \eqref{equ:ILED_gsharp_decomposition} by
    \begin{equation} \label{equ:ILED_gsharp_termI_start}
        \begin{aligned}
            &\bigl\| \jx^{-\frac32-\kappa} e^{i\theta} \bfP_c \bmv^\sharp(t) \bigr\|_{L^2_t([0,T]; H^N_x(\bbR^2))} \\
            &\quad \lesssim \bigl\| \jap{r}^{-\frac32-\kappa} \bfP_c \bmv^\sharp(t) \bigr\|_{L^2_t([0,T]; L^2_{r \ud r})} + \bigl\| \jap{r}^{-\frac32-\kappa} \bfM^{\frac{N}{2}} \bfP_c \bmv^\sharp(t) \bigr\|_{L^2_t([0,T]; L^2_{r \ud r})}.
        \end{aligned}
    \end{equation}
    From the Duhamel formula for $\bfP_c \bmv^\sharp(t)$ and the dual ILED estimate \eqref{equ:dual_ILED} from Proposition~\ref{prop:ILED} for some $0 < \kappa \ll 1$, we then conclude for the first term on the right-hand side of \eqref{equ:HN_gsharp_termI_start} for all $0 \leq t \leq T$ that
    \begin{equation}
        \begin{aligned}
            \bigl\| \bfP_c \bmv^\sharp(t) \bigr\|_{L^2_{r\ud r}} &\lesssim \bigl\| \cos\bigl(t\sqrt{\bfM}\bigr) \bfP_c \bmv^\sharp(0) \bigr\|_{L^2_{r \ud r}} + \bigl\| \sin\bigl(t\sqrt{\bfM}\bigr) \bfM^{-\frac12} \bfP_c \pt \bmv^\sharp(0) \bigr\|_{L^2_{r \ud r}} \\ 
            &\quad + \biggl\| \int_0^t \frac{\sin\bigl((t-s) \sqrt{\bfM}\bigr)}{\sqrt{\bfM}} \bfP_c \bigl( \bfV \bmv^\flat(s) \bigr) \, \ud s \biggr\|_{L^2_{r \ud r}} \\ 
            &\lesssim \bigl\| \bfP_c \bmv^\sharp(0) \bigr\|_{L^2_{r \ud r}} + \bigl\| \bfM^{-\frac12} \bfP_c \pt \bmv^\sharp(0) \bigr\|_{L^2_{r \ud r}} 
            + \bigl\| \jr^{1+\kappa} \bfV \bmv^\flat(t) \bigr\|_{L^2_t([0,T]; L^2_{r \ud r})}.  
        \end{aligned}
    \end{equation}
    Similarly, using the ILED estimate \eqref{equ:ILED} as well as the inhomogeneous ILED estimate \eqref{equ:inhomog_ILED} from Proposition~\ref{prop:ILED}, we also bound the first term on the right-hand side of \eqref{equ:ILED_gsharp_termI_start} by
    \begin{equation}
        \begin{aligned}
            &\bigl\| \jap{r}^{-\frac32-\kappa} \bfP_c \bmv^\sharp(t) \bigr\|_{L^2_t([0,T]; L^2_{r \ud r})} \\
            &\quad \lesssim \bigl\| \bfP_c \bmv^\sharp(0) \bigr\|_{L^2_{r \ud r}} + \bigl\| \bfP_c \pt \bmv^\sharp(0) \bigr\|_{L^2_{r \ud r}} 
            + \bigl\| \jr^{1+\kappa} \bfV \bmv^\flat(t) \bigr\|_{L^2_t([0,T]; L^2_{r \ud r})}. 
        \end{aligned}
    \end{equation}

    Then by the transfer estimate \eqref{equ:transferring_bounds_V_bmvflat_term}, \eqref{equ:new_assumption2}, and Proposition~\ref{prop:HN_g}, we have 
    \begin{equation}
        \begin{aligned}
            \bigl\| \bfP_c \bmv^\sharp(0) \bigr\|_{L^2_{r \ud r}} + \bigl\| \bfP_c \pt \bmv^\sharp(0) \bigr\|_{L^2_{r \ud r}} &\lesssim \bigl\| e^{i\theta} \bmu(0) \bigr\|_{L^2_x} + \bigl\| e^{i\theta} \pt \bmu(0) \bigr\|_{L^2_x} + \sum_{\ast\in\{\real,\imag\}} \bigl\| \jD \bmg^\flat_\ast(0) \bigr\|_{L^2_x} \\
            &\leq C_1 \varepsilon + C_2 \varepsilon^{\frac32-2\delta},
        \end{aligned}
    \end{equation}
    where $C_1 \geq 1$ is some absolute constant, while $C_2 \equiv C_2(C_0)$ is a constant whose size depends on the size of the constant $C_0$ in the bootstrap assumptions \eqref{equ:bootstrap_assumption_weighted_f}, \eqref{equ:bootstrap_assumption_Linfty_f}, \eqref{equ:bootstrap_assumption_HN_f}, \eqref{equ:bootstrap_assumption_HN_ILED_f}.
    Moreover, by the transfer estimate \eqref{equ:transferring_bounds_vflat_Lp}, interpolation, and the bounds from Proposition~\ref{prop:dispersive_decay_g} as well as from Proposition~\ref{prop:HN_g}, we infer the crude but sufficient bound 
    \begin{equation} \label{equ:HN_gsharp_bfV_bmvflat_L2L2}
        \begin{aligned}
            &\bigl\| \jr^{1+\kappa} \bfV \bmv^\flat(t) \bigr\|_{L^2_t([0,T]; L^2_{r \ud r})} \lesssim \Bigl\| \bigl\| \jap{r}^{1+\kappa} \bfV \bigr\|_{L^{\frac{8}{3}}_{r \, \ud r}} \bigl\| \bmv^\flat(t) \bigr\|_{L^{8}_{r \ud r}} \Bigr\|_{L^2_t([0,T])} \\
            &\lesssim \sum_{\ast\in\{\real,\imag\}} \Bigl\| \bigl\| e^{it\jD} \bmg_\ast^\flat(t) \bigr\|_{L^{8}_x} \Bigr\|_{L^2_t([0,T])}  
            \lesssim \sum_{\ast\in\{\real,\imag\}} \Bigl\| \bigl\| \bmg_\ast^\flat(t) \bigr\|_{L^2_x}^{\frac{1}{4}} \bigl\| e^{it\jD} \bmg_\ast^\flat(t) \bigr\|_{L^{\infty}_x}^{\frac{3}{4}} \Bigr\|_{L^2_t([0,T])} \\ 
            &\lesssim \sum_{\ast\in\{\real,\imag\}} \Bigl\| \bigl\| \bmg_\ast^\flat(t) \bigr\|_{L^2_x}^{\frac{1}{4}} \bigl\| e^{it\jD} \bmg_\ast^\flat(t) \bigr\|_{L^{\infty}_x}^{\frac{3}{4}} \Bigr\|_{L^2_t([0,\varepsilon^{-2}])} + \sum_{\ast\in\{\real,\imag\}} \Bigl\| \bigl\| \bmg_\ast^\flat(t) \bigr\|_{L^2_x}^{\frac{1}{4}} \bigl\| e^{it\jD} \bmg_\ast^\flat(t) \bigr\|_{L^{\infty}_x}^{\frac{3}{4}} \Bigr\|_{L^2_t([\varepsilon^{-2},T])} \\
            &\lesssim \Bigl\| \bigl( \varepsilon^{\frac32-2\delta} \bigr)^{\frac{1}{4}} \bigl( \varepsilon^{1-6\delta} \jt^{-1+\delta} \bigr)^{\frac{3}{4}} \Bigr\|_{L^2_t([0,\varepsilon^{-2}])} + \Bigl\| \bigl( \varepsilon^{\frac32-2\delta} \bigr)^{\frac{1}{4}} \bigl( t^{-\frac32+4\delta} \bigr)^{\frac{3}{4}} \Bigr\|_{L^2_t([\varepsilon^{-2}, T])} \\ 
            &\lesssim \varepsilon^{\frac{9}{8}-5\delta} + \varepsilon^{\frac{13}{8}-\frac{13}{2}\delta} \lesssim \varepsilon^{\frac{9}{8}-5\delta}.
        \end{aligned}
    \end{equation}

    For the second term on the right-hand side of \eqref{equ:HN_gsharp_termI_start}, we obtain from the Duhamel formula for $\bfP_c \bmv^\sharp(t)$, the relation \eqref{equ:relation_bfM_bfMnr}, the mapping property \eqref{equ:bfMnr_mapping_Hk}, and the dual ILED estimate \eqref{equ:dual_ILED} from Proposition~\ref{prop:ILED} for some fixed $0 < \kappa \ll 1$ that for all $0 \leq t \leq T$,
    \begin{equation}
        \begin{aligned}
            \bigl\| \bfM^{\frac{N}{2}} \bfP_c \bmv^\sharp(t) \bigr\|_{L^2_{r \ud r}} &\lesssim \bigl\| \cos(t\sqrt{\bfM}) \bfM^{\frac{N}{2}} \bfP_c \bmv^\sharp(0) \bigr\|_{L^2_{r\ud r}} + \bigl\| \sin(t\sqrt{\bfM}) \bfM^{\frac{N-1}{2}} \bfP_c \pt \bmv^\sharp(0) \bigr\|_{L^2_{r \ud r}} \\ 
            &\quad + \biggl\| \int_0^t \sin\bigl((t-s)\sqrt{\bfM}\bigr) \bfM^{-\frac12} \, \bfM^{\frac{N}{2}} \bfP_c \bigl( \bfV \bmv^\flat(s) \bigr) \, \ud s \biggr\|_{L^2_{r\ud r}} \\ 
            &\lesssim \bigl\| e^{i\theta} \bfP_c \bmv^\sharp(0) \bigr\|_{H^N_x} + \bigl\| e^{i\theta} \bfP_c \pt \bmv^\sharp(0) \bigr\|_{H^{N-1}_x} + \bigl\| \jap{r}^{1+\kappa} \bfM^{\frac{N}{2}} \bigl( \bfV \bmv^\flat(t) \bigr) \bigr\|_{L^2_t([0,T]; L^2_{r\ud r})}.
        \end{aligned}
    \end{equation}
    Similarly, we obtain the same bound for the second term on the right-hand side of \eqref{equ:ILED_gsharp_termI_start} using the mapping property \eqref{equ:bfMnr_mapping_Hk}, the ILED estimate \eqref{equ:ILED} and the inhomogeneous ILED estimate \eqref{equ:inhomog_ILED} from Proposition~\ref{prop:ILED},
    \begin{equation} \label{equ:ILED_gsharp_MN}
        \begin{aligned}
            &\bigl\| \jap{r}^{-\frac32-\kappa} \bfM^{\frac{N}{2}} \bfP_c \bmv^\sharp(t) \bigr\|_{L^2_t([0,T]; L^2_{r \ud r})} \\
            &\lesssim \bigl\| e^{i\theta} \bfP_c \bmv^\sharp(0) \bigr\|_{H^N_x} + \bigl\| e^{i\theta} \bfP_c \pt \bmv^\sharp(0) \bigr\|_{H^{N-1}_x} + \bigl\| \jap{r}^{1+\kappa} \bfM^{\frac{N}{2}} \bigl( \bfV \bmv^\flat(t) \bigr) \bigr\|_{L^2_t([0,T]; L^2_{r\ud r})}.
        \end{aligned}
    \end{equation}
    For the first two terms on the right-hand side of \eqref{equ:ILED_gsharp_MN} we have by the transfer estimate \eqref{equ:transferring_bounds_bfPc_vflat},
    \begin{equation}
        \begin{aligned}
            &\bigl\| e^{i\theta} \bfP_c \bmv^\sharp(0) \bigr\|_{H^N_x} + \bigl\| e^{i\theta} \bfP_c \pt \bmv^\sharp(0) \bigr\|_{H^{N-1}_x} \\
            &\lesssim \bigl\| e^{i\theta} \bmu(0) \bigr\|_{H^N_x} + \bigl\| e^{i\theta} \pt\bmu(0) \bigr\|_{H^{N-1}_x} + \bigl\| e^{i\theta} \bfP_c \bmv^\flat(0) \bigr\|_{H^N_x} + \bigl\| e^{i\theta} \bfP_c \pt\bmv^\flat(0) \bigr\|_{H^{N-1}_x} \\
            &\lesssim \bigl\| e^{i\theta} \bmu(0) \bigr\|_{H^N_x} + \bigl\| e^{i\theta} \pt\bmu(0) \bigr\|_{H^{N-1}_x} + \sum_{\ast \in \{\real,\imag\}} \bigl\| \bmg_\ast^\flat(0) \bigr\|_{H^N_x} \\
            &\leq C_1 \varepsilon + C_2 \varepsilon^{\frac32-2\delta}.
        \end{aligned}
    \end{equation}
    For the last term on the right-hand side of \eqref{equ:ILED_gsharp_MN}, we have by the mapping property \eqref{equ:bfMnr_mapping_Hk} and the transfer estimate~\eqref{equ:transferring_bounds_V_bmvflat_term} that
    \begin{equation}
        \begin{aligned}
             \bigl\| \jap{r}^{1+\kappa} \bfM^{\frac{N}{2}} \bigl( \bfV \bmv^\flat(t) \bigr) \bigr\|_{L^2_{r \ud r}} \leq \bigl\| \jx^{1+\kappa} \bfM_{nr}^{\frac{N}{2}} \bigl( e^{i\theta} \bfV \bmv^\flat(t) \bigr) \bigr\|_{L^2_x}  
             &\lesssim \bigl\| \jD^N \jx^{1+\kappa} \bigl( e^{i\theta} \bfV \bmv^\flat(t) \bigr) \bigr\|_{L^2_x} \\ 
             &\lesssim \sum_{\ast\in\{\real,\imag\}} \, \bigl\| \jD^N \jx^{1+\kappa} \bfV \bmv_\ast^\flat(t) \bigr\|_{L^2_x},
        \end{aligned}
    \end{equation}
    whence
    \begin{equation}
        \begin{aligned}
            &\bigl\| \jap{r}^{1+\kappa} \bfM^{\frac{N}{2}} \bigl( \bfV \bmv^\flat(t) \bigr) \bigr\|_{L^2_t([0,T]; L^2_{r\ud r})} 
            \lesssim \sum_{\ast\in\{\real,\imag\}} \sum_{0 \leq |\gamma| \leq N} \, \Bigl\| \bigl\| \partial^\gamma \bigl( \jx^{1+\kappa} \bfV \bmv_\ast^\flat(t) \bigr) \bigr\|_{L^2_x} \Bigr\|_{L^2_t([0,T])} \\ 
            &\lesssim \sum_{\ast\in\{\real,\imag\}} \, \Bigl\| \bigl\| \jx^{1+\kappa} \bfV \bmv_\ast^\flat(t) \bigr\|_{L^2_x} \Bigr\|_{L^2_t([0,T])} \\
            &\quad + \sum_{\ast\in\{\real,\imag\}} \sum_{1 \leq |\gamma| \leq N} \sum_{0 \leq |\beta| \leq N-|\gamma|} \, \Bigl\| \bigl\| \partial^\beta \bigl( \jx^{1+\kappa} \bfV \bigr) \partial^\gamma \bmv_\ast^\flat(t) \bigr\|_{L^2_x} \Bigr\|_{L^2_t([0,T])}.
        \end{aligned}
    \end{equation}
    The first term on the right-hand side in the preceding line enjoys the same bound as \eqref{equ:HN_gsharp_bfV_bmvflat_L2L2}, while for the second term we invoke the ILED bound \eqref{equ:HN_ILED_g} from Proposition~\ref{prop:HN_ILED_g} to conclude that
    \begin{equation}
        \begin{aligned}
            &\sum_{\ast\in\{\real,\imag\}} \sum_{1 \leq |\gamma| \leq N} \sum_{0 \leq |\beta| \leq N-|\gamma|} \, \Bigl\| \bigl\| \partial^\beta \bigl( \jx^{1+\kappa} \bfV \bigr) \partial^\gamma \bmv_\ast^\flat(t) \bigr\|_{L^2_x} \Bigr\|_{L^2_t([0,T])} \\ 
            &\lesssim \sum_{\ast\in\{\real,\imag\}}\sum_{1 \leq k \leq N} \, \bigl\| \jx^{\frac32+2\kappa} \bfV \bigr\|_{W^{N,\infty}_x} \bigl\| \jx^{-\frac12-\kappa} |D|^k e^{it\jD} \bmg_\ast^\flat(t) \bigr\|_{L^2_t([0,T]; L^2_x)} \lesssim \varepsilon^{\frac32 - 2\delta}.
        \end{aligned}
    \end{equation}

    Finally, the second terms on the right-hand side of \eqref{equ:HN_gsharp_decomposition} and of \eqref{equ:ILED_gsharp_decomposition} can be estimated analogously, leading to the same bounds. Collecting all of the preceding estimates, we find that
    \begin{equation} 
        \sup_{0\leq t\leq T} \, \bigl\| \bmg^\sharp_\ast(t)\bigr\|_{H^N_x(\bbR^2)} + \sum_{1 \leq k \leq N} \bigl\| \jx^{-\frac32-\kappa} |D|^k e^{it\jD} \bmg^\sharp_\ast(t)\bigr\|_{L^2_t([0,T]; L^2_x(\bbR^2))} \leq C_1 \varepsilon + C_2 \bigl( \varepsilon^{\frac32-2\delta} + \varepsilon^{\frac98 -5\delta} \bigr),
    \end{equation}
    where $C_1 \geq 1$ is some absolute constant and $C_2 \equiv C_2(C_0)$ is a constant whose size depends on the size of the constant $C_0$ in the bootstrap assumptions \eqref{equ:bootstrap_assumption_weighted_f}, \eqref{equ:bootstrap_assumption_Linfty_f}, \eqref{equ:bootstrap_assumption_HN_f}, \eqref{equ:bootstrap_assumption_HN_ILED_f}.
    Since $C_0 \geq 1$ will be chosen sufficiently large depending only on some absolute constants, we can choose the small constant $0 < \varepsilon_0 \ll 1$ sufficiently small depending on the size of $C_0$ so that $C_2 ( \varepsilon^{\frac32-2\delta} + \varepsilon^{\frac98 -5\delta} ) \leq C_1 \varepsilon$. This concludes the proof of Proposition~\ref{prop:HN_g_sharp}.
\end{proof}

Finally, we obtain weighted energy estimates for the profiles $\bmg_\ast^\sharp(t)$, $\ast \in \{\real,\imag\}$.

\begin{proposition} \label{prop:weighted_energies_g_sharp}
    Suppose the bootstrap assumptions \eqref{equ:bootstrap_assumption_z}, \eqref{equ:bootstrap_assumption_weighted_f}, \eqref{equ:bootstrap_assumption_Linfty_f}, \eqref{equ:bootstrap_assumption_HN_f}, \eqref{equ:bootstrap_assumption_HN_ILED_f} are in place.  
    Then we have for $\ast \in \{\real, \imag\}$, 
    \begin{equation}
        \bigl\| \jxi^2 \jap{\nabla_\xi} \widehat{\bmg}^\sharp_{\ast}(t,\xi) \bigr\|_{L^2_\xi(\bbR^2)} \leq C_2 \left\{ \begin{aligned}
                            &\varepsilon^{1-10\delta} \jt^{\frac12}, \quad &&0 \leq t \leq \varepsilon^{-2}, \\    
                            &t^{5\delta},  &&\varepsilon^{-2} \leq t \leq T,
                         \end{aligned} \right.
    \end{equation}
    where $C_2 \equiv C_2(C_0)$ depends on the size of the constant $C_0$ in the bootstrap assumptions \eqref{equ:bootstrap_assumption_weighted_f}, \eqref{equ:bootstrap_assumption_Linfty_f}, \eqref{equ:bootstrap_assumption_HN_f}, \eqref{equ:bootstrap_assumption_HN_ILED_f}.
\end{proposition}
\begin{proof}
    We have 
    \begin{equation}
        \bigl\| \jxi^2 \jap{\nabla_\xi} \widehat{\bmg}^\sharp_{\ast}(t,\xi) \bigr\|_{L^2_\xi(\bbR^2)} \lesssim \bigl\| \jxi^2 \widehat{\bmg}^\sharp_{\ast}(t,\xi) \bigr\|_{L^2_\xi(\bbR^2)} + \bigl\| \jxi^2 \nabla_\xi \widehat{\bmg}^\sharp_{\ast}(t,\xi) \bigr\|_{L^2_\xi(\bbR^2)}.
    \end{equation}
    For the first term we have by Proposition~\ref{prop:HN_g_sharp} that $\bigl\| \jxi^2 \widehat{\bmg}^\sharp_{\ast}(t,\xi) \bigr\|_{L^2_\xi(\bbR^2)} \lesssim \varepsilon$ for $0 \leq t \leq T$. Using Proposition~\ref{prop:transference} and the evolution equation \eqref{equ:evol_equ_bmv_sharp} for $\bmv^\sharp(t)$, we find that  
    \begin{equation} \label{equ:weighted_energies_g_sharp_black_box_applied}
        \begin{aligned}
            &\sum_{\ast \in \{\real, \imag\}} \bigl\| \jxi^2 \nabla_\xi \widehat{\bmg}^{\sharp}_\ast(t,\xi) \bigr\|_{L^2_\xi} \\
            &\lesssim \bigl\| \jx e^{i\theta} \bfP_c \bmv^\sharp(0) \bigr\|_{H^2_x} + \bigl\| \jx e^{i\theta} \bfP_c (\pt \bmv^\sharp)(0) \bigr\|_{H^1_x} \\ 
            &\quad + \bigl\| \jx e^{i\theta} \bfP_c \bigl( \bfV \bmv^\flat(0,r) \bigr) \bigr\|_{L^2_x} + \bigl\| \jx e^{i\theta} \bfP_c \bigl( \bfV (\pt \bmv^\flat)(0,r) \bigr) \bigr\|_{L^2_x} \\ 
            &\quad + \bigl\| \jx e^{i\theta} \bfP_c \bigl( \bfV \bmv^\flat(t,r) \bigr) \bigr\|_{L^2_x} + \sum_{\ast \in \{\real, \imag\}} \Bigl( \bigl\| \bmg_\ast^\sharp(t) \bigr\|_{H^2_x} + \jt \bigl\| \jD^2 e^{it\jD} \bmg_\ast^\sharp(t) \bigr\|_{L^{\infty-}_x} \Bigr) \\ 
            &\quad + \sum_{j=0,1} \bigl\| s \cdot \jx \bigl( e^{i\theta} \bfP_c \bigl( \bfV \partial_s^j \bmv^\flat(s,r) \bigr) \big) \bigr\|_{L^2_s([0,t]; L^2_x)} + \bigl\| s \cdot \jx \jD \bigl( e^{i\theta}  \bfP_c \bigl( \bfV \bmv^\flat(s,r) \bigr) \bigr) \big) \bigr\|_{L^2_s([0,t]; L^2_x)} \\
            &\quad + \sum_{j=1,2} \bigl\| \jap{x} \bigl( e^{i\theta} \bfP_c \bigl( \bfV \partial_s^j \bmv^\flat(s,r) \bigr) \bigr) \bigr\|_{L^1_s([0,t]; L^2_x)} + \bigl\| \jx \jD \bigl( e^{i\theta} \bfP_c \bigl( \bfV \bmv^\flat(s,r) \bigr) \bigr) \bigr\|_{L^1_s([0,t]; L^2_x)}.
        \end{aligned}
    \end{equation}
    For the first two terms on the right-hand side of \eqref{equ:weighted_energies_g_sharp_black_box_applied} we have by the transfer estimate~\eqref{equ:transferring_bounds_V_bmvflat_term} and by Proposition~\ref{prop:weighted_energies_g},
    \begin{equation}
        \begin{aligned}
            &\bigl\| \jx e^{i\theta} \bfP_c \bmv^\sharp(0) \bigr\|_{H^2_x} + \bigl\| \jx e^{i\theta} \bfP_c (\pt \bmv^\sharp)(0) \bigr\|_{H^1_x} \\
            &\lesssim \bigl\| \jx e^{i\theta} \bfP_c \bmu(0) \bigr\|_{H^2_x} + \bigl\| \jx e^{i\theta} \bfP_c (\pt \bmu)(0) \bigr\|_{H^1_x} + \bigl\| \jx e^{i\theta} \bfP_c \bmv^\flat(0) \bigr\|_{H^2_x} + \bigl\| \jx e^{i\theta} \bfP_c (\pt \bmv^\flat)(0) \bigr\|_{H^1_x} \\ 
            &\lesssim \bigl\| \jx e^{i\theta} \bmu(0) \bigr\|_{H^2_x} + \bigl\| \jx e^{i\theta} (\pt \bmu)(0) \bigr\|_{H^1_x} + \sum_{\ast \in \{\real,\imag\}} \bigl\| \jx \bmg_\ast^\flat(0) \bigr\|_{H^2_x} \lesssim \varepsilon + \varepsilon^{1-6\delta} \lesssim \varepsilon^{1-6\delta}.
        \end{aligned}
    \end{equation}
    For the third and fourth terms on the right-hand side of \eqref{equ:weighted_energies_g_sharp_black_box_applied} we conclude by H\"older's inequality, Sobolev estimates, the transfer estimate \eqref{equ:transferring_bounds_vflat_Lp}, and Proposition~\ref{prop:HN_g} that
    \begin{equation}
        \begin{aligned}
            &\bigl\| \jx e^{i\theta} \bfP_c \bigl( \bfV \bmv^\flat(0,r) \bigr) \bigr\|_{L^2_x} + \bigl\| \jx e^{i\theta} \bfP_c \bigl( \bfV (\pt \bmv^\flat)(0,r) \bigr) \bigr\|_{L^2_x} \\
            &\lesssim \bigl\| \jx \bfV \bigr\|_{L^{2+}_x} \Bigl( \bigl\| e^{i\theta} \bmv^\flat(0,r) \bigr\|_{L^{\infty-}_x} + \bigl\| e^{i\theta} (\pt \bmv^\flat)(0,r) \bigr\|_{L^{\infty-}_x} \Bigr) \lesssim \sum_{\ast \in \{\real,\imag\}} \bigl\| \bmg_\ast^\flat(0) \bigr\|_{H^2_x} \lesssim \varepsilon^{\frac32 - 2\delta}.
        \end{aligned}
    \end{equation}
    Analogously, the fifth term on the right-hand side of \eqref{equ:weighted_energies_g_sharp_black_box_applied} is bounded for all $0 \leq t \leq T$ by
    \begin{equation}
        \bigl\| \jx e^{i\theta} \bfP_c \bigl( \bfV \bmv^\flat(t,r) \bigr) \bigr\|_{L^2_x} \lesssim \sum_{\ast \in \{\real,\imag\}} \bigl\| \bmg_\ast^\flat(t) \bigr\|_{H^2_x} \lesssim \varepsilon^{\frac32-2\delta}.
    \end{equation}
    Next, we turn to the sixth term on the right-hand side of \eqref{equ:weighted_energies_g_sharp_black_box_applied}.
    By Proposition~\ref{prop:HN_g_sharp} we have for all times $0 \leq t \leq T$,  
    \begin{equation}
        \begin{aligned}
            \sum_{\ast \in \{\real, \imag\}} \bigl\| \bmg_\ast^\sharp(t) \bigr\|_{H^2_x} \lesssim \varepsilon.
        \end{aligned}
    \end{equation}
    Using  the interpolation inequality~\eqref{equ:GNS} along with the bounds from Proposition~\ref{prop:dispersive_decay_g_sharp} and Proposition~\ref{prop:HN_g_sharp}, we conclude for all times $0 \leq t \leq T$ that 
    \begin{align*}
        \sum_{\ast\in\{\real,\imag\}} \jt \big\Vert \jD^2 e^{it\langle D\rangle} \bmg_\ast^\sharp(t)\big\Vert_{L^{\infty-}_x} \lesssim \sum_{\ast\in\{\real,\imag\}} \jt \Vert \bmg_\ast^\sharp(t)\Vert_{H^N_x}^{\frac{2-\frac{2}{\infty-}}{N-1}} \Vert e^{it\jD} \bmg_\ast^\sharp(t)\Vert_{L^\infty_x}^{1-\frac{2-\frac{2}{\infty-}}{N-1}} \lesssim \varepsilon^{1-6\delta} \jt^{4\delta}.
    \end{align*}
    For the seventh term on the right-hand side of \eqref{equ:weighted_energies_g_sharp_black_box_applied} we infer by the transfer estimate \eqref{equ:transferring_bounds_vflat_Lp} and the interpolation inequality \eqref{equ:GNS},
    \begin{equation}
        \begin{aligned}
            \sum_{j=0,1} \bigl\| s \cdot \jx \bigl( e^{i\theta} \bfP_c \bigl( \bfV \partial_s^j \bmv^\flat(s,r) \bigr) \big)\bigr\|_{L^2_s([0,t]; L^2_x)} &\lesssim \sum_{j=0,1} \, \Bigl\| s \cdot \bigl\| \jx \bfV \bigr\|_{L^{2+}_x} \bigl\| e^{i\theta} \partial_s^j \bmv^\flat(s,r) \bigr\|_{L^{\infty-}_x} \Bigr\|_{L^2_s([0,t])} \\ 
            &\lesssim \sum_{\ast \in \{\real,\imag\}} \, \Bigl\| s \cdot \bigl\| \jD e^{is\jD} \bmg_\ast^\flat(s) \bigr\|_{L^{\infty-}_x} \Bigr\|_{L^2_s([0,t])} \\ 
            &\lesssim \sum_{\ast \in \{\real,\imag\}} \, \biggl\| s \cdot \bigl\| \bmg_\ast^\flat(s) \bigr\|_{H^N_x}^{\frac{1-\frac{2}{\infty-}}{N-1}} \bigl\| e^{is\jD} \bmg_\ast^\flat(s) \bigr\|_{L^{\infty }_x}^{1-\frac{1-\frac{2}{\infty-}}{N-1}} \biggr\|_{L^2_s([0,t])}.
        \end{aligned}
    \end{equation}
    Using the bounds from Proposition~\ref{prop:dispersive_decay_g} and Proposition~\ref{prop:HN_g}, we conclude for short times $0 \leq t \leq \varepsilon^{-2}$ that
    \begin{equation}
        \begin{aligned}
            \sum_{j=0,1} \bigl\| s \cdot \jx \bigl( e^{i\theta} \bfP_c \bigl( \bfV \partial_s^j \bmv^\flat(s,r) \bigr) \big) \bigr\|_{L^2_s([0,t]; L^2_x)} \lesssim \varepsilon^{1-10\delta} t^{\frac12},
        \end{aligned}
    \end{equation}
    while for long times $\varepsilon^{-2} \leq t \leq T$ we obtain
    \begin{equation}
        \begin{aligned}
            \sum_{j=0,1} \bigl\| s \cdot \jx \bigl( e^{i\theta} \bfP_c \bigl( \bfV \partial_s^j \bmv^\flat(s,r) \bigr) \big) \bigr\|_{L^2_s([0,t]; L^2_x)} \lesssim  t^{5\delta}.
        \end{aligned}
    \end{equation}    
    Similarly, using the transfer estimate \eqref{equ:transferring_bounds_vflat_Lp} and \eqref{equ:GNS}, we estimate the eighth term on the right-hand side of \eqref{equ:weighted_energies_g_sharp_black_box_applied} by 
    \begin{equation}
        \begin{aligned}
            \bigl\| s \cdot \jx \jD \bigl( e^{i\theta} \bfP_c \bigl( \bfV \bmv^\flat(s,r) \bigr) \bigr) \bigr\|_{L^2_s([0,t]; L^2_x)} &\lesssim \Bigl\| s \cdot \bigl\| \jx \bfV \bigr\|_{W^{1,2+}_x} \bigl\| \jD \bigl( e^{i\theta} \bmv^\flat(s,r) \bigr) \bigr\|_{L^{\infty-}_x} \Bigr\|_{L^2_s([0,t])} \\ 
            &\lesssim \sum_{\ast \in \{\real,\imag\}} \, \Bigl\| s \cdot \bigl\| \jD e^{is\jD} \bmg_\ast^\flat(s) \bigr\|_{L^{\infty-}_x} \Bigr\|_{L^2_s([0,t])}.
        \end{aligned}
    \end{equation}
    At this point it is clear that this term enjoys the same bounds as the preceding seventh term. 
    Next, we turn to the ninth term on the right-hand side of \eqref{equ:weighted_energies_g_sharp_black_box_applied}. 
    In the case $j=1$ we infer by the transfer estimate \eqref{equ:transferring_bounds_vflat_Lp}, \eqref{equ:GNS}, and the bounds from Proposition~\ref{prop:dispersive_decay_g} and Proposition~\ref{prop:HN_g} that uniformly for all times $0 \leq t \leq T$,
    \begin{equation}
        \begin{aligned}
            \bigl\| \jap{x} \bigl( e^{i\theta} \bfP_c \bigl( \bfV \partial_s \bmv^\flat(s,r) \bigr) \bigr) \bigr\|_{L^1_s([0,t]; L^2_x)} &\lesssim \sum_{\ast \in \{\real,\imag\}} \int_0^t \bigl\| \jx \bfV \bigr\|_{L^{2+}_x} \bigl\| \jD e^{is\jD} \bmg_\ast^\flat(s) \bigr\|_{L^{\infty-}_x} \, \ud s \lesssim \varepsilon^{1-10\delta}.
        \end{aligned}
    \end{equation}
    In the case $j=2$ we have to insert the evolution equation \eqref{equ:evol_equ_bmv_flat} for $\bmv^\flat(t)$. Along with the transfer estimate \eqref{equ:transferring_bounds_vflat_Lp}, the interpolation inequality \eqref{equ:GNS}, the bounds from Proposition~\ref{prop:dispersive_decay_g} and Proposition~\ref{prop:HN_g}, and the much faster decay of $\bigl\| e^{i\theta} \bmN_c(s) \bigr\|_{L^{\infty-}_x}$, we obtain uniformly for all times $0 \leq t \leq T$ that  
    \begin{equation}
        \begin{aligned}
            &\bigl\| \jap{x} \bigl( e^{i\theta} \bfP_c \bigl( \bfV \partial_s^2 \bmv^\flat(s,r) \bigr) \bigr) \bigr\|_{L^1_s([0,t]; L^2_x)} \\
            &\lesssim \int_0^t \bigl\| \jx \bfV \bigr\|_{L^{2+}_x} \Bigl( \bigl\| \jD^2 \bigl( e^{i\theta} \bmv^\flat(s,r) \bigr) \bigr\|_{L^{\infty-}_x} + \bigl\| e^{i\theta} \bmN_c(s) \bigr\|_{L^{\infty-}_x} \Bigr) \, \ud s \\ 
            &\lesssim \int_0^t \biggl( \sum_{\ast \in \{\real,\imag\}} \bigl\| \bmg^\flat_\ast(s) \bigr\|_{H^N_x}^\frac{2-\frac{2}{\infty-}}{N-1} \bigl\| e^{is\jD} \bmg^\flat_\ast(s) \bigr\|^{1-\frac{2-\frac{2}{\infty-}}{N-1}}_{L^\infty_x} + \bigl\| e^{i\theta} \bmN_c(s) \bigr\|_{L^{\infty-}_x} \biggr) \, \ud s \lesssim \varepsilon^{1-10\delta}.
        \end{aligned}
    \end{equation}
    By the same arguments as in the case $j=1$ for the ninth term, we may bound the tenth and last term on the right-hand side of \eqref{equ:weighted_energies_g_sharp_black_box_applied} by $\varepsilon^{1-10\delta}$ for all times $0 \leq t \leq T$.  

    Combining all of the preceding estimates yields the asserted weighted energy estimate in the statement of Proposition~\ref{prop:weighted_energies_g_sharp} and completes the proof.    
\end{proof}

\section{Bounds for the Bad Part} \label{sec:working_title_w_bounds}

In this section we establish bounds for the profiles $\bmh_\ast(t)$, $\ast \in \{\real, \imag\}$, related to the bad part $\bmw(t)$ of the radiation term. 
We begin with dispersive decay estimates for the evolution $e^{it\jD} \bmh_\ast(t)$.

\begin{proposition} \label{prop:dispersive_decay_h}
    Suppose the bootstrap assumptions \eqref{equ:bootstrap_assumption_z}, \eqref{equ:bootstrap_assumption_weighted_f}, \eqref{equ:bootstrap_assumption_Linfty_f}, \eqref{equ:bootstrap_assumption_HN_f}, \eqref{equ:bootstrap_assumption_HN_ILED_f} are in place.  
    Then we have for $\ast \in \{\real, \imag\}$,
    \begin{equation}
        \bigl\| e^{it\jD} \bmh_\ast(t)\bigr\|_{L^\infty_x(\bbR^2)} 
        \leq \left\{ \begin{aligned}                            
                            &C_2 \varepsilon^2 \log(1+t),  &&0 \leq t \leq \varepsilon^{-2}, \\
                            &C_2 t^{-1} \log(1+t),  &&\varepsilon^{-2} \leq t \leq T,
                         \end{aligned} \right.
    \end{equation}
    where $C_2 \equiv C_2(C_0)$ depends only on the size of the constant $C_0$ in the bootstrap assumptions \eqref{equ:bootstrap_assumption_weighted_f}, \eqref{equ:bootstrap_assumption_Linfty_f}, \eqref{equ:bootstrap_assumption_HN_f}, \eqref{equ:bootstrap_assumption_HN_ILED_f}.    
\end{proposition}
\begin{proof}
    As in the proof of Proposition~\ref{prop:dispersive_decay_g_sharp}, we discuss the details for the dispersive estimate for $e^{it\jD} \bmh_\real(t)$, noting that the case of $e^{it\jD} \bmh_\imag(t)$ can be handled analogously. 
    By the definition \eqref{equ:definition_bmh_profile} of the flat profile $\bmh_\real(t)$ we have 
    \begin{equation} \label{equ:h_dispersive_decay_decomposition}
         e^{it\jD} \bmh_\real(t) = \frac12 \Re\bigl( e^{i\theta} \bfP_c \bmw(t) \bigr) + \frac{1}{2i\jD} \Re\bigl( e^{i\theta} \bfP_c \pt \bmw(t) \bigr) =: I(t) + II(t).
    \end{equation}
    Then the Duhamel formulation of the evolution equation \eqref{equ:evol_equ_bmw} for the bad part (with zero initial data) combined with the dispersive decay estimate \eqref{eq:Linftydispersiveboundnonsharp2} from Proposition~\ref{prop:dispersive_estimate_with_potential} gives for the first term on the right-hand side of \eqref{equ:h_dispersive_decay_decomposition},
    \begin{equation}
        \begin{aligned}
            \|I(t)\|_{L^\infty_x(\bbR^2)} \leq \bigl\| \bfP_c \bmw(t) \bigr\|_{L^\infty_r} &= \biggl\| \int_0^t \frac{\sin((t-s)\sqrt{\bfM})}{\sqrt{\bfM}} \bfP_c \bmN_d(s) \, \ud s \biggr\|_{L^\infty_r} \\ 
            &\lesssim \int_0^t \frac{1}{\jap{t-s}} \, \bigl\| e^{i\theta} \bmN_d(s) \bigr\|_{W^{3,1}_x} \, \ud s. 
        \end{aligned}
    \end{equation}
    Recall that all the nonlinearities in $e^{i\theta} \bmN_d$ are spatially localized. The rate-determining contributions stem from the slowest decaying quadratic terms, where both inputs are given by an internal mode component. These ultimately dictate the dispersive decay estimates for the bad part.
    Specifically, for the contributions of quadratic source terms with two internal mode components of the schematic form $\bmY z_i z_j$ with $i,j \in \{1,2\}$ and $\bmY(x)$ some smooth, rapidly decaying function, we obtain using \eqref{equ:bootstrap_assumption_z} that
    \begin{equation} \label{equ:h_dispersive_decay_proof_z^2_contribution}
        \begin{aligned}
            &\int_0^t \frac{1}{\jap{t-s}} \, \bigl\| \bmY z_i(s) z_j(s) \bigr\|_{W^{3,1}_x} \, \ud s \\ &\quad \lesssim \int_0^t \frac{1}{\jap{t-s}} \frakz(s) \, \ud s \lesssim \int_0^t \frac{1}{\jap{t-s}} \frac{\varepsilon^2}{1+\Gamma_1 \varepsilon^2 s} \, \ud s \leq C_1 \left\{ \begin{aligned}
                            &\varepsilon^2 \log(1+t), \quad &&0 \leq t \leq \varepsilon^{-2}, \\
                            &t^{-1} \log(1+t),  &&\varepsilon^{-2} \leq t \leq T,
                         \end{aligned} \right.  
        \end{aligned}
    \end{equation}
    where $C_1 \geq 1$ is some absolute constant whose size is independent of the constant $C_0$ in the bootstrap assumptions.

    Among the other nonlinearities in $e^{i\theta} \bmN_d$, it is instructive to also discuss the contributions from quadratic terms with one input given by a radiation term and the other input given by an internal mode component, say of the schematic form $\bmY \partial_i u_{k,\real} z_j$ with $k \in \{1,2,3,4\}$, $i, j \in \{1,2\}$, and $\bmY(x)$ some smooth, rapidly decaying function. Moreover, it is useful to consider the contributions from the bad part of the temporal component $\partial_j \pt \eta_{0,d}$, $j \in \{1,2\}$. 
    All other nonlinearities in $e^{i\theta} \bmN_d$ are similar or of higher order, and thus easier to treat. 
    In view of the estimate \eqref{equ:h_dispersive_decay_proof_z^2_contribution}, we observe that for short time scales $0 \leq s \leq \varepsilon^{-2}$, we have by Sobolev embedding and the bootstrap assumptions \eqref{equ:bootstrap_assumption_z}, \eqref{equ:bootstrap_assumption_HN_f} that
    \begin{equation}
        \begin{aligned}
            \bigl\| \bmY \partial_i u_{k,\real}(s) z_j(s) \bigr\|_{W^{3,1}_x} \lesssim \bigl\| \bmY \bigr\|_{W^{3,1}_x} \bigl\| \jD^4 u_{k,\real}(s) \bigr\|_{L^\infty_x} |z_j(s)| \lesssim \bigl\|\bmf_\real(s)\bigr\|_{H^N_x} \frakz(s)^{\frac12} \lesssim C_0 \varepsilon^2,
        \end{aligned}
    \end{equation}
    while for large time scales $\varepsilon^{-2} \leq s \leq T$ the term $\bigl\| \bmY \partial_i u_{k,\real}(s) z_j(s) \bigr\|_{W^{3,1}_x}$ has faster decay than $\frakz(s)$ since $\|\jD^4 u_{k,\real}(s)\|_{L^\infty_x} \lesssim s^{-1+\frac32 \delta}$ enjoys faster $L^\infty_x$-decay by \eqref{equ:GNSell} than $\frakz(s)^{\frac12} \lesssim s^{-\frac12}$. Thus, for all quadratic nonlinearities of the schematic form $\bmY \partial_i u_{k,\real} z_j$ we may use the right-hand side of \eqref{equ:h_dispersive_decay_proof_z^2_contribution} as an upper bound for their contribution to the dispersive decay estimate for $e^{it\jD} \bmh_\real(t)$ upon replacing the absolute constant $C_1$ there by a constant $C_2 \equiv C_2(C_0)$ whose size depends on the size of the constant $C_0$ in the bootstrap assumptions. This bound still suffices to improve the bootstrap assumption~\eqref{equ:bootstrap_assumption_Linfty_f} in the proof of Theorem~\ref{thm:main}.
    The preceding argument also applies to the contributions of $\partial_j \pt \eta_{0,d}$ using interpolation and the corresponding (weighted) high Sobolev norm and $L^\infty_x$-decay estimates from Proposition~\ref{prop:eta0d_bounds} for the bad part of the temporal component.

    Proceeding as at the end of the proof of Proposition~\ref{prop:dispersive_decay_g_sharp}, one sees that the term $II(t)$ on the right-hand side of \eqref{equ:h_dispersive_decay_decomposition} enjoys the same dispersive decay estimates as derived for the term $I(t)$ above.
    This finishes the proof of the proposition.
\end{proof}

Next, we establish high Sobolev norm bounds for the profiles $\bmh_\ast(t)$, $\ast \in \{\real,\imag\}$, along with ILED estimates for the evolutions $e^{it\jD} \bmh_\ast(t)$, $\ast \in \{\real, \imag\}$.

\begin{proposition} \label{prop:HN_h}
    Suppose the bootstrap assumptions \eqref{equ:bootstrap_assumption_z}, \eqref{equ:bootstrap_assumption_weighted_f}, \eqref{equ:bootstrap_assumption_Linfty_f}, \eqref{equ:bootstrap_assumption_HN_f}, \eqref{equ:bootstrap_assumption_HN_ILED_f} are in place.  
    Then we have for $\ast \in \{\real, \imag\}$,
    \begin{equation}
        \sup_{0 \leq t \leq T} \, \bigl\| \bmh_\ast(t)\bigr\|_{H^N_x(\bbR^2)} + \sum_{1 \leq k \leq N} \bigl\| \jx^{-\frac32-\kappa} |D|^k e^{it\jD} \bmh_\ast(t)\bigr\|_{L^2_t([0,T]; L^2_x(\bbR^2))} \leq C_1 \varepsilon,
    \end{equation}
    where $C_1 \geq 1$ is some absolute constant.
\end{proposition}
\begin{proof}
    For $\ast \in \{\real, \imag\}$ we have 
    \begin{equation} \label{equ:HN_h_decomposition}
        \begin{aligned}
            \bigl\| \bmh_\ast(t) \bigr\|_{H^N_x} \lesssim \bigl\| e^{i\theta} \bfP_c \bmw(t) \bigr\|_{H^N_x} + \bigl\| e^{i\theta} \bfP_c \pt \bmw(t) \bigr\|_{H^{N-1}_x}
        \end{aligned}
    \end{equation}
    as well as 
    \begin{equation} \label{equ:ILED_h_decomposition}
        \begin{aligned}
            &\sum_{1 \leq k \leq N} \bigl\| \jx^{-\frac32-\kappa} |D|^k e^{it\jD} \bmh_\ast(t)\bigr\|_{L^2_t([0,T]; L^2_x(\bbR^2))} \lesssim \bigl\| \jx^{-\frac32-\kappa} e^{it\jD} \bmh_\ast(t)\bigr\|_{L^2_t([0,T]; H^N_x(\bbR^2))} \\ 
            &\lesssim \bigl\| \jx^{-\frac32-\kappa} e^{i\theta} \bfP_c \bmw(t) \bigr\|_{L^2_t([0,T]; H^N_x(\bbR^2))} + \bigl\| \jx^{-\frac32-\kappa} e^{i\theta} \bfP_c \partial_t \bmw(t) \bigr\|_{L^2_t([0,T]; H^{N-1}_x(\bbR^2))}.
        \end{aligned}
    \end{equation}        
    Proceeding as at the beginning of Proposition~\ref{prop:HN_g_sharp}, the first term on the right-hand side of \eqref{equ:HN_h_decomposition} is bounded by
    \begin{equation} \label{equ:HN_h_step1}
        \begin{aligned}
            \bigl\| e^{i\theta} \bfP_c \bmw(t) \bigr\|_{H^N_x} \lesssim \bigl\| \bfP_c \bmw(t) \bigr\|_{L^2_{r\ud r}} + \bigl\| \bfM^{\frac{N}{2}} \bfP_c \bmw(t) \bigr\|_{L^2_{r\ud r}}.
        \end{aligned}
    \end{equation}
    Similarly, we can estimate the first term on the right-hand side of \eqref{equ:ILED_h_decomposition} by
    \begin{equation} \label{equ:ILED_h_step1}
        \begin{aligned}
            &\bigl\| \jx^{-\frac32-\kappa} e^{i\theta} \bfP_c \bmw(t) \bigr\|_{L^2_t([0,T]; H^N_x(\bbR^2))} \\
            &\quad \lesssim \bigl\| \jap{r}^{-\frac32-\kappa} \bfP_c \bmw(t) \bigr\|_{L^2_t([0,T]; L^2_{r\ud r} )} + \bigl\| \jap{r}^{-\frac32-\kappa} \bfM^{\frac{N}{2}} \bfP_c \bmw(t) \bigr\|_{L^2_t([0,T]; L^2_{r\ud r})}.
        \end{aligned}
    \end{equation}
    We now only discuss the details for the estimates of the second terms on the right-hand sides of \eqref{equ:HN_h_step1} and of \eqref{equ:ILED_h_step1}. 
    From the Duhamel formula for $\bfP_c \bmw(t)$, the dual ILED estimate \eqref{equ:dual_ILED} from Proposition~\ref{prop:ILED} for some fixed $0 < \kappa \ll 1$, the relation \eqref{equ:relation_bfM_bfMnr}, and the mapping property \eqref{equ:bfMnr_mapping_Hk}, we obtain 
    \begin{equation}
        \begin{aligned}
            \bigl\| \bfM^{\frac{N}{2}} \bfP_c \bmw(t) \bigr\|_{L^2_{r\ud r}} &= \biggl\| \int_0^t \frac{\sin\bigl((t-s)\sqrt{\bfM})}{\sqrt{\bfM}} \bfM^{\frac{N}{2}} \bfP_c \bmN_d(s) \, \ud s \biggr\|_{L^2_{r\ud r}} \\
            &\lesssim \bigl\| \jap{r}^{1+\kappa} \bfM^{\frac{N-1}{2}} \bmN_d(t) \bigr\|_{L^2_t([0,T]; L^2_{r\ud r})} \lesssim \bigl\| \jx^{1+\kappa} \jD^{N-1} \bigl( e^{i\theta} \bmN_d(t) \bigr) \bigr\|_{L^2_t([0,T]; L^2_x)}.
        \end{aligned}
    \end{equation}
    Analogously, using in addition the inhomogeneous ILED estimate \eqref{equ:inhomog_ILED} from Proposition~\ref{prop:ILED}, we conclude that the second term on the right-hand side of \eqref{equ:ILED_h_step1} enjoys the same bound
    \begin{equation}
        \begin{aligned}
            \bigl\| \jap{r}^{-\frac32-\kappa} \bfM^{\frac{N}{2}} \bfP_c \bmw(t) \bigr\|_{L^2_t([0,T]; L^2_{r\ud r})} &\lesssim \bigl\| \jap{r}^{1+\kappa} \bfM^{\frac{N-1}{2}} \bmN_d(t) \bigr\|_{L^2_t([0,T]; L^2_{r\ud r})} \\
            &\lesssim \bigl\| \jx^{1+\kappa} \jD^{N-1} \bigl( e^{i\theta} \bmN_d(t) \bigr) \bigr\|_{L^2_t([0,T]; L^2_x)}.
        \end{aligned}
    \end{equation}
    As in the proof of the dispersive decay estimate for the bad part in the preceding Proposition~\ref{prop:dispersive_decay_h}, we only focus on the most problematic contributions from $e^{i\theta} \bmN_d(t)$ and leave the easier cases to the reader.
    For the contributions of quadratic source terms with two internal mode components of the schematic form $\bmY z_i z_j$ with $i,j \in \{1,2\}$ and $\bmY(x)$ some smooth, rapidly decaying function, we obtain using \eqref{equ:bootstrap_assumption_z} that
    \begin{equation}
        \begin{aligned}
            \bigl\| \jx^{1+\kappa} \jD^{N-1} \bigl( \bmY z_i(t) z_j(t) \bigr) \bigr\|_{L^2_t([0,T]; L^2_x)} &\lesssim \bigl\| \jx^{1+\kappa} \jD^{N-1} \bmY \bigr\|_{L^2_x} \bigl\| \frakz(t) \bigr\|_{L^2_t([0,T])} \\ 
            &\lesssim \bigl\| \varepsilon^2 \bigr\|_{L^2_t([0,\varepsilon^{-2}])} + \bigl\| t^{-1} \bigr\|_{L^2_t([\varepsilon^{-2},T])} \leq C_1 \varepsilon,
        \end{aligned}
    \end{equation}
    where we emphasize that $C_1 \geq 1$ is an absolute constant.
    Next, for the contributions of quadratic terms with one radiation term component and one internal mode component, say of the schematic form $\bmY \partial_i u_{k,\real} z_j$ with $k \in \{1,2,3,4\}$, $i,j \in \{1,2\}$, and $\bmY(x)$ some smooth, rapidly decaying function, we obtain using \eqref{equ:bootstrap_assumption_z}, \eqref{equ:bootstrap_assumption_Linfty_f}, \eqref{equ:bootstrap_assumption_HN_f}, \eqref{equ:bootstrap_assumption_HN_ILED_f} as well as \eqref{equ:GNS} that 
    \begin{equation}     
        \begin{aligned}
            &\bigl\| \jx^{1+\kappa} \jD^{N-1} \bigl( \bmY \partial_i u_{k,\real}(t) z_j(t) \bigr) \bigr\|_{L^2_t([0,T]; L^2_x)} \\
            &\lesssim \bigl\| \jx^{\frac{5}{2}+2\kappa} \bmY \bigr\|_{W^{N-1,\infty}_x} \sum_{1 \leq k \leq N} \bigl\| \jx^{-\frac32-\kappa} |D|^k e^{it\jD} \bmf_\real(t) \bigr\|_{L^2_t([0,T]; L^2_x(\bbR^2))} \bigl\| \frakz(t)^{\frac12} \bigr\|_{L^\infty_t([0,T]; L^\infty_x)} \\
            &\quad + \bigl\| \jx^{1+\kappa} \bmY \bigr\|_{H^{N-1}_x} \Bigl\| \bigl\| \jD e^{it\jD} \bmf_\real(t) \bigr\|_{L^\infty_x} \frakz(t)^{\frac12} \Bigr\|_{L^2_t([0,T])} \\ 
            &\lesssim \varepsilon^2 + \varepsilon^{2-10\delta} \lesssim \varepsilon^{2-10\delta}.
        \end{aligned}
    \end{equation}    
    Finally, we consider the contribution from $\partial_j \pt \eta_{0,d}$ for $j \in \{1,2\}$. Inserting the equation \eqref{equ:pteta0d_equation} for $\pt \eta_{0,d}$, we find that
    \begin{equation} 
        \begin{aligned}
            \bigl\| \jx^{1+\kappa} \jD^{N-1} \bigl( \partial_j \pt \eta_{0,d}(t) \bigr) \bigr\|_{L^2_t([0,T]; L^2_x)} &\lesssim \bigl\| \jx^{1+\kappa} \jD^N (-\Delta+U^2)^{-1} \bigl( \text{RHS of \eqref{equ:pteta0d_equation}} \bigr) \bigr\|_{L^2_t([0,T];L^2_x)}.
        \end{aligned}
    \end{equation}
    Many of the terms on the right-hand side of \eqref{equ:pteta0d_equation} are higher order nonlinearities and easy to treat. Here we only discuss the tightest quadratic terms, beginning with the term $\bigl( \cos(\theta) Y_1 \bigr) z_1 (-\Delta + V_1 + 1) u_{3,\real}$. Using \eqref{equ:bootstrap_assumption_z}, \eqref{equ:bootstrap_assumption_Linfty_f}, \eqref{equ:bootstrap_assumption_HN_ILED_f} we obtain 
    \begin{equation}
        \begin{aligned}
            &\Bigl\| \jx^{1+\kappa} \jD^N (-\Delta+U^2)^{-1} \Bigl( \bigl( \cos(\theta) Y_1 \bigr) z_1 (-\Delta + V_1 + 1) u_{3,\real} \Bigr) \Bigr\|_{L^2_t([0,T];L^2_x)} \\
            &\lesssim \bigl\| \jx^{\frac{5}{2}+2\kappa} \bigl( \cos(\theta) Y_1 \bigr) \bigr\|_{W^{N-2,\infty}_x} \sum_{1 \leq |\gamma| \leq N} \bigl\| \jx^{-\frac{3}{2}-\kappa} \partial^\gamma e^{it\jD} \bmf_\real(t) \bigr\|_{L^2_t([0,T];L^2_x)} \bigl\| \frakz(t)^{\frac12} \bigr\|_{L^\infty_t([0,T]; L^\infty_x)} \\ 
            &\quad + \bigl\| \jx^{1+\kappa} \bigl( \cos(\theta) Y_1 \bigr) \bigr\|_{W^{N-2,\infty}_x} \Bigl\| \bigl\| e^{it\jD} \bmf_\real(t) \bigr\|_{L^\infty_x} \frakz(t)^{\frac12} \Bigr\|_{L^2_t([0,T])} \\ 
            &\lesssim \varepsilon^2 + \varepsilon^{2-10\delta} \lesssim \varepsilon^{2-10\delta}.
        \end{aligned}
    \end{equation}
    Next, we consider the contribution of the term $-2U (\pt \eta_{0,d}) (z_1 Y_1 + u_1)$. Here we have to reinsert the equation \eqref{equ:pteta0d_equation} for $\pt \eta_{0,d}$ once more, at which point all contributions are at least of cubic order and thus harmless to bound. 
    Finally, we discuss the contribution of the term $-2U \bigl(\pt \eta_{0,c}^\flat + \pt \eta_{0,c}^\sharp\bigr) z_1 Y_1$. 
    Using \eqref{equ:bootstrap_assumption_z}, \eqref{equ:pteta0c_HN_bound}, and \eqref{equ:pteta0c_sharp_Hk_bound} we conclude that
    \begin{equation}
        \begin{aligned}
            &\Bigl\| \jx^{1+\kappa} \jD^N (-\Delta+U^2)^{-1} \Bigl( -2U \bigl(\pt \eta_{0,c}^\flat + \pt \eta_{0,c}^\sharp\bigr) z_1 Y_1 \Bigr) \Bigr\|_{L^2_t([0,T];L^2_x)} \\
            &\lesssim \bigl\| \jx^{1+\kappa} U Y_1 \bigr\|_{W^{N-2,\infty}_x} \Bigl\| \bigl\| \jD^{N-2} (\pt \eta_{0,c}^\flat)(t) \bigr\|_{L^2_x} \frakz(t)^{\frac12} \Bigr\|_{L^2_t([0,T])} \lesssim \varepsilon^{3-10\delta}.
        \end{aligned}
    \end{equation}    
    The second terms on the right-hand sides of \eqref{equ:HN_h_decomposition} and of \eqref{equ:ILED_h_decomposition} can be estimated similarly as above, leading to the same bounds.

    Overall, we find that
    \begin{equation}
        \sup_{0 \leq t \leq T} \, \bigl\| \bmh_\ast(t)\bigr\|_{H^N_x(\bbR^2)} + \sum_{1 \leq k \leq N} \bigl\| \jx^{-\frac32-\kappa} |D|^k e^{it\jD} \bmh_\ast(t)\bigr\|_{L^2_t([0,T]; L^2_x(\bbR^2))} \leq C_1 \varepsilon + C_2 \varepsilon^{2-10\delta},
    \end{equation}
    with the size of the constant $C_2 \equiv C_2(C_0)$ depending on the size of the bootstrap constant $C_0$. 
    Since $C_0 \geq 1$ will be chosen sufficiently large depending only on some absolute constants, we can choose the small constant $0 < \varepsilon_0 \ll 1$ sufficiently small depending on the size of $C_0$ so that $C_2 \varepsilon^{2-10\delta} \leq C_1 \varepsilon$. This finishes the proof of Proposition~\ref{prop:HN_h}.
\end{proof}

Finally, we derive weighted energy estimates for the profiles $\bmh_\ast(t)$, $\ast \in \{\real,\imag\}$.

\begin{proposition} \label{prop:weighted_energies_h}
    Suppose the bootstrap assumptions \eqref{equ:bootstrap_assumption_z}, \eqref{equ:bootstrap_assumption_weighted_f}, \eqref{equ:bootstrap_assumption_Linfty_f}, \eqref{equ:bootstrap_assumption_HN_f}, \eqref{equ:bootstrap_assumption_HN_ILED_f} are in place.  
    Then we have for $\ast \in \{\real, \imag\}$,
    \begin{equation} \label{equ:weighted_energies_h_asserted_bound}
        \bigl\| \jxi^2 \jap{\nabla_\xi} \widehat{\bmh}_{\ast}(t,\xi) \bigr\|_{L^2_\xi(\bbR^2)} \leq \left\{ \begin{aligned}
                                &\max\bigl\{ C_1 \varepsilon^2 t^{\frac32}, C_2 \varepsilon^{2-11\delta} \jt^{1+\delta} \bigr\}, \quad &&0 \leq t \leq \varepsilon^{-2}, \\
                                &C_1 t^{\frac12},  &&\varepsilon^{-2} \leq t \leq T.
                              \end{aligned} \right.  
    \end{equation}
    where $C_1 \geq 1$ is some absolute constant and where $C_2 \equiv C_2(C_0)$ depends on the size of the constant $C_0$ in the bootstrap assumptions \eqref{equ:bootstrap_assumption_weighted_f}, \eqref{equ:bootstrap_assumption_Linfty_f}, \eqref{equ:bootstrap_assumption_HN_f}, \eqref{equ:bootstrap_assumption_HN_ILED_f}.
\end{proposition}
\begin{proof}
    We have 
    \begin{equation}
        \bigl\| \jxi^2 \jap{\nabla_\xi} \widehat{\bmh}_{\ast}(t,\xi) \bigr\|_{L^2_\xi(\bbR^2)} \lesssim \bigl\| \jxi^2 \widehat{\bmh}_{\ast}(t,\xi) \bigr\|_{L^2_\xi(\bbR^2)} + \bigl\| \jxi^2 \nabla_\xi \widehat{\bmh}_{\ast}(t,\xi) \bigr\|_{L^2_\xi(\bbR^2)}.
    \end{equation}
    In what follows we just concentrate on estimating the second term on the right-hand side, since the first term is strictly easier to treat.
    Using Proposition~\ref{prop:transference} and the evolution equation \eqref{equ:evol_equ_bmw} for the bad part $\bmw(t)$, we find that 
    \begin{equation} \label{equ:weighted_energies_h_black_box_applied}
        \begin{aligned}
            &\sum_{\ast \in \{\real, \imag\}} \bigl\| \jxi^2 \nabla_\xi \widehat{\bmh}_\ast(t,\xi) \bigr\|_{L^2_\xi} \\
            &\quad \quad \lesssim \bigl\| \jx \bigl( e^{i\theta} \bfP_c \bmN_d(0,r) \bigr) \bigr\|_{L^2_x} + \bigl\| \jx \bigl( e^{i\theta} \bfP_c (\pt \bmN_d)(0,r) \bigr) \bigr\|_{L^2_x} \\ 
            &\quad \quad \quad + \bigl\| \jx \bigl( e^{i\theta} \bfP_c \bmN_d(t,r) \bigr) \bigr\|_{L^2_x} + \sum_{\ast \in \{\real, \imag\}} \Bigl( \bigl\| \bmh_\ast(t) \bigr\|_{H^2_x} + \jt \bigl\| \jD^2 e^{it\jD} \bmh_\ast(t) \bigr\|_{L^{\infty-}_x} \Bigr) \\ 
            &\quad \quad \quad + \sum_{l=0,1} \bigl\| s \cdot \jx \bigl( e^{i\theta} \bfP_c (\partial_s^l \bmN_d)(s,r) \bigr) \bigr\|_{L^2_s([0,t]; L^2_x)} + \bigl\| s \cdot \jx \jD \bigl( e^{i\theta} \bfP_c \bmN_d(s,r) \bigr) \bigr\|_{L^2_s([0,t]; L^2_x)} \\
            &\quad \quad \quad + \sum_{l=1,2} \bigl\| \jap{x} \bigl( e^{i\theta} \bfP_c (\partial_s^l \bmN_d)(s,r) \bigr) \bigr\|_{L^1_s([0,t]; L^2_x)} + \bigl\| \jx \jD \bigl( e^{i\theta} \bfP_c \bmN_d(s,r) \bigr) \bigr\|_{L^1_s([0,t]; L^2_x)} \\ 
            &=: I + II + III + IV + V + VI + VII + VIII.
        \end{aligned}
    \end{equation}
    We now systematically estimate all terms on the right-hand side of \eqref{equ:weighted_energies_h_black_box_applied}.
    Recall that all the nonlinearities in $e^{i\theta} \bmN_d$ are spatially localized. The most problematic contributions stem from those quadratic terms, where both inputs are given by an internal mode. In what follows, we treat the representative contribution of the schematic form $\bmY z_i z_j$ with $i,j \in \{1,2\}$ and $\bmY(x)$ some smooth, rapidly decaying function. As in the proofs of Proposition~\ref{prop:dispersive_decay_h} and of Proposition~\ref{prop:HN_h}, we also treat the other representative quadratic terms of the schematic form $\bmY \partial_i u_{k,\real} z_j$ with $k \in \{1,2,3,4\}$, $i, j \in \{1,2\}$, and $\bmY(x)$ some smooth, rapidly decaying function as well as the contributions from $\partial_j \partial_t \eta_{0,d}$. All other nonlinearities in $e^{i\theta} \bmN_d$ are similar or of higher order, and thus easier to treat. 

    We begin with the first term $I$ on the right-hand side of \eqref{equ:weighted_energies_h_black_box_applied}. For the contribution of the  quadratic nonlinearity $\bmY z_i z_j$ to $I$, we have by \eqref{equ:new_assumption1} that
    \begin{equation}
        \begin{aligned}
            \bigl\| \jx \bigl( e^{i\theta} \bfP_c \bmY z_i(0) z_j(0) \bigr) \bigr\|_{L^2_x} \lesssim \bigl\| \jx \bmY \bigr\|_{L^2_x} \frakz(0) \lesssim \varepsilon^2.
        \end{aligned}
    \end{equation}
    Using \eqref{equ:new_assumption1} and \eqref{equ:new_assumption2}, the contribution of the quadratic nonlinearity $\bmY u_{k,\real} z_j$ to the first term~$I$ can be estimated by
    \begin{equation}
        \begin{aligned}
            \bigl\| \jx e^{i\theta} \bfP_c \bigl( \bmY \partial_i u_{k,\real}(0) z_j(0) \bigr) \bigr\|_{L^2_x} \lesssim \bigl\| \jx \bmY \bigr\|_{L^2_x \cap L^\infty_x} \bigl\|e^{i\theta}\bmu(0)\bigr\|_{H^1_x} \frakz(0)^{\frac12} \lesssim \varepsilon^2.
        \end{aligned}
    \end{equation}
    Finally, invoking the estimate \eqref{equ:pteta0d_Linfty_bound} from Proposition~\ref{prop:eta0d_bounds}, we bound the contribution from $\partial_j \partial_t \eta_{0,d}$ to the first term $I$ by 
    \begin{equation}
        \begin{aligned}
            \bigl\| \jx e^{i\theta} \bfP_c \bigl( \partial_j \partial_t \eta_{0,d}(0) \bigr) \bigr\|_{L^2_x} &\lesssim \bigl\| \jx^3 \jD (\pt \eta_{0,d})(0) \bigr\|_{L^\infty_x} \lesssim \varepsilon^{2-10\delta}.
        \end{aligned}
    \end{equation}

    The second term $II$ on the right-hand side of \eqref{equ:weighted_energies_h_black_box_applied} enjoys the same bounds as the first term~$I$.

    Now we consider the third term $III$ on the right-hand side of \eqref{equ:weighted_energies_h_black_box_applied}. 
    For the contribution of the quadratic nonlinearity $\bmY z_i z_j$ to $III$, we have by \eqref{equ:bootstrap_assumption_z} for all times $0 \leq t \leq T$ that
    \begin{equation}
        \begin{aligned}
            \bigl\| \jx \bigl( e^{i\theta} \bfP_c \bmY z_i(t) z_j(t) \bigr) \bigr\|_{L^2_x} \lesssim \bigl\| \jx \bmY \bigr\|_{L^2_x} \frakz(t) \lesssim \varepsilon^2.
        \end{aligned}
    \end{equation}
    Using \eqref{equ:bootstrap_assumption_z} and \eqref{equ:bootstrap_assumption_HN_f}, the contribution of the quadratic nonlinearity $\bmY u_{k,\real} z_j$ to $III$ can be estimated for all times $0 \leq t \leq T$ by
    \begin{equation}
        \begin{aligned}
            \bigl\| \jx e^{i\theta} \bfP_c \bigl( \bmY \partial_i u_{k,\real}(t) z_j(t) \bigr) \bigr\|_{L^2_x} \lesssim \bigl\| \jx \bmY \bigr\|_{L^2_x \cap L^\infty_x} \bigl\|\bmf_\real(t)\|_{H^1_x} \frakz(t)^{\frac12} \lesssim \varepsilon^2.
        \end{aligned}
    \end{equation}
    Finally, invoking the estimate \eqref{equ:pteta0d_Linfty_bound} from Proposition~\ref{prop:eta0d_bounds}, we bound the contribution from $\partial_j \partial_t \eta_{0,d}$ to $III$ for all times $0 \leq t \leq T$ by 
    \begin{equation}
        \begin{aligned}
            \bigl\| \jx e^{i\theta} \bfP_c \bigl( \partial_j \partial_t \eta_{0,d}(t) \bigr) \bigr\|_{L^2_x} &\lesssim \bigl\| \jx^3 \jD (\pt \eta_{0,d})(t) \bigr\|_{L^\infty_x} \lesssim \varepsilon^{2-10\delta}.
        \end{aligned}
    \end{equation}

    We continue with the fourth term $IV$ on the right-hand side of \eqref{equ:weighted_energies_h_black_box_applied}. 
    By Proposition~\ref{prop:HN_h} we have $\bigl\| \bmh_\ast(t) \bigr\|_{H^2_x} \leq C_1 \varepsilon$ for some absolute constant $C_1 \geq 1$ uniformly for all times $0 \leq t \leq T$ for $\ast \in \{\real,\imag\}$. Moreover, by \eqref{equ:GNS} and the bounds from Proposition~\ref{prop:dispersive_decay_h} and Proposition~\ref{prop:HN_h}, we find that
    \begin{equation}
        \begin{aligned}
            \sum_{\ast \in \{\real, \imag\}} \, \jt \bigl\| \jD^2 e^{it\jD} \bmh_\ast(t) \bigr\|_{L^{\infty-}_x} &\lesssim \sum_{\ast \in \{\real, \imag\}} \, \jt \bigl\| \bmh_\ast(t) \bigr\|_{H^N_x}^{\frac{2-\frac{2}{\infty-}}{N-1}} \bigl\| e^{it\jD} \bmh_\ast(t) \bigr\|_{L^\infty_x}^{1-\frac{2-\frac{2}{\infty-}}{N-1}} \\
            &\lesssim \left\{ \begin{aligned}
                            &\max\bigl\{ \varepsilon^{2-11\delta} \jt^{\frac12-\delta}, \varepsilon^{2-\delta} \jt^{1+\delta} \bigr\}, \quad &&0 \leq t \leq \varepsilon^{-2}, \\
                            &t^{\delta},  &&\varepsilon^{-2} \leq t \leq T.
                         \end{aligned} \right.  
        \end{aligned}
    \end{equation}

    Next, we turn to the fifth term $V$ on the right-hand side of \eqref{equ:weighted_energies_h_black_box_applied}. 
    For the contribution of the quadratic nonlinearity $\bmY z_i z_j$ to the term $V$ we infer using \eqref{equ:bootstrap_assumption_z} that
    \begin{equation}
        \begin{aligned}
            \sum_{j=0,1} \bigl\| s \cdot \jx \partial_s^j \bigl( \bmY(x) z_i(s) z_j(s) \bigr) \bigr\|_{L^2_s([0,t]; L^2_x)} &\lesssim \bigl\| \jx \bmY(x) \bigr\|_{L^2_x} \bigl\| s \cdot \frakz(s) \bigr\|_{L^2_s([0,t])} \\
            &\leq C_1 \left\{ \begin{aligned}
                                &\varepsilon^2 t^{\frac32}, \quad &&0 \leq t \leq \varepsilon^{-2}, \\
                                &t^{\frac12},  &&\varepsilon^{-2} \leq t \leq T,
                              \end{aligned} \right.  
        \end{aligned}
    \end{equation}
    for some absolute constant $C_1 \geq 1$.    
    To estimate the contribution of the schematic quadratic nonlinearity $\bmY \partial_i u_{k,\real} z_j$ to the term $V$, we invoke  \eqref{equ:bootstrap_assumption_z}, \eqref{equ:bootstrap_assumption_Linfty_f}, \eqref{equ:bootstrap_assumption_HN_f} along with the interpolation inequality \eqref{equ:GNS}. This gives    
    \begin{equation}
        \begin{aligned}
            \sum_{l=0,1} \bigl\| s \cdot \jx \partial_s^l \bigl( \bmY(x) \partial_i u_{k,\real}(s) z_j(s) \bigr) \bigr\|_{L^2_s([0,t]; L^2_x)} 
            &\lesssim \bigl\| \jx \bmY(x) \bigr\|_{L^2_x} \Bigl\| s \cdot \bigl\| \jD^2 e^{is\jD} \bmf_\real(s) \bigr\|_{L^\infty_x} \frakz(s)^{\frac12} \Bigr\|_{L^2_s([0,t])} \\
            &\lesssim \left\{ \begin{aligned}
                                &\varepsilon^{2-10\delta} t^{1-3\delta}, \quad &&0 \leq t \leq \varepsilon^{-2}, \\
                                &t^{2\delta},  &&\varepsilon^{-2} \leq t \leq T.
                              \end{aligned} \right.  
        \end{aligned}
    \end{equation}
    Finally, to bound the contribution of $\partial_j \partial_t \eta_{0,d}$ to the term $V$, we invoke \eqref{equ:pteta0d_Linfty_bound} to conclude that
    \begin{equation}
        \begin{aligned}
            \sum_{l=0,1} \bigl\| s \cdot \jx \partial_j \partial_s^{1+l} \eta_{0,d}(s) \bigr\|_{L^2_s([0,t]; L^2_x)} &\lesssim \sum_{l=0,1} \, \Bigl\| s \cdot \bigl\| \jx^3 \jD \partial_s^{1+l} \eta_{0,d}(s) \bigr\|_{L^\infty_x} \Bigr\|_{L^2_s([0,t])} \\
            &\lesssim \left\{ \begin{aligned}
                                &\varepsilon^{2-10\delta} t^{1-3\delta}, \quad &&0 \leq t \leq \varepsilon^{-2}, \\
                                &t^{2\delta},  &&\varepsilon^{-2} \leq t \leq T.
                              \end{aligned} \right.  
        \end{aligned}
    \end{equation}

    It is clear that the sixth term $VI$ on the right-hand side of \eqref{equ:weighted_energies_h_black_box_applied} enjoys the same bounds as the preceding fifth term $V$. Correspondingly, we now consider the seventh term $VII$ on the right-hand side of \eqref{equ:weighted_energies_h_black_box_applied}.
    For the contribution of the quadratic nonlinearity $\bmY z_i z_j$ to the term $VII$, we have in the case $l = 1$ by \eqref{equ:bootstrap_assumption_z} that for $0 \leq t \leq T$,
    \begin{equation} \label{equ:weighted_energies_h_termVII_zsquared}
        \begin{aligned}
            &\bigl\| \jap{x} \partial_s \bigl( \bmY z_i z_j \bigr) \bigr\|_{L^1_s([0,t]; L^2_x)} \lesssim \bigl\| \jx \bmY \bigr\|_{L^2_x} \int_0^t \frakz(s) \, \ud s \lesssim \log\bigl(1+\varepsilon^2 t\bigr).
        \end{aligned}
    \end{equation}
    In the case $l=2$ we also have to reinsert the evolution equations \eqref{equ:system_evol_equations_u_z} for the internal modes $\partial_s^2 z_j = -\lambda^2 z_j + \langle \bmY_j, \bmN \rangle$. Since the nonlinearity $\langle \bmY_j, \bmN \rangle$ is of higher order this produces similar or better terms. Thus, in the case $l=2$ we obtain the same bound \eqref{equ:weighted_energies_h_termVII_zsquared} as for the case $l=1$. 
    Next, we discuss the contribution of the quadratic nonlinearity $\bmY \partial_i u_{k,\real} z_j$ to the term $VII$. 
    In the case $l=1$ we use \eqref{equ:bootstrap_assumption_z}, \eqref{equ:bootstrap_assumption_Linfty_f}, \eqref{equ:bootstrap_assumption_HN_f} and the interpolation inequality \eqref{equ:GNS} to infer
    \begin{equation} \label{equ:weighted_energies_h_termVII_uz}
        \begin{aligned}
            \bigl\| \jap{x} \partial_s \bigl( \bmY \partial_i u_{k,\real} z_j \bigr) \bigr\|_{L^1_s([0,t]; L^2_x)} &\lesssim \bigl\| \jx \bmY \bigr\|_{L^2_x} \int_0^t \bigl\| \jD^2 e^{is\jD} \bmf_\real(s) \bigr\|_{L^\infty_x} \frakz(s)^{\frac12} \, \ud s \\
            &\lesssim \left\{ \begin{aligned}
                                &\varepsilon^{2-10\delta} t^{\frac12-3\delta}, \quad &&0 \leq t \leq \varepsilon^{-2}, \\
                                &\varepsilon^{1-4\delta},  &&\varepsilon^{-2} \leq t \leq T.
                              \end{aligned} \right.  
        \end{aligned}
    \end{equation}
    In the case $l=2$ we again also have to insert the evolution equation $\partial_s^2 \bmu_\real = - \bfM_{\mathrm{nr}} \bmu_\real + \Re\bigl( e^{i\theta} \bfP_c \bmN \bigr)$ obtained from \eqref{equ:system_evol_equations_u_z} for the component $u_{k,\real}$ of the radiation term as well as the evolution equation for the internal mode. Since the nonlinearities give higher order contributions, this again only yields similar or better terms, whence in the case $l=2$ we obtain the same bound \eqref{equ:weighted_energies_h_termVII_uz} as for the case $l=1$. 
    Finally, we estimate the contribution of $\partial_j \partial_t \eta_{0,d}$ to the term $VII$. 
    By \eqref{equ:pteta0d_Linfty_bound} we obtain 
    \begin{equation}
        \begin{aligned}
            \sum_{l=1,2} \bigl\| \jx \partial_j \partial_s^{1+l} \eta_{0,d}(s) \bigr\|_{L^1_s([0,t]; L^2_x)} &\lesssim \sum_{l=1,2} \int_0^t \bigl\| \jx^3 \jD \partial_s^{1+l} \eta_{0,d}(s) \bigr\|_{L^\infty_x} \, \ud s \\ 
            &\lesssim \left\{ \begin{aligned}
                                &\varepsilon^{2-10\delta} t^{\frac12-3\delta}, \quad &&0 \leq t \leq \varepsilon^{-2}, \\
                                &\varepsilon^{1-4\delta},  &&\varepsilon^{-2} \leq t \leq T.
                              \end{aligned} \right.  
        \end{aligned}
    \end{equation}
    It is clear that the last eighth term $VIII$ on the right-hand side of   \eqref{equ:weighted_energies_h_black_box_applied} enjoys the same bounds as the preceding seventh term. We omit the details. 

    Finally, by collecting the fastest growing bounds from all of the preceding estimates, and keeping in mind that the small constant $0 < \varepsilon_0 \ll 1$ will be chosen sufficiently small depending on the size of $C_0$, we arrive at the asserted bound \eqref{equ:weighted_energies_h_asserted_bound} in the statement of Proposition~\ref{prop:weighted_energies_h}.
\end{proof}

\section{Internal Mode I: The System of ODEs and Fermi Golden Rule}\label{sec:FGR1}

In this section we first prove the constrained Hamiltonian formulation
stated in Subsection \ref{subsec:Hamiltonintro}, and then expand the perturbation
Hamiltonian to obtain the almost Hamiltonian system of
Lemma~\ref{lemHamQ}.  This structure will be used in the derivation of the
effective equations for the internal modes and in the computation of the
Fermi Golden Rule coefficients.
In Subsection \ref{secevol} we write Duhamel's formula for the profile of the perturbation
and organize the various nonlinear terms according to their homogeneity and dependence on 
the continuous or discrete components.
We are then in a position to study the system of ODEs satisfied by the internal modes.
In Subsection \ref{secIM1}, exploiting the Hamiltonian and block structure of the system, we calculate the relevant Fermi Golden Rule coefficients that provide dissipation for the internal modes' amplitudes; the necessary non-degeneracy assumption is guaranteed by the analysis and numerical work in \cite{LPPSS1}.
Moreover, we also organize all the non-leading order dissipative terms in the ODE, and provide sufficient conditions (see Lemma \ref{lemODEdamp2ndversion}) so that these can be considered as remainders
relative to the dissipative dynamics.
The estimates for all the remainder terms are then carried out in Section~\ref{secIM2}.

\subsection{Hamiltonian structure}\label{Ssec:Ham}
We begin by giving the details of the derivation of the constrained Hamiltonian formulation stated in Subsection \ref{subsec:Hamiltonintro} and
the validity of \eqref{eq:con_ham_4_alt}-\eqref{eq:con_ham_4_alt_constraint}.
We then specialize the perturbation
Hamiltonian to the equivariant variables and separate its quadratic, cubic,
quartic, and constraint contributions to arrive at the main
equations \eqref{equ:sys3}-\eqref{equ:pt_eta0_equation}

\begin{proof}[Proof of \eqref{eq:con_ham_4_alt}-\eqref{eq:con_ham_4_alt_constraint}]
Recall the position and momentum variables introduced in
\eqref{equ:intro_def_pos_phi}--\eqref{equ:intro_def_mom_phi}.  Since the
Lagrangian does not contain $\partial_0A_0$, the temporal component $A_0$
has no conjugate momentum and is treated as an independent,
non-dynamical variable.
The velocities can be expressed in terms of the momenta as
\[
    \partial_0\phi
    =
    -2\mom_\phi+iA_0\pos_\phi,
    \qquad
    \partial_0A_j
    =
    \partial_jA_0-\mom_{A_j},
    \qquad
    j=1,2.
\]
After this substitution, the Lagrangian density becomes
\begin{equation}
\begin{aligned}
    &\mathscrbf{L}
    \bigl(
    \pos_\phi,\overline{\pos}_\phi,\pos_A,
    \mom_\phi,\overline{\mom}_\phi,\mom_A
    \bigr) =
    \frac12B^2
    -\frac12|\mom_A|^2
    -2|\mom_\phi|^2
    +\frac12
    \bigl|(\nabla_x-i\pos_A)\pos_\phi\bigr|^2
    +\frac18\bigl(1-|\pos_\phi|^2\bigr)^2.
\end{aligned}
\end{equation}
Moreover,
\begin{equation}
\begin{aligned}
    &-\mom_\phi\partial_0\overline{\pos}_\phi
    -\overline{\mom}_\phi\partial_0\pos_\phi
    -\sum_{j=1}^2\mom_{A_j}\partial_0\pos_{A_j} =
    4|\mom_\phi|^2
    -2A_0\im\bigl(\overline{\pos}_\phi\mom_\phi\bigr)
    +|\mom_A|^2
    -\mom_A\cdot\nabla_xA_0.
\end{aligned}
\end{equation}
It follows that the Hamiltonian density is
\begin{align}\label{eq:Hamden}
\begin{split}
    \mathscrbf{H}
    \bigl( \pos_\phi,\overline{\pos}_\phi,\pos_A, \mom_\phi,\overline{\mom}_\phi,\mom_A,A_0
    \bigr) & =
    \frac12B^2
    +\frac12|\mom_A|^2
    +2|\mom_\phi|^2
    +\frac12
    \bigl|(\nabla_x-i\pos_A)\pos_\phi\bigr|^2
    \\ & +\frac18\bigl(1-|\pos_\phi|^2\bigr)^2 - 2A_0\im\bigl(\overline{\pos}_\phi\mom_\phi\bigr)
    -\mom_A\cdot\nabla_xA_0.
\end{split}
\end{align}
The Euler--Lagrange equations
\eqref{equ:EL_equations_selfdual_rectangular} are therefore equivalent to
the constrained Hamiltonian system
\begin{align}\label{eq:con_ham_1'}
    \partial_t\pos
    =
    -\frac{\delta\mathcal H}{\delta\overline{\mom}},
    \qquad
    \partial_t\mom
    =
    \frac{\delta\mathcal H}{\delta\overline{\pos}},
    \qquad
    \frac{\delta\mathcal H}{\delta A_0}
    =
    0,
\end{align}
where
\[
    \pos
    =
    \bigl(
    \pos_\phi,\overline{\pos}_\phi,\pos_{A_1},\pos_{A_2}
    \bigr),
    \qquad
    \mom
    =
    \bigl(
    \mom_\phi,\overline{\mom}_\phi,\mom_{A_1},\mom_{A_2}
    \bigr).
\]
The last equation in \eqref{eq:con_ham_1'} is Gauss' law, 
equation \eqref{equ:EL_equations_selfdual_rectangular}.

We next expand the Hamiltonian about the degree-one vortex.  For general
perturbations, write
\[
    \phi=\underline{\Phi}+\varphi,
    \qquad
    A_x=\underline A_x+\eta_x,
    \qquad
    A_0=\eta_0,
\]
and let $\mathcal H_0$ be the perturbation Hamiltonian defined in
\eqref{eq:Hamper}, which we recall here for convenience is given by
\begin{equation}\label{eq:Hamper'}
\begin{aligned}
    &\mathcal{H}_0\bigl[\varphi,\overline{\varphi},\eta,\mom_\varphi,\overline{\mom}_\varphi,\mom_\eta,\eta_0\bigr] \\
    &\quad := \mathcal{H}\big[\underline{\Phi} + \varphi, \overline{\underline{\Phi}} + \overline{\varphi}, \underline{A}+\eta_x,\mom_{\varphi},\overline{\mom}_{\varphi},\mom_{\eta},\eta_0 \big] - \mathcal{H}\big[ \underline{\Phi},\overline{\underline{\Phi}},\underline{A},0,0,0,0 \big],
\end{aligned}
\end{equation}
where
\[
\mom_\varphi := \mom_\phi, \qquad 
\mom_\eta := \mom_A, \]
since the background is static and $\underline A_0=0$, its conjugate momenta vanish, that is, $\mom_{\underline{\Phi}} = \mom_{\underline A_j} = 0$. 
 
Introduce the covariant derivative associated with the
fixed background,
\[
    \bfD_{\underline A}
    :=
    \nabla-i\underline A,
    \qquad
    (\bfD_{\underline A})_j
    :=
    \partial_j-i\underline A_j.
\]
Expanding \eqref{eq:Hamden} about
$(\underline{\Phi},\underline A)$ gives
\begin{align}\label{eq:calH0def1}
\begin{split}
    &\mathcal H_0
    \bigl[
    \varphi,\overline{\varphi},\eta,
    \mom_\varphi,\overline{\mom}_\varphi,\mom_\eta,\eta_0
    \bigr]
    \\
    &=
    \int_{\R^2}
    \Bigl(
    2|\mom_\varphi|^2
    +\frac12|\mom_\eta|^2
    +\frac12(\partial_1\eta_2-\partial_2\eta_1)^2
    +\frac12|\bfD_{\underline A}\varphi|^2
    +\frac12|\underline{\Phi}|^2|\eta|^2
    \Bigr)\,\ud x
    \\
    &\quad
    +
    \int_{\R^2}
    \Bigl(
    -\eta\cdot
    \im\bigl(
    \overline{\underline{\Phi}}\bfD_{\underline A}\varphi
    +
    \overline{\varphi}\bfD_{\underline A}\underline{\Phi}
    \bigr)
    -\eta\cdot
    \im\bigl(\overline{\varphi}\bfD_{\underline A}\varphi\bigr)
    +|\eta|^2\re\bigl(\overline{\underline{\Phi}}\varphi\bigr)
    \Bigr)\,\ud x
    \\
    &\quad
    +
    \int_{\R^2}
    \Bigl(
    \frac12|\eta|^2|\varphi|^2
    -\frac14(1-|\underline{\Phi}|^2)|\varphi|^2
    +\frac12
    \bigl(\re(\overline{\underline{\Phi}}\varphi)\bigr)^2
    +\frac12
    \re(\overline{\underline{\Phi}}\varphi)|\varphi|^2
    \Bigr)\,\ud x
    \\
    &\quad
    +
    \int_{\R^2}
    \Bigl(
    \frac18|\varphi|^4
    -2\eta_0
    \im\Bigl(
    (\overline{\underline{\Phi}}+\overline{\varphi})
    \mom_\varphi
    \Bigr)
    +\eta_0\nabla\cdot\mom_\eta
    \Bigr)\,\ud x.
\end{split}
\end{align}
In obtaining \eqref{eq:calH0def1}, we have used the static vortex equations
\begin{align*}
    (\bfD_{\underline A})^j
    (\bfD_{\underline A})_j\underline{\Phi}
    +
    \frac12
    (1-|\underline{\Phi}|^2)\underline{\Phi}
    =
    0, \quad \ 
    \partial_2\underline B
    -
    \im\bigl(
    \overline{\underline{\Phi}}
    (\bfD_{\underline A})_1\underline{\Phi}
    \bigr)
    =
    0, \quad \ 
    \partial_1\underline B
    +
    \im\bigl(
    \overline{\underline{\Phi}}
    (\bfD_{\underline A})_2\underline{\Phi}
    \bigr)
    &=
    0,
\end{align*}
where
\(
\underline B
=
\partial_1\underline A_2-\partial_2\underline A_1
\).
These identities cancel the terms that are linear in
$(\varphi,\eta)$.

We now impose the equivariant ansatz \eqref{eq:hamper}.  For the scalar
field, introduce the equivariant position and momentum variables by
\begin{align}\label{eq:posmomchangephi}
    \varphi
    &=
    \frac12e^{i\theta}
    \bigl(\pos_\alpha+i\pos_\beta\bigr),
    &
    \mom_\varphi
    &=
    \frac12e^{i\theta}
    \bigl(\mom_\alpha+i\mom_\beta\bigr),
    \qquad
    \pos_\alpha,\pos_\beta,\mom_\alpha,\mom_\beta\in\R.
\end{align}
For the spatial connection, define
\begin{align}\label{eq:posmomchangeA}
    \pos_\zeta
    &:=
    \frac{\eta_\theta}{r}
    =
    \zeta,
    &
    \pos_\mu
    &:=
    -\eta_r
    =
    \mu,
    &
    \mom_\zeta
    &:=
    \frac1r\mom_{A_\theta},
    &
    \mom_\mu
    &:=
    -\mom_{A_r}.
\end{align}
Since $\eta_0$ is radial, the polar components of $\mom_A$ satisfy
\[
    \mom_{A_\theta}
    =
    -\partial_tA_\theta,
    \qquad
    \mom_{A_r}
    =
    \partial_r\eta_0-\partial_tA_r.
\]

The definitions
\eqref{equ:intro_def_pos_phi}--\eqref{equ:intro_def_mom_phi},
together with \eqref{eq:hamper},
\eqref{eq:posmomchangephi}, and \eqref{eq:posmomchangeA}, give the
position and momentum identities
\eqref{eq:con_ham_4.0}--\eqref{eq:con_ham_4_alt0}.
Since the background is fixed, replacing $\mathcal H$ by
$\mathcal H_0$ does not change the corresponding variational derivatives
with respect to the perturbation variables.  Under the preceding changes
of variables, the constrained system \eqref{eq:con_ham_1'} therefore becomes
\eqref{eq:con_ham_4_alt}, while its $A_0$ equation becomes the constraint
\eqref{eq:con_ham_4_alt_constraint}.
\end{proof}

Up to this point no gauge condition has been imposed.  For the Fermi Golden
Rule analysis, we now impose Stuart's gauge and derive a more explicit form
of the perturbation Hamiltonian in the equivariant variables.  We therefore
decompose \eqref{eq:calH0def1} into its quadratic, cubic, quartic, and
constraint contributions.  Recalling
\eqref{eq:con_ham_4.0}--\eqref{eq:con_ham_4_alt0} and the polar form
\eqref{equ:intro_Stuart_gauge_polar} of Stuart's gauge, we define 
\begin{equation}\label{eq:con_ham_2}\begin{aligned}
    \mathcal{H}_1 & : =
    \dfrac{1}{2} \int_0^\infty\Big( \mom_\alpha^2 + \mom_\beta^2 + \mom_\zeta^2 + \mom_\mu^2 \Big) \, r\,\ud r
    \\ & \  + \dfrac{1}{2}\int_0^\infty \Big(   \alpha L_1\alpha + \beta L_1\beta + \zeta L_2 \zeta + \mu L_2 \mu\Big) \,r\,\ud r - \int_0^\infty 2U'(\alpha \zeta + \beta \mu) \, r\,\ud r,
    \\ \mathcal{H}_2 & := \int_0^\infty \Big( \dfrac{1}{2}U \alpha^3 + \dfrac{1}{2}U \alpha\beta^2 + U \alpha\zeta^2+U\alpha \mu^2 - b \alpha^2\zeta - b\beta^2\zeta + \alpha \mu \partial_r\beta - \beta \mu \partial_r\alpha \Big) \, r\,\ud r,
    \\ \mathcal{H}_3 & := \int_0^\infty\Big( \dfrac{1}{2}(\alpha^2+\beta^2)(\zeta^2+\mu^2)  +\dfrac{1}{8}(\alpha^2+\beta^2)^2 \Big) \, r\,\ud r,
\end{aligned}\end{equation}
and let $\mathcal{H}_{\mathcal{C}}$ denote the \emph{constraint functional} that comes from the terms $-2A_0\im(\overline{\pos}_\phi \mom_\phi) + A_0 \nabla \cdot \mom_A$ in
\eqref{eq:Hamden}, that is,
\begin{align}\label{eq:con_ham_3}
\mathcal{H}_{\mathcal{C}}:= -\int_0^\infty \bigg(  (U+\alpha)\mom_\beta - \beta \mom_\alpha + \dfrac{1}{r}\partial_r(r\mom_\mu)\bigg)\eta_0 \, r\,\ud r.
\end{align} 
We then define
\begin{align}\label{eq:calHsdef1}
\mathcal{H}_S  := \mathcal{H}_1
	+ \mathcal{H}_2+ \mathcal{H}_3+ \mathcal{H}_{\mathcal{C}}.
\end{align}
\begin{lem}\label{lemHamQ}
Suppose $(\phi,A)$ is an equivariant solution of \eqref{equ:EL_equations_selfdual} in Stuart's gauge . Decompose $(\phi,A)= (\underline{\Phi},\underline{A})+(\fy,\eta)$, where $(\underline{\Phi},\underline{A})$ denotes the vortex defined in \eqref{equ:intro_1vortex_selfdual}. Let $\bmr=(\alpha,\zeta,\beta,\mu)^T$ be as in \eqref{syscompact2} and let $\dot{\bmr}=\partial_t \bmr$.
Then, the equations for $(\bmr,\dot{\bmr})$ can be written in 
the (approximate Hamiltonian) form 
\begin{equation}\label{eq:Hamrrdot}
\partial_t\bmr= \nabla_{\dot{\bmr}}\calH_S,\qquad \partial_t\dot{\bmr}=-\nabla_{\bmr}\calH_S+   \big(-\partial_t(\eta_0\beta),\, 0, \, \partial_t(\eta_0(U+\alpha)), \, -\partial_t\partial_r\eta_0 \big)^T.
\end{equation}
Here $\calH_S$, defined in \eqref{eq:calHsdef1}, is viewed as a function of the independent variables $(\bmr,\dot\bmr,\eta_0)$ through the relations \eqref{eq:con_ham_4.0} and \eqref{eq:con_ham_4_alt0}. 
Equivalently, matching \eqref{equ:sys3}, we have the system
\begin{align}\label{eq:lemHamQeq}
\begin{split}
    \partial_t^2 \alpha + L_1 \alpha - 2U' \zeta & = - \dfrac{\delta \mathcal{H}_2}{\delta \alpha} - \dfrac{\delta \mathcal{H}_3}{\delta \alpha} - \dfrac{\delta \mathcal{H}_{\mathcal{C}}}{\delta \alpha } -\partial_t(\eta_0\beta),
    \\ \partial_t^2 \zeta + L_2 \zeta - 2U' \alpha & = - \dfrac{\delta \mathcal{H}_2}{\delta \zeta}  - \dfrac{\delta \mathcal{H}_3}{\delta \zeta}  - \dfrac{\delta \mathcal{H}_{\mathcal{C}}}{\delta \zeta } ,
    \\ \partial_t^2 \beta + L_1 \beta - 2U' \mu & =- \dfrac{\delta \mathcal{H}_2}{\delta \beta} - \dfrac{\delta \mathcal{H}_3}{\delta \beta}  - \dfrac{\delta \mathcal{H}_{\mathcal{C}}}{\delta \beta } +\partial_t\big(\eta_0(U+\alpha)\big),
    \\ \partial_t^2 \mu + L_2 \mu - 2U' \beta  & = - \dfrac{\delta \mathcal{H}_2}{\delta \mu} - \dfrac{\delta \mathcal{H}_3}{\delta \mu} - \dfrac{\delta \mathcal{H}_{\mathcal{C}}}{\delta \mu }   - \partial_t\partial_r \eta_0.
\end{split}
\end{align}
Here $\eta_0$, $\mom_\alpha$, and $\mom_\beta$ in $\calH_\calC$ are viewed as independent variables from $(\alpha,\zeta,\beta,\mu)$ in computing the variational derivatives of $\calH_\calC$.
\end{lem}

Before presenting the proof we make a remark about the notation used in the statement of Lemma~\ref{lemHamQ}.
\begin{remark}
Strictly speaking in the variational derivatives in \eqref{eq:Hamrrdot} and \eqref{eq:lemHamQeq}, the variable $\eta_0$ is regarded as an independent variable separate from $\bmr=(\alpha,\zeta,\beta,\mu)^T$.  Later on, when using \eqref{eq:lemHamQeq} we will view $\eta_0$ as a function of $\bmr$ (by means of the defining elliptic equation \eqref{equ:eta0_equation}) and view the variational derivatives as also acting on $\eta_0$. This abuse of notation is justified because $\eta_0$ appears in $\mathcal{H}_S$ only through $\mathcal{H}_\calC$, and the constraint equation is precisely $\frac{\delta \mathcal{H}_\calC}{\delta \eta_0}=0$. However, in \eqref{eq:lemHamQeq} we still view $\mom_\alpha$ and $\mom_\beta$ in $\calH_\calC$ as independent variables that do not get differentiated with respect to~$\bmr$. After computing the variational derivative we may use \eqref{eq:con_ham_4_alt0} to express them in terms of $\bmr$ and $\dot{\bmr}$.
\end{remark}

\begin{proof}[Proof of Lemma~\ref{lemHamQ}]
We seek to expand $\mathcal{H}_0$ in terms of $(\alpha,\zeta,\beta,\mu)$ and recover the system of PDEs  from Hamilton's equation with respect to  $(\alpha,\zeta,\beta,\mu)$. First note that, using Stuart's gauge, we can rewrite the curl term in $\mathcal{H}_0$ in a more convenient way, namely, \begin{align*}
\int_{\R^2} (\partial_1\eta_2-\partial_2\eta_1)^2 \, dx & = \int_{\R^2} \Big(\vert \nabla \eta_x\vert^2-(\nabla \cdot \eta_x)^2 \Big) \, dx = \int_{\R^2} \Big( \vert \nabla \eta_x\vert^2 - \big(\im(\overline{\underline{\Phi}}\varphi)\big)^2 \Big)\, dx,
\end{align*}
where $\eta_x=(\eta_1,\eta_2)^T$  denotes the spatial components of $\eta$. Similarly, to recover the constraint equation, we rewrite
\[
\nabla \cdot \mom_{\eta} = \Delta \eta_0 - \partial_t \big( \nabla \cdot \eta\big) = \Delta \eta_0 - \im( \overline{\underline{\Phi}} \partial_t\varphi).
\]
At this point it is worth recalling that \[
b := (1-a_\theta)/r, \qquad L_1 := -\Delta + \frac{(1-a_\theta)^2}{r^2} - \frac{\pr a_\theta}{r} + U^2, \qquad L_2 := -\Delta + \frac{1}{r^2} + U^2.
\]
Then, collecting the quadratic terms in $(\varphi,\eta)$ from $\mathcal{H}_0$, and expanding in terms of $(\alpha,\zeta,\beta,\mu)$, we obtain that 
\begin{align*}  
\int_{\R^2} \bigg(  \dfrac{1}{2}\vert \bfD_{\underline{A}}\varphi\vert^2 - \dfrac{1}{4}(1-\vert \underline{\Phi}\vert^2) \vert \varphi\vert^2\bigg) \, dx  & = \int_{\R^2}\bigg( \dfrac{1}{2} \alpha\Big( -\Delta\alpha+b^2 \alpha\Big)
\\ & + \dfrac{1}{2}\beta \Big(-\Delta \beta + b^2 \beta\Big) - \dfrac{1}{4} \big(1-U^2\big) (\alpha^2+\beta^2)\bigg) \, dx ,
    \\  \int_{\R^2}\bigg( \dfrac{1}{2}(\partial_1\eta_2-\partial_2\eta_1)^2  +  \dfrac{1}{2}\vert \underline{\Phi}\vert^2\vert\eta\vert^2  \bigg) \, dx & = \int_{\R^2} \bigg( \dfrac{1}{2}\mu\Big(-\Delta \mu +\dfrac{1}{r^2}\mu\Big) 
    \\ & + \dfrac{1}{2} \zeta\Big(-\Delta \zeta + \dfrac{1}{r^2}\zeta\Big) - \dfrac{1}{2}U^2\beta^2 + \dfrac{1}{2}U^2 (\mu^2+\zeta^2) \bigg) \, dx ,
\end{align*}
and 
\begin{align*}
\int_{\R^2} \bigg( - \eta \cdot \im\big( \overline{\underline{\Phi}}\bfD_{\underline{A}}\varphi  +   \overline{\varphi}\bfD_{\underline{A}}\underline{\Phi}\big) + \dfrac{1}{2}\big(\re( \overline{\underline{\Phi}}\varphi)\big)^2\bigg) \, dx 
    \\ = \int_{\R^2} \bigg( - 2 U'(\alpha\zeta+\beta\mu) + U^2 \beta^2 + \dfrac{1}{2}U^2 \alpha^2 \bigg) \, dx,
\end{align*}
having used $U'=bU$ and Stuart's gauge condition $\nabla\cdot\eta = U\beta$ in the last identity with 
\[
\im( \overline{\underline{\Phi}} \bfD_{\underline{A}}\varphi + \overline{\varphi}\bfD_{\underline{A}}\underline{\Phi}) = \nabla \im(\overline{\underline{\Phi}}\varphi) + 2 \im(\overline{\varphi}\bfD_{\underline{A}}\underline{\Phi}),
\]
which, in particular, produces the term $U^2 \beta^2$. 
Similarly, we collect cubic terms
\begin{align*}
    & \int_{\R^2} \Big(   - \eta \cdot \im\big( \overline{\varphi} \bfD_{\underline{A}}\varphi\big) + \vert \eta\vert^2 \re ( \overline{\underline{\Phi}}\varphi) + \dfrac{1}{2}\re( \overline{\underline{\Phi}}\varphi) \vert \varphi\vert^2 \Big) \, dx 
    \\ & =   \int_{\R^2} \Big( -b \alpha^2 \zeta - b \beta^2 \zeta + \alpha\mu\partial_r\beta - \beta \mu \partial_r\alpha + U\alpha\zeta^2+U\alpha\mu^2  + \dfrac{1}{2}U\alpha^3 + \dfrac{1}{2}U\alpha \beta^2\Big) \, dx , 
\end{align*}
and quartic terms \begin{align*}
    \int_{\R^2}  \Big( \dfrac{1}{2}\vert \eta\vert^2 \vert \varphi\vert^2 +  \dfrac{1}{8}\vert \varphi\vert^4 \Big) \, dx  =  \int_{\R^2}  \Big( \dfrac{1}{2}(\alpha^2+\beta^2)(\zeta^2+\mu^2)+\dfrac{1}{8} (\alpha^2+\beta^2)^2\Big) \, dx.
\end{align*}
%
Therefore using also that $r^{-1}\partial_ra_\theta=\frac{1}{2}(1-U^2)$ we can write the perturbation Hamiltonian in \eqref{eq:Hamper} as
\begin{align*}
\mathcal{H}_0
=: 2\pi\big(\mathcal{H}_1 + \mathcal{H}_2 + \mathcal{H}_{3}
    + \mathcal{H}_{\mathcal{C}}\big),
\end{align*}
where $\calH_\ast$, $\ast\in\{1,2,3,\calC\}$, are defined in \eqref{eq:con_ham_2}--\eqref{eq:con_ham_3}.

Now to prove \eqref{eq:Hamrrdot} we start from \eqref{eq:con_ham_4_alt} and change variables from $(\pos_\ast,\mom_\ast,\eta_0)$, $\ast\in\{\alpha,\beta,\zeta,\mu\}$, to $(\bmr,\dot{\bmr},\eta_0)$. It follows from \eqref{eq:con_ham_4.0} and \eqref{eq:con_ham_4_alt0} that the variational derivatives transform according to
\begin{align*}
&\frac{\delta}{\delta \pos_\alpha}\mapsto \frac{1}{2}\frac{\delta}{\delta \alpha}+\frac{1}{2}\eta_0\frac{\delta}{\delta \dot\beta},\qquad \frac{\delta}{\delta \pos_\beta}\mapsto \frac{1}{2}\frac{\delta}{\delta \beta}-\frac{1}{2}\eta_0\frac{\delta}{\delta \dot\alpha},\quad \frac{\delta}{\delta \pos_\zeta}\mapsto \frac{\delta}{\delta \zeta},\quad {\frac{\delta}{\delta \pos_\mu}\mapsto \frac{\delta}{\delta \mu}},\\
&\frac{\delta}{\delta \eta_0}\mapsto \frac{\delta}{\delta \eta_0}-{\beta}\frac{\delta}{\delta \dot\alpha}+(\alpha+U)\frac{\delta}{\delta \dot\beta}{-\partial^\ast} \frac{\delta}{\delta \dot\mu},\qquad \frac{\delta}{\delta \mom_\ast}=-\frac{\delta}{\delta \dot{\ast}},~\ast\in\{\alpha,\beta,\zeta,\mu\}.
\end{align*}
Then \eqref{eq:Hamrrdot} follows from \eqref{eq:con_ham_4_alt} and the observation that
\begin{equation}
\big(\frac{1}{2}\frac{\delta}{\delta \alpha}+\frac{1}{2}\eta_0\frac{\delta}{\delta \dot\beta}\big)\mom_\ast=\big(\frac{1}{2}\frac{\delta}{\delta \beta}-\frac{1}{2}\eta_0\frac{\delta}{\delta \dot\alpha}\big)\mom_\ast=\big(\frac{\delta}{\delta \eta_0}-\frac{1}{2}\frac{\delta}{\delta \dot\alpha}+(\alpha+U)\frac{\delta}{\delta \dot\beta}+\partial^\ast \frac{\delta}{\delta \dot\mu}\big)\mom_\ast=0
\end{equation}
for $\ast\in\{\alpha,\beta,\zeta,\mu\}$.

For \eqref{eq:lemHamQeq} note that in view of \eqref{eq:con_ham_4_alt0}, 
\begin{equation}\label{eq:con_ham_4}
     \partial_t^2 \alpha = -\partial_t \mom_\alpha - \partial_t(\eta_0\beta),
    \quad \partial_t^2 \zeta = - \partial_t \mom_\zeta ,
    \quad \partial_t^2 \beta = -\partial_t \mom_\beta + \partial_t(\eta_0(U+\alpha)),
    \quad \partial_t^2 \mu =    - \partial_t \mom_\mu - \partial_t\partial_r\eta_0.
\end{equation}
The equations \eqref{eq:lemHamQeq}  then follow by expanding the left-hand side of the second line of \eqref{eq:con_ham_4_alt} using the relations \eqref{eq:con_ham_4_alt0}, and replacing the variational derivatives  with respect to $\mom_\ast$ on the right-hand side by the ones with respect to $\bmr$.
\end{proof}





\medskip
\subsection{Evolution equation and decomposition of the profile}\label{secevol}
In this section we rewrite the equation and organize the nonlinear terms in a more structured way, in order to streamline the analysis in the subsequent sections. We begin by defining the profile of $\bmu$, defined in \eqref{equ:spectral_decomposition}, with respect to $\bfM$: 
\begin{equation}\label{defprof}
\bfrf(t) := e^{-it\bfsM} (2i \bfsM )^{-1} (\partial_t + i\bfsM) \bmu(t),  
\qquad  \bmu  = e^{it\bfsM} \bfrf + e^{-it\bfsM} \overline{\bfrf}
\end{equation}
Here in what follows we abuse notation by writing only $\bfsM$
in place of $\bfP_c \bfsM$. Recall from \eqref{eq:Zjprofiles1} the definitions
\begin{equation}\label{defZj'}
    Z_j(t)  := e^{-i\lambda t} (2i\lambda)^{-1} (\pt + i \lambda) z_j,
\qquad  z_j  = e^{i \lambda t} Z_j + e^{-i\lambda t} \overline{Z_j}.
\end{equation}
In view of \eqref{equ:system_evol_equations_u_z}, the profiles $\bfrf$ and $Z_j$ satisfy the equations
\begin{align}\label{eq:systemZj0}
\begin{split}
\partial_t \bfrf = e^{-it\bfsM} (2i \bfsM)^{-1} \bfP_c \bmN,
\qquad 
\partial_t Z_j = e^{-i\lambda t} (2i\lambda)^{-1} \langle \bmN, \bmY_j\rangle,
\end{split}
\end{align}
for $j = 1, 2$. 
Our goal is to arrange the nonlinear terms in Duhamel's formula for $\bfrf$ so to isolate the main singular contribution in the PDE coming from the feedback of the internal mode.
This will also be useful later on in Sections \ref{secIM1} and \ref{secIM2}.
Before that, 
we arrange the nonlinear terms in the system \eqref{eq:lemHamQeq} (or,  equivalently \eqref{equ:sys3} and \eqref{sysN1}),
into quadratic, cubic, and higher order.
Starting with $\eta_0$ we decompose (to be precise $\eta_0^{(j)}$ below are functions of $\bmr$ and $\dot\bmr$ not just $\bmr$)
\begin{equation}
    \eta_0=\eta_{0}^{(2)}(\bmr)+\eta_{0}^{(3)}(\bmr)+\eta_{0}^{(4)}(\bmr).
\end{equation}
Here, for $\ubmr=(\ualpha,\uzeta,\ubeta,\umu)^T$, $\eta_{0}^{(j)}(\ubmr)$ are defined as the solutions of the elliptic equations
\begin{align}
    &(-\Delta+U^2)\eta_{0}^{(2)}(\ubmr)= -\ubeta \pt \ualpha + \ualpha \pt \ubeta,\\
    &(-\Delta+U^2)\eta_{0}^{(3)}(\ubmr)=- 2 U \eta_{0}^{(2)}(\ubmr) \ualpha,\\
    &(-\Delta+U^2)\eta_{0}^{(4)}(\ubmr)=- 2 U (\eta_{0}^{(3)}(\ubmr)+\eta_{0}^{(4)}(\ubmr)) \ualpha  - (\eta_{0}^{(2)}(\ubmr)+\eta_{0}^{(3)}(\ubmr)+\eta_{0}^{(4)}(\ubmr)) (\ualpha^2 + \ubeta^2).
\end{align}
Note that the nonlinear elliptic equation defining $\eta_0^{(4)}$ has a unique $L^2$ solution as long as $e^{i\theta}\ubmr$ is small in say $H^2_x(\bbR^2)$ (which will be the case for $\ubmr=\bmr$ in our applications). When there is no risk of confusion, we simply write $\eta_{0}^{(j)}$ for $\eta_{0}^{(j)}(\bmr)$ and $\ueta_0^{(j)}$ for $\eta_0^{(j)}(\ubmr)$. We will also write $\eta_0$ and $\ueta_0$ for $\sum_{j=2}^{4}\eta_{0}^{(j)}$ and $\sum_{j=2}^{4}\ueta_{0}^{(j)}$ respectively. For $\partial_t\eta_0$ we repeat the same procedure using \eqref{equ:pt_eta0_equation}. Substituting the expressions for $\calN_\alpha$ and $\calN_\beta$ from \eqref{sysN1} in $\calC_0$ in \eqref{equ:pt_eta0_equation} we write
\begin{equation}
    \calC_0=-U\alpha\partial_t\eta_0+\calC_{3}(\bmr)+\calC_4(\bmr),
\end{equation}
where
\begin{align}
    &\calC_3(\ubmr):=-2U\ueta_0^{(2)}\partial_t\ualpha-\frac{1}{2}U\ubeta(\ualpha^2+\ubeta^2)+U\ubeta(\uzeta^2+\umu^2)+2\umu(\ubeta\partial_r\ubeta+\ualpha\partial_r\ualpha),\\
    &\calC_{4}(\ubmr):=-2U(\ueta_0^{(3)}+\ueta_0^{(4)})\partial_t\ualpha
    -U\eta_0^2\ubeta.
\end{align}
With this notation we decompose $\partial_t\eta_0$ as
\begin{equation}\label{eq:Theta0def1}
    \partial_t\eta_0=\Theta_0^{(2)}(\bmr)+\Theta_0^{(3)}(\bmr)+\Theta_0^{(4)}(\bmr),
\end{equation}
where with $\ubmr$ as above, $\Theta_0^{(j)}(\ubmr)$ are defined as solutions of the elliptic equations
\begin{align}\label{eq:Theta0jdef1}
\begin{split}
    &(-\Delta+U^2)\Theta_0^{(2)}(\ubmr)=\ubeta L_1\ualpha-2U'\ubeta\uzeta-\ualpha L_1\ubeta+2U'\ualpha\umu,\\
    &(-\Delta+U^2)\Theta_0^{(3)}(\ubmr)=-U\ualpha\Theta_0^{(2)}(\ubmr)+\calC_3(\ubmr),\\
    &(-\Delta+U^2)\Theta_0^{(4)}(\ubmr)=-U\ualpha (\Theta_0^{(3)}(\ubmr)+\Theta_0^{(4)}(\ubmr))+\calC_4(\ubmr).
\end{split}
\end{align}
 Note that to be precise $\Theta_0^{(3)}$ and $\Theta_0^{(4)}$ are functions of $\ubmr$ and $\partial_t\ubmr$ not just $\ubmr$. When there is no risk of confusion we will use the notation $\Theta_0^{(j)}=\Theta_0^{(j)}(\bmr)$, $\uTheta_0^{(j)}=\Theta_0^{(j)}(\ubmr)$, as well as $\Theta_0:=\sum_{j=2}^{(4)}\Theta_0^{(j)}$ and $\uTheta_0:=\sum_{j=2}^{(4)}\uTheta_0^{(j)}$. 
 
 Next turning to the nonlinearity $\bmN$ we write
\begin{align}\label{splitN}
\bmN = \bmQ(\bmr) + \bmC(\bmr) + \bmN^{(4)}(\bmr)  
\end{align}
where we define the quadratic terms by
\begin{align}\label{Nrquad}
\begin{split}
\bmQ(\ubmr) & = (Q_1,Q_2,Q_3,Q_4)^T(\ubmr),
\\
Q_1(\ubmr) & = - \frac{3}{2} U \ualpha^2 + \frac12 U \ubeta ^2 
  - 2 \umu \pr \ubeta    + 2 b \ualpha \uzeta - U (\uzeta^2 + \umu ^2),
\\
Q_2(\ubmr) & = b (\ualpha^2 + \ubeta ^2) - 2 U \ualpha \uzeta, 
\\
Q_3(\ubmr) & = - 2U  \ualpha \ubeta + 2\umu  \pr \ualpha  + 2 b \uzeta\ubeta 
   + U \uTheta_0^{(2)}, 
\\
Q_4(\ubmr) & = -\pr  \uTheta_0^{(2)} - \ualpha \pr \ubeta  + \ubeta  \pr \ualpha - 2U \ualpha \umu,
\end{split}
\end{align}
the cubic terms by
\begin{align}\label{Nrcub}
\begin{split}
\bmC(\ubmr) & = (C_1,C_2,C_3,C_4)^T(\ubmr),
\\
C_1(\ubmr) & = - \frac12 (\ualpha^2+\ubeta ^2) \ualpha - (\uzeta^2 + \umu ^2) \ualpha - 2\ueta_0^{(2)} \pt \ubeta  - \uTheta_0^{(2)} \ubeta,
\\
C_2(\ubmr) & = - \uzeta (\ualpha^2 + \ubeta ^2), 
\\
C_3(\ubmr) & = - \frac12 (\ualpha^2+\ubeta ^2)\ubeta  - (\uzeta^2+\umu ^2) \ubeta + 2\ueta_0^{(2)} \pt \ualpha + \uTheta_0^{(2)} \ualpha +U\uTheta_0^{(3)}, 
\\
C_4(\ubmr) & = -\partial_r\uTheta_0^{(3)} - \umu  (\ualpha^2+\ubeta ^2)  ,
\end{split}
\end{align}
and the quartic and higher order $\bmN^{(4)}$ by the remaining contributions
\begin{align}\label{Nr4}
\begin{split}
&\bmN^{(4)}(\ubmr) = (N^{(4)}_1,N^{(4)}_2,N^{(4)}_3,N^{(4)}_4)^T(\ubmr),\\
&N_1^{(4)}(\ubmr)=U \ueta_0^2 + \ueta_0^2 \ualpha- 2(\ueta_0^{(3)}+\ueta_0^{(4)}) \pt \ubeta  - (\uTheta_0^{(3)}+\uTheta_0^{(4)}) \ubeta,\\
&N^{(4)}_2(\ubmr)  = 0,\\
&N^{(4)}_3(\ubmr)  = \ueta_0^2 \ubeta+2(\ueta_0^{(3)}+\ueta_0^{(4)}) \pt \ualpha + (\uTheta_0^{(3)}+\uTheta_0^{(4)}) \alpha + U\uTheta_0^{(4)},
\\
&N^{(4)}_4(\ubmr) = -\partial_r\uTheta_0^{(4)}.
\end{split}
\end{align}
Note that to be precise $\bmQ$, $\bmC$, and $\bmN^{(4)}$ are functions of $\bmr$ and their space and time derivatives.  We use $\bmQ_s$, respectively $\bmC_s$, to denote
the symmetrized quadratic, respectively, cubic form, associated to \eqref{Nrquad}, respectively \eqref{Nrcub}. That is, we denote, for $\ubmr,\ubmr',\ubmr'' \in\R^4$, 
\begin{align}\label{QCforms}
\bmQ_s(\ubmr,\ubmr') & = \frac{1}{2}[\bmQ(\ubmr+\ubmr') - \bmQ(\ubmr)-\bmQ(\ubmr')\big],
\\
\bmC_s(\ubmr,\ubmr',\ubmr'') & = \frac{1}{6}[\bmC(\ubmr+\ubmr'+\ubmr'') 
    - \bmC(\ubmr+\ubmr')-\bmC(\ubmr+\ubmr'')-\bmC(\ubmr'+\ubmr'') + \bmC(\ubmr) + \bmC(\ubmr') + \bmC(\ubmr'')\big].
\end{align}

The next lemma introduces a decomposition the profile $\bfrf$ that will be useful 
for our estimates later on. It identifies the leading order contribution, denoted $\bfrf^F$, 
that leads to the Fermi Golden Rule  damping mechanism. We will use the notation $\bmQ_s$ introduced above, as well as the decomposition of the nonlinearity $\bmN$
as $\bmN = \bmN_c + \bmN_d$ according to the definitions in \eqref{equ:N_cd}-\eqref{equ:decomposition_bmN_cont_disc}.

\begin{lemma}\label{lemnlsetup}
Let $\bfrf$ be the profile defined in \eqref{defprof}.
Then, we have the following decomposition into a `Fermi' and a `Regular' component:
\begin{align}\label{lemnlsetupf}
\bfrf & = \bfrf^F + \bfrf^R  
\end{align}
where
\setlength{\leftmargini}{1.5em}
\begin{itemize}
\smallskip
\item The Fermi component is defined as
\begin{align}\label{deffF}
\bfrf^F := \sum_{j,\ell \in\{1,2\}} \bfrf^F_{j\ell} \qquad  
  \bfrf^F_{j\ell}(t) & := Z_j(t)Z_\ell(t) \big(  \Lambda_{j\ell} - i \Gamma_{j\ell} \big),
\end{align}
with
\begin{equation}
\begin{aligned}\label{defLamGam}
  \Lambda_{j,\ell} & := e^{-it(\bfsM-2\lambda)}   \dfrac{1}{2\bfsM}  \mathrm{p.v.}\frac{1}{\bfsM-2\lambda}
  \bfP_c \bmQ_s(\bfY_j,\bfY_\ell),
\\
\Gamma_{j,\ell} & := \dfrac{\pi}{ 4\lambda } \delta(\bfsM-2\lambda) \bfP_c\bmQ_s(\bfY_j,\bfY_\ell).
\end{aligned}\end{equation}
\smallskip
\item The `Regular' component is given by
\begin{align}\label{deffR}
\bfrf^R & = \bfrf^{R}_c  + \bfrf^{R}_d
\end{align}
where
\begin{align}\label{deffRc}
\begin{split}
\bfrf^{R}_c(t) = \bfrf(0) + \int_0^t e^{-is\bfsM} 
  \dfrac{1}{2i \bfsM} \bfP_c \bmN_c(\bmu) \, \ud s, 
\end{split}
\end{align}
and
\begin{align}\label{deffRd}
\begin{split}
\bfrf^{R}_d& := \bfrf^{R,1}_d + \bfrf^{R,2}_d, 
\\
\bfrf^{R,1}_d(t) & :=  \int_0^t e^{-is\bfsM} 
  \dfrac{1}{2i \bfsM} \big[ \bfP_c \bmN_d - \sum_{j,\ell \in\{1,2\}} 
  e^{2is\lambda} \bfP_c\bmQ_s\big( Z_j(s)\bfY_j, Z_\ell(s)\bfY_\ell \big) \big] \, \ud s,
\\
\bfrf^{R,2}_d(t) & := \sum_{j,\ell \in\{1,2\}} \frac{1}{i(\bfsM-2\lambda + i0^+)} \Big[
 \dfrac{1}{2i\bfsM}  \mathbf{P}_c\bmQ_s\big(Z_j(0)\bmY_j,Z_\ell(0)\bmY_\ell\big)  
\\
& \phantom{:= \sum_{j,\ell \in\{1,2\}} \frac{1}{i(\bfsM-2\lambda + i0^+)} \Big[} +  \int_0^t e^{-is(\bfsM-2\lambda)}  
  \dfrac{1}{i\bfsM} \mathbf{P}_c\bmQ_s\big( \dot{Z}_j(s)\bmY_j,Z_\ell(s)\bmY_\ell\big)  \, ds  \Big].
\end{split}
\end{align}
\end{itemize}
\smallskip
In particular, if we define the `solutions' corresponding to the above profiles by
\begin{align}\label{uFR}
\bmu^\ast
  :=  e^{it\bfsM} \bfrf^\ast   +  e^{-it\bfsM} \overline{\bfrf^\ast}, 
    \qquad \bfrf^\ast \in \{\bfrf^F,\bfrf^R, \bfrf^R_c, \bfrf^{R}_d\}, 
\end{align}
we then have the decomposition of the continuous part of the solution
\begin{align}\label{uFRdec}
\bfP_c\bmr = \bmu = \bmu^F + \bmu^{R} =  \bmu^F + 
  \bmv + \bmu^{R}_d
\end{align}
where $\bmv$ is the good part of $\bmu$, first introduced in \eqref{equ:good_bad_decomposition},
that satisfies the equation 
\eqref{equ:evol_equ_bmv}.
\end{lemma}

\begin{proof}
Starting from the definition of $\bfrf$ and the equation \eqref{equ:system_evol_equations_u_z},
according to the definition \eqref{deffRc}-\eqref{deffRd}, we can write Duhamel's formula as
\begin{align}\label{lemnlsetuppr1}
\begin{split}
\bfrf(t) & = \bfrf(0) + \int_0^t e^{-is\bfsM} 
  \dfrac{1}{2i \bfsM} \bfP_c \bmN(\bmr) \, \ud s
  \\ & = \bfrf^{R}_c(t) + \bfrf^{R,1}_d (t) 
  + \sum_{j,\ell \in\{1,2\}}\int_0^t e^{-i(s\bfsM-2\lambda)} \dfrac{1}{2i \bfsM} \bfP_c\bfQ_s\big(Z_j(s)\bfY_j,Z_\ell(s)\bfY_\ell \big) \, \ud s.
\end{split}
\end{align}
Here we have used that $\bmN = \bmN_c + \bmN_d$ 
according to the definitions in \eqref{equ:N_cd}-\eqref{equ:decomposition_bmN_cont_disc}.
The statement of the lemma is then equivalent to
\begin{align}\label{lemnlsetuppr2}
\sum_{j,\ell \in\{1,2\}} \int_0^t e^{-is(\bfsM-2\lambda)} \dfrac{1}{2i \bfsM} \bfP_c\bmQ_s\big( Z_j(s)\bfY_j, Z_\ell(s)\bfY_\ell \big) \, \ud s
 = \bfrf^{R,2}_d(t) + \bfrf^F(t). 
\end{align}
To verify \eqref{lemnlsetuppr2} we regularize the exponential on the left-hand side
and integrate by parts:
\begin{align}\label{lemnlsetuppr3}
\begin{split}
& \int_0^t e^{-is(\bfsM-2\lambda + i0^+)} \dfrac{1}{2i \bfsM} 
  \bfP_c\bmQ_s\big(Z_j(s)\bfY_j, Z_\ell(s)\bfY_\ell \big) \, \ud s
 \\
 & = \frac{1}{-i(\bfsM-2\lambda + i0^+)} \Big[ e^{-it(\bfsM-2\lambda)}  \dfrac{1}{2i\bfsM}  
 \mathbf{P}_c\bmQ_s\big(Z_j(t)\bmY_j,Z_\ell(t)\bmY_\ell\big)   
  \\ & -  \dfrac{1}{2i\bfsM}  \mathbf{P}_c\bmQ_s\big(Z_j(0)\bmY_j,Z_\ell(0)\bmY_\ell\big)   
  - \int_0^t e^{-is(\bfsM-2\lambda)}  
  \dfrac{1}{i\bfsM} \mathbf{P}_c\bmQ_s\big( \dot{Z}_j(s)\bmY_j,Z_\ell(s)\bmY_\ell\big)  \, ds  \Big];
\end{split}
\end{align}
the choice of sign in $+i0^+$ is motivated by the fact that the operator $e^{it\bfsM}$
acting on the last term in \eqref{lemnlsetuppr3} decays in time. 
The desired conclusion follows from applying the formula 
\begin{align}\label{deltapvidentity}
\lim_{\eta \rightarrow 0^+} \frac{1}{x + 
  i\eta} = \mathrm{p.v.} \frac{1}{x} - 
  i \pi \delta_x
\end{align}
to the first term on the right-hand side of the \eqref{lemnlsetuppr3} to obtain
\begin{align}\label{lemnlsetuppr4}
\begin{split}
\frac{1}{-i(\bfsM-2\lambda + i0^+)} \Big[ e^{-it(\bfsM-2\lambda)}  \dfrac{1}{2i\bfsM}  
 \mathbf{P}_c\bmQ_s\big(Z_j(t)\bmY_j,Z_\ell(t)\bmY_\ell\big) \Big] = \bfrf^F_{j\ell}(t).
\end{split}
\end{align}

Finally, \eqref{uFRdec} follows from \eqref{defprof},
the definitions of the profiles \eqref{lemnlsetupf} and \eqref{deffR},
and the notation \eqref{uFR},
once we observe that $\bmv = \bmu^{R}_c$ in view of \eqref{deffRc}.
\end{proof}


\subsection{Fermi Golden Rule}\label{secIM1}
We now analyze the ODE for the profile of the internal mode.
We first write out the systems of ODE for the internal modes $Z_j$
in \S \ref{secODE}, and organize the various nonlinear terms.
Then, in \S \ref{secFGRweak} we prove that the ODE dynamics lead to damping of the 
amplitudes of the internal modes provided certain coefficients are proven to be non-zero,
and that we have suitable estimates on sub-leading order terms.
The verification of the necessary non-vanishing conditions \eqref{propODEdampFGR} is guaranteed 
by \cite[Proposition~4.3]{LPPSS1}. 
The estimates needed on all remainder terms are proven in Section \ref{secIM2}.

\subsubsection{System of ODEs for the internal mode}\label{secODE}
Recall the relations \eqref{defZj'} and \eqref{eq:systemZj0} for $Z_j$.
Using the splitting of $\bmN$ in \eqref{splitN} we can write
\begin{align}
\pt Z_j & =
 \dfrac{1}{2i\lambda} e^{-i\lambda t}
  \Big[ \big\langle \bmQ(\bmr-\bmu), \bfY_j \big\rangle + \big\langle \bmC(\bmr-\bmu), \bfY_j \big\rangle  \Big]\label{eq:Z12}
\\ 
  & + \dfrac{1}{2i\lambda} e^{-i\lambda t} \Big[ \big\langle \bmQ(\bmr) - \bmQ(\bmr-\bmu), \bfY_j \big\rangle 
    + \big\langle \bmC(\bmr) - \bmC(\bmr-\bmu) , \bfY_j \big\rangle +  \big\langle \bmN^{(4)}(\bmr), \bfY_j \big\rangle \Big].\nonumber
\end{align}
In particular, the first line of \eqref{eq:Z12} contains only nonlinear terms 
(quadratic and cubic) in $Z_1,Z_2$.
%
%
%
%
%
Next, we utilize the formulas and definitions from Lemma \ref{lemnlsetup},
plug them into the ODE \eqref{eq:Z12} and organize all the nonlinear terms
in the following lemma.

\begin{lemma}[The ODE]\label{lemODEsetup}
For $j=1,2$, 
\begin{align}\label{FGRODE1}
\begin{split}
\dot{Z}_j(t) & = \frac{1}{i\lambda} e^{-i\lambda t} \sum_{k=1,2}
  \Big\langle \bmQ_s\Big(e^{-i\lambda t}\overline{Z}_k(t)\bmY_k, \,  
   e^{it\bfsM} \bfrf^{\mathrm{F}}\Big), \, \bmY_j \Big\rangle + R_j(t)
\\
R_j(t) & := Q_{Z,j} + C_{Z,j} + R_{Z,j}^{\mathrm{disc}} + R_{Z,j}^\mathrm{field} 
  + \dfrac{1}{2i\lambda} e^{-i\lambda t} \big[\big\langle \bmC(\bmr) - \bmC(\bmr-\bmu) , \bmY_j \big\rangle + \big\langle \bmN^{(4)}(\bmr), \bfY_j \big\rangle\big]
\end{split}
\end{align}
where:

\setlength{\leftmargini}{1.5em}
\begin{itemize}

\smallskip
\item $Q_{Z,j}$ are quadratic terms in $z_j$,
\begin{align}\label{ODEQZ}
Q_{Z,j} := \dfrac{1}{2i \lambda} e^{-i\lambda t} \big\langle \bmQ_s \big(z_1(t)\bmY_1+z_2(t)\bmY_2\big),\, \bmY_j \big\rangle,
\end{align}
where we recall that $\bmQ_s$ is defined in \eqref{Nrquad}.

\smallskip
\item $C_{Z,j}$ are cubic terms in $z_j$:
\begin{align}\label{ODECZ}
C_{Z,j} :=  \dfrac{1}{2i \lambda}  e^{-i\lambda t} \big\langle \bmC_s \big(z_1(t)\bmY_1+z_2(t)\bmY_2\big),\, \bmY_j \big\rangle,
\end{align}
where we recall that $\bmC_s$ is defined in \eqref{Nrcub}. 

\smallskip
\item $R_{Z,j}^{\mathrm{disc}}$ are quadratic remainder terms containing the 
interactions between the discrete components and the continuous component
that are not the \emph{resonant} ones explicitly written out on the right-hand side of \eqref{FGRODE1}: 
\begin{align}\label{ODERdisc}
\begin{split}
R_{Z,j}^{\mathrm{disc}} & :=  \dfrac{1}{2i \lambda}  e^{-i\lambda t} 
  \sum_{\ell=1,2} 
  \Big[  2 \Big\langle \bmQ_s\Big(e^{-i\lambda t}\overline{Z}_\ell(t)\bmY_\ell,  \, 
   e^{-it\bfsM}\,\overline{\bfrf^{ \mathrm{F}}}\Big), \bmY_j \Big\rangle
  \\
  & \phantom{:=  \dfrac{1}{2i \lambda}  e^{-i\lambda t} 
  \sum_{\ell=1,2} 
  \Big[} + 2 \big\langle \bmQ_s \big(e^{i\lambda t}Z_\ell(t)\bmY_\ell, \bmu^F\big), \bmY_j \big\rangle \Big]
  +    \dfrac{1}{2i \lambda}   e^{-i\lambda t} \big\langle \bmQ_s\big(\bmu^F, \bmu^F\big), \bmY_j \big\rangle.
\end{split}
\end{align}

\smallskip
\item $R_{Z,j}^\mathrm{field}$ are quadratic terms containing all other interactions
that have not been accounted for so far:
\begin{align}\label{ODERfield}
\begin{split}
R_{Z,j}^\mathrm{field} & :=  \dfrac{1}{2i \lambda}  e^{-i\lambda t} \Big[ 
   2 \sum_{\ell=1,2}  { \big\langle \bmQ_s\big(z_\ell(t)\bmY_\ell, \bmu^{R}\big), \bmY_j \big\rangle }
   \\
   & \phantom{:=  \dfrac{1}{2i \lambda}  e^{-i\lambda t} \Big[} + 2\big\langle \bmQ_s\big(\bmu^{F}, \bmu^{R}\big), \bmY_j \big\rangle 
   + \big\langle \bmQ_s\big(\bmu^{R}, \bmu^{R}\big), \bmY_j \big\rangle \Big].
\end{split}
\end{align}

\smallskip
\item The last term in \eqref{FGRODE1} contains cubic and higher order remainders with at least one $\bmu$,
according to the definitions of $\bmC$ in \eqref{Nrcub} and of $C_{Z,j}$ in \eqref{ODECZ}.

\end{itemize}
\end{lemma}

\begin{proof}
The proof of the above lemma is a direct re-arrangement of all the terms from \eqref{eq:systemZj0}
according to the definitions and formulas in Lemma \ref{lemnlsetup}.
\end{proof}

\subsubsection{Dissipative structure and Fermi Golden Rule}\label{secFGRweak}
In this subsection we use the (approximate Hamiltonian) structure of the PDE, the block structure of the operator $\bfM$, and the structure of the quadratic nonlinearities in the PDE,
to identify the leading order damping dynamics in the ODE \eqref{FGRODE1}.
This dissipative structure will arise from the first terms on the right-hand side of \eqref{FGRODE1};
in particular, we will see that the ODE
\begin{align}
\begin{split}
\dot{Z}_j(t) & =
  \frac{1}{i\lambda}e^{-i\lambda t} \sum_{k=1,2}
  \Big\langle \bmQ_s \Big(e^{-i\lambda t}\overline{Z}_k(t)\bmY_k, \,  
  e^{it\bfsM} \bfrf^{\mathrm{F}}\Big), \, \bmY_j \Big\rangle
\end{split}
\end{align}
leads to time-decay of $Z_j$ provided a non-degeneracy condition (the Fermi Golden Rule) is satisfied. 

We begin with an abstract lemma. In the statement we use the notation introduced in \eqref{eq:frakzfirstdef1}. Later in applications $\zed$ in this lemma corresponds to $|Z_1|^2+|Z_2|^2$.

\begin{lem}\label{lemODEdamp2ndversion}
Let $0<\varepsilon\ll1$, and let $\zed :[0,T] \rightarrow [0,\infty)$ be a solution of
\begin{align}\label{lemODEdampeq}
    \dot\zed(t)
    =
    \mathbb{D}(t) + R(t),
    \qquad
    0\leq\zed(0)\leq\varepsilon^2,
\end{align} 
where the operator $\mathbb{D}$ satisfies,
\begin{align}\label{lemODE_damp_2ndv_dissipation}
    -\Gamma_0\zed(t)^2
    \leq
    \mathbb{D}(t)
    \leq
    -\Gamma_1\zed(t)^2,
    \qquad
    0\leq t < T,
\end{align}
for some $0<\Gamma_1\leq\Gamma_0$. For $\Gamma>0$, define
\begin{align}\label{defzG}
    \zed_{\Gamma}(t)
    :=
    \frac{\varepsilon^2}{1+\Gamma\varepsilon^2t}.
\end{align}
Assume that
\[
    R=R_a+R_b+R_h,
\]
and that, for every $0\leq t_1\leq t_2\leq T$,
\begin{equation}\label{lemODE_damp_2ndv_1}
\begin{aligned}
    \left|
    \int_{t_1}^{t_2}R_a(s)\,\ud s
    \right|
    &\lesssim
    \zed_{\Gamma_1}^{1+\delta_1}(t_1),
    \\
    \left|
    \int_{t_1}^{t_2}R_b(s)\,\ud s
    \right|
    &\lesssim
    \epsilon\zed_{\Gamma_1}(t_1),
\end{aligned}
\end{equation}
for some $\delta_1>0$, where $\epsilon>0$ is sufficiently small, and
\begin{align}\label{lemODE_damp_2ndv_2}
    R_h=\frac{\ud}{\ud t}h,
    \qquad
    |h(t)|\lesssim\zed_{\Gamma_1}^{3/2}(t).
\end{align}
Then, for $\varepsilon$ sufficiently small,
the following conclusions hold:
\begin{enumerate}[leftmargin=*]
\item[\rm(i)]
For $0\leq t\leq T$,
\begin{align}\label{lemODE_damp_2ndv_3}
    0
    \leq
    \zed(t)
    \leq
    \frac54\zed_{\Gamma_1}(t).
\end{align}

\item[\rm(ii)]
If, in addition, we assume that
\begin{align}\label{lemODE_damp_2ndv_initial_equality}
    \zed(0)=\varepsilon^2,
\end{align}
then
\begin{align}\label{lemODE_damp_2ndv_3_twosided}
    \frac34\zed_{\Gamma_0}(t)
    \leq
    \zed(t)
    \leq
    \frac54\zed_{\Gamma_1}(t),
    \qquad
    0\leq t\leq T.
\end{align}
\end{enumerate}
\end{lem}

\begin{proof}
    Let us define \[
    W(t) := \zed(t)+\int_t^T \Big( R_a(s)+R_b(s)\Big)\, \ud s - h(t),
    \]
    and observe that, under the given assumptions
    \begin{align}\label{eqw00}
    \big\vert W(t)-\zed(t)\big\vert \lesssim \zed_{\Gamma_1}^{1+\delta_1}(t)+\epsilon \zed_{\Gamma_1}(t)+\zed_{\Gamma_1}^{3/2}(t)\leq \eta \zed_{\Gamma_1}(t),
    \end{align}
    where \[
    \eta := C(\varepsilon^{2\delta_1}+\epsilon+\varepsilon).
    \]
    We choose $\varepsilon$ and $\epsilon$ small enough so that $\eta \ll1$. 

\medskip
{\it Proof of (i)}.
By definition of $W$ and $\zed_{\Gamma_1}$, we have that \begin{equation}\label{lemODE_damp_2ndv_4}\begin{aligned}
        \dot{W}(t) &= \mathbb{D}(t)
        \leq-\Gamma_1\zed(t)^2,
       \qquad  \dot{\zed}_{\Gamma_1}(t)& =-\Gamma_1\zed_{\Gamma_1}(t)^2.
    \end{aligned}\end{equation}
At the initial time we have that
\begin{align}\label{eqw01}
    W(0) \leq \zed(0) + \eta \zed_{\Gamma_1}(0) \leq (1+\eta)\zed_{\Gamma_1}(0) < \dfrac{6}{5}\zed_{\Gamma_1}(0).
    \end{align}
Then, suppose there is a first time $t_*\in(0,T]$ at which \[
    W(t_*)=\dfrac{6}{5}\zed_{\Gamma_1}(t_*).
    \]
    In particular, from \eqref{eqw00}
    it follows that 
    \begin{align}\label{lemODE_damp_2ndv_5}
    \zed(t_*) \geq W(t_*)-\eta \zed_{\Gamma_1}(t_*) = \Big( \dfrac{6}{5}-\eta\Big) \zed_{\Gamma_1}(t_*).
    \end{align}
    Therefore, from the equations \eqref{lemODE_damp_2ndv_4}, along with \eqref{lemODE_damp_2ndv_5}, we conclude that \begin{align*}
        \dfrac{\ud }{\ud t} \Big(W-\dfrac{6}{5}\zed_{\Gamma_1}\Big)(t_*) \leq -\Gamma_1 \zed(t_*)^2 + \dfrac{6}{5}\Gamma_1 \zed_{\Gamma_1}(t_*)^2 \leq -\Gamma_1 \bigg(\Big(\dfrac{6}{5}-\eta\Big)^2-\dfrac{6}{5}\bigg) \, \zed_{\Gamma_1}(t_*)^2 <0,
    \end{align*}
    having used that $\eta\ll1$. We conclude that, for all $t\in[0,T]$, 
    \[
    W(t)<\dfrac{6}{5}\zed_{\Gamma_1}(t),
    \]
    and hence, using again \eqref{eqw00},
    \eqref{lemODE_damp_2ndv_3} follows.

\medskip
{\it Proof of (ii)}.
Assume in addition that \eqref{lemODE_damp_2ndv_initial_equality} holds.
Then $\zed(0)
= \zed_{\Gamma_0}(0) = \zed_{\Gamma_1}(0) = \varepsilon^2.$
Consequently, \eqref{eqw00} gives
\[
    W(0)
    \geq
    \zed(0)-\eta\zed_{\Gamma_1}(0)
    > \frac45\zed_{\Gamma_0}(0).
\]
Suppose that there is a first time $t_*\in(0,T]$ at which
\[
    W(t_*)
    =
    \frac45\zed_{\Gamma_0}(t_*).
\]
At this time, \eqref{eqw00} 
implies that
\begin{align}\label{lemODE_damp_2ndv_5_lower}
\begin{aligned}
    \zed(t_*)
    &\leq
    W(t_*)+\eta\zed_{\Gamma_1}(t_*) \leq
    \left(
    \frac45+\frac{\Gamma_0}{\Gamma_1}\eta
    \right)
    \zed_{\Gamma_0}(t_*).
\end{aligned}
\end{align}
having used that $\zed_{\Gamma_1}(t) \cdot \Gamma_1 \leq 
\zed_{\Gamma_0}(t) \cdot \Gamma_0$. On the other hand, by \eqref{lemODE_damp_2ndv_dissipation},
\begin{align*}
    \dot W(t)
    &=
    \mathbb D(t)
    \geq
    -\Gamma_0\zed(t)^2,
    \qquad \dot{\zed}_{\Gamma_0}(t)
    = -\Gamma_0\zed_{\Gamma_0}(t)^2.
\end{align*}
It follows from \eqref{lemODE_damp_2ndv_5_lower} that
\begin{align*}
    \frac{\ud}{\ud t}
    \left(
    W-\frac45\zed_{\Gamma_0}
    \right)(t_*)
    &\geq
    -\Gamma_0\zed(t_*)^2
    +\frac45\Gamma_0\zed_{\Gamma_0}(t_*)^2 \geq
    \Gamma_0
    \left[
    \frac45
    -
    \left(
    \frac45+\frac{\Gamma_0}{\Gamma_1}\eta
    \right)^2
    \right]
    \zed_{\Gamma_0}(t_*)^2
    >0,
\end{align*}
having used that $\eta\ll1$.  This contradicts the fact that, at a first
crossing from above, the derivative of
$W-(4/5)\zed_{\Gamma_0}$ must be nonpositive and, therefore,
\[
    W(t)
    >
    \frac45\zed_{\Gamma_0}(t),
    \qquad
    0\leq t \leq T.
\]
Using again \eqref{eqw00}
we conclude that
\[
\begin{aligned}
    \zed(t)
    &\geq
    W(t)-\eta\zed_{\Gamma_1}(t) \geq
    \left(
    \frac45-\frac{\Gamma_0}{\Gamma_1}\eta
    \right)
    \zed_{\Gamma_0}(t)
    >
    \frac34\zed_{\Gamma_0}(t).
\end{aligned}
\]
Together with item ${\rm(i)}$, this proves
\eqref{lemODE_damp_2ndv_3_twosided} and completes the proof.
\end{proof}

The following proposition is the main statement about damping of the internal mode. We will use the notation
\begin{equation}
    \zed(t):=|Z_1(t)|^2+|Z_2(t)|^2.
\end{equation}
Recall that we are assuming $\zed(0) \leq \varepsilon^2$.

\begin{proposition}\label{propODEdamp} 
Consider the ODE system \eqref{FGRODE1} for $(Z_1,Z_2)$ for $t\in[0,T]$.
Assume that the following non-degeneracy conditions hold:
\begin{equation}\label{propODEdampFGR} 
\begin{aligned}
\big\langle \delta(\bfsM-2\lambda)\bfP_c\bmQ_s(\bfY_1,\bfY_1), \, \bmQ_s(\bmY_1, \bmY_1) \big\rangle \neq 0,
\\ 
\big\langle \delta(\bfsM-2\lambda)\bfP_c\bmQ_s(\bfY_2,\bfY_2), \, \bmQ_s(\bmY_2, \bmY_2) \big\rangle \neq 0,
\\
\big\langle \delta(\bfsM-2\lambda)\bfP_c\bmQ_s(\bfY_1,\bfY_2), \, \bmQ_s(\bmY_1, \bmY_2) \big\rangle \neq 0.
\end{aligned}
\end{equation}
Let $\Gamma_0$ and $\Gamma_1$ be as in \eqref{eq:Gamma01def} below. Suppose that
\begin{align}\label{propODEshort}
|\dot{Z}_1(t)|+|\dot{Z}_2(t)| \lesssim \varepsilon^2 , \qquad \mbox{for} \quad 0 \leqslant t \leqslant (C\varepsilon)^{-1}
\end{align}
for some absolute constant $C$ large enough,
and that the remainders can be written as
$$R_j = R_{j,a} + R_{j,b} + R_{j,h},$$
with 
\begin{align}\label{propODEdampR1a}
\begin{split}
\Big| \sum_{j=1,2} \int_{t_1}^{t_2} \Re \big( R_{j,a}(s)\overline{Z_j}(s) \big) \, \ud s \Big| 
  \lesssim \frz_{\Gamma_1}^{1+\delta_1}(t_1) , 
  \qquad (C\varepsilon)^{-1} \leqslant t_1\leqslant t_2 \leqslant T,
\end{split}
\end{align}
for some $0 < \delta_1 < 1/4$, and
\begin{align} \label{propODEdampR1b}
&  \bigg\vert \sum_{j=1,2} \int_{t_1}^{t_2} \Re \big( R_{j,b}(s)\overline{Z_j}(s) \big) \, \ud s \bigg\vert  
    \lesssim \epsilon \, \frz_{\Gamma_1}(t_1),  \qquad  (C\varepsilon)^{-1} \leqslant t_1\leqslant t_2 \leqslant T, 
\end{align}
for an absolute constant $\epsilon$ small enough,
and
\begin{align} \label{propODEdampR1h}
\sum_{j=1,2} \Re \big( R_{j,h}(t)\overline{Z_j}(t) \big) = \frac{d}{dt} h(t), \qquad
& |h(t)| \lesssim \frz_{\Gamma_1}^{3/2}(t),\qquad t\in[0,T],
\end{align}
where $\frz_{\Gamma_1}$ is defined as in \eqref{defzG}. Then, 

\begin{enumerate}[leftmargin=*]

\item[\rm(i)]
For $0\leq t\leq T$,
\begin{align}\label{propODEdampconc0}
|Z_1(t)|^2 + |Z_2(t)|^2 \leqslant \frac{5}{4}\frz_{\Gamma_1}(t).
\end{align}

\item[\rm(ii)]
If, in addition, we assume that
\begin{align}\label{lemODE_damp_2ndv_initial_equality2}
    \zed(0)=\varepsilon^2,
\end{align}
then
\begin{align}\label{propODEdampconc}
\frac{3}{4}\frz_{\Gamma_0}(t) \leqslant |Z_1(t)|^2 + |Z_2(t)|^2 \leqslant \frac{5}{4}\frz_{\Gamma_1}(t),\qquad t\in[0,T].
\end{align}

\end{enumerate}


\end{proposition}

\begin{rem}[About short times]\label{remshort}
Notice that if the constant $C$ in \eqref{propODEshort} is sufficiently large, then this smallness condition 
guarantees the desired conclusion \eqref{propODEdampconc}
up to the local time of existence of $O(\varepsilon^{-1})$.
Therefore, the conditions \eqref{propODEdampR1a}-\eqref{propODEdampR1h} on the remainders 
only need to be imposed (and hence later verified)
for times beyond $O(\varepsilon^{-1})$. 
\end{rem}

\medskip
For the proof of the above result we need some intermediary results.
First, we need the following lemma, whose proof exploits the 
(approximate) Hamiltonian structure of the PDE.

\begin{lem}\label{lem1}
Let $\bmQ_s$ be the symmetric quadratic form defined in \eqref{Nrquad}, $\bmY_1,\bmY_2$ be the internal modes, and denote
\begin{align}\label{lem1not}
\bmB^{jk} := \bmQ_s \big( \bmY_j, \bmY_k).
\end{align}
Then, for all $a,b,c,d\in\{1,2\}$ we have
\begin{align}\label{lem1conc}
\begin{split}
\Big\langle \bmQ_s\Big( \bmY_a,  \,  \delta(\bfsM-2\lambda) \bmB^{bc} \Big), \,  \bmY_d \Big\rangle  
 & = \big \langle \delta(\bfsM-2\lambda) \bmB^{bc},  \,  \bmQ_s\big(\bmY_a,\bmY_d) \big\rangle,
\\
\Big\langle \bmQ_s\Big( \bmY_a, \,  \frac{1}{2\bfsM}\mathrm{p.v.}\frac{1}{\bfsM-2\lambda}  \bmB^{bc} \Big),  \, \bmY_d \Big\rangle  
 & = \Big \langle \frac{1}{2\bfsM}\mathrm{p.v.}\frac{1}{\bfsM-2\lambda}  \bmB^{bc}, \,  \bmQ_s\big(\bmY_a,\bmY_d) \Big\rangle.
\end{split}
\end{align}
\end{lem}
\begin{proof}
Let us introduce the `dual quadratic form' of a quadratic symmetric form $Q_s: \R^4\times \R^4 \rightarrow \R^4$, 
denoted $Q_s^\ast$ (which is not necessarily symmetric), by the identity
\begin{align}\label{FGRprQ}
\big\langle Q_s(f,g), h \big\rangle = \big\langle f, Q_s^\ast(g,h) \big\rangle =
\big\langle g, Q_s^\ast(f,h)\big\rangle, 
\end{align}
where $f,g,h$ are Schwartz functions and $\langle\cdot,\cdot\rangle$ denotes the $L^2_{r\ud r}$ inner product.

We want to use the Hamiltonian structure to prove \eqref{lem1conc}.
To account for the quadratic terms that are not in Hamiltonian form on the right-hand side of \eqref{eq:lemHamQeq} we define 
\begin{align}
\mathcal{H}_G(\bmr) := \int \big(-U\beta + \partial^\ast\mu\big) \Theta_0 \, r \ud r,
 \qquad -\dfrac{\delta \mathcal{H}_G}{\delta \bmr}
 = \big( 0,0, U\Theta_0,- \partial_r \Theta_0\big)^{T},
\end{align}    
where we recall that $\Theta_0$ represent $\partial_t\eta_0$ viewed as a function of $\bmr$ and $\dot\bmr$. See \eqref{eq:Theta0def1}. Here we have used the fact that we are working in Stuart's gauge, $U\beta=\partial^*\mu$, 
so that in particular the contribution from the variational derivative of $\Theta_0$ to that of $\calH_G$ vanishes. 
From Lemma \ref{lemHamQ} and the definition of $\mathcal{H}_2$
in \eqref{eq:con_ham_2}, we then write
\begin{align}\label{FGRprH2}
\mathcal{H}_2 + \mathcal{H}_G 
    =\calH_2+\calH_{G,2}+\calH_{G,3}+\calH_{G,4}:= -\int \bmr \cdot \bmQ^{\mathcal{H}_2}_s(\bmr,\bmr) \, r dr+\calH_{G,3}+\calH_{G,4},
\end{align}
where
\begin{equation}\label{eq:calGH34def1}
    \calH_{G,3}:=\int \big(-U\beta + \partial^\ast\mu\big) \Theta_0^{(3)} \, r \ud r,\qquad \calH_{G,4}:=\int \big(-U\beta + \partial^\ast\mu\big) \Theta_0^{(4)} \, r \ud r,
\end{equation}
and $\bmQ_s^{\mathcal{H}_2}$ is the symmetrization of $\bmQ^{\calH_2}$ according to \eqref{QCforms}.
Note that the choice of $\bmQ^{\calH_2}$ is not unique, and one possible choice, the one we use, is 
\begin{align}
\bmQ^{\mathcal{H}_2}(\bmr) = \Big(-\frac12 U \alpha ^2 - \frac12 U \beta ^2 + \partial^\ast(\beta \mu )- \partial_r\beta \mu  
  + b \alpha  \zeta - U \zeta^2 - U \mu ^2, b \beta ^2, U\Theta_0^{(2)}, 
  -\partial_r\Theta_0^{(2)} \Big)^{T}.
\end{align}
By \eqref{FGRprH2}, \eqref{eq:Hamrrdot}, and \eqref{splitN}, we have that
\begin{equation}\label{FGRpr2}
\bmQ(\bmr)  = -\frac{\delta}{\delta \bmr}
    \big( \mathcal{H}_2(\bmr) + \mathcal{H}_{G,2}(\bmr) \big) 
 = 
  \bmQ_s^{\mathcal{H}_2}(\bmr,\bmr) + (\bmQ_s^{\mathcal{H}_2})^\ast(\bmr,\bmr) + 
  (\bmQ_s^{\mathcal{H}_2})^\ast(\bmr,\bmr).
\end{equation}
Then, with the notation \eqref{lem1not}, we calculate
\begin{align*}
& \big\langle \bmQ_s\big( \bmY_a, \delta(\bfsM-2\lambda) \bmB^{bc} \big), \bmY_d \big\rangle  
  \\ & \quad =  \big\langle \bmQ_s^{\mathcal{H}_2}(\bmY_a,\delta(\bfsM-2\lambda)\bmB^{bc}),  \bmY_d\big\rangle
  + \big\langle (\bmQ_s^{\mathcal{H}_2})^\ast(\bmY_a,\delta(\bfsM-2\lambda)\bmB^{bc}),  \bmY_d\big\rangle
\\ 
&  \quad  + \big\langle (\bmQ_s^{\mathcal{H}_2})^\ast(\delta(\bfsM-2\lambda)\bmB^{bc},\bmY_a),  \bmY_d\big\rangle
\\ & \quad  =  \big\langle \delta(\bfsM-2\lambda)\bmB^{bc},(\bmQ_s^{\mathcal{H}_2})^\ast (\bmY_a,\bmY_d)\big\rangle   
  +  \big\langle \delta(\bfsM-2\lambda)\bmB^{bc},\bmQ_s^{\mathcal{H}_2}(\bmY_a,\bmY_d)\big\rangle
\\ &  \quad  + \big\langle \delta(\bfsM-2\lambda)\bmB^{bc},(\bmQ_s^{\mathcal{H}_2})^\ast (\bmY_d,\bmY_a)\big\rangle
\\ & \quad  = \big\langle \delta(\bfsM-2\lambda)\bmB^{bc}, \bmQ_s(\bmY_a,\bmY_d)\big\rangle
\end{align*}
The same calculation applies when substituting the $\delta$ with the $\mathrm{p.v.}$. The proof is complete.
\end{proof}

We also need a simple lemma about orthogonality. In the statement we will use $\{e_1,\dots,e_4\}$ to denote the standard basis of $\bbR^4$.

\begin{lem}\label{lem3}
Define the orthogonal subspaces $V_1 = \mathrm{span}(e_1,e_2)$ and $V_2 = \mathrm{span}(e_3,e_4)$.
Let $\bmQ$ be the quadratic terms in the nonlinearity,
defined through \eqref{splitN}--\eqref{Nrquad}, and let $\bmQ_s$ denote the corresponding symmetric bilinear form defined by \eqref{QCforms}.
Then for all $\bmr_j \in V_j$, $j=1,2$, we have
\begin{align}\label{lem3conc1}
\bmQ_s\big(\bmr_1,\bmr_1), \,\, \bmQ_s\big(\bmr_2,\bmr_2) \in V_1,
  \qquad \bmQ_s\big(\bmr_1,\bmr_2) \in V_2.
\end{align}
Moreover
\begin{align}\label{lem3conc2}
\delta(\bfsM-2\lambda) \bmr_j, \ \frac{1}{2\bfsM}\mathrm{p.v.}\frac{1}{\bfsM-2\lambda} \bmr_j \, \in V_j, \qquad j=1,2.
\end{align}
In particular, using the same notation in \eqref{lem1not}, it follows that
\begin{align}\label{lem3conc3}
& \big\langle\delta(\bfsM-2\lambda) \bmB^{bc}, \, \bmQ_s\big(\bmY_a,\bmY_d) \big\rangle = 0
\\
& \mbox{unless 
$a=d$ and $b=c$, or $\{a,d\}  = \{b,c\} = \{1,2\}$.}
\end{align}
\end{lem}

\begin{proof}
To verify the first property in \eqref{lem3conc1} we inspect of the formulas \eqref{Nrquad}.
By plugging $\bma,\bmb \in V_1$ (that is, $a_3=a_4=b_3=b_4=0$) in \eqref{Nrquad}
one immediately verifies that $Q_3(\bma,\bmb) = U \Theta^{(2)}(\bma,\bmb)$
and $Q_4(\bma,\bmb) = -\pr  \Theta^{(2)}(\bma,\bmb)$.
Then, from the formula for $\Theta^{(2)}$ in \eqref{eq:Theta0jdef1} we see that
these last two expressions also vanish when  $\bma,\bmb \in V_1$.
Similarly, one sees that $Q_3(\bma,\bmb)=Q_4(\bma,\bmb) = 0$ for any $\bma, \bmb$ such that 
$a_1=a_2=b_1=b_2=0$.
A similar direct verification can be made for the last property in \eqref{lem3conc1}.
The property \eqref{lem3conc2} follows from the definition of the distorted Fourier transform 
and the block structure of $\bfM$, see \eqref{sysN1compact}.
The property \eqref{lem3conc3} is a direct consequence of the first two since $\bmY_j \in V_j$, $j=1,2$, 
see \eqref{eq:bmYsdef1}.
\end{proof}

As a  consequence of Lemma \ref{lem3}, we see that each of the terms below vanishes:
\begin{align*}
\big\langle \delta(\bfsM-2\lambda)\bfP_c\bmQ_s(\bfY_1,\bfY_1), \, \bmQ_s(\bmY_1, \bmY_2) \big\rangle = 0,
\\
\big\langle \delta(\bfsM-2\lambda)\bfP_c\bmQ_s(\bfY_2,\bfY_2), \, \bmQ_s(\bmY_1, \bmY_2) \big\rangle = 0,
\\
\big\langle \delta(\bfsM-2\lambda)\bfP_c\bmQ_s(\bfY_1,\bfY_2), \, \bmQ_s(\bmY_1, \bmY_1) \big\rangle = 0,
\\
\big\langle \delta(\bfsM-2\lambda)\bfP_c\bmQ_s(\bfY_1,\bfY_2), \, \bmQ_s(\bmY_2, \bmY_2) \big\rangle = 0,
\end{align*}
and similarly for the $\mathrm{p.v.}$ terms. We are now ready to prove Proposition \ref{propODEdamp}.

\begin{proof}[Proof of Proposition \ref{propODEdamp}]
From the definition $\bfrf^F$ in \eqref{deffF}-\eqref{defLamGam} we see that, for $j=1,2$,
\begin{align}\label{ODEdamppr1}
\begin{split}
\dot{Z}_j(t) & =
  \frac{1}{i\lambda}e^{-i\lambda t} \sum_{a=1,2}
  \Big\langle \bmQ_s\Big(e^{-i\lambda t}\bmY_a,  \, 
    e^{it\bfsM} \sum_{b,c \in\{1,2\}}  \Lambda_{bc}\Big), \, \bmY_j \Big\rangle \overline{Z_a}(t) Z_b(t)Z_c(t) 
  \\
  & + 
  \frac{1}{i\lambda} e^{-i\lambda t} \sum_{a=1,2}
  \Big\langle \bmQ_s\Big(e^{-i\lambda t} \bmY_a,  \, 
    -i  e^{it\bfsM} \sum_{b,c \in\{1,2\}}  \Gamma_{bc}\Big), \, \bmY_j \Big\rangle \overline{Z_a}(t) Z_b(t)Z_c(t)
  + R_j(t)
  \\
  & = \sum_{a,b,c=1,2} D^{a,b,c}_j \, \overline{Z_a}(t) Z_b(t)Z_c(t) 
  + i \sum_{a,b,c=1,2} H^{a,b,c}_j \, \overline{Z_a}(t) Z_b(t)Z_c(t) + R_j(t)
\end{split}
\end{align}
where
\begin{align}
D^{a,b,c}_j :=  -  \frac{\pi}{4\lambda^2}    
  \Big\langle \bmQ_s\Big(\bmY_a, \,  \delta(\bfsM-2\lambda)\bfP_c\bmQ_s(\bfY_b,\bfY_c)\Big), \, \bmY_j \Big\rangle 
\end{align}
and
\begin{align}
H^{a,b,c}_j := - \frac{1}{2\lambda}
  \Big\langle \bmQ_s\Big(\bmY_a, \,  \frac{1}{\bfsM}\mathrm{p.v.}\frac{1}{\bfsM-2\lambda} 
  \bfP_c \bmQ_s(\bfY_b,\bfY_c)\Big), \, \bmY_j \Big\rangle. 
\end{align}
Using Lemma \ref{lem1} we have
\begin{align}
D^{a,b,c}_j & =  -   \frac{\pi}{4\lambda^2}    \big\langle \delta(\bfsM-2\lambda)\bfP_c\bmQ_s(\bfY_b,\bfY_c), 
  \, \bmQ_s\big(\bmY_a, \bmY_j) \big\rangle, 
\\
H^{a,b,c}_j & = -   \frac{1}{2\lambda}   
  \Big\langle \frac{1}{\bfsM}\mathrm{p.v.}\frac{1}{\bfsM-2\lambda} 
  \bfP_c \bmQ_s(\bfY_b,\bfY_c), \, \bmQ_s\big(\bmY_a,\bmY_j) \Big\rangle. 
\end{align}
Then, using Lemma \ref{lem3} we see that the only non-vanishing terms are
\begin{align*}
&D_{1}^{1,1,1}=:D_1,\qquad D_{2}^{2,1,1}=D_{1}^{1,2,2}=:C_{12},\qquad D_2^{2,2,2}=:D_2,\\
&D^{2,1,2}_1 = D^{2,2,1}_1 = D^{1,1,2}_2 = D^{1,2,1}_2 =: D_{12}, \\
&H_{1}^{1,1,1}=:H_1,\qquad H_{2}^{2,1,1}=H_{1}^{1,2,2}=:C_{12}',\qquad H_2^{2,2,2}=:H_2,\\
&H^{2,1,2}_1 = H^{2,2,1}_1 = H^{1,1,2}_2 = H^{1,2,1}_2 =: H_{12},
\end{align*}
where
\begin{align}\label{FGRODE14}
\begin{split}
D_1& := -   \frac{\pi}{4\lambda^2}    \big\langle \delta(\bfsM-2\lambda)\bfP_c\bmQ_s(\bfY_1,\bfY_1),\bmQ_s(\bmY_1, \bmY_1) 
  \big\rangle,
\\
C_{12} & := -  \frac{\pi}{4\lambda^2}   \big\langle \delta(\bfsM-2\lambda)\bfP_c\bmQ_s(\bfY_1,\bfY_1),\bmQ_s(\bmY_2, \bmY_2) 
  \big\rangle ,
\\
D_2 & := - \frac{\pi}{4\lambda^2}   \big\langle \delta(\bfsM-2\lambda)\bfP_c\bmQ_s(\bfY_2,\bfY_2),\bmQ_s(\bmY_2, \bmY_2) 
  \big\rangle,
\\
D_{12}& := - \frac{\pi}{4\lambda^2}  \big\langle \delta(\bfsM-2\lambda)\bfP_c\bmQ_s(\bfY_1,\bfY_2),\bmQ_s(\bmY_2, \bmY_1) 
  \big\rangle,
\end{split}
\end{align}
and
\begin{align}\label{FGRODE14'}
\begin{split}
H_1 & := -   \frac{1}{2\lambda}   \big\langle\frac{1}{\bfsM}\mathrm{p.v.}\frac{1}{\bfsM-2\lambda} 
  \bfP_c\bmQ_s(\bfY_1,\bfY_1),\bmQ_s(\bmY_1, \bmY_1)
  \big\rangle,
\\
C_{12}' & := - \frac{1}{2\lambda}    \big\langle \frac{1}{\bfsM}\mathrm{p.v.}\frac{1}{\bfsM-2\lambda} 
  \bfP_c\bmQ_s(\bfY_1,\bfY_1),\bmQ_s(\bmY_2, \bmY_2) 
  \big\rangle,
\\
H_2 & := -  \frac{1}{2\lambda}   \big\langle \frac{1}{\bfsM}\mathrm{p.v.}\frac{1}{\bfsM-2\lambda} 
  \bfP_c\bmQ_s(\bfY_2,\bfY_2),\bmQ_s(\bmY_2, \bmY_2)
  \big\rangle,
\\
H_{12} & := - \frac{1}{2\lambda}   \big\langle \frac{1}{\bfsM}\mathrm{p.v.}\frac{1}{\bfsM-2\lambda} 
  \bfP_c\bmQ_s(\bfY_1,\bfY_2),\bmQ_s(\bmY_2, \bmY_1) 
  \big\rangle.
\end{split}
\end{align}
Note that all of the coefficients above are real valued and that 
\begin{align}
D_1,D_2,D_{12}\leqslant 0, \qquad |C_{12}| \leqslant \sqrt{D_1D_2}. 
\end{align}
With the above notation, the ODEs read
\begin{equation}
\label{FGRODE13}
\left\{
\begin{array}{ll}
\dot{Z}_1 & = (D_1+iH_1) |Z_1|^2 Z_1 +  2(D_{12}+iH_{12}) |Z_2|^2 Z_1 + (C_{12}+iC_{12}') \overline{Z_1} (Z_2)^2 + R_1,
\\
\\
\dot{Z}_2 & = (D_2+iH_2) |Z_2|^2 Z_2 +  2(D_{12}+iH_{12}) |Z_1|^2 Z_2 + (C_{12}+iC_{12}') \overline{Z_2} (Z_1)^2 + R_2.
\end{array}
\right.
\end{equation}
With $\frz(t) = |Z_1|^2 + |Z_2|^2$, we then calculate
\begin{align}\label{prODEdamp5}
\begin{split}
\frac{d}{dt} \frz(t) & = 2\Re \big( \dot{Z}_1 \overline{Z_1} + \dot{Z}_2 \overline{Z_2} \big)
 \\
 & = 2D_1 |Z_1|^4 + 4D_{12}|Z_2|^2 |Z_1|^2 + 2\Re\big[(C_{12}+iC_{12}') \overline{Z_1}^2 (Z_2)^2 \big] + 
 2\Re( R_1 \overline{Z_1})
 \\ 
 & \quad + 2D_2 |Z_2|^4 + 4D_{12}|Z_1|^2 |Z_2|^2 + 2\Re\big[(C_{12}+iC_{12}') \overline{Z_2}^2 (Z_1)^2 \big] + 
 2\Re( R_2 \overline{Z_2})
 \\
 & = -2 \big(\sqrt{|D_1|} |Z_1|^2 - \sqrt{|D_2|}|Z_2|^2 \big)^2 + 8D_{12}|Z_1|^2 |Z_2|^2
  \\
  & \quad + 4 \big[  C_{12} \Re(\overline{Z_1}^2 (Z_2)^2) - \sqrt{D_1D_2}|Z_1|^2 |Z_2|^2 \big]
  + 2\Re( R_1 \overline{Z_1}) + 2\Re( R_2 \overline{Z_2}).
\end{split}
\end{align}
Next, we use the fact that $|C_{12}| \leqslant \sqrt{D_1D_2}$ to get
\begin{align}\label{prODEdamp6}
\begin{split}
\frac{d}{dt} \frz(t) \leqslant  - 2 
  \Big[ \big(\sqrt{|D_1|} |Z_1|^2 - \sqrt{|D_2|}|Z_2|^2 \big)^2 + 4|D_{12}||Z_1|^2 |Z_2|^2 \Big]
  \\
  + 2\Re( R_1 \overline{Z_1}) + 2\Re( R_2 \overline{Z_2}).
\end{split}
\end{align}
Also, it is not hard to see that 
\begin{align}\label{prODEdamp6.5}
\begin{split}
\frac{d}{dt} \frz(t) &\geqslant  - 2 
  \Big[ \big(\sqrt{|D_1|} |Z_1|^2 - \sqrt{|D_2|}|Z_2|^2 \big)^2 + 4|D_{12}||Z_1|^2 |Z_2|^2
  + 4 \sqrt{D_1D_2}|Z_1|^2 |Z_2|^2 \Big]
  \\
  &\quad + 2\Re( R_1 \overline{Z_1}) + 2\Re( R_2 \overline{Z_2}),
\end{split}
\end{align}
and
\begin{align}\label{prODEdamp10}
\begin{split}
- \Gamma_0 \frz^2(t)
  \leqslant \frac{d}{dt} \frz(t) - 2\Re( R_1 \overline{Z_1}) - 2\Re( R_2 \overline{Z_2}) \leqslant - \Gamma_1 \frz^2(t),
\end{split}
\end{align}
where 
\begin{equation} \label{eq:Gamma01def}
    \begin{aligned} 
    \Gamma_0 &:= 4 \max \bigl\{ |D_1|, |D_2|, |D_{12}| \bigr\} \geq 0, \\
    \Gamma_1 &:= 2 \min \bigl\{ |D_{12}|,\sqrt{D_1D_2} \bigr\} \min \big\{ \sqrt{D_1/D_2}, \sqrt{D_2/D_1} \bigr\} \geq 0.
    \end{aligned}
\end{equation}
To prove the last (right-most) inequality in \eqref{prODEdamp10}, we have used that, for $a,b,c\neq 0$, $xy\geq0$,
\begin{align*}
(ax - by)^2 + 4c xy \geq (ax+by)^2 \qquad \text{if} \quad c>ab, 
\\
(ax - by)^2 + 4c xy \geq c \big(\sqrt{\tfrac{a}{b}}x + \sqrt{\tfrac{b}{a}}y\big)^2  \qquad \text{if}  \quad c<ab.
\end{align*}
The non-degeneracy assumptions \eqref{propODEdampFGR} guarantee that $\Gamma_0 > \Gamma_1 > 0$ and give the desired damping.

Indeed, we can apply Lemma \ref{lemODEdamp2ndversion} with the equation \eqref{prODEdamp5} and the inequality \eqref{prODEdamp10}
playing the role of \eqref{lemODEdampeq} and
\eqref{lemODE_damp_2ndv_dissipation} respectively,
and the assumptions on the remainder terms \eqref{propODEdampR1a}, \eqref{propODEdampR1b}, and \eqref{propODEdampR1h}
matching the required assumptions from 
\eqref{lemODE_damp_2ndv_1} and \eqref{lemODE_damp_2ndv_2}.
From the conclusions of Lemma \ref{lemODEdamp2ndversion} (see also Remark~\ref{remshort}) we obtain the conclusions of Proposition \ref{propODEdamp} as desired.
\end{proof}


\bigskip
\section{Internal Mode II: Estimates}\label{secIM2}

In this section we prove estimates for all the terms in the ODE \eqref{FGRODE1}.
In particular, we prove the estimates needed for Proposition \ref{propODEdamp}.
Here is our main result:

\begin{proposition}[Bounds for remainders]\label{propODERbounds}
With the same notation of Lemma \ref{lemODEsetup} and Proposition \ref{propODEdamp},
assume that the bootstrap assumptions \eqref{equ:bootstrap_assumption_z} and \eqref{equ:bootstrap_assumption_weighted_f}--\eqref{equ:bootstrap_assumption_HN_ILED_f} are satisfied. Then
\begin{equation}\label{eq:dotZshorttimebootstrap1}
    |\dot{Z}_1(t)|+|\dot{Z}_2(t)| \lesssim \varepsilon^2 , \qquad \mbox{for} \quad 0 \leqslant t \leqslant (C\varepsilon)^{-1}.
\end{equation}
Moreover, for $(C\varepsilon)^{-1}\leq t,t_1,t_2\leq T$ we have:

\setlength{\leftmargini}{2em}

\smallskip
\begin{itemize}
 \item The following estimates 
for the remainder terms involving quadratic terms in $Z$:
\begin{align}\label{IM2est1}
\big| Q_{Z,j}(t) \big| & \lesssim \frz_{\Gamma_1}(t),
\qquad
\sum_{j=1,2}\Re \big( Q_{Z,j}(t) \,\overline{Z_j}(t) \big) 
=: \frac{d}{dt}h_2(t), \, \quad |h_2(t)| \lesssim \frz_{\Gamma_1}^{3/2}(t).
\end{align}
\smallskip
\item The following estimates for the remainders involving cubic terms in $Z$:
\begin{align}\label{IM2est2}
\begin{split}
& \big| C_{Z,j}(t) \big| \lesssim \frz_{\Gamma_1}^{3/2}(t),
\\ 
& \sum_{j=1,2}\Re \big( C_{Z,j}(t) \,\overline{Z_j}(t) \big) 
=: \frac{d}{dt}h_3(t) + O\big(\frz_{\Gamma_1}^{2+}(t)\big), 
    \, \quad |h_3(t)| \lesssim \frz_{\Gamma_1}^{2}(t).
\end{split}
\end{align}

\item The following estimates for the (non-resonant) remainder terms involving $Z_j$ and the Fermi component $\bmu^F$:
\begin{align}
\label{IM2estint3}
\Big| \sum_{j=1,2}\int_{t_1}^{t_2} \Re \big( R_{Z,j}^{\mathrm{disc}}(s) \,\overline{Z_j}(s) \big) \, \ud s \Big| 
  & \lesssim \frz_{\Gamma_1}^{1+3\delta}(t_1).
\end{align}

\item The following estimates for the (quadratic) remainders involving $\bmu^R$ where $c_0$ is a small absolute constant:
\begin{align}
\label{IM2estint4}
\Big| \sum_{j=1,2}\int_{t_1}^{t_2} \Re \big( R_{Z,j}^{\mathrm{field}}(s) \,\overline{Z_j}(s) \big) \, \ud s \Big| 
  & \leqslant c_0 \, \frz_{\Gamma_1}(t_1).
\end{align}

\smallskip
\item The following estimates for cubic and higher order terms:
\begin{subequations}\label{IM2estimatesC}
\begin{align}\label{IM2estC}
\Big| \big\langle \bmC(\bmr(t)) - \bmC(\bmr(t)-\bmu(t)) , \bmY_j \big\rangle \Big| 
  \lesssim \frz_{\Gamma_1}^{5/3}(t), 
\end{align}
and
\begin{align}\label{IM2est4}
\Big| \big\langle \bmN^{(4)}(\bmr), \bfY_j \big\rangle 
\Big| \lesssim \frz_{\Gamma_1}^{2}(t).
\end{align}
\end{subequations}

\end{itemize}

\end{proposition}

As a corollary of the above proposition, we obtain that the assumptions \eqref{propODEdampR1a}-\eqref{propODEdampR1h} 
hold for all the remainder terms
on the right-hand side of \eqref{FGRODE1}.

\smallskip
\begin{proof}[Proof of Proposition \ref{propODERbounds}]
Estimate \eqref{eq:dotZshorttimebootstrap1} follows simply from the evolution equation \eqref{eq:systemZj0} for $Z_j$ and the bootstrap assumptions \eqref{equ:bootstrap_assumption_z} and \eqref{equ:bootstrap_assumption_HN_f}. We now turn to the other estimates in the statement of the proposition.
Recall that all the terms that need to be estimated are defined in Lemma \ref{lemODEsetup} 
where the nonlinearities $\bmQ$ and $\bmC$ are given in \eqref{Nrquad} and \eqref{Nrcub}, with \eqref{QCforms},
and with $\bfrf^F$, $\bmu^F$ and $\bmu^R$ defined in Lemma \ref{lemnlsetup}.
with the same notation from Lemma \ref{lemnlsetup}.

Recall, see \eqref{equ:bootstrap_assumption_Linfty_f}-\eqref{equ:bootstrap_assumption_HN_f}, that we have the a priori assumptions
\begin{align}
 \bigl\| e^{it\jD} \bmf_\ast(t)\bigr\|_{L^\infty_x(\bbR^2)} &\leq 2C_0 \left\{ \begin{aligned}
                            &\varepsilon^{1-10\delta} t^{-\frac12-4\delta},  &&1 \leq t \leq \varepsilon^{-2}, \\
                            &t^{-1+\delta},  &&\varepsilon^{-2} \leq t \leq T,
                         \end{aligned} \right. \\ 
    \bigl\| \bmf_\ast(t)\bigr\|_{H^N_x(\bbR^2)} &\leq 2C_0 \varepsilon, \quad 0 \leq t \leq T, 
    \end{align}
where
\begin{equation} \label{equ:definition_bmf_flat_profile_recalled_in_ODE_part}
    \begin{aligned}
        \bmf_\ast(t) &:= e^{-it\jD} (2i\jD)^{-1} (\pt+i\jD) \bmu_\ast(t), \qquad \ast \in \{\real,\imag\},
    \end{aligned}
\end{equation}
and
\begin{equation}
    \begin{aligned}
        \bmu(t) = e^{-i\theta} 
        \big( \bmu_\real(t) + i \bmu_\imag(t) \big), 
        \quad \bmu_\real(t) := \Re \bigl( e^{i\theta} \bmu(t) \bigr), 
        \quad \bmu_\imag(t) := \Im \bigl( e^{i\theta} \bmu(t) \bigr).
    \end{aligned}
\end{equation}
From these we can infer
\begin{equation}\label{ODERestu} 
\sum_{*\in\{\real,\imag\}}\sum_{\vert\alpha\vert\leq 5} 
\Big( {\| \nabla_x^\alpha \bmu_*(t) \|}_{L^\infty(\R^2)} 
+ 
 {\| \nabla_x^\alpha \partial_t \bmu_*(t) \|}_{L^\infty(\R^2)} \Big) 
 \lesssim  
\left\{ \begin{aligned}
                            &\varepsilon^{1-20\delta} t^{-\frac12-8\delta},  &&1 \leq t \leq \varepsilon^{-2}, \\
                            &t^{-1+2\delta},  &&\varepsilon^{-2} \leq t \leq T,
                         \end{aligned} \right. 
\end{equation}
using the $H^N$ norm to handle the derivatives by Gagliardo-Nirenberg-Sobolev interpolation.

\medskip
{\bf Proof of \eqref{IM2est1}.}
The first estimate \eqref{IM2est1} is immediate since $Q_{Z,j}$
is quadratic in $Z$. To see the other property in \eqref{IM2est1} 
we start from the definition \eqref{ODEQZ}. Using that $\bmQ$ is real-valued 
and 
\begin{align}\label{dtzj}
\begin{split}
& \qquad \partial_t z_j = -4\lambda^2\Re\Big( 
    \frac{e^{i\lambda t}}{2i\lambda} Z_j\Big).
\end{split}
\end{align}
we write
\begin{align}
\sum_{j=1,2}
\Re \big( Q_{Z,j}(t) \,\overline{Z_j}(t) \big)  
& =  \sum_{j=1,2} \Re \Big[ \dfrac{1}{2i \lambda}  
e^{-i\lambda t}\overline{Z_j}(t)\big\langle \bmQ\big(z_1(t)\bmY_1+z_2(t)\bmY_2\big), 
   \bmY_j \big\rangle \Big]
    \\
    & = \frac{1}{4\lambda^2}\sum_{j=1,2} \big\langle \bmQ \big(z_1(t)\bmY_1+z_2(t)\bmY_2\big), 
    \partial_t z_j \bmY_j \big\rangle
    \\
    \label{IM2est1pr1}
    & = \frac{1}{4\lambda^2} \big\langle \bmQ \big(z_1(t)\bmY_1+z_2(t)\bmY_2\big), 
    \partial_t (z_1 \bmY_1 + z_2 \bmY_2) \big\rangle .
\end{align}
Recall from \eqref{FGRpr2} that
\begin{align}\label{IM2est1pr2}
\bmQ(\bmr) = - \dfrac{\delta }{\delta \bmr}(\calH_2+\calH_{G,2}).
\end{align}
Then
\begin{equation}\label{IM2est1pr3}
    \big\langle \dfrac{\delta (\mathcal{H}_2+\calH_{G,2})}{\delta \bmr}\big(z_1(t)\bmY_1+z_2(t)\bmY_2\big),
    \partial_t (z_1 \bmY_1 + z_2 \bmY_2) \big\rangle 
    = 
    \frac{\ud}{\ud t}(\mathcal{H}_2+\calH_{G,2})\big(z_1 \bmY_1 + z_2 \bmY_2\big).
\end{equation}
Since $\mathcal{H}_2$ and $\calH_{G,2}$ are cubic polynomials in their arguments
(and spatial derivatives) the property \eqref{IM2est1} is immediately verified.
%
%
%
%

\medskip
{\bf Proof of \eqref{IM2est2}.}
From the definition in \eqref{ODECZ} we have $C_{Z,j} = O(|Z|^3) = O(\frz^{3/2})$.
We then use again the (approximate) Hamiltonian structure to obtain a cancellation as above.
Using that $\bmC$ is real valued and \eqref{dtzj}, we write
\begin{align}
\sum_{j=1,2} \Re \big( C_{Z,j}(t) \,\overline{Z_j}(t) \big)
& = \sum_{j=1,2}\Re \Big[   \dfrac{1}{2i \lambda}  
e^{-i\lambda t}\overline{Z_j}(t)\big\langle \bmC \big(z_1(t)\bmY_1+z_2(t)\bmY_2\big), 
   \bmY_j \big\rangle \Big] 
    \\
    & = \frac{1}{4\lambda^2}\sum_{j=1,2}\big\langle \bmC \big(z_1(t)\bmY_1+z_2(t)\bmY_2\big), 
    \partial_t z_j \bmY_j \big\rangle
    \\
      \label{IM2est2pr1}
    & = \frac{1}{4\lambda^2} \big\langle \bmC \big(z_1(t)\bmY_1+z_2(t)\bmY_2\big), 
    \partial_t (z_1\bmY_1 + z_2\bmY_2) \big\rangle 
\end{align}
From Lemma \ref{lemHamQ} we have 
\begin{align}
\bmC(\bmr) &= - \dfrac{\delta \mathcal{H}_3}{\delta \bmr}
 + \big( -\Theta_{0}^{(2)}\beta-\eta_0^{(2)}\partial_t\beta,0, \Theta_0^{(2)}\alpha+\eta_0^{(2)}\partial_t \alpha,0\big)^T + \big(-\eta_0^{(2)}\partial_t\beta,0, \eta_0^{(2)}\partial_t \alpha,0\big)^T\nonumber\\
 &\quad+\big( 0,0, U\Theta_0^{(3)},-\partial_r\Theta_{0}^{(3)}\big)^T.\label{IM2est2pr2}
\end{align}
Here to compute $\frac{\delta\calH_\calC}{\delta r}$ in the notation convention from Lemma \ref{lemHamQ} we have used that 
\begin{align*}
   \dfrac{ \delta \mathcal{H}_{\mathcal{C}} }{ \delta \alpha }   =  \eta_0 \partial_t\beta - \eta_0^2(U+\alpha),  \qquad\dfrac{\delta \mathcal{H}_{\mathcal{C}}}{\delta \zeta} =0,\qquad
     \dfrac{\delta \mathcal{H}_{\mathcal{C}}}{\delta \beta} = - \eta_0\partial_t\alpha -  \eta_0^2\beta , \qquad  \dfrac{\delta \mathcal{H}_{\mathcal{C}}}{\delta \mu} = 0.
\end{align*}
We can then handle the contribution from $\mathcal{H}_3$ in \eqref{IM2est2pr1}
exactly as in \eqref{IM2est1pr3}, and obtain a more than sufficient bound
since $\mathcal{H}_3$ is a quartic polynomial in its argument.
To deal with the  second term on the right-hand side of \eqref{IM2est2pr2}
we recall from  \eqref{eq:Theta0jdef1} that
\begin{align}\label{dteta0quad}
\Theta_0^{(2)}(\alpha,\zeta,\beta,\mu) & = (-\Delta+U^2)^{-1}
     \big(\beta L_1\alpha - 2U' \beta \zeta - \alpha L_1 \beta + 2U' \alpha \mu\big).
\end{align}
Recalling that $\bmY_1=(Y_1,Y_2,0,0)^t$ and $\bmY_2=(0,0,Y_1,Y_2)^t$,
we can check that 
\begin{align}\label{dteta0quad=0}
& \Theta_0^{(2)}(z_1(t)\bmY_1+z_2(t)\bmY_2) = \Theta_0^{(2)}(z_1Y_1,z_1Y_2,z_2Y_1,z_2Y_2) 
\\
& = (-\Delta+U^2)^{-1}
     \big(z_2Y_1 L_1(z_1Y_1) - 2U' z_2Y_1 z_1Y_2 - z_1Y_1 L_1(z_2Y_1) + 2U' z_1Y_1z_2Y_2\big) = 0.
\end{align}
Then, we use this to write the contribution of the second term on the right-hand side of \eqref{IM2est2pr2} 
to \eqref{IM2est2pr1} as
\begin{align}
& 
\int  \Big[ -\big(\Theta_0^{(2)}\beta+\eta_0^{(2)}\partial_t\beta\big)\Big|_{\bmr=z_1(t)\bmY_1+z_2(t)\bmY_2} \partial_t z_1 Y_1
    +  \big(\Theta_0^{(2)}\alpha+\eta_0^{(2)}\partial_t\alpha\big)\Big|_{\bmr=z_1(t)\bmY_1+z_2(t)\bmY_2} \partial_t z_2Y_1 \Big] 
    \, r\ud r
\\
\label{IM2est2pr5}
& = \int Y_1 \eta_0 \Big[ -\partial_t(z_2Y_1) (\partial_t z_1)
    +  (\partial_t z_1 Y_1) \partial_t z_2 \Big] \, r\ud r
\end{align}
which vanishes. 
Similarly, regarding the third term on the right-hand side of \eqref{IM2est2pr2}, it is enough to notice that \begin{align*}
\Big\langle (-\eta_0^{(2)}\partial_t\beta,0,\eta_0^{(2)}\partial_t\alpha,0)\Big|_{\bmr=z_1(t)\bmY_1+z_2(t)\bmY_2}, \,  \partial_t (z_1\bmY_1 + z_2\bmY_2) \Big\rangle =0.
\end{align*} 
For the last term on the right-hand side of \eqref{IM2est2pr2} the argument is similar to \eqref{IM2est1pr3} with the difference that now the variational derivative with respect to $\dot\bmr$ also enters. Concretely, recalling \eqref{FGRprH2}, the contribution of this term to \eqref{IM2est2pr1} can be written up to a multiplicative constant as (by an abuse of notation we have evaluated $\calH_{G,3}$ at $z_1\bmY_1+z_2\bmY_2$, while to be precise it should be evaluated at $(z_1\bmY_1+z_2\bmY_2,\dot z_1\bmY_1+\dot z_2\bmY_2)$ because $\calH_{G,3}$ is a function of $\bmr$ and $\dot\bmr$)
\begin{align}
    &\big\langle\frac{\delta }{\delta\bmr}\calH_{G,3}\big(z_1(t)\bmY_1+z_2(t)\bmY_2\big),
    \partial_t (z_1 \bmY_1 + z_2 \bmY_2) \big \rangle\\
    &=\frac{\ud}{\ud t}\calH_{G,3}\big(z_1 \bmY_1 + z_2 \bmY_2\big) -\big\langle \dfrac{\delta }{\delta \dot\bmr}\calH_{G,3}\big(z_1(t)\bmY_1+z_2(t)\bmY_2\big),
    \partial_t (\dot z_1 \bmY_1 + \dot z_2 \bmY_2) \big\rangle.
\end{align}
The first term on the right is already of the right form. For the second term, we recall from \eqref{eq:calGH34def1} that in view of our gauge condition, $\frac{\delta}{\delta\dot\bmr}\calH_{G,3}(\bmr,\dot\bmr)=0$. It follows that $\dfrac{\delta }{\delta \dot\bmr}\calH_{G,3}\big(z_1(t)\bmY_1+z_2(t)\bmY_2\big)$ is a cubic expression in $(\bmr,\dot\bmr)$ containing at least one factor of $\bmu$ or $\dot\bmu$. As such, its contribution to \eqref{IM2est1pr3} is of the form $O(\frz_{\Gamma_1}^{2+})$ which is acceptable in \eqref{IM2est2}. See the proof of \eqref{IM2estimatesC} below for more details in a similar contribution.

\medskip
{\bf Proof of \eqref{IM2estint3}.}
From the definition of $R_{Z,j}^{\mathrm{disc}}$ in \eqref{ODERdisc} we can write (we omit the unimportant dependence on $j=1,2$)
\begin{align}\label{IM2est31}
\begin{split}
R_{Z,j}^{\mathrm{disc}} & = R_{1}^{\mathrm{disc}} + R_{2}^{\mathrm{disc}} + R_{3}^{\mathrm{disc}},
\\
R_{1}^{\mathrm{disc}} & :=  \dfrac{1}{i\lambda}  e^{-i\lambda t} \sum_{\ell=1,2} 
   \big\langle \bmQ_s\big(e^{-i\lambda t}\overline{Z}_\ell(t)\bmY_\ell, 
   e^{-it\bfsM}\overline{\bfrf^{\mathrm{F}}}), \bmY_j \big\rangle,
\\
R_{2}^{\mathrm{disc}} & :=   \dfrac{1}{i\lambda}    e^{-i\lambda t} \sum_{\ell=1,2}    
\big\langle \bmQ_s \big(e^{i\lambda t}Z_\ell(t)\bmY_\ell, \bfu^F\big), \bmY_j \big\rangle ,
\\
R_{3}^{\mathrm{disc}} & := \dfrac{1}{2i\lambda}  e^{-i\lambda t} \big\langle \bmQ_s\big(\bfu^F, \bfu^F\big), \bmY_j \big\rangle.
\end{split}
\end{align}
Recall that $\bfu^F$ and $\bfrf^F$ as defined in \eqref{uFR}  and \eqref{deffF}, respectively (see also \eqref{lemnlsetuppr4}),
are given through the identities
\begin{equation} \label{equ:ODE_part_uFR_deffF_recalled}
    \begin{aligned}
    \bfu^F(t) & =  e^{it\bfsM} \bfrf^F(t)   +  e^{-it\bfsM} \overline{\bfrf^F}(t), \\ 
    \frf^{F}_{jk}(t) & = \frac{1}{2\bfsM(\bfsM-2\lambda + i0^+)} e^{-it(\bfsM-2\lambda)} \mathbf{P}_c\bmQ_s\big(Z_j(t)\bmY_j,Z_k(t)\bmY_k\big).
    \end{aligned}
\end{equation}
Recall also from Lemma~\ref{lem:regresbd1} that
\begin{align}\label{eq:resest'}
\begin{split}
& \sum_{q=0,1} \|\jap{r}^{-2}\partial_r^q e^{it\sqrt{\bfM}}(\sqrt{\bfM}-2\lambda+i0^+)^{-1}\bmf\|_{L^2_{r\ud r}}\lesssim \frac{1}{\jap{t}}\|\jap{r}^2\bmf\|_{L^2_{r\ud r}}. 
\end{split}
\end{align}
Let us then record that, using \eqref{eq:resest'} with $t=0$, and the exponential decay of $\bmY$ and its derivatives, 
we see that
\begin{align} \label{estuF}
\begin{split}
& \sum_{q=0,1} {\big\| \langle r\rangle^{-2} \partial_r^q  \bmu^F(t)  \big\|}_{L^2}
    + {\big\| \langle r\rangle^{-2} \partial_r^q e^{-it\bfsM} \overline{\bfrf^F}  \big\|}_{L^2} 
      \\ 
      &\lesssim \sum_{j,\ell=1,2}
     {\Big\| \langle r\rangle^{2} \mathbf{P}_c\bmQ_s\big(\bmY_j,\bmY_\ell\big) \Big\|}_{L^2} 
     \big| Z_j(t)Z_\ell(t) \big| \lesssim \frz_{\Gamma_1}(t). 
\end{split}
\end{align}

We now show the validity of \eqref{IM2estint3}.
From \eqref{IM2est31} we see that $R_{1}^{\mathrm{disc}}$ 
is a linear combination of non-resonant interactions so that we can integrate by parts in time.
More precisely, recalling also \eqref{equ:ODE_part_uFR_deffF_recalled}, we have, for $j=1,2$,
\begin{align}
&\int_{t_1}^{t_2} R_{1}^{\mathrm{disc}}(s) \,\overline{Z_j}(s) \, \ud s  
\\
& = \dfrac{1}{i\lambda} \sum_{k,j_1,j_2=1,2}  
    \int_{t_1}^{t_2} e^{-4i\lambda s} \big\langle \bmQ_s\big(\bmY_k, 
   \frac{1}{2\bfsM(\bfsM-2\lambda + i0^+)} \mathbf{P}_c\bmQ_s
   \big(\bmY_{j_1},\bmY_{j_2}\big)\big) , \bmY_j \big\rangle \,
   \\
   & \qquad \qquad \qquad \qquad \qquad \qquad \qquad \qquad \qquad \times 
   \overline{Z}_{j_1}(s)\overline{Z}_{j_2}(s)\overline{Z}_k(s)\overline{Z_j}(s) \, \ud s.
\end{align}
Using the oscillating factor $e^{-4i\lambda s}$ we perform an integration by parts
in time and, using also \eqref{eq:resest'} with $t=0$ and $\dot{Z}(s) = O(\frz_{\Gamma_1}(s))$, 
we estimate 
\begin{align}\label{IM2est32}
\Big| \sum_{j=1,2}\int_{t_1}^{t_2} \Re \big( R_{1}^{\mathrm{disc}}(s) \,\overline{Z_j}(s) \big) \, \ud s \Big|
 & \lesssim 
  \frz_{\Gamma_1}(t_1)^2 + \int_{t_1}^{t_2} O(\frz_{\Gamma_1}^{3/2}(s)) \cdot O(\dot{Z})(s) \, \ud s
  \lesssim O(\frz_{\Gamma_1}^{3/2}(t_1))
\end{align}
this is stronger than the desired bound. 

The term $R_{2}^{\mathrm{disc}}$ also contains only non-resonant interactions 
$(Z_{j_1}, Z_{j_2}, Z_{k})$ 
and therefore can be handled as above, using integration by parts in the same manner. 

Finally, for the contribution from the term $R_{3}^{\mathrm{disc}}$ we use
the exponential decay of $\bmY$ and \eqref{estuF} to obtain
\begin{align}
\Big| \big\langle \bmQ_s\big(\bfu^F(s), \bfu^F(s) \big), \bmY_j \big\rangle \Big| \lesssim \frz_{\Gamma_1}^2(s)
\end{align}
so that the corresponding integrand in \eqref{IM2estint3} is $O(\frz_{\Gamma_1}^{5/2})$ which is more than sufficient.

\medskip
{\bf Proof of \eqref{IM2estint4}.}
From the definition of $R_{Z,j}^\mathrm{field}$ in \eqref{ODERfield} we write 
(again omitting the dependence on $j$)
\begin{align}
R_{Z,j}^\mathrm{field} & = R_{1}^\mathrm{field} + R_{2}^\mathrm{field} + R_{3}^\mathrm{field},   
\\
\label{ODERfield1}
R_{1}^\mathrm{field} & := \dfrac{1}{i\lambda} e^{-i\lambda t} 
  \sum_{\ell=1,2}  \big\langle \bmQ_s\big(z_\ell(t)\bmY_\ell, \bmu^{R}), \bmY_j \big\rangle, 
   \\
\label{ODERfield2}
R_{2}^\mathrm{field} & :=  \dfrac{1}{i\lambda}  e^{-i\lambda t} \big\langle \bmQ_s\big(\bmu^{F}, \bmu^{R}), \bmY_j \big\rangle, 
   \\
\label{ODERfield3} 
R_{3}^\mathrm{field} & :=   \dfrac{1}{2i\lambda}  e^{-i\lambda t} \big\langle \bmQ_s\big(\bmu^{R}, \bmu^{R}), \bmY_j \big\rangle.
\end{align}
Recall, from Lemma \ref{lemnlsetup}, see \eqref{uFR}, that we have the decomposition
\begin{align}\label{uRvd}
\bmu^{R} = \bmv + \bmu^{R}_d, \qquad \bmv & = e^{it\bfsM}\bfrf^{R}_c + e^{-it\bfsM}\overline{\bfrf^{R}_c},  
\qquad \bmu^{R}_d = e^{it\bfsM}\bfrf^{R}_d + e^{-it\bfsM}\overline{\bfrf^{R}_d},
\end{align}
where the profiles are defined in \eqref{deffRc} and \eqref{deffRd}.
We claim that we have the following bounds, for some absolute $c_0 \ll 1$ that can be chosen sufficiently small:

\begin{itemize}

\item[(i)] For $j,\ell=1,2$, $t_1<t_2 \in [(C\varepsilon)^{-1},T]$,
\begin{subequations}\label{ODERfieldest}
\begin{align}
\label{ODERfield1est1}
\Big| 
  \int_{t_1}^{t_2} e^{-i\lambda s}
  \big\langle \bmQ_s\big(z_\ell(s)\bmY_\ell, \bmv\big), \bmY_j \big\rangle \,\overline{Z_j}(s) 
  \, \ud s \Big| & \leqslant 
  c_0 \, \frz_{\Gamma_1}(t_1), 
\\
\label{ODERfield1est2}
\Big| 
  \int_{t_1}^{t_2} e^{-i\lambda s} 
  \big\langle \bmQ_s\big(z_\ell(s)\bmY_\ell, \bmu^R_d\big), \bmY_j \big\rangle \,\overline{Z_j}(s) \, \ud s \Big| 
  & \leqslant 
  c_0 \, \frz_{\Gamma_1}(t_1);  
\end{align}

\item[(ii)] For all $t\in[(C\varepsilon)^{-1},T]$,
\begin{align}
\label{ODERfield2est}
\big| \big\langle \bmQ_s\big(\bmu^{F}(t), \bmu^{R}(t)\big), \bmY_j \big\rangle \big| 
  & \leqslant  \dfrac{\varepsilon^{0+}
  \frz_{\Gamma_1}^{1/2}(t)}{\langle t\rangle}  
\end{align}

\item[(iii)] For all $t\in[(C\varepsilon)^{-1},T]$,
\begin{align}
\label{ODERfield3est}
\big| \big\langle \bmQ_s\big(\bmu^{R}(t), \bmu^{R}(t)\big), \bmY_j \big\rangle \big| 
  & \leqslant  \dfrac{\varepsilon^{0+}
  \frz_{\Gamma_1}^{1/2}(t)}{\langle t\rangle}    
\end{align}
\end{subequations}

\end{itemize}

The bounds \eqref{ODERfield1est1}-\eqref{ODERfield1est2} give the estimate \eqref{IM2estint4}
with $R_{1}^\mathrm{field}$ instead of $R_{Z,j}^\mathrm{field}$ on the left-hand side,
while the estimate \eqref{ODERfield2est}, respectively \eqref{ODERfield3est},
imply the desired conclusion \eqref{IM2estint4} for the terms $R_{2}^\mathrm{field}$,
respectively $R_{3}^\mathrm{field}$.
To obtain \eqref{IM2estint4} it then suffices to prove the bounds \eqref{ODERfieldest}.

\medskip
{\it Proof \eqref{ODERfield1est1}}.
We first write $z_\ell$ in terms of its profile and split the integral in \eqref{ODERfield1est1} as the sum of the terms
\begin{align}\label{ODERfield1pr1}
I[\bmv]:= \int_{t_1}^{t_2} e^{-2i\lambda s} 
  \big\langle \bmQ(\bmY_\ell, \bmv), \bmY_j \big\rangle \,\overline{Z_\ell(s)}\overline{Z_j}(s) \, \ud s,
\end{align}
and
\begin{align}\label{ODERfield1pr1'}
II[\bmv]:= \int_{t_1}^{t_2}
  \big\langle \bmQ(\bmY_\ell, \bmv), \bmY_j \big\rangle \, Z_\ell(s)\overline{Z_j}(s) \, \ud s. 
\end{align}
We concentrate on estimating the first term, and it will be apparent to the 
reader how a similar argument can handle the easier term $II$ since it has non-vanishing oscillations.

To estimate \eqref{ODERfield1pr1} we first recall 
the flat-sharp decomposition $\bmv(t) =\bfP_c \bmv^\flat(t) + \bfP_c\bmv^\sharp(t)$ of the good part introduced 
in Subsection~\ref{subsec:flat_sharp_decomposition}.
For the contribution of the flat component $\bmv^\flat$ to \eqref{ODERfield},
by  Proposition \ref{prop:dispersive_decay_g} and Lemma~\ref{lem:transferring_bounds_vflat_variables}
\begin{align}\label{equ:definition_bmg_flat_profile22}
\Vert \bmv^\flat \Vert_{L^\infty_r} \lesssim 
\sum_{\ast \in \{\real, \imag\}}\Vert e^{it\langle D\rangle} \bmg_*^\flat(t)\Vert_{L^\infty_x} 
\lesssim \left\{ \begin{aligned}
                            &\varepsilon^{1-6\delta} \jt^{-1+\delta}, \quad  &&0 \leq t \leq \varepsilon^{-2}, \\
                            &t^{-\frac32+4\delta},  &&\varepsilon^{-2} \leq t \leq T.
                         \end{aligned} \right.
\end{align}
Similarly, using the fact that $\Vert \bmg_*^\flat(t)\Vert_{H^N_x(\R^2)}\lesssim \varepsilon^{2-11\delta}$ (see Proposition \ref{prop:HN_g}) along with Gagliardo-Nirenberg-Sobolev, we infer that \[
\Vert \partial_r \bmv^\flat(t)\Vert_{L^\infty_r} \lesssim 
\sum_{\ast \in \{\real, \imag\}} \Vert \nabla e^{it\langle D\rangle} \bmg_*^\flat(t)\Vert_{L^\infty_x} \lesssim \Vert e^{it\langle D\rangle} \bmg_*^\flat(t)\Vert_{L^\infty_x}^{(N-2)/(N-1)} \Vert \bmg_*^\flat(t)\Vert_{H^N_x}^{1/(N-1)}.
\]
Therefore, since $N$ is large, we have
\begin{equation}\label{equ:definition_bmg_flat_profile33}
\begin{aligned}
\Vert \partial_r \bmv^\flat(t)\Vert_{L^\infty_r}
%
%
%
%
%
& \lesssim   \left\{ \begin{aligned}
& \varepsilon^{1-6\delta} \langle t\rangle^{-1+2\delta}  , \quad  &&0 \leq t \leq \varepsilon^{-2}, 
\\ &   \langle t\rangle^{-3/2+5\delta} &&\varepsilon^{-2} \leq t \leq T,
                         \end{aligned} \right.
\end{aligned}\end{equation}
and we can estimate 
\begin{align}\label{ODERfield1pr2}
\begin{split}
\big| I[\bfP_c\bmv^\flat] \big| 
    & \lesssim \int_{t_1}^{t_2} \Big( {\| 
    \bmv^\flat(s) \|}_{L^\infty_r} 
    +  {\| 
    \partial_r 
    \bmv^\flat(s)
    \|}_{L^{\infty}_r} \Big) \, \frz(s) \, \ud s 
    \\
    & \lesssim 
    \int_{t_1}^{t_2} \min\Big\{ \varepsilon^{1-6\delta}\langle s\rangle^{-1+2\delta}, \, \langle s\rangle^{-3/2+5\delta} \Big\} \, \frz_{\Gamma_1}(s) \, \ud s \lesssim \varepsilon^\delta \frz_{\Gamma_1}(t_1).
\end{split}\end{align}
The last inequality holds since, if $t_1<t_2\lesssim \varepsilon^{-2}$, then 
\[
\int_{t_1}^{t_2} \varepsilon^{1-6\delta}\langle s\rangle^{-1+2\delta}\, \frz_{\Gamma_1}(s) \, \ud s \lesssim \varepsilon^{1-10\delta}\frz_{\Gamma_1}(t_1),
\]
whereas, if $t_1 \gtrsim \varepsilon^{-2}$, 
\[
\int_{t_1}^{t_2} \langle s\rangle^{-3/2+5\delta}
    \, \frz_{\Gamma_1}(s) \, \ud s \lesssim \frz_{\Gamma_1}(t_1) \int_{t_1}^{t_2} \langle s\rangle^{-3/2+5\delta} \, \ud s \lesssim \varepsilon^{1-10\delta} \frz_{\Gamma_1}(t_1);
\]
the general case is simply a combination of both cases above.

Next, for the contribution of  $\bfP_c\bmv^\sharp$, we separate the contribution of the initial data and the Duhamel term. Let
\begin{equation}
    \bmv^\sharp_{h}(t):=\cos(t\sqrt{\bfM})\bfP_c\bmv^\sharp(0)+\frac{\sin(t\sqrt{\bfM})}{\sqrt{\bfM}}\bfP_c\partial_t\bmv^\sharp(0),\qquad \bmv^\sharp_p(t):=\int_0^t\frac{\sin((t-s)\sqrt{\bfM})}{\sqrt{\bfM}}\bfP_c(\bfV\bmv^\flat(s))\ud s.
\end{equation}
Then by Proposition~\ref{prop:dispersive_decay_g_sharp}, specifically equations~\eqref{eq:vsharpforFGR1} and~\eqref{eq:vsharpforFGR2} in its proof, 
\begin{equation}\label{eq:bridgebmabound12}
\big\|\bmv^\sharp_p
\big\|_{L^\infty_r} \lesssim
\left\{
\begin{aligned}
&\varepsilon^{1-6\delta}\jt^{-1+3\delta},
&&0\leq t\leq \varepsilon^{-2},
\\
&\varepsilon^{1-12\delta}t^{-1},
&&\varepsilon^{-2}\leq t\leq T.
\end{aligned}
\right.
\end{equation}
Moreover, using  integration by parts to move derivatives away from $\bmv^\sharp$, we infer that
\begin{align}\label{ODERfield1pr3}
\begin{split}
\big|
\big\langle \bmQ_s\big(\bmY_\ell,\bfP_c\bmv_p^\sharp(s)\big),\bmY_j\big\rangle
\big|
&\lesssim \|\bmv^\sharp_p\|_{L^\infty_r}.
\end{split}
\end{align}
Combining \eqref{eq:bridgebmabound12} and \eqref{ODERfield1pr3}, we obtain
\begin{align}\label{ODERfield1pr7}
\begin{split}
\big|I[\bmv_p^\sharp]\big|
&\lesssim
\int_{t_1}^{t_2}
\big|
\big\langle \bmQ_s\big(\bmY_\ell,\bfP_c\bmv_p^\sharp(s)\big),\bmY_j\big\rangle
\big|\,\frz(s)\,\ud s
\\
&\lesssim
\int_{t_1}^{t_2}
\min\big\{ \varepsilon^{1-6\delta}\js^{-1+3\delta} ,
\, \varepsilon^{1-12\delta}s^{-1}\big\} \, \frz(s) \,\ud s
\\
&\leq c_0\,\frz_{\Gamma_1}(t_1),
\end{split}
\end{align}
provided $\varepsilon_0>0$ is chosen sufficiently small. In fact, if
$t_2\leq \varepsilon^{-2}$, then
\[
\int_{t_1}^{t_2}\varepsilon^{1-6\delta}\js^{-1+3\delta}\frz_{\Gamma_1}(s)\,\ud s
\lesssim \varepsilon^{1-12\delta}\frz_{\Gamma_1}(t_1).
\]
If $t_1\geq \varepsilon^{-2}$, then using
$\frz_{\Gamma_1}(s)\lesssim \frz_{\Gamma_1}(t_1)t_1/s$ for $s\geq t_1\geq \varepsilon^{-2}$, we get
\[
\int_{t_1}^{t_2}\varepsilon^{1-12\delta}s^{-1}\frz_{\Gamma_1}(s)\,\ud s
\lesssim \varepsilon^{1-12\delta}\frz_{\Gamma_1}(t_1).
\]
The general case follows by splitting the integral at $\varepsilon^{-2}$.  Next, note that by Proposition~\ref{prop:dispersive_decay_g_sharp},  $\bmv^\sharp_h$ satisfies the same estimate as \eqref{eq:bridgebmabound12} when $t\leq \varepsilon^{-2}$. Therefore the contribution of $\bmv^\sharp_h$ to the $s$ integral in \eqref{ODERfield1est1} for $s\leq \varepsilon^{-2}$ can be treated as above. Therefore, writing $z_\ell=2\Re\big(e^{i\lambda t}Z_\ell\big)$, for \eqref{ODERfield1est1} it suffices to estimate the integral
\begin{equation}
    \int_{t_1}^{t_2} e^{-2i\lambda s}\langle \bmQ_s(\bmY_\ell,\bmv^\sharp_h),\bmY_j\rangle \overline{Z_j}(s)\overline{Z_\ell}(s)\ud s,\qquad t_1\geq \varepsilon^{-2},
\end{equation}
and a similar integral with $e^{-2i\lambda s}\overline{Z_\ell}(s)$ replaced by $Z_\ell(s)$. We ignore this  latter integral because it can be treated more simply in view of the more favorable phase as will be clear from the argument below. We also write
\begin{equation}
    \bmv^\sharp_h(s)=e^{is\sqrt{\bfM}}\bma_++e^{-is\sqrt{\bfM}}\bma_-,\qquad \bma_\pm=\frac{1}{2}\big(\bfP_c\bmv^\sharp(0)\mp i \bfM^{-1/2}\bfP_c\partial_t\bmv^\sharp(0)\big).
\end{equation}
Again  $\bma_-$ contributes a more favorable phase function, so we concentrate on $\bma_+$. For this we write the corresponding integral as
\begin{equation}
    \int_{t_1}^{t_2} \langle \bmQ_s(\bmY_\ell,e^{is(\sqrt{\bfM}-2\lambda)}\bma_+),\bmY_j\rangle \overline{Z_j}(s)\overline{Z_\ell}(s)\ud s =\int_{t_1}^{t_2}(I_h(s)+II_h(s))\overline{Z_j}(s)\overline{Z_\ell}(s)\ud s
\end{equation}
with 
\begin{align}
    &I_h:=\langle \bmQ_s\bigl( \bmY_\ell, \chi_{\leq 10^{-2}}(\sqrt{\bfM-\bfI})e^{is(\sqrt{\bfM}-2\lambda)}\bma_+ \bigr), \bmY_j \bigr\rangle , \\ &II_h:=\langle \bmQ_s\bigl( \bmY_\ell, \chi_{\geq 10^{-2}}(\sqrt{\bfM-\bfI})e^{is(\sqrt{\bfM}-2\lambda)}\bma_+ \bigr), \bmY_j \bigr\rangle.
\end{align}
For $I$ we note that the phase is non-vanishing so we can write
\begin{equation}
    I_h=\frac{\ud}{\ud s}\langle \bmQ\bigl( \bmY_\ell, (i\sqrt{\bfM}-2i\lambda)^{-1}\chi_{\leq 10^{-2}}(\sqrt{\bfM-\bfI})e^{is(\sqrt{\bfM}-2\lambda)}\bma_+ \bigr), \bmY_j \bigr\rangle.
\end{equation}
The contribution of this term can then be handled by integration by parts in $s$ and then estimating the evolution $e^{is\sqrt{\bfM}}\bma_+$ in $L^\infty_{r}$ using \eqref{eq:Linftyweighteddispersivebound1} and estimating $\partial_s\big(\overline{Z_j}(s)\overline{Z_\ell}(s)\big)$ by $s^{-3/2}$ (recall that $s\geq t_1\geq \varepsilon^{-2}$). For $II$ we instead use the local decay estimate from Lemma~\ref{lem:linpaperlocaldecaycutoff1} to bound
\begin{equation}
    \int_{t_1}^{t_2} II_h(s)\overline{Z_j}(s)\overline{Z_\ell}(s)\ud s\lesssim \|\jap{r}\bma_{+}\|_{L^2_{r\ud r}}\int_{t_1}^{t_2}\frac{\ud s}{s^2}\lesssim \|\jap{r}\bma_{+}\|_{L^2_{r\ud r}}\frz_{\Gamma_1}(t_1),
\end{equation}
and the desired smallness is furnished by $\|\jap{r}\bma_{+}\|_{L^2_{r\ud r}}$.

For future reference we also note that combining Proposition~\ref{prop:HN_g_sharp} with
Gagliardo-Nirenberg-Sobolev interpolation,
\begin{equation}\label{vsharpnltouse}
\sum_{|\alpha|\leq 1}
\Big(
\big\| 
\nabla_x^\alpha
\Re\big(e^{i\theta}\bfP_c\bmv^\sharp(t)\big)\big\|_{L^\infty_x(\bbR^2)}
+
\big\| 
\nabla_x^\alpha
\Im\big(e^{i\theta}\bfP_c\bmv^\sharp(t)\big)\big\|_{L^\infty_x(\bbR^2)}
\Big) \lesssim \varepsilon^{0+}\frz_{\Gamma_1}^{1/4}(t)
    \jt^{-1/2}.
\end{equation}

\medskip
{\it Proof of \eqref{ODERfield1est2}}.
Recall the definition 
\begin{align}\label{ufRd}
\bmu^{R}_d & = e^{it\bfsM}\bfrf^{R}_d + e^{-it\bfsM}\overline{\bfrf^{R}_d}, 
\end{align}
where $\bfrf^{R}_d$ is given in \eqref{deffRd} by $\bfrf^{R}_d := \bfrf^{R,1}_d + \bfrf^{R,2}_d$,
and observe that \eqref{ODERfield1est2} is implied by the following inequalities for a small constant $c_0$: 
\begin{align}
\label{ODERfield1est2.1}
& \Big| \int_{t_1}^{t_2} 
  \big\langle \bmQ_s\big(\bmY_\ell, e^{is(\bfsM-2\lambda)}\bfrf^{R,n}_d\big), \bmY_j \big\rangle \, \overline{Z_j}(s)\overline{Z_\ell}(s) \, \ud s \Big| 
  \lesssim c_0 \frz_{\Gamma_1}(t_1), \qquad n=1,2,
\\
\label{ODERfield1est2.1-}
& \Big| \int_{t_1}^{t_2}  
  \big\langle \bmQ_s\big(\bmY_\ell, e^{-is(\bfsM+2\lambda)}\overline{\bfrf^{R,n}_d}\big), \bmY_j \big\rangle \, \overline{Z_j}(s)\overline{Z_\ell}(s) \, \ud s \Big| 
  \lesssim c_0 \frz_{\Gamma_1}(t_1),  \qquad n=1,2,
\\
\label{ODERfield1est2.2}
& \Big| \int_{t_1}^{t_2}  
  \big\langle \bmQ_s\big(\bmY_{\ell'}, e^{is\bfsM}\bfrf^{R,n}_d\big), \bmY_{j'} \big\rangle \,\overline{Z_j}(s) Z_\ell(s)\, \ud s \Big| 
  \lesssim c_0 \frz_{\Gamma_1}(t_1), \qquad n,\ell',j',\ell,j = 1,2.
\end{align}
As in the proof of the estimate \eqref{ODERfield1est1} above, cfr. \eqref{ODERfield1pr1}-\eqref{ODERfield1pr1'},
the estimate \eqref{ODERfield1est2.1} is strictly harder to prove than \eqref{ODERfield1est2.1-}
and \eqref{ODERfield1est2.2}, due to the vanishing of the oscillations $\bfsM-2\lambda$.
We then concentrate on \eqref{ODERfield1est2.1}, and it will be clear to the reader how to adapt the proof 
to obtain \eqref{ODERfield1est2.1-} and \eqref{ODERfield1est2.2}.

\medskip
\noindent
{\it Proof of \eqref{ODERfield1est2.1} for $n=1$}.
To organize our estimates we recall that $\bmN_d = \bmN - \bmN_c$ 
and $\bmN = \bmQ + \bmC + \bmN^{(4)}$, see \eqref{splitN} with \eqref{Nrquad}-\eqref{Nr4}.
According to \eqref{deffRd}, we rewrite $\bfrf^{R,1}_d$ as
\begin{align}\label{fR1ddec}
\bfrf^{R,1}_d = \bmQ^{\mathrm{mix}}_s(\bmu,Z) + \bmQ^{\mathrm{disc}}(Z) + \bmR(\bmr).
\end{align}
Here  $\bmQ^{\mathrm{mix}}_s(\bmu,Z)$, $\bmQ^{\mathrm{disc}}(Z)$ and  $\bmR(\bmr)$ are defined as follows:

\setlength{\leftmargini}{1.5em}
\begin{itemize}
\smallskip
\item The `continuous-discrete' interactions are (with $Z=(Z_1,Z_2)$)
\begin{align}\label{Qmix}
\begin{split}
\bmQ^{\mathrm{mix}}_s(\bmu,Z) & 
  := 2 \sum_{j\in\{1,2\}} \int_0^t e^{-is\bfsM} \dfrac{1}{2i \bfsM}  \bfP_c 
  \bmQ_s \big(\bmu(s), z_j(s) \bfY_j \big) \, \ud s.
\end{split}
\end{align}

\smallskip
\item The `non-resonant discrete' interactions are 
\begin{align}\label{Qdisc2}
\begin{split}
& \bmQ^{\mathrm{disc}}(Z) =  \bmQ^{\mathrm{disc},1}(Z) + \bmQ^{\mathrm{disc},2}(Z),
\\
& \bmQ^{\mathrm{disc},1}(Z) := 2 \sum_{j,\ell\in \{1,2\}} \int_0^t e^{-is \bfsM}  
\dfrac{1}{2i \bfsM}   \bfP_c\bmQ_s\big(\overline{Z}_j(s)\bmY_j, Z_\ell(s)\bmY_\ell\big)  \, \ud s,
  \\ 
& \bmQ^{\mathrm{disc},2}(Z) := \sum_{j,\ell\in \{1,2\}} \int_0^t e^{-is(\bfsM+2\lambda)} 
 \dfrac{1}{2i \bfsM}   \bfP_c\bmQ_s\big(\overline{Z}_j(s)\bmY_j, \overline{Z}_\ell(s)\bmY_\ell\big) \, \ud s.
\end{split}
\end{align}

\smallskip
\item The cubic (and higher order terms) are
\begin{align}
\nonumber
\bmR & = \bmR^{(3)} + \bmR^{(4)}, 
\\
\label{lemnlsetupcontC}
& \bmR^{(3)}(\bmr) := \int_0^t e^{-is \bfsM} \dfrac{1}{2i \bfsM} 
  \big[ \bfP_c \bmC(\bmr(s)) - \bfP_c \bmC(\bfP_c\bmr(s))  \big] \, \ud s,
\\
\label{lemnlsetupcontC+}
& \bmR^{(4)}(\bmr) := \int_0^t e^{-is \bfsM} \dfrac{1}{2i \bfsM} 
  \big[ \bfP_c \bmN^{(4)}(\bmr(s)) - \bfP_c \bmN^{(4)}_c \big] \, \ud s,
\end{align}
\end{itemize}
where $\bmN^{(4)}_c =\bmN_c-\bmQ(\bfP_c\bmr)-\bmC(\bfP_c\bmr)$.

We now proceed to prove the bound \eqref{ODERfield1est2.1} 
for all the contributions in the decomposition \eqref{fR1ddec} above.
First, using Proposition  \ref{prop:dispersive_estimate_with_potential}
and \eqref{ODERestu}, 
 we get
\begin{align}
&  {\big\| e^{is(\bfsM-2\lambda)}\bmQ_s^{\mathrm{mix}}(\bmu,Z)(s) \big\|}_{L^\infty_r} 
\\
& \lesssim {\Big\| \sum_{j=1,2}\int_0^s e^{i(s-\tau)\bfsM} \dfrac{1}{2i \bfsM}  \bfP_c 
  \bmQ_s \big(\bmu(\tau), z_j(\tau) \bfY_j \big) \, \ud \tau \Big\|}_{L^\infty_r} 
  \\
  & \lesssim \sum_{j=1,2} \int_0^s \frac{1}{\langle s -\tau \rangle}
  {\big\| e^{i\theta}  \bmQ_s \big(\bmu(\tau), z_j(\tau) \bfY_j \big) \big\|}_{W^{3,1}_x} \, \ud \tau 
  \\
  & \lesssim \int_0^s \frac{1}{\langle s -\tau \rangle}
  \cdot \frz^{1/2}(\tau) \cdot \frac{\varepsilon^{0+}}{\langle \tau \rangle^{1-C\delta}}
  \, \ud \tau\lesssim \varepsilon^{1/2} \js^{-1},
\end{align}
Integrating by parts to move derivatives only on $\bmY$, it follows that
\begin{align}\label{ODERfield1est2.1pr1}
\begin{split}
\Big| \int_{t_1}^{t_2} 
  \big\langle \bmQ_s\big(\bmY_\ell, e^{is(\bfsM-2\lambda)} \bmQ_s^{\mathrm{mix}}(\bmu,Z)(s) \big) , \bmY_j \big\rangle 
  \, \overline{Z_j}(s)\overline{Z_\ell}(s) \, \ud s \Big| 
  &\lesssim \int_{t_1}^{t_2} \varepsilon^{1/2} \js^{-1} \cdot \frz(s) \, \ud s\\
  &\leqslant c_0 \frz_{\Gamma_1}(t_1),
\end{split}
\end{align}
which gives the desired bound. A similar bound holds for the cubic and higher order remainder terms: using the bootstrap assumption \eqref{equ:bootstrap_assumption_Linfty_f}, the high Sobolev regularity \eqref{equ:bootstrap_assumption_HN_f}, and the bootstrap assumption for the discrete component \eqref{equ:bootstrap_assumption_z}, it is not hard to verify that 
\begin{align}
\begin{split}
& {\Big\| \langle x\rangle\jap{D}^2\big( e^{i\theta} \big[ \bfP_c \bmC(\bmr(s)) - \bfP_c \bmC(\bfP_c\bmr(s)) \big]\big) \Big\|}_{L^2_x}
\\
& + {\Big\| \langle x\rangle\jap{D}^2\big( e^{i\theta}\big[ \bfP_c \bmN^{(4)}(\bmr(s)) - \bfP_c \bmN^{(4)}_c(s) \big]\big) \Big\|}_{L^2_x}  
 \lesssim \varepsilon^{0+} \langle s\rangle^{-1-\delta}.
\end{split}
\end{align}
Therefore, after integration by parts so that only $\bmY$ is differentiated, and using \eqref{eq:Linftydispersiveboundnonsharp2},
\begin{align}
\big\langle \bmQ_s\big(\bmY_\ell, e^{is(\bfsM-2\lambda)}  \bmR(\bmr)(s) \big), \bmY_j \big\rangle
  \lesssim \varepsilon^{0+} \js^{-1},
\end{align}
and hence,
\begin{align}\label{ODERfield1est2.1pr2}
\begin{split}
\Big| \int_{t_1}^{t_2} 
  \big\langle \bmQ_s\big(\bmY_\ell, e^{is(\bfsM-2\lambda)}  \bmR(\bmr)(s) \big), \bmY_j \big\rangle 
  \, \overline{Z_j}(s)\overline{Z_\ell}(s) \, \ud s \Big| 
  \lesssim c_0 \frz_{\Gamma_1}(t_1).
\end{split}
\end{align}

To conclude the proof of \eqref{ODERfield1est2.1} for $n=1$ we need to estimate the 
contribution from $\bmQ^{\mathrm{disc}}(Z)$. The same arguments above do not apply directly,
but we can exploit the oscillations to integrate by parts.
We concentrate on the first term in \eqref{Qdisc2} since the second can be estimated identically.
Integrating by parts we can write 
(we omit the dependence on $j,\ell=1,2$)
\begin{align}\label{Qdisc123}
\begin{split}
\bmQ^{\mathrm{disc},1}(Z) & = 
I + II + III, 
\\ 
I(s) & := \frac{1}{\bfM}e^{-is \bfsM} \bfP_c\bmQ_s\big(\overline{Z}_j(s)\bmY_j, Z_\ell(s)\bmY_\ell\big),
\\
II(s) & := -\frac{1}{\bfM}\bfP_c\bmQ_s\big(\overline{Z}_j(0)\bmY_j, Z_\ell(0)\bmY_\ell\big),
  \\
III(s) & :=- \int_0^s e^{-i\tau \bfsM} \dfrac{1}{\bfM}  
  \frac{d}{d\tau}\bfP_c\bmQ_s\big(\overline{Z}_j(\tau)\bmY_j, Z_\ell(\tau)\bmY_\ell\big)  \, \ud \tau.
\end{split}
\end{align}
Using Proposition  \ref{prop:dispersive_estimate_with_potential}   we obtain
\begin{align}
\begin{split}
\big\langle \bmQ_s\big(\bmY_\ell, e^{is(\bfsM-2\lambda)} \big(II(s)+III(s)\big), \,  \bmY_j \big\rangle  \lesssim \Vert  e^{is(\bfsM-2\lambda)} (II+III)\Vert_{L^\infty_r} \lesssim \varepsilon^{0+} \frz_{\Gamma_1}(s),
\end{split}
\end{align}
(recall that $s \gtrsim \varepsilon_0^{-1}$) and it follows that
\begin{align}
\begin{split}
\Big| \int_{t_1}^{t_2} 
  \big\langle \bmQ_s\big(\bmY_\ell, e^{is(\bfsM-2\lambda)} [II(s) + III(s)] \big), \bmY_j \big\rangle 
  \, \overline{Z_j}(s)\overline{Z_\ell}(s) \, \ud s \Big| 
  \lesssim \Big| \int_{t_1}^{t_2} \varepsilon^{0+} \frz_{\Gamma_1}(s) \cdot \frz(s) \, \ud s \Big| 
  \lesssim c_0 \frz_{\Gamma_1}(t_1).
\end{split}
\end{align}

To handle the contribution from the first term in \eqref{Qdisc123} we first write it out explicitly
and then integrate by parts using the $e^{-2i\lambda s}$ factor:
\begin{align}
\begin{split}
& \Big| \int_{t_1}^{t_2} 
  \big\langle \bmQ_s\big(\bmY_\ell, e^{is(\bfsM-2\lambda)} I(s)\big), \bmY_j \big\rangle 
  \, \overline{Z_j}(s)\overline{Z_\ell}(s) \, \ud s \Big|
  \\
  & = \Big| \int_{t_1}^{t_2} 
  \big\langle \bmQ_s\big(\bmY_\ell, e^{-2i\lambda s} 
  \frac{1}{\bfM} \bfP_c\bmQ_s\big(\overline{Z}_m(s)\bmY_m, Z_k(s)\bmY_k\big)\big), \bmY_j \big\rangle 
  \, \overline{Z_j}(s)\overline{Z_\ell}(s) \, \ud s \Big| 
  \\ & \lesssim  \vert  A(t_1)  \vert + \vert A(t_2)   \vert + B( t_1,t_2),
\end{split}
\end{align}
where
\begin{align}
\begin{split}
& A(s):= \big\langle \bmQ_s\big(\bmY_\ell, e^{-2i\lambda s} 
  \frac{1}{\bfM}\bfP_c\bmQ_s\big(\overline{Z}_m(s)\bmY_m, Z_k(s)\bmY_k\big)\big), \bmY_j \big\rangle 
  \, \overline{Z_j}(s)\overline{Z_\ell}(s), 
\\
& B( t_1, t_2 ):= \Big| \int_{t_1}^{t_2} 
  \big\langle \bmQ_s\big(\bmY_\ell, e^{-2i\lambda s} 
  \frac{1}{\bfM} \, \frac{d}{ds} \Big[ \bfP_c\bmQ_s\big(\overline{Z}_m(s)\bmY_m, Z_k(s)\bmY_k\big)\big), \bmY_j \big\rangle 
  \, \overline{Z_j}(s)\overline{Z_\ell}(s) \Big] \, \ud s \Big|.
\end{split}
\end{align}
We can then verify, using also $\dot{Z}(s) = O(\frz_{\Gamma_1}(s))$, that 
\begin{align*}
& |A(s)| \lesssim \frz_{\Gamma_1}^2(s), \quad |B( t_1, t_2 )| \lesssim \frz_{\Gamma_1}^{3/2}( t_1 ),
\end{align*}
which is more than sufficient to obtain
\begin{align}\label{ODERfield1est2.1pr3}
& \Big| \int_{t_1}^{t_2} 
  \big\langle \bmQ_s\big(\bmY_\ell, e^{is(\bfsM-2\lambda)} \bmQ^{\mathrm{disc},1}(Z)(s) \big), 
  \bmY_j \big\rangle \, \overline{Z_j}(s)\overline{Z_\ell}(s) \, \ud s \Big| 
  \lesssim c_0 \frz_{\Gamma_1}(t_1).
\end{align}
Together with the estimates \eqref{ODERfield1est2.1pr1}, \eqref{ODERfield1est2.1pr2} as well as
the decomposition \eqref{fR1ddec}, we conclude the proof of \eqref{ODERfield1est2.1} for $n=1$.

Note that the arguments above also yield the following bounds,
which will be helpful later:
\begin{align}\label{fR1touse}
\begin{split}
\sum_{\vert \alpha\vert\leq 1} {\| 
\nabla_x^\alpha\big( e^{i\theta} \,  e^{is\bfsM} \bfrf^{R,1}_d(s)\big) \|}_{L^\infty_x} \lesssim \varepsilon^{0+}\js^{-1+C\delta},
\end{split}
\end{align}
and
\begin{align}\label{fR1touse'}
\begin{split}
\sum_{\vert \alpha\vert\leq 1} {\|  
\nabla_x^\alpha\big( e^{i\theta} \, e^{is\bfsM} \bfrf^{R,1}_d(s)\big) \|}_{L^\infty_x} \lesssim     \varepsilon^{0+}\frz_{\Gamma_1}^{1/4}(s)\js^{-1/2}.
\end{split}
\end{align}

\medskip
\noindent
{\it Proof of \eqref{ODERfield1est2.1} for $n=2$}.
Recalling the definition \eqref{deffRd}, that is,
\begin{align}\label{fR2ddec}
\begin{split}
\bfrf^{R,2}_d(t) & := \sum_{j,\ell \in\{1,2\}} \frac{1}{i(\bfsM-2\lambda + i0^+)} \Big[
  \dfrac{1}{2i\bfsM}  \mathbf{P}_c\bmQ_s\big(Z_j(0)\bmY_j,Z_\ell(0)\bmY_\ell\big)  
\\
& \phantom{:= \sum_{j,\ell \in\{1,2\}} \frac{1}{i(\bfsM-2\lambda + i0^+)} \Big[}+  \int_0^t e^{-is(\bfsM-2\lambda)}  
  \dfrac{1}{i\bfsM} \mathbf{P}_c\bmQ_s\big( \dot{Z}_j(s)\bmY_j,Z_\ell(s)\bmY_\ell\big)  \, ds  \Big],
\end{split}
\end{align}
and using \eqref{eq:resest'} and $\dot{Z}_j(s) = O(\frz_{\Gamma_1}(s))$, we see that
\begin{align}\label{fR2ddec2}
\begin{split}
\sum_{q=0,1} {\big\|  \langle r\rangle^{-2} \partial_r^q
e^{is(\bfsM-2\lambda)} \bfrf^{R,2}_d(s) 
\big\|}_{L^2_{r\ud r}} \lesssim \varepsilon^{1-2\delta}\js^{-1}.
\end{split}
\end{align}
It follows that
\begin{align}
& \Big| \int_{t_1}^{t_2} 
  \big\langle \bmQ_s\big(\bmY_\ell, e^{is(\bfsM-2\lambda)}\bfrf^{R,2}_d\big),
  \bmY_j \big\rangle \, \overline{Z_j}(s)\overline{Z_\ell}(s) \, \ud s \Big| 
  \lesssim \Big| \int_{t_1}^{t_2} \varepsilon^{1-2\delta} \js^{-1} \frz(s)\, \ud s \Big| 
  \leqslant c_0 \, \frz_{\Gamma_1}(t_1),
\end{align}
as desired. The proof of \eqref{ODERfield1est2.1} is complete, and \eqref{ODERfield1est2} follows.

\medskip
{\it Proof of \eqref{ODERfield2est}}.
Recall the definitions of $\bmu^{R}$, $\bmu^{R}_d$ and $\bmu^F$ in \eqref{uRvd}, \eqref{ufRd}, and \eqref{uFR} (see also \eqref{deffF}), respectively. Using the estimate \eqref{estuF} on $\bmu^{F}$ we see that
(recall the notation \eqref{deffRd})
\begin{align}
\big| \big\langle \bmQ_s\big(\bmu^{F}(t), \bmu^{R}(t)\big), \bmY_j \big\rangle \big| 
  &  \lesssim \frz_{\Gamma_1}(t) \cdot 
  \Big( {\| \bmu^{R}(t) - \bmu_d^{R,2}(t) \|}_{L^\infty_r}
  + {\| \jap{r}^{-2} \bmu_d^{R,2}(t) \|}_{L^2_{r\ud r}} \Big)
\end{align}
having used again integration by parts to move the derivative away from $\bmu^R$.
To obtain the desired upperbound 
it is then sufficient to show 
\begin{align}
\label{estuR0} 
{\big\| \bmu^{R}(t)-\bmu^{R,2}_d(t) \big\|}_{L^\infty_r} + {\| \jap{r}^{-2} \bmu_d^{R,2}(t) \|}_{L^2_{r\ud r}} \lesssim \varepsilon^{0+}\frz_{\Gamma_1}^{1/4}(t)\jt^{-1/2}, 
\end{align}
which, in view of \eqref{uRvd}, is implied by 
\begin{align}\label{estuR} 
\begin{split}
{\|  
\bmv(t) 
\|}_{L^\infty}
+ 
{\| 
\bmu^{R,1}_d(t)
\|}_{L^\infty_r} + {\| \jap{r}^{-2} \bmu_d^{R,2}(t) \|}_{L^2_{r\ud r}}
\lesssim \varepsilon^{0+}\frz_{\Gamma_1}^{1/4}(t)\jt^{-1/2}. 
\end{split}
\end{align}
The estimate \eqref{estuR} follows by combining:

\noindent
- \eqref{equ:definition_bmg_flat_profile22}-\eqref{equ:definition_bmg_flat_profile33} to bound $\bmv^\flat$, 

\noindent
- the estimate \eqref{eq:bridgebmabound12}-\eqref{vsharpnltouse} to control $\bmv^\sharp$,  

\noindent
- the estimates  \eqref{fR1touse'} and \eqref{fR2ddec2} to bound $e^{it\bfsM}\bfrf_d^{R,1}$ and $e^{it\bfsM}\bfrf_d^{R,2}$, respectively.

\medskip
{\it Proof of \eqref{ODERfield3est}}.
This case follows by noticing that the same estimates listed after \eqref{estuR} that we used to prove \eqref{estuR0}, also imply that \begin{align*}
     \sum_{\vert \alpha\vert\leq1} {\Big\| \nabla_x^\alpha\Big( e^{i\theta} \big( \bmu^{R}(t)-\bmu^{R,2}_d(t) \big)\Big) \Big\|}_{L^\infty_x} + \sum_{q=0,1} {\| \jap{r}^{-2} \partial_r^q\bmu_d^{R,2}(t) \|}_{L^2_{rdr}} \lesssim \varepsilon^{0+}\frz_{\Gamma_1}^{1/4}(t)\jt^{-1/2}.
\end{align*}
Then,  
\begin{align}
& \big| \big\langle \bmQ_s\big(\bmu^{R}(t), \bmu^{R}(t)\big), \bmY_j \big\rangle \big|  
  \\ & \lesssim\Big(  \sum_{\vert \alpha\vert\leq1} {\Big\| \nabla_x^\alpha\Big( e^{i\theta} \big( \bmu^{R}(t)-\bmu^{R,2}_d(t) \big)\Big) \Big\|}_{L^\infty_x} + \sum_{q=0,1} {\| \jap{r}^{-2} \partial_r^q\bmu_d^{R,2}(t) \|}_{L^2_{rdr}} \Big)^2   
  \lesssim \dfrac{\varepsilon^{0+}\frz_{\Gamma_1}^{1/2}(t)}{\langle t\rangle}.
\end{align}
The proofs of the inequalities \eqref{ODERfieldest} are complete,
and \eqref{IM2estint4} follows.

\medskip
{\bf Proof of \eqref{IM2estimatesC}.}
%
According to \eqref{Nrcub}, and $\bmr(t) = \bmu(t) + z_1 \bmY_1 + z_2 \bmY_2$, 
the quantity on the left-hand side of \eqref{IM2estC}
is a (localized) cubic expression in $\bmu$ and $z_j \bmY_j$ with at least one factor of $\bmu$
and at most one derivative on one of the three inputs.
Using \eqref{ODERestu}, it follows that
\begin{align*}
& \Big| \big\langle \bmC(\bmr(t)) - \bmC(\bmr(t)-\bmu(t)) , \bmY_j \big\rangle \Big| 
  \\
  & \lesssim  \sum_{\vert \alpha\vert, \vert \beta\vert \leq1} {\big\|  \nabla_x^\alpha \partial_t^\beta \big(e^{i\theta} \bmu(t)\big) \big\|}_{L^\infty_x}
  \Big( \sum_{\vert \alpha\vert, \vert \beta\vert \leq1} 
  {\big\|  \nabla_x^\alpha \partial_t^\beta \big(e^{i\theta} \bmu(t)\big) \big\|}_{L^\infty_x}^2 + \frz(t) \Big)
  \lesssim \frz_{\Gamma_1}^{5/3}(t)
\end{align*}
for all times $t\geqslant (C\varepsilon)^{-1}$.
Similarly, for the quartic and higher order terms on the left-hand side of \eqref{IM2est4} we have
\begin{align*}
\Big| \big\langle \bmN^{(4)}(\bmr), \bfY_j \big\rangle \Big| 
  \lesssim \big( {\|  \bmu(t) \|}_{L^\infty_r} +  \big{\| \nabla_x \big(e^{i\theta} \bmu(t)\big) \big\|}_{L^\infty_x}
  + |Z_1(t)| + |Z_2(t)| \big)^4 \lesssim \frz_{\Gamma_1}^{2}(t)
\end{align*}
The proof of Proposition \ref{propODERbounds} is complete.
\end{proof}


\section{Conclusion of the Proof of Theorem~\ref{thm:main}} \label{sec:conclusion_proof_thm_main}

In this final section we detail how the results from the preceding sections combine to furnish a proof of Theorem~\ref{thm:main}.

\begin{proof}[Proof of Theorem~\ref{thm:main}]
We begin by explaining how the local and global existence result for the abelian Yang-Mills-Higgs equation in the Lorenz gauge from \cite{BurMon1} applies in the present gauge to furnish a local-in-time solution to \eqref{equ:EL_equations_selfdual} at $H^N_x$-regularity on some maximal interval of existence $[0,T_\ast)$. \cite{BurMon1} establishes that for initial data
\begin{equation}
    \bigl( \eta_0,\partial_t\eta_0,\eta_{x^1},\partial_t\eta_{x^2},\eta_{x^2},\partial_t\eta_{x^2},\varphi,\bfD_0\varphi \big) \big\vert_{t=0}\in \mathscr{H}^{(s)} := \big(H^{s+1}_x(\mathbb{R}^2)\times H^{s}_x(\mathbb{R}^2)\big)^4, \quad s \geq 2,
\end{equation}
satisfying the Lorenz gauge condition $\nabla^\mu \eta_\mu=0$, the abelian Yang-Mills-Higgs system of equations \eqref{equ:EL_equations_selfdual} with $\phi=\underline{\Phi}+\varphi$ and $A=a+\eta$, admits a unique solution defined globally in time.
Note that in the Lorenz gauge $A_\mu$, $\mu=0,1,2,3$, and $\phi$ all satisfy wave-type equations (see for instance \cite[equation (2.19)]{BurMon1}), and the initial data are subject to the constraint $\Delta A_0=\Im(\overline{\phi}\bfD_0\phi)$ (which in our case is just the elliptic equation \eqref{equ:eta0_equation} for $\eta_0$; see also \cite[equation (2.11)]{BurMon1}).
Now given sufficiently regular data in a different gauge, such as Stuart's gauge \eqref{equ:intro_orthogonal_gauge_condition} in our case, one can obtain data satisfying the Lorenz gauge condition by means of a gauge transformation. Indeed, suppose we are given initial data satisfying Stuart's gauge condition \eqref{equ:intro_orthogonal_gauge_condition}. These data determine $A_0$ 
by means of the elliptic equation \eqref{equ:eta0_equation}. 
Then the modified initial data
\begin{align*}
    &\mathring\phi(0)=\phi(0), \quad\mathring{A}_\mu(0)=A_\mu(0),\\
    &\partial_t\mathring\phi(0)=\partial_t\phi(0),\quad \partial_t\mathring{A_j}(0)=\partial_tA_j(0),\quad \partial_t\mathring{A_0}(0)=\nabla^jA_j(0),
\end{align*}
satisfy the Lorenz gauge condition $\nabla^\mu \mathring{A}_\mu = 0$ as well as the constraint \eqref{equ:eta0_equation} at time $t=0$. It follows from \cite{BurMon1} that there is a unique global evolution of this data satisfying equations \cite[equation (2.19)]{BurMon1} in Lorenz gauge\footnote{Incidentally, the modified data still satisfy Stuart's gauge condition at $t=0$ as well as \eqref{equ:eta0_equation}. But since the evolution equations are not the ones derived from \eqref{equ:EL_equations_selfdual} according to Stuart's gauge, this condition will not necessarily be satisfied for later times.} and the constraint \cite[equation (2.11)]{BurMon1}. Moreover, the resulting solution is equivariant if the original (and hence the modified) data are equivariant. To return to a solution in Stuart's gauge, we decompose $\mathring\phi$ as $\mathring\phi=\underline{\Phi}+\psi$. If $\chi$ satisfies the elliptic equation (see \cite[p.73]{Stu})
\begin{equation}\label{eq:LoretnztoStuartchi1}
    -\Delta\chi+|\underline{\Phi}|^2\sin\chi+\Im\big(e^{i\chi}\psi\overline{\underline{\Phi}}\big)-\nabla^j\mathring{A}_j=0,
\end{equation}
then $(e^{i\chi}\mathring{\phi},\mathring{A}+\ud\chi)$ are a solution to \eqref{equ:EL_equations_selfdual} in Stuart's gauge. To see that \eqref{eq:LoretnztoStuartchi1} has a solution for short time, observe that if $\varphi$ and $\eta$ are small, say in $H^N_x(\bbR^2)$, then $\nabla^j\mathring{A}_j$ and $\psi$ remain small in $H^N_x(\bbR^2)$ for small time. Existence of a solution to \eqref{eq:LoretnztoStuartchi1} then follows from standard perturbative methods. We also need to verify that the initial data for $(e^{i\chi}\mathring{\phi},\mathring{A}+\ud \chi)$ agree with our original data. This follows if $\chi(0)=\partial_t\chi(0)=0$, which in turn follows from evaluating  \eqref{eq:LoretnztoStuartchi1} and its time derivative at time $t=0$, and comparing the initial constraints \eqref{equ:eta0_equation} and \cite[equation (2.11)]{BurMon1}. 

Having thus obtained a local-in-time solution\footnote{Local existence for data close to $(\underline{\Phi},a)$ can also be proved directly in Stuart's gauge by setting up an iteration as in \cite[Theorem~5.1]{Stu}. However, an additional argument as in \cite{BurMon1} would be required to extend this solution globally.}, we can start our bootstrap \eqref{equ:bootstrap_assumption_z}--\eqref{equ:bootstrap_assumption_HN_ILED_f} on some interval $[0,T]$. 
Once these bootstrap assumptions have been improved, the resulting uniform bounds allow us to continue the solution and obtain global existence. 
It therefore remains to explain how the estimates established in the preceding sections improve \eqref{equ:bootstrap_assumption_z} and \eqref{equ:bootstrap_assumption_weighted_f}--\eqref{equ:bootstrap_assumption_HN_ILED_f}. This closes the bootstrap and proves part (i) of Theorem~\ref{thm:main}. Under the additional assumption \eqref{equ:new_assumption_frakz_equal} of part (ii), no further bootstrap argument is needed: the lower bound in \eqref{equ:decay_zed_lower_bound_mainthm} will follow directly from part (ii) of Proposition~\ref{propODEdamp}.
We also note that the pointwise decay estimate \eqref{eq:decay_radiation_mainthm} for the radiation follows from the improved dispersive bound \eqref{equ:bootstrap_assumption_Linfty_f}, while for very short times it follows by Sobolev embedding from the high Sobolev bound \eqref{equ:bootstrap_assumption_HN_f}.


First, the bootstrap bound \eqref{equ:bootstrap_assumption_z} about the decay estimate for the internal mode components is improved as a consequence of Propositions~\ref{propODEdamp} and~\ref{propODERbounds}. More precisely, under our bootstrap assumptions, Proposition~\ref{propODERbounds} implies that the hypotheses \eqref{propODEshort}, \eqref{propODEdampR1a}, \eqref{propODEdampR1b}, \eqref{propODEdampR1h} of Proposition~\ref{propODEdamp} are satisfied. The latter then implies that the bootstrap bound \eqref{equ:bootstrap_assumption_z} can be improved.
Moreover, part (ii) of Proposition~\ref{propODEdamp} shows that, under the additional assumption \eqref{equ:new_assumption_frakz_equal} in part (ii) of the statement of Theorem~\ref{thm:main}, the lower bound in \eqref{equ:decay_zed_lower_bound_mainthm} for the decay of the internal mode components follows directly.

Next we turn to the remaining bootstrap bounds \eqref{equ:bootstrap_assumption_weighted_f}, \eqref{equ:bootstrap_assumption_Linfty_f}, \eqref{equ:bootstrap_assumption_HN_f}, \eqref{equ:bootstrap_assumption_HN_ILED_f} concerning the profiles $\bmf_\real(t)$ and $\bmf_\imag(t)$. We recall that the decomposition of the radiation term
\begin{equation}
    \bmu(t) = \bfP_c \bmu(t) = \bfP_c \bmv^\flat(t) + \bfP_c \bmv^\sharp(t) + \bfP_c \bmw(t)
\end{equation}    
implies the decompositions of the profiles
\begin{equation} \label{equ:proof_theorem_decomposition_profiles}
    \begin{aligned}
        \bmf_\real(t) &= e^{-it\jD} (2i\jD)^{-1} (\pt + i\jD) \, \Re\bigl( e^{i\theta} \bfP_c \bmv^\flat(t) \bigr) + \bmg^\sharp_\real(t) + \bmh_\real(t), \\
        \bmf_\imag(t) &= e^{-it\jD} (2i\jD)^{-1} (\pt + i\jD) \, \Im\bigl( e^{i\theta} \bfP_c \bmv^\flat(t) \bigr) + \bmg^\sharp_\imag(t) + \bmh_\imag(t).
    \end{aligned}
\end{equation} 
It suffices to consider $\bmf_\real(t)$.
For what follows, recall that $C_1 \geq 1$ denotes an absolute constant, while $C_2 \equiv C_2(C_0) \geq 1$ denotes a constant whose size depends on the size of the constant $C_0$ in the bootstrap assumptions.
We begin with the bootstrap assumption \eqref{equ:bootstrap_assumption_weighted_f} about the weighted energy of $\bmf_\real(t)$. The decomposition \eqref{equ:proof_theorem_decomposition_profiles} implies
\begin{equation} \label{equ:proof_theorem_weighted_energy_decomp}
    \begin{aligned}
        \bigl\| \jxi^2 \jap{\nabla_\xi} \widehat{\bmf}_\real(t,\xi) \bigr\|_{L^2_\xi} &\leq \Bigl\| \jD^2 \jx \Bigl( e^{-it\jD} (2i\jD)^{-1} (\pt + i\jD) \Re \bigl( e^{i\theta} \bfP_c \bmv^\flat(t) \bigr) \Bigr) \Bigr\|_{L^2_x} \\ 
        &\quad + \bigl\| \jxi^2 \jap{\nabla_\xi} \widehat{\bmg}^\sharp_\real(t,\xi) \bigr\|_{L^2_\xi} + \bigl\| \jxi^2 \jap{\nabla_\xi} \widehat{\bmh}_\real(t,\xi) \bigr\|_{L^2_\xi}.
    \end{aligned}
\end{equation}
In the first term on the right-hand side of \eqref{equ:proof_theorem_weighted_energy_decomp} we insert the definition \eqref{equ:definition_bfP_c} of the projection to the continuous spectral subspace. Using the transfer estimate \eqref{equ:transferring_bounds_vflat_weighted_energy} as well as the weighted energy bounds from Proposition~\ref{prop:weighted_energies_g}, the dispersive decay estimate from Proposition~\ref{prop:dispersive_decay_g}, the Sobolev norm bound from Proposition~\ref{prop:HN_g} along with the interpolation inequality~\eqref{equ:GNS}, this gives uniformly for all times $0 \leq t \leq T$,
\begin{equation}
    \begin{aligned}
        &\Bigl\| \jD^2 \jx \Bigl( e^{-it\jD} (2i\jD)^{-1} (\pt + i\jD) \Re \bigl( e^{i\theta} \bfP_c \bmv^\flat(t) \bigr) \Bigr) \Bigr\|_{L^2_x} \\ 
        &\leq \Bigl\| \jD^2 \jx \Bigl( e^{-it\jD} (2i\jD)^{-1} (\pt + i\jD) \Re \bigl( e^{i\theta} \bmv^\flat(t) \bigr) \Bigr) \Bigr\|_{L^2_x} \\ 
        &\quad + \sum_{j=1,2} \Bigl\| \jD^2 \jx \Bigl( e^{-it\jD} (2i\jD)^{-1} (\pt + i\jD) \Re \bigl( \langle \bmY_j, \bmv^\flat(t) \rangle e^{i\theta} \bmY_j \bigr) \Bigr) \Bigr\|_{L^2_x} \\ 
        &\lesssim \sum_{\ast\in\{\real,\imag\}} \, \bigl\| \jD^2 \jx \bmg_\ast^\flat(t) \bigr\|_{L^2_x} + \sum_{j=1,2} \Bigl( \bigl\|\bmv^\flat(t)\bigr\|_{L^\infty_r} + \bigl\| \pt \bmv^\flat(t) \bigr\|_{L^\infty_r} \Bigr) \Bigl\| \jD^2 \jx \Bigl( e^{\pm i t \jD} \bigl( e^{i\theta} \bmY_j \bigr) \Bigr) \Bigr\|_{L^2_x} \\ 
        &\lesssim \sum_{\ast\in\{\real,\imag\}} \Bigl\| \jxi^2 \jap{\nabla_\xi} \Bigl( \widehat{\bmg}^\flat_{\ast}(t,\xi) - \widehat{\calF}\bigl[ \bmB_\ast[\bmf, \bmf](t) \bigr](\xi) \Bigr) \Bigr\|_{L^2_\xi(\bbR^2)} \\
        &\quad + \sum_{\ast\in\{\real,\imag\}} \Bigl\| \jxi^2 \jap{\nabla_\xi} \widehat{\calF}\bigl[ \bmB_\ast[\bmf, \bmf](t) \bigr](\xi) \Bigr\|_{L^2_\xi(\bbR^2)} + \sum_{\ast\in\{\real,\imag\}} \, \bigl\| \jD e^{it\jD} \bmg_\ast^\flat(t) \bigr\|_{L^\infty_x} \cdot \jt \\
        &\leq C_2 \left\{ \begin{aligned}
                            &\varepsilon^{1-10\delta} \jt^{\frac12-3\delta}, \quad &&0 \leq t \leq \varepsilon^{-2}, \\
                            &t^{2\delta},  &&\varepsilon^{-2} \leq t \leq T.
                         \end{aligned} \right. 
    \end{aligned}
\end{equation}
Inserting the preceding bound back into \eqref{equ:proof_theorem_weighted_energy_decomp} along with the weighted energy bound for $\widehat{\bmg}^\sharp_\real(t,\xi)$ from Proposition~\ref{prop:weighted_energies_g_sharp} and for $\widehat{\bmh}_\real(t,\xi)$ from Proposition~\ref{prop:weighted_energies_h}, it follows that by choosing the bootstrap constant $C_0 \geq 1$ sufficiently large depending only on the size of the absolute constants, and by then choosing the small constant $0 < \varepsilon_0 \ll 1$ sufficiently small depending on the size of~$C_0$, we can conclude that
\begin{equation}
        \bigl\| \jxi^2 \jap{\nabla_\xi} \widehat{\bmf}_\ast(t,\xi) \bigr\|_{L^2_\xi(\bbR^2)} \leq C_0 \left\{ \begin{aligned}
                            &\varepsilon^{1-20\delta}, &&0 \leq t \leq 1, \\
                            &\varepsilon^{1-20\delta} t^{\frac12}, &&1 \leq t \leq \varepsilon^{-1-20\delta}, \\
                            &\varepsilon^2 t^{\frac32},  &&\varepsilon^{-1-20\delta} \leq t \leq \varepsilon^{-2}, \\
                            &t^{\frac12},  &&\varepsilon^{-2} \leq t \leq T,
                         \end{aligned} \right.
\end{equation}
thus improving the bootstrap assumption \eqref{equ:bootstrap_assumption_weighted_f}.

Next, we consider the bootstrap assumption \eqref{equ:bootstrap_assumption_Linfty_f} about the dispersive decay of the evolution of the profile $e^{it\jD} \bmf_\real(t)$. 
Here we have
\begin{equation} \label{equ:proof_theorem_dispersive_decay_decomp}
    \begin{aligned}
        \bigl\| e^{it\jD} \bmf_\real(t) \bigr\|_{L^\infty_x} &\leq \bigl\| (2i\jD)^{-1} (\pt + i\jD) \Re \bigl( e^{i\theta} \bfP_c \bmv^\flat(t) \bigr) \bigr\|_{L^\infty_x} \\ 
        &\quad + \bigl\| e^{it\jD} \bmg^\sharp_\real(t) \bigr\|_{L^\infty_x} + \bigl\| e^{it\jD} \bmh_\real(t) \bigr\|_{L^\infty_x}.
    \end{aligned}
\end{equation}
For the first term on the right-hand side of \eqref{equ:proof_theorem_dispersive_decay_decomp} we obtain using the transfer estimate \eqref{equ:transferring_bounds_vflat_Lp} along with the dispersive decay estimates from Proposition~\ref{prop:dispersive_decay_g}, the Sobolev bound from Proposition~\ref{prop:HN_g}, and the interpolation inequality \eqref{equ:GNS},
\begin{equation}
    \begin{aligned}
        &\bigl\| (2i\jD)^{-1} (\pt + i\jD) \Re \bigl( e^{i\theta} \bfP_c \bmv^\flat(t) \bigr) \bigr\|_{L^\infty_x} \\
        &\leq \bigl\| (2i\jD)^{-1} (\pt + i\jD) \Re \bigl( e^{i\theta} \bmv^\flat(t) \bigr) \bigr\|_{L^\infty_x} 
        + \sum_{j=1,2} \bigl\| (2i\jD)^{-1} (\pt + i\jD) \Re \bigl( \langle \bmY_j, \bmv^\flat(t) \rangle e^{i\theta} \bmY_j \bigr) \bigr\|_{L^\infty_x} \\ 
        &\lesssim \bigl\| \bmv^\flat(t) \bigr\|_{L^\infty_r} + \bigl\| \partial_t \bmv^\flat(t) \bigr\|_{L^\infty_r} \lesssim \sum_{\ast \in \{\real,\imag\}} \bigl\| \jD e^{it\jD} \bmg^\flat_\ast(t) \bigr\|_{L^\infty_x} \leq C_2 \left\{ \begin{aligned}
                            &\varepsilon^{1-6\delta} \jt^{-1+2\delta}, \quad &&0 \leq t \leq \varepsilon^{-2}, \\    
                            &t^{-\frac32+5\delta},  &&\varepsilon^{-2} \leq t \leq T.
                         \end{aligned} \right.
    \end{aligned}
\end{equation}
Inserting the preceding bound for the first term on the right-hand side of \eqref{equ:proof_theorem_dispersive_decay_decomp} along with the dispersive decay estimates for $e^{it\jD} \bmg^\sharp_\real(t)$ from Proposition~\ref{prop:dispersive_decay_g_sharp} and for $e^{it\jD} \bmh_\real(t)$ from Proposition~\ref{prop:dispersive_decay_h}, we infer, arguing as above, that 
    \begin{equation}
        \bigl\| e^{it\jD} \bmf_\ast(t)\bigr\|_{L^\infty_x(\bbR^2)} \leq C_0 \left\{ \begin{aligned}
                            &\varepsilon^{1-10\delta}, &&0 \leq t \leq 1, \\    
                            &\varepsilon^{1-10\delta} t^{-\frac12-4\delta},  &&1 \leq t \leq \varepsilon^{-2}, \\
                            &t^{-1+\delta},  &&\varepsilon^{-2} \leq t \leq T,
                         \end{aligned} \right.  
    \end{equation}
thus improving the bootstrap assumption \eqref{equ:bootstrap_assumption_Linfty_f}.

Now we turn to the bootstrap assumption \eqref{equ:bootstrap_assumption_HN_f} about the high Sobolev norm of $\bmf_\real(t)$.
From \eqref{equ:proof_theorem_decomposition_profiles} we obtain
\begin{equation} \label{equ:proof_theorem_HN_decomp}
    \begin{aligned}
        \bigl\| \bmf_\real(t) \bigr\|_{H^N_x} &\leq \bigl\| (2i\jD)^{-1} (\pt + i\jD) \Re\bigl( e^{i\theta} \bfP_c \bmv^\flat(t) \bigr) \bigr\|_{H^N_x} + \bigl\| \bmg_\real^\sharp(t) \bigr\|_{H^N_x} + \bigl\| \bmh_\real(t) \bigr\|_{H^N_x}.
    \end{aligned}
\end{equation}
Then for the first term on the right-hand side of \eqref{equ:proof_theorem_HN_decomp} we use the transfer estimate \eqref{equ:transferring_bounds_vflat_HN} along with the high Sobolev norm bound from Proposition~\ref{prop:HN_g} to conclude  
\begin{equation}
    \begin{aligned}
        &\bigl\| (2i\jD)^{-1} (\pt + i\jD) \Re\bigl( e^{i\theta} \bfP_c \bmv^\flat(t) \bigr) \bigr\|_{H^N_x} \\
        &\leq \bigl\| (2i\jD)^{-1} (\pt + i\jD) \Re\bigl( e^{i\theta} \bmv^\flat(t) \bigr) \bigr\|_{H^N_x} + \sum_{j=1,2} \bigl\| (2i\jD)^{-1} (\pt + i\jD) \Re\bigl( \langle \bmY_j, \bmv^\flat(t) \rangle  e^{i\theta} \bmY_j \bigr) \bigr\|_{H^N_x} \\ 
        &\lesssim \sum_{\ast \in \{\real,\imag\}} \bigl\| (2i\jD)^{-1} (\pt + i\jD) \bmv_\ast^\flat(t) \bigr\|_{H^N_x} + \bigl\| \bmv^\flat(t) \bigr\|_{L^2_{r\ud r}} + \bigl\| \pt \bmv^\flat(t) \bigr\|_{L^2_{r\ud r}} \\ 
        &\lesssim \sum_{\ast \in \{\real,\imag\}} \bigl\| \bmg_\ast^\flat(t) \bigr\|_{H^N_x} \leq C_2 \varepsilon^{\frac32 - 2\delta}.
    \end{aligned}
\end{equation}
Inserting this bound for the first term on the right-hand side of \eqref{equ:proof_theorem_HN_decomp} along with the high Sobolev norm bound for $\bmg_\real^\sharp(t)$ from Proposition~\ref{prop:HN_g_sharp} and for $\bmh_\real(t)$ from Proposition~\ref{prop:HN_h}, we conclude, arguing as before, that
\begin{equation}
    \bigl\| \bmf_\real(t) \bigr\|_{H^N_x} \leq C_0 \varepsilon,
\end{equation}
whence improving the bootstrap assumption \eqref{equ:bootstrap_assumption_HN_f}.

It remains to consider the bootstrap assumption \eqref{equ:bootstrap_assumption_HN_ILED_f} about the ILED estimate at high Sobolev norm for the evolution of the profile $e^{it\jD} \bmf_\real(t)$. From \eqref{equ:proof_theorem_decomposition_profiles} we have 
\begin{equation} \label{equ:proof_theorem_ILED_decomp}
    \begin{aligned}
        &\sum_{1 \leq k \leq N} \bigl\| \jx^{-\frac32-\kappa} |D|^k e^{it\jD} \bmf_\real(t) \bigr\|_{L^2_t([0,T]; L^2_x)} \\
        &\leq \sum_{1 \leq k \leq N} \bigl\| \jx^{-\frac32-\kappa} |D|^k (2i\jD)^{-1} (\pt + i\jD) \Re\bigl( e^{i\theta} \bfP_c \bmv^\flat(t) \bigr) \bigr\|_{L^2_t([0,T]; L^2_x)} \\ 
        &\quad  + \sum_{1 \leq k \leq N} \bigl\| \jx^{-\frac32-\kappa} |D|^k e^{it\jD} \bmg^\sharp_\real(t) \bigr\|_{L^2_t([0,T]; L^2_x)} 
        \\
        &\quad + \sum_{1 \leq k \leq N} \bigl\| \jx^{-\frac32-\kappa} |D|^k e^{it\jD} \bmh_\real(t) \bigr\|_{L^2_t([0,T]; L^2_x)}.
    \end{aligned}
\end{equation}
For the first term on the right-hand side of \eqref{equ:proof_theorem_ILED_decomp} we obtain using the transfer estimate \eqref{equ:transferring_bounds_vflat_ILED}, the interpolation inequality \eqref{equ:GNS}, and the bounds from Proposition~\ref{prop:dispersive_decay_g}, Proposition~\ref{prop:HN_g}, as well as Proposition~\ref{prop:HN_ILED_g},
\begin{equation}
    \begin{aligned}
        &\sum_{1 \leq k \leq N} \bigl\| \jx^{-\frac32-\kappa} |D|^k (2i\jD)^{-1} (\pt + i\jD) \Re\bigl( e^{i\theta} \bfP_c \bmv^\flat(t) \bigr) \bigr\|_{L^2_t([0,T]; L^2_x)} \\ 
        &\leq \sum_{1 \leq k \leq N} \bigl\| \jx^{-\frac32-\kappa} |D|^k (2i\jD)^{-1} (\pt + i\jD) \Re\bigl( e^{i\theta} \bmv^\flat(t) \bigr) \bigr\|_{L^2_t([0,T]; L^2_x)} \\ 
        &\quad + \sum_{j=1,2} \sum_{1 \leq k \leq N} \bigl\| \jx^{-\frac32-\kappa} |D|^k (2i\jD)^{-1} (\pt + i\jD) \Re\bigl( \langle \bmY_j, \bmv^\flat(t) \rangle e^{i\theta} \bmY_j \bigr) \bigr\|_{L^2_t([0,T]; L^2_x)} \\
        &\lesssim \sum_{1 \leq k \leq N} \bigl\| \jx^{-\frac32-\kappa} |D|^k e^{it\jD} \bmg^\flat_\real(t) \bigr\|_{L^2_t([0,T]; L^2_x)} + \sum_{j=1,2} \Bigl( \bigl\| \langle \bmY_j, \bmv^\flat(t) \rangle \bigr\|_{L^2_t([0,T])} + \bigl\| \langle \bmY_j, \pt \bmv^\flat(t) \rangle \bigr\|_{L^2_t([0,T])} \Bigr) \\ 
        &\lesssim \sum_{1 \leq k \leq N} \bigl\| \jx^{-\frac12-\kappa} |D|^k e^{it\jD} \bmg^\flat_\real(t) \bigr\|_{L^2_t([0,T]; L^2_x)} + \sum_{l=0,1} \Bigl\| \bigl\| \pt^l \bmv^\flat(t)\bigr\|_{L^2_{r\ud r}}^{\frac12} \bigl\| \pt^l \bmv^\flat(t)\bigr\|_{L^\infty_{r}}^{\frac12} \Bigr\|_{L^2_t([0,T])}  \\ 
        &\lesssim \sum_{1 \leq k \leq N} \bigl\| \jx^{-\frac12-\kappa} |D|^k e^{it\jD} \bmg^\flat_\real(t) \bigr\|_{L^2_t([0,T]; L^2_x)} + \sum_{\ast \in \{\real,\imag\}} \Bigl\| \bigl\| \jD e^{it\jD} \bmg_\ast^\flat(t) \bigr\|_{L^2_x}^{\frac12} \bigl\| \jD e^{it\jD} \bmg_\ast^\flat(t) \bigr\|_{L^\infty_x}^{\frac12} \Bigr\|_{L^2_t([0,T])} \\ 
        &\leq C_2 \varepsilon^{\frac32 - 2\delta} + C_2 \varepsilon^{\frac54-5\delta}.
    \end{aligned}
\end{equation}
Together with the ILED estimates for $e^{it\jD} \bmg^\sharp_\real(t)$ from Proposition~\ref{prop:HN_g_sharp} and for $e^{it\jD} \bmh_\real(t)$ from Proposition~\ref{prop:HN_h}, we conclude
\begin{equation}
    \sum_{1 \leq k \leq N} \bigl\| \jx^{-\frac32-\kappa} |D|^k e^{it\jD} \bmf_\real(t) \bigr\|_{L^2_t([0,T]; L^2_x)} \leq C_0 \varepsilon,
\end{equation}
thus improving the bootstrap assumption \eqref{equ:bootstrap_assumption_HN_ILED_f}.
This finishes the proof of Theorem~\ref{thm:main}.
\end{proof}


\bigskip

\begin{refcontext}[sorting=nyt]
  \printbibliography
\end{refcontext}

\end{document}